\documentclass[10pt,twoside,openright]{book}
\pdfoutput=1
\usepackage[nofoot]{geometry}
\usepackage{graphicx}
\usepackage{xcolor}
\usepackage{fancyhdr}

\usepackage[english,dutch]{babel}

\usepackage{relsize}
\usepackage{ifthen}
\usepackage{tikz}
\usepackage{mathrsfs}
\usepackage{amssymb}
\usepackage{amsfonts}
\usepackage{marvosym}
\usepackage{amsmath, latexsym, amscd}
\usepackage[T1]{fontenc} 
\usepackage[center,small]{caption} 
\usepackage{epsfig}
\usepackage{pict2e}

\iffalse 
\definecolor{mydarkblue}{rgb}{0,0,0.7}
\definecolor{mydarkgreen}{rgb}{0,0.7,0}
\definecolor{mydarkred}{rgb}{0.7,0,0}
\definecolor{fileblue}{rgb}{0,0,1}
\else 
\definecolor{mydarkblue}{rgb}{0,0,0}
\definecolor{mydarkgreen}{rgb}{0,0,0}
\definecolor{mydarkred}{rgb}{0,0,0}
\definecolor{fileblue}{rgb}{0,0,0}
\fi

\usepackage[pdftex,plainpages=false,hyperfootnotes=false,pdfpagelabels,
            linkcolor=mydarkblue,menucolor=mydarkblue,%
            citecolor=mydarkred,urlcolor=fileblue,filecolor=fileblue,%
            linktocpage,colorlinks]{hyperref}

\input{defs.tex}
\begin{document}
\selectlanguage{english}

\newcommand{\clearplaindoublepage}{\newpage{\pagestyle{plain}\cleardoublepage}}

\renewcommand{\headrulewidth}{0pt}%
\renewcommand{\footrulewidth}{0pt}%

\pagestyle{fancy}
\fancyhead{}
\fancyfoot{}
\fancyhead[LE,RO]{\thepage}%
\fancyhead[CE]{\leftmark} 

\renewcommand{\chaptermark}[1]{\markboth{\textsc{Chapter} \thechapter}{}}%
\renewcommand{\sectionmark}[1]{\markright{\textsc{\thesection{}. #1}}{}}%

\fancypagestyle{plain}{%
\fancyhead[LE,RO,LC]{}%
\fancyhead[LE,RO]{\thepage}
}

\fancypagestyle{empty}{%
\fancyhead[LE,RO,LC]{}%
}

\makeatletter
\renewcommand\chapter{\if@openright\clearplaindoublepage\else\clearpage\fi
                    \thispagestyle{empty}%
                    \global\@topnum\z@
                    \@afterindentfalse
                    \secdef\@chapter\@schapter}
                    
\renewcommand\backmatter{%
  \if@openright
    \clearplaindoublepage
  \else
    \clearpage
  \fi
  \@mainmatterfalse}                    
                    
\makeatother                   

\frontmatter

\pagestyle{empty}

	\newlength{\cacheparindent}
\setlength{\cacheparindent}{\parindent}
\setlength{\parindent}{0pt}

\def\Fulltitle{Transformation \& Uncertainty
Some Thoughts on 
Quantum Probability Theory,\\ 
Quantum Statistics, and Natural Bundles}
\def\Title{Transformation \& Uncertainty}
\def\Subtitle{Some Thoughts on Quantum Probability Theory,
\\ Quantum Statistics, and Natural Bundles}
\def\NLTitle{Transformatie \& Onzekerheid}
\def\NLSubtitle{Idee\"{e}n over Kwantumkansrekening,\\
Kwantumstatistiek
en 
Natuurlijke Bundels}

\title{\Fulltitle}
\author{Bas Janssens}

	\makeatletter
	\begin{titlepage}%
	\null\vfil
	\begin{center}%
	{\bf \LARGE \Title \par}%
	\vskip 2em%
	{\large \Subtitle \par}%
	\end{center}
	\vfil\null
	\end{titlepage}

	\makeatother

\newpage
\thispagestyle{empty}
\noindent\emph{Thesis committee}:\\
Prof.~dr.~E.\,P.~van~den~Ban, Universiteit Utrecht\\
Dr.~N.~Datta, University of Cambridge\\
Prof.~dr.~M.~Fannes, Katholieke Universiteit Leuven\\
Prof.~dr.~E.\,J.\,N.~Looijenga, Universiteit Utrecht\\
Prof.~dr.~C.~Wockel, Universit\"at Hamburg\\

\vfill
	%
ISBN 978-90-393-5409-4\\[1mm]\medskip
Copyright~\copyright~2010 by B. Janssens.

\newpage
\thispagestyle{empty}
\pdfbookmark[0]{Title page}{tit}
\noindent
	\null\vskip10mm
	{\Huge \Title}\\
	\\
	{\Large \Subtitle}

	\vfill

	{\Huge \NLTitle}\\
	\\
	{\Large \NLSubtitle}
	\\
	\\
	\noindent\normalsize (met een samenvatting in het Nederlands)

	\vfill
	\vfill
	\vfill

	\noindent\Large Proefschrift
	\par
	\bigskip
	\par\noindent\normalsize  ter verkrijging van de graad van doctor aan
	de Universiteit Utrecht op gezag van de rector magnificus,
	prof.~dr.~J.\,C.~Stoof, ingevolge het besluit van het college voor
	promoties in het openbaar te verdedigen op maandag 4 oktober 2010 des
	middags te 2.30 uur

	\vspace*{1ex}
	\noindent
	door

	\vspace*{4ex}

	Bas Janssens
	
	\vspace*{3ex}
	
	\par\noindent\normalsize geboren op 20 mei 1981 te Maastricht

\newpage
\thispagestyle{empty}
	\noindent
	\begin{tabular}{ll}
{}	Promotor: & Prof.~dr.~R.~Fern\'andez
{}	\\
{}	Copromotoren: & Dr.~J.\,W.~van~de~Leur
{}	\\
{}	& Dr.~J.\,D.\,M.~Maassen
	\end{tabular}%

%
%

\newpage
\thispagestyle{empty}
	\null
	\vfill
	\begin{flushright}
	\itshape\selectfont
  \parbox{7.4cm}{
    \ldots and to him the meaning of an episode was not inside like a kernel 
    but outside, enveloping the tale which brought it out only as a glow
    \mbox{brings out} \mbox{a haze}\ldots
    }\\\medskip	
    \normalfont
    Joseph Conrad, Heart of Darkness
\end{flushright}

\newpage
\thispagestyle{empty}

\setlength{\parindent}{\cacheparindent}




	\tableofcontents

\mainmatter
\pagestyle{fancy}
\fancyhead[LE,RO]{\thepage}%
\fancyhead[CE]{\leftmark} 

\chapter{Introduction}
\label{ch:introductie}
\setsubdir{hoofdstukken/intro/}

This thesis aims to investigate mathematically two of the most robust 
phenomena in modern physics. 
The first is the `uncertainty principle' in quantum mechanics.
It states that no measurement can be performed
without disturbing the system.
The second is the `equivariance principle' in classical field theory,
which states that any transformation of space-time must also 
transform the fields. 

The word `robust' should be understood in the sense 
of `independent of the dynamics of the particular model'. 
A physical theory can usually be separated into a kinematic and a
dynamical part. The kinematic part of a theory describes
the basic
physical objects
in terms of mathematics.
The dynamical part then prescribes
the precise way in which these variables
evolve. Our results will largely be independent of the dynamics,
making them applicable to a wide range of different models.

In the case of quantum mechanics for example, the kinematic part 
of the theory states that a system is described by a Hilbert space,
observables by self-adjoint operators, states by vectors and 
transformations by unitary operators. This alone will suffice for us 
to prove quantitative versions of the uncertainty principle, 
as well as asymptotic
optimality of our many-qubit state estimation scheme.
The dynamics of the theory, embodied by a model-specific Hamilton
operator, only come into play when we test these ideas
in specific situations.

In the case of classical field theory, the kinematic part of the theory
entails that 
space-time is modelled by a smooth manifold,
fields are modelled by sections
of a bundle over this manifold, and
transformations are described by bundle automorphisms.  
This alone will suffice for us to classify the relevant bundles.
The dynamical part of the theory, in the form of a model-specific
Lagrangian, will have no part in this. 

Although I would warmly recommend it, I am aware that the reader may not
wish to read this thesis from cover to cover. Chapters \ref{ch:ITISC}
through \ref{ch:BLID} can be read independently, with perhaps
an occasional glance at the `background' chapter \ref{ch:intro}.
We give a brief description of each chapter.
\paragraph{Chapter \ref{ch:intro}. Background}
We give a short introduction to quantum probability theory and
classical field theory, providing background and motivation for
the rest of the thesis.
\paragraph{Chapter \ref{ch:ITISC}. Information Transfer Implies State Collapse}
We prove that within the framework of quantum mechanics,
a transfer of information from a system to the outside
necessarily causes decoherence on that system. 
We also prove a quantitative version of this, linking
the amount of decoherence to the quality of information transfer.
Finally, we shed some light on the fact that coherence 
is very hard to observe in macroscopic systems. 
This chapter is an adapted version of \cite{JM}, written 
together with Hans Maassen. 
\paragraph{Chapter \ref{ch:UDHP}. Unifying Decoherence and the Heisenberg Principle}
\hspace{-2mm}With\-in the framework of quantum probability theory, 
we prove quantitative versions of the Heisenberg principle, 
the joint measurement theorem, the no-cloning theorem, and
the collapse of the wave function. 
An important difference with chapter \ref{ch:ITISC} is that 
all bounds will be in terms of the measurement procedure alone;
they will not depend
on the particular state of the measured system.
This is a good thing because the state of the measured system
is generally unknown-- hence the need for measurement.
A part of this chapter appeared earlier in
\cite{JanssensClassicalCoding}.
\paragraph{Chapter \ref{ch:OPJM}. Optimal Pointers for Joint Measurement
}
We apply 
the joint measurement inequality of chapter \ref{ch:UDHP}
to the specific situation of
a 2-level atom, i.e.\ a qubit, coupled to the quantized electromagnetic field. 
Information on $\sigma_{x}$ and $\sigma_{z}$ leaks into the field
and affects its quadrature, 
which can then be probed using homodyne detection.
The joint measurement inequality provides us with a sharp upper bound
on the quality of any simultaneous measurement of 
$\sigma_{x}$ and $\sigma_{z}$.  
Using quantum stochastic calculus, we investigate how close
one can come to this bound using the above setup.
%
Somewhat to our surprise, we found that it cannot be reached,
although one can come as close as 5.6\%.
This chapter is an adapted version of \cite{Luc}, written together 
with Luc Bouten.
\paragraph{Chapter \ref{ch:OEQS}. Optimal Estimation of Qubit
States 
}
We propose an asymptotically optimal scheme to estimate the state 
of $n$ identically prepared qubits.
The qubits are coupled to the quantized electromagnetic field
for some period of time, during which 
continuous time measurements on the field are performed.
%
The proof of optimality
relies on the concept of local asymptotic normality for qubits: for large $n$,
the ensemble of qubits becomes statistically equivalent to the
tensor product of a classical Gaussian distribution 
and a Gaussian state on a quantum harmonic oscillator.
The Gaussian states are well understood, 
which means that we can use
the statistical correspondence to better understand the qubits.

In a way, this is not unlike the classical problem of determining
the bias of an unfair coin that is tossed $n$ times.
In this case, 
the $n$ coin tosses become 
statistically equivalent to a single Gaussian
in the limit of large $n$.
One should guard against taking this analogy too far though:
in the quantum case, there is a difference between 
the act of handing someone $n$ identical qubits to probe, 
and the act of handing someone $n$ times a single qubit to probe. 
At least for mixed states,
a joint measurement on an ensemble of $n$ qubits performs strictly better
than the combined information of $n$ single qubit measurements.
The whole is more than the sum of its parts, so to speak. 
This chapter is an adapted version of \cite{GutaJanssensKahn}, 
which was written together with M\u{a}d\u{a}lin Gu\c{t}\u{a}
and Jonas Kahn.
\paragraph{Chapter \ref{ch:BLID}. Bundles with a Lift of Infinitesimal 
Diffeomorphisms}
A bundle is called `natural' if its fibres transform in a local fashion 
under diffeomorphisms of the base. Examples include 
tangent bundles, frame bundles, and most other bundles
of a local geometric nature.
%
We extend the notion of a `natural' bundle to that of an
`infinitesimally natural' one by requiring the fibres 
to transform
under \emph{infinitesimal} diffeomorphisms only.
We classify exactly what is gained by this extension.
Bundles of spin type are never natural, but it turns out that 
some of them are
indeed natural in the above infinitesimal sense.
In classical field theory, fields are sections of a bundle
over space-time. 
Spin and energy of a field are well-defined 
precisely when the corresponding bundle is 
infinitesimally natural. 
The fact that bundles that correspond to fermionic fields 
(such as quarks, electrons, etc.) are never natural
shows that the generalization is really called for.

\paragraph{Samenvatting} We give a
non-rigorous Dutch summary of the main results,
containing no formulae and as many pictures as possible. 
If you don't have much experience with (this kind of) mathematics, 
and if you understand Dutch, 
then this may be a good place to start reading.

\chapter{Background}
\label{ch:intro}
\setsubdir{hoofdstukken/intro/}
The chapters to come can be read independently, and are 
self-sufficient to a certain degree.
Complete autarky however is a dream, and, when realized in a PhD thesis,
probably a nightmare.
From time to time then, the reader is referred elsewhere, asked to
trust a statement bluntly placed before him, or simply assumed to 
`know what is going on'.
This is inevitable. 
The flip side of the progress of science is
the overspecialization that comes with it. 

However, as details and applications become more and more complex, 
the foundations of science have a tendency
to simplify
with time.
This makes it feasible for us to sketch the fundamentals
of quantum probability theory
and classical field theory in the next forty-odd pages.
We aim to 
explain some of the more vital ingredients that will enter 
our reasoning, but the main goal here is to provide
background and motivation
for the problems that will be addressed later on.



We will start by giving a short introduction to quantum mechanics,
focusing first on the way in which the theory gives rise to probability measures,
and then on matters of transformation 
and symmetry.
An even shorter introduction to 
quantum probability theory will follow, with an emphasis on transformations
of the theory.  
Everything will be done first in the finite dimensional setting, 
and only
then in the general case. 
The purpose of this is to 
prevent the elegant, yet somewhat involved
functional analytic side of the theory
from overshadowing its simpler algebraic counterpart. 

An important tool for us will be the quantum stochastic
calculus of Hudson and Parthasarathy.
This is 
a generalization 
of the It\^{o} stochastic calculus on Wiener space,
which is essentially calculus with stochastic
infinitesimal increments.
After a primer on Fock spaces and Weyl operators, 
we will spend some words
on the Wiener measure and It\^{o} calculus.
We will then see how the field operators in a free quantum
field theory form 
a quantum version of
the stochastic infinitesimal increments of It\^{o} calculus,
giving rise to 
quantum stochastic calculus.

The last part of this chapter will be a brief introduction 
to classical field theory from a differential geometric point of view.
In a nutshell, this means that we will consistently treat fields as sections 
of a fibre bundle.
The emphasis is once again on transformation and symmetry.
We formulate the requirement that fields should transform in a definite fashion
under infinitesimal space-time transformations.
The ramifications of this will become clear in chapter \ref{ch:BLID}.

As already stated,
we mean to sketch only the most elemental
aspects of the theory. For a more complete picture, at perhaps a more
responsible pace, the reader is referred to 
\cite{Par92}, \cite{Frankel}, \cite{FR},
\cite{vdBanGQT}, and
\cite{Hans}, on which we based the following account.


\section{Quantum Mechanics}
%
In quantum mechanics, a system is modelled by a Hilbert space $\CH$.
%
A state of the system is modelled by a ray $\mathbb{C}\psi$ in $\CH$,
and observables are represented by self-adjoint, possibly unbounded
operators $A$ on $\CH$. 
Transformations of a closed system are described by unitary (or antiunitary)
operators $U$ on $\CH$, and there is a distinguished continuous 
one-parameter group
$t \mapsto U_t$ of unitary operators that describes time evolution.  

Quantum mechanics can only be interpreted in terms of probabilities.
There is a canonical 
map $(\mathbb{C}\psi,A) \mapsto \mathbb{P}_{\psi}$
that associates to each ray $\mathbb{C}\psi$ and self-adjoint operator $A$
a probability measure $\pp$ on the spectrum of $A$. It is interpreted as
follows.
\begin{center}
\emph{If the system is in state $\mathbb{C}\psi$, then observation of 
an observable $A$
yields an outcome $\lambda \in S$ with probability $\pp(S)$.}
\end{center}
Everything else follows from these postulates. 

\subsection{Probabilities and the Spectral Theorem}
We will first describe the map $(\mathbb{C}\psi,A) \mapsto \mathbb{P}_{\psi}$ 
in the special case of
finite dimensional Hilbert space, 
and then in the general case.

The 
finite dimensional case is
simple but interesting.
It is simple because
a choice of orthonormal basis identifies $\CH$ 
with $\mathbb{C}^{n}$, 
reducing many of our problems to linear algebra.
It is interesting because the results in chapter 
\ref{ch:ITISC} and \ref{ch:UDHP} 
are as nontrivial
here as they are in the general case, and
--if the reader will forgive me my patronizing remark-- 
it might be helpful to keep the finite dimensional case in mind 
when reading these chapters.
  
Let us take the opportunity to fix some notation.
All Hilbert spaces are over $\mathbb{C}$ unless specifically mentioned
otherwise, and we will take
the   
inner product $\inp{\,\,\cdot\,\,}{\,\,\cdot\,\,}$ 
to be linear on the right. The adjoint of $A$ is denoted $A^{\dagger}$,
and $\CB(\CH)$ denotes the algebra of bounded 
linear operators
on $\CH$. For $n$-dimensional $\CH$, this is just the algebra of 
$n \times n$-matrices over $\mathbb{C}$.
 

\subsubsection{Hilbert Spaces of Finite Dimension}
An orthogonal projection $E$
is in particular a self-adjoint operator.
Because its spectrum is $\{0,1\}$,
it should be interpreted as an \emph{event}. The outcome
1 means that the event occurs, the outcome 0 means that it does not.

If $E$ projects onto $\CH_{1} \subset \CH$, then we have the orthogonal
decomposition 
$\CH = \CH_{0} \oplus \CH_{1}$. A unit vector $\psi$ splits as 
$\psi = \psi_0 + \psi_1$, with $\|\psi_{0}\|^2 + \|\psi_{1}\|^2 = 1$.
This leads us to define the
probability distribution $\pp$ on $\{0,1\}$ by 
$$
\pp(0) = \|\psi_0\|^{2} \,, \quad 
\pp(1) = \|\psi_1\|^{2} \,.
$$
This depends only on the ray $\mathbb{C}\psi$, so that we have our map
$(\mathbb{C}\psi , A) \mapsto \pp$ in the special case that $A$ is a projection.

We can play the same game if we split $\CH$ into any number of mutually
orthogonal spaces. Such a decomposition determines a projection valued measure.
\begin{basdef}Let $\Omega$ be a finite set, and $\CH$ a finite dimensional
Hilbert space. Then a \emph{projection valued measure} (or PVM) is a map
$E : \Omega \rightarrow \CB(\CH)$ such that:
\begin{itemize}
\item[-]$E(\omega)$ is an orthogonal projection for all $\omega \in \Omega$.
\item[-]If $\omega \neq \omega'$, then $E(\omega)E(\omega') = 0$.
\item[-]$\sum_{\Omega} E(\omega) = \one$.
\end{itemize}
\end{basdef}
If $\CH_{\omega}$ is the image of $E(\omega)$, it is clear that
$\CH = \bigoplus_{\Omega} \CH_{\omega}$. Accordingly, a unit vector
$\psi$ decomposes as $\psi = \sum_{\Omega} \psi_{\omega}$, with 
$\|\psi_{\omega}\|^{2} = \inp{\psi}{E(\omega)\psi}$ and 
$\sum_{\Omega} \inp{\psi}{E(\omega)\psi} = 1$.
This means that a unit vector $\psi$ gives rise to a measure
$\pp$ on $\Omega$, defined by
$
\pp(\omega) = \inp{\psi}{E(\omega)\psi}
$.

For any unit vector $\psi$, we will define the map
$
\rho_{\psi} : \CB(\CH) \rightarrow \mathbb{C}\,,
$
called a \emph{pure state}, by
$$
\rho_{\psi}(A) = \inp{\psi}{A\psi}\,.
$$
With this notation, we have
$$
\pp = \rho_{\psi} \circ E \,.
$$ 
Note that $\rho_{\psi} = \rho_{\psi'}$ if and only if 
$\psi' = c \psi$ for some $c \in \mathbb{C}^{\times}$, so that 
a pure state corresponds precisely to a 
ray $\mathbb{C}\psi$ in the Hilbert space.
We will identify the two, and think of the \emph{pure state space}
as either
$$
\mathcal{P}(\CH) = \{\rho_{\psi} \,;\, \psi \in \CH,\, \|\psi\| = 1\}
$$
or as the space of rays
$$
\mathcal{P}(\CH) =  \CH - \{0\}\, /\, \mathbb{C}^{\times}\,.
$$ 
We leave it up to the reader to decide whether  
$\mathcal{P}$ stands for 
`pure' or `projective'.

The spectral theorem says that each self-adjoint operator $A$ gives rise to a 
projection valued measure on its spectrum. 
\begin{theorem}
Let $A$ be a self-adjoint operator on a finite dimensional Hilbert space $\CH$,
and let $\spec(A)$ be its spectrum of eigenvalues. Then there exists
a unique projection valued measure $E_{A} : \spec(A) \rightarrow \CB(\CH)$ such that
$$
A = \sum_{\lambda \in \spec(A)} \lambda E_{A}(\lambda)\,.  
$$
\end{theorem}    
The $E_{A}(\lambda)$ are the projections 
onto the eigenspaces $\CH_{\lambda}$, and 
$\CH = \bigoplus_{\spec(A)} \CH_{\lambda}$
is the eigenspace decomposition.
 
The map $(\mathbb{C}\psi , A) \mapsto \pp$ is now defined as follows. 
\begin{itemize}
\item[-] The ray $\mathbb{C}\psi$ is identified with the pure state
$\rho_{\psi} : \CB(\CH) \rightarrow \mathbb{C}$.
\item[-] An observable $A$ gives rise to a PVM 
$E_{A} : \spec(A) \rightarrow \CB(\CH)$.
\item[-]The measure $\pp$ on $\spec(A)$ is given by
$\pp = \rho_{\psi} \circ E$.
\end{itemize} 
We interpret $\pp(\lambda)$ 
as the probability that if the system is in state $\mathbb{C}\psi$,
a measurement of $A$ yields the outcome $\lambda$. 


\subsubsection{Infinite Dimensional Hilbert Spaces}

Although there do exist physically relevant systems where the Hilbert space
is finite dimensional, this is not typically the case.
We adapt our framework to encompass general Hilbert spaces.

We equip the algebra $\CB(\CH)$ of bounded operators on $\CH$ with the
\emph{weak topology}, the coarsest topology that makes
$\rho_{\psi}(A) = \inp{\psi}{A\psi}$ continuous for all $\psi$ in $\CH$.
We then introduce the notion of a projection valued measure for general Hilbert
spaces as follows.


%
\begin{basdef}\label{pvem} 
Let $(\Omega,\Sigma)$ be a measurable space, and $\CH$ a Hilbert space.
Then a 
\emph{projection valued measure} or PVM
is a map   
$E : \Sigma \rightarrow \CB(\CH)$ with the following properties.
\begin{itemize}
\item[-]$E(S)$ is a projection for all $S \in \Sigma$.
\item[-]$E(\emptyset) = 0$, and $E(\Omega) = \one$. 
\item[-]If $S_1 \cap S_2 = \emptyset$, then $E(S_1)E(S_2) = 0$.
\item[-]If $S_i$ is a sequence of disjoint sets in $\Sigma$, then
$E(\bigcup_{i=0}^{\infty}S_i ) = \sum_{i=0}^{\infty}E(S_i)$, where 
the limit is in the weak topology. 
\end{itemize}
\end{basdef}
Each ray $\mathbb{C}\psi$
induces a probability measure $\pp$ on $(\Omega,\Sigma)$ by
applying the pure state
$\rho_{\psi}(A) := \inp{\psi}{A\psi}$ to $E$, i.e. 
$$
\pp(S) = (\rho_{\psi} \circ E)(S)\,.
$$ 

If $f$ is a measurable function on $\Omega$, then 
an integral 
\begin{equation}\label{oetkhijp}
I_E(f) := \int_{\Omega} f(\omega) E(d\omega)
\end{equation}
can be defined by the weak limit of Riemann sums.
The resulting operator $I_{E}(f)$ on $\CH$ 
is bounded if $f$ is bounded,
which invites us to consider 
$I_{E}$ as a linear map
$I_E : L^{\infty}(\Omega,\Sigma) \rightarrow \CB(\CH)$.
It is a homomorphism of $*$-algebras, i.e.
it respects multiplication 
$I_E(fg) = I_E(f)I_E(g)$ and adjunction $I_E(f^{\dagger}) = I_E(f)^{\dagger}$.
In particular, $I_E(f)$ is self-adjoint if $f$ is real, and 
positive if $f$ is positive.
With respect to $\pp$, we have
$\mathbb{E}(f) = \rho_{\psi}(I_E(f))$. 

Observables correspond to self-adjoint (i.e. symmetric, closed 
and possibly unbounded)
operators on $\CH$. 
The spectral theorem says that any self-adjoint operator
gives rise to a projection valued measure on its spectrum, called the 
spectral measure.
\begin{theorem}[spectral theorem]
Let $A$ be a self-adjoint operator on $\CH$.
Then there exists a unique PVM 
$E_{A} : \Sigma(\mathbb{R}) \rightarrow \CB(\CH)$ 
on the Borel sigma-algebra of $\mathbb{R}$
such that,
with $\lambda$ the identity function on $\mathbb{R}$,
$$
A \psi = \int_{\mathbb{R}} \lambda E(d\lambda) \psi\,
$$
for all $\psi$ in the domain of $A$. 
The PVM 
is concentrated on the spectrum 
of $A$. 
\end{theorem}
\proof See for example \cite{Par92}. \qed
Any pure state $\rho_\psi$ therefore induces a probability measure 
$\pp$ on $\spec(A)$ by $\pp = \rho_{\psi} \circ E_{A}$.
Under this probability measure, the expectation
of the observable $A$ 
of a system in a state $\psi$ 
is
$\mathbb{E}(\lambda) = \rho_{\psi}(A)$.
%
Its variance is $\rho_{\psi}(A^2) - \rho_{\psi}(A)^{2}$.

\subsection{Transformation, Symmetry and Conservation Laws}
Let us devote some attention to the way 
in which transformations and symmetries occur in
quantum mechanics, and then give examples of relevant quantum systems.
Although the proper treatment of symmetry in quantum mechanics will not 
play a major role in this thesis,
it will help us
to properly interpret the 
examples. 
 
\subsubsection{Hilbert Spaces of Finite Dimension}
As before, we will first focus attention on finite dimensional Hilbert spaces,
before looking at the general case. Throughout, we will denote by $U(\CH)$
the group of unitary operators on $\CH$.



%

\paragraph{Transformations}
A unitary operator $U \in U(\CH)$
induces an automorphism 
$\alpha_U : A \mapsto U^{\dagger}AU$
of $\CB(\CH)$. This in turn
acts on $\mathcal{P} (\CH)$ by its dual, 
$\alpha_U^* : \rho_{\psi} \mapsto \rho_{\psi} \circ \alpha_U$.
The action of the unitary operators on the pure state space
is therefore governed by the group of
inner automorphisms
$\mathcal{P}U(\CH) = \{\alpha_U \,;\, U \in U(\CH)\}$.

Recall that the pure state space has the equivalent
description
$\mathcal{P}(\CH) = \CH  - \{0\}/\mathbb{C}^{\times}$, on which $\alpha_U$
acts by $\mathbb{C}\psi \mapsto \mathbb{C}U\psi$.  
This shows that 
$\alpha_U = \alpha_{U'}$ if and only if 
$U' = c U$ for some $c \in \mathbb{C}^{\times}$.
We therefore have 
the equivalent description
$\mathcal{P}(U(\CH)) = U(\CH) / \mathbb{C}^{\times}$.
Apparently, the $\mathcal{P}$ stands for `projective'.

The action of an abstract Lie group $G$ on a quantum mechanical system
should therefore be described by a
\emph{projective} unitary representation, i.e.
a continuous group homomorphism
$\pi : G  \rightarrow \mathcal{P}U(\CH)$.

Recall that a \emph{linear} unitary representation
is a continuous group homomorphism 
$\pi : G  \rightarrow U(\CH)$.
Every linear unitary representation gives rise to a 
projective unitary representation, 
but it is not true that
every projective representation arises in this way. 

Since the theory of linear unitary representations
is well understood, we would like to reduce projective representations
to linear ones. This is achieved by the following theorem,
at the expense of slightly enlarging the group.
\begin{theorem}\label{symprojhond}
Let $G$ be a connected Lie group, with universal cover $\tilde{G}$.
Let $\CH$ be a finite dimensional Hilbert space.
Then every projective unitary representation $\pi$ of $G$ on 
$\CH$ comes from a linear
unitary representation $\tilde{\pi}$ of $\tilde{G}$ on $\CH$, in the sense that 
the following diagram commutes.
\begin{center}
\begin{tikzpicture}
\node (linksboven) at (0, 1.2, 0) {$\tilde{G}$};
\node (rechtsboven) at (1.5, 1.2, 0) {$U(\CH)$};
\node (linksonder) at (0, 0, 0) {$G$};
\node (rechtsonder) at (1.5,0,0) {$\mathcal{P}U(\CH)$};
\draw [->] (linksboven) -- node[above] {$\tilde{\pi}$}(rechtsboven);
\draw [->] (linksonder) -- node[above]{$\pi$}(rechtsonder);
\draw [->] (linksboven) -- (linksonder);
\draw [->] (rechtsboven) -- (rechtsonder);
\end{tikzpicture}
\end{center}
\end{theorem}
\proof Since $\CH \simeq \mathbb{C}^{n}$ for some finite $n$, we 
can consider $\pi$ as a homomorphism 
$G \rightarrow \mathcal{P}(U(\mathbb{C}^{n}))$.
The derived Lie algebra homomorphism 
$\mathrm{Lie}(G) \rightarrow \mathfrak{su}(\mathbb{C}^{n})$
extends to a group homomorphism 
$\tilde{G} \rightarrow \mathrm{SU}(\mathbb{C}^{n})$ by Lie's 
second theorem. 
\qed

\noindent In what follows, we will assume that $G$ has been replaced by
the universal cover of its connected component of unity $\tilde{G}_{0}$,
and $\pi$ by a linear unitary representation of $\tilde{G}_{0}$.



\paragraph{Conservation Laws}
It is a remarkable fact of life that in every fundamental physical theory,
infinitesimal symmetries correspond to conservation laws.
In quantum mechanics, this correspondence is particularly straightforward. 


Infinitesimal transformations are elements of the Lie algebra 
$\mathrm{Lie}(G)$ of $G$. 
The following theorem describes how 
they give rise to observables.
\begin{theorem}
Let $\pi$ be a continuous unitary representation of $G$ on $\CH$. 
Then for each $X \in \mathrm{Lie}(G)$, we have a well defined 
skew-symmetric 
operator $\dot{\pi}(X)$ defined by
\begin{equation}\label{schroot}
\dot{\pi}(X)  = \pi(e^{-tX})\frac{d}{dt} \pi(e^{tX})\,.
\end{equation}
As $\dot{\pi}([X,Y]) = [\dot{\pi}(X),\dot{\pi}(Y)]$, this defines a 
Lie algebra homomorphism 
$\dot{\pi}$ 
of $\mathrm{Lie}(G)$ into the Lie algebra $\mathfrak{su}(\CH)$
of skew-symmetric operators on $\CH$.
\end{theorem}
%
The operator $i \dot{\pi}(X)$ is self-adjoint, 
and should be considered as the observable corresponding to 
the infinitesimal transformation $X$. We denote it $J(X)$. 
If one has to pass to the
universal cover, then $\dot{\pi}$ and $J$ are still well-defined, 
since the Lie algebras of $G$ and $\tilde{G}_{0}$ are the same. 

Time evolution is given by a one-parameter group of unitary
transformations,
i.e. a unitary representation $t \mapsto U_t$
of $(\mathbb{R},+)$. Its Lie algebra is again
$\mathbb{R}$, and the self-adjoint operator $J(1)$
is called the \emph{Hamiltonian}, denoted $H$.  
The observable associated to this operator is the energy.

Equation (\ref{schroot}) for the particular case of time evolution
reads 
$$
i \frac{d}{dt}U_t = H U_t\,.
$$
In the `Heisenberg picture', where $U$ acts by $\alpha_U$, one 
writes $A_t = U_t^{\dagger} A U_t$
to obtain \emph{Heisenberg's equation of motion}
\begin{equation}
\frac{d}{dt}A_t = i[H,A_t]\,.
\end{equation} 
In the dual `Schr\"odinger picture', where $U$ acts by $\alpha_{U}^{*}$, one writes $\psi_t = U_t \psi$ to 
obtain the \emph{Schr\"odinger equation}
\begin{equation}\label{schrodfin}
i \frac{d}{dt}\psi_t = H \psi_t\,.
\end{equation}  

A transformation is called a \emph{symmetry} if it commutes with the
time evolution $t \mapsto U_t$,  
i.e.~if $[U_t, \pi(g)] = 0$ for all $t\in \mathbb{R}$.
A group $G$ is a group of symmetries if $\pi(g)$ is a symmetry 
for every $g \in G$.  
Because this implies that $\dot{\pi}(X)$ commutes with $H$, we immediately
obtain from Heisenberg's equation of motion the conservation law
$$
\frac{d}{dt}J(X)_t = 0\,. 
$$
If $X$ is an infinitesimal symmetry, then the 
probability measure on the spectrum of $J(X)$ is constant in time.
For $J(X) = H$, this is called conservation of energy. 

The fundamental problem in
quantum mechanics is to reconstruct
the time evolution $t \mapsto U_t$ from the given Hamiltonian $H$.
This is usually done by finding 
the eigenspace 
decomposition  
$\CH = \bigoplus_{\spec(H)}\CH_{E}$ for $H$,
and then computing $U_t = e^{-itH}$.
This is greatly facilitated by the presence of symmetries.
 

\paragraph{The Qubit}
Time for an example.
The internal degrees of freedom of an 
electron 
are 
described 
by the Hilbert space $\CH = \mathbb{C}^{2}$. 
It transforms under the rotation group
$\mathrm{SO}(3,\mathbb{R})$ by a projective unitary representation.
 
Consider the map $\mathcal{P}U(\mathbb{C}^{2}) \rightarrow 
\mathrm{SO}(3,\mathbb{R})$ defined as follows. 
In terms of the basis
$$
\one = 
\begin{pmatrix}
1 & 0\\
0 & 1
\end{pmatrix}
,\quad
\sigma_{x} = 
\begin{pmatrix}
0 & 1\\
1 & 0
\end{pmatrix}
, \quad
\sigma_{y} = 
\begin{pmatrix}
0 & -i\\
i & 0
\end{pmatrix}
, \quad
\sigma_{z} = 
\begin{pmatrix}
1 & 0\\
0 & -1
\end{pmatrix}
$$
of $\CB(\mathbb{C}^{2})$, one readily checks that
$\alpha_U (\lambda \one) = \lambda \one$,
and that
$$
\alpha_{U} (\sum_{i=1}^{3} x_i \sigma_i) 
=
\sum_{i,j=1}^{3} R_{ij} x_j \sigma_i
$$
for some $R \in \mathrm{SO}(3,\mathbb{R})$. 
The map $\alpha_u \mapsto R$ 
turns out to be an isomorphism 
$\mathcal{P}U(\mathbb{C}^{2}) 
\stackrel{\sim}{\rightarrow} \mathrm{SO}(3,\mathbb{R})$,
and its inverse $\pi$ is the projective representation of $\mathrm{SO}(3,\mathbb{R})$
on $\mathbb{C}^{2}$.
If we choose the basis 
$$
\tau_{x} = 
\begin{pmatrix}
0 & 0 & 0\\
0 & 0 & 1\\
0 & -1 & 0
\end{pmatrix}
,\quad
\tau_{y} = 
\begin{pmatrix}
0 & 0 & 1\\
0 & 0 & 0\\
-1 & 0 & 0
\end{pmatrix}
,\quad
\tau_{z} = 
\begin{pmatrix}
0 & 1 & 0\\
-1 & 0 & 0\\
0 & 0 & 0
\end{pmatrix}
$$ 
of $\mathfrak{so}(3,\mathbb{R})$, 
then the Lie algebra homomorphism $\dot{\pi} : \mathfrak{so}(3,\mathbb{R}) 
\rightarrow \mathfrak{su}(\mathbb{C}^{2})$ gives rise to three observables
$J(\tau_x) = \half \sigma_x$,
$J(\tau_y) = \half \sigma_y$ and
$J(\tau_z) = \half \sigma_z$.
They should therefore be considered as the 
\emph{internal angular momenta} or \emph{spin} of the electron.

If the particle is placed in a vertically aligned magnetic field 
of strength $B$, then the effective 
Hamiltonian of the system is $H = B \sigma_{z}/2$.
The group of symmetries is then \emph{not} all of 
$\mathrm{SO}(3,\mathbb{R})$, but only
the subgroup $\mathrm{SO}(2,\mathbb{R})$ of 
rotations around the $z$-axis. The only conserved quantity is therefore 
$H \sim \half \sigma_{z}$.
Heisenberg's equation of motion $i \dt A_t = [H,A_t]$ can be solved to yield 
$$
t \mapsto a \Big( \cos(Bt + \phi) \sigma_x + \sin(Bt + \phi) \sigma_y \Big) + 
b \sigma_z\,,
$$
which says that the angular momentum rotates around the $z$-axis 
with a frequency proportional to the magnetic field.
In particular, the angular momentum in the $z$-direction
is preserved, as was predicted by symmetry considerations.

The system $\CH = \mathbb{C}^{2}$ is called a \emph{qubit}, and one can
use it to model any quantum system with only two (relevant) states,
not just the electron.
For example, 
in quantum optics, the interaction of an atom with a laser beam is 
often such that only two energy-eigenstates of the atom 
are affected. It is then sufficient to model this atom by
a qubit $\mathbb{C}^{2}$, rather than by its full Hilbert space.
We will make extensive use of this in chapters \ref{ch:OPJM} and \ref{ch:OEQS}.   

%
%
%
%
%
%
%
%
%
%

\subsubsection{Infinite Dimensional Hilbert Spaces}\label{gramsmid}
In infinite dimensional Hilbert space, the situation
is similar, but slightly more involved.
The only qualitative difference is that for some --but certainly not for all--
Lie groups $G$, the universal cover $\tilde{G}$ no longer suffices
to describe all projective unitary representations, in which case 
a central deformation of the algebra of symmetries occurs.

\paragraph{Transformations}
A Lie group $G$ of transformations acts on $\CH$ by a 
projective
unitary representation. This is a continuous group
homomorphism $G \rightarrow \mathcal{P}U(\CH)$
as before, 
but as we are no longer in a finite setting, we must
pay some attention to matters of topology.

We equip $\mathcal{P}(\CH)$ with the \emph{trace distance}\label{loosspoor}
or Kolmogorov-distance,
which is the worst case difference in probability on a 
single event (projection) $E$, i.e.
$$d(\rho_{\psi} , \rho_{\phi}) = \sup_{E} (\rho_{\psi}(E) - \rho_{\phi}(E))\,.$$
The \emph{strong topology}
on $\mathcal{P}U(\CH)$ is the one defined  
by pointwise convergence on $\mathcal{P}(\CH)$.
If we define the strong topology on $U(\CH)$
by pointwise convergence on $\CH$, then the strong topology on  
$\mathcal{P}U(\CH)$ is the strongest one that makes
$U(\CH) \rightarrow \mathcal{P}U(\CH)$ continuous.
It is with respect to these strong topologies that we 
require our (projective) representations to be continuous.  

If the Hilbert space is infinite dimensional, one can still
reduce projective unitary representations of $G$
to linear representations of its universal cover $\tilde{G}$,
provided that the second Lie algebra cohomology of 
$\mathrm{Lie}(G)$ vanishes.
\begin{theorem}\label{symprojhondinf}
Let $G$ be a connected Lie group such that 
$H^{2}(\mathrm{Lie}(G),\mathbb{R})=0$, and let 
$\tilde{G}$ be its universal cover.
Then every projective unitary representation of $G$ comes from a 
unitary representation of $\tilde{G}$.
\end{theorem}
Lie algebras with vanishing second cohomology 
include the real semisimple Lie algebras, but also
for example
the Abelian Lie algebra $\mathbb{R}$, and the (infinite dimensional) 
Lie algebra of vector fields on any smooth manifold of dimension
bigger than one.

If $H^{2}(\mathrm{Lie}(G),\mathbb{R})$ is nonzero, then not all is lost.
For every projective representation of $G$, there exists a central 
extension $\hat{G}$ of $G$ by the circle $S^{1}$, 
and a linear unitary
representation of $\hat{G}$ that induces $\pi$.
Since the Lie algebra
of observables then follows $\mathrm{Lie}(\hat{G})$
rather than $\mathrm{Lie}(G)$, a `central deformation'
is said to have occurred.

\paragraph{Conserved Quantities}
Let us consider a group of transformations as a unitary representation
$\pi$ of $G$, having replaced the group by a bigger one
if necessary. 
Then the observables corresponding to infinitesimal transformations
$X \in \mathrm{Lie}(G)$ 
are still precisely 
the derived Lie algebra elements $J(X) = i \dot{\pi}(X)$. 
 

Recall that $\pi$ is 
strongly continuous by definition, so that 
$\pi_{\psi}(g) = \pi(g)\psi$ is always a continuous map.
A vector $\psi$ is called \emph{smooth} if 
$\pi_{\psi}$
is smooth, 
and we denote the Fr\'{e}chet space of smooth vectors by $\CH^{\infty}$. 
\begin{theorem}[Stone, G{\aa}rding] \label{StoGa}
Let $\pi$ be
a unitary representation of a Lie group $G$ on 
a Hilbert space $\CH$. Then the smooth vectors $\CH^{\infty}$ 
lie dense in $\CH$.
For each $X \in \mathrm{Lie}(G)$, the expression
$$
\dot{\pi}(X) = \exp(-tX)\dt \pi(\exp(tX))
$$ 
defines a closed, densely defined, 
skew-symmetric operator $\dot{\pi}(X)$.  
Its domain contains $\CH^{\infty}$, and 
on this space, we have 
$$\dot{\pi}([X,Y]) = [\dot{\pi}(X) , \dot{\pi}(Y)].$$
\end{theorem}
\proof See \cite[ch.\ III]{KnRevEx}, but also 
\cite[sec.\ 12,13]{Par92}. \qed

\noindent In particular, the generator of the 1-parameter group
$t \mapsto U_t$ of time translations is again
the Hamiltonian $H$. In realistic models,
its spectrum is usually bounded from below,
but not from above.

\paragraph{A Point Particle on a Line}
Let us consider a single particle on the 1-dimensional
Euclidean space $\mathbb{R}$.
We should have a self-adjoint operator $X$ on $\CH$
to represent its position. The group of Euclidean motions $(\mathbb{R},+)$
should have a projective unitary representation $\pi$ on $\CH$, and its 
effect on $X$ ought to be be a shift, $\pi(s)^{\dagger}X \pi(s) = X + s \one$.
Since $H^{2}(\mathbb{R} , \mathbb{R}) = 0$ and $\tilde{\mathbb{R}} = \mathbb{R}$,
we may take $\pi$ to be linear rather than projective.
The observable $J(1) = i \dot{\pi}(1)$ is called the \emph{momentum},
and denoted $P$. The above relation then differentiates to 
the \emph{canonical commutation relation}
\begin{equation}\label{ccrel}
[X,P] = i\one\,.
\end{equation}
One way to realize this is by setting $\CH = L^{2}(\mathbb{R})$,
$(X\psi)(x) = x\psi(x)$, and $(\pi(s)\psi)(x) = \psi(x-s)$. The associated
momentum $P = -i\frac{d}{dx}$ does indeed satisfy $[X,P] = i\one$.
A theorem of Stone and von Neumann says that this is essentially
the only way to realize (the integrated version of) equation (\ref{ccrel}).

The PVM for $X$ is not overly complicated; $E_X(S)$ is given by 
multiplication with the indicator function $\one_S$. That is, 
$(E_X(S)\psi)(x) = \one_S(x)\psi(x)$. This means that
the probability distribution $\pp^{X}$ induced by $\psi$ 
on the spectrum $\mathbb{R}$ of $X$ is just
$\pp^{X}(S) = \int_{S}|\psi(x)|^2 dx$. The function
$|\psi|^2$ is the probability density of $\pp^{X}$ w.r.t.
the Lebesgue measure.

The PVM for $P$ can be derived from this. 
If $F$ is the
unitary Fourier transform
on $L^{2}(\mathbb{R})$, then one has
$P = F^{\dagger}(-x)F$, so that $E_P(S) = F^{\dagger}E_X(-S)F$.
This means that
$\pp^{P}(S) = \int_S |F\psi|^2 dx$, and the squared Fourier transform
of $\psi$ yields the probability density of $\pp^{P}$.

In  
sections \ref{hofock} and \ref{hogauss},
we will find two other ways to describe the point particle on a line.


\paragraph{A Point Particle in $\mathbb{R}^{3}$}

Let us consider a single particle in the 3-dimensional
Euclidean space $\mathbb{R}^{3}$. We require 3 
commuting self-adjoint operators $X_1$, $X_2$ and $X_3$ to represent its position.
The Euclidean motion group $\mathbb{R}^{3} \rtimes \mathrm{O}(3)$ should
have a projective unitary representation $\pi$ on $\CH$ such that
each $v \in \mathbb{R}^{3}$ shifts $X_i$ to $X_i + v_i \one$, i.e.
$\pi(v)^{\dagger} X_i \pi(v) = X_i + v_i \one$,
and each rotation $R \in \mathrm{O}(3)$ rotates $X_i$ to $\sum_{j=1}^{3} R_{ij}X_j$, 
i.e. $\pi(R)^{\dagger}X_i \pi(R) = \sum_{j=1}^{3} R_{ij}X_j$.
 
The Lie algebra $\mathbb{R}^{3} \rtimes \mathfrak{so}(3)$ of the Euclidean
motion group has zero second cohomology, so $\pi$ induces a linear 
representation $\dot{\pi}$ on the level of Lie algebras.
The standard basis $e_x, e_y, e_z$ for 
$\mathbb{R}^{3}$ and $\tau_x, \tau_y, \tau_z$ for 
$\mathfrak{so}(3)$ gives rise to the momenta
$P_k := J(e_k)$ and the angular momenta $L_k := J(\tau_k)$. 
The commutation relations between the $P_k$ and $L_l$ are then
easily obtained from the Lie bracket on 
$\mathbb{R}^{3} \rtimes \mathfrak{so}(3)$,
and their commutator with the $X_k$ follows from the requirement that the 
$\pi(g)$
shift or rotate the $X_k$.
This yields the following extension of the canonical commutation relations:
\begin{eqnarray}\label{Euclied}
{}[X_k,X_l] = 0\,, & 
{}[X_k,P_l] = i\delta_{kl} \one\,, & 
{}[P_k,P_l]=0 \\ \nonumber
{}[X_k,L_l] = i\epsilon_{klm}X_m \,, &
{}[L_k,L_l] = i\epsilon_{klm}L_m\,, &  
{}[L_k,P_l] = i\epsilon_{klm}P_m \,.
\end{eqnarray}

 
This can be achieved by taking the Hilbert space $\CH = L^{2}(\mathbb{R}^{3})$, 
with the observables $(X_k \psi)(x) = x_k \psi(x)$ and the 
unitary representation $(\pi(g)\psi) (x) = \psi(g^{-1}x)$.
The corresponding observables
$$P_k = -i \frac{\del}{\del x_k} \,, \quad 
L_k = -i\sum_{l,m}\epsilon_{klm} x_l \frac{\del}{\del x_m}$$
can easily be seen to obey \ref{Euclied}. 

Note that we have established all of this without any reference 
to the dynamics whatsoever.
In other words, everything
up to this point is independent of the particular Hamiltonian.  
A Hamiltonian usually takes
the shape $H = \sum_{k=1}^{3} P_{k}^{2}/2m + V(x)$, with $V$ the 
\emph{potential function}. This gives rise to the 
Schr\"odinger equation
$$
i \frac{d}{dt} \psi(x) = -\frac{1}{2m}\sum_{k=1}^{3}\frac{{\del}^2}{\del x_k^2} 
\psi(x) + V(x)\psi(x)\,.
$$

It is the choice of Hamiltonian that determines 
which transformations are symmetries.
For example,
the hydrogen atom is described by the potential 
$V(x) = -1/\sqrt{x_1^2 + x_2^2 + x_3^2}$. In this case, the group
of Euclidean motions $\mathbb{R}^{3} \rtimes \mathrm{O}(3)$ is
\emph{not} a group of symmetries, but its subgroup $\mathrm{O}(3)$
of orthogonal transformations is. This means that the angular momenta
$L_k$ are conserved quantities, whereas the momenta $P_k$ are not.

\section{Quantum Probability Theory}
 
Quantum probability theory (QPT) is a 
mathematical framework which generalizes 
both quantum mechanics and classical probability theory. It is
ideally suited to describe \emph{open} quantum systems, and 
in particular quantum measurement. 

We first describe quantum probability theory on a 
finite dimensional Hilbert space $\CH$. This allows us to 
focus on the algebra, postponing the functional analysis 
involved in the infinite dimensional case to later.

\subsection{QPT on Finite Dimensional Hilbert Spaces.}

A subalgebra $\CA$ of $\CB(\CH)$ is called a $*$-algebra if it is closed
under the adjoint.
In quantum probability theory, 
a system is modelled by a $*$-algebra 
$\CA \subset \CB(\CH)$.
The observables of the system are described 
by self-adjoint operators $A$ in $\CA$, and
the states of the system are modelled as follows.
\begin{basdef}
Let $\CH$ be a finite dimensional Hilbert space, 
and $\CA$ a $*$-sub\-algebra of $\CB(\CH)$.
Then a 
\emph{state} on $\CA$ is a positive, normalized, 
linear functional $\rho : \CA \rightarrow \mathbb{C}$. 
The space of states is denoted $\CS(\CA)$.
\end{basdef}
One should think of a state as a map that assigns to each 
observable $A$ its expectation $\rho(A)$. 
Since the expectation should be $\mathbb{R}$-linear on the observables,
we may as well extend it to a $\mathbb{C}$-linear
map on all of $\CA$. 
Positivity means that $A \geq 0$ implies
$\rho(A) \geq 0$. This is reasonable; 
if $A$ has nonnegative spectrum, then its expectation
should also be nonnegative.
Normalization,
$\rho(\one) = 1$,
means that the 
expectation of $\one$ is 1.


\subsubsection{Quantum Mechanics}
We have already seen that in our garden-variety quantum mechanics, every
self-adjoint operator corresponds to an observable, 
so that $\CA = \CB(\CH)$.
Any unit vector $\psi \in \CH$ gives rise to a pure state
$\rho_{\psi}(A) = \inp{\psi}{A\psi}$, 
but it turns out that not all states are pure.
%
\begin{theorem}
If $\CH$ is finite dimensional, then any state 
\mbox{$\rho : \CB(\CH) \rightarrow \mathbb{C}$} can be written 
$$
\rho(A) = \tr(RA)
$$
for a unique positive operator $R$ with $\tr(R) = 1$, called the
\emph{density matrix}. 
\end{theorem}
\proof
Because 
$(A,B)\mapsto \tr(AB)$ is a perfect pairing, any linear functional 
$\rho: \CB(\CH) \rightarrow \mathbb{C}$ can be uniquely written as
$\rho(A) = \tr(RA)$. If $E = \ketbra{\psi}$ is the projection onto
a 1-dimensional subspace $\mathbb{C}\psi$, then
$\rho(E) = \inp{\psi}{R\psi}$. Positivity of $\rho$ therefore implies 
$R \geq 0$, and normalization corresponds to $\tr(R) = 1 $. \qed

\noindent The pure states $\rho_{\psi}$ correspond to 1-dimensional projections
$R = \ketbra{\psi}$. They are the extreme points of the convex
state space $\CS(\CB(\CH))$. 

Indeed, if we choose an orthonormal basis $\psi_i$ of eigenvectors
of $R$, then we can write
$R = \sum_{i=1}^{n} p_i \ketbra{\psi_i}$.
Positivity implies $p_{i} \geq 0$, and normalization implies 
$\sum_{i=1}^{n} p_{i} = 1$. 
Any state on $\CB(\CH)$ 
can therefore be written
as a convex combination of vector states
$\rho = \sum_{i} p_i \rho_{\psi_i}$, albeit not necessarily in a unique fashion. 

This invites us to interpret a system in state $\rho$ as 
one that is with probability $p_i$ in 
state $\rho_{\psi_i}$.
We will do so, and even elevate this to the more general
\emph{stochastic equivalence principle}, which holds for
convex linear combinations of any two states 
$\rho_1$ and $\rho_2$ on any system $\CA$.
\begin{center}
\emph{
A system in state $p \rho_1 + (1-p)\rho_2$ cannot be distinguished from one
that is in state $\rho_1$ with probability $p$, and in state $\rho_2$
with probability $1-p$.
}
\end{center}


\subsubsection{Classical Probability Space}
Let $(\Omega,\mathbb{\mathbb{P}})$ be a finite probability space.
A \emph{random variable} is a map $f : \Omega \rightarrow \mathbb{C}$.
It makes sense to identify two random variables $f$ and $f'$ if 
they are equal almost surely,
$\mathbb{P}([f \neq f']) = 0$,
and we write $L^{\infty}(\Omega,\mathbb{P})$ for the resulting algebra. 
The expectation defines a linear map 
$\mathbb{E} : L^{\infty}(\Omega,\mathbb{P}) \rightarrow \mathbb{C}$.

The Hilbert space $L^{2}(\Omega,\mathbb{P})$ is the same vector space, 
equipped with the inner product 
$\inp{\phi}{\psi} = \mathbb{E}(\bar{\phi}\psi)$.
%
Because $L^{\infty}(\Omega,\mathbb{P})$ acts on $L^{2}(\Omega,\mathbb{P})$ 
by multiplication, we can
consider it as a commutative 
$*$-subalgebra of $\CB(L^{2}(\Omega,\mathbb{P}))$.
This places us in the framework of quantum probability theory, with 
$\CA = L^{\infty}(\Omega,\mathbb{P})$,
$\CH = L^{2}(\Omega,\mathbb{P})$, and $\rho = \mathbb{E}$.
Let us describe the state space $\CS(\CA)$.


\begin{theorem}
Any finite probability space $(\Omega,\mathbb{P})$ gives rise to 
a finite dimensional Hibert space
$L^{2}(\Omega,\mathbb{P})$, 
and a $*$-subalgebra $L^{\infty}(\Omega,\mathbb{P})$ of 
$\CB(L^{2}(\Omega,\mathbb{P}))$. 
If $\mathbb{P}'$ is absolutely continuous w.r.t. $\mathbb{P}$, then
it defines a state $\mathbb{E}' : 
L^{\infty}(\Omega,\mathbb{P}) \rightarrow \mathbb{C}$.
All states arise in this way.  
\end{theorem}
\proof The set 
$\{\delta_{\omega} \, ; \, \omega \in \Omega,\, \mathbb{P}(\omega) \neq 0 \}$ 
constitutes a basis of $L^{\infty}(\Omega,\mathbb{P})$, so that 
any state $\rho$ can be written
$\rho = \sum p'_{\omega} \delta_{\omega}^{*}$.
Positivity and normalization amount to $p'_{\omega} \geq 0$ and
$\sum p'_{\omega} = 1$ respectively. The numbers $p'_{\omega}$ 
thus define a probability measure $\mathbb{P}' \ll \mathbb{P}$,
and $\rho = \mathbb{E}'$ is its associated expectation. \qed
 



\subsubsection{Operations}
We wish to describe transformations of a system
$\CA \subset \CB(\CH)$ to a system $\CB \subset \CB(\CH')$.
In the context of quantum probability theory, these are called
\emph{operations}.
Each operation is determined by a map
$\tau : \CS(\CA) \rightarrow \CS(\CB)$.
We investigate which maps can reasonably be interpreted as
operations.

We interpret the convex combination $p \rho_1 + (1-p)\rho_2 \in \CS(\CA)$ 
as a system that is in state $\rho_1$ with probability $p$, and in 
state $\rho_2$ with probability $(1-p)$. This means that after the
operation, the system $\CB$ must be described as being
in state $\tau(\rho_1)$ with 
probability $p$, and in state $\tau(\rho_2)$ with probability
$1-p$. Thus
$$
\tau(p \rho_1 + (1-p)\rho_2) = p \tau(\rho_1) + (1-p)\tau(\rho_2)\,.
$$  
This means that $\tau$ extends to a \emph{linear} map 
$\tau : \CA^* \rightarrow \CB^*$. We can consider it as 
the dual $\tau = T^{*}$ of a linear map $T : \CB \rightarrow \CA$, and
will do so from now on.
The requirement that $T^*$ map states to states implies that
$T(B) \geq 0$ if $B \geq 0$ (it preserves positivity), and that
$T(\one) = 1$ (normalization).

\begin{basdef}
Let $\CA$ and $\CB$ be $*$-subalgebras of $\CB(\CH)$ and $\CB(\CH')$
respectively. A map $T : \CB \rightarrow \CA$ is called \emph{positive}
if it is linear, if $B \geq 0$ implies $T(B) \geq 0$, and
if $T(\one) = 1$.
\end{basdef}
One would be tempted to interpret each positive map as an operation,
but 
this is \emph{not} what one should do. 
There is one more requirement that we must impose on $T$
in order for it to qualify as an operation.

If $T : \CB \rightarrow \CA$ describes an operation from 
$\CA$ to $\CB$, then
the map 
$T \otimes \id_n : \CB \otimes \CB(\mathbb{C}^{n}) 
\rightarrow \CA \otimes \CB(\mathbb{C}^{n})$
describes the act of transforming $\CA$ to $\CB$, while doing absolutely
nothing on the system $\CB(\mathbb{C}^{n})$.
Surely this should map states on $\CA \otimes \CB(\mathbb{C}^{n})$
to states on 
$\CB \otimes \CB(\mathbb{C}^{n})$,
and we therefore require that $T \otimes \id_n$ be a positive map for all 
$n \in \mathbb{N}$. 
\begin{basdef}
A map $T : \CB \rightarrow \CA$ is called \emph{completely positive}
if 
$$
T \otimes \id_n : \CB \otimes \CB(\mathbb{C}^{n}) 
\rightarrow \CA \otimes \CB(\mathbb{C}^{n})
$$
is positive for all $n \in \mathbb{N}$.
\end{basdef}
Perhaps surprisingly, this is not automatic. There exist maps which are
positive, but not completely positive.
We shall model operations from  $\CA$ to $\CB$ by completely positive
maps $T : \CB \rightarrow \CA$. (CP-maps for short.)

\paragraph{Operations in Classical Probability}
Let $(\Omega, \mathbb{P})$ and $(\Omega', \mathbb{P}')$ be 
finite probability spaces. Then a \emph{transition kernel}
is a positive function $P$ on $\Omega' \times \Omega$
such that $\sum_{\Omega'} P(\omega', \omega) \mathbb{P}'(\omega') = 1$
for all $\omega \in \Omega$. One should think of 
$P(\omega', \omega) \mathbb{P}'(\omega')$ as the probability 
that a system in state $\omega$ makes a transition to $\omega'$. 

An operation $T : L^{\infty}(\Omega' , \mathbb{P}') 
\rightarrow L^{\infty}(\Omega , \mathbb{P})$ is related to such
a kernel. Indeed, if we write
$Tf (\omega) = \sum_{\Omega'}  P(\omega', \omega) 
f(\omega') \mathbb{P}'(\omega')$, then 
positivity of $T$ entails
$P(\omega', \omega)\geq 0$, 
and normalization $T(\one) = \one$ translates to 
$\sum_{\Omega'} P(\omega', \omega) \mathbb{P}'(\omega') = 1$. 
An operation is precisely an equivalence class of kernels,
if we deem 
two kernels
to be equivalent when 
their difference is nonzero only on those $(\omega' , \omega)$
with
either $\mathbb{P}(\omega) = 0$ or $\mathbb{P}'(\omega') = 0$.


%
%
\paragraph{Operations from Quantum to Classical Systems} 
A completely positive map 
$T : L^{\infty}(\Omega,\mathbb{P}) \rightarrow \CB(\CH)$
is determined by its values on $\delta_{\omega}$, i.e. by the map
$E : \Omega \rightarrow \CB(\CH)$ defined as $E(\omega) = T(\delta_{\omega})$.
\begin{basdef} \label{povem}
A \emph{positive operator valued measure} or POVM is a map
$E : \Omega \rightarrow \CB(\CH)$ that satisfies
\begin{itemize}
\item[-]$E(\omega) \geq 0$ for all $\omega \in \Omega$
\item[-]$\sum_{\Omega} E(\omega) = \one$.
\end{itemize}
\end{basdef}
The map $E(\omega) = T(\delta_{\omega})$ is a POVM;  
positivity of $T$ implies the first requirement, 
and normalization the second.
%
Conversely, any POVM that is measurable w.r.t. $\mathbb{P}$
determines a CP-map $T : L^{\infty}(\Omega,\mathbb{P}) \rightarrow \CB(\CH)$.
A CP-map 
is 
a $*$-homomorphism if and only if its POVM is 
projection valued.

\paragraph{Operations between Quantum Systems}
Any $*$-homomorphism of algebras is a CP-map. In particular,
the transformations $\alpha_{U}(A) = U^{\dagger} A U$ 
of $\CB(\CH)$ are operations in the sense of quantum probability.

Another type of operation is 
considering a subsystem 
as part of a bigger whole. 
If $\CH$ is a Hilbert subspace of $\CK$, then it comes 
with the inclusion $V : \CH \hookrightarrow \CK$, 
and the adjoint projection
$V^{\dagger} : \CK \twoheadrightarrow \CH$. The 
CP-map $T : \CB(\CK) \rightarrow \CB(\CH)$
corresponding to this inclusion
is given by $T(A) = V^{\dagger} A V$. 

According to the following classification theorem, 
any CP-map is a combination of a $*$-homomorphism
and an inclusion.
\begin{theorem}[Stinespring]
Let $\CH$ and $\CH'$ be
Hilbert spaces of finite dimension, 
let $\CA$ be a $*$-subalgebra of $\CB(\CH')$, and 
let $T$ be a 
CP-map $\CA \rightarrow \CB(\CH)$.
Then there exists a finite dimensional Hilbert space $\CK$, 
an inclusion $V : \CH \hookrightarrow \CK$ and a $*$-homomorphism
$\pi : \CA \rightarrow \CB(\CK)$ such that $T(A) = V^{\dagger}\pi(A)V$.
\begin{center}
\begin{tikzpicture}
\node (linksboven) at (0, 1.3, 0) {$\CB(\CK)$};
\node (linksonder) at (0, 0, 0) {$\CA$};
\node (rechtsonder) at (1.7,0,0) {$\CB(\CH)$};
\draw [->] (linksboven) -- node[above] {$\quad \quad \quad V^{\dagger} \,\cdot\,V$}(rechtsonder);
\draw [->] (linksonder) -- node[above]{$T$}(rechtsonder);
\draw [->] (linksonder) -- node[left]{$\pi$} (linksboven);
\end{tikzpicture}
\end{center}
\end{theorem}
\proof
Let $\CK_{0} = \CA \otimes \CH$. Then the inclusion $V_0 : \CH \rightarrow \CK_{0}$
is given by $\psi \mapsto \one \otimes \psi$, and if $A \in \CA$, then 
$\pi_0(A) :\CK_{0} \rightarrow \CK_0$ is defined by
$B \otimes \psi \mapsto AB \otimes \psi$.
We equip $\CK_{0}$ with the sesquilinear form defined by 
$\inp{B \otimes \psi}{B' \otimes \psi'} = 
\inp{\psi}{T(B^{\dagger}B') \psi'}_{\CH}$.

Let $F = \sum_{i=1}^{k} B_{i} \otimes \psi_{i}$ be 
an arbitrary element of $K_{0}$.
In terms of the map 
$T \otimes \id_{k} : 
\CA \otimes \CB(\mathbb{C}^{k}) \rightarrow 
\CB(\CH \otimes \mathbb{C}^{k})
$, we can write
\begin{equation}\label{spronkje}
\inp{F}{F} = 
\left\langle
\Psi
\,|\,
T \otimes \id_{k} 
\left(
\underline{B}^{\dagger}\underline{B}
\right) 
\Psi
\right\rangle\,,
\end{equation}
with $\Psi = \sum_{i=1}^{k} \psi_i \otimes e_i$ in $\CH \otimes \mathbb{C}^{k}$
and $\underline{B} = \sum_{i=1}^{k} B_i \otimes \ket{e_1}\bra{e_i}$
in $\CA \otimes \CB(\mathbb{C}^{k})$.
Because $T\otimes \id_{k}$ is a positive map, we conclude that 
$\inp{F}{F}\geq 0$. Our sesquilinear form is positive semidefinite.

We can therefore define the Hilbert space
$\CK = \CK_{0}/\CN$, with 
$\CN$
the kernel
of $\inp{\,\cdot\,}{\,\cdot\,}$.
The map $V_0$ is an isometry w.r.t. $\inp{\,\cdot\,}{\,\cdot\,}$,  
and therefore
induces an inclusion 
$V : \CH \hookrightarrow \CK$ of Hilbert spaces. Its adjoint projection
is given by
$V^{\dagger} ( B \otimes \psi) = T(B) \psi$.

The map
$\pi_0 : \CA \rightarrow \CB(\CK_{0})$ is a contraction.
Indeed, if we replace $F$ by $\pi_{0}(A)F$ in equation
(\ref{spronkje}),
then on the right hand side, $\underline{B}^{\dagger}\underline{B}$
changes into 
the expression 
$\underline{B}^{\dagger} (A^{\dagger} A\otimes \one)\underline{B}$. 
Since this is dominated by
$
\|A\|^2 \underline{B}^{\dagger}\underline{B}$, we have 
$\inp{\pi_{0}(A)F}{\pi_{0}(A)F} \leq \|A\|^2 \inp{F}{F}$, and the
operator norm of
$\pi_0(A)$ does not exceed that of $A$.

In particular, 
$\CN$ is invariant under $\pi_0(A)$, so that we have
a map  
$\pi(A) : \CK \rightarrow \CK$.
It is then straightforward to check that $\pi$ is a 
$*$-homomorphism, and that $V^{\dagger}\pi(A)V = T(A)$.
\qed

%
%


\subsection{QPT on General Hilbert Spaces}
We define quantum probability theory in the infinite dimensional setting, which requires 
a little bit of functional analysis.

\subsubsection{Von Neumann Algebras}
We have seen that in quantum mechanics, all the structure of a system is 
encoded in its algebra $\CB(\CH)$ of bounded operators.  
A subalgebra of $\CB(\CH)$ is called a \emph{von Neumann algebra} if it is 
closed not only under the adjoint $A \mapsto A^{\dagger}$, but also 
in the weak topology.
In quantum probability theory, a system is modelled by a 
von Neumann algebra $\CA$.

A self-adjoint, possibly unbounded operator $A$ on $\CH$ 
is said to be \emph{affiliated} to $\CA$
if all its spectral projections are in $\CA$.
Observables of the system are modelled by
self-adjoint operators affiliated to $\CA$.

If $\CA$ and $\CB$ are von Neumann algebras, then we will
denote by $\CA \otimes \CB$ the von Neumann algebra generated
by the algebraic tensor product.

\subsubsection{Normal States}

A state is again described by a 
normalized positive linear map $\rho : \CA \rightarrow \mathbb{C}$, 
but it is wise to impose a continuity condition.
If $A$ is affiliated to $\CA$, 
then we would like its spectral measure
$E : \spec(A)\rightarrow \CA$ to induce a probability measure
$\mathbb{P}_{\rho} = \rho \circ E$ on $\spec(A)$.
But although $E$ is $\sigma$-additive, i.e.
$\sum_{\mathbb{N}}E(S_i) = E(\cup_{\mathbb{N}}S_i)$ in the weak topology,
this need not hold for $\mathbb{P}_{\rho} = \rho \circ E$, unless we  
impose the following continuity requirement.  
\begin{basdef}
A state $\rho : \CA \rightarrow \mathbb{C}$ is called \emph{normal}
if it is weakly continuous on the unit ball $\{A \in \CA \,; \, \|A\| \leq 1\}$. 
\end{basdef}
A normal state $\rho \in \CS(\CA)$ induces a probability measure $\mathbb{P}_{\rho}$ on the spectrum of 
any self-adjoint operator affiliated to $\CA$, so that we have our beloved map
$(\rho,A) \mapsto \mathbb{P}_{\rho}$. 
With respect to $\mathbb{P}_{\rho}$, we have of course
$\mathbb{E}(\lambda) = \rho(A)$.

\paragraph{Normal States in Quantum Mechanics}
The normal states of an ordinary quantum system $\CB(\CH)$ can again be described 
by a density matrix.
\begin{theorem}\label{sporenstel}
Let $\CH$ be a separable Hilbert space. Then
any normal state $\rho$ on $\CB(\CH)$ can be written
$$
\rho(A) = \tr(RA)
$$
for a unique positive trace-class operator $R$ on $\CH$ with $\tr(R)=1$.
\end{theorem}
\proof See for example chapter 7 of \cite{KaRiII}.\qed

\noindent A positive normalized trace-class operator admits a basis of eigenvectors $\psi_i$ 
with nonnegative eigenvalues $p_i$ summing to $1$. 
In other words, 
any state can be written as a countably infinite convex
combination of pure states,
$\rho = \sum_{i=1}^{\infty} p_i \rho_{\psi_i}$.

\paragraph{Normal States in Classical Probability}
A classical probability space, denoted $(\Omega , \Sigma , \mathbb{P})$, 
gives rise to the algebra $L^{\infty}(\Omega , \Sigma, \mathbb{P} )$ of
essentially bounded random variables on $\Omega$, 
i.e. bounded measurable functions 
up to the equivalence $\sim$,
with $f \sim g$ if $f = g$ almost surely.

Since $L^{\infty}(\Omega , \Sigma, \mathbb{P} )$ 
acts on the Hilbert space $L^2(\Omega , \Sigma , \mathbb{P})$
by multiplication, we can regard it as a 
commutative $*$-subalgebra of $\CB(L^2(\Omega , \Sigma , \mathbb{P}))$,
which turns out to be weakly closed.
\begin{theorem}\label{comvn}
The algebra
$L^{\infty}(\Omega , \Sigma, \mathbb{P})$ is a commutative von Neumann subalgebra of
$\CB(L^2(\Omega , \Sigma , \mathbb{P}))$. 
If $\mathbb{P}'$ is absolutely continuous w.r.t. $\mathbb{P}$, then the expectation
$\mathbb{E}'(f) = \int_{\Omega} f(\omega) \mathbb{P}'(d\omega)$ constitutes
a normal state on $L^{\infty}(\Omega , \Sigma, \mathbb{P})$, and
every normal state arises in this way.
\end{theorem}
\proof See e.g.\ \cite{Hans2} and the references therein. \qed


\noindent In short, a classical probability space $(\Omega,\Sigma,\mathbb{P})$ gives rise to a 
commutative von Neumann algebra, and a 
normal state corresponds to a probability measure $\mathbb{P}' \ll \mathbb{P}$.
Essentially every\footnote{ 
To be more precise, a state $\rho$ is called \emph{faithful} if 
$A \geq 0$ and $\rho(A) = 0$ implies $A=0$.  
Any commutative von Neumann algebra $\CA$ with a faithful normal state $\rho$ 
is isomorphic to $L^{\infty}(\Omega, \Sigma, \mathbb{P})$, and
the state $\rho$ corresponds to the expectation w.r.t. $\mathbb{P}$.
}
commutative von Neumann algebra arises in this way. 
%

\subsubsection{Operations}
The definition of a completely positive map remains unaltered.
A linear map $T : \CB \rightarrow \CA$ between von Neumann algebras is 
called 
positive if $T(\one) = \one$ and if
$B \geq 0$ implies $T(B)\geq 0$. It is called 
completely positive
if $T \otimes \id_n : 
\CB \otimes \CB(\mathbb{C}^{n}) \rightarrow \CA \otimes \CB(\mathbb{C}^{n})$ 
is positive for all $n \in \mathbb{N}$.
This already implies that 
$T \otimes \id : \CB \otimes \CB(\CH) \rightarrow \CA \otimes \CB(\CH)$
is positive for arbitrary $\CH$, so that we need not
add this as an extra requirement.


We do however require $T^*$ to map normal states to normal states.
We therefore model operations by 
\emph{weakly continuous} CP-maps.

\paragraph{Operations in Classical Probability}
Let $(\Omega,\Sigma,\mathbb{P})$
and 
$(\Omega',\Sigma',\mathbb{P}')$
be two probability spaces. A \emph{transition kernel} is 
a positive measurable function $P$ on $\Omega' \times \Omega$ 
such that $\int_{\Omega'} P(\omega', \omega) \mathbb{P}'(d\omega') = 1$
for all $\omega \in \Omega$. 
We will identify two transition kernels if 
they are equal almost surely w.r.t. $\mathbb{P}'\times \mathbb{P}$.
One should think of 
$P(\omega', \omega) \mathbb{P}'(d\omega')$ as the probability
that the system ends up in $d\omega'$, provided that it starts in $\omega$.

Completely positive maps between classical probability spaces
correspond precisely to transition kernels.
\begin{theorem}
A transition kernel $P$ gives rise to the weakly continuous
CP-map $T_P : L^{\infty}(\Omega',\Sigma',\mathbb{P}')
\rightarrow  L^{\infty}(\Omega,\Sigma,\mathbb{P})$
by $T_P(f) (\omega) = \int_{\Omega'} f(\omega')p(\omega,\omega')
\mathbb{P}'(d\omega)$. 
Every weakly continuous CP-map $T : L^{\infty}(\Omega',\Sigma',\mathbb{P}')
\rightarrow  L^{\infty}(\Omega,\Sigma,\mathbb{P})$
is of this form.
\end{theorem}
\proof
We prove the second statement. The map $T$
gives rise to the normal state $\rho$ on
$L^{\infty}(\Omega',\Sigma',\mathbb{P}')\otimes 
L^{\infty}(\Omega,\Sigma,\mathbb{P})$
defined by $\rho(f' \otimes f) = \mathbb{E}_{\mathbb{P}}(f T(f'))$.
According to theorem \ref{comvn}, this corresponds to a 
probability measure on $\Omega' \times \Omega$ which is
absolutely continuous w.r.t. $\mathbb{P}' \times \mathbb{P}$.
Its Radon-Nikodym derivative, well defined up to equivalence,
is a transition kernel $P$. The requirement
$\int_{\Omega'} P(\omega', \omega) \mathbb{P}'(d\omega') = 1$
corresponds to the fact that
$\mathbb{P}$ is the marginal probability distribution
on $\Omega$, as is clear from $\rho(\one \otimes f) = \mathbb{E}_{\mathbb{P}}(f)$.  
\qed

\paragraph{Operations from Quantum to Classical Systems}
We have seen that a projection valued measure $E : \Sigma \rightarrow \CB(\CH)$
gives rise to a $*$-algebra homomorphism $I_E : L^{\infty}(\Omega,\Sigma) \rightarrow \CB(\CH)$, by
way of
\begin{equation}\label{oetkhijpgu}
I_E(f) := \int_{\Omega} f(\omega) E(d\omega)\,.
\end{equation}
If $E \ll \mathbb{P}$, this factors through a weakly continuous homomorphism 
$I_E : L^{\infty}(\Omega,\Sigma, \mathbb{P}) \rightarrow \CB(\CH)$
of von Neumann algebras.
Conversely, any weakly continuous homomorphism
$T : L^{\infty}(\Omega,\Sigma, \mathbb{P}) \rightarrow \CB(\CH)$
of von Neumann algebras
gives rise to a PVM by way of $E(S) = T(\one_S)$. 

In order to capture all CP-maps, not just the homomorphisms, we introduce
the notion of a Positive Operator Valued Measure or POVM 
(cf.~def.~\ref{povem}).
\begin{basdef}
Let $(\Omega,\Sigma)$ be a measurable space, and $\CH$ a Hilbert space.
Then a 
POVM is a map
$E : \Sigma \rightarrow \CB(\CH)$ with the following properties.
\begin{itemize}
	\item[-]$E(S)\geq 0$ for all $S \in \Sigma$.
	\item[-]$E(\Omega) = \one$.
	\item[-]If $S_i$ is a sequence of disjoint sets in $\Sigma$, then 
				  $\sum_{i=1}^{\infty} E(S_i) = E(\bigcup_{i=1}^{\infty} S_i)$
				  in the weak topology.
\end{itemize}
\end{basdef}
A weakly continuous CP-map 
$T : L^{\infty}(\Omega,\Sigma, \mathbb{P}) \rightarrow \CB(\CH)$
corresponds to a POVM $E \ll \mathbb{P}$ by 
$T(f) = \int_{\Omega}f(\omega)E(d\omega)$. The POVM
is a PVM (cf.~def.~\ref{pvem}) if and only if $T$ is a homomorphism.  

\paragraph{Operations in Quantum Systems}
Although
every homomorphism of von Neumann algebras is a CP-map, 
it is certainly not true that every CP-map is a homomorphism.
We do however have the following theorem.
\begin{theorem}\label{geelsnuitje} 
Every CP-map $T : \CB \rightarrow \CA$ with a completely positive inverse
is an isomorphism.
\end{theorem}\vspace{-1mm}
See page \pageref{huuuk} for the proof.
The group of invertible operations of a system
$\CA$ is thus simply $\mathrm{Aut}^{\mathrm{ct}}(\CA)$,
the group of weakly continuous automorphisms.
Its action on the state space $\CS(\CA)$ equips it  
with the topology of pointwise 
convergence in trace distance. (See page \pageref{loosspoor}.)
From the point of view of quantum probability theory, 
the action of a Lie group $G$ should therefore be 
described by a continuous homomorphism
$\pi : G \rightarrow \mathrm{Aut}^{\mathrm{ct}}(\CA)$.
The following theorem shows that
$\mathrm{Aut}^{\mathrm{ct}}(\CB(\CH)) = \mathcal{P}U(\CH)$.
\begin{theorem}[Wigner]
Let $\CH$ be a separable Hilbert space. Then every
weakly continuous 
automorphism $\alpha$ of $\CB(\CH)$ is inner, and
can be written
$\alpha (A) = U^{\dagger}AU$ for some unitary $U$. 
\end{theorem}\vspace{-2mm}
\proof
Because $\alpha^{*}$ maps normal states to normal states, 
we can consider it as a
continuous trace-preserving
$\mathbb{C}$-linear map on the  
trace-class operators $\mathcal{T}(\CH)$,
cf.\ theorem \ref{sporenstel}.
Because
$\alpha^{*}$ maps $R$ to $\alpha^{-1}(R)$, 
it restricts to a map on 
the pure state space $\mathcal{P}(\CH)$,
and $\alpha$ is determined 
by this restriction because of continuity.
The map on $\mathcal{P}(\CH)$ preserves the pairing 
$(\mathbb{C}\phi , \mathbb{C}\psi) = |\inp{\phi}{\psi}|/\|\phi\|\|\psi\|$
because it preserves the trace.
According to Wigner's classification of 
the automorphisms of $\mathcal{P}(\CH)$, any such map must be given by
$\alpha^{*} (\ketbra{\psi}) = U \ketbra {\psi} U^{\dagger}$ with $U$
either unitary or anti-unitary. The map $U \, \cdot \, U^{\dagger}$ 
can only be extended to a 
$\mathbb{C}$-linear map on $\mathcal{T}(\CH)$ if $U$ is unitary, 
or of course if $\CH$ is one-dimensional, in which case the theorem is
trivially true. 
\qed
%
%
%
%
An action of $G$ by weakly continuous CP-maps
on $\CB(\CH)$ is therefore precisely a projective unitary representation. 
In particular, \emph{invertible} time evolution on a \emph{closed} quantum system 
is still
described by a Hamiltonian, and we do not gain anything new. 

\emph{Non-invertible} time evolution on an \emph{open} quantum system however should 
be described by a continuous homomorphism of the semigroup 
$(\mathbb{R}^{\geq 0} , +)$ into
the semigroup of weakly continuous CP-maps on $\CB(\CH)$, which will in 
general \emph{not} be automorphic. 
CP-maps are classified by the following theorem.
\begin{theorem}[Stinespring]
Let $\CA$ be a von Neumann algebra, $\CH$ a Hilbert space, and 
let $T$ be a 
CP-map $\CA \rightarrow \CB(\CH)$.
Then there exists a Hilbert space $\CK$, 
an inclusion $V : \CH \hookrightarrow \CK$ and a $*$-homomorphism
$\pi : \CA \rightarrow \CB(\CK)$ such that $T(A) = V^{\dagger}\pi(A)V$.
\vspace{-5mm}
\begin{center}
\begin{tikzpicture}
\node (linksboven) at (0, 1.3, 0) {$\CB(\CK)$};
\node (linksonder) at (0, 0, 0) {$\CA$};
\node (rechtsonder) at (1.7,0,0) {$\CB(\CH)$};
\draw [->] (linksboven) -- node[above] {$\quad \quad \quad V^{\dagger} \,\cdot\,V$}(rechtsonder);
\draw [->] (linksonder) -- node[above]{$T$}(rechtsonder);
\draw [->] (linksonder) -- node[left]{$\pi$} (linksboven);
\end{tikzpicture}
\end{center}
\vspace{-2mm}
If we impose that $\CK$ be the closure of $\pi(\CA)\CH$, then
$(\pi,V,\CK)$ is unique up to unitary transformation, and
$T$ is weakly continuous if and only if $\pi$ is.
\end{theorem}
\proof See Stinespring's paper \cite{St}.\qed
\subsubsection{Summarizing}
We arrive at the following picture of quantum probability theory.
A system is modelled by a von Neumann algebra $\CA$
of bounded operators on a Hilbert space $\CH$. 
An observable is modelled by a self-adjoint operator affiliated to $\CA$.
A state on this system is described by a 
normalized positive linear map $\rho : \CA \rightarrow \mathbb{C}$
that is weakly continuous on the unit ball,
and
an operation from $\CA$ to $\CB$ is described by
a weakly continuous completely positive map $T : \CB \rightarrow \CA$.

\section{Quantum Stochastic Differential Equations}\label{introqsdesec}

This concludes our description of quantum probability spaces in 
general. 
We will now focus on a single model that will be
of particular relevance to this thesis: 
the second quantized electromagnetic field.
Its Hilbert space $\CF(L^{2}(\mathbb{R}^{+}))$
is the Fock space of
$L^{2}(\mathbb{R}^{+})$.


\subsection{Fock Space}\label{focksec}

%
%
If $\CH$ is a Hilbert space, then so is
its $n$-fold symmetric tensor power $\CH^{\otimes_{s}n}$. 
We set $\CH^{\otimes_{s}0} := \mathbb{C}$, and 
define the Fock space over $\CH$ as follows.
\begin{basdef}
The (symmetric) \emph{Fock space} $\CF(\CH)$ over a Hilbert space $\CH$
is defined as
$$
\CF(\CH) = \bigoplus_{n=0}^{\infty} \CH^{\otimes_{s}n} \,,
$$
the Hilbert space of sequences $n \mapsto v_n$ with 
$v_n \in \CH^{\otimes_{s}n}$
and $\sum_n \inp{v_n}{v_n} < \infty$.
\end{basdef}
In many respects, taking $\CH$ into $\CF(\CH)$ can be seen as exponentiating
a Hilbert space.
For example, the set of \emph{exponential vectors} 
$\{\vexp(v) \,;\, v \in \CH\}$, with
$$
\vexp(v) := \bigoplus_{n=0}^{\infty} \frac{1}{\sqrt{n!}} \, v^{\otimes n} \,,
$$
is linearly independent 
and dense in 
$\CF(\CH)$ (see \cite{Par92}).
They are easily seen to satisfy 
$$
\inp{\vexp(v)}{\vexp(w)}_{\CF(\CH)} = \exp(\inp{v}{w}_{\CH}).
$$
We will denote the vectors $v^{\otimes n}$ by $\ket{v}_{n}$. 
The normalized exponential vectors are called
\emph{coherent states}. We denote them by
$\ket{v} := e^{-\frac{1}{2}\inp{v}{v}} \vexp(v)$.

\begin{proposition}
There is a 
natural isomorphism of Hilbert spaces 
$$
\CF(\CH_1 \oplus \CH_2) \simeq \CF(\CH_1) \otimes \CF(\CH_2)\,,
$$
under which $\vexp(v_1 \oplus v_2) \simeq \vexp(v_1)\otimes \vexp(v_2)$.

\end{proposition}
\proof
Since the exponential vectors are dense and linearly independent,
the expression 
$\vexp(v_1 \oplus v_2) \mapsto \vexp(v_1)\otimes \vexp(v_2)$
yields a densely defined linear map.
Since both $\inp{\vexp(v_1 \oplus v_2)}{\vexp(w_1 \oplus w_2)}$ 
and $\inp{\vexp(v_1)\otimes \vexp(v_2)}{\vexp(w_1)\otimes \vexp(w_2)}$
equal $\exp(\inp{v_1}{w_1} + \inp{v_2}{w_2})$, it extends to a unitary
isomorphism.
\qed
By the same token, it is clear that an isometry $U : \CH_1 \rightarrow \CH_2$
induces the isometry $\CF(U) : \CF(\CH_1) \rightarrow \CF(\CH_2)$ defined by 
$\CF(U) \vexp(v) = \vexp(Uv)$. Because $\CF(U' \circ U) = \CF(U')\circ \CF(U)$, 
one can view $\CF$ as a functor from the category of Hilbert spaces
to itself, taking direct sums into tensor products\footnote{
In the words of E.~Nelson: 
`Quantization is a mystery, but second quantization is a functor.'
}. 
In particular, $\CF(\CH)$ carries a unitary representation of $U(\CH)$.

\subsection{Weyl Operators, Fields and Momenta}
The unitary representation $U\mapsto \CF(U)$ of $U(\CH)$ extends to a 
\emph{projective} unitary representation of the group $U(\CH)\ltimes \CH$ 
with product \mbox{$(U,v)(U'\!,v') \!=\! (UU'\!, Uv'\!+v)$.}
\begin{proposition}
The \emph{Weyl operators} $W(U,v)$, defined by 
\begin{equation}\label{WeylOp}
W(U,v)\vexp(w) = \exp(-\half \|v\|^{2} - \inp{v}{Uw}) \vexp(Uw + v)\,,
\end{equation}
constitute a projective unitary representation of $U(\CH) \ltimes \CH$, satisfying
\begin{equation}\label{WeylCoc}
W(U , v)W(U' , v') = \exp(-i \basim \inp{v}{U v'})W(UU' , Uv' + v )\,.
\end{equation}
\end{proposition}
\proof
Since $\inp{W(U,v)\vexp{w}}{W(U,v)\vexp(w')} = \inp{\vexp{w}}{\vexp{w'}}$, 
equation (\ref{WeylOp}) does indeed define a unitary operator. Equation
(\ref{WeylCoc}) is then readily verified on the dense set of exponential vectors.
Continuity essentially follows from continuity of the map
$v \mapsto \vexp(v)$, see \cite{Par92} for details.
\qed
We will denote $W(\one,v)$ by $W(v)$,
and designate $W(U,0)$ by its old name $\CF(U)$.
%
If $\CH$ is finite dimensional, we find a projective Lie algebra homomorphism 
$\dot{W} : \mathfrak{su}(\CH) \ltimes \CH \rightarrow \mathfrak{su}(\CF(\CH))$,
cf. section \ref{gramsmid}.
Rather than making this precise for infinite dimensional $\CH$, 
we will restrict attention to 
Lie subgroups of $U(\CH)\ltimes \CH$ of dimension 1 or 2.
%
\subsubsection{Associated Observables}
Every self-adjoint operator $H$ gives rise to
a one-parameter Lie group in $U(\CH) \ltimes \CH$, 
as does every vector $v \in \CH$.
The restriction of $W$ to these one-parameter groups 
is a linear unitary representation. We can therefore define
the skew-symmetric closed operators  
\begin{equation}
\dot{\CF}(H) = \CF(e^{itH})\frac{d}{dt} \CF(e^{-itH})\,
\end{equation}
and
\begin{equation}\label{gruffel}
\dot{W}(v) = W(-tv)\frac{d}{dt}W(tv)\,.
\end{equation}
The observable $\lambda(H) := i\,\dot{\CF}(H)$ is called the 
differential second quantization of $H$. 
We will not be needing commutation relations involving it. 
(See however \cite{Par92}.)
The operators $J(v) = i \dot{W}(v)$ are called  
fields and momenta.

They satisfy the \emph{canonical commutation relation}
\begin{equation}\label{CanComRel}
[J(u) , J(v)] = 2i\,\basim\inp{u}{v}\one\,,
\end{equation}
as can be seen by differentiating (\ref{WeylCoc}).
Since every exponential vector is in the domain of every product
$J(v_1) \ldots J(v_n)$, the l.h.s. of  
equation (\ref{CanComRel}) is densely defined.


\subsubsection{Fields and Momenta}
As is the case in Hamiltonian mechanics, the division of the $J(v)$
into fields and momenta is
to a certain extent arbitrary. 
If we choose a 
completely real subspace $\CH_{\mathbb{R}}$,
then we define for $f \in \CH_{\mathbb{R}}$ the \emph{fields} to be
$$
\Phi(f) = - J(f)\,,  
$$
and their \emph{conjugate momenta}
$$
\Pi(f) = -J(if)\,. 
$$
They satisfy the commutation relations 
$$
[\Phi(f) , \Phi(g)] = 0 \,, \quad
[\Phi(f) , \Pi(g)] = 2i \inp{f}{g}\,, \quad
[\Pi(f) , \Pi(g)] = 0\,.
$$

\subsubsection{Creation and Annihilation Operators}
Even without a completely real subspace, we can 
introduce for any $u \in \CH$ the creation and annihilation operators
\begin{equation}\label{crean}
a^{\dagger}(u) = - \half (J(iu) + i J(u)) , \quad 
a(u) = -\half (J(iu) - i J(u))\,,
\end{equation}
respectively, satisfying
$$
[a(u) , a(v)] = 0 \,, \quad
[a(u) , a^{\dagger}(v)] = \inp{u}{v}\one\,, \quad
[a^{\dagger}(u) , a^{\dagger}(v)] = 0\,.
$$
By sandwiching (\ref{WeylOp}) between exponential 
vectors and applying equation (\ref{gruffel}), we obtain:
\begin{eqnarray}
\inp{\vexp{v}}{a(u)\vexp(v')} &=& \inp{u}{v'}\inp{\vexp{v}}{\vexp(v')}\,,\\
\inp{\vexp{v}}{a^{\dagger}(u)\vexp(v')} &=& \inp{v}{u}\inp{\vexp{v}}{\vexp(v')}\,.
\end{eqnarray}
We will also have use for the second order correlation functions
\begin{eqnarray}
\inp{\vexp{v}}{a^{\dagger}(u_1)a^{\dagger}(u_2)\vexp(v')} &=&
\inp{v}{u_1}\inp{v}{u_2}\inp{\vexp(v)}{\vexp{v'}}\label{aadaad}\\
\inp{\vexp{v}}{a(u_1)a(u_2)\vexp(v')} &=&  
\inp{u_1}{v'}\inp{u_2}{v'} \inp{\vexp{v}}{\vexp(v')}\label{aaaa}\\
\inp{\vexp{v}}{a^{\dagger}(u_1)a(u_2)\vexp(v')} &=&
\inp{v}{u_1}\inp{u_2}{v'}\inp{\vexp{v}}{\vexp(v')}\label{aadaa}\\
\inp{\vexp{v}}{a(u_1)a^{\dagger}(u_2)\vexp(v')} &=&
\inp{u_1}{v'}\inp{v}{u_2} \inp{\vexp{v}}{\vexp(v')}\label{aaaad}\\
& & + \,\,\,\,\quad\inp{u_1}{u_2} \inp{\vexp{v}}{\vexp(v')}\,. \nonumber 
\end{eqnarray}
The seemingly innocuous hiccup at (\ref{aaaad}), essentially due to the
fact that $a$ and $a^{\dagger}$ do not commute, is what 
will eventually make the quantum stochastic calculus 
of section \ref{qsdesec} work.  

\subsection{The Harmonic Oscillator \`{a} la Fock}\label{hofock}
It is instructive to see how all of this 
works out in the simple case $\CF(\mathbb{C})$.
This describes the 
harmonic oscillator, i.e.~a particle moving in a potential
$V(x) \sim x^2$.

The $n$-fold symmetric tensor product 
$\mathbb{C}^{\otimes_{s}n}$ of $\mathbb{C}$ is one dimensional, 
and we choose the unit basis element 
$\ket{\,1\,}_{n} = 1 \otimes \ldots \otimes 1$.
We denote it $\ket{\,n\,} := \ket{\,1\,}_{n}$ for brevity.
These combine to form the
\emph{particle basis}
$\{\,\ket{\,n\,} \,;\, n \in \mathbb{N} \,\} $
of the Fock space
$$
\CF(\mathbb{C}) = \bigoplus_{n=0}^{\infty}\mathbb{C}^{\otimes_{s}n}  \,.
$$
In terms of this particle basis, we have the exponential vectors
$$
\vexp(z) = \sum_{n=0}^{\infty} \frac{z^{n}}{\sqrt{n!}} \,\ket{\,n\,}\,.
$$
Conversely, we can recover the particle basis from the exponential vectors by
\begin{equation}\label{huilen}
\ket{\,n\,} = \frac{1}{\sqrt{n!}} \frac{d^{n}}{dz^{n}}\Big|_{z=0} \vexp(z)\,.
\end{equation}

Let us investigate what the observables $\lambda(h) = i \dot{\CF}(h)$ and 
$J(v) = i\dot{W} (\xi)$
look like in the particle basis. 
(They combine to form a projective representation of
the three-dimensional real Lie algebra 
$\mathfrak{su}(\mathbb{C}) \ltimes \mathbb{C}$.)

The 1-dimensional unitary group $U(\mathbb{C})$ acts on the exponential 
vectors by
$
W(e^{-ih}) \vexp(z) = \vexp(e^{-ih} z)
$.
From equation (\ref{huilen}), we see that 
$W(e^{-ith}) \ket{\,n\,} = e^{-inth} \ket{\,n\,}$, so that
\begin{equation}\label{snuuuk}
\lambda (h) \,\ket{\,n\,} = nh \ket{\,n\,}\,.
\end{equation}
We denote the second quantized Hamiltonian by $H = \lambda(1)$.
The particle basis is a basis of eigenvectors of $H$. 

Using equation (\ref{huilen}) and the fact that 
$J(v) = i \frac{d}{dt}|_{0} W(tv)$, a short calculation 
starting from
$$
W(tv) \vexp(z) = \exp(-\half t^{2}|v|^{2} - t\bar{v}z) \vexp(z+tv)
$$
yields
$$
J(v) \, \ket{\,n\,} = iv\sqrt{n+1}\,\, \ket{\,n+1\,} - i 
\bar{v} \sqrt{n} \,\, \ket{\,n-1\,}
$$
for $n \geq 1$, and $J(v) \,\ket{\,0\,} = iv\,\ket{\,1\,}$.
From equation (\ref{crean}), we then have 
\begin{equation}\label{fruuuk}
a(1) \,\ket{\,n\,} = \sqrt{n} \,\,\ket{\,n-1\,} , 
\quad a^{\dagger}(1) 
\, \ket{\,n\,} = \sqrt{n+1} \,\,\ket{\,n+1\,}\,.
\end{equation}

From equation \ref{CanComRel}, we see that
$
[J(1) , J(i)] = 2i
$.
If we define 
$X = -J(1)$ and  $P = -J(i)/2$,
we have the canonical commutation relation $[X,P] = i$.
The creation and annihilation operators can be expressed as
$$
a(1) = P - \half iX , \quad a^{\dagger}(1) = P + \half iX\,
$$
with $[a(1) , a^{\dagger}(1)] = \one$.

We see from equations (\ref{snuuuk}) and (\ref{fruuuk}) that
$\lambda(1) = a(1)a^{\dagger}(1) - \one$, so that
$$
H =   P^{2} +  \fourth X^{2} - \half \one\,.
$$
This means that $\CF(\mathbb{C})$ serves as a Hilbert space for the harmonic
oscillator in which $H$ is diagonal w.r.t. the particle basis.
The Schr\"{o}dinger equation 
$i\frac{d}{dt} \psi(t) = H \psi(t)$
is therefore trivially solved, 
the fundamental solutions are $\psi(t) = e^{-int} \ket{\,n\,}$.

The fact that $H$ is diagonal means that
$X$ and $P$ are not, because they do not commute with $H$. 
We will now consider a representation 
in which $X$ is diagonal, but $H$ and $P$ are not.

\subsection{The Harmonic Oscillator \`{a} la Gau\ss}\label{hogauss}
Consider the Hilbert space $L^{2}(\mathbb{R} , \mathbb{P}_{G})$, 
with $\mathbb{P}_{G}$ is the Gaussian probability measure
$$
\mathbb{P}_{G}(dx) = \frac{1}{\sqrt{2\pi}} \, e^{-\frac{1}{2} x^2} dx\,.
$$
It has the characteristic function 
$\chi_{G}(z) = \mathbb{E}(e^{izx}) = e^{- \frac{1}{2}z^2}$.
Because the functions $e_{z}(x) = e^{izx}$ are in 
$L^{2}(\mathbb{R} , \mathbb{P}_{G})$ not only for 
$z \in \mathbb{R}$, but
for all $z \in \mathbb{C}$, we may consider 
$\chi_{G}(z)$ as a holomorphic function. 

The inner product $\inp{e_z}{e_w} = 
\mathbb{E}(e^{i(w- \bar{z})x})$ equals 
$\chi_{G}(w-\bar{z}) = e^{- \frac{1}{2}(w- \bar{z})^{2}}$.
The linear span of the $e_z$ is dense in $L^{2}(\mathbb{R} , \mathbb{P}_{G})$,
and if we define $E_{z} = e^{\frac{1}{2} z^2} e_{z}$, then
$\inp{E_z}{E_w} = e^{\bar{z}w}$.


Recall that similarly, the exponential vectors $\vexp(z)$ span $\CF(\mathbb{C})$,
and that $\inp{\vexp(z)}{\vexp(w)} = e^{\bar{z}w}$.
The map $\vexp(z) \mapsto E_{z}$ therefore yields an isomorphism 
$\gamma : \CF(\mathbb{C}) \stackrel{\sim}{\rightarrow} 
L^{2}(\mathbb{R} , \mathbb{P}_{G})$. 
We can now transport the Weyl operators to 
$L^{2}(\mathbb{R} , \mathbb{P}_{G})$, as well as $X$, $P$
and $H$. 

On exponential vectors, the Weyl operators are given by
$\gamma W(v) \gamma^{-1} 
E_z (x) = e^{\frac{1}{2}v(v - \bar{v})} e^{ivx} E_z(x + v - \bar{v})$.
The $E_z$ are dense, so we must have 
$\gamma W(v)\gamma^{-1} 
\psi (x) = e^{\frac{1}{2}v(v - \bar{v})} e^{ivx} \psi(x + v - \bar{v})$
for all $\psi \in L^{2}(\mathbb{R} , \mathbb{P}_{G})$.
The corresponding observables 
are then obtained from
$J(v) = i\frac{d}{dt}|_{0}\gamma W(tv) \gamma^{-1}$.
We get
$J(1) = -x$ and $J(i) = i( 2 \frac{d}{dx} - x)$, so that
$X = x$, $P = - i\frac{d}{dx} + \half i x$ and
$H = -\frac{d^{2}}{dx^{2}} + x\frac{d}{dx}$.
We conclude that indeed,
in the
Gaussian picture, $X$
is diagonal. 

Although $H$ is not diagonal here,
we can use $\gamma$ 
to transport its basis $\ket{\,n\,}$ of eigenvectors
to $L^{2}(\mathbb{R},\mathbb{P}_{G})$.
The particle basis then transforms into the Hermite polynomials
%
$$
\psi_n(x) = 
\frac{1}{\sqrt{n!}} \frac{d^{n}}{dz^{n}}\Big|_{0} e^{\frac{1}{2} z^2 + izx},
$$
which apparently satisfy 
$( -\frac{d^{2}}{dx^{2}} + x\frac{d}{dx})\psi_{n} = n \psi_{n}$.
Schr\"{o}dinger's equation 
$i \frac{d}{dt}\psi = ( -\frac{d^{2}}{dx^{2}} + x\frac{d}{dx}) \psi$
therefore
has the fundamental solutions $e^{-int}\psi_n(x)$.
\subsection{The Wiener Measure}
We wish to obtain a Gaussian representation
of the Fock space $\CF(L^{2}(\mathbb{R^{+}}))$, 
in which the fields $\Phi(f)$ for all real functions
$f$
are \emph{simultaneously} diagonalized. 
The germane probability distribution 
turns out to be the so-called Wiener measure, 
obtained from the random walk in the limit of 
small step size.

\subsubsection{The Random Walk}
The random walk is the following stochastic process.
Divide the time line $\mathbb{R}^{+}$ into equal intervals
of length $\Delta t$. At the start of each 
time interval, flip a coin to step either left or right by an amount
$\sqrt{\Delta t}$.
If we let $\Delta t$ tend to zero, then  
according to the central limit theorem, the 
position $B_t$ at time $t$ will have a Gaussian
distribution with $\mu = 0$ and $\sigma = \sqrt{t}$.

\begin{center}
\setlength{\unitlength}{1 cm}
\begin{picture}(9,4)
\put(0,1){\line(1,0){9}}
\put(0,0){\line(0,1){3.5}}
\put(-0.4,2){\mbox{\footnotesize B}}
\put(-0.4,2.5){\mbox{$\uparrow$}}
\put(4,0.6){\mbox{\footnotesize t} \quad $\longrightarrow$}
\put(0,1){\line(1,1){0.5}}
\put(0.5,1.5){\line(1,1){0.5}}
\put(1,2){\line(1,-1){0.5}}
\put(1.5,1.5){\line(1,1){0.5}}
\put(2,2){\line(1,1){0.5}}
\put(2.5,2.5){\line(1,1){0.5}}
\put(3,3){\line(1,-1){0.5}}
\put(3.5,2.5){\line(1,-1){0.5}}
\put(4,2){\line(1,1){0.5}}
\put(4.5,2.5){\line(1,-1){0.5}}
\put(5,2){\line(1,-1){0.5}}
\put(5.5,1.5){\line(1,1){0.5}}
\put(6,2){\line(1,-1){0.5}}
\put(6.5,1.5){\line(1,1){0.5}}
\put(7,2){\line(1,1){0.5}}
\put(7.5,2.5){\line(1,1){0.5}}
\end{picture}\\
{\small \fig Random walk, produced with a 0,05\EUR \,\,random bit generator.}
\end{center}

For each finite step size, the random walk places a 
probability measure $\mathbb{P}_{\Delta t}$ on the set of possible paths
starting at $0$. 
The Wiener measure $\mathbb{P}_{W}$ is 
its limit as $\Delta t$ tends to zero. 
It lives on
$C_{S}(\mathbb{R^{+}})$, the space of
continuous paths $B:  \mathbb{R}^{+} \rightarrow \mathbb{R}$
that start at $0$ and grow at most polynomially. 
We briefly sketch its construction.

\subsubsection{Measures on the Space of Tempered Distributions}

Because the increments $\Delta B_t$  
rather than the endpoints $B_T$ are independent,
it is natural 
to view $B_T$ as the `noise' $N_t = \frac{d}{dt}B_t$
integrated from $0$ to $T$, and then consider $\mathbb{P}_{W}$
as a probability measure on the set
$\Omega = \{N = \frac{d}{dt} B \,;\, B \in C_{S}(\mathbb{R}^{+})\}$ 
of noises. 

Unfortunately, not every continuous function is differentiable.
Because 
one cannot expect the erratic ($\Delta B_t \sim \sqrt{\Delta t}$) 
paths to be differentiable, we cannot solve this problem by restricting
the space of paths. 
%
We therefore define the 
noise $N_t$ as the tempered distribution 
$$
N(\phi) =  -\int_{0}^{\infty} B_t \frac{d}{dt}\phi dt\,,
$$
a continuous linear functional on
the space
$C_{S}^{\infty}(\mathbb{R}^{+})$ of rapidly decreasing functions
vanishing at $0$. 
Using highly illegal 
partial integration, we write this 
symbolically as 
$N(\phi) = \int_{0}^{\infty} \phi(t) dB_{t}$.  

Having realized our noises 
as a subset $\Omega \subset \mathcal{D}_{S}(\mathbb{R}^{+})$ of 
the tempered distributions,
we can (and will) consider the Wiener measure 
as a measure on 
$\mathcal{D}_{S}(\mathbb{R}^{+})$ 
with $\mathbb{P}_{W}(\Omega) = 1$.

Now for any probability measure $\mathbb{P}$ on $\mathcal{D}_{S}(\mathbb{R}^{+})$, 
we define the characteristic function 
$\chi : C_{S}^{\infty}(\mathbb{R}^{+}) \rightarrow \mathbb{C}$
to be 
$$
\chi(\phi) = \int_{\mathcal{D}_{S}(\mathbb{R}^{+})} 
\exp(iF(\phi))
\mathbb{P}(dF)\,.
$$
According to the following theorem, 
a probability measure is completely characterized by its characteristic 
function. (Hence the name.)
%
\begin{theorem}[Bochner-Minlos]\label{Minlos}
A functional $\chi : C_{S}^{\infty}(\mathbb{R}^{+}) \rightarrow \mathbb{C}$ 
is the
characteristic function of a probability measure on the space
$\mathcal{D}_{S}(\mathbb{R}^{+})$ of tempered distributions
if and only if it is continuous, if $\chi(0) = 1$ and 
if it is positive semidefinite, i.e. 
$\sum_{i,j} \bar{z}_i z_j\chi(\phi_i - \phi_j) \geq 0$ 
for all $z_1 \ldots z_m \in \mathbb{C}$, 
$\phi_1 \ldots \phi_m \in C_{S}^{\infty}(\mathbb{R}^{+})$. 
\end{theorem}
\proof
Any characteristic function is 
continuous, and $\chi(0) = \mathbb{P}(\mathcal{D}_{S}(\mathbb{R}^{+}))$ 
is equal to $1$.
It is also positive semidefinite, as
$$
\sum_{i,j} \bar{z}_i z_j \int_{\mathcal{D}_{S}(\mathbb{R}^{+})} 
\exp(iF(\phi_i - \phi_j)) 
\mathbb{P}(dF)
=
\int_{\mathcal{D}_{S}(\mathbb{R}^{+})} 
\big| \sum_i z_i \exp(iF (\phi_i)) \big|^{2} \mathbb{P}(dF)
$$
is certainly nonnegative. 
%
For the converse statement, which is of course the hard part,
see e.g. \cite[p.~322]{GV64}. 
\qed
Our strategy is now to
calculate the characteristic function $\chi_{\Delta t}$ 
of the random walk process with step size $\Delta t$,
and show that it tends to a limit $\chi_{W}$ as $\Delta t$
tends to zero. 
We then use theorem \ref{Minlos} to define $\mathbb{P}_{W}$
as the unique measure characterized by $\chi_{W}$.

\subsubsection{The Wiener Measure}

Calculating $\chi_{\Delta t}$ is rather straightforward. 
If we think of $\omega \in \{\pm 1\}^{\mathbb{N}}$ as an an infinite coin toss,
then 
$B_t(\omega) = \int_{0}^{t}N_s(\omega) ds$, with 
$
N(\omega) = 
\sum_{k=0}^{\infty} \frac{\omega_{k}}{\sqrt{\Delta t}}
\one_{[k\Delta t, (k+1)\Delta t)}\,.
$  
As a distribution, this reads
$$N(\phi) = 
\sum_{k=0}^{\infty} \frac{\omega_{k}}{\sqrt{\Delta t}}
\int_{k\Delta t}^{(k+1)\Delta t} \phi(s) ds\,.
$$  
Then 
\begin{eqnarray*}
\chi_{\Delta t}(\phi) &=& \lim_{n \rightarrow \infty}
2^{-n} \sum_{\omega \in \{\pm 1\}^{n}} 
\exp\left(
\sum_{k=0}^{n-1} \frac{i \omega_{k}}{\sqrt{\Delta t}}
\int_{k\Delta t}^{(k+1)\Delta t} \phi(s) ds
\right)\\
&=& \lim_{n \rightarrow \infty}\prod_{k=0}^{n-1}
\cos\left(\frac{1}{\sqrt{\Delta t}} 
\int_{k\Delta t}^{(k+1)\Delta t} \phi(s) ds \right)\,.
\end{eqnarray*}
This means that as the step size $\Delta t$ tends to zero, the characteristic 
function of the random walk process 
approaches the limit
$\chi_{\Delta t} (\phi) \rightarrow 
\exp(- \frac{1}{2} \int_{0}^{\infty} \phi^{2}(s) ds)$. 
%
%
The function 
$\chi_{W} (\phi) = \exp(- \frac{1}{2} \int_{0}^{\infty} \phi^{2}(s) ds)$
is indeed continuous and positive semidefinite, so that we may finally define 
the Wiener measure. 
\begin{basdef}\label{wieien}
The Wiener measure $\mathbb{P}_{W}$ is 
the unique probability measure
$\mathbb{P}_{W}$
on $\mathcal{D}_{S}(\mathbb{R}^{+})$
with characteristic function 
$\chi_{W} (\phi) = \exp(- \frac{1}{2} \int_{0}^{\infty} \phi^{2}(s) ds)$.
\end{basdef}
%
One can show that
the Wiener measure concentrates on the set of noises,
$\mathbb{P}_{W}(\Omega) = 1$.
It therefore induces 
a probability measure on $C_{S}(\mathbb{R}^{+})$, which we
denote again by 
$\mathbb{P}_{W}$.  
It is the unique one with the following properties.
\begin{itemize}
\item[-]
Every increment $B_{t_2} - B_{t_1}$
has a Gaussian distribution with $\mu = 0$ and 
$\sigma = \sqrt{t_2 - t_1}$.
\item[-]The increments are independent, i.e.~for any sequence 
$t_1 \leq \ldots \leq t_k$ of times,
the steps $B_{t_{i+1}} - B_{t_{i}}$ are mutually independent.
\item[-]The increments are stationary, i.e.~the 
joint probability distribution of
the steps 
$(B_{t_1} - B_{s_1}), \ldots , (B_{t_k} - B_{s_k})$ does not change
if all $t_i$ and $s_i$ are shifted by the same amount.   
\end{itemize}

\subsection{Stochastic Differential Equations}
Consider $B_s$ as a random variable on the space 
$\Omega \subset \mathcal{D}_{S}(\mathbb{R}^{+})$ of noises.
Let $\Sigma_{t}$ be the $\sigma$-algebra generated by
$\{B_s \,;\, 0 \leq s \leq t\}$. The $\sigma$-algebras
$\Sigma_{s}$ with $s \in \mathbb{R}^{+}$ filter $\Omega$, 
in the sense that $\Sigma_t$ is finer than $\Sigma_s$ if $t > s$.
 
A \emph{process} is a measurable map 
$Y : \Omega \times \mathbb{R}^{+} \rightarrow \mathbb{R}$.
A process is called \emph{adapted} 
if $Y_t$ is measurable w.r.t. $\Sigma_{t}$.
One should think of an adapted process as one in which 
$Y_t$ depends only on the paths $B_{s}$ up to time $t$.
We denote by $L^{2}(\Omega \times \mathbb{R}^{+}, \Sigma_t)$
the Hilbert space of adapted processes 
that are square integrable w.r.t. the product of
$\mathbb{P}_{W}$ and the Lebesgue measure.

An adapted step function is an adapted process of the form
$$
Y_{t}(\omega) = \sum_{k=0}^{\infty} c_{k}(\omega) \one_{[t_k , t_{k+1})}(t) 
$$
for some 
increasing sequence $t_k$ with $t_0 = 0$ and $t_k \rightarrow \infty$.
Adaptedness implies that the $c_{k}$ are measurable w.r.t. $\Sigma_{t_k}$.
If $Y$ is an adapted step function, and $T = t_n$, then
we define the stochastic integral
$$
\int_{0}^{T} Y_t dB_t = 
\sum_{k=0}^{n-1} c_{k}(\omega)(B_{t_{k+1}}(\omega) -
B_{t_{k}}(\omega))\,. 
$$
If $T$ is not one of the $t_k$, we simply make a subdivision to
fix this.

\begin{theorem}[It\^{o}]
The adapted step functions are dense in 
$L^{2}(\Omega \times \mathbb{R}^{+}, \Sigma_t)$,
and the map $Y \mapsto \int_{0}^{\infty}Y_s dB_s$ on step functions 
extends to an isometry 
$I :
L^{2}(\Omega \times \mathbb{R}^{+}, \Sigma_t) 
\rightarrow L^2(\Omega ,\mathbb{P}_{W})$.
We denote $I(Y \one_{[0,t]})$ by $\int_{0}^{t}Y_{s}dB_s$. 
The process $Z_t = \int_{0}^{t}Y_s dB_s$ has a version with 
$t \mapsto Z_t(\omega)$ continuous a.s.
\end{theorem}

\noindent Now that we know how to integrate adapted processes, 
we can define the \emph{stochastic differential equation}
\begin{equation}\label{Stooch}
dX_t = A(X_t , t) dt + S(X_t , t) dB_t\,
\end{equation}
to be just shorthand for the integral equation
$$
X_T = X_0 + \int_{0}^{T} A(X_t , t) dt + \int_{0}^{T} S(X_t , t) dB_t\,.
$$
If $A$ and $S$ both satisfy $|f(x,t)| \leq C (1 + |x|)$
and $|f(x,t) - f(x',t)| \leq D |x - x'|$ uniformly for $t \in [0,T]$,
then a Picard iteration scheme shows that 
equation (\ref{Stooch}) has a unique continuous adapted solution $X_t$
for any given $X_0$.

\begin{theorem}
If $X_t$ and $X'_t$ satisfy $dX_t = A(X_t,t) dt + S(X_t,t) dB_t$
and $dX'_t = A'(X'_t,t) dt + S'(X'_t,t) dB_t$, then 
$X_t X'_t$ satisfies
$d(X_t X'_t) = X_t dX'_t + dX_t X_t + dX_t dX'_t$,
where $dX_t dX'_t$ should be evaluated according to the following It\^{o} table.
\begin{center}
\begin{tabular}{c | c  c}
 $\times$ & $dt$ & $dB_t$\\
 \hline
 $dt$ & $0$ & $0$\\
 $dB_t$ & $0$ & $dt$
\end{tabular}
\end{center}
In other words,
$$
d(X_t X'_t) = \left( X_t A'_t + A_t X'_t + S_t S'_t \right) dt
+ \left( X_tS'_t + SX'_t \right) dB_t\,.
$$
\end{theorem}
This is what ultimately makes stochastic differential equations tractable.
We can just do calculus, as long as we keep track of the second order terms.

For example, let us solve the stochastic differential equation
\begin{equation}\label{vreugde}
dX_t = X_t dB_t
\end{equation}
with initial condition $X_0 = 1$. One would expect exponential behaviour.
For any $C^2$ function $f$, we have
$df(X_t) = f'(X_t) dX_t + \half f''(X_t)(dX_t)^{2}$.
This means that $\exp(B_t)$ will not quite be a solution; 
$d\exp(B_t) = \exp(B_t) dB_t + \half \exp(B_t) dt$.
We correct this, and try $X_t = \exp(B_t - \half t)$ instead.
This does the trick; 
$d \exp(B_t - \half t) = (dB_t - \half dt)\exp(B_t - \half t) +
\half \exp(B_t - \half t) dt = \exp(B_t - \half t) dB_t$,
so that $X_t = \exp(B_t - \half t)$ is the unique solution of 
equation (\ref{vreugde}) with $X_0 = 1$.

\subsection{Fock Space and the Wiener Measure}
We can now realize the fields $\Phi(f)$ as multiplication
operators 
on the Hilbert space
$L^{2}(C_{S}(\mathbb{R}^{+}) , \mathbb{P}_{W})$.
The construction of the required isomorphism
$
\gamma : 
\CF(L^{2}(\mathbb{R}^{+})) \stackrel{\sim}{\rightarrow}
L^{2}(C_{S}(\mathbb{R}^{+}) , \mathbb{P}_{W})
$
is analogous to 
that of the homonymous map
of section
\ref{hogauss}.

The characteristic function 
$\chi_{W} : 
C_{S}^{\infty}(\mathbb{R}^{+} ) \rightarrow \mathbb{C}$
of the Wiener measure
extends to a holomorphic map
$\chi_{W} : 
C_{S}^{\infty}(\mathbb{R}^{+} , \mathbb{C}) \rightarrow \mathbb{C}$,
given by the same formula $\chi_{W}(\phi) = 
\exp(-\frac{1}{2} \int_{0}^{\infty} \phi^{2}(s) ds)$.
The functions
$e_{\phi} : F \mapsto \exp(iF(\phi))$ therefore
have finite expectation $\chi_{W}(\phi) = \mathbb{E}(e_{\phi})$
for all
$\phi \in C_{S}^{\infty}(\mathbb{R}^{+} , \mathbb{C})$.
They are even square integrable.
%
Indeed,
\begin{eqnarray*}
\inp{e_{\phi}}{e_{\psi}} 
&=& 
\int_{\mathcal{D}_{S}(\mathbb{R}^{+})} \exp(-iF(\bar{\phi}) ) \exp(iF(\psi)) 
\mathbb{P}_{W}(dF)\\ 
&=& 
\int_{\mathcal{D}_{S}(\mathbb{R}^{+})} \exp(iF(\psi - \bar{\phi}))\,,
\end{eqnarray*}
so that $\inp{e_{\phi}}{e_{\psi}} = \chi_{W}(\psi - \bar{\phi})$,
and in particular $\inp{e_{\phi}}{e_{\phi}} < \infty$.
Moreover, in terms of the renormalized functions 
\begin{equation}\label{smoelie}
E_{\phi}(F) = \exp\left(\half \int_{0}^{\infty} \phi^{2}(t) dt\right) 
\exp(iF(\phi))\,,
\end{equation}
we have 
the coveted relation
$\inp{E_{\phi}}{E_{\psi}} = \exp(\inp{\phi}{\psi})$.
Since the $E_{\phi}$ are dense in 
$L^{2}(\mathcal{D}_{S}(\mathbb{R}^{+}) , \mathbb{P}_{W})$,
the map $\vexp(\phi) \mapsto E_{\phi}$ extends to 
an isomorphism. 
It is this map, considered 
as an isomorphism
$\gamma :
\CF(L^{2}(\mathbb{R}^{+})) \stackrel{\sim}{\rightarrow} 
L^{2}(C_{S}(\mathbb{R}^{+}) , \mathbb{P}_{W})
$, 
that realizes the fields 
$\Phi(f)$ as multiplication operators.
\begin{theorem}\label{diagveld}
The isomorphism 
$\gamma :
\CF(L^{2}(\mathbb{R}^{+})) \stackrel{\sim}{\rightarrow} 
L^{2}(C_{S}(\mathbb{R}^{+}) , \mathbb{P}_{W})
$ sends $\vexp(\phi)$ to 
\begin{equation}\label{hoelie}
\gamma \vexp(\phi) : B \mapsto 
\exp\left(\half {\int_{0}^{\infty} \phi(t)^2 dt}\right)
\exp\left( {i\int_{0}^{\infty} \phi(t) dB_t} \right)\\
\end{equation}
for all $\phi \in L^{2}(\mathbb{R}^{+}, \mathbb{C})$, and
$\Phi(f)$ to
\begin{equation}\label{troelie}
\gamma \Phi(f)\gamma^{-1}  : L(B) \mapsto  
\left( \int_{0}^{\infty} f dB_t \right) L(B)
\end{equation}
for all $f \in L^{2}(\mathbb{R}^{+},\mathbb{R})$. Moreover, the 
probability measure induced on $C_{S}(\mathbb{R}^{+})$
by the vacuum state 
$\vexp(0)$
is precisely the Wiener measure $\mathbb{P}_{W}$.
\end{theorem}
\proof
If we consider the Wiener measure to be based on 
$C_{S}(\mathbb{R}^{+})$, then
equation (\ref{smoelie}) 
and the definition of $\gamma$ yield
(\ref{hoelie}). 
Using equation (\ref{WeylOp}), one calculates that
$$
\gamma W(v) \gamma^{-1} E_{\phi}(B) = 
\exp\left(i \int_{0}^{\infty} v (dB + \basim (v) ds)\right)
E_{\phi}(B_t + 2 \int_{0}^{t} \basim (v) ds)\,.  
$$
Because the $E_{\phi}$ are dense,
this continues to hold if one replaces $E_{\phi}$
by any $L \in L^{2}(C_{S}(\mathbb{R}^{+}) , \mathbb{P}_{W})$. 
Differentiation then yields\footnote{The conjugate momenta are given by
$\smash[b]{
\Pi(f) L(B) = 
\left(i \int_{0}^{\infty} f dB_t - 2 i \frac{\del}{\del \int f} \right) L(B)\
}$, but we will not need them here.}
equation (\ref{troelie}). 

The vacuum state $\vexp(0)$ induces a probability measure 
$\mathbb{P}_{(f_1, \ldots, f_n)}$
on the joint spectrum $\mathbb{R}^{n}$ of 
$
(\Phi(f_1), \ldots, \Phi(f_n))
$. 
We show that this coincides with the joint 
probability distribution $\mathbb{P}^{W}_{(f_1, \ldots, f_n)}$ of 
$\left(\int_{0}^{\infty} f_1(t) dB_t, \ldots, 
\int_{0}^{\infty} f_n(t) dB_t\right)
$ under the Wiener measure. 
This suffices to prove the theorem,
because the cylinder sets 
$\{F \,;\, (F(f_1), \ldots F(f_n)) \in S \}$
with open $S\subset \mathbb{R}^{n}$ generate the 
Borel $\sigma$-algebra of $\mathcal{D}_{S}(\mathbb{R}^{+})$.  

The 
characteristic function
of $\mathbb{P}_{(f_1, \ldots, f_n)}$ is
$$
\chi(f_1 , \ldots ,f_n \,;\, k_1 ,\ldots k_n) = 
\inp{\vexp(0)}{
\prod_{r=1}^{n}\exp(i \Phi(k_r f_r))
\vexp(0)}\,.
$$
Because 
$\exp(i\Phi(kf)) = W(kf)$, 
and $\prod_{i} W(k_i f_i) = W(\sum_i k_i f_i)$, 
we see that 
$
\chi(f_1 , \ldots ,f_n \,;\, k_1 ,\ldots k_n) =
\chi(\sum_{i} k_i f_i)
$
with
$\chi(f) = 
\inp{\vexp(0)}{
W(f)
\vexp(0)}
$. It equals $\chi(f) = \exp(-\half \inp{f}{f})$
by equation (\ref{WeylOp}).

Similarly, $\mathbb{P}^{W}_{(f_1, \ldots, f_n)}$
has characteristic function
$$
\chi_{W}(f_1 , \ldots ,f_n \,;\, k_1 ,\ldots k_n) =
\int  
\prod_{r=1}^{n} \exp(ik_r \int_0^{\infty} f_r(t) dB_t) \mathbb{P}_{W}(dB)\,.
$$
This is equal to $\chi_{W}(\sum_{i}k_i f_i)$,
with 
$\chi_{W}(f) = 
\int  
\exp(i\int_0^{\infty} f(t) dB_t) \mathbb{P}_{W}(dB)$.
We have
$\chi_{W}(f) = 
\exp(-\half \inp{f}{f})
$ by definition \ref{wieien}.
The probability measures $\mathbb{P}_{(f_1, \ldots, f_n)}$
and $\mathbb{P}^{W}_{(f_1, \ldots, f_n)}$ now agree because 
their characteristic functions do.
\qed

\subsection{Quantum Stochastic Differential Equations}\label{qsdesec}

We have defined stochastic calculus in terms of the random variables 
$B_{t}$. 
According to theorem (\ref{diagveld}), they correspond 
to the field operators
$\Phi(\one_{[0,t]}) = i(a(\one_{[0,t]}) - a^{\dagger}(\one_{[0,t]}))$
through the isomorphism $\gamma$.
%

However, we may as well choose to diagonalize\footnote{This 
is easily done by composing
$\gamma : \CF(L^{2}(\mathbb{R}^{+})) 
\stackrel{\sim}{\rightarrow} L^{2}(C_{S}(\mathbb{R}^{+}) , \mathbb{P}_{W})$ 
with the `Fourier transform'
$\CF(i) \in U(\CF(L^{2}(\mathbb{R}^{+})))$.
} 
the conjugate momenta
$\Pi(\one_{[0,t]}) = a(\one_{[0,t]}) + a^{\dagger}(\one_{[0,t]})$,
and integrate w.r.t. those instead.
We then have two \emph{noncommuting} notions of stochastic integration
living on the same Fock space. 
 
This train of thought is brought to its logical conclusion by the
theory of quantum stochastic integration,  
due to Hudson and Parthasarathy.

\subsubsection{Adapted Processes}

We take $\CH = L^{2}(\mathbb{R}^{+})$, and write 
$\CH = \CH_{t} \oplus \CH^{t}$ with
$\CH_{t} = L^{2}([0,t))$ and $\CH^{t} = L^{2}([t,\infty))$. 
If we write $\CF = \CF(\CH)$, $\CF_t = \CF(\CH_t)$
and $\CF^{t} = \CF(\CH^{t})$, then this translates into
$\CF = \CF_{t} \otimes \CF^{t}$.

We take $S$ to be the dense set of locally bounded functions in $\CH$.
The set $\vexp(S)$ of exponential vectors then spans $\CF$.

We choose an auxiliary Hilbert space $\CK$. 
One may think of $\CF$ as the Hilbert space of the
electromagnetic field, and of $\CK$ as the Hilbert space of some system 
(an atom, say) in 
interaction with this field.
Our dynamics will take place
in $\CK \otimes \CF$, the Hilbert space of the combined system.

The von Neumann algebra $\CN = B(\CK \otimes \CF)$ of observables 
possesses the natural filtration $\{\CN_{t} \,;\, t \in \mathbb{R}^{+}\}$,  
with $\CN_{t} = B(\CH \otimes \CF_{t}) \otimes \one_{\CF^{t}}$
the von Neumann algebra of observables up to time $t$.
We denote its commutant by $\CN^{t}$. 
It is the algebra of future events, 
$\CN^{t} = \one_{\CK \otimes \CF_t} \otimes B(\CF^t)$.

\begin{basdef}
An \emph{adapted process} based on $S$
is a family $\{F_t \,;\, t \in \mathbb{R}^{+}\}$ of operators
such that $F_t$ is affiliated to $\CN_{t}$, and $\CK \otimes \vexp(S)$
is in the domain of each member. 
\end{basdef}
The operator $F_t$ does not `know' what happens after time $t$.
The common core is a convenient and not overly restrictive technical
requirement. 

\begin{basdef}
An adapted process $F_t$ is called \emph{simple} if there is an
increasing sequence
$t_k$ with $t_0 = 0$ and $t_k \rightarrow \infty$ such that
$F_t = \sum_{k=0}^{\infty} F_k \one_{[t_k , t_{k+1})}(t)$.
It is called \emph{continuous} if the map $\gamma : t \mapsto F_t u 
\otimes \vexp(f)$
is continuous for all $u \in \CK$, $f \in S$, and 
it is called \emph{locally square
integrable} if $\gamma$ is measurable, and $\int_{0}^{t} \|\gamma(s)\|^{2} ds$
is finite for each $t \in \mathbb{R}^{+}$. 
\end{basdef}
The locally square integrable processes are sufficiently general 
to cater to our needs, and the simple ones are, as the name suggests,
easy to handle. 
Luckily, locally square integrable processes can be approximated by simple ones.
\begin{proposition}[Muthuramalingam, \cite{HuP84}]\label{multitroela} 
Any locally square in\-te\-gra\-ble process $F_t$ can be approximated by a 
sequence of simple processes $F_t^{n}$, in the sense that
for each $t \in \mathbb{R}^{+}$, $u \in \CK$ and $f \in S$, we have
$$
\lim_{n \rightarrow \infty} \int_{0}^{t}
\|(F_s - F^{n}_s) u \otimes \vexp(f)\|^{2} ds = 0\,.
$$
\end{proposition}
\proof
One first shows that the simple processes approximate continuous ones.
From this, one sees that
it suffices to approximate $F_t$ by the convolution processes 
$F_t^{k} = (F * k \one_{[0,1/k]})(t)$, which are continuous. 
See \cite{HuP84} for details.\qed  

\subsubsection{Quantum Stochastic Integration}

Stochastic integration is integration against the distinguished
Wiener process $B_t$. In the noncommutative version, there are 
two separate processes which take its role:
the annihilation process $A(t) = a(\one_{[0,t]})$
and the creation process $A^{\dagger}(t) = a^{\dagger}(\one_{[0,t]})$.

\begin{basdef}
Let $F_t = \sum_{k=0}^{\infty} F_k \one_{[t_k , t_{k+1})}(t)$
be a simple process, and let $T = t_n$. (Refine the sequence $t_k$
by inserting $T$ if necessary.)
%
We then
define the \emph{stochastic integral}
$$
\int_{0}^{t} F_t dA_t = \sum_{k=0}^{n-1} F_k (A(t_{k+1}) - A(t_k))\,.
$$
The expression $\smash{\int_{0}^{t}} F_t dA^{\dagger}_t$ is 
defined in a similar fashion,
with $A^{\dagger}$ in stead of $A$, and of course we have 
$\int_{0}^{t}F_t dt = \sum_{k=0}^{n-1} F_k (t_{k+1} - t_k)$ 
as usual.
%
\end{basdef}
Note that $F_k$ is affiliated to $\CN_{t_k}$, 
acting essentially on
$\CK \otimes \CF_{t_k}$.
On the other hand, 
$A(t_{k+1}) - A(t_k) = \smash{a(\one_{[t_k , t_{k+1}]})}$ 
is affiliated to its commutant $\CN'_{t_{k}}$, 
acting essentially on $\smash{\CF^{t_k}}$.
In particular, $F_t$ commutes with $\Delta A_{t_k}$ and 
$\Delta \smash{A^{\dagger}_{t_k}}$. 
We have used (and will use from now on) 
the notation $\Delta f_t$ for $f(t_{k+1}) - f(t_k)$.

We can now define the stochastic integral of any locally square integrable
process $F_t$ by approximating it with simple processes.
If $F_t^n \rightarrow F_t$ in the sense of proposition \ref{multitroela}, 
then one can show, (and this is not entirely trivial \cite{HuP84}),
that
$$
 \lim_{n \rightarrow \infty} F_t^n dA_t\,
$$
converges on vectors of the type $u \otimes \vexp(f)$.
We will take this to be our stochastic integral
$\int_{0}^{t} F_t dA_t$ (and similarly for $\int_{0}^{t} F_t dA^{\dagger}_t$).
We introduce the `stochastic differential equation'
$dM_t = F_t dA_t + G_t dA^{\dagger}_t + H_t dt$
as shorthand notation for the integral equation 
$M_t = M_0 + 
\int_{0}^{t} F_s dA_s + 
\int_{0}^{t} G_s dA^{\dagger}_s + 
\int_{0}^{t}H_s ds$.

\subsubsection{A Quantum Version of It\^{o}'s Formula}

We can now prove a noncommutative version of It\^{o}'s formula for stochastic 
integrals of continuous processes.
\begin{theorem}[Hudson, Parthasarathy]\label{qiet}
Let $M_t$ and $M'_t$ be integrals of continuous processes, i.e. 
$dM_t = F_t dA_t + G_t dA^{\dagger}_t + H_t dt$
and
$dM'_t = F'_t dA_t + G'_t dA^{\dagger}_t + H'_t dt$.
Suppose that the
products between the operators $M_t , F_t , G_t, H_t$
and $M'_s ,F'_s , G'_s, H'_s$
are well defined on $\CK \otimes \vexp(S)$.
Then the adapted process $t \mapsto M_t M'_t$ satisfies 
$
d (M_t M'_t) = d(M_t) M'_t + M_t d(M'_t) + dM_t dM'_t \,,
$
where $dM_t dM'_t$ should be evaluated according to the following
quantum It\^{o} table.
\begin{center}
\begin{tabular}{c | c  c  c}
 $\times$ & $dt$ & $dA_t$ & $dA^{\dagger}_t$ \\
 \hline
 $dt $ & $0$ & $0$ & $0$ \\
 $dA_t$ & $0$ & $0$ & $dt$\\
 $dA^{\dagger}_t$ & $0$ & $0$ & $0$\\  
\end{tabular}
\end{center}
In other words,
$$
d(M_tM'_t) = 
(F_t M'_t +M_t F'_t) dA_t + (G_t M'_t + M_t G'_t ) dA^{\dagger}_t + 
(H_t M'_t + M_t H'_t  + F_t G'_t)dt\,.
$$
\end{theorem}
\proof
We expand $\Delta (M_t M'_t) = (\Delta M_t) M'_t + M_t (\Delta M'_t) +
\Delta M_t \Delta M'_t$, and focus on the last term.
In the limit $\Delta t \rightarrow 0$, this equals
\begin{eqnarray*}
\Delta M_t \Delta M'_t &=& 
F_t F'_t (\Delta A_t)^2  
+ G_t G'_t (\Delta A^{\dagger}_t)^{2} 
+ H_t H'_t (\Delta t)^2\\
& &
+ F_t G'_t \Delta A_t \Delta A^{\dagger}_t
+ (F_t H'_t + H_t F'_t) \Delta A_t \Delta t\\
& &
+ G_t F'_t \Delta A^{\dagger}_t \Delta A_t 
+ (G_t H'_t + H_t G'_t) \Delta A^{\dagger}_t \Delta t\,.
\end{eqnarray*}
Each one of these 7 terms is a product of 
a term in $\CN_t$ and one in $\CN^t$, 
so that we can split their contributions.
For example, take
$F_t G'_t \Delta A_t \Delta A^{\dagger}_t$. Consider it as
$F_t G'_t \otimes \one_{\CF^t} \times 
\one_{\CK \otimes \CF_t} \otimes \Delta A_t \Delta A^{\dagger}_t$, 
and 
split $u \otimes \vexp(h) = 
(u \otimes \vexp(h \one_{[0,t]})) \otimes \vexp(h \one_{[t , \infty]})$
accordingly. We see that
$$
\inp{u \otimes \vexp(h)}
{ F_t G'_t \Delta A_t \Delta A^{\dagger}_t u' \otimes \vexp(h')}
$$
equals 
$$
\inp{u \otimes \vexp(h \one_{[0,t]})}
{ \,F_t G'_t \,\, u' \otimes \vexp(h'\one_{[0,t]})}
$$
times 
$$
\inp{\vexp(h \one_{[t,\infty)})}
{ a(\one_{[t_k , t_{k+1}]})  a^{\dagger}(\one_{[t_k , t_{k+1}]})
 \vexp(h'\one_{[t , \infty)})}\,.
$$
A look at equation (\ref{aaaad}) reveals that 
this second factor equals 
$$\left \{ \left(\int_{t_k}^{t_{k+1}} h'(s) ds \right)
\left( \int_{t_k}^{t_{k+1}}\bar{h}(s) ds \right) +
\Delta t \right\}
\inp{\vexp(h \one_{[t,\infty)})}{\vexp(h' \one_{[t,\infty)})}\,.$$
As $h'$ and $h$ are locally bounded, the first term
is quadratic in $\Delta t$. If we take $\Delta t \rightarrow 0$, 
we are left with
$
\inp{u \otimes \vexp(h)}
{ \,F_t G'_t \,\, u' \otimes \vexp(h')}\Delta t
$.
The other 6 terms are of second order, as 
equations (\ref{aadaad}), (\ref{aaaa}), (\ref{aadaa})
lack the extra $\inp{u_1}{u_2}$ term of equation (\ref{aaaad}).
\qed

\subsubsection{Quantum Stochastic Differential Equations}

Let $V_1$, $V_2$ and $V_3$ be bounded operators on $\CK$.
We can then define the quantum stochastic differential equation
\begin{equation}\label{stog}
dM_t = M_t (V_1 dA_t + V_2 dA^{\dagger}_t + V_3 dt)
\end{equation}
as shorthand for the quantum integral equation
$$
M_t = M_0 + \int_{0}^{t} M_t (V_1 dA_t + V_2 dA^{\dagger}_t + V_3 dt)\,.
$$

\begin{theorem}
Equation (\ref{stog}) has a unique continuous solution with $M_0 = \one$. 
It is unitary
if and only if $(V_1 , V_2 , V_3)$ takes the shape
$(V^{\dagger} , -V , iH - \half V^{\dagger}V)$.
\end{theorem}
\proof
The existence of a unique solution can be shown by means of a 
Picard-type iteration procedure,
see \cite{HuP84}. The solution is unitary, $M^{\dagger}_t M_t = \one$,
if and only if $d (M^{\dagger}_t M_t) = 0$.
Using theorem \ref{qiet}, we calculate this by adding  
$M^{\dagger}_t d M_t = V_1 dA_t + V_2 dA^{\dagger}_t + V_3 dt$, 
$d(M^{\dagger}_t) M_t = V_{1}^{\dagger} dA^{\dagger}_t + 
V^{\dagger}_2 dA_t + V^{\dagger}_3 dt$ and
$d(M^{\dagger}_t) d (M_t) = V^{\dagger}_{2} V_2 dt$. 
The result is identically zero if and only if 
$V_1 + V_{2}^{\dagger} = 0$ and $V_3 + V^{\dagger}_3 + V^{\dagger}_2 V_2 =0$.
\qed
If we think of $\CF(L^{2}(\mathbb{R}^{+}))$ as the Hilbert space
of the electromagnetic field, and of $\CK$ as the Hilbert space
of an atom coupled to this field, then time evolution
is described by the Schr\"{o}dinger equation
$$
i \frac{d}{dt}U^{\lambda}(t) = 
(H_S\otimes \one + \one \otimes H_F + \lambda H_I)U^{\lambda}(t)\,
$$
on $\CK \otimes \CF(L^{2}(\mathbb{R}^{+}))$. In this expression, 
$H_S$ is the Hamiltonian
of the system $\CK$, $H_F$ is the the free field Hamiltonian on
$\CF(L^{2}(\mathbb{R}^{+}))$, and 
$H_I = i(D \otimes a^{\dagger}(g) - D^{\dagger} \otimes a(g))$
is the interaction 
Hamiltonian on $\CK \otimes \CF(L^{2}(\mathbb{R}^{+}))$ with
coupling constant $\lambda$.

Since $U^{0}(t)$ can often be easily solved, it suffices to find
the interaction picture time evolution 
$\hat{U}^{\lambda}(t) := U^{\lambda}(t)U^{0}(-t)$. 
The `van Hove' limit of
$U^{\lambda}(t/\lambda^{2})$ for $\lambda \downarrow 0$ 
exists under suitable conditions.
According to a theorem of Accardi, Frigerio and Lu, see \cite{AFL},
it is given
by a unitary
quantum stochastic differential equation of the form
$$
dU_t = U_t(V^{\dagger}dA_t - V dA_{t}^{\dagger} - 
\half V^{\dagger} V dt) \,.
$$
This equation will play a pivotal role in chapters 
\ref{ch:OPJM} and \ref{ch:OEQS}.

\section{Bundles and Classical Fields}\label{bunklasveld}
We have reviewed the foundations of quantum mechanics, 
and scratched the surface of a tiny part of
quantum field theory. We will now take a step back, 
and submit the basic framework of classical field theory
to a closer look.

In classical field theory, a system is modelled by a 
smooth fibre bundle
$\pi : F \rightarrow M$, and fields are described
by sections of $\pi$.
The base manifold $M$ represents space-time.
The dynamics are prescribed by a Lagrangian density
$\CL : J^{1}F \rightarrow \wedge^{n}T^{*}M$, with
$n$ the dimension of $M$.
Transformations of the theory correspond to automorphisms 
of $F$.

\subsection{Bundles, Sections and Jets}
Intuitively speaking, 
a fibre bundle over $M$ with fibre $F_0$
is a smoothly varying family $\{F_x \,;\, x \in M\}$ of copies of
$F_0$. A section is a function $\phi : M \rightarrow F$ 
such that its value 
on $x \in M$ lies in the fibre $F_{x}$ over $x$.
The $k$-jet $j^{k}_{x}\phi$ of $\phi$ at $x$ is its value $\phi(x)$ at $x$,
together with its derivatives at $x$ up to order $k$.    
Let us make this a bit more precise. 

\subsubsection{Fibre Bundles}
If $F$, $F'$, $M$ and $M'$ are smooth manifolds, then an \emph{isomorphism} 
between the maps 
$\pi : F \rightarrow M$
and $\pi' : F' \rightarrow M'$ is by definition
a pair $(\phi , \phi_{M})$ of diffeomorphisms  
$\phi : F \rightarrow F'$ and $\phi_{M} : M \rightarrow M'$
such that $\pi' \circ \phi = \phi_{M} \circ \pi$.
Since $\phi_{M}$ is completely determined by $\phi$, we will
suppress it in the notation, and  
write $\phi$ for $(\phi , \phi_{M})$.
If $M = M'$, then $\phi$ is called
\emph{vertical} if $\phi_{M} = \id_{M}$. 
 
\begin{basdef}\label{vezelrijk}
Let $F$, $M$ and $F_0$ be smooth manifolds. 
Then a \emph{fibre bundle} $\pi : F \rightarrow M$ with fibre $F_0$
is a smooth map with the property that each $x \in M$ has an open neighbourhood
$U$ such that the restriction $ \pi : \pi^{-1}(U) \rightarrow U$ 
is vertically isomorphic to the projection $\pi' : U \times F_0 \rightarrow U$
on the first factor.
\end{basdef}
In short, $F$ is locally the product of the base by $F_{0}$.
A vertical isomorphism $\phi : \pi^{-1}(U) \rightarrow U \times F_0$ 
is called a \emph{local trivialization}. 
The manifold $F_{x} := \pi^{-1}(\{x\})$ is called the \emph{fibre} 
over $x$. 
Each local trivialization restricts to a diffeomorphism 
$\phi|_{F_{x}} : F_{x} \stackrel{\sim}{\rightarrow} \{x\} \times F_{0}$,
so that
in particular, 
any fibre $F_{x}$ is diffeomorphic to $F_0$.

For any fibre bundle $F \rightarrow M$, we can produce 
\emph{adapted co-ordinates} as follows.
Choose
co-ordinates $x^\mu : M \supset U \rightarrow \mathbb{R}$ on $M$, 
co-ordinates $v^{a} : F_0 \supset V \rightarrow \mathbb{R}$ on $F_0$,
and local trivializations 
$\phi_{U} : \pi^{-1}(U) \rightarrow U \times F_0$.
The adapted co-ordinates on $F$ are then
$(x^{\mu} \circ \phi_{U} , v^{a} \circ \phi_{U})$.
We will simply denote them
$(x^\mu , v^{a})$. They have the convenient property that 
$x^{\mu}(f) = x^{\mu}(f')$ for all $\mu$ if and only if 
$f$ and $f'$ are in the same
fibre.

Finally, if $F \rightarrow M$ is any fibre bundle, and
$f : N \rightarrow M$ is a smooth map, then we can define 
the pullback
bundle $f^*F \rightarrow N$ as the manifold
$f^*F = \{(n,f) \in N \times F \,;\, f(n) = \pi(f)\}$,
with projection $(n,f) \mapsto n$.
The idea is that the map $f$ allows one to
regard the fibre over $f(n) \in M$ as a fibre over $n \in N$,
i.e.\
$f^*F_{n} = F_{f(n)}$.

\subsubsection{Vector Bundles and Principal Fibre Bundles}
It is useful to study fibre bundles with additional structure
on the fibres.  
A fibre bundle $\pi : E \rightarrow M$ is called a 
\emph{vector bundle} if $E_0$ and
each of the fibres $E_{x}$ are vector spaces, 
and if one can choose the local trivializations $\phi$ over $U$ 
in such a way that for any $x\in M$, the map
$\phi|_{E_x} : E_{x} \stackrel{\sim}{\rightarrow} E_0$ is an 
isomorphism of vector spaces.
%
Similarly, 
a fibre bundle $\pi : P \rightarrow M$ is called a 
\emph{principal fibre bundle} with Lie group $G$ if 
$P_0$ and the fibres $P_{x}$ 
carry a free transitive
right action of $G$, and if the local trivializations 
can be chosen to be $G$-equivariant.

The tangent bundle is a natural
example of a vector bundle.
The \emph{tangent space} $T_{x}M$ of $M$ at $x \in M$
is by definition the
space of vectors tangent to $M$ at $x$.
If $M$ is $n$-dimensional, then $T_{x}M$ is 
a vector space
isomorphic to $\mathbb{R}^{n}$. 
This makes the
\emph{tangent bundle}
$TM = \cup_{x \in M}T_{x}M$ 
into a vector bundle, with projection 
$\pi : TM \rightarrow M$ assigning to each
tangent vector $\xi_{x} \in T_{x}M$ 
its base point $\pi(\xi_{x}) = x$.

Starting from the tangent bundle, one can construct many other
interesting fibre bundles.
If $\pi : F \rightarrow M$ is a any fibre bundle, then the 
kernel $T_{f}^{V}F$ of the pushforward 
$\pi_{*} : T_{f}F \rightarrow T_{\pi(f)}M$ is a linear 
subspace of $T_{f}F$, and we obtain the vertical
tangent bundle $T^{V}F \rightarrow F$, a vector subbundle
of $TF \rightarrow F$.
One can think of it as the union of all the separate 
tangent spaces $T(F_{x})$. 

The \emph{cotangent bundle} $\pi : T^{*}M \rightarrow M$ 
is the vector bundle defined by letting $T^{*}_{x}M$
be the dual of $T_{x}M$, i.e.\ 
the vector space of linear maps $T_{x}M \rightarrow \mathbb{R}$.

Similarly, a \emph{volume form} on the vector space $T_{x}M$ is
an alternating multilinear map 
$T_{x}M \times \ldots \times T_{x}M \rightarrow \mathbb{R}$
on $n$ copies of $T_{x}M$.
We denote the 1-dimensional space of volume forms 
on $T_{x}M$ by
$\wedge^{n}T_{x}^{*}M$. Their union 
$\wedge^{n}T^{*}M = \cup_{x \in M} \wedge^{n}T_{x}^{*}M$
is a line bundle; a vector bundle over $M$ with 1-dimensional fibre. 

 
The \emph{frame bundle} $F(M) \rightarrow M$
is a natural example of a principal fibre bundle. 
A frame at $x$ is by definition a linear
isomorphism $f : \mathbb{R}^{n} \rightarrow T_{x}(M)$. 
If we fix the standard basis $(e_1, \ldots, e_{n})$ of
$\mathbb{R}^{n}$, then frames $f$ at $x$ correspond precisely
to bases $(f(e_1), \ldots, f(e_n))$ of $T_{x}(M)$. 
The space of frames at $x$, denoted $F_{x}(M)$,
has a free transitive right action $ g : f \mapsto f \circ g$ of 
$\mathrm{GL}(\mathbb{R}^{n})$, making 
$F(M) = \cup_{x\in M}F_{x}(M)$ into a principal fibre bundle.

\subsubsection{Sections}
A \emph{section} of $\pi : F \rightarrow M$ is a smooth map 
$\phi : M \rightarrow F$ such that $\pi \circ \phi = \id_{M}$,
and we
denote the space of sections by $\Gamma(F)$. 

If $F$ happens to be the trivial bundle $M \times F_0$, then 
sections correspond precisely to functions $M \rightarrow F_0$.
One way to view sections and bundles is as follows. 
Take the base manifold $M$, and cut it into parts $U_{i}$
that overlap only at their boundaries, and that are small enough to
fit into a trivializing neighbourhood.
A section $\phi$ of $F$ is then a collection of functions
$\phi_i : U_{i} \rightarrow F_{0}$, and the bundle $F \rightarrow M$ encodes
the boundary conditions that $\phi_i$ and 
$\phi_j$ have to satisfy on $\del U_i \cap \del U_j$. 

A section of the tangent bundle $TM \rightarrow M$ is 
precisely a smooth vector field on $M$, so that the 
commutator bracket $[\xi,\chi]$ of vector fields endows $\Gamma(TM)$
with a structure of Lie algebra.
A section of $T^*M$ is called a 1-form, and a section of $\wedge^{n}T^*M$
a volume form.

The tangent bundle $TS^{2} \rightarrow S^{2}$ of the
$2$-sphere is a nice nice example of a vector bundle
that is not trivial.
Every vector field on 
$S^{2}$ must have at least one zero,\footnote{%
The theorem of P.\ Bohl, J.\ Hadamard and L.\ E.\ J.\ Brouwer that every 
(continuous) vector field
on a sphere of even dimension must have at least one zero is known
colloquially 
as the `hairy ball theorem': one cannot 
comb the hair on a $2n$-sphere without creating a crown. 
} showing
that a section of $TS^{2}\rightarrow S^{2}$ is not quite the same as
a function $S^{2}\rightarrow \mathbb{R}^{2}$. 


\subsubsection{Jet Bundles}
Intuitively speaking, the $k$-jet of a section $\phi$ at $x$ is
its value $\phi(x)$, together with its derivatives at $x$ up to order $k$. 

To put this in a more formal setting, we introduce 
the relation $\sim$ on the local sections around $x$, where
$\phi \sim \phi'$ if $\phi$ and $\phi'$ have the same value at $x$, 
and if their derivatives agree up to order $k$.
The $k$-jet of $\phi$ at $x$, denoted $j^{k}_{x}\phi$,
is now simply the equivalence class of $\phi$
modulo $\sim$. 

We denote the set of all $k$-jets at $x$ by $J_{x}^{k}F$, and 
we define $J^{k}F := \cup_{x \in M} J_{x}^{k}F$
to be the \emph{jet bundle}.
It is a fibre bundle in two different ways.
The map $\pi_{k} : J^{k}F \rightarrow M$ defined by 
$j_{x}^{k}(\phi) \mapsto x$ makes $J^{k}F$ into a 
bundle over $M$, and the map
$\pi_{k,0} : J^{k}F \rightarrow F$ that maps 
$j_{x}^{k}(\phi)$ to $\phi(x)$ makes it into
a bundle over $F$. We have $\pi \circ \pi_{k,0} = \pi_{k}$.
Any section $\phi$ of $\pi : F\rightarrow M$ yields a section $j^{1}\phi$
of $\pi_{k} : J^{k}F \rightarrow M$, but 
not all sections of $J^{k}F$ are of this form.

Adapted co-ordinates $(\tilde{x}^{\mu} , \tilde{v}^{a})$ on $F$ give rise to 
the co-ordinates 
$(x^{\mu} , v^{a} , v^{a}_{\mu})$ on $J^{1}F$ by 
$x^{\mu}(j^{1}_{x}\phi) = \tilde{x}^{\mu}(x) $,
$v^{a}(j^{1}_{x}\phi) = \tilde{v}^{a}(\phi(x))$ and
$v^{a}_{\mu}(j^{1}_{x}\phi) = \del_{\mu} \tilde{v}^{a}(\phi(x))$.
If $\phi$ is a section of $F$, then we denote 
$\phi^{a} := v^{a} \circ \phi$, so that 
$\del_{\mu}\phi^{a} = v^{a}_{\mu}(j^{1}\phi)$.

\subsection{The Euler-Lagrange Equation}
Having described its basic ingredients, we now have a closer look 
at classical field theory itself.
A physical system is described by a fibre bundle $\pi : F \rightarrow M$, 
the base manifold $M$ representing space-time.
Fields are modelled by sections of $\pi$, and the
dynamics are prescribed 
by a bundle map $\CL : J^{1}F \rightarrow \wedge^{n}(T^{*}M)$ 
called\footnote{Generalization to Lagrangian densities on  
$J^{k}F$ with $k>1$ is straightforward, and the Euler-Lagrange 
equation then becomes a PDE of order $k+1$.  
We focus on the case $k=1$ for notational convenience, but also because
most fundamental field equations are of order two.}
the \emph{Lagrangian density}, or Lagrangian for short.
It singles out a class of preferred sections by means of the 
\emph{action principle}. These are considered to be the 
fields of physical relevance.

\subsubsection{The Action Principle}

Let $U$ be an open submanifold of $M$ with smooth boundary 
and
compact closure.  We define the
\emph{action functional} $S_{U} : \Gamma(F) \rightarrow \mathbb{R}$
by $S_{U}(\phi) = \int_{U} j^{1}\phi^{*} \CL$. The action principle then 
states that the physical fields are precisely the critical points of the action
$S_{U}$ under infinitesimal variations with support inside $U$.

An infinitesimal variation of a section $\phi \in \Gamma(F)$ is precisely  
a section
$\delta \phi$ of $\phi^{*}T^{V}F$, 
the pullback of the vertical tangent bundle of 
$F$ along $\phi$. 
Indeed, if $\phi_{\varepsilon}$ is a one-parameter family of sections of $F$ with 
$\phi_{0} = \phi$, then  
$\delta \phi(x) = \frac{d}{d\varepsilon}|_{0} \phi_{\varepsilon}$
yields a vertical vector in $T_{\phi(x)}F$
for each $x \in M$.

If we choose adapted co-ordinates 
$(x^{\mu} , v^{a} , v^{a}_{\mu})$, then we can write 
the Lagrangian density as
$\CL = L(x^{\mu} , v^{a} , v^{a}_{\mu}) \vol$, with 
$\vol = dx^{1} \wedge \ldots \wedge dx^{n}$
the volume form assigning unit volume to the frame
$(\frac{\del}{\del x^{1}}, \ldots, \frac{\del}{\del x^{n}})$, 
and $L$
the so-called Lagrangian function.
The action then reads 
$S_{U}(\phi) = \int_{U} L(x^{\mu} , \phi^{a} , \del_{\mu} \phi^{a}) \vol$.

In order to avoid a jungle of sigmas, we will adopt the convention that
repeated indices indicate summation. 
We also use the shorthand notation \smash{$\frac{\del L}{\del \phi^{a}}$}
for 
$x \mapsto 
\frac{\del L}{\del v^{a}}(x^{\mu}(x) , \phi^{a}(x) , \del_{\mu} \phi^{a}(x))$,
and similarly, we write 
\smash[b]{$\frac{\del L}{\del (\del_{\mu} \phi^{a})}$} for the function
$x \mapsto 
\frac{\del L}{\del v^{a}_{\mu}}(x^{\mu}(x) , \phi^{a}(x) , 
\del_{\mu} \phi^{a}(x))$.

The following theorem shows that
critical sections 
correspond to solutions of a second order P.D.E.\ called the
\emph{Euler-Lagrange equation}.
\begin{theorem}[Euler-Lagrange]
A section $\phi$ of $F$ is a stationary point of
the action functionals
$S_{U}$ under variations $\delta \phi$ with support 
strictly contained in $U$ if and only if it satisfies
the second order P.D.E.\ 
\begin{equation}\label{eulagr}
[\CL\,]^{\phi} = 0\,,
\end{equation} 
called the Euler-Lagrange equation.
The map 
$[\CL\,]^{\phi} : \phi^{*}T^{V}F \rightarrow \wedge^{n}T^{*}M$
of vector bundles
is well
defined by 
\begin{equation}\label{preel}
[\CL\,]^{\phi}(\delta \phi) = 
\left(
	\left(\del_{\mu} \frac{\del L}{\del (\del_{\mu}\phi^{a})}\right) 
	-
	\frac{\del L}{\del \phi^{a}} 
\right)
\vol \cdot \delta\phi^{a}\,.
\end{equation}
The r.h.s. of (\ref{preel}) does not depend on the 
choice of co-ordinates, whether $\phi$ is critical or not.
\end{theorem}
\proof
Since $\frac{\del L}{\del \phi^{a}}(x)$ and 
$\frac{\del L}{\del (\del_{\mu}\phi^{a})}(x)$ only 
depend on $\phi$ through $j^{1}_{x}\phi$, we see from (\ref{preel}) 
that
$[\CL \, ]^{\phi}(x)$ depends only on $\phi$ through 
$j^{2}_{x}\phi$, so that (\ref{eulagr}) 
is indeed a P.D.E.\ of order 2.

The variation of $S_{D}$ along $\delta \phi = \delta \phi^{a}\del_{a}$ reads
$$
\delta S_{U}(\phi) = \int_{U} 
\left(
\frac{\del L}{\del \phi^{a}} \delta \phi^{a}
+
\frac{\del L}{\del (\del_{\mu} \phi^{a})} \delta \del_{\mu} \phi^{a}
\right)
\vol\,.
$$
Note that 
$
\delta \del_{\mu} \phi^{a} 
= 
\del_{\mu} \delta \phi^{a}
$, because $\frac{d}{d\varepsilon}|_{0} \del_{\mu}\phi_{\varepsilon}^{a} 
= 
\del_{\mu} \frac{d}{d\varepsilon}|_{0} \phi_{\varepsilon}^{a} 
$.
We can then write 
$$
\frac{\del L}{\del (\del_{\mu} \phi^{a})} \delta \del_{\mu} \phi^{a}\vol
=
d\left(
\frac{\del L}{\del (\del_{\mu} \phi^{a})} \delta \phi^{a}
i_{\del_{\mu}} \vol \right)
- \left(\del_{\mu} \frac{\del L}{\del (\del_{\mu}\phi^{a})}\right) 
\delta \phi^{a}\vol\,.
$$
If we now define
$[\CL\,]^{\phi} : \phi^{*}T^{V}F \rightarrow \wedge^{n}T^{*}M$
by 
\begin{equation}
[\CL\,]^{\phi}(\delta \phi) = 
\left(
	\left(\del_{\mu} \frac{\del L}{\del (\del_{\mu}\phi^{a})}\right) 
	-
	\frac{\del L}{\del \phi^{a}} 
\right)
\vol \cdot \delta\phi^{a}
\end{equation}
and 
$
\{\CL\}^{\phi} : \phi^{*}T^{V}F \rightarrow \wedge^{n-1}T^{*}M
$
by
\begin{equation}\label{prestroom}
\{\CL\}^{\phi}(\delta \phi) = 
\frac{\del L}{\del(\del_{\mu}\phi^{a})}
i_{\del_{\mu}} \vol \cdot \delta \phi^{a}\,,
\end{equation}
then we have 
\begin{equation}\label{preding}
\delta S_{U}(\phi) = -\int_{U} [\CL\,]^{\phi}(\delta \phi) 
+ \int_{\del U} \{\CL \}^{\phi}(\delta \phi)\,.
\end{equation}
It can be seen from (\ref{preding}) that $[\CL\,]^{\phi}$ 
and $\{\CL\}^{\phi}$
do not depend on our choice of co-ordinates, a fact not at all 
obvious from equations (\ref{preel}) and (\ref{prestroom}). 

The argument runs as follows. Let $[\CL\,]_{A}^{\phi}$ and $[\CL\,]_{B}^{\phi}$
be two versions of (\ref{preel}), defined w.r.t\ different co-ordinate systems.
Let $\delta \phi$ be any 
section of $\phi^{*}T^V(F)$ 
that vanishes on $\del U$. 
Then the second term on the r.h.s.\ of (\ref{preding}) vanishes, and 
we have both
$\delta S_{U} (\phi) = - \int_{U} [\CL\,]_{A}^{\phi}(\delta \phi)$
and $\delta S_{U} (\phi) = - \int_{U} [\CL\,]_{B}^{\phi}(\delta \phi)$. 
This means that
$\int_{U} ([\CL\,]_{A}^{\phi} - [\CL\,]_{B}^{\phi})(\delta \phi) = 0$
for all infinitesimal variations vanishing at 
$\del U$, so that $[\CL\,]_{A}^{\phi} = [\CL\,]_{B}^{\phi}$ holds
over the interior of $U$. But since $U$ was arbitrary, equality
must hold over all of $M$.

We can now show that also $\{\CL\}^{\phi}$ does not depend
on the choice of co-ordinates. 
Let $\delta \phi$ be any section of $\phi^{*}T^V(F)$, not necessarily
vanishing at $\del U$.
According to equation 
(\ref{preding}), both $\int_{\del U}\{\CL\}_{A}^{\phi}(\delta \phi)$
and $\int_{\del U}\{\CL\}_{B}^{\phi}(\delta \phi)$ are equal to
$\delta S_{U}(\phi) + \int_{U} [\CL\,]^{\phi}(\delta \phi)$.
(We just established that $[\CL\,]^{\phi}(\delta \phi)$ is the same
for $A$ and $B$.) Once again, 
$\int_{\del U} \{\CL\}_{A}^{\phi}(\delta \phi) - \{\CL\}_{B}^{\phi}(\delta \phi)
= 0$ for arbitrary $\delta \phi$ and $U$ implies that 
$\{\CL\}_{A}^{\phi} = \{\CL\}_{B}^{\phi}$ must hold over all of $M$.

It is clear from (\ref{preding}) that $S_{U}$ is stationary at $\phi$
under variations $\delta \phi$ with support inside $U$ if and only 
if $\int_{U}[\CL\,]^{\phi}(\delta \phi)$ is zero
for all such variations. This is the case if and only if
$[\CL\,]^{\phi}$ is identically zero. \qed

\subsubsection{Example: the Klein-Gordon Equation}

One of the simplest examples of an Euler-Lagrange equation the following.
Let $M = \mathbb{R}^{4}$, and equip it with the Minkowski
metric $\eta$ of signature $(-,+,+,+)$.
In co-ordinates $(x^{0},x^{1},x^{2},x^{3})$, this reads 
$\eta (\xi,\chi) = \eta_{\mu\nu} \xi^{\mu} \chi^{\nu}$,
with $\eta_{\mu\nu}$ the diagonal matrix $\mathrm{diag}(-1,1,1,1)$.
We denote $\xi_{\mu} = \eta_{\mu\nu}\xi^{\nu}$, so that
$\xi^{\mu} = \eta^{\mu\nu}\xi_{\nu}$.

Let $F$ be the trivial bundle
$\mathbb{R}^{4} \times \mathbb{C} \rightarrow \mathbb{R}^{4}$, which 
comes with co-ordinates $(x^{\mu};v^{a})$. The index $\mu$ runs from $0$
to $3$, and $a$ runs from $0$ to $1$. (We write $z \in \mathbb{C}$ as 
$z = v^{0} + iv^{1}$.) 
We define the Lagrangian 
$$
\CL(j^{1}\phi) = \left(
\frac{1}{2m}\del^{\mu}\phi_{a} \del_{\mu} \phi^{a} - 
\frac{m}{2} \phi_{a}\phi^{a}\right) \vol\,.
$$
We have $\frac{\del L}{\del \phi^{a}} = m \phi_{a}$
and $\frac{\del L}{\del(\del_{\mu} \phi^{a})} = m^{-1} \del^{\mu} \phi_{a}$, so that
the Euler-Lagrange equation in this case is the \emph{Klein Gordon} equation
$$
\del^{\mu}\del_{\mu} \phi = m^{2} \phi\,.
$$
It allows for plane wave solutions $\phi(x) = \exp(k_{\mu}x^{\mu})$
with $k^{\mu}k_{\mu} = m^{2}$.

Note that certainly not all Euler-Lagrange equations
are linear. For instance, if 
$F$ is the bundle $\mathbb{R}^2\times \mathbb{R} \rightarrow \mathbb{R}^2$,
$\eta$ is of signature $(+,-)$, and 
$\CL(j^1\phi) = \half \del_\mu \phi \del^{\mu} \phi - (1 - \cos(\phi))$,
then  
the resulting Euler-Lagrange equation 
$\del_{\mu}\del^{\mu}\phi = - \sin(\phi)$,
rather childishly called the sine-Gordon equation, 
is a nonlinear P.D.E. with some bearing on surfaces 
of constant negative curvature.


\subsubsection{Example: Electromagnetism}
The most classical of all field theories is electromagnetism.
The base manifold $M$ is again $\mathbb{R}^{4}$ with the Minkowski metric,
but the bundle $F$ is now the cotangent bundle 
$T^{*}M \rightarrow M$. 
This is a trivial vector bundle for $M = \mathbb{R}^{4}$.
If $dx^{\nu}$ is the dual of $\del_{x^{\nu}}$ in $T^*_{x}M$,
then a global trivialization is given by 
describing 
$p_{\nu} dx^{\nu} \in T_{x}^*M$ 
with the co-ordinates $(x^{\mu},p_{\nu})$. 

The basic field in electromagnetism is a 1-form
$A \in \Gamma(T^*M)$ called the electromagnetic four-potential.
We set $p_{\nu}(A(x)) = A_{\nu}(x)$, and write
$A = A_{\mu}(x)d x^{\mu}$. We define the
antisymmetric tensor
$F_{\mu\nu} := \del_{\mu}A_{\nu} - \del_{\nu} A_{\mu}$, so that
$F = 2dA$, with $F = F_{\mu\nu} dx^{\mu}\wedge dx^{\nu}$. 
The (vacuum) EM-Lagrangian is then
$$
\CL = -\frac{1}{4} F^{\mu\nu}F_{\mu\nu} \vol\,.
$$
We have $\frac{\del L}{\del \del_{\mu}A_{\nu}} = -F^{\mu\nu}$,
and
$\frac{\del L}{\del A_{\nu}} = 0$. 
The Euler-Lagrange equation is then the (vacuum) \emph{Maxwell equation}
$$
\del_{\mu} F^{\mu\nu} = 0\,.
$$
The 3 components $E_i = F_{i0}$ are identified with the electric field,
and the 3 components $B_{i} = \half \epsilon_{ijk}F_{jk}$ with the magnetic field.
Light is described by the plane wave solutions
$A_{\mu}(x) = a_{\mu}\exp(ik_{\mu}x^{\mu})$ with $k_{\mu}k^{\mu} = 0$,
$a_{\mu} k^{\mu} = 0$ and, at least in the Coulomb gauge, $a_{0}=0$.

\subsection{Transformations and Symmetries}\label{trafosymsym}
Transformations of the system correspond to automorphisms of 
the fibre bundle.
An automorphism of $\pi : F \rightarrow M$ 
is a pair $(\alpha,\alpha_{M})$ with 
$\alpha$ a diffeomorphism of the total space $F$, $\alpha_{M}$
a diffeomorphism of the base $M$, and 
$\pi \circ \alpha = \alpha_{M} \circ \pi$.
If $F$ is a vector bundle, 
we additionally require that
$\alpha$ be linear on the fibres,
and if $F$ is a principal fibre bundle, 
we require $\alpha$ to be equivariant w.r.t.\ the group action. 
Because $\alpha_{M}$ is completely determined by $\alpha$,
we drop it from the notation, and write
$\alpha$ instead of $(\alpha,\alpha_{M})$.  
An automorphism is called \emph{vertical} if $\alpha_{M} = \id$.

The automorphisms of $F$ form the group $\mathrm{Aut}(F)$, and the vertical
automorphisms $\mathrm{Aut}^{V}(F)$ form a normal subgroup.
If we denote the diffeomorphism group of $M$ by $\mathrm{Diff}(M)$,
then $\alpha \mapsto \alpha_{M}$ is a homomorphism
from $\mathrm{Aut}(F)$ to $\mathrm{Diff}(M)$ with kernel 
$\mathrm{Aut}^{V}(F)$. Its image $\mathrm{Diff}^{F}(M)$,
the group of liftable diffeomorphisms, is an open subgroup of 
$\mathrm{Diff}(M)$.
We have the exact sequence
\begin{equation}\label{snurkluit}
1 
\rightarrow 
\mathrm{Aut}^{V}(F) 
\rightarrow
\mathrm{Aut}(F)
\rightarrow 
\mathrm{Diff}^{F}(M)
\rightarrow 
1\,.
\end{equation}

The group $\mathrm{Aut}(F)$ of transformations acts on 
the space $\Gamma(F)$ of fields by 
$\phi \mapsto \alpha \circ \phi \circ \alpha_{M}^{-1}$.
Note however that without any
additional structure, the diffeomorphism group
$\mathrm{Diff}(M)$ does \emph{not} 
have a natural action on the fields $\Gamma(F)$.
The extra structure that is required to have $\mathrm{Diff}^{F}(M)$
act on the fields $\Gamma(F)$ is precisely
a homomorphism 
$\Sigma : \mathrm{Diff}^{F}(M) \rightarrow \mathrm{Aut}(F)$
splitting (\ref{snurkluit}).

An automorphism $\alpha$ of $F$ acts on $J^{1}F$ by
$j^{1}\alpha : 
j^{1}_{x}\phi \mapsto 
j_{x}^{1} (\alpha \circ \phi \circ \alpha^{-1}_{M})$.
This means that $\alpha$ takes the Lagrangian 
$\CL : J^{1}F \rightarrow \wedge^{n}T^{*}M$
into $\alpha^{*}_{M} \circ \CL \circ j^{1}\alpha$.
\begin{basdef}\label{echtsym}
An automorphism $\alpha$ of $F$ is called a \emph{symmetry}
of the Lagrangian $\CL : J^{1}F \rightarrow \wedge^{n}T^{*}M$
if $\CL = \alpha^{*}_{M} \circ \CL \circ j^{1}\alpha$.
\end{basdef}
This is a sensible definition. If $\phi$ is a local section
of $F$ over $U$, then $\alpha \circ \phi \circ \alpha^{-1}_{M}$
is a local section over $\alpha(U)$. One readily checks that 
if $\alpha$ is a symmetry of $\CL$, then 
$S_{U}(\phi) = S_{\alpha(U)}(\alpha \circ \phi \circ \alpha^{-1}_{M})$.
In particular, $\alpha \circ \phi \circ \alpha^{-1}_{M}$ is a critical
section if and only if $\phi$ is, and the action of 
$\alpha$ on $\Gamma(F)$ restricts to the solution space
$\Gamma_{EL}(F)$ of the Euler-Lagrange
equation.  

\subsubsection{Infinitesimal Symmetries}
It will be expedient to look at this picture from the infinitesimal 
point of view. The Lie algebra of $\mathrm{Diff}(M)$ is, at least 
morally speaking,
the Lie algebra $\Gamma(TM)$ of vector fields on $M$.
The group $\mathrm{Aut}(F)$ corresponds to the Lie algebra 
$\Gamma^{P}(TF)$ of projectable vector fields on $F$, i.e.\ 
vector fields $\xi \in \Gamma(TF)$ such that 
$\pi_{*} \xi_{f} = \pi_{*} \xi_{f'}$ whenever $\pi(f) = \pi(f')$.
The pushforward $\pi_{*} : \Gamma^{P}(TF) \rightarrow
\Gamma(TM)$ is not only well defined, but also a homomorphism of Lie algebras. 
Its kernel, the ideal $\Gamma^{V}(TF)$ of vertical vector fields, 
corresponds to the group $\mathrm{Aut}^{V}(F)$, and we have the exact sequence
of Lie algebras
\begin{equation}\label{rijworst}
0 
\longrightarrow
\Gamma^{V}(TF)
\longrightarrow
\Gamma^{P}(TF)
\stackrel{\pi_{*}}{
\longrightarrow}
\Gamma(TM)
\longrightarrow
0\,.
\end{equation}

The infinitesimal variation $\delta \phi \in \Gamma(\phi^{*} T^{V}F)$ 
of $\phi$ under $\xi \in \Gamma^{P}(TF)$
is given by 
$\delta\phi_{x} = 
\frac{d}{d\varepsilon}|_{0} \alpha_{\varepsilon} \circ \phi \circ 
\alpha^{-1}_{\varepsilon,M}(x)$, 
with
$\alpha_{\varepsilon} = \exp(\varepsilon \xi)$ and 
$\alpha_{\varepsilon, M} = \exp(\varepsilon \pi_{*}\xi)$
the induced flows on $F$ and $M$ respectively.  
It is readily seen to be
$\delta \phi_{x}
 = 
d^{v}\phi (\xi_{\phi(x)})$,
with 
$d^{v}\phi (\xi) =
\xi - \phi_{*} \pi_{*}\xi$
the vertical projection 
$T_{\phi(x)}F \rightarrow T^{V}_{\phi(x)}F$ 
along the image of $\phi_* : T_{x}M \rightarrow T_{\phi(x)}F$.

In order to obtain the infinitesimal version of definition \ref{echtsym},
we require that $\frac{d}{d\varepsilon}|_{0} 
\alpha^*_{\varepsilon , M} \CL \circ j^{1}\alpha_{\varepsilon}$
vanish. %
Variation of $\alpha^*_{\varepsilon , M}$ gives rise to a term
$L_{\pi_{*} \xi} \CL(j^{1}(\phi)) = d (i_{\pi_* \xi}\CL(j^{1}\phi))$.
Variation of $\CL$ through the variation 
$\frac{d}{d \varepsilon}
\alpha_{\varepsilon} \circ \phi \circ \alpha^{-1}_{\varepsilon,M}
$ of the field gives rise to an additional term 
$ - [\CL\,]^{\phi}(\delta \phi) + d(\{\CL\}^{\phi}(\delta \phi))$,
cf.\ equation (\ref{preding}), with $\delta \phi_x = d^{v}\phi(\xi_{\phi(x)})$. 
This leads us to the following definition of infinitesimal symmetries.
\begin{basdef} \label{zooloog}
A projectable vector field $\xi \in \Gamma^{P}(TF)$ is called an 
\emph{infinitesimal symmetry} of the Lagrangian 
$\CL : J^{1}F \rightarrow \wedge^{n}T^*M$
if 
$$-[\CL\,]^{\phi}(d^{v}\phi(\xi)) +
d(i_{\pi_{*} \xi} \CL + \{\CL\}^{\phi}(d^{v}\phi(\xi)) 
 = 0\,.$$
\end{basdef}

From this definition, we immediately see that infinitesimal symmetries
give rise to conserved quantities. If we define the \emph{Noether current}
associated to the infinitesimal symmetry $\xi$ 
and the field $\phi$
to be the $(n-1)$-form
\begin{equation}\label{kanoniekestroom}
J^{\phi}(\xi) = \{\CL\}^{\phi}(d^{v}\phi(\xi)) + i_{\pi_{*}\xi} \CL(j^{1}(\phi))\,,
\end{equation}
then the following theorem follows straight from the definitions. 
\begin{theorem}[Noether]
If $\phi$ is a critical point of the action, and $\xi$ is an infinitesimal
symmetry, then $d J^{\phi}(\xi) = 0$.
\end{theorem}
\proof We have $\delta \CL = - [\CL\,]^{\phi}(d^{v}\phi(\xi)) + d J^{\phi}(\xi)$.
Now $\delta \CL = 0$ because $\xi$ is an infinitesimal symmetry, and 
$[\CL\,]^{\phi} = 0$ by the Euler-Lagrange equation.\qed
The rationale behind considering closed $(n-1)$-forms as conserved
quantities is the following.
Let $t$ be a (local) time co-ordinate, and let $\Sigma$ be a 
compact $(n-1)$-dimensional submanifold with boundary of the (local)
time-slice $t = t_{0}$.
If we interpret $Q_{\Sigma} = \int_{\Sigma} J^{\phi}(\xi)$
as the `charge' of $\Sigma$, then $\frac{d}{dt}Q_{\Sigma}
= \int_{\Sigma}\CL_{\del_{t}}J^{\phi}(\xi)$, with
$\CL_{\del_{t}}$ the Lie derivative along $\del_{t}$. Now 
$\CL_{\del_{t}}J^{\phi}(\xi) = d i_{\del_{t}}J^{\phi}(\xi)$
because $d J^{\phi}(\xi) = 0$, and 
$\int_{\Sigma} d i_{\del_{t}}J^{\phi}(\xi) = 
\int_{\del \Sigma}i_{\del_{t}}J^{\phi}(\xi)$ by Gau\ss'
theorem. If we now interpret
$
\Phi_{\del \Sigma} = \int_{\del \Sigma}i_{\del_{t}}J^{\phi}(\xi)\,
$
as the flux through $\del{\Sigma}$,
then 
%
$\frac{d}{dt} Q_{\Sigma} = \Phi_{\del \Sigma}$
says that the rate of change
of the charge of $\Sigma$ is equal to
the flux through its boundary. 

Strictly speaking then, $Q_{\Sigma}$ is 
not a conserved quantity at all; we merely know that its rate of change 
can be attributed to a flux. However, 
if the current goes to zero fast enough at spatial infinity,
then the total `charge' of the time-slice $t = t_0$ is conserved.

\subsubsection{Example: the Current of a Klein-Gordon Field}

The Klein-Gordon Lagrangian 
allows the one-parameter group of vertical rotations
$\alpha_{\theta} : (x^{\mu}; v^{0}, v^{1}) \mapsto
(x^{\mu} ; \cos(\theta)v^{0} + \sin(\theta)v^{1} , -\sin(\theta)v^{0}
 + \cos(\theta)v^{1} )$
as a group of symmetries. The corresponding infinitesimal
symmetry is the vector field 
$\xi = v^{1}\frac{d}{dv^{0}} - v^{0}\frac{d}{dv^{1}}$.
We calculate the corresponding current from equation (\ref{kanoniekestroom}). 
Because $\xi$ is vertical, $\pi_{*}\xi = 0$
and $d^{v}\phi(\xi) = \xi$. 
Since $\{\CL\}_{a}^{\phi} = \frac{1}{m}\del^{\mu}\phi_{a} i_{\del_{\mu}}\vol$,
we get $J^{\phi}(\xi) = \frac{1}{m}
(\phi^{0}\del^{\mu}\phi^{1} - \phi^{1}\del^{\mu}\phi^{0})i_{\del_{\mu}}\vol$.
It is not hard to check `by hand' that $\del_{\mu}J^{\mu}(\xi) = 0$ for solutions
of the Klein-Gordon equation, as predicted by Emmy Noether.

\subsubsection{Transformations and Symmetries of Space-Time}

We would like to associate
conserved quantities to
infinitesimal symmetries of space-time.
However, in order to even define what a space-time symmetry is, 
we need to know how $\mathrm{Diff}(M)$ acts on the fields.
In other words, we need a homomorphism 
$\Sigma : \mathrm{Diff}^{F}(M) \rightarrow \mathrm{Aut}(F)$ that splits
(\ref{snurkluit}).
A transformation $\alpha_{M} \in \mathrm{Diff}(M)$ is then called a 
space-time symmetry
of $\CL$
if $\Sigma (\alpha_{M})$ is a symmetry in the sense of definition 
\ref{echtsym}. 

An infinitesimal transformation of space-time is a vector field
$\xi_{M} \in \Gamma(TM)$. 
In order to describe its action on the fields, we need a splitting
$\sigma : \Gamma(TM) \rightarrow \Gamma^{P}(TF)$ of 
the exact sequence (\ref{rijworst}) of Lie-algebras. A vector field
$\xi_{M} \in \Gamma(TM)$ is called an infinitesimal space-time
symmetry of $\CL$ if $\sigma(\xi_{M})$ is a symmetry in the sense of
definition \ref{zooloog}.
    
This Lie algebra homomorphism 
$\sigma : \Gamma(TM) \rightarrow \Gamma^{P}(TF)$ will be 
the object of our attention in chapter \ref{ch:BLID}.
We will argue there that it is of more fundamental importance than
its counterpart 
$\Sigma : \mathrm{Diff}(M) \rightarrow \mathrm{Aut}(F)$
at the level of groups.

Although every flat connection
gives rise to a splitting of (\ref{rijworst}),
it is certainly not true that every splitting arises
in this way.
However, in chapter \ref{ch:BLID}, proposition \ref{alhetvet} and 
corollary \ref{structuurcor}, we will prove that under 
reasonable assumptions\footnote{Namely, we assume that
the fibre bundle $F$ is associated to a principal fibre bundle, 
and that the infinitesimal transformations respect this association.},
$\sigma$ must be a differential operator of finite order $k$. 

This means that 
the map of sections
$\sigma : \Gamma(TM) \rightarrow \Gamma^{P}(TF)$
comes from a map
$\nabla : \pi^{*}J^{k}TM \rightarrow TF$
of bundles over $F$.
In local co-ordinates, 
there then exist functions $C^{a; \rho_1 \ldots \rho_j}_{\tau}(x^{\mu} , v^{a})$
such that
$$\sigma(\xi^{\mu}\del_{\mu}) 
= \xi^{\mu}\del_{\mu} + 
\left( C^{a}_{\tau}\xi^{\tau} +
C^{a; \rho}_{\tau} \del_{\rho}\xi^{\tau} +
\cdots
+ 
C^{a; \rho_1 \ldots \rho_{k}}_{\tau} \del_{\rho_{1}} 
\ldots \del_{\rho_{k}}\xi^{\tau}
\right) \del_{a}\,.
$$
Note that the case $k=0$ corresponds to a flat connection.

\subsubsection{The Stress-Energy-Momentum Tensor}
Once we have a splitting $\sigma$, it is tempting to define 
the associated conserved current as $J^{\phi}(\sigma(\xi_{M}))$.
Appealing as this may seem, it is \emph{not} 
the proper course of action in the presence of a sufficiently 
large group of space-time symmetries, as noted first in \cite{No}. 

Instead, one can use $\sigma$ in order to construct a 
\emph{Stress-Energy-Momentum (SEM) tensor} as follows.
%
Choose adapted co-ordinates, and expand
$$
\delta \CL = \left(
\del_{\mu}(\xi^{\mu} L) + 
\frac{\del L}{\del \phi^{a}} \delta \phi^{a} 
+
\frac{\del L}{\del \del_{\mu}\phi^{a}} \del_{\mu} \delta \phi^{a} 
\right)\vol
$$
into 
\begin{equation}\label{smoelwerk}
\delta \CL = T_{\tau} \xi^{\tau}\vol + 
{T^{\rho}}_{\tau}\del_{\rho}\xi^{\tau}\vol + \ldots
+ T^{\rho_1 \ldots \rho_{k+1}}_{\tau} \del_{\rho_{1}}
\ldots \del_{\rho_{k+1}}\xi_{\tau}\vol\,,
\end{equation} 
using 
$\delta \phi = d^{v}\phi(\sigma(\xi_{M}))$, i.e.\
$$\delta \phi^{a} = (C^{a}_{\tau} - \del_{\tau}\phi^{a})\xi^{\tau} 
+
C^{a;\rho}_{\tau}\del_{\rho}\xi^{\tau} + \ldots +
C^{a;\rho_1,\ldots,\rho_{k}}_{\tau}\del_{\rho_1} \ldots
\del_{\rho_{k}}\xi^{\tau}\,.
$$

The tensor density ${T^{\rho}}_{\tau}(\phi)\vol$ is 
then the desired SEM-tensor. 
If $\Sigma$ is a compact $(n-1)$-dimensional submanifold with boundary
of the equal-time slice $t = t_{0}$, then for $\tau$ equal to zero,
the quantity
$P_{\tau} = \int_{\Sigma} {T^{\rho}}_{\tau}(\phi)i_{\del_{\rho}}\vol$
is interpreted as the energy of the fields
in $\Sigma$.
For $\tau$
equal to $1$ through $n-1$,
it is interpreted as the momentum.

The problem with this definition is of course that the tensor 
density ${T^{\rho}}_{\tau}(\phi)\vol$ may well depend on the choice
of co-ordinates. However, if $\CL$ is diffeomorphism invariant, 
then $\delta \CL = 0$ in equation $(\ref{smoelwerk})$ for any choice
of $\xi \in \Gamma(TM)$, forcing all the 
$T^{\rho_{1}\ldots\rho_{s}}_{\tau}(\phi)$ to vanish identically.  

\subsubsection{Example: the SEM-Tensor of Electromagnetism}

The tangent bundle comes with a natural lift of 
$\mathrm{Diff}(M)$. For any diffeomorphism $\alpha$, 
the pushforward $\alpha_{*}$ maps $T_{x}M$ to $T_{\alpha(x)}M$, 
and is therefore an automorphism of $TM$ covering $\alpha$. The lift 
$\Sigma(\alpha) = \alpha_{*}$ is a homomorphism of groups.

Similarly, the exact sequence (\ref{snurkluit}) for the 
cotangent bundle $T^*M$ is split by
the inverse pullback $\Sigma(\alpha) = (\alpha^{-1})^{*}$,
considered as a map
$\mathrm{Diff(M)} \rightarrow \mathrm{Aut}(T^*M)$. 
Its infinitesimal version is the Lie algebra 
homomorphism 
$\sigma(\xi) = \frac{d}{d\varepsilon}|_{0} \exp(-\varepsilon \xi)^{*}$,
which
reads
$\sigma(\xi^{\mu} \del_{\mu}) = 
\xi^{\mu} \del_{\mu} - p_{\sigma} \del_{\nu}\xi^{\sigma} \frac{\del}{\del p_{\nu}}$
in local co-ordinates.

The infinitesimal variation $\delta A \in \Gamma(A^{*}(T^{V}T^{*}M))$
of $A \in \Gamma(T^*M)$
along $\sigma(\xi)$ 
reads $\delta A = d^{v}A(\sigma(\xi)) = \sigma(\xi)_{A(x)} - A_* \xi$. 
In local co-ordinates, this is $\delta A = 
(\xi^{\mu}\del_{\mu} - A_{\sigma}\del_{\nu}\xi^{\sigma}\frac{\del}{\del p_{\nu}})
-
(\xi^{\mu}\del_{\mu} + \xi^{\sigma}\del_{\sigma} A_{\nu} \frac{\del}{\del p_{\nu}})
$,
i.e.\
$
\delta A = - 
(
A_{\sigma}\del_{\nu}\xi^{\sigma}
+
\xi^{\sigma}\del_{\sigma} A_{\nu}) \frac{\del}{\del p_{\nu}}
$.
If we identify $A^{*}(T^{V}T^{*}M)$ with $T^*M$, which we can do
because $T^*M$ happens to be a vector bundle, then we simply have 
$\delta A = -\CL_{\xi} A$.

Now that we know that action of $\Gamma(TM)$ on $\Gamma(T^*M)$
is the Lie derivative \mbox{-- What} else could it have \mbox{been? --} 
calculating the SEM-tensor is straightforward.
We expand  
$\delta \CL = -(F^{\mu\nu}\del_{\mu}\delta A_{\nu} - 
\frac{1}{4}\del_{\tau}(F^{\mu\nu}F_{\mu\nu}\xi^{\tau}))\vol$
into derivatives of $\xi$.
\begin{eqnarray*} 
\delta \CL &=& 
\left(
F^{\mu\nu}\del_{\mu}(\xi^{\tau}\del_{\tau}A_{\nu} + 
A_{\tau} \del_{\nu}\xi^{\tau}) - \fourth\del_{\tau}
(F^{\mu\nu}F_{\mu\nu}\xi^{\tau})
\right) \vol \\
&=&
(
F^{\mu\nu}\del_{\mu}\del_{\tau}A_{\nu} - \fourth \del_{\tau}
(F^{\mu\nu}F_{\mu\nu})
) \xi^{\tau} \vol \\
& & +\,\,
\left(
F^{\mu\nu}(\del_{\mu} \xi^{\tau} \del_{\tau}A_{\nu} 
+ \,\, \del_{\mu} A_{\tau} \del_{\nu} \xi^{\tau})
-
\fourth F^{\mu\nu}F_{\mu\nu} \del_{\tau}\xi^{\tau} 
\right)
\vol\\
& & 
+\,\, A_{\tau}F^{\mu\nu}\del_{\mu}\del_{\nu}\xi^{\tau}\vol\\
&=&
- \,\,\fourth \del_{\tau}
(F^{\mu\nu}F_{\mu\nu})
\xi^{\tau} \vol \\
& &+\,\,
\left(F^{\mu\nu}F_{\tau\nu} 
- \fourth F^{\mu\nu}F_{\mu\nu}{\delta^{\mu}}_{\tau}
\right)\del_{\mu}\xi^{\tau} \vol
\end{eqnarray*}
In the last step, we used that $F_{\mu\nu}$ is antisymmetric,
and that $F^{\mu\nu}\del_{\mu}A_{\tau}\del_{\nu}\xi^{\tau}
= - F^{\mu\nu}\del_{\nu}A_{\tau}\del_{\mu} \xi^{\tau}$.
According to our definition, the SEM-tensor of electromagnetism
is 
\begin{equation}\label{emstress}
{T^{\mu}}_{\tau}(A) = 
F^{\mu\nu}F_{\tau\nu} 
- \fourth F^{\mu\nu}F_{\mu\nu}{\delta^{\mu}}_{\tau}\,.
\end{equation}
In particular, the energy density $T^{00}$ of the EM-field
is $\half(|E|^{2} + |B|^{2})\vol$. 

\basremark{The above tensor density is co-ordinate independent, 
which, in our limited setup, is a small miracle. 
As ever so often, the miracle disappears if we broaden our view a little.
The transformation behaviour $\Sigma : \mathrm{Diff}(M) \rightarrow 
\mathrm{Aut}(F)$ of the EM-field allows us to couple it to a metric
$g^{\mu\nu}$, resulting in a Lagrangian that equals the one 
above for $g = \eta$, but is $\mathrm{Diff}(M)$-invariant because
the transformation of the metric cancels that of the fields. 
The SEM-tensor then vanishes identically. It is the sum of a term 
resulting from $\delta A$, and an equal term of opposite sign
resulting from $\delta g$. The former term was calculated above, 
and the latter term can be seen to represent a bona fide tensor density
from general principles. See e.g.\ \cite{FR} for details.}
\basremark{Our definition of the SEM-tensor differs from the one in
\cite{GM} and \cite{FR} by a term $[\CL\,]_{a}C^{a;\mu}_{\tau}$
which vanishes on shell, 
i.e.\ for solutions of
the Euler-Lagrange equations.
Consequently, the above expression for the SEM-tensor is correct
on shell as well as off shell, whereas the derivation in \cite{GM}
is valid only on shell.
The term $[\CL\,]_{a}C^{a;\mu}_{\tau}$ does appear in Emily Noether's 
seminal paper \cite{No}, 
but her aim is slightly different from ours. She expresses 
$d(J^{\mu}(\xi) - [\CL\,]_{a}C^{a\mu}_{\tau}\xi^{\tau}\vol )$ 
in terms of $[\CL\,]^{\phi}$
and its derivatives, showing that the Euler-Lagrange
equations are linearly dependent in the presence of symmetry.}

\subsection{Infinitesimally Natural Bundles}

The point of deriving (\ref{emstress})
is not to acquaint the reader of the precise form of the
electromagnetic SEM-tensor. It will not occur
in this thesis again, and could be found
in a myriad of introductory books on electrodynamics if it did.   
The point is to stress the crucial role that the 
splitting $\sigma : \Gamma(TM) \rightarrow \Gamma^{P}(TF)$
of (\ref{rijworst}) plays in its derivation. 

Note that in order to derive equation (\ref{emstress}), we needed to use 
the nontrivial splitting derived from the pullback, 
even though 
$T^{*}\mathbb{R}^{4} \rightarrow \mathbb{R}^{4}$ is trivial 
as a vector bundle!  
Had we used the trivial splitting 
$\hat{\sigma}(\xi^{\mu}\del_{\mu}) = \xi^{\mu}\del_{\mu}$ that comes from 
the global trivialization of $T^*M$, we would have obtained the 
\emph{incorrect} SEM-tensor
$$
{T^{\mu}}_{\tau} = F^{\mu\nu}\del_{\tau}A_{\nu} - \fourth F^{\mu\nu}F_{\mu\nu}
{\delta^{\mu}}_{\tau}\,,
$$  
which differs from the correct one by a term 
$d(F^{\mu\nu}A_{\tau}i_{\del_{\nu}}\vol) 
- (\del_{\nu}F^{\mu\nu})A_{\tau}\vol$.
This would increase the electromagnetic energy of a
$3$-volume $\Sigma \subset \mathbb{R}^{3}$ by a 
boundary term $\int_{\del \Sigma} A_{0}\vec{E} \cdot d\vec{S}$,
which is
accessible to experiment. It is known to be incorrect. 

We conclude that a classical field theory is not fully
specified by a fibre bundle $F \rightarrow M$
and a Lagrangian density $\CL : J^{1}F \rightarrow \wedge^{n}T^{*}M$.
In order to derive the correct SEM-tensor, one has to specify 
the reaction of the fields to infinitesimal space-time
transformations, that is one has to specify a Lie algebra homomorphism
$\sigma : \Gamma(TM) \rightarrow \Gamma^{P}(TF)$ splitting
the exact sequence (\ref{rijworst}). 

Even if one has no interest 
in the SEM-tensor whatsoever, 
this is a very natural thing to ask for.
It is hard to imagine 
fields on space-time without a well
defined transformation behaviour under infinitesimal 
space-time transformations.  
For one thing,
such fields would have ill-defined spin. 

We therefore consider the
splitting
$\sigma : \Gamma(TM) \rightarrow \Gamma^{P}(TF)$
as part of the input data of a classical field theory,
on equal footing with the bundle $\pi : F \rightarrow M$
and the Lagrangian $\CL : J^{1}F \rightarrow \wedge^{n}T^*M$. 
This poses an obvious restriction on the fibre bundle $F \rightarrow M$,
namely that 
it has to admit at least one such splitting. 
We will call such fibre bundles \emph{infinitesimally natural}.  
In chapter
\ref{ch:BLID}, we will classify these infinitesimally natural
bundles-- or at least the ones with sufficiently 
structured fibres.

\chapter{Information Transfer Implies State Collapse}
\label{ch:ITISC}
\setsubdir{hoofdstukken/hoofdstuk1/}



We attempt to clarify certain puzzles concerning 
state collapse and decoherence.
In open quantum systems decoherence is shown to be a necessary
consequence of the transfer of information to the outside;
we prove an upper bound for the amount of coherence which can survive
such a transfer.
We claim that in large closed systems decoherence has never been observed,
but we will show that it is usually harmless to assume its occurrence.
An independent postulate of state collapse over and above
Schr\"odinger's equation and the probability interpretation
of quantum states, is shown to be redundant.

\section{Introduction}
In its most basic formulation,
quantum theory encodes the preparation of a system in a pure quantum state,
a unit vector $\psi$ in a Hilbert space $\CH$. 
Observables are modelled by (say, nondegenerate) self-adjoint operators
on $\CH$.
The expectation value of an observable $A$ in a state $\psi$
is given by $\inp{\psi}{A \psi}$. 
If $a$ is an eigenvalue of $A$ and $\psi_a$ a unit eigenvector,
and information concerning $A$ is somehow extracted from the system,
then the probability for the value $a$ to be observed
is $|\inp{\psi_a}{\psi}|^2$.
If this observation is indeed made,
then the subsequent behaviour of the system is predicted using
the pure state $\psi_a$.
This is called \emph{state collapse}.
It follows that,
if the information extraction has taken place 
but the information on the value of $A$ is disregarded,
then the subsequent behaviour can be described optimally using a mixture
of eigenstates.
This is called \emph{decoherence}. 
In this chapter we substantiate the following claim concerning
decoherence and state collapse.
\begin{center}
\emph{Decoherence is only observed in open systems, where it is a necessary\\
consequence of the transfer of information to the outside.\\}
\end{center}

So the observed occurrence of decoherence does not contradict the unitary
time evolution postulated by quantum mechanics,
since open systems do not evolve unitarily.
Decoherence can be explained in quantum theory
by embedding the quantum system into a larger, closed whole,
which in itself evolves unitarily.
This is well-known (see e.g. \cite{Ne}).
We add the observation that
decoherence is not only a {\em possibility} for an open system,
but a {\em necessary consequence}
of the leakage of information out of the system.
We prove an inequality relating the decoherence between two pure states
to the degree in which a decision between the two is possible
by a measurement outside.
This is the content of theorem~\ref{cor5} in section~\ref{open}.

Also,
we have claimed that one {\it does not} actually observe decoherence
in closed macroscopic systems.
First of all, 
most of the systems that are ever observed are actually open,
since it is extremely difficult to shield large systems from
interaction.
But more to the point,
the difference between coherence and decoherence can only 
be seen by measuring some highly exotic `stray observables'
which are almost always forbiddingly hard to observe.
And indeed,
in those rare cases where experimenters have succeeded in measuring them,
ordinary unitary evolution was found, not decoherence. 
(See \cite{Arn}, \cite{Fri}, \cite{Wal}.)

We illustrate the latter point in section~\ref{oester},
where we show that the measurement of two classes of observables 
\emph{can not} reveal the difference between coherence and decoherence:
a class of microscopic observables and a large class of macroscopic
observables.
Take as an example a volume of gas.
Microscopic observables such as the position of one particular atom in a gas, 
only relate to a small fraction of the system.
Macroscopic observables like the center of mass of the gas,
are the average over a large number of microscopic
observables.
Belonging neither to the macroscopic nor to the microscopic class, 
the `stray observables'
referred to above describe detailed correlations between 
large numbers of atoms in the gas. 
This kind of information is  experimentally almost inaccessible.

Driving home our point concerning decoherence in closed systems:
coherent superpositions of macroscopically distinguishable states
are not the strange monsters produced by a quantum theory applied
outside its domain.
They are, on the contrary, everyday occurrences which,
however, {\it can not} be distinguished
from the more classical incoherent superpositions in practice,
and can therefore always be regarded as such.


\section{Abstract Information Extraction}\label{perse}
Quantum phenomena are inherently stochastic.
This means that, if quantum systems are prepared in identical ways,
then nevertheless different events may be observed.
A quantum state describes an ensemble of
physical systems, e.g. a beam of particles,
and is modelled by a normalized 
trace-class operator $\rho$ on the Hilbert space.
The expectation value of an observable $A$ in the state $\rho$ 
is then $\tr(\rho A)$.

An information extraction or measurement on a quantum state is to be
considered as the partition of such an ensemble into subensembles,
each subensemble corresponding to a measurement outcome.
Let us, in the present section,
not wonder {\it how} the splitting of ensembles can be described 
by quantum theory,
but let us see what such an information extraction, if it can be done,
will entail for the subsequent behaviour of the subensembles.
Note that this process may serve as part of the preparation
for further experiments on the system, so that it must again lead to a state.

\subsection{Information Extraction}\label{secIESC}
For simplicity let us assume that only two outcomes can occur, labelled 
0 and 1, say with probabilities $p_0$ and $p_1$.
The ensemble is then split in two parts,
described by their respective states $\rho_0$ and $\rho_1$.
The map
\begin{equation}
M:\rho\mapsto  p_0 \rho_0\oplus p_1 \rho_1
\label{Measure}
\end{equation}
must be normalized, affine and positive.
Indeed, normalization is the property 
that $p_0 + p_1 = 1$, and 
positivity is the requirement that states must be mapped to states.
The affine property entails that for all states $\rho$ and $\theta$
on the original system, and for all $\lambda\in[0,1]$,
   $$M(\lambda\rho+(1-\lambda)\vartheta)=\lambda M(\rho)+(1-\lambda)M(\vartheta)\;.$$
This follows from the physical principle that a system which is prepared
in the state $\rho$ with probability $\lambda$ and in the state 
$\vartheta$ with probability
$1-\lambda$, say by tossing a coin,
can not be distinguished from a physical system in the state
$\lambda\rho+(1-\lambda)\vartheta$.
We emphasize that indeed this is a {\it physical} principle,
not a matter of definitions.
It states, for instance, that a bundle of particles
having 50\% spin up and 50\% spin down can not be distinguished from
a bundle having 50\% spin left and 50\% spin right.
This is a falsifyable statement.

\subsection{State Collapse}
The above elementary observations are sufficient to prove that information
extraction implies state collapse.
If $M$ distinguishes perfectly between
the pure states $\psi_0$ and $\psi_1$, then of course
$p_0=1$ in case $\rho=\ketbra{\psi_0}$, and
$p_1=1$ if $\rho=\ketbra{\psi_1}$.
\begin{proposition}\label{abscol}
Let $\CT(\CH)$ \hspace{-0.5mm}denote\hspace{-0.15mm} the\hspace{-0.15mm} space\hspace{-0.15mm} 
of\hspace{-0.15mm} trace\hspace{-0.15mm} class\hspace{-0.15mm} operators\hspace{-0.15mm} on\hspace{-0.2mm} 
a\hspace{-0.15mm} Hilbert space
$\CH$,
and let the map
$M:\CT(\CH)\to\CT(\CH)\oplus\CT(\CH):\rho\mapsto M_0(\rho)\oplus M_1(\rho)$
be the linear extension of some normalized, affine and positive map
on the states.
Suppose that unit vectors $\psi_0$ and $\psi_1$ exist such that
\begin{eqnarray*}
     M\bigl(\ketbra{\psi_0}\bigr) & = & M_0\bigl(\ketbra{\psi_0}\bigr)\oplus0\\
     M\bigl(\ketbra{\psi_1}\bigr) & = & 0\oplus M_1\bigl(\ketbra{\psi_1}\bigr)\;.
\end{eqnarray*}
Then we have $M\bigl( \ket{\psi_0} \bra{\psi_1} \bigr)
= M\bigl( \ket{\psi_1} \bra{\psi_0} \bigr) = 0$.
\end{proposition}
\proof 
The positivity of $M$ yields 
$M\bigl( \ketbra{\varepsilon e^{i\varphi} \psi_0 + \psi_1} \bigr) \geq 0$
as an operator inequality.
In particular, the 0-th component must be positive.
As $M_0 (\ketbra{\psi_1}) =0$, it follows that for all
$\varepsilon, \phi \in \RR$, we have
$$\varepsilon^2 M_0 \bigl( \ketbra{\psi_0}\bigr) + 
\varepsilon \left( e^{i\varphi} M_0 \bigl( \ket{\psi_0}\bra{\psi_1}\bigr) + 
e^{-i\varphi} M_0 \bigl( \ket{\psi_1}\bra{\psi_0}\bigr) \right) 
\geq 0\,.$$
Taking the limit $\varepsilon \downarrow 0$ yields  
$\left( e^{i\varphi} M_0 \bigl( \ket{\psi_0}\bra{\psi_1}\bigr) + 
e^{-i\varphi} M_0 \bigl( \ket{\psi_1}\bra{\psi_0}\bigr) \right) 
\geq 0$ for all $\varphi \in \RR$. In particular for 
$\varphi = 0,\frac{\pi}{2},\pi,\frac{3\pi}{2}$, implying that
both
$M_0 \bigl( \ket{\psi_0}\bra{\psi_1} \bigr)$ and 
$M_0 \bigl( \ket{\psi_1}\bra{\psi_0} \bigr)$ are equal
to zero.
Exchanging the roles of $\psi_0$ and $\psi_1$ in the argument above
shows that also
$M_1 \bigl( \ket{\psi_0}\bra{\psi_1} \bigr)$ and 
$M_1 \bigl( \ket{\psi_1}\bra{\psi_0} \bigr)$ vanish, 
proving the proposition.
\qed
We may draw two conclusions from proposition~\ref{abscol}.
The first is that, 
for all $\ket{\psi}=\alpha_0\ket{\psi_0}+\alpha_1\ket{\psi_1}$, we have 
\begin{equation}\label{ITISC:collapse}
(M_0 + M_1)\bigl(\ketbra{\psi}\bigr) = (M_0 + M_1)\bigl(\Col \bigr)\,.
\end{equation}
In words: for the prediction of events {\it after} the splitting
of the ensemble in two,
it no longer matters whether {\it before} the splitting
the system was in the pure state
$\ketbra {\alpha_0 \psi_0 + \alpha_1 \psi_1}$ 
or in the mixed state 
$\abss{\alpha_0} \ketbra{\psi_0} + \abss{\alpha_1} \ketbra{\psi_1}$.
This phenomenon, which is a direct consequence of the structure
(\ref{Measure}) of the measurement process,
we will call {\it decoherence}.
   
The second conclusion from proposition~\ref{abscol} is the following.
For all $\ket{\psi}=\alpha_0\ket{\psi_0}+\alpha_1\ket{\psi_1}$, we have 
   \begin{equation}M \bigl(\ketbra\psi\bigr)=
        \abss{\alpha_0} M_0 \bigl(\ketbra{\psi_0} \bigr) \oplus 
		\abss{\alpha_1} M_1 \bigl(\ketbra{\psi_1}\bigr)\;.
   \end{equation}
In words: if an ensemble is split in two parts, then the `0-ensemble'
will further behave as if the system had been in state $\psi_0$ instead of
$\psi$ prior to splitting, and the 
`1-ensemble' as if it had been in state $\psi_1$ instead of $\psi$.   
This phenomenon will be called {\it collapse}.

Throughout this chapter, we will maintain a
sharp distinction between the collapse $M : \CT(\CH) \to \CT(\CH) \oplus \CT(\CH)$ 
and the decoherence $(M_0 + M_1) : \CT(\CH) \to \CT(\CH) $.
The former represents the splitting of an ensemble in two parts by means of
measurement, whereas the latter represents the splitting and subsequent 
recombination of this ensemble.  

\section{Open Systems} \label{open}
A decoherence-mapping $(M_0 + M_1)\,\colon\,\CT(\CH)\to\CT(\CH)$
maps the pure (vector) state $\ketbra{\alpha_0 \psi_0 + \alpha_1 \psi_1}$
and the mixed state 
$|\alpha_0|^2\ketbra{\psi_0} + |\alpha_1|^2\ketbra{\psi_1}$ to the same final state.
Since unitary maps preserve purity, 
there can not exist 
a unitary map $U:\CH\to \CH$ such
that for all $\rho\in\CT(\CH)$:
   $$(M_0 + M_1)(\rho)=U\rho U^*\;.$$
However, according to Schr\"odinger's equation the development 
of a closed quantum system is given by a unitary operator.
We conclude that the 
decoherence (\ref{ITISC:collapse}) 
is impossible in a closed system.
On the other hand decoherence is
a well known and experimentally confirmed phenomenon. 

We will therefore consider open systems,
i.e. quantum systems which do not obey the Schr\"odinger equation, but are part of a
larger system which does. 
It has often been pointed out (e.g. \cite{Ne}, \cite{Zu})
that decoherence can well occur in this situation,
provided that states are only evaluated on the observables of the smaller system.
We are more ambitious here: 
we shall prove that this form of `local' decoherence is not just a
\emph{possible}, but an an \emph{unavoidable} consequence
of information-transfer out of the open system.


\subsection{Unitary Information Transfer and Decoherence}\label{ITSC}
We assume that the open system has Hilbert space $\CH$,
and that its algebra of observables is given by $B(\CH)$,
the bounded operators on $\CH$.
We may then assume that the larger system has Hilbert space $\CK \ten \CH$,
since the only way to represent $B(\CH)$ on a Hilbert space
is in the form $A\mapsto\one\ten A$ \cite{Ta}.
We may think\footnote{Sometimes it may happen,
as for instance in fermionic systems,
that the observables of the ancilla do not all commute
with those of the open system.
Also the observable algebra on $\CK$ may be smaller than
$B(\CK)$,
but we will neglect these complications here.} of 
$B(\CK)$ as the observable algebra of some ancillary system 
in contact with our open quantum system. 
In this context, $\CH$ will be referred to
as the `open system', $\CK$ as the `ancilla' and $\CK \ten \CH$
as the `closed system'.

We couple the system to the ancilla during a finite time interval $[\, 0, t \,]$.
Let $\tau \in \CT(\CK)$ denote the state of the ancilla at time 0,
and $\rho  \in \CT(\CH)$ that of the small system. 
The effect of the interaction is described by a unitary operator
$U:\CK\ten\CH\to\CK\ten\CH$,
and the state of the pair at time $t$ is given by
the density operator 
$U(\tau \ten \rho)U^* \in \CT(\CK \ten \CH)$.
For convenience, we will define the information transfer map $T : \CT(\CH) \to \CT(\CK\ten\CH)$ by
$ T(\rho) := U(\tau\ten\rho) U^*\,.$

\subsubsection{Decoherence}\label{zwaan}
In the above setup, we are interested in distinguishing whether the open system $\CH$
was in state $\ket{\psi_0}$ or $\ket{\psi_1}$ at time 0.
This can be done if
there exists a `pointer observable' $B\ten \one$ in the ancilla $B(\CK)$
which takes average value $b_0$ in state $T(\ketbra{\psi_0})$ 
and $b_1$ in state $T(\ketbra{\psi_1})$. 
By looking only at the ancilla $\CK$ at time $t$, we are then able to 
gain information on the state of the open system $\CH$ at time 0. 
We say that information is \emph{transferred} from $\CH$ to $\CK$.

Under these circumstances, we wish to prove that 
decoherence occurs on the open
system.  
We prepare the ground by proving the following lemma.

\begin{lemma}\label{lem4}.
Let $\vartheta_0,\vartheta_1$ be unit vectors in a Hilbert space $\CL$,
and let $A$ and $B$ be bounded self-adjoint operators on $\CL$ satisfying
$\norm{[A,B]}\le\delta\norm{A}\cdot\norm{B}.$
For $j=0$ or $1$,
let $b_j := \inp{\theta_j}{B\theta_j}$ denote the expectation and 
$\s_j^2 := \inp{\theta_j}{B^2 \theta_j} - \inp{\theta_j}{B\theta_j}^2$ 
the variance of $B$ in the state $\vartheta_j$.
Then, if $b_0\ne b_1$,
 $$\bigl|\inp{\vartheta_0}{A\vartheta_1}\bigr|
   \le \frac{\delta\norm{B}+\s_0+\s_1}{|b_0-b_1|}\norm A\;.$$
\end{lemma}
\proof 
We write $|(b_0-b_1)\inp{\vartheta_0}{A\vartheta_1}|$
as 
$|\inp{\vartheta_0}{\bigl(A(B-b_1)-(B-b_0)A+[B,A]\bigr)\vartheta_1}|$\,.
Since
$\norms{(B-b_j)\vartheta_j}=\inp{\vartheta_j}{(B-b_j)^2\vartheta_j}=\s_j^2$,
we can then use the Cauchy-Schwarz inequality to show that
$$\bigl|(b_0-b_1)\inp{\vartheta_0}{A\vartheta_1}\bigr|
     \le\norm A(\s_1+\s_0)+\delta\norm A\cdot\norm B\;.\vspace{-3mm}$$
\qed
Note that, for $\delta = \sigma_0 = \sigma_1 = 0$, lemma~\ref{lem4}
merely states that commuting operators respect each other's eigenspaces.
We proceed to prove that information transfer causes 
decoherence on the open system. (See \cite{scriptie}.)
\begin{theorem}\label{cor5}
Let $\psi_0$ and $\psi_1$ be mutually orthogonal unit vectors in a Hilbert
space $\CH$,
and let $\tau \in \CT(\CK)$ be a state on a Hilbert space $\CK$.
Let $U : \CK \ten \CH \rightarrow \CK \ten \CH$ be unitary and define
$T : \CT(\CH) \rightarrow \CT(\CK \ten \CH)$ by
$T(\rho) = U (\tau \ten \rho) U^* $.
Let $B$ be a bounded self-adjoint operator on $\CK\ten\CH$, 
and denote by $b_j$ 
and $\s_{j}^{2}$  its expected value  and variance 
in the state $T(\ketbra{\psi_j})$ for $j = 0,1$.
Suppose that $b_0 \neq b_1$. 
Then for all $\psi=\alpha_0\psi_0+\alpha_1\psi_1$ with $|\alpha_{0}|^2 + 
|\alpha_{1}|^2 = 1$
and for all bounded self-adjoint operators $A$ on $\CK \ten \CH$ such that
$\norm{[A,B]}\le\delta\norm{A}\cdot\norm{B}$, we have
 $$\Bigl| \tr \bigl( T(\ketbra{\psi}) A \bigr) -
 \tr \Bigl( T\bigl( \Col ) \bigr) A \Bigr) \Bigr| \le$$
\begin{equation}\label{margel}
 \frac{\delta\norm B+\s_0+\s_1}{|b_0-b_1|}\norm A \,.
\end{equation}
 
\end{theorem}
\proof 
First, we prove (\ref{margel}) in the special case that $\tau = \ketbra{\varphi}$ 
for some vector $\varphi \in \CK$. 
We introduce the notation $\theta_j := U(\varphi\ten\psi_j)$. 
Recall that the expectation of $B$ is given by
$ b_j = \tr \bigl( T(\ketbra{\psi_j}) B \bigr)$, and its variance by
$ \s_{j}^{2} = \tr \bigl( T(\ketbra{\psi_j}) B^2 \bigr) - 
{\tr}^2 \bigl( T(\ketbra{\psi_j}) B \bigr)$.
In terms of
$\theta_j$, this reduces to 
$b_j =\inp{\theta_j}{B \theta_j}$ 
and
$\s_{j}^{2} =\inp{\theta_j}{B^2 \theta_j} - \inp{\theta_j}{B \theta_j}^2$. 
Similarly, the l.h.s. of
(\ref{margel}) equals 
$| \overline{\alpha_0}\alpha_1 \inp{\theta_0}{A\theta_1} + 
\alpha_0\overline{\alpha_1}
\inp{\theta_1}{A\theta_0} |$, a quantity bounded by $| \inp{\theta_0}{A\theta_1} |$
since $2|\alpha_0|\cdot|\alpha_1| \le 1$. Formula (\ref{margel}) is then a direct application
of lemma~\ref{lem4}.

To reduce the general case to the case above, we note that
a non-pure state $\tau$ can always be represented as a
vector state. Explicitly, suppose that $\tau$ decomposes as 
$\tau = \sum_{i \in \NN} |\beta_i|^2 \ketbra{\varphi_i}$. Then define the Hilbert space 
$\tilde{\CK} := \bigoplus_{i \in \NN} \CK_i$, where each $\CK_i$ is a copy of $\CK$.
Now since 
$(\bigoplus_{i \in \NN} \CK_i) \otimes \CH \cong \bigoplus_{i \in \NN} (\CK_i
\ten \CH)$,  
we may define, for each $X \in B(\CK\ten\CH)$, the operator  
$\tilde{X} \in B(\tilde{\CK}\ten\CH)$ by 
diagonal action on the components of the sum, i.e. 
$\tilde{X} \bigr( \bigoplus_{i \in \NN}(k_i \ten h_i)\bigr) := 
\bigoplus_{i \in \NN} X(k_i \ten h_i)$. 
If we now define the vector $\tilde{\varphi} \in \tilde{\CK}$ by 
$\tilde{\varphi} = \bigoplus_i \beta_i \varphi_i$, then we have
for all $X \in \CK\ten\CH$ and $\chi \in \CH$:
\begin{eqnarray*}
\tr \bigl( \tilde{U} ( \ketbra{\tilde{\varphi}} \otimes \ketbra{\chi} ) \tilde{U}^*
\tilde{X}\bigr)\!\!\!
&=&
\!\!\!
\inp{{\textstyle \bigoplus}_{i \in \NN}(\beta_i \varphi_i \ten \chi)}
{\tilde{U}^*\tilde{X}\tilde{U}{\textstyle \bigoplus}_{j \in \NN}
(\beta_j \varphi_j \ten \chi)}_{\tilde{\CK}\ten \CH}\\
&=&
\!\!\!
\inp{{\textstyle \bigoplus}_{i \in \NN} (\beta_i\varphi_i \ten \chi)}
{{\textstyle \bigoplus}_{j \in \NN} U^*XU(\beta_j\varphi_j \ten
\chi)}_{\tilde{\CK}\ten \CH}\\
&=&
\!\!\!
\sum_{i \in \NN} |\beta_i|^2 \inp{(\varphi_i \ten \chi)}{U^*XU(\varphi_i \ten
\chi)}_{\CK\ten\CH}\\
&=&
\!\!\!
\sum_{i \in \NN} |\beta_i|^2 \tr \big( U (\ketbra{\varphi_i} \ten \ketbra{\chi})
U^* X \big)
\\
&=&
\!\!\!
\tr \bigl( U (\tau \otimes \ketbra{\chi}) U^*X\bigr)
\end{eqnarray*}
The second step is due to the diagonal action of the operators on
$\tilde{\CK}\ten\CH$.
The problem is now reduced to the vector-case by
applying the above to $\chi = \psi$, $\chi = \psi_0$ or
$\chi = \psi_1$ and on the other hand $X = A$, $X = B$ or $X = B^2$.
\qed
The backbone of theorem~\ref{cor5} is 
formed by the special case that $\s_0 = \s_1 =0$, 
$[A,B] = 0$ and $\tau = \ketbra{\phi}$, which allows for a short and
transparent proof. 

In order to arrive at a physical interpretation of theorem~\ref{cor5}, we 
focus on the case $B = \tilde{B} \ten \one$,
when information is transferred from $\CH$ to $\CK$. Indeed, examining $\CK$ at time
$t$ yields information about $\CH$ at time 0.

\subsubsection{Quality of Information Transfer}
A small ratio $\frac{\s_0 + \s_1}{| b_0 - b_1 |}$ indicates a good quality of
information transfer.
The ratio equals 0 in the perfect case, when $\sigma_0 = \sigma_1 =
0$. Thus $\tilde{B} \ten \one$ takes a definite value of either $b_0$
or $b_1$,  depending on whether the initial state of $\CH$ was $\ket{\psi_0}$
or $\ket{\psi_1}$. 
In this case, one can infer the initial state of $\CH$
with certainty by inspecting only the ancilla $\CK$. 
More generally, 
it is still possible to reliably determine from 
the ancilla $\CK$ whether the open system $\CH$ was initially in state
$\ket{\psi_0}$ or $\ket{\psi_1}$ 
as long as the standard deviations are
small compared to the difference in mean, 
$\sigma_0 , \sigma_1 \ll |b_0 - b_1|$. 
\begin{center}
\setlength{\unitlength}{0.9 cm}
\begin{picture}(8,4.5) 
\put(0,1){\line(1,0){8}}
\put(7.2,0.7){\mbox{\footnotesize $b$}}
\put(4,0.6){\line(0,1){3.9}}
\put(4.1,4){\mbox{\footnotesize $p(b)$}}
\qbezier(0,1.2)(0.7,1.3)(1,2.5)
\qbezier(1,2.5)(1.5,5)(2,2.5)
\qbezier(2,2.5)(2.3,1.3)(3,1.2)
\qbezier(3,1.2)(3.03 , 1.205)(8,1.05)
\qbezier(8 ,1.2)(7.3,1.3)(7,2.5)
\qbezier(7,2.5)(6.5,5)(6,2.5)
\qbezier(6,2.5)(5.7,1.3)(5,1.2)
\qbezier(5,1.2)(4.97 , 1.201)(0,1.05)
\put(1.5,2.5){\vector(1,0){0.5}}
\put(1.4,2.7){\mbox{\footnotesize $\sigma_0$}}
\put(1.5,2.5){\vector(-1,0){0.5}}
\put(6.5,2.5){\vector(1,0){0.5}}
\put(6.4,2.7){\mbox{\footnotesize $\sigma_1$}}
\put(6.5,2.5){\vector(-1,0){0.5}}
\put(4,1.5){\vector(1,0){2.5}}
\put(4.2,1.7){\mbox{\footnotesize $| b_1 - b_0 |$} }
\put(4,1.5){\vector(-1,0){2.5}}
\put(-0.5,0.2){\mbox{\small \fig Probability densities $p$ of $B$ for
input $\ket{\psi_0}$ and $\ket{\psi_1}$}}
\end{picture}
\end{center}
As the ratio increases, the restriction (\ref{margel})
gets less severe, reaching triviality at 
$\s_0 + \s_1 = 2|b_0 - b_1|$.

\subsubsection{Decoherence on the Commutant of the Pointer}
Assume perfect information transfer, i.e. $\s_0 = \s_1 =
0$. If $[A,B] = 0$, then theorem~\ref{cor5} 
says that 
coherent and mixed initial states yield identical distributions of $A$
at time $t$. 
In order to distinguish, at time $t$, whether or not $\CH$ was in a pure state 
at time 0, we will have to use observables $A$ which do not commute
with $B$. But then $A$ and $B$ cannot be observed simultaneously. Summarizing:
\begin{quote}\label{quote2}
\emph{At time $t$, it is possible to distinguish whether $\CH$ was in state
$\psi_0$ or $\psi_1$ at time 0. It is also possible to distinguish  whether $\CH$ was in state
$\psi$ or $\Col$ at time 0. But it is {\it not} possible to do both.
}
\end{quote}
We emphasize that this holds even when one has all observables of 
the entire closed system $\CK \ten \CH$ at one's disposal.

\subsubsection{Decoherence on the Open System}
We consider the final state of the open system $\CH$, obtained from the
final state of the closed system $\CK \ten \CH$ by tracing out the degrees of
freedom of the ancilla $\CK$: an initial state $\rho \in S(\CH)$
yields final state $\tr_{\CK}(T(\rho)) \in S(\CH)$.

Suppose that information is transferred to a pointer $B = \tilde{B} \ten \one$ in the
ancilla $\CK$ with perfect quality, that is $\s_0 = \s_1 = 0$. 
Since $[\one\ten \tilde{A}, \tilde{B}\ten\one] = 0$, we see 
from theorem~\ref{cor5} that
$\tr ( T(\ketbra{\psi}) (\one\ten \tilde{A}) ) $ 
must equal
$\tr ( T ( \Col ) ( \one\ten \tilde{A} ) )$
for {\it all} $\tilde{A} \in B(\CH)$,
or equivalently
\begin{equation}\label{schildpad}
   \tr_{\CK} \bigl( T(\ketbra{\psi}) \bigr) =
   \tr_{\CK} \Bigl( T \bigl( \Col  \bigr) \Bigr) \; .
\end{equation}
In words:
\begin{quote}\emph{Suppose that at time $t$, by making a hypothetical 
measurement of $\tilde{B}$ on the ancilla,
it would be possible to distinguish perfectly
whether the open system had been in state $\psi_0$ or $\psi_1$ at time 0.
Then, by looking only at the observables of the open system, 
it is not possible to distinguish whether $\CH$ 
had been in
the pure state $\psi=\alpha_0\psi_0+\alpha_1\psi_1$ or the collapsed state
$|\alpha_0|^2\ketbra{\psi_0}+|\alpha_1|^2\ketbra{\psi_1}$ at time 0.}
\end{quote}
This statement holds true, regardless whether $\tilde{B}$
is actually measured or not.
(So we do not assume here that such a measurement is physically possible.)
We have shown that the map
$M_0 + M_1 = \tr_{\CK} \circ T$,
with $T : \CT(\CH) \to \CT(\CK\ten\CH)$ the
information-transfer operation defined by
$ T(\rho) := U(\tau\ten\rho) U^*$,
constitutes a physical realization of  
the
abstract decoherence mapping $(M_0 + M_1)$ of section~\ref{perse}.

All in all, 
we have proven that decoherence is an unavoidable consequence of
information transfer out of an open system.

\subsubsection{Example}\label{spinflip}
The simplest possible example of 
unitary information transfer
is the following.
Let $\CK\sim\CH\sim{\mathbb{C}}^2$ be the Hilbert space of a qubit;
let $\psi_0=(1,0)$ and $\psi_1=(0,1)$ be the `computational basis',
and let $U:{\mathbb{C}}^2 \otimes {\mathbb{C}}^2\to{\mathbb{C}}^2 \otimes {\mathbb{C}}^2$ 
be the `controlled-not gate'.
Explicitly, $U$ is defined by 
$U\ket{\psi_1 \ten \psi_1} = \ket{\psi_0 \ten \psi_1}$,
$U\ket{\psi_0 \ten \psi_1} = \ket{\psi_1 \ten \psi_1}$,
$U\ket{\psi_1 \ten \psi_0} = \ket{\psi_1 \ten \psi_0}$, and
$U\ket{\psi_0 \ten \psi_0} = \ket{\psi_0 \ten \psi_0}$.
That is, it flips the first qubit whenever the second qubit is set to 1.
Let $\tau$ be the 0 state of the first qubit.

Since the initial state of the second qubit can be read off from the first, 
this situation satisfies the hypotheses of theorem~\ref{cor5} with 
$B = \s_z \ten \one$ and
$\s_0 = \s_1 = 0$.
We verify equation (\ref{schildpad}). For any state 
$\ket{\psi} = \alpha_0 \ket{\psi_0} + \alpha_1 \ket{\psi_1}$:
\begin{eqnarray*}
U \ket{\psi_0 \ten \psi} &=& \alpha_0 \ket{\psi_0 \ten \psi_0} 
+ \alpha_1 \ket{\psi_1 \ten \psi_1} \, := \, \ket{\theta} ;\\
\tr_{\CK} \left(\ketbra{\theta}\right) & = & \Col\,. 
\end{eqnarray*}
Thus we have
$\tr_{\CK} \bigl( T(\ketbra{\psi}) \bigr) = \Col$. 
This agrees with equation (\ref{schildpad}), as one can easily check that
$\tr_{\CK}\bigl( T(\Col)\bigr) $ equals $\Col$ as well.

\subsection{Unitary Information Transfer and State Collapse}\label{realcol}
We have derived that,
in the context of  
information transfer to an ancillary system, the initial states 
$\ketbra{\psi}$
and $|\alpha_0|^2 \ketbra{\psi_0} + |\alpha_1|^2
\ketbra{\psi_1}$ lead to the same final state.
This is decoherence.

State collapse is a much stronger statement:  
if outcome `0' is observed, 
then the system will further behave as if its initial state
had been $\psi_0$ instead of $\psi$.
Similarly, if outcome `1' is observed, then the system 
will behave as if its initial state had been $\psi_1$.
Now suppose that we ignore the outcome.
Since `0' happens with probability $|\alpha_{0}|^2$
and `1' with probability $|\alpha_{1}|^2$, 
the system will behave as if its initial state had been 
$|\alpha_0|^2 \ketbra{\psi_0} + |\alpha_1|^2
\ketbra{\psi_1}$. 
We see that collapse implies decoherence. 

The converse does not hold however: imagine a Stern-Gerlach experiment, 
in which a beam of particles in a $\s_x$-eigenstate is split in two 
according to spin in the $z$-direction. 
State collapse is the statement that 
one beam consists of particles with positive spin, the other 
of particles with negative spin and that both beams  
have equal intensity.
Decoherence is the statement that both outgoing beams \emph{together}
consist for 50\% of positive-spin particles and for 50\% of negative-spin
particles. The former statement is strictly stronger than the latter, 
and deserves separate investigation.   

We will therefore answer the following question: suppose that we
transfer information
to an ancilla $\CK$, and then separate $\CK$ from $\CH$,
dividing $\CH$ into subensembles according to outcome.
What states do we use to describe these subensembles?
 

\subsubsection{Joint Probability Distributions}
A special case of an observable is an \emph{event} $p$, which 
in quantum mechanics is represented by 
a projection $P$.
The relative frequency of occurrence of $p$ 
is given by $ {\mathbb P}(p=1) = \tr(\rho P)$. 

The projection $\one - P$ is interpreted as `not $p$'.
Furthermore, if a projection $Q$ corresponding to an observable $q$ 
commutes with $P$, then $PQ$ is again a projection.
According to quantum mechanics, $p$ and $q$ can then be observed
simultaneously, and the projection $PQ$ 
is interpreted as the event `$p$ and $q$ are both observed'.

A state $\rho$ therefore induces a joint probability distribution on
$p$ and $q$:
\begin{center}
\begin{tabular}{rcl}
$\tr(\rho PQ) = {\mathbb P}(p=1,q=1)\!\!\!$&
,&
$\!\!\!{\mathbb P}(p=0,q=1) = \tr(\rho (\one -P)Q) $ \\
$\tr(\rho P(\one -Q)) = {\mathbb P}(p=1,q=0)\!\!\!$ &
,&
$\!\!\!{\mathbb P}(p=0,q=0) = \tr(\rho (\one -P)(\one -Q))$
\end{tabular}
\end{center}  
Particularly relevant is the case in which $\rho$ is a state on 
a combined space $\CK \ten \CH$, and
the projections are of the form $Q \ten \one$ and $\one \ten P$.
(The commuting projections are properties of different systems.)
We then see that 
${\mathbb P}(p=1,q=1)$ 
equals 
$\tr( (\one \ten P) (Q \ten \one) \rho )$,
which in turn is equal to 
$\tr \left( P \, \tr_{\CK} ((Q\ten\one) \rho)  \right)$.
This holds for all projections $P$ on $\CH$, so that
the normalized version of 
$\tr_{\CK} ((Q\ten\one) \rho) \in \CT(\CH)$ 
must be interpreted as 
the state of $\CH$, given that $q = 1$.
Similarly, the normalized version of $\tr_{\CK} ((\one - Q)\ten\one \rho) \in \CT(\CH)$
is the state of $\CH$, given that $q=0$ is observed. 

\subsubsection{Collapse}
Let $T : \rho \mapsto U (\tau \ten \rho) U^*$ from $\CT(\CH)$ to $\CT(\CK \ten
\CH)$ be an information transfer from $\CH$ to a pointer-projection $Q \in
B(\CK)$. That is, $\tr((Q\ten\one) T(\ketbra{\psi_0})) = 0$ and 
$\tr((Q\ten\one) T(\ketbra{\psi_1})) = 1$, so that at time $t$, one can 
see from $\CK$ whether $\CH$ was in state $\psi_0$ or $\psi_1$ at time $0$.

Since $Q \ten \one$ 
commutes with all of $\one \ten B(\CH)$, it is possible to separate
$\CH$ from $\CK$, and divide $\CH$ into subensembles according to 
the outcome of $Q$.
This is done as follows: 
with any measurement 
on $\CH$, a simultaneous measurement of
$Q$ on $\CK$ is performed to determine in which ensemble 
this particular system should fall. 
It follows from the above that the $1$-ensemble should be
described by the normalized version of  
$M_1 (\rho) := \tr_{\CK} ((Q\ten \one) T(\rho))$, and the $0$-ensemble by
the normalized version of $M_0(\rho) := \tr_{\CK} ((\one-Q\ten \one) T(\rho))$.
Since $Q$ commutes with $B(\CH)$,
this is just conditioning on a classical probability space
at time $t$.
We have arrived at an interpretation of the map 
$M(\rho) := M_0(\rho) \oplus M_1(\rho)$ of section~\ref{perse}.

We will now prove that $\hspace{-0.4mm}M\hspace{-0.2mm}$ takes the form
$M(\ketbra{\psi}) = |\alpha_0|^2 \tr_{\CK}T(\ketbra{\psi_0}) \oplus 
|\alpha_1|^2 \tr_{\CK}T(\ketbra{\psi_1})$.
This is a strong physical statement. For instance,
any spin-system 
$\alpha_0 \ket{\psi_0} + \alpha_1 \ket{\psi_1}$ that is
found to have spin 1 in the $z$-direction 
may subsequently be treated as if it had been in state $\psi_1$ 
at time 0. This is nontrivial: a priori, it is perfectly conceivable
that the different initial states $\psi_0$ and $\psi$ result
in different final states, even though they yield the same $Q$-output.  

One could alternatively, (and more traditionally), arrive at the `collapse of the wavefunction'
$M(\ketbra{\psi}) = |\alpha_0|^2 \tr_{\CK}T(\ketbra{\psi_0}) \oplus 
|\alpha_1|^2 \tr_{\CK}T(\ketbra{\psi_1})$ by assuming that, 
at time 0, 
the quantum system 
makes either the jump $\ketbra{\psi} \mapsto \ketbra{\psi_0}$ 
or the jump $\ketbra{\psi} \mapsto \ketbra{\psi_1}$.
Since we arrive at the same conclusion, namely the above 
`collapse of the wavefunction', using only 
open systems,
unitary transformations and
the probabilistic
interpretation of quantum mechanics, such an  
assumption of `jumps' at time $0$ is made redundant.

\begin{proposition}\label{prop6}
Let $T : \rho \mapsto U (\tau\ten\rho) U^*$ from $\CT(\CH)$ to $\CT(\CK \ten
\CH)$ satisfy $\tr((Q\ten\one) T(\ketbra{\psi_0})) = 0$ and
$\tr((Q\ten\one) T(\ketbra{\psi_1})) = 1$ for some `pointer-projection' $Q$
on $\CK$. 
Define a map $M : \CT(\CH) \mapsto \CT(\CH)\oplus \CT(\CH)$ by 
$M(\rho) := \tr_{\CK} ((\one -Q\ten\one) T(\rho)) \oplus \tr_{\CK} ((Q\ten\one)
T(\rho))$.
Then for $\psi = \alpha_0\psi_0 + \alpha_1 \psi_1$ we have
$M(\ketbra{\psi}) = |\alpha_0|^2 \tr_{\CK}T(\ketbra{\psi_0})
\oplus |\alpha_1|^2 \tr_{\CK}T(\ketbra{\psi_1})$.
\end{proposition}
This can be seen almost directly from proposition~\ref{abscol}:

\proof 
Since $M_{1}(\ketbra{\psi_0}) \geq 0$ is a positive operator, we may
conclude from $\tr \left(M_{1}(\ketbra{\psi_0})\right) = 0$ that
$M_{1}(\ketbra{\psi_0})$ itself must vanish. 
Similarly, reversing the role of $0$ and $1$, we see that 
$M_{0}(\ketbra{\psi_1}) = 0$.
From proposition~\ref{abscol}, 
we then obtain
$M(\ketbra{\psi}) = |\alpha_0|^2 M_{0}(\ketbra{\psi_0}) 
\oplus |\alpha_1|^2 M_{1}(\ketbra{\psi_1})$.
The proof is completed by noting from 
$\tr_{\CK}((\one - Q\ten \one) T(\ketbra{\psi_1}))=0$ 
that 
$\tr_{\CK}(T(\ketbra{\psi_1})) =   
\tr_{\CK}((Q\ten \one) T(\ketbra{\psi_1})) =
M_{1}(\ketbra{\psi_1})$,
and similarly, again reversing the role of $0$ and $1$, that
$\tr_{\CK}(T(\ketbra{\psi_0})) =   M_{0}(\ketbra{\psi_0})$.
\qed
\newpage
We summarize: 
\begin{quote}
\emph{Consider an ensemble of systems of type $\CH$ in state $\psi$. 
	Suppose that information is transferred   
	to a pointer-projection 
	$Q$ on an ancillary system $\CK$. 
	Subsequently, the ensemble is divided into two subensembles according
	to outcome. 
	Then all observations on $\CH$ made afterwards, conditioned on the 
	observation that
	the measurement outcome was 0, will be as if the system had 
	originally been in 
	the collapsed state $\psi_0$ instead of $\psi$.
	No independent `collapse postulate' is needed to arrive at this
	conclusion.}
	
\end{quote}

\subsubsection{Example} 
In the simple model of information transfer introduced in 
Section~\ref{ITSC}, we will now demonstrate why 
repeated spin-measurements 
yield identical outcomes. 

The probed system is once again a single spin $\CH = {\mathbb{C}}^2$,
whereas the ancillary system now consists of two spins, 
$\CK = {\mathbb{C}}^2 \otimes {\mathbb{C}}^2$ in initial state $\ket{\psi_0 \ten \psi_0}$. 
Repeated information-transfer, first to pointer 
$\s_{z,1}$ and then to $\s_{z,2}$, is then represented by the unitary 
$U := U_2 U_1$ on $\CK \ten \CH$. In this expression, $U_1$ is the 
controlled not-gate flipping the first qubit of $\CK$ if $\CH$ is 
set to 1,
and $U_2$ flips the second qubit of $\CK$ if $\CH$ is set to 1. 

Since $U  \ket{\psi_0 \ten \psi_0 \ten (\alpha_0\psi_0 + \alpha_1\psi_1)} = 
\ket{\alpha_0 \psi_0 \ten \psi_0 \ten \psi_0}
+ \ket{\alpha_1 \psi_1 \ten \psi_1 \ten \psi_1}$,
we can explicitly calculate the joint probability distribution 
on the two pointers $\s_{z,1}$ and $\s_{z,2}$ in the final state: 
\begin{center}
\begin{tabular}{rcl}
$ {\mathbb P}(s_{z,1}= \,\,\,\,\,1,s_{z,2}=1) = |\alpha_1|^2$ &,&  
$\hspace{2mm}0\hspace{3,6mm} = {\mathbb P}(s_{z,1}=\,\,\,\,\,1,s_{z,2}=-1)$ \\
$ {\mathbb P}(s_{z,1}=-1,s_{z,2}=1)=\hspace{2mm}0\hspace{3,6mm}$ &,& 
$|\alpha_0|^2={\mathbb P}(s_{z,1}=-1,s_{z,2}=-1)$
\end{tabular}
\end{center}
In particular, we see that if the first outcome is $1$ (which happens with
probability $|\alpha_1|^2$), then so is the second.
proposition~\ref{prop6} shows that this is the general situation, 
independent of the (rather simplistic) details of this particular model. 

\subsection {Information Leakage to the Environment}\label{leakage}  
On closed systems decoherence does not occur, because unitary time
evolution preserves the purity of states.  
However, macroscopic systems are almost never closed.

Imagine, for example, that $\CH = { \mathbb C}^2$ represents a two-level atom, 
and $\CK$ some large measuring device.
Information about the energy $ \one \ten \sigma_z$ of the atom is
transferred to the apparatus, where it is stored as the position 
$\tilde{B} \ten \one$ of a pointer. 
Then as soon as information on the pointer-position $ \tilde{B} \ten \one$
leaves the system, collapse on the combined atom-apparatus system takes place. 
For example, a ray of light may reflect on the pointer, revealing its
position to the outside world. (See \cite{JZ}.)
It is of course immaterial whether or not someone is actually {\it looking} at the photons. 
If even the smallest speck of light were to fall on the pointer, the information
about the pointer position would 
already be encoded in the light, causing full collapse on the atom-apparatus system.
(See \cite{Zu} for an example.)

The quality of this information transfer will not be perfect.
If a macroscopic system is interacting normally with the outside world, (the occasional photon happens to 
scatter on it, for instance),
then a number of macroscopic observables $X$ will leak
information continually, with a macroscopic uncertainty $\sigma$. 
This enables us to apply theorem~\ref{cor5}. 
It says that all coherences between eigenstates $\smash{\psi_{x_1}}$ and $\smash{\psi_{x_2}}$ of 
macroscopic observables 
$X$ are continually vanishing on the macroscopic system $\smash{\CL}$,  
provided that their eigenvalues $x_1$ and $x_2$ satisfy the inequality 
$\smash{| x_1 - x_2 | \gg 2\sigma}$. 
(The pointer, e.g. a
beam of light, is outside the system, so that $\delta = 0$.)

Take for example a collection of $N$ spins, with Hilbert space 
$\smash{\CL = \bigotimes_{i = 1}^{N}}
{\mathbb C}^2$.
Suppose that the average spin-observables 
$\smash[t]{S_{\alpha} = \frac{1}{N}\sum_{i = 1}^{N} \sigma_{\alpha}^{i}}$ 
(for $\alpha = x,y,z$)
are continually being measured
with an accuracy\footnote{Since $[S_{x}, S_{y}] \neq 0$, they cannot be
simultaneously measured with complete accuracy, see e.g.\cite{We}.
However, this problem disappears if the accuracy satisfies 
$\s^2 \geq \frac{1}{2} \|[S_{x} , S_{y}]\| = \frac{1}{N}$, see \cite{scriptie}.
For large $N$, (typically $N \sim 6 \times 10^{23}$),
this allows for extremely accurate measurement. 
} $N^{-1/2} \ll \sigma \ll 1$. Then between 
macroscopically different eigenstates of $S_{\alpha}$, i.e. states for which the eigenvalues
satisfy 
$|s_{\alpha} - \smash{s'_{\alpha}}| \gg \sigma$, coherences 
are constantly disappearing. 
However, the information leakage need not have any effect 
on states which only differ on a microscopic scale.
Take for instance $\rho \otimes \ketbra{+} $ 
and $\rho \otimes \ketbra{-} $, with $\rho$ an arbitrary state on $N-1$ spins.
Indeed, $|s_{\alpha} - s'_{\alpha}| \leq 2/N \ll \sigma$, so 
theorem~\ref{cor5} is 
vacuous in this case: no decoherence occurs. 

We see how the variance $\sigma^2$ produces a smooth boundary between the macroscopic 
and the microscopic world: macroscopically distinguishable states
(involving $S_{\alpha}$-differences $\gg \sigma$) continually 
suffer from loss of coherence, while states that only differ microscopically 
(involving $S_{\alpha}$-differences $\ll \sigma$) 
are unaffected. 

In case of a system monitored by a macroscopic measurement apparatus, we are 
interested in coherence between eigenstates of the macroscopic
pointer. By definition, these eigenstates are macroscopically distinguishable. 
We may then give the following answer 
to the question why
it is so hard, in practice, to witness coherence: 
\begin{quote}
\emph{If information leaks from the pointer into the outside world,
decoherence 
takes place on the 
combination of system and measurement apparatus. In practice, macroscopic
pointers constantly leak information.}
\end{quote}\vspace{-3mm}

\section{Closed Systems}\label{oester}\vspace{-1mm}
Closed systems evolve according to unitary time evolution, so that
coherence which is present initially will still be there at later times. 
Yet on macroscopic systems, coherent superpositions are almost never observed. 
Why is this so?\vspace{-1mm} 

\subsection {Macroscopic Systems} \label{macro}\vspace{-1mm}
Because of the direct link that it provides between
the scale of a system on the one hand, and on the other hand
the difficulties in witnessing coherence, we feel that 
the following line of reasoning, essentially due to Hepp \cite{Hp},
is the most important mechanism hiding coherence.

Let us first define what we mean by macroscopic and microscopic
observables. We consider a system consisting of $N$ distinct
subsystems, i.e. $\CK = \bigotimes_{i=1}^{N} \CK_i$.
If one thinks of $\CK_i$ as the atoms out of which a macroscopic system
$\CK$ is constructed, $N$ may well be in the order of $10^{23}$.

We will define the \emph{microscopic} observables to be the ones that refer
only to one particular subsystem $\CK_i$:\vspace{-1mm}  
\begin{basdef}
An observable $X \in B(\CK)$ is called \emph{microscopic} if it is of the
form $X = \one \otimes \ldots \otimes \one \otimes X_i \otimes \one 
\otimes \ldots \otimes \one$ for some $i \in \{ 1,2,\ldots,N \}$
and some $X_i \in B(\CK_i)$.
\end{basdef}\vspace{-1mm} 
In this situation we will identify $X_i \in B(\CK_i)$ with $X \in B(\CK)$. 
We take macroscopic
observables to be averages of microscopic observables `of the same size':\vspace{-5mm}
\begin{basdef}
An observable $Y \in B(\CK)$ is called \emph{macroscopic} if it is of the
form $Y = \frac{1}{N} \sum_{i=1}^{N} Y_i$, 
with $Y_i \in B(\CK_i)$ such that $\|Y_i\| \leq \|Y\|$. 
\end{basdef}\vspace{-1mm}
We will only use the term `macroscopic' in this narrow
sense from here on, even though there do exist observables which are 
called `macroscopic' in daily life, but do not fall under the above
definition.

Now suppose that we transfer information from a system $\CH$ to a 
macroscopic system $\CK = \bigotimes_{i=1}^{N} \CK_i$, using a 
macroscopic pointer $\tilde{B} \in B(\CK)$. As explained before, we then have a map
$T : \CT(\CH) \rightarrow \CT(\CK \otimes \CH)$ such that the pointer
$\tilde{B} \otimes \one$ has 
different expectation values $b_0$ and $b_1$ in
the states
$T(\ketbra{\psi_0})$ and $T(\ketbra{\psi_1})$.

Since the pointer $\tilde{B}$ is macroscopic, it is unrealistic to require 
$T(\ketbra{\psi_0})$ and $T(\ketbra{\psi_1})$
to be eigenstates of $\tilde{B}$. Instead, we will require their standard 
deviations in $\tilde{B}$ to be negligible compared to their difference in mean,
that is $\s_0 \ll |b_0 - b_1|$ and $\s_1 \ll |b_0 - b_1|$.  

After this information transfer, we try to distinguish whether the 
system $\CH$ had initially been in the coherent state 
$\alpha_0 \ket{\psi_0} + \alpha_1 \ket{\psi_1}$
or in the incoherent mixture $|\alpha_0|^2 \ketbra{\psi_0} + |\alpha_1|^2 \ketbra{\psi_1}$. 
We have already shown that this cannot be done by measuring 
observables in $\one \otimes B(\CH)$. 
The following adaptation of theorem~\ref{cor5} shows that it 
is also impossible to do this
by measuring
macroscopic or microscopic observables on the closed system $\CK \otimes \CH$.

\begin{corollary}\label{cor6}
Let $\psi_0$ and $\psi_1$ be orthogonal unit vectors in a Hilbert space $\CH$
and let $\tau \in \CT(\CK)$ be a state on the Hilbert space 
$\CK = \bigotimes_{i=1}^{N} \CK_i$.
Let $U : \CK \ten \CH \rightarrow \CK \ten \CH$ be unitary and define
$T : \CT(\CH) \rightarrow \CT(\CK \ten \CH)$ by
$T(\rho) = U (\tau \ten \rho) U^* $.
Let $\tilde{B}$ be a macroscopic observable in 
$B(\CK)$, and define $B := \tilde{B} \ten \one$. 
Denote by $b_j$ 
and $\s_{j}^{2}$  its expected value  and variance 
in the state $T(\ketbra{\psi_j})$ for $j = 0,1$.
Suppose that $b_0 \neq b_1$. 
Then for all $\psi=\alpha_0\psi_0+\alpha_1\psi_1$ with 
$|\alpha_{0}|^2 + |\alpha_{1}|^2 = 1$
and for all microscopic and macroscopic
observables $A \in B(\CK \otimes \CH)$, we have
$$ \label{mac}
   \Bigl| \tr \bigl( T(\ketbra{\psi}) A \bigr) -
   \tr \Bigl( T\bigl( |\alpha_0|^2 \ketbra{\psi_0} + 
   |\alpha_1^2| \ketbra{\psi_1}\bigr) A \Bigr) \Bigr| 
   \le
   \frac{\frac{2}{N} \norm B\!+\!\s_0\!+\!\s_1}{|b_0-b_1|}\norm A \,.
$$
\end{corollary}
\proof 
If A is microscopic, we have the inequality
$$ \big\|[A,B]\big\| = \Big\|\Big[A_i , {\textstyle \frac{1}{N}}
\sum_{j=1}^{N} B_j\Big]\Big\| = {\textstyle \frac{1}{N}} \big\|\big [A_i, B_i\big]\big\| 
\leq \frac{2\|A\|
\|B\|}{N}\,.$$
If A is macroscopic, we have 
$$ \big\|[A,B]\big\| = \Big\|\Big[ {\textstyle  \frac{1}{N+1}} 
\sum_{i=0}^{N} A_i , {\textstyle \frac{1}{N}} \sum_{j=1}^{N} B_j\Big]\Big\| = 
{\textstyle \frac{1}{N(N+1)}} \sum_{i=1}^{N} \big\|[A_i, B_i]\big\| \leq \frac{2\|A\|
\|B\|}{N}\,.$$
Either way, we can now apply theorem~\ref{cor5}. \qed

\subsection{Examples}
In order to illustrate the above,  
we discuss four examples of information transfer to
a macroscopic system.

\subsubsection{The Finite Spin-Chain}\label{spin-chain}
We study a single spin $\CH = {\mathbb{C}}^2$ in interaction with 
a large but finite spin-chain $\CK = \bigotimes_{i=1}^{N} {\mathbb{C}}^2$,
the latter acting as a measurement apparatus.
Once again, let $\psi_0=(1,0)$ and $\psi_1=(0,1)$ be the `computational
basis'. Initially, all spins in the spin-chain are down: 
$\tau = \ketbra{\psi_0 \otimes\ldots \otimes \psi_0}$.
Let $U_i:\CK\ten\CH\to\CK\ten\CH$ be the `controlled-not gate',
which flips spin number $i$ in the chain whenever the single qubit is set to 1.
(We define $U_j = \one$ for $j \notin \{1,2,\ldots,N\}$.)
$$U_i = \one \otimes P_{-} + \s_{x,i} \otimes P_{+} \quad \mbox{with} \quad
	P_+ =  \begin{pmatrix}	1&0 \\
				0&0 \\
		\end{pmatrix}, \quad
	P_- =  \begin{pmatrix} 	0&0 \\
				0&1 \\
		\end{pmatrix}	\,.
$$
In discrete time $n \in {\mathbb Z}$, the unitary evolution is given 
by $n \mapsto U_n U_{n-1} \ldots U_2 U_1$. (See \cite{Hp}.)
This represents a single spin flying over a spin-chain from 1 to $N$,
interacting with spin $n$ at time $n$. 

Obviously, we have
$U_N  \ket{\psi_0 \otimes\ldots\otimes\psi_0} \otimes \ket{\psi_0}=
     \ket{\psi_0 \otimes\ldots\otimes\psi_0} \otimes \ket{\psi_0} $
and     
$U_N \ket{\psi_0 \otimes\ldots\otimes\psi_0} \otimes \ket{\psi_1}=
    \ket{\psi_1 \otimes\ldots\otimes\psi_1} \otimes \ket{\psi_1}$.
We consider the average spin of the spin-chain as pointer, 
$B = \frac{1}{N} \sum_{i=1}^{N} \s_{z,i}$.
This makes the map $T : \rho \mapsto U_N \tau \otimes \rho U_N^{*}$ an
information transfer to a macroscopic system.
Applying Corollary~\ref{cor6} with $b_0 = -1$, $b_1 = 1$ and $\s_0 = \s_1 = 0$,
we see that
$$
   \Bigl| \tr \bigl( T(\ketbra{\alpha_0\psi_0+\alpha_1\psi_1}) A \bigr) -
   \tr \Bigl( T\bigl( |\alpha_0|^2 \ketbra{\psi_0} + 
   |\alpha_1^2| \ketbra{\psi_1}\bigr) A \Bigr) \Bigr| 
$$
is bounded by $\frac{1}{N} \norm A$
for all microscopic and macroscopic $A \in B(\CK \ten \CH)$. 
It is not hard to see that in this particular model, the estimated quantity is even identically
zero. Indeed, $\inp{\psi_0 \otimes \ldots \psi_0}
{X_i \,\psi_1 \otimes\ldots\otimes\psi_1 } = 
\inp{\psi_0}{\psi_1}^{N-1}\inp{\psi_0}{X_i \psi_1} = 0$ for all 
microscopic $X_i$.

Of course coherence \emph{can} be detected on the closed system $\CK \ten \CH$,
but only using observables that are neither macroscopic nor microscopic, 
such as the `product of all spins' $\s_x \ten \ldots \ten \s_x$.
\subsubsection{Finite Spin-Chain at Nonzero Temperature}
A more realistic initial state for the spin-chain is the
nonzero-temperature state $\tau_{\beta} = \frac{e^{-\beta H}}{\tr e^{-\beta H}}$.
For the spin-chain Hamiltonian we will take 
$H = \sum_{i} \s_{z,i} = NB$, so that $\tau_{\beta}$ 
becomes the tensor product of
$N$ copies of the ${\mathbb{C}}^2$-state 
$$\hat{\tau}_{\beta} = \frac{1}{e^{\beta} + e^{-\beta}}  
		\begin{pmatrix}	e^{-\beta}	&0 \\
				0		&e^{\beta} \\
		\end{pmatrix}\,.
$$
With the same time evo\-lu\-tion as before, we have
$T(\ketbra{\psi_0}) = \ketbra{\psi_0} \otimes \rho_{\beta}$ and 
$T(\ketbra{\psi_1}) = \ketbra{\psi_1} \otimes \rho_{-\beta}$. 
Again we choose the mean energy $B$ as our pointer.
\mbox{A brief} calculation shows that 
$\tr(B\tau_\beta) = \frac{e^{-\beta} - e^{\beta}}{e^{\beta} +
e^{-\beta}} =: \varepsilon(\beta)$ and that
$\tr(B^2 \tau_\beta) - \tr(B \rho_\beta)^2 = \frac{1}{N}(1 - \varepsilon^2(\beta)) $.
Corollary~\ref{cor6} now gives us, for microscopic and macroscopic $A$,
$$
\Bigl| \tr \bigl( T(\ketbra{\alpha_0\psi_0+\alpha_1\psi_1}) A \bigr) -
   \tr \Bigl( T\bigl( |\alpha_0|^2 \ketbra{\psi_0} + 
   |\alpha_1^2| \ketbra{\psi_1}\bigr) A \Bigr) \Bigr|
$$\vspace{-4mm}
$$
\le
   \left( \frac{1}{\varepsilon(\beta) N} + \frac{\sqrt{1 - \varepsilon^2(\beta)}}
   {\varepsilon(\beta) \sqrt{N}}\right) 
   \norm A \,. 
$$
For large $N$, we see that the term $\sim \frac{1}{N}$ due to
the fact that $[A,B] \neq 0$ is dominated by the thermodynamical
fluctuations, which of course go as $\sim \frac{1}{\sqrt{N}}$. 
In statistical physics, it is standard practice to neglect even the latter.

\subsubsection{Energy as a Pointer}
Hamiltonians often fail to be macroscopic in our narrow sense of the word,
since they are generically unbounded and contain interaction terms.
However, this does not imply failure of our scheme to estimate coherence.

    
For example, consider an $N$-particle system with Hilbert space  $\CK =
\bigotimes_{i=1}^{N} \CK_{i}$ 
and Hamiltonian $H = \sum_{i=1}^{N} \frac{p_{i}^{2}}{2m_{i}} + V(x_1, x_2,
\ldots, x_{N} )$.   
Information is transferred from $\CH$ to $\CK$ with $H$ as pointer, that is 
the two states $\tr_{\CH}(T(\ketbra{\psi_0}))$ and
$\tr_{\CH}(T(\ketbra{\psi_1}))$ have
different energies $E$ and $E'$.
Without loss of generality, assume that they are vector states:
$\tr_{\CH}(T(\ketbra{\psi_0})) = \ketbra{\psi}$ and
$\tr_{\CH}(T(\ketbra{\psi_1})) =
\ketbra{\psi'}$. (Density matrices can always be represented as vectors
on a different Hilbert space, cf. the proof of theorem~\ref{cor5}.)

We thus have two vector states $\ket{\psi}$ and $\ket{\psi'}$ with two different
energies, $E := \inp{\psi}{H\psi}$ and $E' := \inp{\psi'}{H\psi'}$. 
We estimate the coherence between $\ket{\psi}$ and $\ket{\psi'}$ on $x_{n}$, 
the position of particle $n$.
\begin{eqnarray*}
(E \!-\! E')\inp{\psi}{x_{n} \psi'} \hspace{-1mm}&=&\hspace{-1mm} \inp{E \psi}{x_{n} \psi'} - \inp{x_{n}\psi}
{E' \psi'}\\
&=&\hspace{-1mm} \inp{H \psi - (H - E) \psi}{x_{n} \psi'} - \inp{x_{n} \psi}{ H \psi' - (H
- E')\psi'}\\
&=&\hspace{-1mm} \inp{[H,x_{n}] \psi}{ \psi'} - \inp{(H - E) \psi}{x_{n} \psi'} + 
\inp{x_{n} \psi}{(H - E')\psi'}
\end{eqnarray*}
Now since 
$
[H,x_{n}] = \frac{1}{2m_{n}} [p_{n}^{2}, x_{n}] = \frac{-i\hbar
p_{n}}{m_{n}} 
$, we can apply the Cauchy-Schwarz inequality in each term to obtain

\begin{eqnarray*}
|E - E'| \, |\inp{\psi}{x_{n} \psi'}| &\leq& 
\frac{\hbar}{m_{n}} \sqrt{\inp{\psi}{p_{n}^{2} \psi}} +
\sqrt{\inp{\psi}{x_{n}^2 \psi}} \sqrt{\inp{\psi'}{(H - E')^2 \psi'}} \\[-1mm]
&+&
\sqrt{\inp{\psi'}{x_{n}^2 \psi'}} \sqrt{\inp{\psi}{(H - E)^2 \psi}} \;.
\end{eqnarray*}
If we define the characteristic speed
$V_{n} := \sqrt{ \inp{\psi}{({p_{n}}/{m_{n}})^2 \psi} }$,
the characteristic positions $ X_{n} := \sqrt{\inp{\psi}{x_{n}^2 \psi}}$  
and $ X'_{n} := \sqrt{\inp{\psi'}{x_{n}^2 \psi'}}$, and the 
standard deviations $\sigma := \sqrt{\inp{\psi}{(H - E)^2 \psi}}$ and
$\sigma' := \sqrt{\inp{\psi'}{(H - E')^2 \psi'}}$, we obtain
$$
|\inp{\psi}{x_{n} \psi'}| \leq \frac{\hbar V_{n} + \sigma X'_{n} +
\sigma' X_{n}}{|E - E'|}\;.
$$ 
As such, this doesn't tell us very much. We will have to make some
physically plausible assumptions on the state of the system
in order to obtain results.
First, we assume that the
system is encased in an $L \times L \times L$ box so that $X_{n},
X'_{n} \leq L$. Also, we assume $V_{n} < c$. This yields   
$
|\inp{\psi}{x_{n} \psi'}| \leq \frac{\hbar c + L( \sigma + \sigma' )}{|E - E'|}
$. 
Secondly, we assume that scaling the system in any meaningful way will 
produce $|E - E'| \sim N$
and $\sigma + \sigma' \sim  \sqrt{N}$, so that the coherence on $x_n$
approaches zero as $\sim\frac{1}{\sqrt{N}}$.
Notice the almost thermodynamic lack of detail required for this estimate.

\subsubsection{Schr\"odinger's Cat}
Let us finally analyze the rather drastic extraction of information from a
radioactive particle that has become known\footnote{
Actually, Schr\"odinger's proposal was slightly different.
In the original thought experiment, death of the cat was correlated 
with decay of the atom at time $t$ instead of $0$, which wouldn't
make it an information transfer in our sense of the word. 
} 
as `Schr\"odinger's cat'. (See \cite{Sch}.) 
The experiment is performed as follows. We are interested in a radioactive
particle. Is it in a decayed state $\psi_0$ or in a non-decayed state
$\psi_1$? 

In order to determine this, we set up the following experiment. 
A Geiger counter is placed next to the radioactive particle. If the particle
decays, then the Geiger counter clicks. A mechanism then releases a hammer,
which smashes a vial of hydrocyanic acid, killing a cat.
All of this happens in a closed box no higher than $1 m$, and
completely impenetrable to
information.  
A measurement of the atom is done as follows: 
first, place it inside the box. Then wait for
a period of time that is long compared to the decay time of the atom.
Finally, open the box, 
and inspect whether the cat has dropped dead or is still standing upright. 

The atom is described by a Hilbert space $\CH$, the combination of
Geiger counter, mechanism, hammer, vial and cat by a Hilbert space $\CK$.
Initially, the latter is prepared in a state $\ket{\theta}$.
As a pointer, we take the center of mass of the cat, 
$Z := \frac{1}{N} \sum_{i=1}^{N} z_i$. In this expression,
$N$ is the amount of atoms out of which the cat is constructed,
and $z_i$ is the $z$-component of particle number $i$.
(It is a harmless assumption that all atoms in the cat have the same mass.)  
Since the 
box only measures $1 m$ in height, we may take $\|Z\| = 1$.
The unitary evolution $U \in B(\CK \otimes \CH)$ then  produces
$U \ket{\psi_0 \otimes \theta} := \ket{\gamma_0}$ and
$U \ket{\psi_1 \otimes \theta} := \ket{\gamma_1}$, which are 
eigenstates\footnote{
As discussed before, it would be more realistic to allow for a nonzero 
variance $0 < \s_j \ll 1$ instead of requiring $\theta_j$ to be
eigenstates of $Z$. We use $\s_j = 0$ for clarity,
leaving the argument essentially unchanged .   
} of $Z$ with different eigenvalues.

Suppose that, initially, the atom is either in the decayed
state $\psi_0$ with probability $|\alpha_0|^2$ or in the
non-decayed state $\psi_1$ with probability $|\alpha_1|^2$.
That is, the initial state is the incoherent mixture 
$|\alpha_0|^2 \ketbra{\psi_0} + |\alpha_1|^2 \ketbra{\psi_1}$.
By linearity, the final state is then the incoherent
state $|\alpha_0|^2 \ketbra{\gamma_0} + |\alpha_1|^2 \ketbra{\gamma_1}$.

On the other hand,
if the atom starts out in the coherent superposition 
$\alpha_0 \ket{\psi_0} + \alpha_1 \ket{\psi_1}$, then 
the combined system ends up in the coherent state
$U \ket{(\alpha_0 \psi_0 + \alpha_1 \psi_1) 
\otimes \theta} = \alpha_0 \ket{\gamma_0} + \alpha_1 \ket{\gamma_1}$.

The question is now this: why do we not notice the difference 
between these two situations if we open the box?
First of all, according to theorem~\ref{cor5} (and the observations following 
it in section~\ref{quote2}), 
it is impossible to detect coherence between $\gamma_0$
and $\gamma_1$ \emph{and} ascertain the position of the cat.
Upon opening the black box, we must make a choice. 

Secondly, according to the discussion in section~\ref{leakage}, the 
coherences between the macroscopically different states $\gamma_0$
and $\gamma_1$ are extremely volatile. Any speck of light falling 
on the cat will reveal its position with reasonable accuracy, 
causing the coherence 
to disappear according to theorem~\ref{cor5}.

Yet even if we were able to open the box without any information on the
position of the cat leaking out, even then would we be unable to detect
coherence between $\gamma_0$ and $\gamma_1$.
Apply Corollary~\ref{cor6} to the transfer of information from atom to cat.
We have $\s_0 = \s_1 = 0$, and with pointer $Z$ we have $\|Z\| = 1$
(the height of the box is 1 m) and $z_1 - z_0 = 0.1$ (the difference 
between a cat that is standing up and one that has dropped dead is 10
$cm$). We then obtain for all macroscopic and microscopic $A$:
$$
	\left| \inp{\alpha_0 \gamma_0 + \alpha_1 \gamma_1}
	{A \alpha_0 \gamma_0 + \alpha_1 \gamma_1} -
	\left( |\alpha_0|^2 \inp{\gamma_0}{A \gamma_0} +
	|\alpha_1|^2 \inp{\gamma_1}{A \gamma_1} \right) \right|
	\leq {\textstyle\frac{20}{N}} \|A\|\;.
$$
On the subset of observables we are normally able to measure,
the distinction between coherent and incoherent mixtures 
practically vanishes for $N \sim 10^{23}$. 
For all practical intents and purposes, it is completely 
harmless to assume that the 
final state of the cat is 
$|\alpha_0|^2 \ketbra{\gamma_0} + |\alpha_1|^2 \ketbra{\gamma_1}$
instead of
$\alpha_0 \ket{\gamma_0} + \alpha_1 \ket{\gamma_1}$.
But it would be false to state that the former has actually been observed.

\section{Conclusion}
In open systems, we have proven that decoherence is a necessary consequence
of information transfer to the outside. More in detail, 
we have reached the following conclusions:
\begin{itemize}
\item[-]Suppose that an open system $\CH$ interacts with an ancillary system
$\CK$ in such a way, that it is possible, in principle, to determine from
$\CK$ whether $\CH$ had been in state $\psi_0$ or $\psi_1$ before the
interaction. If $\CH$ started out in a coherent state 
$\alpha_0 \ket{\psi_0} + \alpha_1 \ket{\psi_1}$, then it will behave after
the information transfer as if it had started out in the incoherent mixture
$|\alpha_0|^2 \ketbra{\psi_0} + |\alpha_1|^2 \ketbra{\psi_1}$ instead.
This is called `decoherence'.
\item[-]Suppose again that the information whether $\CH$ was in state
$\psi_0$ or $\psi_1$ is transported to an ancillary system $\CK$.
This is done with an ensemble of $\CH$-systems described by the state 
$\alpha_0 \ket{\psi_0} + \alpha_1 \ket{\psi_1}$.
The ensemble is then split into subensembles, according to outcome.
The `0-ensemble' then behaves as if it had been in state $\psi_0$
at the beginning of the procedure,
and the `1-ensemble' as if it had started in state $\psi_1$.
This is called `state collapse'.
\item[-]These results were obtained entirely within the framework of 
traditional quantum mechanics and unitary time evolution on a
larger, closed system containing $\CH$. No `reduction-postulate' 
is needed. From proposition~\ref{abscol}, we see that any information
extraction causes collapse, quite independent of its particular mechanism.  
\item[-]On the closed system containing the smaller, open one no 
decoherence occurs in principle. 
In practice however, closed systems are very hard to achieve. 
We have argued that 
information transfer from a macroscopic observable $A$, performed with 
macroscopic precision $\s$, causes decoherence between eigenstates 
of $A$ if their values satisfy $\s \ll |a_1 - a_0|$. 
Since information on macroscopic observables tends to leak out, 
coherence between macroscopically different states tends to vanish.   
\end{itemize}
Still, even if the combined system $\CK \ten \CH$ is considered perfectly
closed, there are some results to be obtained.
Again, we investigated the case that a system $\CH$ interacts 
unitarily with a system $\CK$ 
in such a way that the information whether $\CH$ was in state $\psi_0$
or $\psi_1$ can be read off from a pointer in $\CK$.
We have reached the following conclusions concerning the closed system $\CK
\ten \CH$:
\begin{itemize}
\item[-]Using\hspace{-0.3mm} only\hspace{-0.3mm} observables\hspace{-0.3mm} 
on\hspace{-0.3mm} the\hspace{-0.3mm} closed\hspace{-0.3mm} system\hspace{-0.3mm} 
that\hspace{-0.3mm} commute\hspace{-0.3mm} with\hspace{-0.3mm} the\hspace{-0.3mm} 
pointer, 
it is impossible to detect whether $\CH$ had started out in state 
$\alpha_0 \ket{\psi_0} + \alpha_1 \ket{\psi_1}$ or
$|\alpha_0|^2 \ketbra{\psi_0} + |\alpha_1|^2 \ketbra{\psi_1}$.
Physically, this means that it is impossible to distinguish between 
coherent and incoherent initial states while at the same time 
distinguishing between $\psi_0$ and $\psi_1$.  
\item[-]Suppose that the closed system $\CK \ten \CH$ is macroscopic,
and that one has access to its macroscopic and microscopic
observables only. Then it is almost impossible to distinguish 
whether $\CH$ had started out in state 
$\alpha_0 \ket{\psi_0} + \alpha_1 \ket{\psi_1}$ or
$|\alpha_0|^2 \ketbra{\psi_0} + |\alpha_1|^2 \ketbra{\psi_1}$.
We have obtained upper bounds on the coherences
$\inp{\psi_0}{A \psi_1}$, evaluated on microscopic or macroscopic $A$.
Assuming perfect information transfer ($\s_0 \!=\! \s_1 \!= 0$),
they approach zero as $\sim \frac{1}{N}$, where $N$ is the size 
of the system.
\end{itemize}
In short: no decoherence ever occurs on perfectly closed systems,
even if they are macroscopic.
It is just very hard to distinguish coherent from incoherent states, 
creating the false impression that it does.

\newpage
The link between decoherence and macroscopic systems was brought 
forward by Klaus Hepp in his 
fundamental paper \cite{Hp},
where he considered infinite closed systems,
displaying decoherence in infinite time.
In infinite systems, the microscopic observables form a 
non-commutative C$^*$-algebra $\cal{A}$.
Its weak closure ${\cal A}''$ is considered as the (von Neumann-)algebra of
all observables.
The macroscopic observables form a commutative algebra ${\cal C}$ which is
contained in the centre of ${\cal A}''$, i.e. ${\cal C} \subset {\cal Z} = 
\{ Z \in {\cal A}'' | [Z,A] = 0 \, \, \forall A
\in {\cal A}'' \}$, yet is almost disjoint from the microscopic
observables:
${\cal C} \cap {\cal A} = {\mathbb{C}} \one$.  
Transfer of information to a macroscopic observable 
therefore implies perfect decoherence on all microscopic and macroscopic
observables (cf. section~\ref{quote2}).

Unfortunately, this transfer cannot be done by any automorphic 
time evo\-lu\-tion, since the
macroscopic observables are central. Hepp proposed information transfer by
a $t \rightarrow \infty$ limit of automorphisms. 
He was able to show that this 
causes decoherence in the weak-operator sense.   
That is, on each fixed microscopic observable, the coherence
becomes arbitrarily small for sufficiently large time $t$.   

The paper was criticized by John Bell a few years later \cite{Bell},
on the grounds that, for each fixed time $t$, there are
observables to be found on which coherence is not small.
Since Bell was of the opinion that 
a `wave packet reduction', even on closed systems, 
`takes over from the Schr\"odinger equation', this was not to his
satisfaction.
He did agree however that these observables would become arbitrarily 
difficult to observe in practice for large $t$.

By considering large but finite closed systems subject to unitary
time evolution, we hope to clarify the role that macroscopic systems 
play in making us mistake coherent superpositions for classical mixtures. 
It seems striking that the same, simple mathematics can also
be used to understand why open systems do undergo decoherence as soon as
they lose information.

\chapter{Unifying Decoherence and the Heisenberg Principle}
\label{ch:UDHP}
\setsubdir{hoofdstukken/hoofdstuk2/}



The laws of quantum mechanics impose severe restrictions on the
quality of measurement. In this chapter, we investigate 
some of these restrictions
from a quantum 
probabilistic point of view.
We introduce state independent figures of merit for
measurement performance, and then use them to prove  
sharp bounds.

\section{Introduction}
In quantum probability theory, an observable is modelled by a self-adjoint 
operator $A$ in a von Neumann algebra $\CA$, and a state by a positive,
normalized, linear map $\rho : \CA \rightarrow \mathbb{C}$. 
As discussed in chapter \ref{ch:intro}, a (normal) state induces
a probability measure on the spectrum of $A$.
It is the objective of a \emph{quantum measurement}
to portray this probability measure as 
faithfully as possible.

According to the uncertainty relation of Heisenberg, 
Kennard and Robertson 
$$
\sigma_{X} \sigma_{Y}
\geq \half | \rho( [X,Y])|\,,
$$
(see \cite{He, Ke, Ro}),
there is an inherent variance 
in the quantum state.
Furthermore, quantum theory puts severe 
re\-stric\-tions on the performance of measurement.
These restrictions, which come on top of the measurement restrictions
implied by the above uncertainty relation, fall into four distinct classes.
\begin{itemize}
\item[I]The impossibility of perfect \emph{joint measurement}.
It is not possible to perform a simultaneous measurement
of two noncommuting observables in such a way that 
both measurements have perfect quality. 
\item[II]The \emph{Heisenberg principle}, (see \cite{He}).
This states that 
quantum information cannot be extracted from a system without disturbing 
that system. 
\item[III]The impossibility of \emph{classical coding}.
One cannot extract classical information from a quantum system 
and then use this information to reconstruct the quantum state. 
\item[IV]The \emph{collapse of the wave function}.
When information is extracted from a quantum system, 
a so-called decoherence is experimentally known to 
occur on this system. 
\end{itemize}
We aim to find a quantitative description for each of these fundamental 
measurement restrictions. 

The subject of uncertainty relations in quantum measurement is already 
endowed with an extensive literature.
For example, 
the Heisenberg principle and the impossibility of joint measurement
are quantitatively illustrated in
\cite{AK,Oz, Is, Ha}.
However, 
the inequalities in these papers depend on the state $\rho$,
which somewhat limits their practical use.
Indeed, the bound on the measurement quality 
can only be calculated if the state $\rho$
is known, in which case there is no need for a 
measurement in the first place.

In section \ref{trans}, we introduce two state independent 
figures of merit: the \emph{maximal added variance}
and the
\emph{measurement infidelity}.
These will lead us quite naturally to 
state independent bounds on the performance of measurement.

In section \ref{secCS}, we introduce our main tool in the proof of
these inequalities: on operator valued Cauchy-Schwarz inequality.

%

In section \ref{JM}, we prove a sharp, state independent bound on 
the performance of jointly unbiased measurement. This generalizes the
impossibility of perfect joint measurement.

In section \ref{Heis}, we prove a sharp, state independent bound on the 
performance of a measurement in terms of the 
maximal disturbance that it causes.
This generalizes the Heisenberg principle.

In section \ref{sec:CCCSI}, using a result of Keyl and Werner, we prove
a sharp, state independent bound on the performance of classical coding.
  
In section \ref{collapse}, we investigate the collapse of the wave function.  
In contrast with the Heisenberg principle 
and its abundance of inequalities,
this phenomenon has mainly been 
investigated in specific examples (see e.g. \cite{Hp,Zu,JZ}).
Although there are some bounds on the 
remaining coherence in terms of the measurement quality
(for example \cite{Se}, but also theorem \ref{cor5} in the previous chapter),
a sharp, information-theoretic inequality does not yet appear to exist.
We will provide such an inequality, and prove
a sharp upper bound 
on the amount of coherence 
that can survive information transfer.

First of all, this generalizes the collapse of the wave function, 
and shows that decoherence is a mathematical consequence of 
information extraction.
But it also shows that no information can be extracted
if all coherence is left perfectly intact, which is
precisely the Heisenberg principle.
It is therefore a \emph{unified} description of both 
the Heisenberg principle and the collapse of the wave function.


\section{Information Transfer}\label{trans}

We wish to model information transfer between quantum systems. 
Recall from chapter \ref{ch:intro} that a quantum system
is modelled by a \emph{quantum probability space}, a pair
$(\CA,\rho)$ of a von Neumann-algebra $\CA$ and a normal state
$\rho \in \CS(\CA)$. Information transfer then proceeds as follows.
We are interested in the distribution of 
the observable $A \in \CA$, 
with the system $\CA$ in some unknown state $\rho$. We perform the
operation $T^* : \CS(\CA) \rightarrow \CS(\CB)$, and then observe 
the `pointer' $B$ in the system $\CB$ in order to obtain information on $A$. 


\subsection{General Information Transfer} \label{geninftranf}
We  
take the point of view (see \cite{Ha})
that \emph{any} CP-map 
$T : \CB \rightarrow \CA$
is an information transfer 
from \emph{any} observable $A \in \CA$ to \emph{any}
pointer $B \in \CB$.
Undeniably though, for the purpose of transferring information, 
some operations are better than others.
What does it mean for a CP-map $T$ to be `good'
at transferring information?

\subsubsection{Perfect Information Transfer}

The initial state $\rho$ defines a probability measure 
$\mathbb{P}_{i}$ on the spectrum of $A$
by $\mathbb{P}_{i}(S) := \rho( \one_{S}(A))$,
where $\one_{S}(A)$ denotes the spectral projection of $A$
associated to the set $S$. 
Similarly, the final state $T^*\rho$ defines a probability distribution 
$\mathbb{P}_{f}$ on the spectrum of $B$.

The object of information transfer is to gain information on 
$\mathbb{P}_{i}$, the probability measure of the observable $A$
\emph{before} the operation, by looking at
the probability measure $\mathbb{P}_{f}$ of the pointer $B$
\emph{afterwards}. 
We will therefore consider an information transfer to be perfect
if $\mathbb{P}_{i} = \mathbb{P}_{f}$ for any initial state $\rho$.
Translated to algebra, this reads as follows.
\begin{basdef} An information transfer $T : \CB \rightarrow \CA$ 
from $A \in \CA$ to $B \in \CB$ is called \emph{perfect} if 
$T(B) = A$ and if
the restriction of $T$ to $B''$, the von Neumann algebra generated by $B$, 
is 
an isomorphism $B'' \rightarrow A''$ of commutative non Neumann algebras.
\end{basdef}
Indeed, $\mathbb{P}_{i} = \mathbb{P}_{f}$ if and only if  
$\mathbb{E}_{i}(g) = \mathbb{E}_{f}(g)$ for any measurable function 
$g$. In other words, $T^*\rho(g(B)) = \rho(g(A))$
for all $\rho \in \CS(\CA)$ and $g \in L^{\infty}(\mathbb{R})$. 
But if $\rho(T(g(B)))$ equals $\rho(g(A))$ for any normal state $\rho$,
we must have $T(g(B)) = g(A)$ and vice versa. Since 
$X'' = \{g(X) | g \in L^{\infty}(\mathbb{R})\}$, the two notions of perfection 
coincide. 

\subsubsection{Measurement Infidelity}
Not unlike a man, an information transfer which is interesting  
is not necessarily perfect. 
The aim of information transfer is to infer 
$\mathbb{P}_{i}$ from $\mathbb{P}_{f}$, so surely
$T$ is \emph{good} if 
$\mathbb{P}_{i}$ and $\mathbb{P}_{f}$
are \emph{close}. 
 
Recall that both $\mathbb{P}_{i}$ and $\mathbb{P}_{f}$ depend on 
the state $\rho$.
If we choose a metric $d$ on the space of probability measures, 
then we can define the \emph{measurement infidelity} $\delta$
to be the worst case distance between the input and the output measure,
$$
\delta := \sup_{\rho \in \CS(\CA)} d(\mathbb{P}_{i} , \mathbb{P}_{f}) \,. 
$$
We consider $T$ to be a good information transfer from $A$
to $B$ if 
$\delta$ is small. Clearly, $T$ is perfect if and only if $\delta = 0$. 
 
As a metric on the space of probability measures, we employ the
\emph{trace distance} (a.k.a. variational distance or Kolmogorov distance) 
$$d( \mathbb{P}_{f} , \mathbb{P}_{i} ) := \sup_{S \subset \mathbb{R}} \{| \mathbb{P}_{i}(S) 
- \mathbb{P}_{f}(S) |\}\,.$$
It is
the difference between the probability that 
the event $S$
occurs in the distribution $\mathbb{P}_{i}$ and the 
probability that it occurs in the distribution 
$\mathbb{P}_{f}$, for the worst case 
Borel set $S \subset \mathbb{R}$. 
 
Let us expand this definition, and write it in a more algebraic fashion. 
We have $\mathbb{P}_{i}(S) = \rho (\one_{S}(A))$
and $\mathbb{P}_{f}(S) = \rho (T(\one_{S}(B)))$, so that
$
\sup_{\rho} d(\mathbb{P}_{i} , \mathbb{P}_{f} )$ 
can be written
$
\sup_{\rho, S} |\rho (\one_{S}(A) - T(\one_{S}(B)))|
$.
As $\sup_{\rho} |\rho(X)| = \|X\|$, this is just
$\sup_{S} \| \one_{S}(A) - T(\one_{S}(B)) \|$, and we arrive 
at the following less intuitive yet cleaner expression.
\begin{basdef}\label{measurementinfidelity} 
Let $T : \CB \rightarrow \CA$ be a CP-map. Its
\emph{measurement infidelity} $\delta$ in transferring information 
from the observable $A$ to the pointer $B$ is defined as
$ \delta := \sup_{S} \| \one_{S}(A) - T(\one_{S}(B)) \|$,
where $S$ runs over the Borel subsets of $\mathbb{R}$. 
\end{basdef}
The following is a standard example of a perfect information transfer. 
 
\subsubsection{Example: von Neumann Qubit Measurement} \label{Neuqubit} 
Let $\CA = M_2(\mathbb{C})$ and $A = \sigma_z$. 
%
Let $\OO = \{ +1 , -1 \}$ be the spectrum of $\sigma_z$. 
Denote by $L^{\infty}(\OO)$ the (commutative von Neumann) algebra
of $\mathbb{C}$-valued random variables on $\OO$.
(This is just $\mathbb{C}^{2}$ with pointwise multiplication.) 
A state $L^{\infty}(\OO) \rightarrow \mathbb{C}$ is precisely 
the expectation $\mathbb{E}$ with respect to
some probability measure $\mathbb{P}$ on $\OO$.
Define the probability measures $\mathbb{P}_{\pm}$ 
to assign probability 1 to $\pm 1$, so that 
$\mathbb{E}_{\pm}(f) = f(\pm 1)$.
We take $\CB = M_2(\mathbb{C})\otimes L^{\infty}(\Omega)$
and $B = \one \ten (\delta_{+1} - \delta_{-1})$.
 
The von Neumann-measurement 
$T : M_2(\mathbb{C}) \ten L^{\infty}(\Omega) \rightarrow M_2(\mathbb{C})$
is then defined as $T(X \ten f) := f(+1) P_+ X P_+ + f(-1) P_- X P_-$, with 
$P_+ = \up$ and $P_- = \down$.
The dual
$T^* : \CS(M_2(\mathbb{C})) \rightarrow 
\CS(M_2(\mathbb{C})) \ten \CS(L^{\infty}(\Omega))$
reads 
$T^* \rho = \rho (P_+) \up \ten \mathbb{P}_+ + 
\rho (P_-) \down \ten \mathbb{P}_-$.

In words: with probability $\rho (P_+)$, the output $+1$
occurs and the qubit is left in \mbox{state $\ket{\!\uparrow\,}$}.  
With probability $\rho (P_-)$, the output $-1$
occurs, leaving the qubit in state $\ket{\!\downarrow\,}$. 
The von Neumann-measurement $T$ is a perfect 
information 
transfer from $\s_z \in M_2(\mathbb{C})$ to 
$\one \ten (\delta_{+1} - \delta_{-1}) \in 
M_2(\mathbb{C}) \ten L^{\infty}(\OO)$. 

\subsubsection{POVMs}

Quantum measurements are often (e.g. \cite{Ho,Ha}) modelled by 
Positive Operator Valued Measures or POVMs.
From a CP-map $\CB \rightarrow \CA$ and a pointer $B \in \CB$, we may distill the 
POVM $\mu : \Sigma(\spec(B)) \rightarrow \CA$ by way of
$\mu(S) := T(\one_{S}(B))$. It will only be projection valued if the information transfer 
is perfect.

For instance, the above von Neumann qubit measurement gives rise to the
POVM 
$\mu : \OO \rightarrow M_2$ defined by 
$\mu(+1) = P_+$ and $\mu(-1) = P_-$. 

Note that from the POVM, the CP-map can be reconstructed only on $B''$.
A CP-map can thus be seen as an extension of a POVM that
keeps track of the system output as well as the measurement output.
Since we will be interested in disturbance of the system,
it is imperative that we consider the full CP-map rather than 
merely its POVM.

\subsection{Unbiased Information Transfer}

Not unlike a man, it does an information transfer credit to be \emph{unbiased}.
This means that the \emph{expectation} 
of the observable $A$ before the operation is equal to the \emph{expectation}
of the pointer $B$ afterwards, for any initial state $\rho$.
In other words, if $T^*\rho(B) = \rho(A)$ for any $\rho$.
This is equivalent to $\rho(T(B) - A) = 0$ for any $\rho$, and thus to 
$T(B) = A$.  
\begin{basdef}\label{geenvoorkeur} A CP-map $T : \CB \rightarrow \CA$ is called an \emph{unbiased}
information transfer from $A \in \CA$ to $B \in \CB$ if $T(B) = A$.
\end{basdef}
Any perfect information transfer is unbiased, but the converse does not hold.
Schematically, we have
\begin{center}
\setlength{\unitlength}{1cm}
\begin{picture}(10,3)(0,0)
\put(4.5,1.6){\mbox{\footnotesize Perfect}}
\put(5,0.85){\mbox{\footnotesize Unbiased}}
\put(2.5,0.25){\mbox{\footnotesize General}}
\put(5,1.5){\oval(8 , 2.88)}
\put(5,1.6){\oval(4.5 , 2)}
\put(5,1.7){\oval(2 , 1)}
\end{picture}\\
{\small \fig Kinds of information transfer.}
\end{center}

\subsubsection{Maximal Added Variance}\label{mav}
For unbiased information transfer, there is a natural figure of merit 
other than the measurement infidelity.
The variance of an observable $X$ is defined as 
$\var(X,\rho) = \mathbb{E}_{\rho}(X^2) - \mathbb{E}_{\rho}(X)^2$, 
or equivalently $\rho (X^2) - \rho (X)^2$.
An unbiased information transfer is good if the variance 
$\var(B,T^*(\rho))$ of the pointer in the
final state is as low as possible.


The output variance can be split in two parts.
One part $\var(A,\rho)$ is the variance of the input, which
is intrinsic to the quantum state $\rho$.
The other part $\var(B,T^*(\rho))$ $-$
$\var(A,\rho) \geq 0$ is \emph{added}
by the measurement procedure.

It would be unfair to hold the information transfer responsible for
uncertainty that was already there to begin with, so we will take 
this added variance 
to measure its performance.
The \emph{maximal added variance}, \label{maximaaa}
$$
\Sigma^2 := \sup_{\rho \in \CS(\CA)} \var(B,T^*(\rho)) -
\var(A,\rho)\,,
$$
will be our figure of (de)merit.

For example, in the case of perfect information 
transfer from $A$ to $B$, we have that
$\var(B,T^*(\rho)) = \var(A,\rho)$;
the maximal added variance is 0. 
There is uncertainty in the measurement outcome,
but all uncertainty `comes from' the quantum state, 
and none is added by the measurement procedure.

To obtain a cleaner definition, we expand 
$\var(B,T^*(\rho)) - \var(A,\rho)$ into
$
( \rho(T(B^2)) - \rho(T(B))^2 )
$ 
$-$ 
$
( \rho(A^2) - \rho(A)^2)
$. Using the fact that $T(B) = A$, we obtain
$\Sigma^2 = \sup_{\rho} \rho(T(B^2) - T(B)^2)$,
and thus $\Sigma^2 = \| T(B^2) - T(B)^2\|$.

\begin{basdef}\label{MAVdef} The \emph{maximal added variance} of an 
unbiased information transfer $T: \CB \rightarrow \CA$
from an observable $A = T(B)$ to a pointer $B$
is defined as
\begin{equation}\label{oelsmufknut}
\Sigma^2 := \| T(B^2) - T(B)^2\|.
\end{equation}
\end{basdef}

\section{An Operator Valued Form}\label{secCS}
Inspired by formula (\ref{oelsmufknut}), we define an 
operator valued 
form which allows us to consider 
$\Sigma$
as the `length' of the pointer;  
$\Sigma = \sqrt{\|(B,B)\|}$. 
\begin{basdef}\label{indef} Let $T : \CB \rightarrow \CA$ be a CP-map. 
We define the operator-valued 
sesquilinear form 
$( \,\cdot \,,\,\cdot \,)_{T} \, : \,  \CB \times \CB \rightarrow \CA$ by
\begin{equation}\label{inproduct}
(X,Y)_{T} := T(X^\dag Y) - T(X)^\dag T(Y)\,.
\end{equation} 
\end{basdef}
If there is no danger of confusion, we will drop the subscript $T$.

As we will see, the form
$( \,\cdot \,,\,\cdot \,)_{T}$
is 
sesquilinear and positive semidefinite.
%
Its most useful feature however is that it satisfies an operator valued 
Cauchy-Schwarz inequality.

\subsection{A Cauchy-Schwarz Inequality}
Before proving the operator valued Cauchy-Schwarz inequality, let us first
check that $( \,\cdot \,,\,\cdot \,)_{T}$
is indeed sesquilinear and positive semidefinite.
\begin{proposition}\label{posinp}
The form $( \,\cdot \,,\,\cdot \,)_{T}$
from definition \ref{indef} 
is linear in the second argument, it satisfies 
$(X,Y)^{\dagger} = (Y,X)$, and 
$(X,X) \geq 0\,$ for all $X \in \CB$.
\end{proposition}
\proof
Linearity in the second argument follows from linearity of $T$,
and the fact that $(X,Y)^{\dagger} = (Y,X)$ is a consequence of
$T(X^{\dagger}) = T(X)^{\dagger}$. 
The following proof of positivity is taken from \cite{Hans}.
 
As $T$ is completely positive, 
$T \otimes \id : \CB \otimes M_2(\mathbb{C}) 
\rightarrow \CA \otimes M_2(\mathbb{C})$ 
is positive, so that
$$
T \otimes \id 
\left(
{\begin{pmatrix}
X & \one\\
0 & 0
\end{pmatrix}}^{\dagger}
\begin{pmatrix}
X & \one\\
0 & 0
\end{pmatrix}
\right)
=
\begin{pmatrix}
T(X^{\dagger}X) & T(X^{\dagger})\\
T(X) & \one
\end{pmatrix}
$$
is a positive matrix.
Using the fact that $A \geq 0$
implies $B^{\dagger}AB \geq 0$, we find that
$$
{\begin{pmatrix}
\one & 0\\
-T(X) & 0
\end{pmatrix}}^{\dagger}
\begin{pmatrix}
T(X^{\dagger}X) & T(X^{\dagger})\\
T(X) & \one
\end{pmatrix}
\begin{pmatrix}
\one & 0\\
-T(X) & 0
\end{pmatrix}
$$
is also positive,
and hence $T(X^{\dagger}X) - T(X)^{\dagger}T(X) \geq 0$
as required.\qed

\noindent We are now ready to 
prove the operator valued Cauchy-Schwarz inequality.
Simple as it may be, it will be crucial to
the remainder of this chapter. All our bounds on the performance 
of information transfer will depend upon it. 
\begin{lemma}[Cauchy-Schwarz] \label{causch}
Let $T: \CB \rightarrow \CA$ be a CP-map, and define
$(X,Y) := T(X^\dag Y) - T(X)^\dag T(Y)$.  
Then for all $X,Y \in \CB$:
\begin{equation}\label{alsch}
(X,Y)(Y,X) \leq \|(Y,Y)\| (X,X).
\end{equation}
\end{lemma}
\proof 
If $T : \CB \rightarrow \CA$ is completely positive\footnote{%
We state the result for completely positive maps, but it is clear from the proof
that $4$-positivity suffices.
}, 
then so is 
$T \otimes \id : 
\CB\otimes M_2(\mathbb{C}) \rightarrow \CA \otimes M_2(\mathbb{C})$.
According to proposition \ref{posinp}, 
$$
\left( 
\begin{pmatrix}
X & Y\\
0 & 0
\end{pmatrix}
,
\begin{pmatrix}
X & Y\\
0 & 0
\end{pmatrix}
\right)_{T \otimes \id}
= 
\begin{pmatrix}
(X,X)_{T} & (X,Y)_{T} \\
(Y,X)_{T} & (Y,Y)_{T}
\end{pmatrix}
$$
is a positive matrix, so that also 
\begin{equation}\label{hakzaad2}
\begin{pmatrix}
\one & 0 \\
-(Y,X)_{T} & 0
\end{pmatrix}^{\dagger}
\begin{pmatrix}
(X,X)_{T} & (X,Y)_{T} \\
(Y,X)_{T} & (Y,Y)_{T}
\end{pmatrix}
\begin{pmatrix}
\one & 0 \\
-(Y,X)_{T} & 0
\end{pmatrix}
\end{equation}
is positive, and therefore 
$
(X,X) - 2 (X,Y)(Y,X) + (X,Y)(Y,Y)(Y,X)\geq 0
$.
Rewrite this as 
$
(X,Y)(Y,X) \leq (X,X) + (X,Y)((Y,Y) - \one)(Y,X)\,,
$
and substitute $\lambda Y$ for $Y$. If $\lambda$ 
is such that 
\begin{equation}\label{alikruikje}
|\lambda|^2 \|(Y,Y)\| \leq 1\,,
\end{equation} 
then 
$|\lambda|^{2}(X,Y)(|\lambda|^{2}(Y,Y) - \one)(Y,X)$ 
is negative, so that
\begin{equation}\label{kniepertje}
|\lambda|^2 (X,Y)(Y,X) \leq (X,X)\,.
\end{equation}
If $(Y,Y) = 0$, then (\ref{alikruikje}) holds for arbitrary $\lambda$,
so that (\ref{kniepertje})
forces $(X,Y)(Y,X)$ to be zero.
If $(Y,Y) \neq 0$, then taking $|\lambda|^2 =  \|(Y,Y)\|^{-1}$ in 
(\ref{kniepertje})
yields the inequality $(X,Y)(Y,X) \leq \|(Y,Y)\|(X,X)$.
Either way, we have (\ref{alsch}), 
which was to be demonstrated. 
\qed\vspace{-4mm}
\basremark{The special case that $(X,X)=0$ implies 
$(X,Y) =0$ is due to Reinhard Werner, see e.g. \cite{We}.}\vspace{-2mm}

\noindent We will have occasion for the following corollary.
\begin{corollary}\label{gramcolmuts}
Let $C : \CA \rightarrow \CN$ and $T : \CB \rightarrow \CA$
be completely positive maps. Define 
$( \,\cdot \,,\,\cdot \,)_{C,T} \, : \,  \CB \times \CB \rightarrow \CN$
by 
$(X,Y)_{C,T} = C \left( (X,Y)_{T}\right)$. Then
$$
(X,Y)_{C,T}(Y,X)_{C,T} \leq \|(Y,Y)_{C,T}\| (X,X)_{C,T}\,.
$$
\end{corollary}
\proof This is a minor adaptation of the proof of lemma \ref{causch}. 
The map $C \otimes \id : 
\CA \otimes M_{2}(\mathbb{C}) \rightarrow \CN \otimes M_{2}(\mathbb{C})$
is positive. If we apply it to each of the three factors in 
equation (\ref{hakzaad2}) separately, the result will still be positive.
This has the effect of replacing every $(\,\cdot\,,\,\cdot\,)_{T}$
that occurs in the entries by $(\,\cdot\,,\,\cdot\,)_{C,T}$. 
The rest of the proof then goes through uninterrupted. \qed

%
%
%

\subsection{Perfect Information Transfer}

If an information transfer is perfect, then its maximal added variance 
$\Sigma^2 = \|(B,B)\|$ will be zero. Indeed,
the initial and final probability measures are
identical, so that no variance is added. 
Using our newly acquired Cauchy-Schwarz inequality, 
we will show that the converse is also true. 
\begin{theorem}\label{struc}
Let $T: \CB \rightarrow \CA$ be a CP-map, let $B \in \CB$ be Hermitean.
Then among the statements
\begin{itemize}
\item[1] $(B,B) = 0$.
\item[2] The restriction of $T$ to $B''$, the von Neumann algebra generated
by $B$, is a ${}^*$-homomorphism $B''
\rightarrow T(B)''$. 
\item[3] $(f(B), f(B)) = 0$ for all measurable functions $f$ on the
spectrum of $B$.
\item[4] T maps the relative commutant $B' = \{X \in \CA ; [X,B] = 0\}$ 
into $T(B)'$.
\end{itemize}
the following relations hold: 
$(1) \Leftrightarrow (2) \Leftrightarrow (3) \Rightarrow(4)$. 
\end{theorem}
\proof  First, we prove $(1) \Rightarrow (2)$. 
By Cauchy-Schwarz (lemma \ref{causch}), we have
$T(B^n) - T(B)T(B^{n-1}) \leq \|(B,B)\| (B^{n-1}, B^{n-1}) =0$.
Induction on $n$ shows that $T(B^n) = T(B)^n$, so we have $T(f(B)) = f(T(B))$
for all polynomials $f$ by linearity.  Thus $T$ is a ${}^*$-homomorphism from the algebra of polynomials 
on the spectrum of $B$ to that on $T(B)$. Since $T$ is weakly continuous, this
statement extends to the algebras of measurable functions on the spectra of
$B$ and $T(B)$, isomorphic to $B''$ and $T(B)''$ respectively. 
For $(2)\Rightarrow(3)$, note that $T(f(B)^2) = T(f(B))^2$. 
For $(3)\Rightarrow(1)$, one simply takes \mbox{$f(x)=x$}. 
Finally we prove the implication $(1)\Rightarrow(4)$. Suppose that $[A,B] = 0$. Then 
$[T(B), T(A)] = T([A,B]) - [T(A), T(B)] = (A^\dag,B) - (B^\dag , A)$. 
($B$ is Hermitean.)  By Cauchy-Schwarz, the last term equals 
zero if $(B,B)$ does.
\qed
So the maximal added variance $\Sigma^2$ is zero if
and only if $T$ is a perfect information transfer,
i.e. a
${}^*$-homomorphism
$B'' \rightarrow \CA$.  

We can now give an easy proof of theorem \ref{geelsnuitje}
in chapter \ref{ch:intro}. 
\begin{theorem}
Every CP-map $T : \CB \rightarrow \CA$ with a completely
positive inverse
is an isomorphism. 
\end{theorem}
\proof \label{huuuk}
Let $S$ be the inverse of $T$. 
According to proposition~\ref{posinp}, $(A,A)_{S}\geq 0$,
so that also $T ((A,A)_{S})\geq 0$.
Expanding this as 
$T(S(A^{\dagger}A)\! -\! S(A^{\dagger})S(A))\geq 0$
and using the fact that $T$ is the inverse of $S$, we see that
$A^{\dagger}A \geq T(S(A^{\dagger})S(A))$. On the other hand,
$T(S(A^{\dagger})S(A)) - T(S(A^{\dagger}))T(S(A)) \geq 0$ implies 
$A^{\dagger}A \leq T(S(A^{\dagger})S(A))$, so that 
$A^{\dagger}A = T(S(A^{\dagger})S(A))$ for all $A \in \CA$. 
For $A = T(B)$, this reads $T(B)^{\dagger} T(B) = T(B^{\dagger}B)$
for all $B \in \CB$, or $(B,B)_{T}=0$. The Cauchy-Schwarz inequality hen
implies that also $(B'^{\dagger},B)_{T} = 0$ for all $B',B\in \CB$, so that
$T(B'B) = T(B')T(B)$. 
\qed
%

\section{Joint Measurement}\label{JM}
In a jointly unbiased measurement, information on two 
observables $A$ and $\tilde{A}$ 
is transferred to two \emph{commuting} pointers $B$ and $\tilde{B}$.  
Lemma \ref{causch} will provide a sharp bound on its performance.
\subsection{Joint Measurement Inequality} 

If $A$ and $\tilde{A}$ do not commute, then it is not possible for both
information transfers to be perfect. (See \cite{Ne}, \cite{We}.)
The following inequality states that the degree of imperfection is 
determined by the amount of
noncommutativity. 
\begin{theorem}\label{TJM}
Let $T: \CB \rightarrow \CA$ be a CP-map, let $B$,$\tilde{B}$ be commuting
Hermitean observables in $\CB$, and define $A := T(B)$, $\tilde{A} :=
T(\tilde{B})$, $\Sigma_{B}^2 := \|(B,B)\|$ and 
$\Sigma_{\tilde{B}}^2 := \|(\tilde{B},\tilde{B})\|$. Then 
\begin{equation}\label{jointje}
\Sigma_B \Sigma_{\tilde{B}} \geq \half \|[A,\tilde{A}]\| \,.
\end{equation}
\end{theorem}
\proof Since $[B, \tilde{B}] = 0$, we have 
$[\tilde{A}, A] = T([B,\tilde{B}]) - [T(B), T(\tilde{B})] = 
(B,\tilde{B}) - (\tilde{B} ,B)$. By Cauchy-Schwarz, the 
latter is at most $2 \Sigma_{B} \Sigma_{\tilde{B}}$ in norm.\qed 
The following application will show that this bound is sharp.
\subsection{The Beamsplitter as a Joint Measurement}
A \emph{beamsplitter} is a device which takes two beams of light as input. 
A certain fraction of each incident beam is refracted and the rest is
reflected, in such a way that the refracted part of the first beam
coincides with the reflected part of the second and vice versa. 
\begin{center}
\setlength{\unitlength}{1 cm}
\begin{picture}(8,3.2)
\put(0.1,2.2){\mbox{\footnotesize air}}
\put(0.1,1.7){\mbox{\footnotesize glass}}
\put(0,2){\line(1,0){8}}
\put(4,2){\vector(3,1){3}}
\put(4,2){\vector(1,-1){1.5}}
\put(2.5,0.5){\vector(1,1){1.5}}
\put(1,3){\vector(3,-1){3}}
\end{picture}\\
{\small \fig Beamsplitter.}
\end{center}
We will show that the beamsplitter serves as an optimal joint 
unbiased measurement.

\subsubsection{The Maximal Added Variance of the Beamsplitter}

In cavity QED, a single mode in the field is described by a Hilbert space
$\CH$ of a harmonic oscillator, with creation and annihilation
operators $a^{\dag}$ and $a$ satisfying $[a,a^{\dag}] = 1$, as well as 
$x = \frac{a + a^{\dag}}{\sqrt{2}}$ and 
$p = \frac{a - a^{\dag}}{\sqrt{2}i}$.
The coherent states $\ket{\alpha} = \smash{e^{-|\alpha|^2/2} \sum_{n=0}^{\infty}
\frac{\alpha^n}{\sqrt{n!}} \ket{n}}$ are dense in $\CH$, and satisfy
$a\ket{\alpha} = \alpha\ket{\alpha}$.

Quantummechanically, a beamsplitter is described by the unitary operator 
$U = \exp( \theta ( a^{\dag} \ten{a} - a \ten a^{\dag} ))$
on 
$\CH \ten \CH$. 
In terms of the coherent vectors, 
we have
$U \ket{\alpha} \ten \ket{\beta} = 
\ket{\alpha \cos(\theta)  + \beta \sin(\theta) } \ten 
\ket{-\alpha\sin(\theta)  + \beta \cos(\theta) }
$. Note that 
$U^{\dag} a\ten\one U = \cos(\theta) a\ten\one   + \sin(\theta) \one\ten a  $
and that
$U^{\dag} \one\ten a U = - \sin(\theta) a\ten\one   + \cos(\theta) \one\ten a $.
(This can be seen by sandwiching both sides between coherent vectors.)
Since the map $Y \mapsto U^{\dag}YU$ respects $+$, $\cdot$ and
${}^{\dag}$, we readily calculate 
\begin{eqnarray*}
U^{\dag} x\ten\one U &=& \cos(\theta) x\ten\one + \sin(\theta) \one\ten x,\\
U^{\dag} x^2\ten\one U &=& \cos^2(\theta) x^2\ten\one + 
2\sin(\theta)\cos(\theta) x\ten x + \sin^2(\theta) \one\ten x^2,\\
U^{\dag} \one\ten p U &=& -\sin(\theta) p\ten\one + \cos(\theta) \one\ten
p,\\ 
U^{\dag} \one\ten p^2 U &=& \sin^2(\theta) p^2\ten\one - 
2\cos(\theta)\sin(\theta) p\ten p + \cos^2(\theta) \one \ten p^2.
\end{eqnarray*}
Let us identify normal states $\rho$ on $B(\CH)$ with normalized trace class
operators $R \in \CT(\CH)$ by way of $\rho(X) = \tr(RX)$.
We are now interested in the map $ R \mapsto U R \ten \ketbra{0}
U^{\dag}$, from $\CS(\CH)$ to $\CS(\CH) \ten \CS(\CH)$.
In other words, we feed the beamsplitter only one beam of light in a state
$\rho$, the other input being the vacuum. The dual of this is the CP-map
$T : \CB(\CH) \ten \CB(\CH) \rightarrow \CB(\CH)$ defined by 
$T(Y) := id \ten \Phi ( U^{\dag}YU )$, with $\Phi$ the vacuum state
$\Phi(X) = \bra{0} X \ket{0}$.

Take for instance the pointer $B = \cos^{-1}(\theta) x \ten \one$. Then 
$T(B) = x\bra{0}\one\ket{0} + \tan(\theta) \one \bra{0}x \ket{0} = x$. 
Starting with the pointer $\tilde{B} = -\sin^{-1}(\theta) \one \ten p$ instead, 
we end up with $T(\tilde{B}) = p$. Apparently, splitting a beam of light in two parts,
measuring $x\ten\one$ in the first beam and $\one \ten p$ in the second,
and then compensating for the loss of intensity provides a 
\emph{simultaneous unbiased measurement} of $x$ and $p$ in the 
original beam. 
Since $[x,p] = i$, we must\footnote{We neglect the technical complication of
$x$ and $p$ being unbounded operators.
} have $\Sigma_{B} \Sigma_{\tilde{B}} \geq \half$.

We now calculate $\Sigma_{B}$ and $\Sigma_{\tilde{B}}$ explicitly.
From $\bra{0}x^2\ket{0} = \half$, we see that
$T(B^2) = x^2 + \smash{\half} \tan^2(\theta) \one$. Thus
$\Sigma_{B}^2 = \| (B,B) \| = \smash{\half} \tan^2 (\theta)$.
Similarly $\Sigma_{\tilde{B}}^2 = \smash{\half} \tan^{-2}(\theta)$.
We see that $\Sigma_{B}\Sigma_{\tilde{B}} = \smash{\half}$, so that 
the beamsplitter is indeed an optimal jointly unbiased measurement.

\subsection{Sharpness and Comparison}

The beamsplitter is an optimal joint
measurement in the sense that it minimizes the product 
$\Sigma^{2}_{B} \Sigma^{2}_{\tilde{B}}$ of the maximal added 
variances.
By scaling $B$, optimal joint measurements can be
found for arbitrary values of $\Sigma_{B}$ and $\Sigma_{\tilde{B}}$, 
which shows the bound in
theorem \ref{TJM} to be sharp.

The beamsplitter also performs well with other figures of merit. 
For example, if the quality of joint measurement is judged by the
state-dependent cost  
$R(T):= \var(B,T^*(\rho)) + \var(\tilde{B},T^*(\rho))$,
then at least for Gaussian $\rho$, the optimal measurement is 
the above beamsplitter with $\theta = \pi/4$. (See \cite{Ho}.)

Choosing a figure of merit is --to a degree-- a moral decision. Is it 
fair to prefer the angle $\theta = \pi/4$ over the others, or should
should we deem all angles to be of equal value? The former point of view 
leads to the figure of merit $R(T)$, the latter to 
$\Sigma_{B} \Sigma_{\tilde{B}}$. 
   
We will take the second point of view, and propose the product of 
the two maximal added variances as a figure of merit for 
joint unbiased measurement.
In chapter \ref{ch:OPJM}, we will
show that homodyne detection of 
the spontaneous decay 
of a two-level atom constitutes a joint measurement with
$\Sigma_{B} \Sigma_{\tilde{B}} = 1.056$, slightly above the bound 
$\Sigma_{B} \Sigma_{\tilde{B}} \geq 1$ 
provided by theorem \ref{TJM}.

%
%

\section{The Heisenberg Principle}\label{Heis}
The Heisenberg Principle may be stated as follows:
\begin{quote}
\begin{center}
\emph{If all states are left intact, then no 
quantum information can be extracted from a system.}
\end{center}
\end{quote}
This alludes to information transfer from an initial system 
$\CA$ to a final system consisting of two parts: the system 
$\CA$ and an ancilla $\CB$, containing the pointer~$B$. 
We thus have an information
transfer $T: \CA \ten \CB \rightarrow \CA$ from $A$ to $\one \ten B$. 

An initial state $\rho \in \CS(\CA)$ gives rise to a final state 
$T^*\rho \in \CS(\CA \ten \CB)$. Restricting this final state to the
system $\CA$ (i.e. taking the partial trace over $\CB$)
yields a `residual' state $R^*\rho \in \CS(\CA)$, 
whereas taking the partial trace over
$\CA$ yields the final state $Q^*\rho \in \CS(\CB)$ of the ancilla.
We define the CP-maps $R: \CA \rightarrow \CA$ by 
$R(A) := T(A \ten \one)$ and $Q: \CB \rightarrow \CA$ by 
$Q(B) := T(\one \ten B)$. The map $R$ describes what happens to $\CA$
if we forget about the ancilla $\CB$, and $Q$ describes the ancilla, 
neglecting the original system $\CA$.\vspace{-1mm} 
\begin{center}
\setlength{\unitlength}{1 cm}
\begin{picture}(10,3.5)
\put(2,1){\oval(1.5,1.5)}
\put(1.8,0.85){\mbox{\Large $\CA$}}
\put(5,1.75){\oval(1.5,3)}
\put(4.8,0.85){\mbox{\Large $\CA$}}
\put(4.8,1.65){\mbox{\Large $\otimes$}}
\put(4.8,2.35){\mbox{\Large $\CB$}}
\put(8,1){\oval(1.5,1.5)}
\put(8,2.5){\oval(1.5,1.5)}
\put(7.8,0.85){\mbox{\Large $\CA$}}
\put(7.8,2.35){\mbox{\Large $\CB$}}
\put(2.9,1){\vector(1,0){1.2}}
\put(3.2,1.2){\mbox{\Large $T^{*}$}}
\put(5.9,1){\vector(1,0){1.2}}
\put(6.2,1.2){\mbox{\Large $\tr_{\CB}$}}
\put(5.9,2.5){\vector(1,0){1.2}}
\put(6.2,2.7){\mbox{\Large $\tr_{\CA}$}}
\end{picture}\\
{\small \fig Constructing $R^*$ and $Q^*$ from $T^*$.}\vspace{1.5mm}
\end{center}

\noindent We wish to find a quantitative version of the Heisenberg principle,
i.e. we want to relate the imperfection of the extracted quantum-information  
to the amount of state disturbance.
\subsection{A Figure of Merit}
We have already established a figure of merit for the quality of 
information transfer.
For general information transfer, we use 
the measurement infidelity
$\delta$ of definition \ref{measurementinfidelity}, 
and for unbiased information transfer the maximal added variance 
$\Sigma^2$ of definition \ref{MAVdef}. 
We must now find a figure of merit to quantify how well 
states on $\CA$ are preserved by $R^*$.

For any initial state $\rho \in \CS(\CA)$, we want the 
final, residual state $R^*\rho$ to be as close to $\rho$ as
possible. We therefore define the 
\emph{maximal disturbance} $\Delta$ to be the maximum trace distance 
between the initial and the final state,
\begin{equation}\label{gronsbalksmurf}
\Delta := \sup_{\rho \in \CS(\CA)} d( \rho , R^* \rho)\,.
\end{equation}

The trace distance or Kolmogorov distance 
between two states $\rho$ and $\sigma$ is defined as the worst case 
difference in probability that an event $P$ occurs, 
$$
d(\rho,\sigma) := \sup_{P} (\rho(P) - \sigma(P))\,,
$$
where $P$
runs over the projections (or `events') in $\CA$.  
It is shown in \cite{NC} that one may equivalently have
$P$ run
over all positive operators $0 \leq P \leq \one$, rather than
just the projections.

The name `trace distance' stems from the fact that
if $\rho = \tr(R\,\,\cdot\,\,)$ and $\sigma = \tr(S\,\,\cdot\,\,)$ 
for trace class operators 
$R$ and $S$, 
then \label{treesafstand}  
$d(\rho,\sigma) = \frac{1}{2} \tr(|R - S|)$. 

In order to arrive at a cleaner definition, replace 
$R^*(\rho) (P)$ in
$$
\Delta = \sup
\{|R^*(\rho) (P) - \rho (P) |  \,\, ; \,\, 
\rho \in \CS(\CA) \,,\, 0 \leq P \leq \one\}
$$ 
by $\rho(R(P))$ in order to obtain
$\Delta = \sup_{0 \leq P \leq \one} \|R(P) - P\|$.
Summarizing, we arrive at the following quantification of state disturbance.
\begin{basdef}\label{hoedjekapoetje} 
The \emph{maximal disturbance} $\Delta$ 
of a map $R : \CA \rightarrow \CA$
is given by 
$$
\Delta := \sup \{ 
\mathrm{ \mbox{$\| R(P) - P \|$} } \,;\,
P \in \CA \,,\, 0 \leq P \leq \one 
\}\,.
$$
\end{basdef}
One has $\Delta = 0$ if and only if the operation $R$ leaves all
states of the original system perfectly in place, or equivalently
iff $R = \id$.

%

\subsection{HP for the Unbiased Case}
We first turn our attention to unbiased information transfer. 
The imperfection of the transfer is then captured in the maximal
added variance $\Sigma^2$.

The Heisenberg principle only holds for quantum-information.
Classical observables are contained in the centre 
$\CZ = \{A \in \CA \,; \, [A,X] = 0 \,\,\, \forall X \in \CA \}$, 
whereas quantum observables are not. 
The degree in which an observable $A$ is `quantum'
is given by its distance to the centre 
$d(A,\CZ) = \inf_{Z \in \CZ} \|A - Z\|$.
In the following, we will take the algebra of observables to be
$\CA = B(\CH)$ for some Hilbert space $\CH$.
The centre is then simply $\mathbb{C} \one$.

\begin{theorem}\label{HP}
Let $T : B(\CH) \otimes \CB \rightarrow B(\CH)$ be a CP-map, let $B \in \CB$
be Hermitean. Define $A := T(\one \otimes B)$, and 
$\Sigma^2 := \|(\one \ten B,\one \ten B)\|$. 
Furthermore, define 
$\Delta := \sup_{P} \{ \| R(P) - P \| \}$,
with $R$ the restriction of $T$ to $B(\CH) \ten \one$.  
Then 
\begin{equation}\label{bound}
\Sigma \geq d(A,\CZ) \frac{\half -\Delta}{\sqrt{\Delta(1 - \Delta)}}\,.
\end{equation}
This bound is sharp in the sense that for all $\Delta \in [0,\half]$,
there exist $T$ and $A$ for which (\ref{bound}) attains equality. 
\end{theorem}
We prepare the ground for theorem \ref{HP} by means of the following lemma.
\begin{lemma}\label{schatje}
Let $P$ be a projection, and let $0 \leq X \leq \one$ be such that
$\|P - X\| \leq \Delta$. Then $\|X(\one - X)\| \leq \Delta(1 - \Delta)$.
\end{lemma}
\proof
Let $x$ be in $\spec(X)$, the spectrum of $X$. 
As $0 \leq X \leq 1$, we have 
$0 \leq x \leq 1$.
Without loss of generality, assume that there
exists a normalized eigenvector $\psi$ such that $X \psi = x
\psi$. (If this is not the case, one may complete the proof using
approximate eigenvectors.)
Decompose $\psi$ over the eigenspaces of $P$, i.e. write
$\psi = \chi_1 + \chi_0 $, with $\chi_1$ in the image and
$\chi_0$ in the kernel of $P$.
Then $(X - P) \psi = (x - 1) \chi_1 + x \chi_0$.
Since $\|\chi_1\|^2 + \|\chi_0\|^2 = 1$, the inequality
$\Delta^2 \geq \|(X-P) \psi \|^2 = 
(x-1)^2 \|\chi_1\|^2 + x^2 \|\chi_0\|^2$ implies that
either $(1-x)^2 \leq \Delta^2$ or $x^2 \leq \Delta^2$. 
We conclude that $\spec(X) \subseteq [0, \Delta] \cup
[(1 - \Delta),1]$.
This implies that 
$\spec(X - X^2) 
\subseteq [0,\Delta (1 - \Delta )]$, as desired.
\qed
We proceed with the proof of theorem \ref{HP}.

\proof
We may assume $\Delta < \half$, as inequality (\ref{bound})
is trivially satisfied otherwise. 
Let $x := \sup(\spec(A))$ and $y := \inf(\spec(A))$, so that 
$d(A,\CZ) = (x - y)/2$.
Without loss of generality, assume that there exist normalized 
eigenvectors $\psi_x$ and $\psi_y$ satisfying 
$A \psi_x = x \psi_x$ and $A \psi_y = y \psi_y$.
(If this is not the case, one uses 
approximate eigenvectors.)
Define $\psi_{-} := (\psi_x + \psi_y)/\sqrt{2}$,  
and $\tilde{B} := \ketbra{\psi_-}$.  

On the one hand, we have 
$\|R(\tilde{B}) - \tilde{B}\| \leq \Delta$, so that
$\|[A, R(\tilde{B}) - \tilde{B}]\| \leq 2 \Delta d(A,\CZ)$.
On the other hand, we have $\|[A,\tilde{B}]\| = d(A,\CZ)$.
We then use the triangle inequality to see that 
\begin{eqnarray*}
\|[T(\one \ten B), T(\tilde{B} \ten \one)]\| &=& \|[A,R(\tilde{B})]\| \\
 &=& \|[A,\tilde{B}] + [A, R(\tilde{B}) - \tilde{B}]\|\\
& \geq & d(A,\CZ) (1 - 2\Delta)\,.
\end{eqnarray*}
%
This brings us in a position to apply theorem \ref{TJM}
to the commuting pointers $\tilde{B} \ten \one$ and $\one \ten B$,
yielding 
\begin{equation}\label{goudvis}
2 \Sigma \sqrt{\|(\tilde{B} \ten \one , \tilde{B} \ten \one)_{T}\|} \geq 
d(A,\CZ) (1 - 2 \Delta)\,.
\end{equation}

In order to estimate 
$\|(\tilde{B} \ten \one, \tilde{B} \ten \one)_{T}\| = 
\|(\tilde{B} , \tilde{B})_{R}\|$,
we use the fact that 
$\tilde{B}^2 = \tilde{B}$
to write
$(\tilde{B}, \tilde{B})_{R} = 
R(\tilde{B}) - R(\tilde{B})^2$.
Lemma \ref{schatje} with $X = R(\tilde{B})$ and $P = \tilde{B}$
then yields $\|R(\tilde{B}) - R(\tilde{B})^2\| \leq \Delta(1 - \Delta)$.
Inserting
$$
\sqrt{\Delta(1 - \Delta)} \geq
\sqrt{\|(\tilde{B} \ten \one, \tilde{B} \ten \one)_{T}\|}
$$
into inequality (\ref{goudvis}) yields
$2\Sigma \sqrt{\Delta(1 - \Delta)} \geq d(A,\CZ) (1 - 2\Delta)$,
which was to be demonstrated. 
For sharpness of the bound, see section \ref{sharp}.
\qed 


\noindent In the case of no disturbance, $\Delta = 0$, we see that 
$\Sigma \rightarrow \infty$.
No information
transfer from $\CA$ is allowed if all states on $\CA$ are left intact. 
This is Werner's (see \cite{We})
formulation of the Heisenberg principle.

In the opposite case of perfect information transfer, $\Sigma = 0$, 
inequality \ref{bound} shows that $\Delta$ must equal at least one half.
We shall see in section \ref{collapse} that this 
corresponds with a so-called `collapse of the wave
function'. \newpage
These two extreme situations are connected by theorem \ref{HP} in a
continuous fashion,
as indicated in the graph below.\\[-1.1cm]
\begin{center}
\begin{tabular}{p{8cm}}
\includegraphics[width = 6.5 cm, viewport = 0 0 500 500]{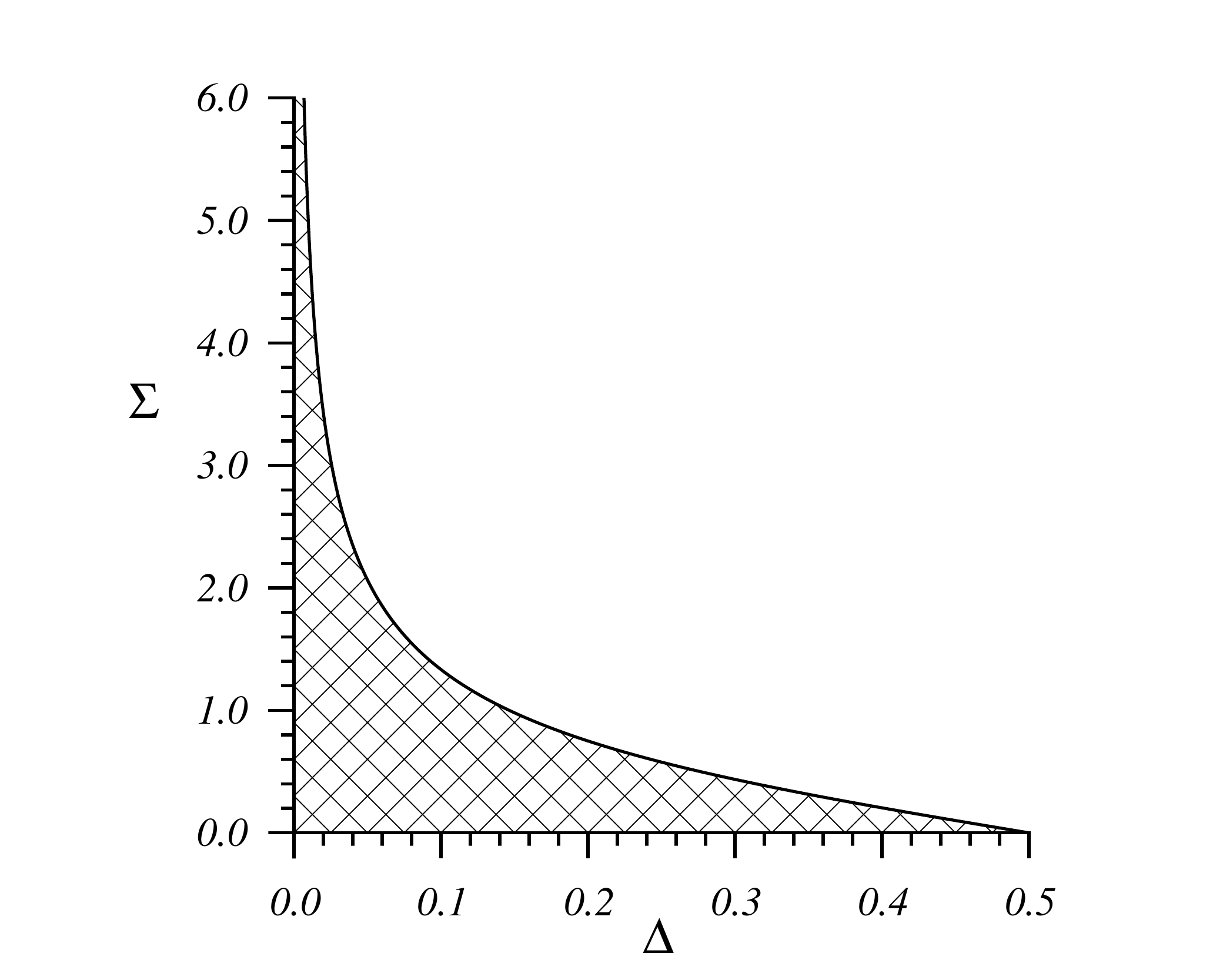}\\[-3 mm]
{\small \fig 
The combinations $(\Delta,\Sigma)$ below the curve 
are forbidden, those above are allowed. 
(With $d(A,\CZ)=1$.)}
\end{tabular}
\end{center}
The upper left corner of the curve illustrates 
the Heisenberg principle,
whereas in the the lower right corner, we can see the collapse of the
wave function.

\subsection{HP for General Information Transfer}

We now prove a version of the Heisenberg Principle 
for general information transfer. We formulate it as a trade-off
between the maximal disturbance $\Delta$ and the measurement 
infidelity $\delta$.
\begin{corollary}\label{hpdel}
Let $T : B(\CH) \otimes \CB \rightarrow B(\CH)$ be a CP-map, let 
$A \in B(\CH)$ and
$B \in \CB$ be Hermitean, $A \notin \CZ = \mathbb{C} \one$. 
Define 
$\Delta := \sup_{P} \{ \| R(P) - P \| \}$,
with $R$ the restriction of $T$ to $B(\CH) \ten \one$.  
Define
$\delta := \sup_{S} \{ \|T(\one \ten \one_{S}(B)) - \one_{S}(A)\| \}$.
Then, for $\delta$ and $\Delta$ in $[0,\half]$, we have
\begin{equation}\label{bound2}
(\half - \delta)^2 + (\half - \Delta)^2 \leq \fourth\,.
\end{equation}
This bound is sharp in the sense that for all 
$\Delta \in [0,\half]$, there exists a $T$
for which (\ref{bound2}) attains equality.
\end{corollary}
\proof
Choose a nontrivial subset $S$
of $\spec(A)$ and put $P := \one \ten \one_{S}(B)$. 
Since $ \| T(P) - \one_{S}(A) \| \leq \delta$,
lemma \ref{schatje} yields 
$\Sigma^2 = \|T(P) - T(P)^2 \|
\leq \delta(1-\delta)$.
From the proof of the lemma, we see that
$\spec(T(P))$ is contained in $[0,\delta] \cup [1-\delta, 1]$.
Since $\spec(T(P))$ contains points in both
$[0,\delta]$ and $[1 - \delta,1]$,
we have $d(T(P), \CZ) \geq \half- \delta$.
Then apply theorem \ref{HP} to the pointer $P$
to obtain 
$
\sqrt{\delta(1-\delta)} \geq
(\half - \delta)(\half - \Delta)/\sqrt{\Delta(1-\Delta)}
$,
or equivalently 
$(\half - \delta)^2 + (\half - \Delta)^2 \leq \fourth$.
For sharpness, see section \ref{sharp}. 
\qed

\noindent A measurement which does not disturb any state ($\Delta = 0$)
cannot yield information ($\delta \geq \half$).
This is the Heisenberg principle.
On the other hand, perfect information ($\delta = 0$)
implies full disturbance ($\Delta \geq \half$), corresponding
to the collapse of the wave function.
Both extremes are connected in a continuous fashion, 
as depicted below.\vspace{-5mm}
\begin{center}
\begin{tabular}{p{8cm}}
\includegraphics[width = 6.5 cm, viewport = 0 0 500 500]{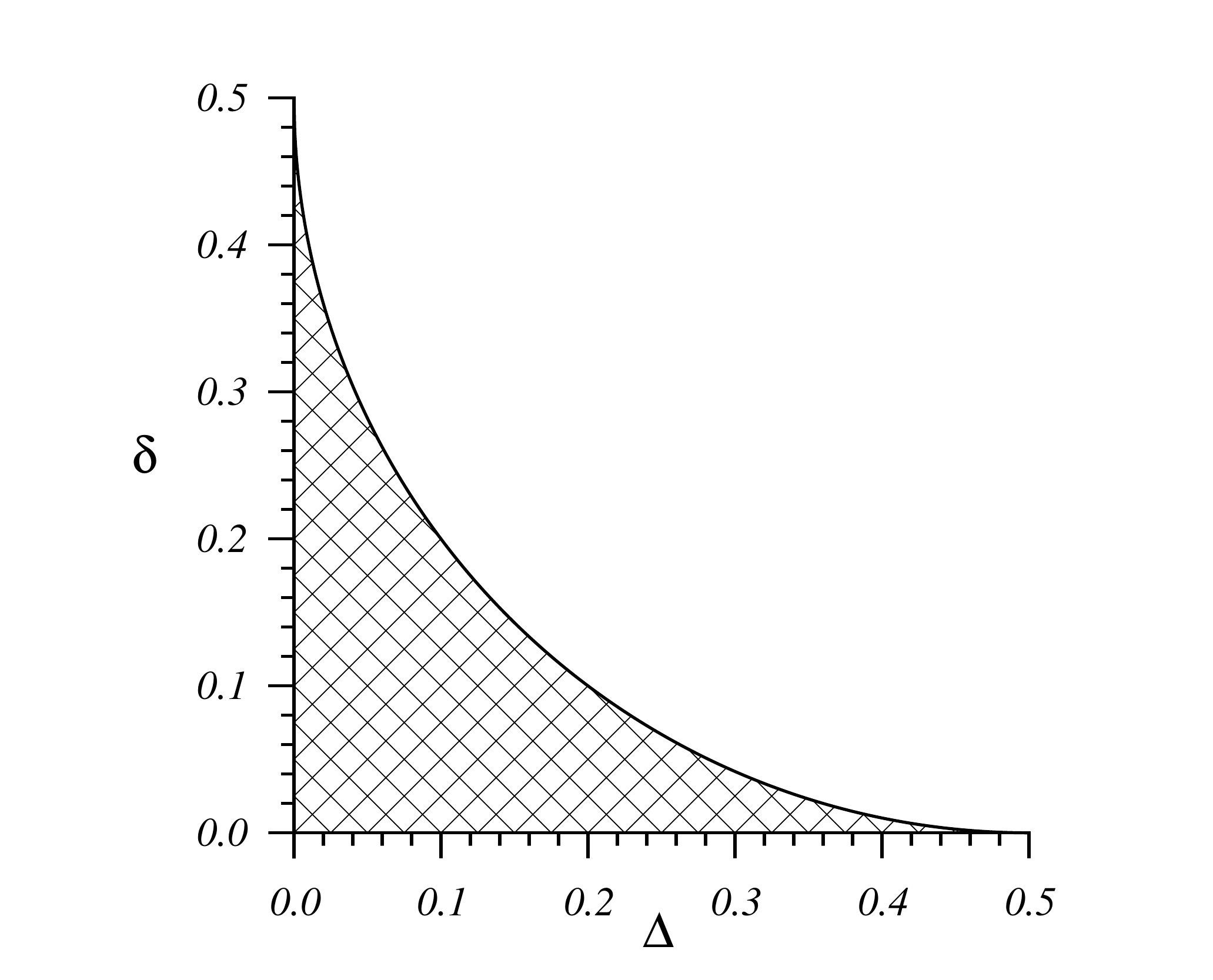}\\[-4 mm]
{\small \fig 
The combinations $(\Delta,\delta)$ below the curve 
are forbidden, those above are allowed.}
\end{tabular}
\end{center}
 
\subsection{Application: Resonance Fluorescence} \label{Resflu}
Corollary \ref{hpdel} may be used to determine the minimum amount of 
disturbance if the quality
of the measurement is known. Alternatively, if the system is only
mildly disturbed, one may  
find a bound on the attainable measurement quality. 
Let us concentrate on the latter option.

We investigate the radiation emission of a laser-driven two-level atom.
The emitted EM-radiation yields information on the atom.
A two-level atom (i.e. a qubit) only has three independent
observables: $\s_x$, $\s_y$ and $\s_z$. 
There are various ways to probe the EM field: photon counting, 
homodyne detection, heterodyne 
detection, et cetera. 
For a strong ($\OO \gg 1$) resonant 
($\omega_{\mathrm{laser}} = \omega_{\mathrm{atom}}$)
laser, we will use corollary \ref{hpdel} to prove that \emph{any} 
EM-measurement of $\s_x$, $\s_y$ or 
$\s_z$ will have a 
measurement infidelity of at least
$$
\delta \geq 
\half - \half \sqrt{1 - e^{-\frac{3}{2}\lambda^{2} t}}\,,
$$
with $\lambda$ the coupling constant.
For a measurement with two outcomes,
$\delta$ is the maximal probability of getting the
wrong outcome.

\subsubsection{Unitary Evolution on the Closed System} 
The atom is modelled by the Hilbert-space $\mathbb{C}^2$ (only two
energy-levels are deemed relevant). 
In the field, we discern a 
\emph{forward} and a \emph{side} channel, each described by a 
bosonic Fock-space $\CF$.  
The laser is put on the forward channel, which is thus initially in the 
state $\Phi_\OO$, the coherent state with frequency $\omega$ and strength
$\OO$. 
(The field strength is parametrized by the frequency of the induced Rabi-oscillations).
The side channel starts in the vacuum state $\Phi_0$.

If we permit ourselves the luxury of identifying a state with its
trace class operator, then    
the state at time $t$ is given by 
$$
T_{t}^* (\rho) =  U( t ) (\rho \ten \Phi_{\OO} \ten \Phi_0 )U^{\dag}(t)\,,
$$
with time evolution 
$$
\frac{d}{dt}  U_t = -i(H_S + H_F + \lambda H_I) U_t \,.
$$
$H_S \in B(\mathbb{C}^2)$ is the Hamiltonian of the two-level atom,
$H_F \in B(\CF \otimes \CF)$ that of the field and 
$\lambda H_I \in B(\mathbb{C}^2) \ten B(\CF \ten \CF)$ is the 
interaction-Hamiltonian.
Define the interaction-picture time evolution by 
$$
\hat{T}^*_t (\rho):= U_1 (t)^{\dag}U_2 (t)^{\dag} T^*_t (\rho) 
U_2 (t) U_1 (t)\,,
$$
where
$U_1 (t) := e^{-iH_S t}$ and
$U_2 (t) := e^{- i H_F t}$
are the 1-parameter groups of unitaries for
the \smash{`unperturbed'} time evolution. 

We now investigate $\hat{T}_t$ instead of $T_t$.
We are looking for a bound on the measurement infidelity 
$\delta = \sup_{S}\{ \|T(\one \ten \one_S(B)) - \one_S(A)\| \}$ of $T$, but
if we move to the interaction picture and define $\hat{B} := U_2\smash{^{\dag}}B U_2$, then 
$\hat{T}(\one_S (\hat{B})) = T(\one_S(B))$, so that
$\smash{\hat{\delta}} = 
\sup_{S}\{ \|\smash{\hat{T}}(\one \ten \one_S(\smash{\hat{B}})) - \one_S(A)\| \}
= \delta$.
If we find the interaction-picture disturbance $\smash{\hat{\Delta}}$, then
corollary \ref{hpdel} will yield a bound on 
$\smash{\hat{\delta}}$, and therefore also on $\delta$.

In the weak coupling limit $\lambda \downarrow 0$, the interaction picture 
time evolution $\hat{T}_t$ is given by 
$\hat{T}_{t}^* (\rho) = 
\hat{U}( t/\lambda^2 ) 
(\rho \ten \Phi_{\OO} \ten \Phi_0 )
\hat{U}^{\dag}(t/\lambda^2)$,
where the evolution of the unitary cocycle
$t \mapsto \hat{U}_t$ is described (see \cite{AFL}) by a 
Quantum Stochastic Differential Equation 
or QSDE. 
Explicitly calculating the maximal added variances $\Sigma^2$ by 
solving the
QSDE is in general rather nontrivial, if indeed possible at all. 

In chapter \ref{ch:OPJM}, we will do this for the
case of spontaneous decay, i.e. $\OO = 0$, with the map $\hat{T}_t$ 
restricted to the commutative algebra of homodyne 
measurement results. But here we take an other approach.    

\subsubsection{Master Equation for the Open System}

Fortunately, in contrast to the somewhat complicated time evolution
$\hat{T}_t$ of the 
combined system, the evolution \emph{restricted to the 
two-level system} is both well-known and uncomplicated.
If we use $\lambda^{-2}$ as a unit of time, then
the restricted evolution 
$\hat{R}^*_{t}(\rho) := \tr_{\CF \ten \CF} \hat{T}^{*}_{t} (\rho)$
of the two-level system
is known (see \cite{Bo}) to satisfy the Master equation 
\begin{equation} \label{bitrev}
\frac{d}{dt} \hat{R}_{t}^* (\rho) = L(\hat{R}_{t}^* (\rho))\,,
\end{equation}
with Liouvillian
$
L(\rho) := 
\half i \OO [e^{-i(\omega - E) t} V + e^{i(\omega - E)t} V^{\dag}, \rho] -
\half \{ V^\dag V , \rho  \} + V \rho V^{\dag}
$. In this expression,
$E$ is the energy-spacing of the two-level atom 
and $V^\dag = \s_+$, $V = \s_-$ are its raising and lowering operators.
In the case $\omega = E$ of \emph{resonance fluorescence}, 
we obtain  
$$
L(\rho) = 
\half i \OO [ V + V^{\dag}, \rho] -
\half \{ V^\dag V , \rho  \} 
+ V \rho V^{\dag}\,.
$$
If we parametrize a state by its Bloch-vector 
$\hat{R}^*_t (\rho) = \half(\one + x \s_x + y \s_y + z \s_z)$, 
then equation \ref{bitrev}
is simply the following differential equation on $\mathbb{R}^3$:
$$
\frac{d}{dt} \begin{pmatrix} x \\ y \\ z \end{pmatrix} =
\begin{pmatrix}		-\half & \,0& \, 0 \\ 
			0 &-\half& \, \OO \\
			0 & -\OO & -1
\end{pmatrix}
\begin{pmatrix} x \\ y \\ z
\end{pmatrix}
- \begin{pmatrix} 0 \\ 0 \\ 1
  \end{pmatrix}			
$$
This can be solved explicitly.
For $\OO \gg 1$, the solution approaches
$$\begin{pmatrix} x \\ y \\ z \end{pmatrix} = 
\begin{pmatrix} e^{-\half t} & 0 & 0 \\
		0 & \, e^{-\frac{3}{4} t} \cos(\OO t) & e^{-\frac{3}{4} t}	
		\sin( \OO t)	\\
		0 & - e^{-\frac{3}{4} t} \sin(\OO t)& 
		e^{-\frac{3}{4} t} \cos(\OO t)
\end{pmatrix}		
\begin{pmatrix}		
x_0 \\ y_0 \\ z_0
\end{pmatrix}\,.
$$
If we move to the interaction picture 
once more to counteract the Rabi oscillations, 
i.e. with 
$U_1(t) = e^{\frac{i}{2} \OO t \s_x}$ and $U_2 = \one$,
we see that the time evolution is transformed to
$$
\begin{pmatrix}
x \\ y \\ z 
\end{pmatrix}
= 
\begin{pmatrix}e^{-\half t} & 0 & 0 \\
		0 & \, e^{-\frac{3}{4} t}  & 0	\\
		0 &  0 & e^{-\frac{3}{4} t}
\end{pmatrix}		 
\begin{pmatrix}		
x_0 \\ y_0 \\ z_0
\end{pmatrix}\,.
$$
Since the trace distance $D( \rho , \tau )$ is exactly half the Euclidean distance between the
Bloch vectors of $\rho$ and $\tau$, (see \cite{NC}), we see that $\Delta = \half(1- e^{-\frac{3}{4}t})$.
For any measurement of $\s_x$, $\s_y$ or $\s_z$, we therefore have
$
\delta \geq \half - \half \sqrt{1 - e^{-\frac{3}{2} t}}
$ by corollary \ref{hpdel}
(remember that $t$ is in units of $\lambda^{-2}$).\vspace{-5mm}
\begin{center}
\begin{tabular}{p{8cm}}
\includegraphics[width = 6.5 cm, viewport =  0 0 500 500]{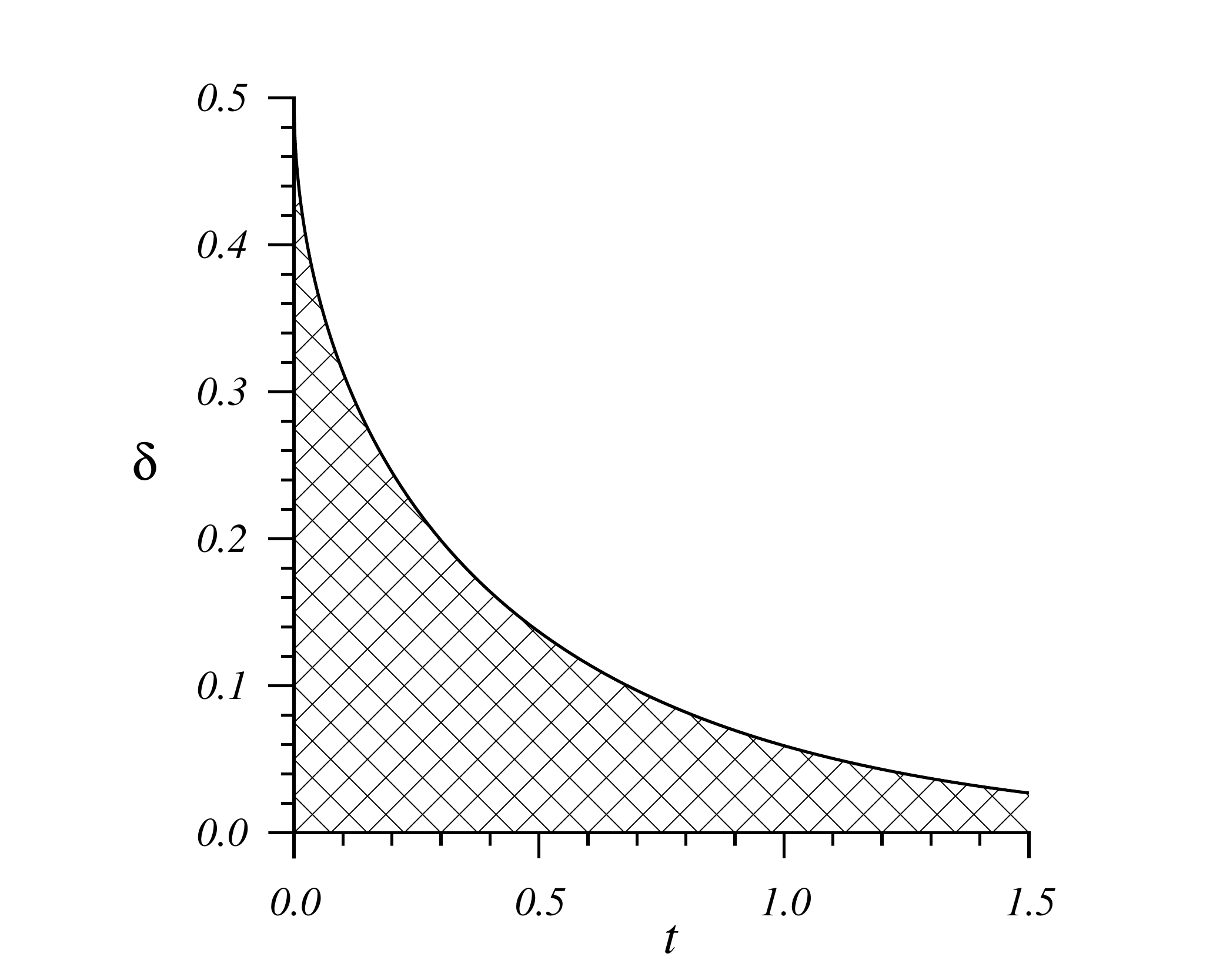}\\[-4 mm]
{\small \fig Lower bound on $\delta$ in terms of $t$ 
(in units of $\lambda^{-2}$).}
\end{tabular}
\end{center}

\section{Classical Coding}
\label{sec:CCCSI}

In classical coding, a single quantum state is encoded into 
classical information. Decoding this classical
information in order to regain the original quantum state is
known to be impossible.
However, one can attempt to construct a state which comes as close as 
possible.
We give bounds on the smallest possible trace distance 
between the original and the decoded state
which can be reached.
 
We start with an approach based on the 
operator-valued Cauchy-Schwarz inequality. 
It will yield a bound, but in contrast to the
case of joint measurement and the Heisenberg principle, 
it will not be sharp.
We will then prove a sharp bound, 
starting from Keyl and Werner's
no-cloning theorem \cite{KW}. 

\subsection{Classical Coding}\label{klaskonijn}

We are interested in operations that take 
quantum states into classical probability measures 
and vice versa.
As we have seen in chapter \ref{ch:intro}, 
a classical probability space $(\Omega, \Sigma, \mathbb{P})$ 
corresponds to a \emph{commutative} quantum probability space
$(\CA,\rho)$.
The algebra is given by the random variables `up to measure zero',
$\CA := L^{\infty}(\Omega , \Sigma , \mathbb{P})$, and
the state is simply the expectation 
$\rho(f) = \int_{\Omega} f(\omega) \mathbb{P}(d\omega)$.

%
An operation that maps quantum states in $\CS(B(\CH))$ to 
probability measures in $\CS(\CA)$ must therefore
be the dual of a CP-map $C : \CA \rightarrow B(\CH)$, with 
$\CA := L^{\infty}(\Omega , \Sigma , \mathbb{P})$.
According to section \ref{geninftranf}, such a map is precisely a 
POVM $\Sigma \rightarrow B(\CH)$.
Similarly, an operation that maps classical probability distributions
into quantum states must be the dual of a CP-map $D : B(\CH) \rightarrow \CA$.

\subsubsection{Classical Coding}

By \emph{classical coding}, we mean the following procedure.
First, classical information is extracted from a quantum system.
This is described by the coding map $C : \CA \rightarrow B(\CH)$,
with $\CA = L^{\infty}(\Omega, \Sigma)$.
Then, on the basis of this classical information, 
the original state is reconstructed as well as possible 
by means of some decoding procedure. This is described by a CP-map  
$D : B(\CH) \rightarrow \CA$. In the dual picture, this gives
$$
\CS(B(\CH)) \stackrel{C^*}{\rightarrow} \CS(\CA) 
\stackrel{D^*}{\rightarrow} \CS(B(\CH))\,.
$$
The coding procedure is flawless iff every state is reconstructed
perfectly, i.e. iff $C\circ D : B(\CH) \rightarrow B(\CH)$ is the identity.  
 
It is well known \cite{We, Hans} that perfect classical coding is impossible.
So let us investigate how close we can come to perfection.

\subsubsection{A Figure of Merit}
To define a figure of merit for classical coding, we proceed in close 
analogy to the maximal 
disturbance $\Delta$. 

For any initial state $\rho \in \CS(\CB(\CH))$, we want the 
final, reconstructed state $D^*C^* \rho$ to be as close to $\rho$ as
possible. 
We therefore define the 
\emph{coding imperfection} $\Delta$ to be the maximum trace distance 
between the initial state and its reconstruction,
$$
\Delta := \sup_{\rho \in \CS(\CB(\CH))} d( \rho , D^*C^* \rho)\,.
$$

Comparing this with 
equation (\ref{gronsbalksmurf}), we see that
this is precisely the maximal disturbance 
for the map $R^* = D^*C^*$. 
We can follow the reasoning leading up to definition
\ref{hoedjekapoetje} word by word to arrive at the following reformulation.
\begin{basdef}
A \emph{classical coding procedure} consists of 
a commutative von Neumann algebra $\CA$,
and a pair of CP-maps 
$C : \CA \rightarrow \CB(\CH)$ and $D : \CB(\CH) \rightarrow \CA$. 
Its \emph{coding imperfection} is
defined as
$$
\Delta = \sup_{0 \leq B \leq \one} \|B - CD(B)\|\,.
$$
\end{basdef}

The coding imperfection quantifies the quality of the
coding procedure: the better the procedure, the smaller the 
imperfection.
One would only have $\Delta = 0$ if all states could be encoded
perfectly.
We now investigate how close to zero $\Delta$ can come. 

\subsection{Classical Coding Inequality}

We use the Cauchy-Schwarz inequality to obtain a
the bound on $\Delta$ for any classical coding procedure.
(The proof is inspired by the `impossibility of classical coding' version
in \cite{Hans}.)
\begin{theorem}\label{trinklet}
Any classical coding procedure has coding imperfection
$$\Delta \geq (3 - \sqrt{5})/4\,.$$
\end{theorem} 
\proof
Take two orthogonal vectors $\psi, \phi \in \cal{H}$,
and define $X$ to be the projection on $\psi$,
and $Y$ the projection on $(\psi + \phi)/\sqrt{2}$.
We have $\|[X,Y]\| = 1/2$. Since $\CA$ is Abelian, we have
$D(X)D(Y) = D(Y)D(X)$, and we can write
\begin{eqnarray}\label{mandofiraat}
[X,Y] &=& [X,Y] - CD([X,Y]) + \label{mangtor} \\
  & &  C\left( D(XY) - D(X)D(Y) \right) - \label{koptig}\\
     & & C \left( D(YX) - D(Y)D(X) \right)\,. \nonumber
\end{eqnarray}
We will bound the r.h.s. in terms of $\Delta$. Remembering 
that the l.h.s. is at least $1/2$ in norm will then yield
a minimum value on $\Delta$.  

We start with (\ref{mangtor}).
Like any antihermitean operator, $[X,Y]$ can be
written as $[X,Y] = i(A_+ - A_-)$, with 
$0 \leq A_\pm \leq \|[X,Y]\| I = \frac{1}{2} I$. Therefore,
we have $\|[X,Y] - CD([X,Y])\| \leq 
\Delta \|A_+\| + \Delta \|A_-\| \leq \Delta$.

We then proceed with (\ref{koptig}).
Consider the positive semidefinite sesquilinear form
$(X,Y)_{C,D} := C\left( D(X^\dagger Y) - D(X)^\dagger D(Y)\right)$,
in terms of which the expression (\ref{koptig}) equals
$(X,Y)_{C,D} - (Y,X)_{C,D}$.
According to corollary \ref{gramcolmuts}, we then have
$$
\|(X,Y)_{C,D} - (Y,X)_{C,D}\| \leq 2 \sqrt{\|(X,X)\|\|(Y,Y)\|}\,.
$$

Now $(X,X)_{C,D} = (X,X)_{C D} -(D(X),D(X))_C \leq 
(X,X)_{C D}$, which in turn can be written
$(X,X)_{C D} = C D (X) (I - C D (X))$
because $X^\dagger X = X$.

As $X$ is a projection, and $\|X - CD(X)\| \leq \Delta$, 
we can apply 
lemma \ref{schatje} 
to obtain
$\|C D (X)  - C D (X)^2)\| \leq \Delta (1 - \Delta)$,
or equivalently $\|(X,X)_{C D}\| \leq \Delta(1-\Delta)$.
Of course, we similarly have
$\|(Y,Y)_{C D}\| \leq \Delta(1-\Delta)$.
All in all, we can bound (\ref{koptig})
as $\|(X,Y)_{C,D} - (Y,X)_{C,D}\| \leq 2 \Delta(1-\Delta)$.
  
We conclude that $1/2 = \|[X,Y]\| \leq \Delta + 2\Delta(1-\Delta)$,
or $(\Delta - 3/4)^2 \leq 5/16$. Thus 
$\Delta \geq (3 - \sqrt{5})/4$, which was to be proven. \qed
The bound $\Delta \geq (3 - \sqrt{5})/4$ is unfortunately not
the best possible one.
However, the good news is that  
the best bound follows easily from the no-cloning theorem
of Keyl and Werner \cite{KW}.
\subsection{The No-Cloning method}\label{secclo}
Each classical coding 
procedure yields a `cloning machine'.
All one has to do is to decode the classical information 
an arbitrary amount $M$ of times, rather than just once.
A bound for cloning will therefore induce a bound for 
classical coding. 
\begin{theorem}\label{KeylWerner}
Any classical coding procedure on a Hilbert space $\CH$ 
of dimension $d$ has coding imperfection 
$$
\Delta \geq \frac{d-1}{d+1}\,.
$$
\end{theorem}
\proof
We set
$\CH = \mathbb{C}^{d}$, 
and 
$\CA = L^{\infty}(\Omega)$  with  $\Omega = \{1,2, \ldots ,n\}$.
The `diagonal' map $\Omega \rightarrow \Omega^{M} \,:\, 
i \mapsto (i, \ldots ,i)$ induces the `classical 
cloning' map $K : \CA^{\otimes M} \rightarrow \CA$, i.e.
$(K f)(i) := f(i,\ldots,i)$. 
(Note that for e.g. $\Omega = \mathbb{R}$, a cloning map 
poses difficulties.)

The composition $ T := C \circ K \circ D^{\otimes M}$, mapping
$B(\mathbb{C}^{d})^{\otimes M}$ to $B(\mathbb{C}^{d})$,
is a so-called $1 \rightarrow M$ cloner.
By construction of the operation $T$, it is clear that
$T(I \otimes \ldots \otimes B \otimes \ldots \otimes I)
=CD(B)$.
The main result of \cite{KW} then says that 
$\sup_{0 \leq B \leq I} 
{\|T(I \otimes \ldots \otimes B \otimes \ldots \otimes I) - B\|}
\geq \frac{(M-1)}{M}\frac{d-1}{d+1}$.
Since $M \in \mathbb{N}^+$ was arbitrary,
this implies $\Delta \geq \frac{d-1}{d+1}$.
\qed
The number $(3 - \sqrt{5})/4 \approx 0.19$, obtained from the
Cauchy-Schwarz inequality, is inferior to the $\frac{1}{3} \approx 0.33$
from the no-cloning theorem in the best possible case $d=2$, 
in usefulness as well as in
the standard order on $\mathbb{R}$. Theorem \ref{KeylWerner}
is thus the  better result. 

\section{Collapse of the Wave function}\label{collapse}
Recall that the Heisenberg principle states that
if all states are left intact, then no information can be
extracted from a system.
The `collapse of the wave function' may be seen as the flip side
of the Heisenberg principle. It states that if information is
extracted from a system, then its states undergo a very specific
kind of perturbation, called decoherence.
\subsection{Collapse for Unbiased Information Transfer}
We start out by investigating unbiased information transfer. 
We prove a sharp upper bound on the amount of remaining coherence
in terms of the measurement quality. 
\begin{theorem} \label{offdiag}
Let $T : \CA\ten\CB \rightarrow \CA$ be a CP-map. Let $B \in \CB$ 
be Hermitean, and 
consider $T$ as an unbiased information transfer from 
$A = T(\one \otimes B)$ to $\one \otimes B$
with maximal added variance 
$\Sigma^2 = \|(\one \otimes B,\one \otimes B)\|$.
Suppose that $\psi_x$ and $\psi_y$ are unit eigenvectors of $A$ with 
different eigenvalues 
$x$ and $y$ respectively. Define, 
for $|\alpha|^2 + |\beta|^2 = 1$, the coherent and the incoherent state\vspace{-0.5mm}
$$\rho_{\mathrm{coh}} := \ketbra{\alpha \psi_x + \beta \psi_y}
,\,\,\,\,\,\,\,\,\,
\rho_{\mathrm{inc}} := |\alpha|^2 \ketbra{\psi_x} + |\beta|^2 \ketbra{\psi_y}\,.\vspace{-0.5mm}
$$
If $R: \CA \rightarrow \CA$ is the restriction of 
$T$ to $\CA \ten \one$, then we have the following bound in trace distance:\vspace{-1mm}
\begin{equation}\label{karper}
d \Big(
R^* \big( \rho_{\mathrm{coh}} \big) , 
R^*\big( \rho_{\mathrm{inc}} \big) 
\Big) \leq
\frac{ \Sigma / |x-y| }
{\sqrt{1 + 4 \left( \Sigma / |x-y| \right)^2}} \,.\vspace{-1mm}
\end{equation}  
This bound is sharp in the sense that for all values of $\Sigma/|x-y|$,
there exists a CP-map $T$ 
for which (\ref{karper}) attains equality.\vspace{-6mm}
\end{theorem}
\begin{center}
\begin{tabular}{p{8cm}}
\includegraphics[width = 6 cm, viewport =  0 0 500 500]{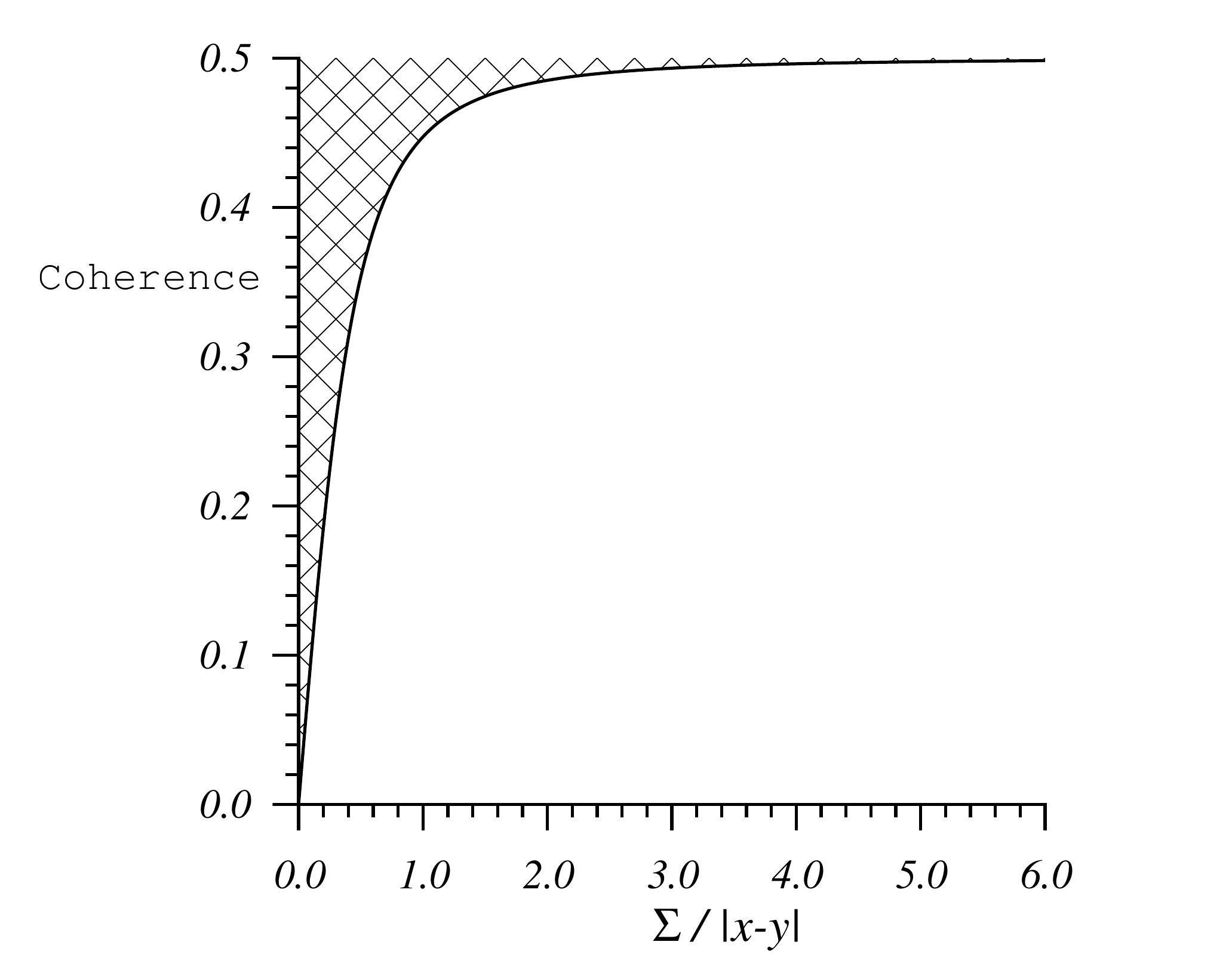}\\[-2mm]
{\small \fig Bound on the coherence as a function of $\Sigma/|x-y|$.
Only points below the curve are allowed.}
\end{tabular}
\end{center}
Consider the ideal case of perfect ($\Sigma = 0$) 
information transfer.
Suppose that the system
$\CA$ is initially 
in the coherent state 
$\rho_{\mathrm{coh}} = \ketbra{\alpha \psi_x + \beta \psi_y}$.
Then 
theorem \ref{offdiag} says that, after the information 
transfer 
to the ancilla $\CB$, the system $\CA$
cannot be distinguished from one that started out
in the incoherent 
state 
$\rho_{\mathrm{inc}} = |\alpha|^2 \ketbra{\psi_x} + |\beta|^2 \ketbra{\psi_y}$
instead.
As far as the behaviour of $\CA$ \emph{after} the information transfer
is concerned, it is therefore 
completely harmless to assume that
a collapse 
$$
\ketbra{\alpha \psi_x + \beta \psi_y}
\mapsto
|\alpha|^2 \ketbra{\psi_x} + |\beta|^2 \ketbra{\psi_y}
$$
will have occurred \emph{at the start} of the procedure. 

\proof 
The l.h.s. of equation (\ref{karper}) can be written
$$ 
\sup \{
\,
\bar{\alpha}\beta \inp{\psi_x}{R(P) \psi_y} + \mathrm{{c.c.}}
\, \,; \,\,
P \in \CA, \, P^{\dagger} = P^2 = P\,\}\,.
$$   
Furthermore, since $ 2|\alpha||\beta| \leq 1$,
it suffices to bound the `coherence' 
$\inp{\psi_x}{R(P) \psi_y}$ on all projections $P$. 
Now  
$(x-y) \inp{\psi_x}{R(P) \psi_y} =
\inp{\psi_x}{[A,R(P)] \psi_y}$, and furthermore
$[A,R(P)] = 
(P \ten \one, \one \ten B) - 
(\one \ten B, P \ten \one)$. Therefore
\begin{equation} \label{garbok}
(x-y) \inp{\psi_x}{R(P) \psi_y} =
\inp{\psi_x}{(P \ten \one, \one \ten B) \psi_y} - 
\inp{\psi_x}{(\one \ten B, P \ten \one) \psi_y}\,,
\end{equation}
and we will bound these last two terms.

By Stinespring's theorem \cite{St, Ta}, we may assume 
without loss of generality that
$T$ is of the form $T(X) = V^\dag X V $ for some contraction $V$.
We then obtain 
$(X,Y) = \smash{V^\dag X^\dag (\one - VV^\dag)} YV$.
Defining $g(X) := \sqrt{\one - VV^\dag} X V$, we have
$(X,Y) = g(X)^\dag g(Y)$, and therefore
$\| g(\one \ten B) \| = \Sigma$. 
%
%
This leads to the estimate\vspace{-1mm}
\begin{eqnarray*}
\inp{\psi_x}{(P \ten \one, \one \ten B) \psi_y}
&=&
\inp{g(P \ten \one) \psi_x}{g(\one \ten B)\psi_y}\\
&\leq&
\|g(P \ten \one) \psi_x\| \|g(\one \ten B)\psi_y\| \\
&\leq&
\Sigma \sqrt{\inp{\psi_x}{(P\ten\one , P\ten\one) \psi_x}}\,.
\end{eqnarray*}
We will bound 
$
\inp{\psi_x}{ (T(P^2 \ten \one) - T(P\ten \one)^2) \psi_x } =
\inp{\psi_x}{ (R(P) - R(P)^2) \psi_x }
$
in terms of the coherence.
For brevity, let us denote $\inp{\psi_x}{X \psi_{x'}}$
by $X_{xx'}$.
Since $\psi_x \perp \psi_y$,  we have 
$
(R(P)^2)_{xx} \geq |R(P)_{xx}|^2 + |R(P)_{xy}|^2
$, so that 
$$ 
(R(P) - R(P)^2)_{xx}
\leq
R(P)_{xx}(1- R(P)_{xx}) - |R(P)_{xy}|^2\,.
$$
Since $x(1-x) \leq \smash{\fourth}$ for all $x \in \mathbb{R}$,
this is at most
$
\smash{\fourth} - |R(P)_{xy}|^2
$. 

All in all, we have obtained
$$
(P \ten \one, \one \ten B)_{xy} \leq 
\Sigma \sqrt{ \fourth - |R(P)_{xy}|^2 }\,,
$$ and of course the same for $x \leftrightarrow y$.
Plugging this into equation \ref{garbok} yields
$$
|x-y| |R(P)_{xy}| \leq 
2 \Sigma  \sqrt{\fourth - |R(P)_{xy}|^2}\,,
$$
or equivalently the desired result
$\smash{
|R(P)_{xy}| \leq 
\frac{ \Sigma / |x-y| }
{\sqrt{1 + 4 \left( \Sigma / |x-y| \right)^2}}
} 
$. For sharpness, see section \ref{sharp}.
\qed

\noindent We have seen that if $\Sigma = 0$, then theorem \ref{offdiag}
is precisely the collapse of the wave function.
Now consider the other extreme of a measurement which leaves all states intact,
i.e. $R^*\rho = \rho$ for all $\rho$. 
Then there exist states for which the l.h.s. of
equation (\ref{karper}) equals $\half$, forcing
$\Sigma \rightarrow \infty$; no information can be obtained. 
This is Werner's formulation of the Heisenberg principle.

Theorem \ref{offdiag} thus unifies the Heisenberg principle
and the collapse of the wave function. 
For $\Sigma = 0$ we have a full decoherence,
whereas if all states are left intact, we have $\Sigma \rightarrow \infty$. 
For all intermediate cases,
the bound \ref{karper}
on the remaining coherence
is an increasing function of 
$\Sigma / |x - y|$.

This agrees with physical intuition: 
decoherence between 
$\psi_x$ and $\psi_y$ is expected to occur in case the information
transfer is able to distinguish between the two. 
This is the case if the added variance is small w.r.t. the differences in mean.


\subsection{Application: Perfect Qubit Measurement}
In section \ref{trans}, we have encountered 
the von Neumann Qubit measurement.
Now consider any perfect measurement $T$ of $\s_z$
with pointer $\one \ten (\delta_+ - \delta_-)$ which
leaves $\up$ and $\down$ in place, that is
$R^*(\up) = \up$ and $R^*(\down) = \down$.
(Such a measurement is often called \emph{nondestructive}.) 
Theorem \ref{offdiag} then reads 
$R^*(\ketbra{\alpha \uparrow + \beta \downarrow}) = 
|\alpha|^2 \ket{\!\uparrow\,}\bra{\,\uparrow\!} + 
|\beta|^2 \ket{\!\downarrow\,}\bra{\,\downarrow\!}$, as
illustrated below.\vspace{-2mm}
\begin{center}
\begin{tabular}{c}
\includegraphics[width = 7.3 cm]{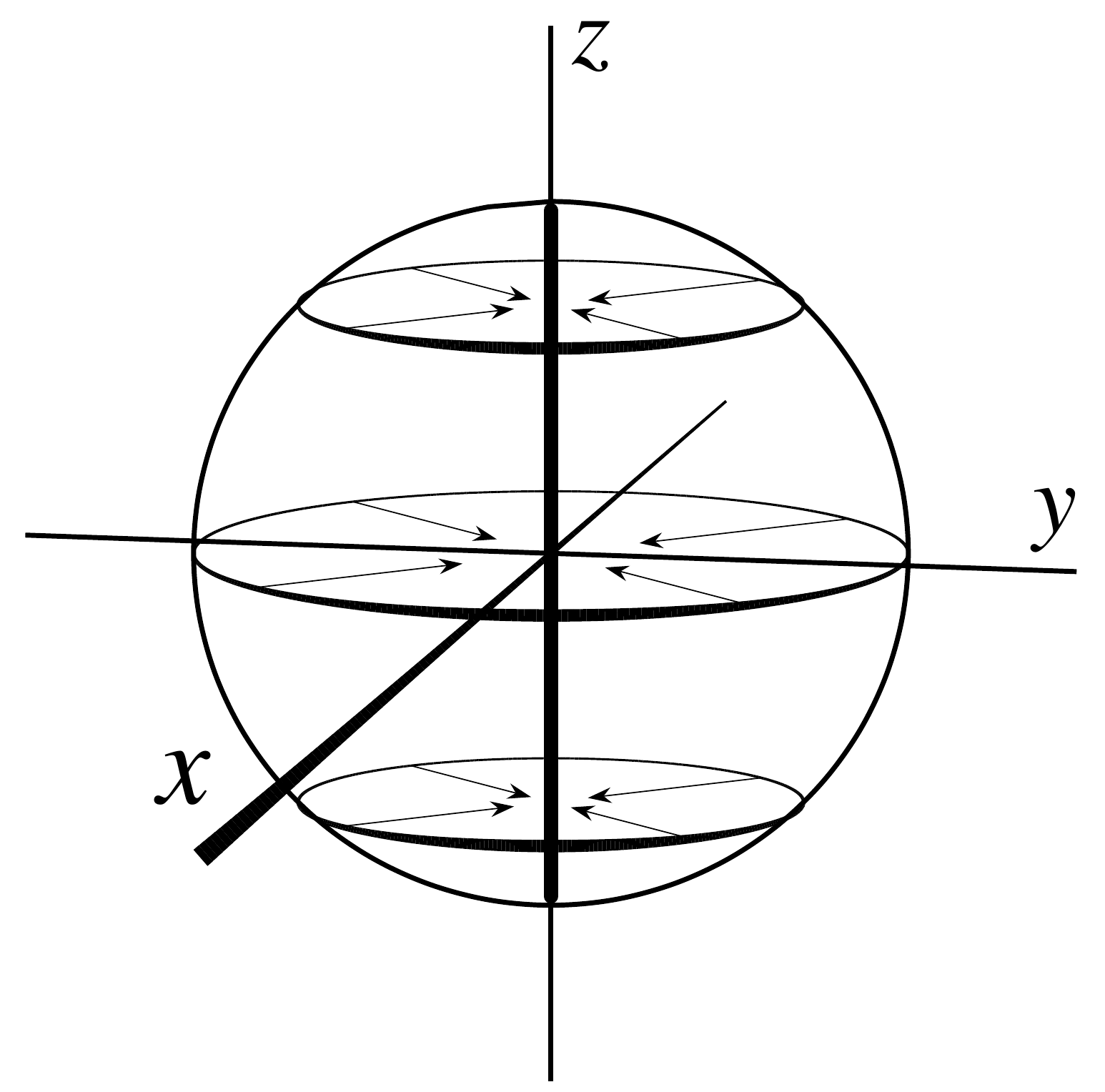}\\[1mm]
{\small \fig Collapse on the Bloch-sphere for perfect measurement.}\\ 
\end{tabular}
\end{center}
Incidentally, the trace distance between the centre of the Bloch sphere and
its surface is $\half$, so that we read off
$\Delta = \sup \{ D(R^*(\rho),\rho) ; \rho \in \CS(M_2) \} = \half$.
This was predicted by theorem \ref{HP}.

\subsection{Collapse for Nondestructive Measurement }

We will prove a sharp bound on the remaining coherence 
in general information transfer. 
The word `general' in the previous sentence means that the measurement
need not be unbiased, but
we will focus attention on 
nondestructive measurements.

A measurement of $A$ is called `nondestructive' 
(or `conserving' or `quantum nondemolition')
if it leaves the eigenstates of $A$ intact,
so that repetition of the measurement will yield the same result.
For example, the measurement in section \ref{Neuqubit}
is nondestructive, the one in section \ref{Resflu} is destructive.

Restricting attention to the nondestructive case is certainly not uncommon  
in quantum measurement theory (see \cite{Per}). In all fairness though, we 
make the nondemolition assumption out of technical convenience, 
not because we believe it to be of fundamental importance.

\begin{corollary}\label{gencol}
Let $T : B(\CH) \otimes \CB \rightarrow B(\CH)$ be a CP-map, and let 
$A \in B(\CH)$ and
$B \in \CB$ be Hermitean. 
Suppose that $A$ has an orthonormal basis 
of eigenvectors $\psi_i$ with eigenvalue $a_i$.
%
If $a_i \neq a_j$, and $|\alpha|^2 + |\beta|^2 = 1$, define
the coherent and the incoherent state
$$\rho_{\mathrm{coh}} := \ketbra{\alpha \psi_i + \beta \psi_j}
,\,\,\,\,\,\,\,\,\,
\rho_{\mathrm{inc}} := |\alpha|^2 \ketbra{\psi_i} + |\beta|^2 \ketbra{\psi_j}\,.
$$
Suppose that $T$ is nondestructive, i.e.
$R^*(\ketbra{\psi_{i}}) = \ketbra{\psi_{i}}$ for all $\psi_{i}$,
with $R$ the restriction of $T$ to $B(\CH) \ten \one$.
Then if 
the measurement infidelity
$\delta := \sup_{S} \{ \|T(\one \ten \one_{S}(B)) - \one_{S}(A)\| \}$
satisfies
$\delta \in [0 , \half]$, then
\begin{equation}\label{bound3}
d \Big(
R^* \big( \rho_{\mathrm{coh}} \big) , 
\rho_{\mathrm{inc}} 
\Big) 
\leq \sqrt{\delta(1-\delta)}\,.
\end{equation}
This bound is sharp in the sense that for all $\delta \in [0,\half]$,
there exists a CP-map $T$ 
for which (\ref{bound3}) attains equality.\vspace{-2mm}
\end{corollary}
\begin{center}
\begin{tabular}{p{7cm}}
\includegraphics[width = 6.5 cm]{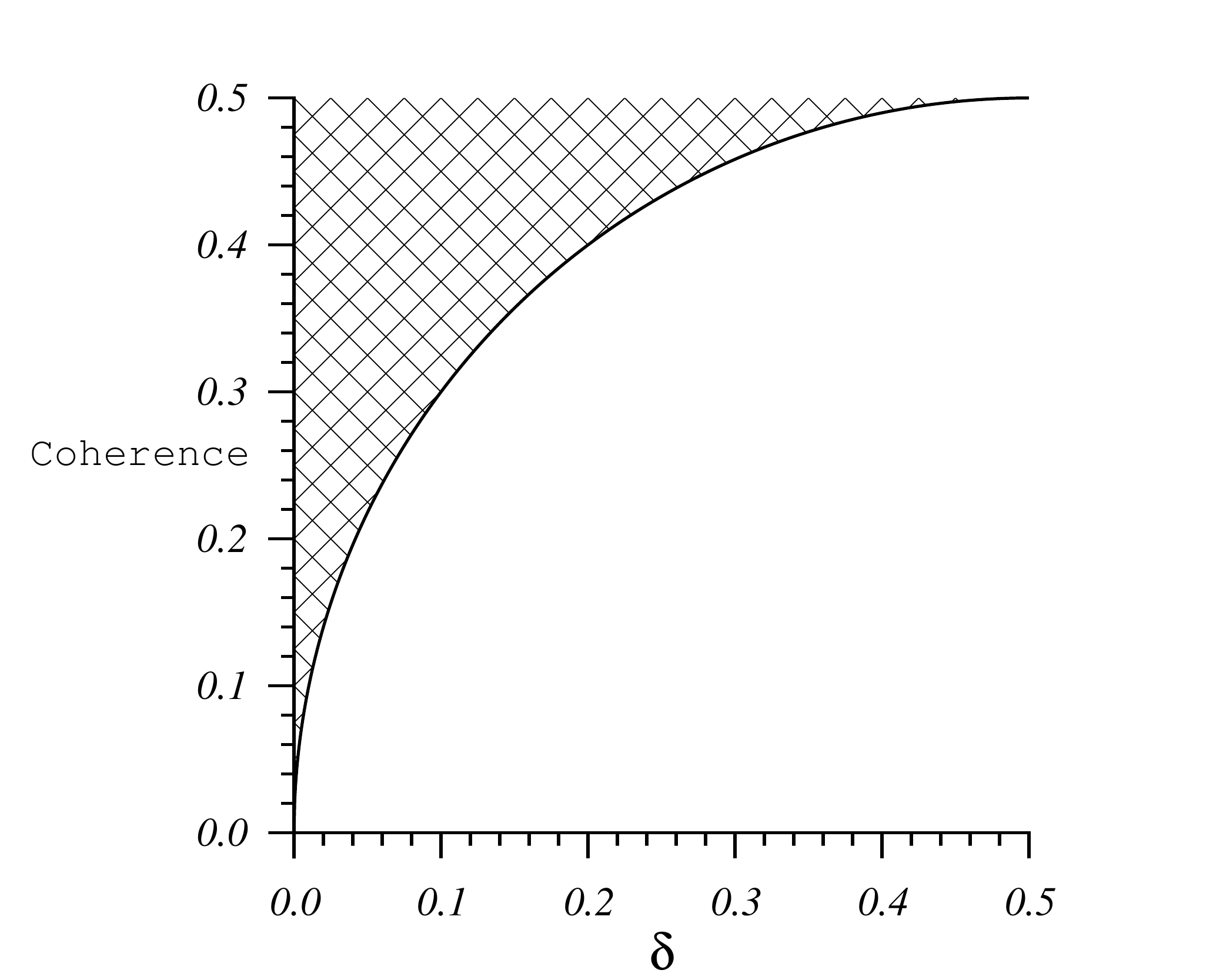}\\[-2 mm]
{\small \fig \label{figr}
Bound on the coherence in terms of $\delta$.
Only points below the curve are allowed.}
\end{tabular}
\end{center}
\proof
For sharpness, see section \ref{sharp}. 
Choose a set $S$ such that $a_i \in S$ and $a_j \notin S$.
$T$ is an unbiased measurement of $T(\one \ten \one_{S}(B))$
with pointer $\one_{S}(B)$ and maximal added variance
$\Sigma^2 \leq \delta(1 - \delta)$ (cf. the proof of corollary
\ref{hpdel}).
We will prove that $\psi_{i}$ and $\psi_{j}$ are
eigenvectors of $T(\one \ten \one_{S}(B))$ with  
eigenvalues $x$ and $y$ which differ by at least $1 - 2\delta$.

Define $P_{i} := \ketbra{\psi_{i}}$. Since $T$ is nondestructive, 
we have
$\inp{\psi_{j}}{ R(P_{i}) \psi_{j}} = 
\inp{\psi_{j}}{ P_{i} \, \psi_{j}}$
for all $j$.
Apparently, $R(P_{i})$ has only one nonzero diagonal element, a $1$ 
at position $(i,i)$. Since $R(P_{i}) \geq 0$, this implies
$R(P_{i}) = P_{i}$. 

From 
$(P_{i} \ten \one , P_{i} \ten \one ) = 0$,
we obtain
$(P_{i} \ten \one , \one \ten \one_{S}(B)) = 0$
by Cauchy-Schwarz.
Since $P_{i} = T( P_{i} \ten \one )$, we have
$
[T(\one \ten \one_{S}(B)) , P_{i} ] =
(P_{i} \ten \one , \one \ten \one_{S}(B)) -
(\one \ten \one_{S}(B) , P_{i} \ten \one)
= 0
$. Therefore $\psi_{i}$ is an eigenvector of
$T(\one \ten \one_{S}(B))$, with eigenvalue $x$, say. 
By a similar reasoning, $\psi_{j}$ is also an eigenvector,
denote its eigenvalue by $y$.

Since $\| T(\one \ten \one_{S}(B)) - \one_S (A) \| \leq \delta$,
we have in particular 
$\| (T(\one \ten \one_{S}(B)) - \one_S (A)) \psi_{i} \| = 
|x-1| \leq \delta$
and
$\| (T(\one \ten \one_{S}(B)) - \one_S (A)) \psi_{j} \| 
= |y| \leq \delta$, so that $|x - y| \geq 1 - 2\delta$.
We can now apply Theorem \ref{offdiag}. 

On the l.h.s.
of the bound (\ref{karper}), we may substitute
$
R^*(\rho_{\mathrm{inc}}) = \rho_{\mathrm{inc}}
$
on account of $T$ being nondestructive.
On the r.h.s., 
we substitute $\Sigma = \sqrt{\delta(1-\delta)}$ 
and $|x - y| = (1 - 2\delta)$. 
Strikingly enough, this yields the bound
$$\big( \sqrt{\delta(1-\delta)} / (1-2\delta) \big)   / 
\sqrt{1 + 4 
\big( \delta(1-\delta) / (1 - 2\delta)^2 \big)} = 
\sqrt{\delta(1-\delta)}\,\vspace{-6mm}
$$ 
\hfill $\Box$\par\medskip
\noindent For perfect measurement ($\delta = 0$), this yields
$R^* (\ket{\alpha \psi_{i} + \beta \psi_{j} \,} \bra{ \, \alpha \psi_{i} +
\beta \psi_{j} } )  =  |\alpha|^2 \ketbra{\psi_{i}} + 
|\beta|^2 \ketbra{\psi_{j}}$; all coherence between $\psi_i$ and $\psi_j$
must vanish. 
This collapse of the wave function is illustrated in the lower left corner of fig. \ref{figr}.
On the other hand, if all states are left intact so that 
$R^* = \id$, then we must have $\delta = \half$;
no information can be gained.
This is illustrated in the upper right corner of fig.\ref{figr}.
Corollary \ref{gencol} is a unified description of the Heisenberg principle
and the collapse of the wave function.

\subsection{Application: Nondestructive Qubit-Measurement}
In quantum information theory, a
$\s_z$-measurement is often taken to yield output $+1$
or $-1$, according to whether the input was 
$\ket{\!\uparrow\,}$ or $\ket{\!\downarrow\,}$.
It is nondestructive if it leaves the states
$\ket{\!\uparrow\,}$ and $\ket{\!\downarrow\,}$ intact,
yet it is only unbiased if it is perfect. 
Corollary \ref{gencol} shows that in the nondestructive case, 
the Bloch-sphere 
collapses to the cigar-shaped region depicted below:\vspace{-2mm}  
\begin{center}
\begin{tabular}{c}
\includegraphics[width = 7.5 cm]{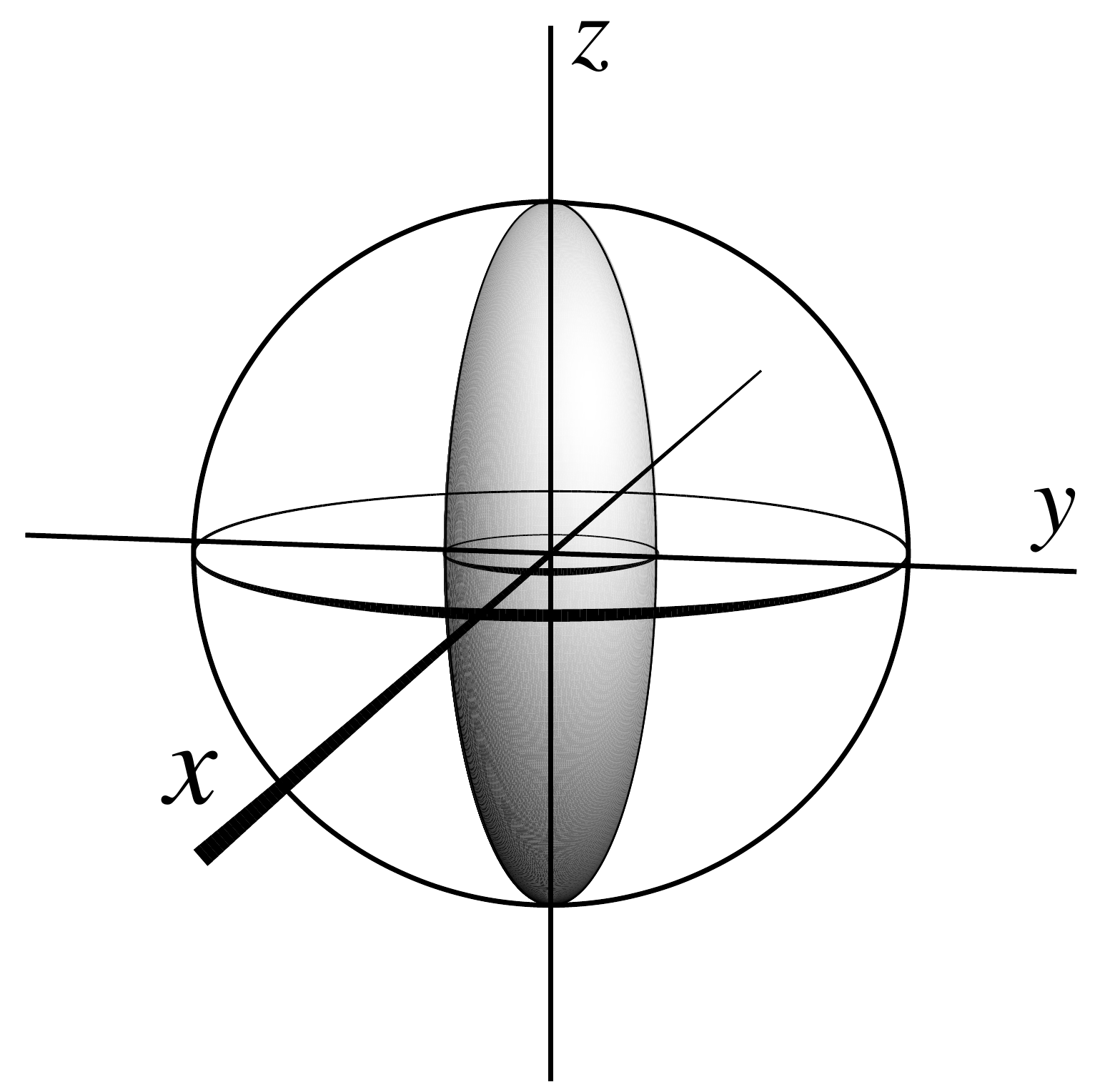}\\[2 mm]
{\small \fig Collapse on the Bloch-sphere with $\delta = 0.01$}.
\end{tabular}
\end{center}
Current single-qubit readout technology is just now moving into 
the regime where the bound (\ref{bound3}) becomes significant. 
in \cite{Lu}, a nondestructive measurement
of a SQUID-qubit was described, with experimentally determined 
measurement infidelity $\delta = 0.13$.
The bound then 
reads $d( R^* \rho_{\mathrm{coh}} , \rho_{\mathrm{inc}}) \leq 0.336$.

\subsection{Sharpness of the Bounds}\label{sharp}

We have yet to prove sharpness of the bounds. As it turns out, 
the following CP-map will do the trick for all the bounds
in a single stroke.  

Let $\Omega = \{+1,-1\}$, and $\CB = L^{\infty}(\Omega)$. Let
$$V_+ := 
\begin{pmatrix}
\sqrt{1-p} & 0 \\ 0 & \sqrt{p}
\end{pmatrix}
,\quad
V_- := \begin{pmatrix}
\sqrt{p} & 0 \\ 0 & \sqrt{1-p}
\end{pmatrix}\,,
$$ and define 
$T :  M_2(\mathbb{C}) \ten \CB \rightarrow M_2(\mathbb{C})$ by
$$
T(X \ten f) := f(+1) V_+ X V_+ + f(-1) V_- X V_-\,.
$$
For $p=0$, this is the von Neumann-measurement. 
As a measurement of $\s_z$
with pointer $B:= (\delta_+ - \delta_-)/(1-2p)$,
we have
$\delta = p$.
This yields bounds on the disturbance and on the coherence.
Corollary \ref{hpdel} and theorem \ref{HP} yield
$\Delta \geq \half - \sqrt{p(1-p)}$,
corollary \ref{gencol} and theorem \ref{offdiag} yield  
$
d (
R^* \big( \rho_{\mathrm{coh}} \big) , 
\rho_{\mathrm{inc}})
\leq \sqrt{p(1-p)}
$
We now explicitly calculate the restriction of $T$ to $M_2$,
and find 
$$R^*(\rho) = \begin{pmatrix} 
\rho_{11} & 2\sqrt{p(1-p)} \rho_{12} \\ 
2\sqrt{p(1-p)} \rho_{21} & \rho_{22}
\end{pmatrix}
\,.$$
The maximal remaining coherence occurs for $\alpha = \beta = 1/\sqrt{2}$,
for which it equals $\sqrt{p(1-p)}$. 
The maximal disturbance equals $\Delta = \half - \sqrt{p(1-p)}$.
This shows the bounds (\ref{bound}), (\ref{bound2}), (\ref{karper})
and (\ref{bound3}) to be sharp.

\section{Conclusion}
Our expedition into the heart of darkness has brought us the following
spoils.
\begin{itemize}
\item[I]
Theorem \ref{TJM} provides a sharp, state independent bound on 
the performance of unbiased joint measurement of noncommuting observables. 
In the case of perfect ($\Sigma=0$)
measurement of one observable, it implies that 
no information whatsoever 
($\Sigma' = \infty$) can be gained on the other. 
\item[II]
Theorem \ref{HP} (for unbiased information transfer)
and corollary \ref{hpdel} (for general information transfer) 
provide sharp, state independent bounds on the 
performance of a measurement in terms of the 
maximal disturbance that it causes.
In the case of zero disturbance, when all states are left intact,
it follows that no information can be obtained.
This is the Heisenberg principle.
\item[III]
Theorem \ref{KeylWerner}, based on a result of Keyl and
Werner, provides a sharp, state independent bound on the quality 
of classical coding. In contrast to I, II and IV, the result
based on the Cauchy-Schwarz inequality 
is not optimal in this case.  \newpage
\item[IV]
Theorem \ref{offdiag} (for unbiased information transfer) and 
corollary \ref{gencol} (for nondestructive information transfer) 
provide a sharp upper bound on the amount of coherence 
which can survive information transfer.
For perfect information transfer, 
all coherence vanishes.
This clearly proves that decoherence on a system is a mathematical
consequence of information transfer out of this system.
If, on the other hand, all states are left intact, then it follows 
that 
no information can be obtained.
This is the Heisenberg principle. 
Theorem \ref{offdiag} and 
corollary \ref{gencol} connect these two extremes 
in a continuous fashion; they form
a unified description of the Heisenberg
principle and the collapse of the wave function.
\end{itemize} 
All these restrictions follow directly from the formalism of 
quantum probability theory as set forth in chapter \ref{ch:intro}.

\chapter{Optimal Pointers for Joint Measurement 
}
\label{ch:OPJM}
\setsubdir{hoofdstukken/hoofdstuk4/}
%
%
%
%
Having concluded our investigation into the theoretical restrictions
of quantum measurement in general, we turn our attention towards
a particular example of a joint measurement procedure.
In this chapter, we 
study a qubit in 
interaction with the electromagnetic field. 
By means of homodyne detection, 
the field-quadrature $A_t+A_t^{\dagger}$ is 
observed continuously in time.
Due to the interaction, 
information about the initial state of the qubit is transferred into 
the field, thus influencing the homodyne measurement results.  
We construct random variables
(pointers) on the probability space of homodyne measurement outcomes having 
distributions close to the initial distributions of $\sigma_x$ and 
$\sigma_z$. Using variational calculus, we find the pointers that 
are optimal. 
These optimal pointers are very close to 
hitting the bound imposed by the joint measurement inequality,
theorem \ref{TJM} in chapter \ref{ch:UDHP}.

\section{Introduction}\label{sec intro}

The implementation of quantum filtering and 
control \cite{Bel88} in recent experiments \cite{AASDM02}, \cite{GSM04} 
has brought new interest to the field of continuous time 
measurement of quantum systems \cite{Dav69}, \cite{Dav76}, \cite{Bel88}, 
\cite{Bel92b}, \cite{Car93}, \cite{WiM93a}, \cite{Bo}.
In particular, homodyne detection has played a considerable role in this
development \cite{Car93}. 
In this chapter, we aim to gain insight into the transfer of information 
about the initial state of a qubit from this qubit, a two-level atom, 
to the homodyne photocurrent, which is observed in actual experiments. 
Our goal is to perform 
a joint measurement of two noncommuting observables in the 
initial system. In order to achieve this, we construct random variables 
(pointers) on the space of possible homodyne measurement results, having 
distributions close (in a sense to be defined) to the 
distributions of these observables in the initial state.  

\pagebreak[4]
The problem of joint measurement of noncommuting 
observables has been studied by several authors before, see 
\cite{Ne}, \cite{Dav76}, \cite{Ho} and the references 
therein. As a measure for the quality of an unbiased measurement, we use the 
\emph{maximal added variance}, i.e. the 
difference between the variance of the pointer in the final state and the 
variance of the observable in the initial state, evaluated in  
the worst case initial state. (See page \pageref{maximaaa}.)
%
Recall that
according to theorem 
\ref{TJM} in chapter \ref{ch:UDHP},
these maximal added variances for two 
pointers, corresponding to two noncommuting observables 
of the initial system, satisfy the Heisenberg-like relation (\ref{jointje}) 
that bounds how well their joint measurement can be 
performed. 

The present chapter concentrates on the example of a qubit 
coupled to the quantized electromagnetic field. We study 
this system in the weak coupling limit \cite{Gou05}, i.e.\ 
the interaction between qubit and field is governed by a 
quantum stochastic differential equation in the sense of 
Hudson and Par\-tha\-sa\-ra\-thy \cite{HuP84}. In the electromagnetic field 
we perform a homodyne detection experiment. Its integrated 
photocurrent is the measurement result for measurement of 
the field-quadrature $A_t+A_t^{\dagger}$ continuously in time. 
Using the characteristic functions introduced by Barchielli 
and Lupieri \cite{BaL85}, we find the probability density 
for these measurement results. In this density the $x$- and 
$z$-component of the Bloch vector of the initial state appear, 
indicating that homodyne detection is in fact a joint measurement 
of $\sigma_x$ and $\sigma_z$ in the initial state. 

Our goal is to construct random variables (pointers) 
on the probability space of homodyne measurement results 
having distributions as close as possible  
to those of the ob\-ser\-va\-bles $\sigma_x$ and $\sigma_z$ in 
the initial state of the qubit.
`As close as possible' is taken to mean that the pointer must give an
unbiased estimate of the observable, with its maximal added variance
as low as possible.
Using an argument due to Wiseman \cite{Wis96},
we first show that optimal random variables will only depend 
on the endpoint of a weighted path of the integrated 
photocurrent. Allowed to restrict our attention to this smaller class of
pointers, we are able to use standard variational calculus to obtain the 
optimal random variables. They do not achieve the bound 
imposed by the joint measurement theorem \ref{TJM}, 
but will be off by less than 5.6\%. 

The remainder of the chapter is organized as follows. In 
section \ref{sec model} we introduce the model of the 
qubit coupled to the field in the weak coupling limit. 
In section \ref{sec quality}, we briefly recall the definition of maximal 
added variance,
as well as the Heisenberg-like relation (\ref{jointje}) 
for joint measurement. In section \ref{sec characteristic function}
we calculate the characteristic function of Barchielli 
and Lupieri for the homodyne detection experiment. 
Section \ref{sec variation} deals with the variational 
calculus to find the optimal pointers. In section 
\ref{sec densities} we calculate the densities of the 
optimal pointers and then capture our main results 
graphically. 
In the last section we discuss our 
results.

\section{The Model}\label{sec model}

We consider a two-level atom, i.e. a qubit, in interaction with 
the quantized electromagnetic field. The qubit is described 
by $\BB{C}^2$ and the electromagnetic field by the \emph{symmetric 
Fock space} $\CF$ over the Hilbert space of quadratically integrable 
functions $L^2(\BB{R})$ (space of one-photon wave functions), i.e.\vspace{-1mm}
  \begin{equation*}
  \CF := \BB{C} \oplus \bigoplus_{k = 1}^\infty L^2(\BB{R})^{\ten_s k}.
  \end{equation*}
With the Fock space $\CF$ we can describe superpositions of field-states 
with different numbers of photons. The joint 
system of qubit and field together is described by the 
Hilbert space $\BB{C}^2\ten\CF$.

The interaction between the qubit and the electromagnetic field is 
studied in the weak coupling limit \cite{Gou04}, \cite{Gou05}, \cite{AFL}. 
This means that in the interaction picture the unitary dynamics of 
the qubit and the field together is given by a quantum 
stochastic differential equation (QSDE) in the sense of Hudson 
and Parthasarathy \cite{HuP84} 
  \begin{equation}\label{eq QSDE}
  dU_t =\Big\{\sigma_{-} dA_t^{\dagger} - \sigma_{+} dA_t - 
  \half\sigma_{+} \sigma_{-} dt\Big\}U_t,\vspace{-1mm}
  \end{equation}
  with\vspace{-1mm}
  $$
  \sigma_{-}  = \begin{pmatrix}0 & 0 \\ 1 & 0 \end{pmatrix}, \ 
  \sigma_{+}  = \begin{pmatrix}0 & 1 \\ 0 & 0 \end{pmatrix}, \
  $$
  and $ 
  U_0 = \one$.   
The operators $\sigma_{-}$ and $\sigma_{+}$ are the annihilator and
creator on the two-level system. 
The field annihilation and creation processes
are denoted $A_t$ and $A^{\dagger}_t$, respectively. 
Keep in mind that the evolution
$U_{t}$ acts nontrivially on the combined system $\BB{C}^2\ten\CF$, 
whereas $\sigma_{\pm}$ and $A_{t}$ are understood to designate
the single system-operators $\sigma_{\pm}\ten \one$ and $\one\ten A_{t}$.
Throughout
this chapter we will remain in the interaction picture. Equation \eqref{eq QSDE} 
should be understood as a shorthand for the integral equation
  \begin{equation*}
  U_t = \one + \int_0^t\sigma_{-} U_{\tau} dA^{\dagger}_{\tau} 
   - \int_0^t\sigma_{+} U_{\tau}dA_{\tau} 
   - \half\int_0^t\sigma_{+} \sigma_{-} U_{\tau}d\tau, 
  \end{equation*} 
where the integrals on the right-hand side are 
stochastic integrals in the sense of Hudson and Parthasarathy \cite{HuP84},
and Picard iterations can be used to show that (\ref{eq QSDE})
has a unique solution.
The value of these integrals does not lie so much in 
their actual definition (see section \ref{qsdesec}), but in 
the It\^o rule satisfied by them, allowing for easy calculations. 
The following is a slight extension of theorem \ref{qiet}.

\begin{theorem}[Quantum It\^o rule \cite{HuP84}, \cite{Par92}]\label{Itorule}
Let $X_t$ and $Y_t$ be quantum stochastic integrals of the form 
\begin{equation*}\begin{split}
&dX_t = C_t dA_t + D_t dA_t^{\dagger} + E_tdt \\
&dY_t = F_t dA_t + G_t dA_t^{\dagger} + H_tdt
\end{split}\end{equation*}
for some stochastically integrable processes $C_t, D_t, E_t, F_t, 
G_t$ and $H_t$ (see \cite{HuP84}, \cite{Par92} for definitions).
Suppose that $X_t$ has an adjoint process, and that all
the necessary products exist and are stochastically integrable. 
Then the process $X_tY_t$ satisfies the relation
  \begin{equation*}
  d(X_tY_t) = X_tdY_t + (dX_t)Y_t + dX_tdY_t,
  \end{equation*}
where $dX_tdY_t$ should be evaluated according to 
the quantum It\^o table:
\begin{center}
{\large \begin{tabular} {l|lll}
 & $dA_t$  & $dA^{\dagger}_t$ & $dt$\\
\hline 
$dA_t$ & $0$ & $dt$ & $0$\\
$dA^{\dagger}_t$ & $0$  & $0$ & $0$\\
$dt$ & $0$ & $0$ & $0$
\end{tabular} }
\end{center}
i.e.\ $dX_tdY_t = C_tG_tdt$. 
\end{theorem}
As a corollary we have that, for any $f \in C^{2}(\BB{R})$, the 
process $f(X_{t})$ satisfies
$\textstyle d (f(X_{t})) = f'(X_{t})dX_{t} + \half f'' (X_{t}) (dX_{t})^2,$
where $(dX_t)^2$ should be evaluated according to the quantum It\^o table.


Let us return to equation \eqref{eq QSDE}. In order to illustrate how the 
quantum It\^o rule will be used,
we calculate the time evolution on the qubit explicitly. 
We choose an orthonormal basis of $\mathbb{C}^{2}$, and identify
the algebra $\CB(\mathbb{C}^{2})$ of qubit-observables with 
$M_2(\BB{C})$,
the algebra of 
$2 \times 2$-matrices.
The algebra of observables in the field is given by $\CB(\CF)$, the bounded operators
on $\CF$.
If $\id :\ 
M_2(\BB{C})\to M_2(\BB{C})$ is the identity map and $\phi :\ \CB(\CF)\to\BB{C}$ is 
the expectation with respect to the vacuum state 
$\Omega := 1 \oplus 0 \oplus 0 \oplus \ldots \in \CF$ 
(i.e.\ $\phi (Y) := \inp{\Omega}{Y\Omega}$),
then time evolution on the qubit $T_t :\ M_2(\BB{C}) \to M_2(\BB{C})$ 
is given by
$ T_t(X) := \id \ten\phi(U_t^{\dagger}X\ten \one U_t)$.
On the combined system, the full time evolution 
$j_t :\ M_2(\BB{C})\ten \CB(\CF) \to M_2(\BB{C})\ten \CB(\CF)$  
is given by $j_t(W) := U^{\dagger}_tWU_t$.
In a diagram this 
reads 
   \begin{equation}\label{dildiag}\begin{CD}
  M_2(\BB{C}) @>T_t>> M_2(\BB{C})              \\
   @V{  \mbox{\footnotesize id} \ten \one}VV  
   @AA{\mbox{\footnotesize id} \ten \phi}A      \\
   M_2(\BB{C}) \ten \CB(\CF) @>j_t>> M_2(\BB{C}) \ten \CB(\CF).            \\
  \end{CD}\end{equation}
In the Schr\"odinger picture the arrows would be reversed.  
A qubit-state $\rho$ would be extended with the vacuum 
to $\rho\ten\phi$, time evolution would change it into 
into $\rho \ten \phi (U_t^{\dagger} \,\cdot \, U_t)$, 
and in the last 
step the partial trace over the field would be taken,
resulting in the state $\rho \circ T_t$.

Using the It\^o rule we can derive a (matrix-valued) differential 
equation for $T_t(X)$ as follows. We start by 
using the It\^{o} rule in order to write
  \begin{equation}\label{eq T}\begin{split}
  dT_t(X) &= \id \ten\phi\big(d(U_t^{\dagger}X\ten \one U_t)\big) \\
  &= \id \ten\phi\big((dU_t^{\dagger})X\ten \one U_t + U_t^{\dagger}X \ten \one(dU_t) + 
  (dU_t^{\dagger})X\ten \one(dU_t)\big)\,. \\ 
  \end{split}\end{equation}
We then use the fact that stochastic integrals with respect to $dA_t$ and 
$dA^{\dagger}_t$ vanish with respect to the vacuum expectation, leaving us only with the 
$dt$ terms. 
The QSDE for $U^{\dagger}_t$ easily 
follows from \eqref{eq QSDE};
  \begin{equation*}
  dU^{\dagger}_t = U_t^{\dagger}\Big\{\sigma_{+} dA_t -\sigma_{-} 
  dA^{\dagger}_t - \half\sigma_{+} \sigma_{-} dt\Big\}, \ \ \ \ 
  U^{\dagger}_0 = \one.
  \end{equation*} 
The first term in equation \ref{eq T} can then be written  
\begin{eqnarray*}
\id \otimes \phi (dU_t^{\dagger} X \ten \one U_t) 
&=& 
\id \otimes \phi (U_t^{\dagger} 
\big\{	\sigma_{+} dA_t -\sigma_{-} dA^{\dagger}_t - 
	\half\sigma_{+} \sigma_{-} dt
\big\}
X \ten \one U_t)\\
&=&
\id \otimes \phi (U_t^{\dagger} 
\big\{	-\half\sigma_{+} \sigma_{-} dt
\big\}
X \ten \one U_t)\\
&=&
T_t(- \half \sigma_{+} \sigma_{-} X)dt\,.
\end{eqnarray*}
Similarly, the second term yields
$$
\id \otimes \phi (U_t^{\dagger} X \ten \one dU_t) 
=
T_t( -\half X \sigma_{+} \sigma_{-} )dt\,.
$$
For the third term, we use the quantum It\^{o} table to see that
$$
U_t^{\dagger} 
\big\{	\sigma_{+} dA_t -\sigma_{-} dA^{\dagger}_t - 
	\half\sigma_{+} \sigma_{-} dt
\big\}
X \otimes \one
\big\{	\sigma_{-} dA^{\dagger}_t -\sigma_{+} dA_t - 
	\half\sigma_{+} \sigma_{-} dt
\big\}
U_t
$$
equals
$
U_t^{\dagger} (\sigma_{+}X\sigma_{-}\otimes \one ) U_t \,dt 
$,
so that
$$
\id \otimes \phi(d U_{t}^{\dagger} \, X \otimes \one \, d U_t ) = 
T_t (\sigma_{+} X \sigma_{-}) dt \,.
$$
Combining these three terms, equation (\ref{eq T}) becomes
\begin{equation}\label{fruthulkje}
\frac{d}{dt} T_t(X) = T_t(L(X))\,,
\end{equation}
with $L$ is the Lindblad generator 
  \begin{equation*}
  L(X) := -\half (\sigma_{+} \sigma_{-} X +X\sigma_{+} \sigma_{-} ) + \sigma_{+} X \sigma_{-} . 
  \end{equation*} 
The matrix-valued differential equation \eqref{fruthulkje} with initial 
condition $T_0(X) = X$ is solved by 
$T_t(X) =  \exp(tL)(X)$, which is exactly the 
time evolution of a two-level system spontaneously decaying to 
the ground state, as it should be.
(Compare this to the discussion in section \ref{Resflu},
with zero field strength $\Omega$.)
Although the arguments above 
are completely standard (cf.\ \cite{HuP84}), they do 
illustrate nicely and briefly some of the techniques used also in following 
sections.

\section{Quality of Information Transfer}\label{sec quality}

Now suppose we do a homodyne detection experiment, enabling 
us to measure the observables $A_t + A_t^{\dagger}$ in the field 
continuously in time \cite{Bar90}.
If initially the qubit is in state 
$\rho$, then at time $t$ the qubit and field together are in  
a state $\rho^t$ on $M_2 (\BB{C}) \ten \CB(\CF)$ given by 
$\rho^t(W) := \rho\ten\phi(U^{\dagger}_tWU_t) = 
\rho\big(\id \ten\phi(U^{\dagger}_tW U_t)\big)$.
We find $\rho^{t}(\one \otimes (A_t + A_t^{\dagger}))$ by solving
the differential equation
  \begin{eqnarray}
  d\Big(\id \ten\phi(U^{\dagger}_t \one\ten (A_{t} + A^{\dagger}_{t}) U_t)\Big) 
  &=& 
  \id \ten\phi\Big(d(U^{\dagger}_t \one\ten (A_{t} + A^{\dagger}_{t}) U_t)\Big)
  \label{gter1} \\
  &=&  
  \id \ten\phi\Big(U_t^{\dagger}(\sigma_{-}  + \sigma_{+} )\ten \one U_t\Big)dt
  \label{gter2} \\
  &=& 
  \exp(t L)(\sigma_{-} +\sigma_{+} )dt\\ 
  &=& e^{-\frac{t}{2}}\sigma_xdt \,. 
  \end{eqnarray}
In order to go from (\ref{gter1}) to (\ref{gter2}), 
we used once again the fact that the ground state kills all terms with 
$dA_t$ and $dA_t^{\dagger}$.
Using the quantum It\^{o} table to isolate the terms with $dt$, we then expand
$d(U^{\dagger}_t \one\ten (A_{t} + A^{\dagger}_{t}) U_t)$
into six terms, three of first order and three of second order. 
The $dt$-terms coming from 
$dU^{\dagger}_t \one\ten (A_{t} + A^{\dagger}_{t}) U_t$, 
$U^{\dagger}_t \one\ten (A_{t} + A^{\dagger}_{t}) dU_t$
and 
$dU^{\dagger}_t \one\ten (A_{t} + A^{\dagger}_{t}) dU_t$
cancel each other out, and obviously
$U^{\dagger}_t \one\ten (dA_{t} + dA^{\dagger}_{t}) U_t$
lacks terms with $dt$.
We are left with
$dU^{\dagger}_t \one\ten (dA_{t} + dA^{\dagger}_{t}) U_t$
and
$U^{\dagger}_t \one\ten (dA_{t} + dA^{\dagger}_{t}) dU_t$,
which contribute
$U_t^{\dagger} ((\sigma_+ + \sigma_-)\otimes \one )U_t$.

This shows that regardless the initial state $\rho$ of the qubit,  
the expectation of $(A_{t} + A^{\dagger}_{t})$ in the final state $\rho^t$ 
will equal the 
expectation of 
$( 2 - 2e^{-\frac{t}{2}} )\sigma_x$ in the initial state $\rho$.

\subsection{Defining the Quality of Information Transfer}

The process at hand is thus a transfer of information about $\sigma_x$ 
to the `pointer' $(A_t + A^{\dagger}_t)$, which can be read off by means of 
homodyne
detection. We recall definition \ref{geenvoorkeur} of
\emph{unbiased measurement} in the particular case 
of a qubit measurement by means of the electromagnetic field. 


\begin{basdef}[Unbiased Measurement]
Let X be an observable of the qu\-bit, i.e.\ a self-adjoint 
element of $M_2(\BB{C})$, and let $Y$ be an observable of 
the field, i.e.\ a self-adjoint operator in (or affiliated to) $\CB(\CF)$. 
Then an \emph{unbiased measurement} 
$M$ of $X$ with \emph{pointer} $Y$ is by definition a completely 
positive map $M:\ \CB(\CF) \to M_2(\BB{C})$ such that $M(Y) = X$. 
\end{basdef}  

\noindent Needless to say, for each fixed point $t$ in time, the map 
$M\,\colon\, \CB(\CF) \to M_2(\BB{C})$ defined by 
$M(B) := \id \ten\phi(\smash{U^{\dagger}_t} \one\ten BU_t)$  
is a measurement of $\sigma_x$ 
with pointer 
$\smash{Y = (2-2e^{-\frac{t}{2}})^{-1}(A_t + A_t^{\dagger})}$. This means that, 
after the measurement procedure of coupling to the field in the 
vacuum state and allowing for interaction with the qubit for $t$ 
time units, the distribution of the measurement 
results of the pointer $Y$ has inherited the expectation of 
$\sigma_x$, regardless of the initial state $\rho$. However, 
we are more ambitious and would like its distribution as a whole 
to resemble that of $\sigma_x$. 
We recall from section \ref{mav} that the performance 
of an unbiased measurement
is gauged by its \emph{maximal added variance} $\Sigma^2$.
We will call its square root $\Sigma$ the \emph{quality}
of measurement.
\begin{basdef}[Quality] 
Let $M:\ \CB(\CF) \to M_2(\BB{C})$ be an unbiased measurement of 
$X$ with pointer $Y$. Then its \emph{quality} $\Sigma$ is defined by
  \begin{equation*}
  \Sigma^2 := \sup \Big\{\var_{\rho\circ \hspace{-0.2 mm} M}(Y) - 
  \var_\rho(X)\big|\ 
  \rho \in \mathcal{S}(M_2 (\BB{C}))\Big\},
  \end{equation*}
where $\mathcal{S}(M_2 (\BB{C}))$ denotes the state space of $M_2(\BB{C})$ (i.e.\ all 
positive normalized linear functionals \mbox{on $M_2(\BB{C})$}).
\end{basdef}

\noindent This means that $\Sigma^2$ is the variance added to the initial 
distribution of $X$ by the measurement procedure $M$ for the worst 
case initial state $\rho$. As we have seen,
  \begin{equation*}
  \var_{\rho\circ \hspace{-0.2 mm} M}(Y) - 
  \var_\rho(X) = \rho\big(M(Y^2)-M(Y)^2\big), 
  \end{equation*}
which implies that $\Sigma^2 = \| M(Y^2) -  M(Y)^2 \| $, where 
$X \mapsto \| X \|$
denotes the operator norm \mbox{on $M_2(\BB{C})$}. In particular this shows 
that $\Sigma^2$ is positive, as one might expect. It follows from 
theorem \ref{struc} that
$\Sigma$ equals zero if and only if the measurement procedure 
$M$ exactly carries over 
the distribution of $X$ to $Y$. In short, $\Sigma$ is a suitable measure for 
how well $M$ transfers information about $X$ to the pointer $Y$. 

\subsection{Calculating the Quality of Information Transfer}

Let us return to the example at hand, i.e.\ 
$M(B) = \id \ten\phi(U_t^{\dagger}\one\ten B U_t)$, with field-observable $Y = (2-2e^{-\frac{t}{2}})^{-1}
(A_t + A_t^{\dagger})$ 
as a pointer for $\sigma_x$. Let us calculate its quality, which amounts to
evaluating $M(Y^2) = 
(2-2e^{-\frac{t}{2}})^{-2}M\big((A_t+A_t^{\dagger})^2\big)$. 
To this aim,
we will first introduce some ideas which will be 
of use to us in later calculations as well. 

\begin{basdef}\label{def F} 
Let $f$ and $h$ be real valued functions, $h$ twice 
differentiable. Let $Y_t$ be given by $dY_t = f(t)(dA_t + dA_t^{\dagger})$, 
$Y_0 = 0$. For $X \in M_2(\BB{C})$ we define 
 \begin{equation*}
 F_h(X,t) := \id \ten\phi\big(U_t^{\dagger} X \ten h(Y_t)U_t\big). 
 \end{equation*}
When no confusion can arise we shall shorten $F_h(X,t)$ to $F_h(X)$.
\end{basdef}

\noindent The homodyne detection experiment has given us a measurement 
result (the integrated photocurrent) which is just the path of 
measurement results of $A_t + A_t^{\dagger}$ continuously in time. Given 
this result, we post-process it by weighting the increments of the 
path with the function $f(t)$ and letting 
$h(y)$ act on the result. 
The following lemma will considerably shorten calculations.

\begin{lemma}\label{lem eq F}
  \begin{equation*}
  \frac{dF_h(X)}{dt} = F_h\big(L(X)\big) + f(t)F_{h'}(\sigma_{+} X +
  X\sigma_{-} ) + \half f(t)^2F_{h''}(X)
  \end{equation*}
\end{lemma} 
\proof
Just as we did in equation (\ref{gter1}),
we use that the vacuum expectation kills all terms with 
$dA_t$ and $dA^*_t$. We then expand 
$d(U_t^{\dagger} X \ten h(Y_t)U_t)$ into three first order 
and three second order terms.
Using theorem \ref{Itorule} we see that after the vacuum expectation, 
the terms $(dU_{t}^{\dagger})X \otimes h(Y_t) U_t$, 
$U_t^{\dagger} X \ten h(Y_t)dU_t$ and 
$(dU_t^{\dagger}) X \ten h(Y_t)dU_t$ 
make up $F_h\big(L(X)\big)dt$. 
From 
  \begin{equation*}
  dh(Y_t) = h'(Y_t)f(t)(dA_t + dA^{\dagger}_t) + \half h''(Y_t)f(t)^2 dt,
  \end{equation*}
we find that, after taking vacuum expectations, the terms 
$d(U_t^{\dagger}) d(X \ten h(Y_t))U_t$ and 
$U_t^{\dagger} d(X \ten h(Y_t))dU_t$ 
make up
the second term 
$f(t)F_{h'}(\sigma_{+} X + X\sigma_{-} )dt$, and 
$U_t^{\dagger} d(X \ten h(Y_t))U_t$ provides the 
last term $\frac{1}{2}f(t)^2F_{h''}(X)dt$.       
%
\qed

\noindent We are now well-equipped to calculate $M\big((A_t+A_t^{\dagger})^2\big)$.
Choose $f(t) = 1$ and $h(x) = x^2$. (The maps $x \mapsto x^n$ will be denoted $\x^n$
hereafter.) Then we have
$M\big((A_t+A_t^{\dagger})^2\big) = F_{\x^2}(\one)$ and by lemma \ref{lem eq F}
  \begin{equation}\label{eq Fx2}
  \frac{dF_{\x^2}(\one)}{dt} = 2 F_{\x}(\sigma_{-} +\sigma_{+} ) + F_{1}(\one) = 
  2 F_{\x}(\sigma_{-} +\sigma_{+} ) + \one, \ \ \ \ 
  F_{\x^2}(\one, 0) = 0. 
  \end{equation}
Applying lemma \ref{lem eq F} to 
$F_{\x}(\sigma_{-} +\sigma_{+} )$, we obtain  
  \begin{equation}\label{eq Fx}
  \frac{dF_{\x}(\sigma_{-} +\sigma_{+} )}{dt} = -\half F_{\x}(\sigma_{-} + \sigma_{+} ) + 
  2F_{1}(\sigma_{+} \sigma_{-} ), \ \ \ \
  F_{\x}(\sigma_{-} +\sigma_{+} , 0) = 0.
  \end{equation}
Finally, $F_{1}(\sigma_{+} \sigma_{-} )$ satisfies
  \begin{equation} \label{eq Fx0}
  \frac{dF_{1}(\sigma_{+} \sigma_{-} )}{dt} = - F_{1}(\sigma_{+} \sigma_{-} ),\ \ \ \ 
  F_{1}(\sigma_{+} \sigma_{-} , 0) = \sigma_{+} \sigma_{-}.   
  \end{equation}
Solving \eqref{eq Fx0}, \eqref{eq Fx} and \eqref{eq Fx2} successively leads
first to  
$F_1(\sigma_{+} \sigma_{-} ) = e^{-t}\sigma_{+} \sigma_{-} $, then to
$F_{\x}(\sigma_{-} +\sigma_{+} ) = 4(e^{-\frac{t}{2}} - e^{-t})\sigma_{+} \sigma_{-} $
and finally to
$F_{\x^2}(\one) = 8(e^{-\frac{t}{2}} -1)^2 \sigma_{+} \sigma_{-}  + tI$. 
Consequently, the quality of the measurement $M$ of $\sigma_x$ 
with pointer $Y = (2-2e^{-\frac{t}{2}})^{-1}(A_t+ A_t^{\dagger})$ 
is given by
  \begin{eqnarray*} 
  \Sigma^2 & = &\| M(Y^2)- M(Y)^2 \| \\
  &=& 
  \left\| \frac{8(e^{-\frac{t}{2}} -1)^2 \sigma_{+} \sigma_{-}  + tI}
  {(2-2e^{-\frac{t}{2}})^{2}} - \one \right\| \\
  &=& 
  \left\| 2\sigma_{+} \sigma_{-}  + 
  \left(\frac{t}{(2-2e^{-\frac{t}{2}})^{2}} - 1\right)\one \right\| \\
  & =& \frac{t}{(2-2e^{-\frac{t}{2}})^{2}} + 1.
  \end{eqnarray*}

This expression takes its minimal value $2.228$ at $t = 2.513$, 
leading to a quality $\Sigma = 1.493$. 

The calculation above has an interesting side product.
The observable $M\big( (A_t + A_t^{\dagger})^2 \big)$ depends linearly on 
$\sigma_{z}$, indicating that in addition to information
on $\sigma_{x}$, also information on $\sigma_z$ in the initial
qubit-state ends up in the measurement outcome. 
Indeed, if we use as a pointer
  \begin{equation}\label{eq tilde Y}
  \tilde{Y} := \frac{(A_t+A_t^{\dagger})^2 -tI}{4 (e^{-\frac{t}{2}}-1)^2} - \one,  
  \end{equation}
then we have $M(\tilde{Y}) = \sigma_z$, so that 
$M$ is also 
a measurement of $\sigma_z$ with pointer~$\tilde{Y}$.   
  
Note that the pointers $Y$ and $\smash{\tilde{Y}}$ commute, i.e.\ measuring 
$A_t+\smash{A_t^{\dagger}}$ via the homodyne detection scheme is an indirect 
\emph{joint measurement} of $\sigma_x$ and $\sigma_z$. If we would 
also like to gain some information about $\sigma_y$, we could 
for example sweep the measured quadrature through $[0, 2\pi)$ in time 
by measuring $e^{i\omega t}A_t + e^{-i\omega t}A^{\dagger}_t$ instead. 
We will however restrict ourselves to continuous time measurement of $A_t + A_t^{\dagger}$,
as additional information on $\sigma_y$ would deteriorate the quality of 
$\sigma_x$- and/or $\sigma_z$-measurement. 
We formulate the joint measurement theorem \ref{jointje} in the particular case of
qubit measurement by means of an electromagnetic field.
%
\begin{theorem}[Joint Measurement]\label{thm joint}
Let $M:\ \CB(\CF) \to M_2(\BB{C})$ be an unbiased measurement 
of self-adjoint observables $X \in M_2(\BB{C})$ and $\tilde{X} \in M_2(\BB{C})$ 
with self-adjoint commuting pointers 
$Y$ and $\tilde{Y}$ in (or affiliated to) $\CB(\CF)$, respectively. Then for 
their corresponding qualities $\Sigma$ and $\tilde{\Sigma}$ 
the following relation holds
  \begin{equation*}
  2\Sigma\tilde{\Sigma} \ge \| [X, \tilde{X}]\|.
  \end{equation*}
\end{theorem}

\noindent Denote by $\tilde{\Sigma}$ the quality of the $\sigma_z$
measurement with the pointer $\tilde{Y}$ defined in \eqref{eq tilde Y}.
Since $[\sigma_x,\sigma_z] = -2i\sigma_y$, the qualities $\Sigma$ and 
$\tilde{\Sigma}$ (corresponding to the pointers $Y$ and $\tilde{Y}$, respectively)
satisfy the inequality 
  \begin{equation}\label{eq ineq} 
  \Sigma\tilde{\Sigma} \ge 1. 
  \end{equation}
Using similar techniques as before, that is recursively calculating 
$F_{\x^4}(\one)$ via lemma \ref{lem eq F}, we find 
  \begin{equation*}
  \tilde{\Sigma}^2 = \frac{t^2}{8(e^{-\frac{t}{2}}-1)^4} + 
  \frac{2t-4(e^{-\frac{t}{2}}-1)^2}{(e^{-\frac{t}{2}}-1)^2}.
  \end{equation*}
This expression takes its minimal value $8.836$ at $t = 2.513$. This 
leads to a quality $\tilde{\Sigma} = 2.973$, which means that 
$\Sigma\tilde{\Sigma} = 4.437$, i.e.\ we are far removed from hitting the 
bound $1$ in \eqref{eq ineq}. However, there is still some room for 
manoeuvring by post-processing 
of the homodyne measurement data.

\section{The Weighted Path}\label{sec characteristic function}

Let us presently return to our homodyne detection experiment. 
We observe $A_{\tau} +A_{\tau}^{\dagger}$ continuously in time, i.e.\ the result of our 
measurement is a path $\omega$ of measurement results $\omega_{\tau}$ (the 
photocurrent integrated up to time $\tau$) for $A_{\tau}+A_{\tau}^{\dagger}$.
This means that we have a space $\Omega$ of all possible measurement 
paths and that we can identify an operator $A_{\tau}+A_{\tau}^{\dagger}$ with the 
map from $\Omega$ to $\BB{R}$ mapping a measurement 
path $\omega \in \Omega$ to the measurement result $\omega_{\tau}$ at time $\tau$.     
That is, we have simultaneously diagonalized the family of commuting 
operators $\{A_{\tau}+A^{\dagger}_{\tau}|\ \tau \ge 0\}$ and viewed them as random 
variables on 
the spectrum $\Omega$. The spectral projectors of the operators 
$\{A_{\tau}+A_{\tau}^{\dagger}|\ 0\le \tau \le t\}$ endow $\Omega$ with a 
filtration 
of $\sigma$-algebras $\{\Sigma_t \, ;\, t \in \mathbb{R}^{+} \}$. Furthermore, the 
states $\rho^{\tau}$, defined by $\rho^{\tau} (W) := \rho\ten\phi(U^{\dagger}_{\tau} W U_{\tau})$ 
provide a family of consistent measures $\BB{P}_{\tau}$ on $(\Omega,\Sigma_{\tau})$, 
turning it into the probability space $(\Omega,\Sigma_t,\BB{P})$. 
(See e.g. \cite{BT4}.)

We aim to find random variables on $(\Omega,\Sigma_t,\BB{P})$ 
having distributions resembling those of $\sigma_x$ and $\sigma_z$ in 
the initial state $\rho$. In the previous section we used the 
random variables 
  \begin{equation} \label{naive choice}
  Y(\omega) = \frac{\omega_{\tau}}{2-2e^{-\frac{\tau}{2}}} \ \ \ \ \mbox{and} \ \ \ \
  \tilde{Y}(\omega) =\frac{\omega_{\tau}^2 -{\tau}}{4 (e^{-\frac{\tau}{2}}-1)^2} - 1,
  \ \ \ \ \tau = 2.513
  \end{equation}
for $\sigma_x$ and $\sigma_z$, respectively. Our next goal is to 
find the optimal  
random variables, in the sense of the previously defined quality. 


\subsection{Restricting the Class of Pointers}

In our specific example, $M$ is given by $M(B) = 
\id \ten\phi(U_{\tau}^{\dagger} \one\ten B U_{\tau})$. Note that stochastic 
integrals with respect to the annihilator $A_{\tau}$ acting 
on the vacuum vector $\Omega$ are zero. Therefore, we can modify 
$U_{\tau}$ to $Z_{\tau}$, given by
  \begin{equation*}
  dZ_{\tau} = \Big\{\sigma_{-} (dA^{\dagger}_{\tau}+dA_{\tau}) - 
  \half \sigma_{+} \sigma_{-} d\tau \Big\}Z_{\tau},\ \ \ \  Z_0 = \one, 
  \end{equation*}
without affecting $M$ \cite{Bel92a}.  
Therefore, for all $B \in \CB(\CF)$, we have $M(B) = \id \ten\phi(U^{\dagger}_{\tau} \one\ten
BU_{\tau}) =$ 
\mbox{$\id \ten\phi(Z^{\dagger}_{\tau} \one\ten B Z_{\tau})$}. The solution 
$Z_t$ can readily be found, it is given by
  \begin{equation*}
  Z_t = \begin{pmatrix}e^{-\frac{1}{2}t} & 0 \\ 
  \int_0^te^{-\frac{1}{2}\tau}(dA_{\tau}+dA_{\tau}^{\dagger}) & 1 \end{pmatrix}.
  \end{equation*}
Note that $Z_t$, as a matrix valued function of the measurement 
path, is an element of $M_2(\BB{C}) \ten \CC_t$, where $\CC_t$ is 
the commutative von Neumann algebra generated by 
$A_{\tau} + A_{\tau}^{\dagger},\ 0\le \tau \le t$.
Moreover we see that $Z_t$ is not a function of all the $(A_{\tau}+A_{\tau}^{\dagger})$'s 
separately, it is only a function of the endpoint of the weighted path 
$Y_{t} = \int_0^t e^{-\frac{1}{2}\tau} (dA_{\tau}+dA_{\tau}^{\dagger})$ \cite{Wis96}.
Therefore if we define $\mathcal{S}_t \subset \CC_t$ to be 
the commutative von Neumann algebra generated by $Y_t$,
then we even have $Z_{t} \in M_2(\BB{C}) \ten \mathcal{S}_t$. 

Denote by $C \mapsto \BB{E}[\, C \,|\CS_t]$ the unique classical 
conditional expectation from $\CC_t$ onto $\CS_t$ that 
leaves $\phi$ invariant, i.e.\ $\phi(\BB{E}[\,C\,|\CS_t]) = \phi(C)$ for 
all $C \in \CC_t$.
We can extend $ \BB{E}[\ \cdot \ |\CS_t]$ by tensoring it 
with the identity map on the $2 \times 2$ matrices 
to obtain a map $\id \ten \BB{E}[\ \cdot \ |\CS_t]$ from 
$M_2(\BB{C}) \ten \CC_t$ onto $M_2(\BB{C})\ten \CS_t$. 
From the positivity of $\BB{E}[\ \cdot \ |\CS_t]$ as 
a map between commutative algebras, it follows that 
$\id  \ten \BB{E}[\ \cdot \ |\CS_t]$ is completely positive.
Since $\BB{E}[\ \cdot \ |\CS_t]$ satisfies $\BB{E}[\,CS\,|\CS_t] = 
\BB{E}[\,C\,|\CS_t] S$ for all $C \in \CC_t$ and $S \in \CS_t$, 
we find that $\id  \ten \BB{E}[\ \cdot \ |\CS_t]$ satisfies 
the \emph{module property}, i.e. 
  \begin{equation*}
  \id  \ten \BB{E}[\ \cdot \ |\CS_t](A_1BA_2) 
  = A_1\Big(\id \ten\BB{E}[\ \cdot \ |\CS_t](B)\Big)A_2, 
  \end{equation*}
for all $A_1, A_2 \in M_2(\BB{C})\ten\CS_t$ and 
$B \in M_2(\BB{C})\ten\CC_t$. Moreover, if $\rho$ is a state on 
$M_2(\BB{C})$, then it follows from the invariance of $\phi$ under 
$\BB{E}[\ \cdot \ |\CS_t]$ that $\id \ten\BB{E}[\ \cdot \ |\CS_t]$ 
leaves $\rho \ten \phi$ invariant. We conclude that, given $\rho$ 
on $M_2(\BB{C})$, the map 
$\id \ten\BB{E}[\ \cdot \ |\CS_t]$ from $M_2(\BB{C})\ten\CC_t$ onto 
$M_2(\BB{C})\ten \CS_t$ is the 
unique conditional expectation in the noncommutative 
sense of \cite{Tak71} that leaves $\rho\ten\phi$ 
invariant. We will use the shorthand $\BB{E}_{\CS_t}$ for 
$\id \ten \BB{E}[\ \cdot \ |\CS_t]$ in the following.

\begin{lemma} \label{essence}
Let $C \in \CC_t $ be a pointer with quality $\Sigma_{C}$ such that $M(C) = X$. Then 
$\tilde{C} := \BB{E}[\,C\, | \mathcal{S}_t]$ is 
also a pointer with $M(\tilde{C}) = X$,
and its quality is at least as good; $\Sigma_{\tilde{C}} \le \Sigma_{C}$.
\end{lemma}
\proof
Note that for all states $\rho$ on $M_2(\BB{C})$ we have
\begin{eqnarray*}
  \rho\big(M(\tilde{C})\big) &=& 
  \rho\ten\phi\big(Z_t^{\dagger}\one \ten \tilde{C} Z_t\big)\\
  &=& 
  \rho\ten\phi\Big(Z_t^{\dagger} \BB{E}_{\CS_t} (\one\ten C) Z_t\Big)\\
  &=&
  \rho\ten\phi\Big(\BB{E}_{\CS_t} \big(Z_t^{\dagger} \one\ten C Z_t\big)\Big)\\ 
  &=& \rho\ten\phi(Z_t^{\dagger} \one \ten C Z_t)\\
  &=& \rho\big(M(C)\big)\,.
\end{eqnarray*}  
In the third step,
we used the module property and the fact that $Z_t$ is an 
element of $M_2 (\BB{C})\ten \CS_t$, and in the fourth 
step, we used 
invariance 
of $\rho \ten \phi$. Since this holds for 
all states $\rho$ on $M_2(\BB{C})$, we conclude
that $M(\tilde{C}) = M(C) = X$. 

As for the variance, we note first that the conditional expectation 
$\smash{\BB{E}_{\CS_t}}$ is a completely positive identity preserving map.
Therefore,
for all self-adjoint $C \in \CC_t$, we have
  \begin{equation*}
  \BB{E}_{\CS_t}(\one\ten C^2) \ge 
  \Big(\BB{E}_{\CS_t}(\one\ten C)\Big)^2
  \end{equation*}
by proposition \ref{posinp}.  

We can now apply the same strategy as before. For all states $\rho$ on 
$M_2(\BB{C})$ we have 
  \begin{eqnarray*}
  \rho\big(M(C^2)\big) &=& \rho \ten \phi (Z_t^{\dagger} \one \ten C^2 Z_t) \\
  &=& 
  \rho \ten \phi \Big(\BB{E}_{\CS_t}
  \big(Z_t^{\dagger} \one \ten C^2  Z_t\big)\Big)\\
  &=&  
  \rho \ten \phi \Big(Z_t^{\dagger}\BB{E}_{\CS_t}(\one \ten C^2) Z_{t}\Big)\\
  &\geq& 
  \rho \ten \phi \Big( Z_t^{\dagger}\Big(\BB{E}_{\CS_t}(\one \ten C)\Big)^2Z_{t}\Big)\\ 
  &=& \rho\big(M(\tilde{C}^2)\big) \,.
  \end{eqnarray*}
Thus $M(C^2) \geq M(\tilde{C}^2)$, so that 
$M(C^2) - M(C)^2 \geq M(\tilde{C}^2) - M(\tilde{C})^2$,
and in particular
$\Sigma_{C}^2 
\geq  
\Sigma_{\tilde{C}}^2$.
\qed

\noindent This has a very useful consequence: 
if we are looking for pointers that record, say, 
$\sigma_{x}$ or $\sigma_{z}$ in an optimal fashion, 
then it suffices to examine 
only pointers in $\mathcal{S}_t$.
Instead of sifting through the collection of all random variables 
on the measurement outcomes, we are thus allowed to 
confine the scope of our search to the 
rather transparent collection of  
measurable functions of $Y_{t}$. 
In the following, we will look at such pointers $h_t(Y_{t})$. 
We will usually drop the subscript $t$ on $h$ to 
make the notation lighter. 

\subsection{Distribution of \texorpdfstring{$Y_{t}$}{Y}}

At this point we are interested in the 
probability distribution of the random variable $Y_t$. 
Its characteristic function \cite{BaL85} is given by
  \begin{equation*}
  E(k) := \BB{E}_{\rho^t}\big[\exp(-ik Y_t)\big] = 
  \rho \ten \phi \Big(U_t^{\dagger} \one \ten \exp(-ik Y_t)U_t\Big) = 
  \rho \left(F_{\exp(-ik \x)}(\one)\right), 
  \end{equation*}  
so that we need only calculate $F_{\exp(-ik \x)}(\one)$. For notational 
convenience we will replace the subscript $\exp(-ik \x)$ by $k$ in the
following. Using lemma \ref{lem eq F}, we find the following system of
matrix valued differential equations:
\begin{equation*} \mathrm{
  \begin{tabular}{l c l}
  $\displaystyle \frac{d}{dt}F_k(\one) $ \rule[-4,2 mm]{0 mm}{5 mm} &=&$
  \displaystyle -ik e^{-\frac{t}{2}}F_k(\sigma_x ) -
  \frac{k^2e^{-t}}{2} F_k(\one),$\\ 
  $\displaystyle \frac{d}{dt} F_k(\sigma_{x})$ \rule[-4,2 mm]{0 mm}{5 mm} 
  &=&$ \displaystyle -\half F_k(\sigma_{x} ) 
  - 2ike^{-\frac{t}{2}} F_k(\sigma_{+} \sigma_{-} ) 
  - \frac{k^2 e^{-t}}{2}F_k(\sigma_{x}), $\\
  $\displaystyle \frac{d}{dt}F_k(\sigma_{+} \sigma_{-} ) $ 
  \rule[-4,2 mm]{0 mm}{5 mm}&=&$
  \displaystyle -F_k(\sigma_{+} \sigma_{-} ) - 
  \frac{k^2e^{-t}}{2}F_k(\sigma_{+} \sigma_{-} ),$\\
   \end{tabular}
 }\end{equation*}
 with initial values 
 $\displaystyle F_k(\one,0) = \one$, 
 $F_k(\sigma_{x} , 0) = \sigma_{x}$ and
 $F_k(\sigma_{+} \sigma_{-} , 0) = \sigma_{+} \sigma_{-}$.
 (Recall that $\sigma_x = \sigma_+ + \sigma_-$.)
Solving this system leads to
  \begin{equation*}
  F_k(\one) = e^{-\frac{k^2(1-e^{-t})}{2}}\Bigg(\one  
  -ik \big(1-e^{-t}\big)\sigma_{x}  - 
  k^2\big(1-e^{-t}\big)^2\sigma_{+} \sigma_{-} \Bigg).
  \end{equation*} 
We define the Fourier transform to be $\CF(f)(x) := 
\frac{1}{\sqrt{2\pi}}\int_{-\infty}^\infty f(k) e^{ikx}dk$. Then 
the probability density of $Y_t$ with respect to the Lebesgue
measure is given by $\frac{1}{\sqrt{2\pi}} \CF(E)(x) = 
\frac{1}{\sqrt{2\pi}} \rho\Big(\CF\big(F_k(\one)\big)(x)\Big)$. 
Defining $p(x):=\frac{1}{\sqrt{2\pi}} \CF\big(F_k(\one)\big)(x)$, 
we can write
  \begin{equation*}
  p(x) =  \frac{e^{-\frac{1}{2}\frac{x^2}{1-e^{-t}}}}{\sqrt{2\pi(1-e^{-t})}}
  \Big(\one + x\sigma_x +
   (x^2-1+e^{-t})\sigma_{+} \sigma_{-} \Big),
  \end{equation*}
i.e.\ $Y_t$ is distributed according to a Gaussian perturbed by the matrix 
elements of the initial state $\rho(\sigma_x)$ and $\rho(\sigma_{+} \sigma_{-} ) = 
\frac{1}{2}\rho(\sigma_z)+\frac{1}{2}$. No information about $\rho$ on 
$\sigma_y$ enters the distribution though. To gain information about $\sigma_y$ 
we would have to change our continuous time measurement setup, 
as we discussed before. If we absorb a constant $(1-e^{-t})^{-\frac{1}{2}}$ in 
the definition of $Y_t$, i.e. if we redefine $Y_t := 
(1-e^{-t})^{-\frac{1}{2}}\int_0^t
e^{-\frac{\tau}{2}}(dA_{\tau}+dA_{\tau}^{\dagger})$,  
then its density becomes
  \begin{equation}\label{eq p}
  p(y) =  \frac{e^{-\frac{y^2}{2}}}{\sqrt{2\pi}}
  \Big(\one + \beta_ty\sigma_x +
   \beta_t^2(y^2-1)\sigma_{+} \sigma_{-} \Big),
  \end{equation}
where $\beta_t : = \sqrt{1-e^{-t}}$.

\section{Variational Calculus} \label{sec variation}

In lemma \ref{essence}, we have shown that it suffices to 
consider only random variables of the form $h(Y_t)$ for some
measurable $h$.
In equation \eqref{eq p}, we have captured the probability 
distribution of $Y_t$.
All that remains now is to 
calculate the optimal $h$, which can be done with variational 
calculus. 

\subsection{Optimal \texorpdfstring{Measurement of $\sigma_x$}{Measurement 
of the Spin in the X-Direction}}

We seek the function $h^*$ for which 
the quality $\Sigma$ of the pointer $h^*(Y_t)$ for 
$\sigma_x$-measurement is optimal.
 In other words, we need  
\begin{equation}\label{eq sigma2}
  \Sigma^2 := \left\| \int_{-\infty}^\infty h^2(y)p(y)dy - 
  \Big(\int_{-\infty}^\infty h(y)p(y)dy\Big)^2\right\|
  := \left\| \begin{pmatrix} d_1 & 0 \\ 0 & d_2 \end{pmatrix}  \right\|
  \end{equation}  
to be minimal under the restriction 
$\int_{-\infty}^\infty h(y)p(y)dy = \sigma_{x}$.

Now $\Sigma^2$ is the norm of a diagonal $2 \times 2$-matrix with entries $d_1$ and $d_2$.
Both depend smoothly on $h$, but $\Sigma^2 = \mbox{max}\{ d_1 , d_2 \}$ does not.
There are three possibilities:
\begin{itemize}
\item[{\footnotesize I)}] $\Sigma^2 = d_1$ in some open neighbourhood of $h^*$. 
	   To find these $h^*$, we must minimize the smooth functional $d_1$ and then
	   check whether $d_1 < d_2$. 	  
\item[{\footnotesize II)}] $\Sigma^2 = d_2$ in some open neighbourhood of $h^*$. To find these $h^*$, we must 
	   minimize $d_2$ and check whether $d_2 < d_1$.  
\item[{\footnotesize III)}] $d_1 = d_2$ for $h^*$. To find these $h^*$, 
we must minimize $d_1$ subject to the condition $d_1 = d_2$.
\end{itemize} 
In principle, we need three different functionals $\Lambda_1$, 
$\Lambda_2$ and $\Lambda_3$ 
for these three distinct cases. 
However, it turns out that we can make due with the following functional
 \begin{equation}\label{eq tildeLambda}\begin{split} 
  \Lambda(h,\kappa,\gamma_1,\gamma_2,\gamma_3) :=\ 
  &\Big( \frac{1}{\sqrt{2\pi}}
  \int_{-\infty}^\infty h^2(y)e^{-\frac{1}{2}y^2}dy - 1\Big)\ + \\
  &\kappa\Big(\frac{\beta_t^2}{\sqrt{2\pi}}
  \int_{-\infty}^\infty h^2(y)(y^2-1)e^{-\frac{1}{2}y^2}dy\Big)\ + \\
  &\gamma_1\Big(\int_{-\infty}^\infty h(y)e^{-\frac{1}{2}y^2}dy\Big) \ + \\   
  &\gamma_2\Big(\int_{-\infty}^\infty h(y)(y^2-1)e^{-\frac{1}{2}y^2}dy\Big)\ + \\
  &\gamma_3\Big(\frac{\beta_t}{\sqrt{2\pi}}
  \int_{-\infty}^\infty h(y)ye^{-\frac{1}{2}y^2}dy - 1\Big).
  \end{split}\end{equation} 
The constants 
$\gamma_1, \gamma_2$ and $\gamma_3$ are the Lagrange multipliers enforcing 
the identity $\int_{-\infty}^\infty h(y)p(y)dy = \sigma_{x}$. 
These are needed in all cases: $\Lambda_1$, $\Lambda_2$ and $\Lambda_3$.
One can readily check that setting $\kappa = 0$ in $\Lambda$ yields $\Lambda_1$,
setting $\kappa = 1$ yields $\Lambda_2$ 
and considering $\kappa$ as a free Lagrange
multiplier forces $d_1 = d_2$, so that one has $\Lambda_3 = \Lambda$. 

All three cases lead to similar optimality conditions.
The requirement that the optimal solution be stable under first order variations 
yields $h^*$ satisfying either
\begin{equation}\label{one}
  h^*(x) = \frac{C_1 x + C_2}{x^2 + \varepsilon} + C_3 
\end{equation}
or
\begin{equation}\label{two}
  h^*(x) = C_{4}x^2 + C_{5}x + C_{6} \,     
\end{equation}
for some real constants $C_1,C_2,C_3,C_4, C_5, C_6$ and $\varepsilon$ depending 
on $\kappa,\gamma_1,\gamma_2,\gamma_3$.

Suppose that $h^*$ takes the form (\ref{two}). The constraint
$\int_{-\infty}^{\infty} h^*(y) p(y) dy = \sigma_{x}$ will then force
$C_{4} = C_{6} = 0$ and $C_{5} = \beta_{t}^{-1}$, so that 
$h^*(y) = \beta_{t}^{-1} y$. The random variable 
we are investigating is simply the observed path, weighted by the function 
$f(\tau) = \beta_t^{-1}e^{-\tau /2}$, with $t$ the final time of measurement.
Since all the integrals we encounter are Gaussian moments, we can readily compute 
$M(h^{*2}(Y_t)) = \int_{-\infty}^{\infty} h^{*2}(y) p(y) dy \ $ 
to be $ \ 2\sigma_{+} \sigma_{-}  + \beta_{t}^{-2} \one$.
Thus $$\Sigma^2 = \|(2\sigma_{+} \sigma_{-}  + \beta_{t}^{-2} \one ) 
- (\sigma_{-} +\sigma_{+} )^2  \| = 1 + \beta_{t}^{-2}\,.$$
For $t \to \infty$, this amounts to $\Sigma \to \sqrt{2}$.
Already, we have improved on the naive result $\Sigma = 1.493$ 
obtained previously.

We proceed with the more involved case (\ref{one}), which will 
provide us with the optimal solution.
Before we continue with the constants 
$C_1,C_2,C_3$ and $\varepsilon$ however,
we calculate some integrals for later use.

\begin{basdef}  
Define the error function $\mathrm{erf}(x)$ and integrals $I(\varepsilon)$ and $J(\varepsilon)$ by
  \begin{equation*}
  \mathrm{erf}(x) := \frac{2}{\sqrt{\pi}}\int_{0}^x e^{- u^2}du,\ \ \ \
  I(\varepsilon) := \int_{-\infty}^\infty \frac{e^{-\frac{x^2}{2}}}{x^2+\varepsilon}dx,
  \ \ \ \
  J(\varepsilon) := \int_{-\infty}^\infty \frac{e^{-\frac{x^2}{2}}}{(x^2+\varepsilon)^2}dx.
  \end{equation*}
\end{basdef}

\begin{lemma}\label{integrals}
  \begin{equation*}
  J(\varepsilon) = \frac{\sqrt{2\pi} + (1-\varepsilon)I(\varepsilon)}{2\varepsilon} \ \ \ \
  \mbox{and} \ \ \ \ 
  I(\varepsilon)= \pi \sqrt{\frac{e^{\varepsilon}}{\varepsilon}}
  \left( 1  - \mathrm{erf}\left(\sqrt{\frac{\varepsilon}{ 2}}\right) \right).
  \end{equation*}
\end{lemma}
\proof
Since the Fourier transform of $e^{-\sqrt{\varepsilon}|k|}$ is equal to 
$\sqrt{\frac{2\varepsilon}{\pi}}\frac{1}{x^2 + \varepsilon}$, we find
  \begin{eqnarray*}
  I(\varepsilon) &=& \sqrt{\frac{\pi}{2\varepsilon}}\int_{-\infty}^\infty 
  \CF\big(e^{-\sqrt{\varepsilon}|k|}\big)\CF\big(e^{-\frac{k^2}{2}}\big)dx \\
  &=& \sqrt{\frac{\pi}{2\varepsilon}}\int_{-\infty}^\infty 
  e^{-\sqrt{\varepsilon}|k|}e^{-\frac{k^2}{2}}dk\,,
  \end{eqnarray*}
in view of the unitarity of the Fourier transform. Then
\begin{eqnarray*}
  I(\varepsilon)&=& 
  \sqrt{\frac{2\pi}{\varepsilon}}\int_{0}^\infty 
  e^{-\sqrt{\varepsilon}k}e^{-\frac{k^2}{2}}dk \\
  &=&
  \sqrt{\frac{2\pi}{\varepsilon}}
  e^{\frac{1}{2}\varepsilon}\int_{\sqrt{\varepsilon}}^\infty
  e^{-\frac{u^2}{2}}du\\
  &=& \pi \sqrt{\frac{e^{\varepsilon}}{\varepsilon}}
  \left( 1 - \mathrm{erf}\left(\sqrt{\frac{\varepsilon}{2}}\right)\right).
\end{eqnarray*}
%
The expression for $J$ then follows from 
  \begin{equation*}\begin{split}
  0 & = \frac{xe^{-\frac{x^2}{2}}}{x^2+\varepsilon}\Bigg|_{-\infty}^\infty = 
  \int_{-\infty}^\infty
  \frac{d}{dx}\left(\frac{xe^{-\frac{x^2}{2}}}{x^2+\varepsilon}\right)dx =
  \int_{-\infty}^\infty 
  \left(\frac{1-x^2}{x^2+\varepsilon} -
  \frac{2x^2}{(x^2+\varepsilon)^2}\right)e^{-\frac{x^2}{2}}dx  \\
  & = \int_{-\infty}^\infty \left(-1 + \frac{\varepsilon -1}{x^2 + \varepsilon} + 
  \frac{2\varepsilon}{(x^2+\varepsilon)^2}\right)e^{-\frac{x^2}{2}}dx = 
  -\sqrt{2\pi} + (\varepsilon -1)I(\varepsilon) + 2\varepsilon J(\varepsilon).
  \end{split}\end{equation*}
\qed
The condition $\int_{-\infty}^\infty h^*(y)p(y)dy = \sigma_{-}  + \sigma_{+}  
=\sigma_x$ implies
  \begin{equation*}
  C_1 = \frac{\sqrt{2\pi}}{\beta_t\big(\sqrt{2\pi}- \varepsilon I(\varepsilon)\big)}, 
  \ \ \ \ C_2 = C_3 = 0
  \end{equation*}
which fixes $C_1$ as a function of $\varepsilon$. 
The next step is to express $d_1$ and $d_2$ in terms of $\varepsilon$:
  \begin{eqnarray*}
  d_2 &=& \frac{C_1^2}{\sqrt{2\pi}}\int_{-\infty}^
  \infty \frac{y^2}{(y^2+\varepsilon)^2}e^{-\frac{y^2}{2}}dy - 1\\
  &=& \frac{C_1^2}{\sqrt{2\pi}}\Big(I(\varepsilon)-
  \varepsilon J(\varepsilon)\Big)-1,\\  
  d_1 & =& \frac{C_1^2\beta_t^2}{\sqrt{2\pi}}\int_{-\infty}^\infty
  \frac{y^2(y^2-1)}{(y^2 + \varepsilon)^2}e^{-\frac{y^2}{2}}dy + d_2  \\
  &=& 
  \frac{C_1^2\beta_t^2}{\sqrt{2\pi}}\Big(\sqrt{2\pi} - (1+2\varepsilon)I(\varepsilon)+ 
  \varepsilon(1+\varepsilon)J(\varepsilon)\Big) + d_2.
  \end{eqnarray*}
%
First, we use lemma \ref{integrals} to express the above in terms of elementary functions
and the error function. Then, using Maple, we find that $\varepsilon \mapsto \mbox{max}\{d_1,d_2\}$ 
has a unique minimum at $\varepsilon = 0.605$, 
for which $d_1 =d_2 = 0.470$.
This leads to 
a $C_1$ that equals 2.359,
and to a quality of 
  \begin{equation*}
  \Sigma = \sqrt{\mbox{max}\{d_1,d_2\}} = 0.685 \,.
  \end{equation*}  

\subsection{Optimal Measurement of \texorpdfstring{$\sigma_z $}{ the Spin in the 
Z-Direction}}

For optimal $\sigma_z$-measurement, we can run the same program.  
We search for the function $\tilde{h}$ that optimizes the quality 
$\tilde{\Sigma}$, under the restriction that $\tilde{h}(Y_t)$ 
be a pointer for $\sigma_z$-measurement. 
That is, we search for a function $\tilde{h}$ minimizing 
the functional of equation \eqref{eq sigma2}, but now 
under the restriction $\int_{\infty}^\infty h(y)p(y)dy = \sigma_z$.
Again there are three cases of interest, $d_1=d_2$, $d_1>d_2$ and $d_2>d_1$, 
which we can treat simultaneously by introducing, analogous to  
equation \eqref{eq tildeLambda}, the functional
  \begin{eqnarray*} 
  \tilde{\Lambda}(h,\kappa,\gamma_1,\gamma_2,\gamma_3) &:=& 
  \Big( \frac{1}{\sqrt{2\pi}}\int_{-\infty}^\infty 
  h^2(y)e^{-\frac{1}{2}y^2}dy - 1\Big) \\
  &+& 
  \kappa\Big(\frac{\beta^2_t}{\sqrt{2\pi}}
  \int_{-\infty}^\infty h^2(y)(y^2-1)e^{-\frac{1}{2}y^2}dy\Big)\\
  &+& 
  \gamma_1\Big(\frac{1}{\sqrt{2\pi}}\int_{-\infty}^\infty h(y)
  e^{-\frac{1}{2}y^2}dy+1\Big)\\ 
  &+&   
  \gamma_2\Big(\frac{\beta^2_t}{\sqrt{2\pi}}\int_{-\infty}^\infty h(y)
  (y^2-1)e^{-\frac{1}{2}y^2}dy-2\Big)\\
   &+& 
  \gamma_3\Big(\int_{-\infty}^\infty h(y)ye^{-\frac{1}{2}y^2}dy\Big).
  \end{eqnarray*}
Indeed, $\gamma_1, \gamma_2$ and $\gamma_3$ are the Lagrange multipliers 
enforcing the restriction $\int_{\infty}^\infty h(y)p(y)dy = \sigma_z$. 
Again, the functional $\Sigma^2$ of equation \eqref{eq sigma2} depends 
nondifferentiably on $h$ when $d_1=d_2$. We then have to search for the 
optimum among the points of nondifferentiability, in which case 
$\kappa$ is the Lagrange multiplier confining us to these points. If $d_1 >d_2$
then $\kappa = 1$ and if $d_2 > d_1$ then 
$\kappa = 0$. Summarizing, wherever $\Lambda$ takes its minimal value, 
optimality implies $\frac{\delta\tilde\Lambda}{\delta h}
(\tilde{h}, \kappa, \gamma_1, \gamma_2, \gamma_3) = 0$ for some 
$\kappa, \gamma_1, \gamma_2$ and $\gamma_3$. Performing the 
functional derivative 
yields either
  \begin{equation}\label{one1}
  \tilde{h}(x) = \frac{D_1x+D_2}{x^2 + \delta} + D_3  
  \end{equation}
or
\begin{equation}\label{two2}
  \tilde{h}(x) = D_{4} x^2 + D_{5}x + D_{6}  
\end{equation}   
for some (time-dependent) constants $D_1, D_2, D_3, D_4, D_5, D_6$ and 
$\delta$ depending on 
$\kappa, \gamma_1, \gamma_2$ and $\gamma_3$. 

Again, we begin with the least demanding case \eqref{two2}, resulting from $\kappa = 0$.
The condition  $\int_{-\infty}^\infty \tilde{h}(y)p(y)dy = \sigma_z$ implies 
$D_{5} = 0$, $D_{4} = \beta_{t}^{-2}$ and $D_{6} = -1 -\beta_{t}^2$.
For $t \to \infty$,
this leads to 
\begin{eqnarray*}
\Sigma^2 & = & \| M(\tilde{h}^2(Y_{t})) - M(\tilde{h}(Y_{t}))^2 \|\\
         & = & \|(4 \sigma_{+} \sigma_{-}  + 3 \one) - \one  \| \\
         & = & 6\,,
\end{eqnarray*}
so that $\Sigma \to \sqrt{6}$.

This improves the result $\tilde{\Sigma} = 2.973$ obtained 
previously, but once again the ultimate bound will be 
reached in the more arduous case \eqref{one1}.
There, the condition 
$\int_{-\infty}^\infty \tilde{h}(y)p(y)dy = \sigma_z$ implies 
  \begin{eqnarray*}
  D_3 &=& - \frac{\sqrt{2\pi} + I(\delta)D_2}{\sqrt{2\pi}},\\ 
  D_2 &=& \frac{2\sqrt{2\pi}}{\beta_t^2\big(\sqrt{2\pi}-(1+\delta)I(\delta)\big)},\\ 
  0 & = & D_1\big(\sqrt{2\pi}-\delta I(\delta)\big) \,.
  \end{eqnarray*}
This leads to expressions for $d_1$ and $d_2$ as a function of $\delta$.
Using lemma \ref{integrals} and Maple once more, we find that the 
function $\delta \mapsto \mbox{max}\{d_1,d_2\}$ 
has a unique minimum at $\delta = 2.701$,
for which $d_1=d_2 = 2.373$. This leads to a quality of
  \begin{equation*}
  \tilde{\Sigma} = \sqrt{\mbox{max}\{d_1,d_2\}} = 1.540,
  \end{equation*}  
  attained for $D_1 = 0$, $D_2 = -21.649$ and $D_3 = 5.391$.
For the joint measurement this leads to
  \begin{equation*}
  \Sigma\tilde{\Sigma} = 1.056\,.
  \end{equation*}
Although we did not achieve the bound of 1 provided by theorem \ref{thm joint},
we have come as close as the measurement setup allows.
We conclude that, 
using the setup investigated in this chapter, no simultaneous measurement 
of $\sigma_x$ and $\sigma_z$ will be able to approach the quantum bound 
by more than 5.6 \%.  
Furthermore, we have identified the unique pointers for this optimal 
measurement in equations \eqref{one} and \eqref{one1}.

\section{Distribution of Pointer Variables}\label{sec densities}

We have designed pointers $h^*(Y_t)$ and 
$\tilde{h}(Y_t)$ in such a way that their distributions in the final state
best resemble the distributions of $\sigma_x$
and $\sigma_z$ in the initial state.
We will now calculate and plot these final densities. 

\subsection{Calculation of \texorpdfstring{$h^*$- and $\tilde{h}$-}{some }Densities}

Let $\rho$ be the initial state of 
the qubit and let it be parameterized by its Bloch vector $(P_x,P_y,P_z)$. 
By equation \eqref{eq p}, the density $q(y)$ of $Y_t$ is given by
  \begin{equation}\label{eq q}
  q(y) = \rho\big(p(y)\big) = \frac{e^{-\frac{1}{2}y^2}}{\sqrt{2\pi}}
  \Big(1 + \beta_t yP_x + \beta_t^2(y^2-1)\frac{P_z+1}{2}\Big).
  \end{equation} 

We are interested in the the distributions $r(x)$ and $s(x)$ 
of $h^*(Y_t)$ and $\tilde{h}(Y_t)$ respectively. 
Let us start with $h^*$. From equation 
\eqref{one}, we first calculate the points $y$ 
where $h^*(y) = x$ for some fixed value of $x$.  
  \begin{equation*}
  y_{\pm} = \frac{C_{1}\pm\sqrt{C_{1}^2-4x^2\varepsilon}}{2x} 
  \end{equation*}
By the Frobenius-Peron equation (see e.g.\ \cite{Ott}), $r(x)$ is given by
  \begin{equation*}
  r(x) = \sum_{+,-} \frac{q(y_\pm)}{|{h^*}'(y_\pm)|}, 
  \end{equation*}
which leads immediately to
  \begin{equation} \label{rofx}
  r(x) = \sum_{+,-} \frac{(y_\pm^2+\varepsilon)^2
  \Big(1 + \beta_t y_\pm P_x + \beta_t^2(y_\pm^2-1)\frac{P_z+1}{2}\Big)}{C_{1}|y^2_\pm -\varepsilon|} 
  \frac{e^{-\frac{1}{2}y_\pm^2}}{\sqrt{2\pi}}, 
  \end{equation}
where it is understood that $r(x) \neq 0$ only for $x \in
[-\frac{C_{1}}{2 \sqrt{\varepsilon}}, \frac{C_{1}}{2 \sqrt{\varepsilon}}]$. 
We run a similar analysis for $s(x)$. The points $y$ in which 
$\tilde{h}(y)= x$ are given by 
  \begin{equation*}
  y_\pm = \pm \sqrt{\frac{(x-D_3)\delta-D_2}{D_3-x}}.
  \end{equation*}
This leads to 
  \begin{equation} \label{sofz}
  s(x) = \sum_{+,-} 
  \frac{(y_\pm^2+\delta)^2
  \Big(1 + \beta_t y_\pm P_x + \beta_t^2(y_\pm^2-1)\frac{P_z+1}{2}\Big)}{2|D_2y_\pm|}
  \frac{e^{-\frac{1}{2}y_\pm^2}}{\sqrt{2\pi}},  
  \end{equation}  
with $s(x) \neq 0$ only for $x \in [D_{3} + \frac{D_{2}}{\delta}, D_{3}]$. 
We proceed with a graphical illustration of the results obtained so far. 

\subsection{Plots of \texorpdfstring{$\sigma_{x}$-}{}Measurement}

According to formula \ref{eq q}, the 
distribution of the endpoint of the weighted path depends on the input qubit-state. 
For instance, the negative $\sigma_{x}$-eigenstate, the tracial state and
the positive $\sigma_{x}$-eigenstate lead to the
distributions below:\vspace{-7mm}

\begin{center}
\begin{tabular}{l c l}
	\begin{minipage}[t]{3.6cm}
		\includegraphics[width = 3.6 cm, viewport = 50 100 500 700]{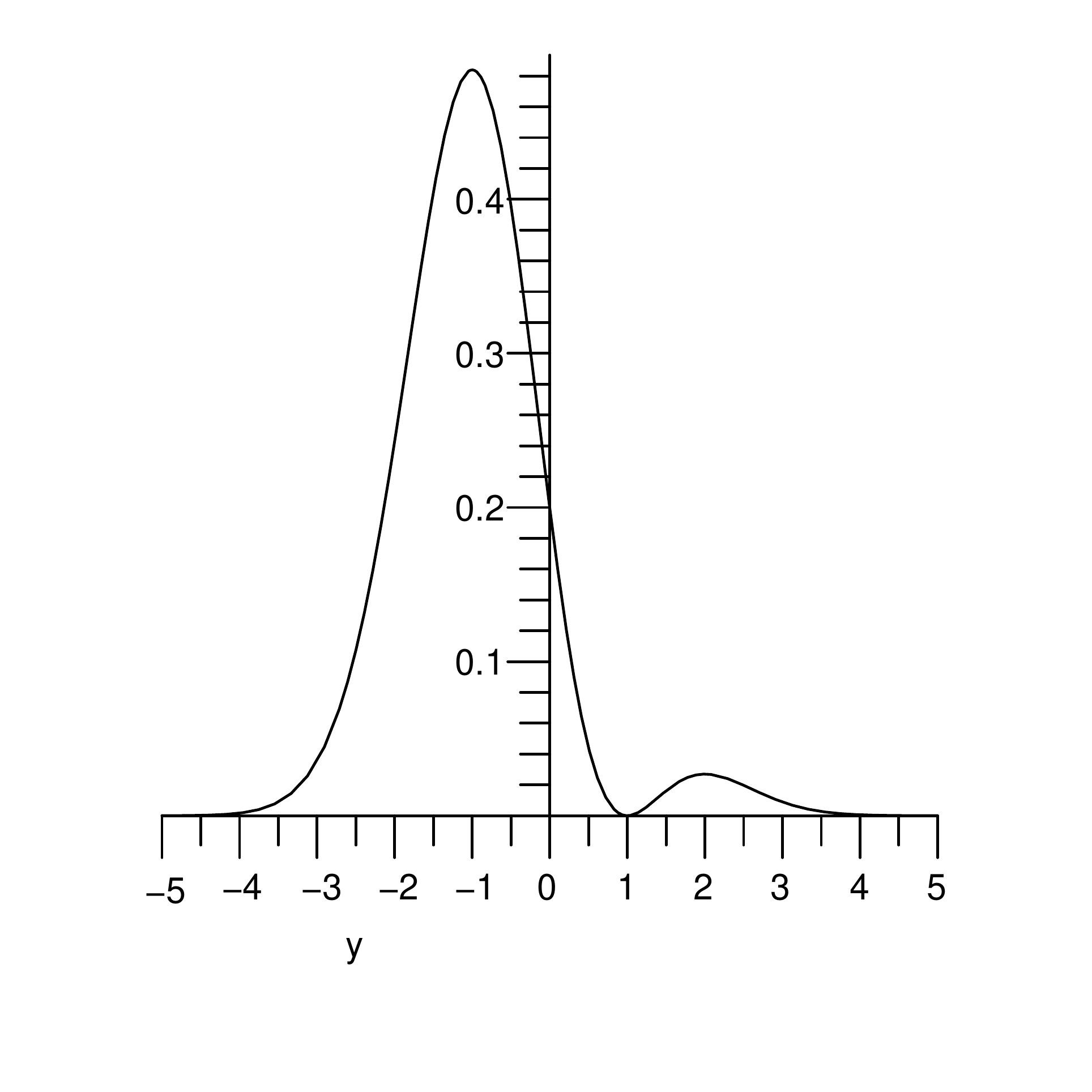}
		\\[5mm]{\footnotesize \fig Probability density of the endpoint 
		of the weighted path 
		 for the input state $|\leftarrow \, \rangle $.}
		\label{fig: g(y)-100}
	\end{minipage}
&  
	\begin{minipage}[t]{3.6cm}
		\includegraphics[width = 3.6 cm, viewport = 50 100 500 700]{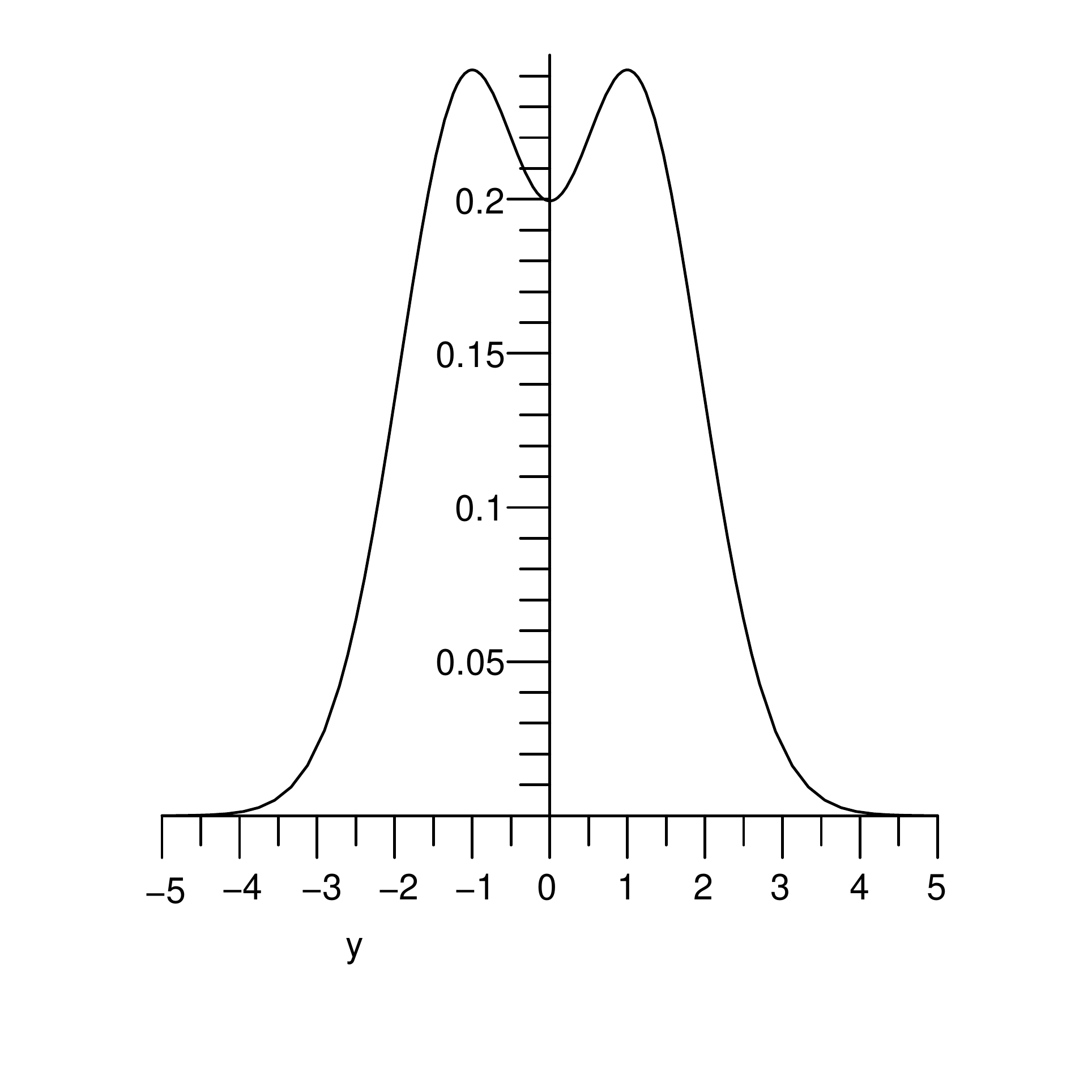}
		\\[5mm]{\footnotesize \fig Probability density of the endpoint of the 
		weighted path 
		 for 
		 the tracial input state.}
		\label{fig: g(y)000}
	\end{minipage}
& 
	\begin{minipage}[t]{3.6cm}
		\includegraphics[width = 3.6 cm, viewport = 50 100 500 700]{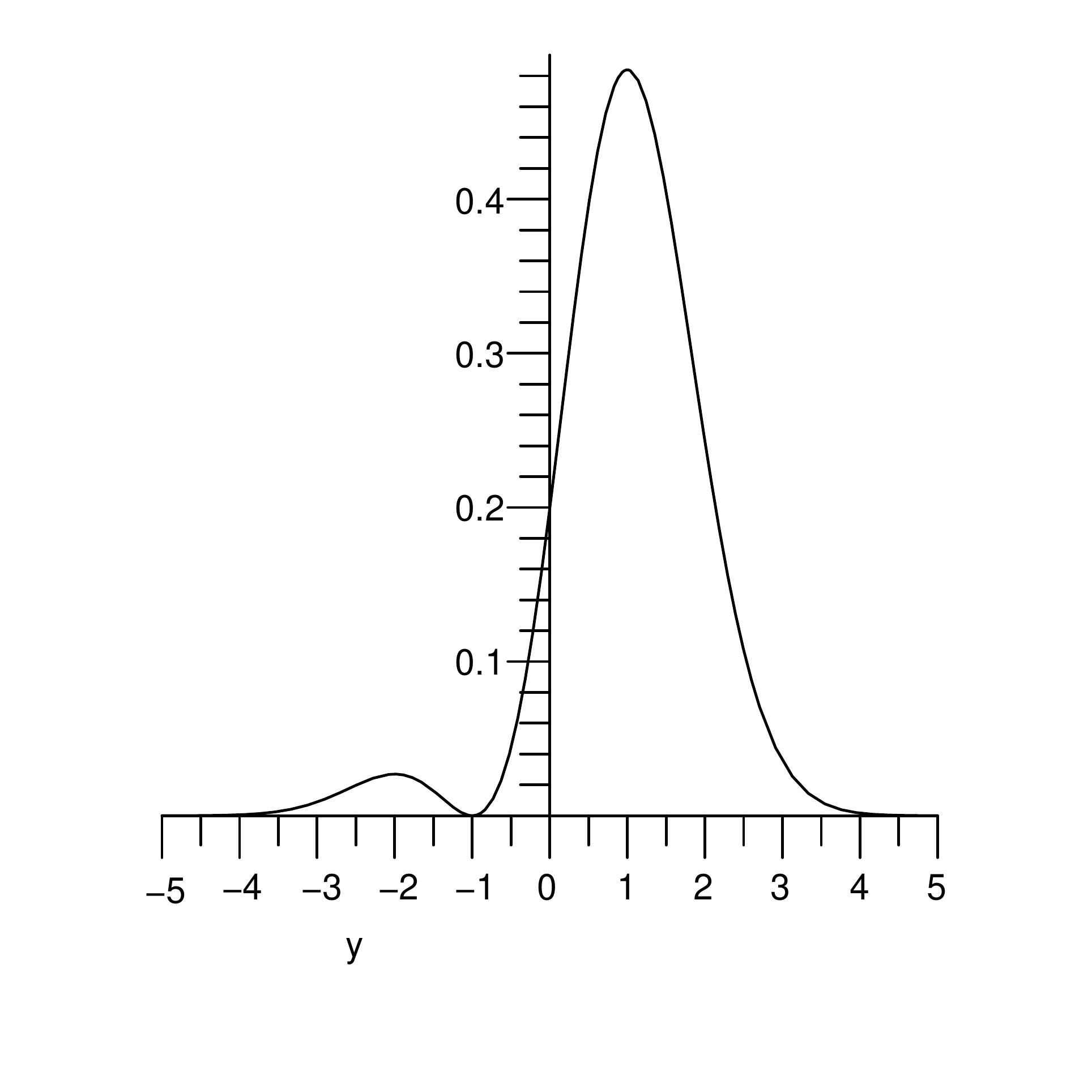}
		\\[5mm]{\footnotesize \fig  Probability density of the endpoint of the 
		weighted path 
		 for the input state $|\rightarrow \, \rangle $.}
	\end{minipage}
\\[2.8cm]
	\multicolumn{2}{p{7.6 cm}} 
	   {	In order to estimate $\sigma_x$, we use the pointer			
	   	of $\sigma_x$ given by\vspace{-1mm} 
		\begin{equation}\label{hstar}
		h^*(x) = \frac{C_{1} x}{x^2 + \varepsilon},
		\end{equation}
		with
		$\varepsilon = 0.605$ 
		and
		$C_{1} = 2.359$. It is illustrated to the right. 
		In formula \eqref{rofx}, we have calculated the 
		probability distribution of this
		pointer under the distribution of 
		the endpoint of the weighted path. The distribution of the weighted path
		is illustrated above, that of the pointer 
		below. 
 	  }	
& 
	\begin{minipage}[t]{3.6cm}
		\mbox{}\\[-2.1mm]
		\includegraphics[width = 3.6 cm, viewport = 0 0 400 400]{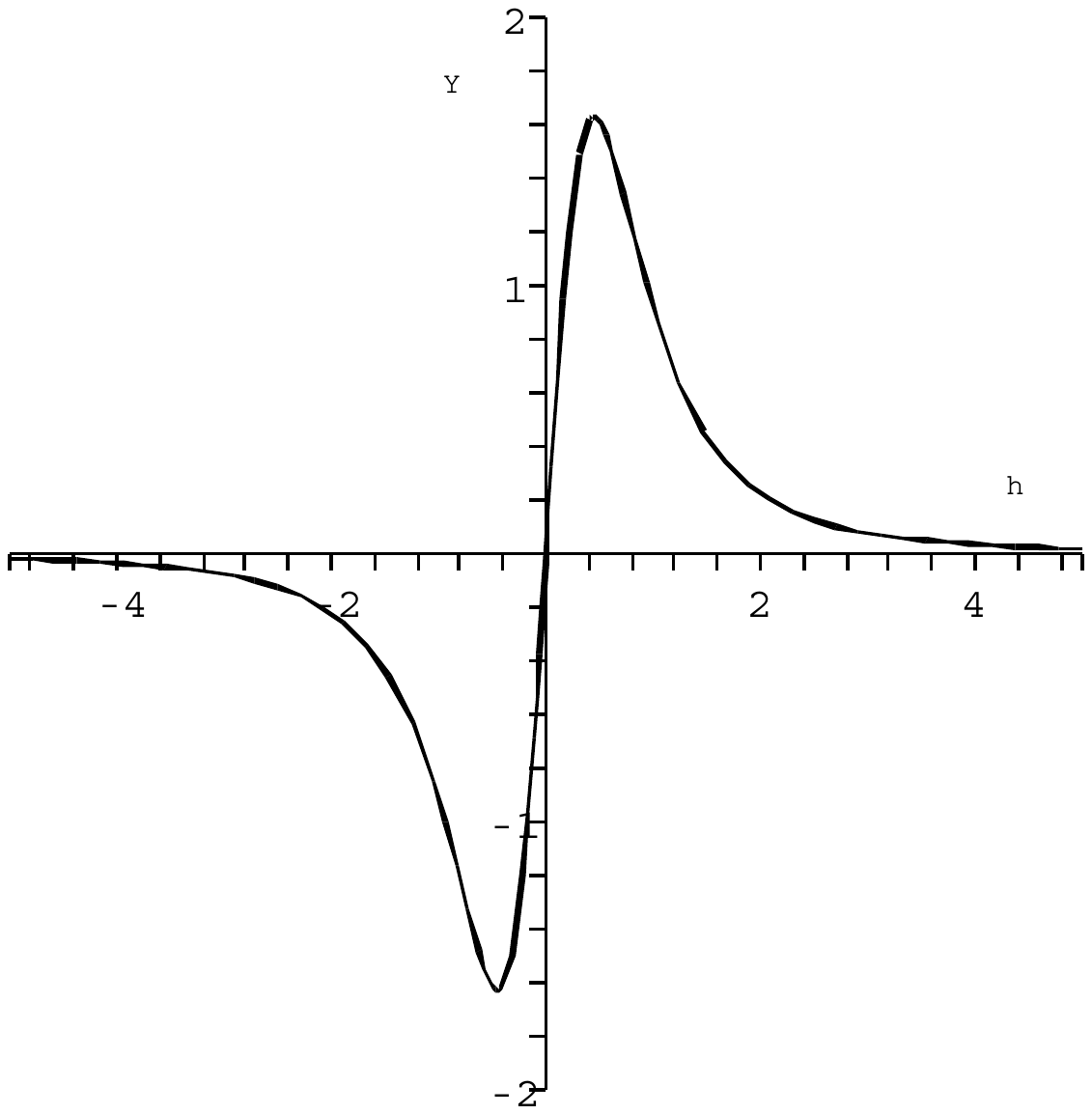}
		\label{Pointer x}
		 \\[-5mm] {\footnotesize \fig Pointer for $\sigma_x$}\\
		\mbox{}  
	\end{minipage}
\\[-8mm]
	\begin{minipage}[t]{3.6 cm}
		\includegraphics[width = 3.6 cm, viewport = 50 100 500 700]{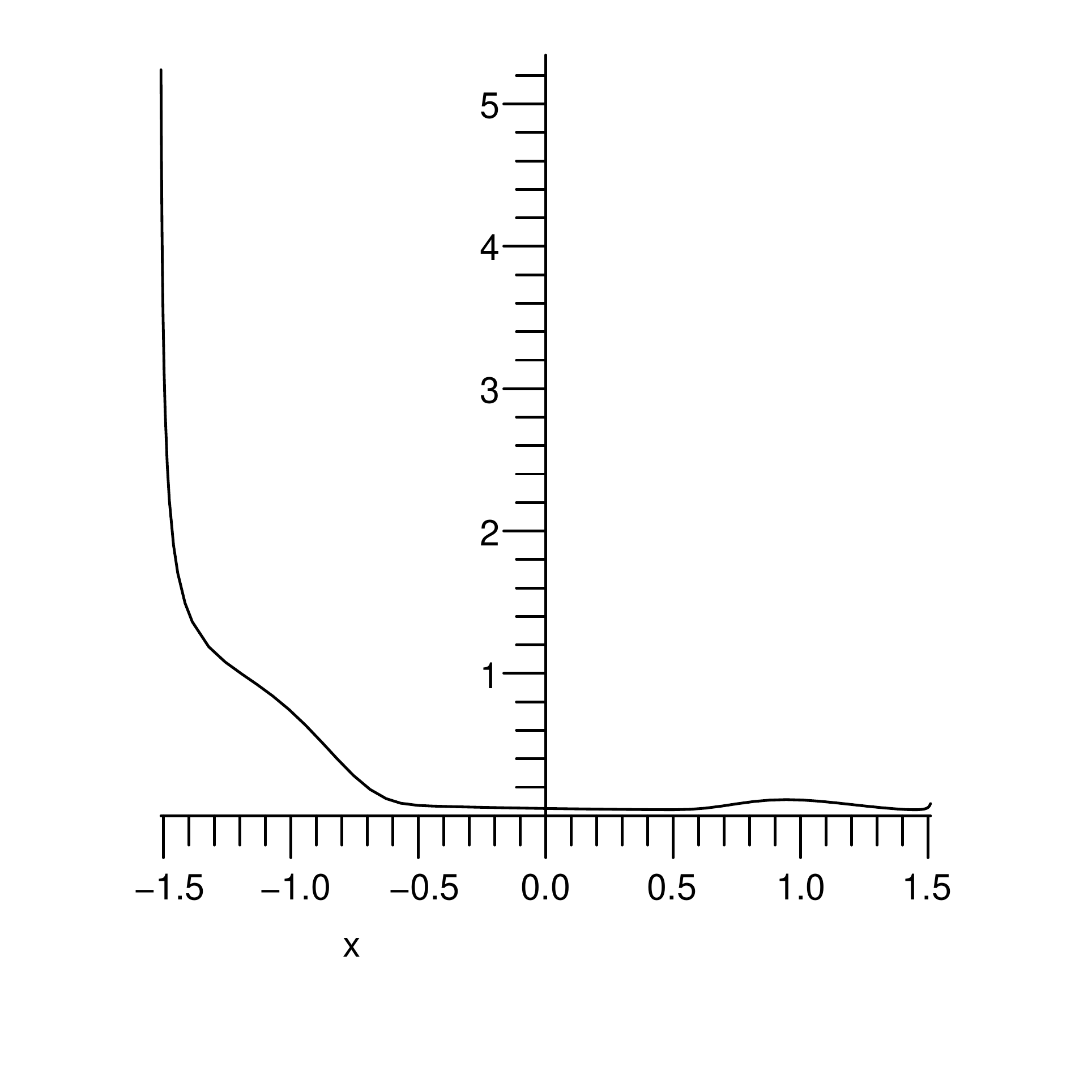}
		\\[5mm] {\footnotesize \fig Probability density of the $\sigma_x$-pointer 
		for the input state $|\leftarrow \, \rangle $.}
	\end{minipage}
&
	\begin{minipage}[t]{3.6 cm}
		\includegraphics[width = 3.6 cm, viewport = 50 100 500 700]{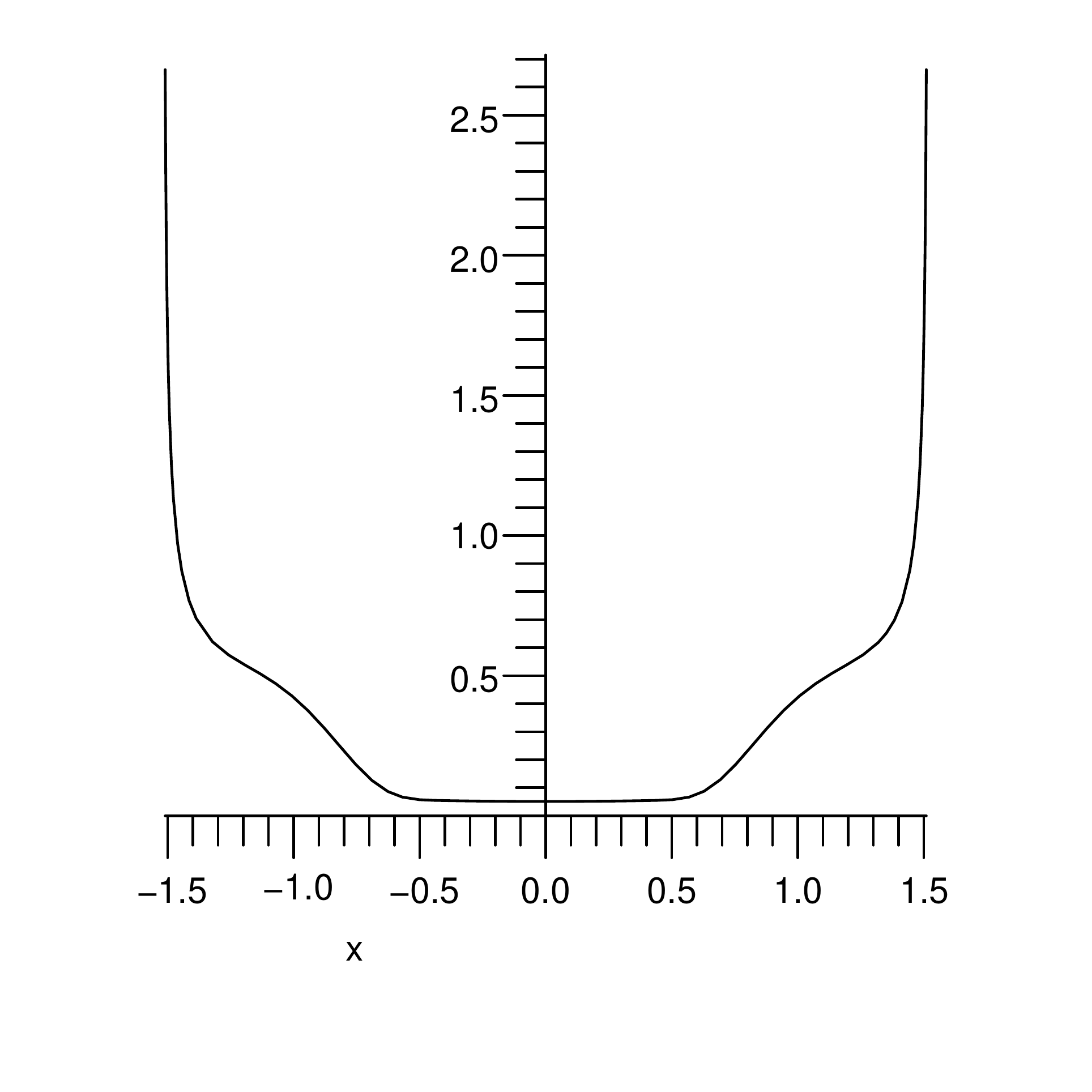} 
		\\[5mm] {\footnotesize \fig 
		Probability density of the $\sigma_x$-pointer 
  		for the tracial input state.}
	\end{minipage}
&
	\begin{minipage}[t]{3.6 cm}
		\includegraphics[width = 3.6 cm, viewport = 50 100 500 700]{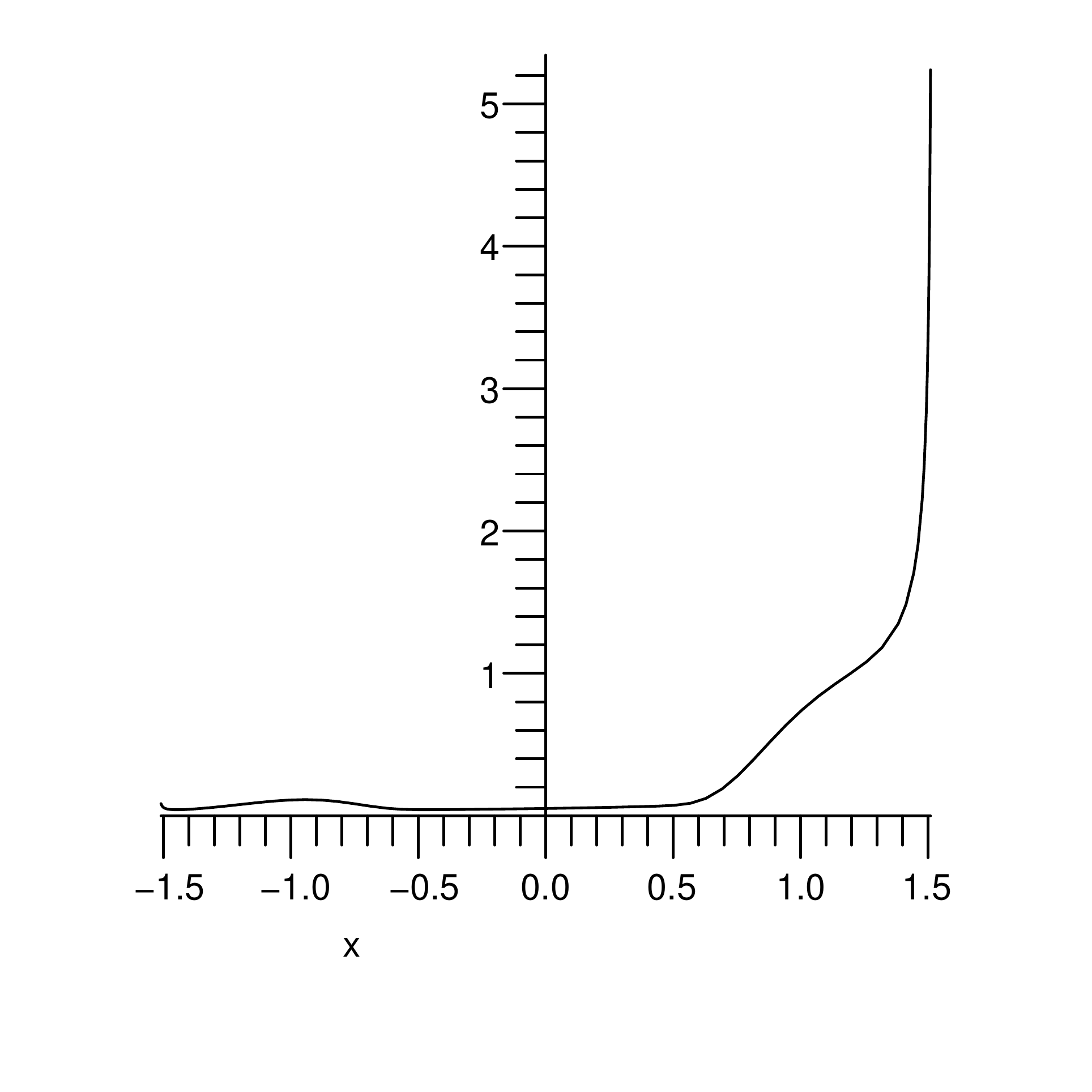}
		\\[5mm] {\footnotesize \fig Probability density of the $\sigma_x$-pointer 
		for the input state $|\rightarrow \, \rangle $.}
		\label{fig: r(x)100}
	\end{minipage}
\\
\end{tabular}
\end{center}

\subsection{Plots of \texorpdfstring{$\sigma_{z}$-}{}Measurement} \label{plotz}

We repeat this for the $\sigma_{z}$-pointer.  
By formula \ref{eq q}, the positive $\sigma_{z}$-eigenstate, the
tracial state and the negative $\sigma_{z}$-eigenstate 
lead to the distributions of the endpoint of the weighted 
path that are shown below:\vspace{-15mm}

\begin{center}
\begin{tabular}{l c l}
	\begin{minipage}[t]{3.6cm}
	  	\includegraphics[width = 3.6 cm, viewport = 50 100 500 700]
		{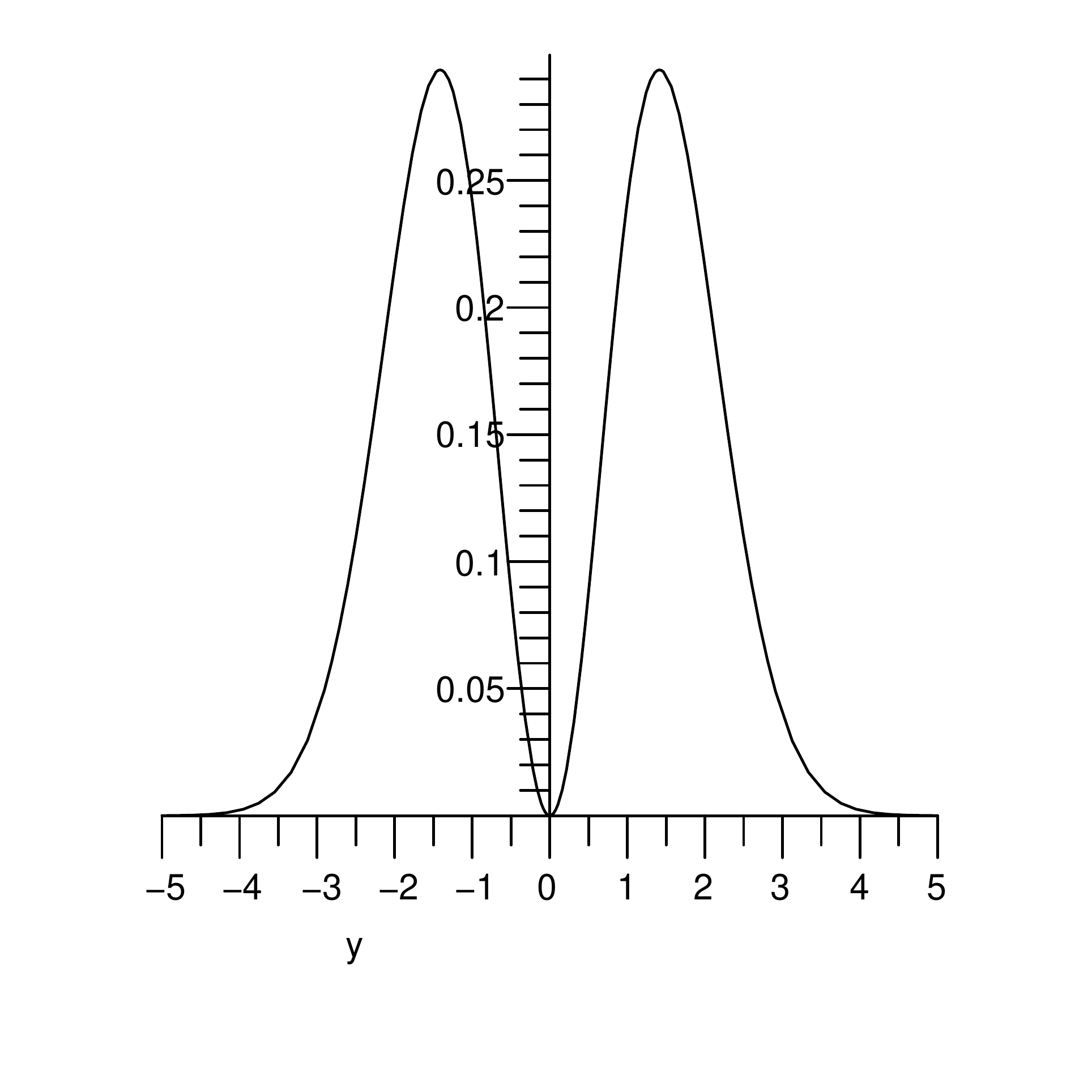} 
	 	\\[4mm] {\footnotesize \fig Probability density of the endpoint of the weighted path 
		for the input state $|\uparrow \, \rangle $.}
	\end{minipage}
&
	\begin{minipage}[t]{3.6cm}
		\includegraphics[width = 3.6 cm, viewport = 50 100 500 700]
		{qy000.pdf}
		\\[4mm]  {\footnotesize \fig Probability 
		density of the endpoint of the weighted path 
		for the tracial input state.}
	\end{minipage}
&
	\begin{minipage}[t]{3.6cm}
		  \includegraphics[width = 3.6 cm, viewport = 50 100 500 700]
		  {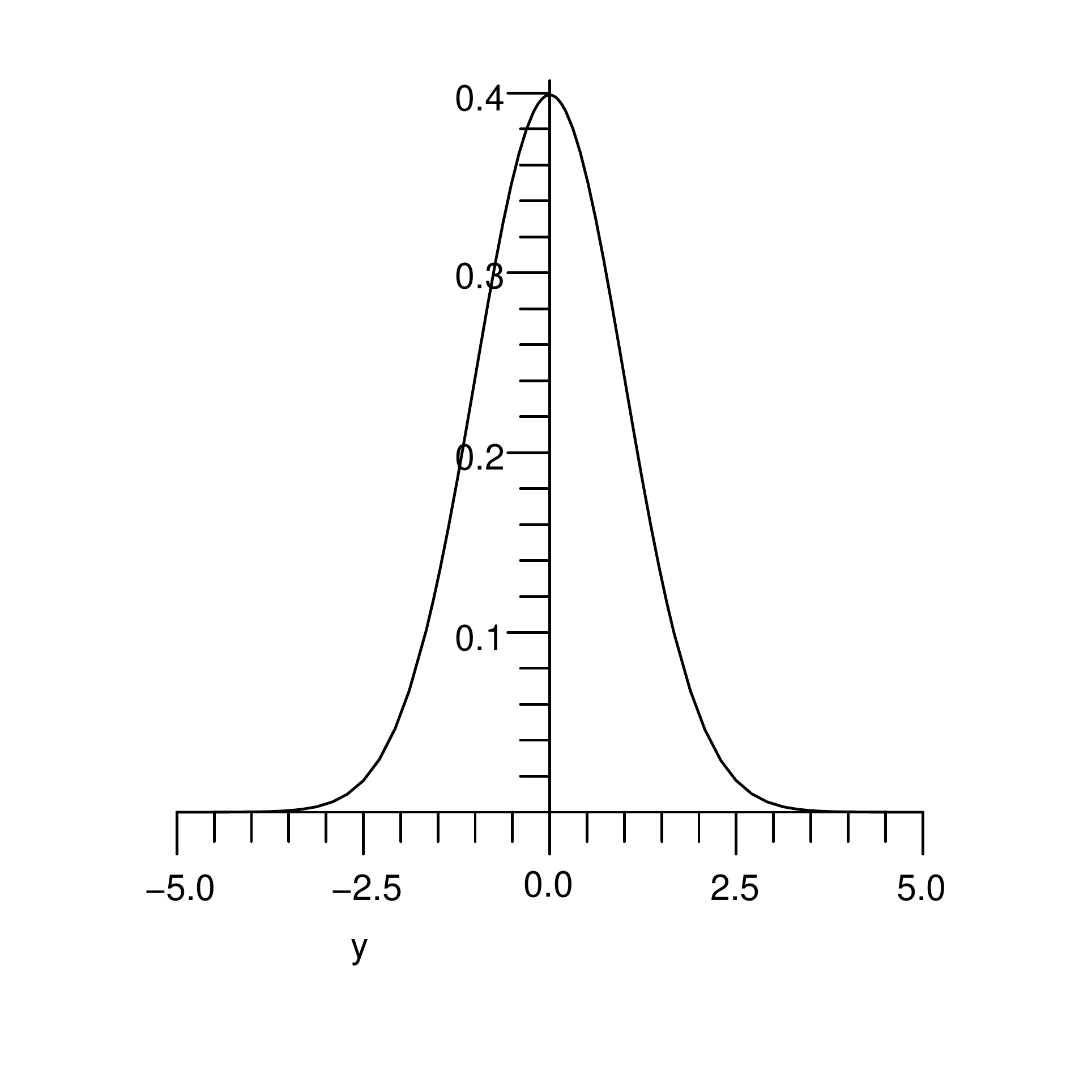}
		  \\[4mm] {\footnotesize \fig Probability 
		  density of the endpoint of the weighted path 
		  for the input state $|\downarrow \, \rangle $.}
	\end{minipage}
\\[2.8cm]
\begin{minipage}[t]{3.6cm}
		\mbox{}\\[-3mm]
		\includegraphics[width = 3.6 cm, viewport = 0 0 430 430 ]{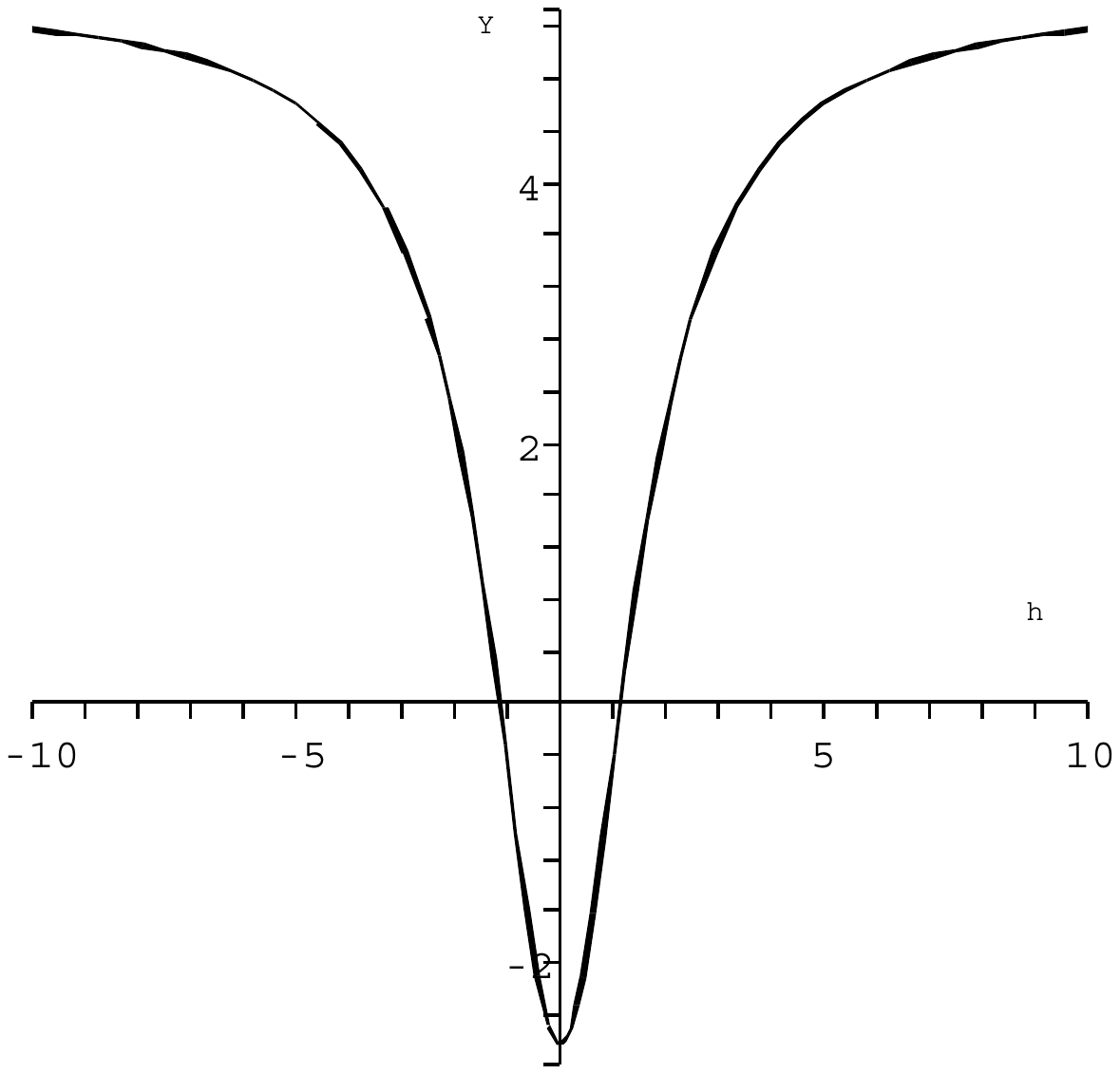}\\[5mm]
  		{\footnotesize \fig Pointer for $\sigma_z$.}\\
	\end{minipage}
	&
	\multicolumn{2}{p{7.6 cm}} 
	   {
	   	In order to estimate $\sigma_z$, we use the pointer
		of $\sigma_z$ illustrated here to the left. It is given by 
		\begin{equation}\label{htilde}\vspace{-1mm}
		\tilde{h}(x) = \frac{D_2}{x^2 + \delta} + D_3,
		\end{equation}
		with
		$\delta = 2.701$, $D_2 = -21.649$ and $D_3 = 5.391$.
		From formula \eqref{sofz}, we read off the probability 
		distributions of this
		pointer under the distributions of the endpoint of the weighted 
		path. 
		The probability distributions of the weighted path are illustrated above, 
		those of the pointer below.
 	  }	
\\[-8mm]
	\begin{minipage}[t]{3.6 cm}
		\includegraphics[width = 3.6 cm, viewport = 50 100 500 700]{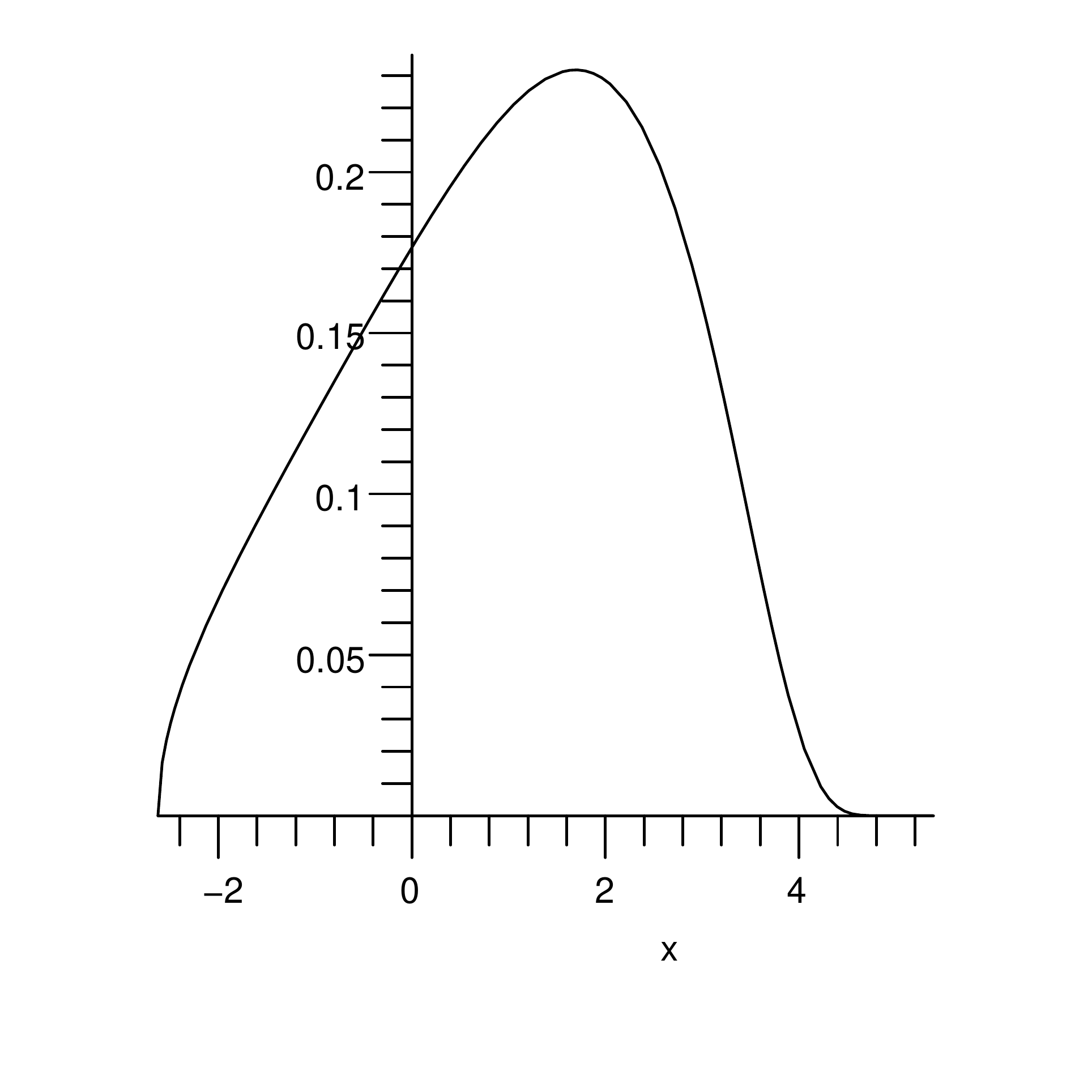}
 		\\[4mm] {\footnotesize \fig Probability density of the 
		$\sigma_z$-pointer  
  		for the input state $|\uparrow \, \rangle $.}
	\end{minipage}
&
	\begin{minipage}[t]{3.6 cm}
		\includegraphics[width = 3.6 cm, viewport = 50 100 500 700]{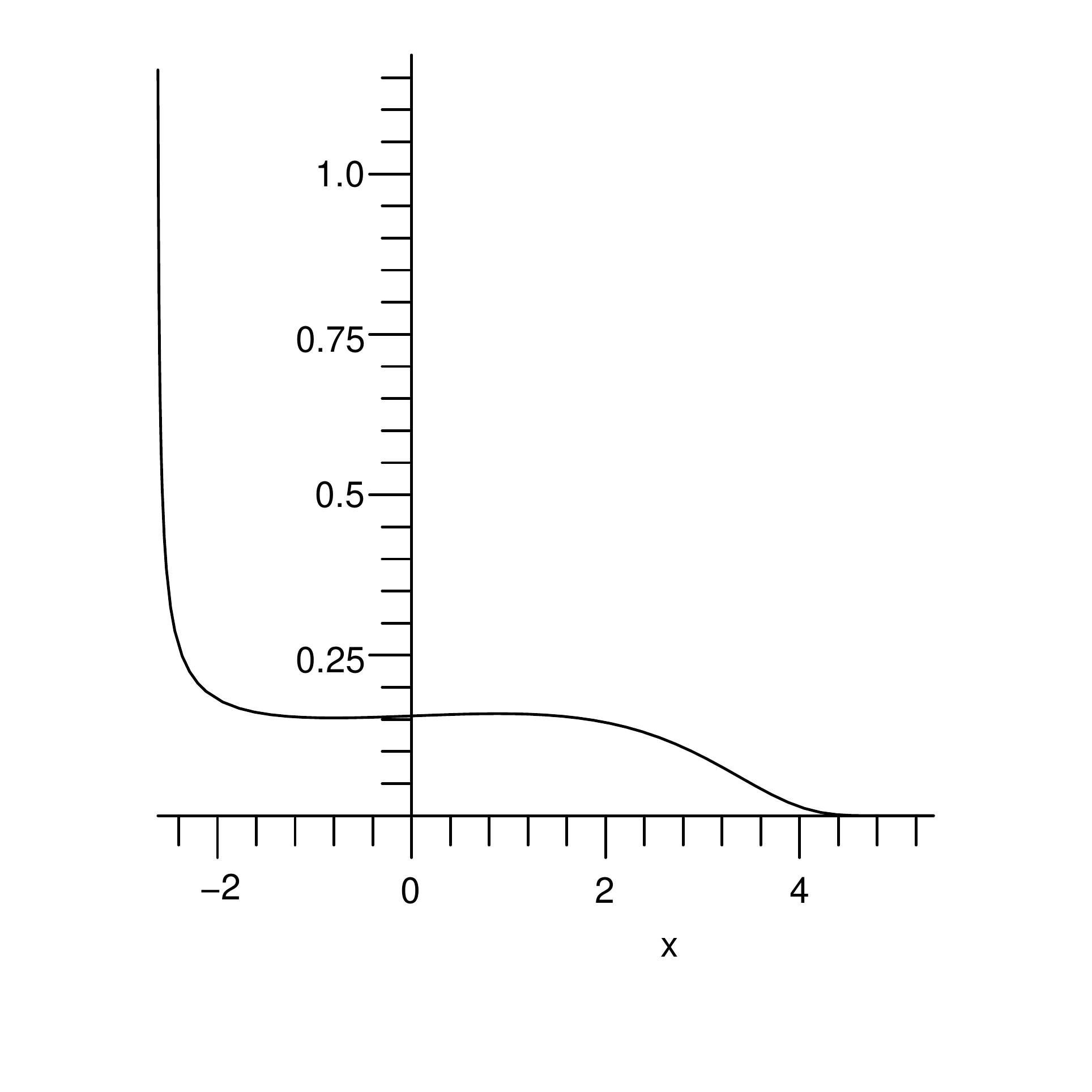}
		\\[4mm]  {\footnotesize \fig Probability density of the 
		$\sigma_z$-pointer  
		  for the tracial input state.}
	\end{minipage}
&
	\begin{minipage}[t]{3.6 cm}
		\includegraphics[width = 3.6 cm, viewport = 50 100 500 700]{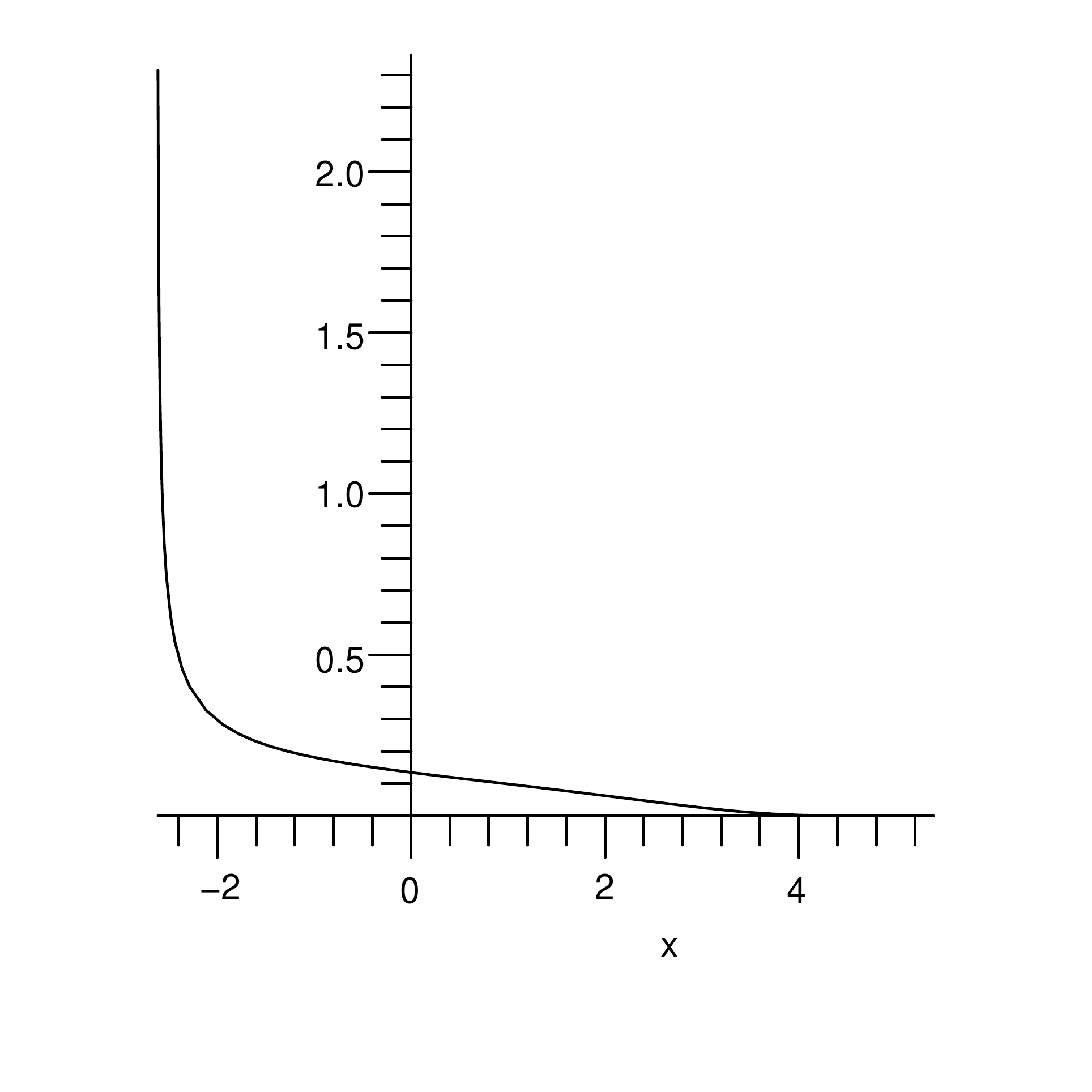}  
		\\[4mm] {\footnotesize \fig Probability density of the 
		$\sigma_z$-pointer 
		for the input state $|\downarrow \, \rangle $.}
	\end{minipage}
\\
\end{tabular}
\end{center}








\section{Discussion}\label{sec discussion}

In this chapter, we have investigated homodyne detection of 
spontaneous decay of a two-level atom into the 
electromagnetic field.    
We have seen how the photocurrent, besides carrying information 
on $\sigma_x$ (which is immediate from the innovations
term in the filtering equation), also carries information 
on $\sigma_z$. 
Homodyne detection can thus be viewed as 
a joint measurement of the noncommuting observables $\sigma_x$ 
and $\sigma_z$ in the initial state of the qubit, and we have identified
the optimal pointers for this procedure in equations \eqref{hstar} and
\eqref{htilde}.

One particular feature of the pointers we constructed might seem 
counter-intuitive at
first: they yield values outside $[-1,1]$ with nonzero probability.
This is a direct result of our requirement that the measurement be 
unbiased. Suppose, for example, that the input state is $| \uparrow \, \,\rangle$,
so that $\sigma_{z}$ has value 1. Since the photocurrent carries
information on $\sigma_{x}$ as well, its information on $\sigma_{z}$ is
certainly flawed, and will yield estimates $\sigma_{z} < 1$ at 
least some of the time. Unbiasedness then implies that also 
estimates $\sigma_{z} > 1$ must occur.   

On the other hand, an unbiased measurement will yield \emph{on average} 
the `true' value of
$\sigma_{z}$ for any possible input state. 
(Not just for the 3 possibilities
sketched on page \pageref{plotz}.) 
In repeated experiments, optimality of our pointers guarantees fast
convergence to these averages.  

Theorem \ref{thm joint} provides a theoretical bound for the quality of joint
measurement of $\sigma_{x}$ and $\sigma_{z}$. 
No conceivable measurement procedure can ever achieve 
$\Sigma \tilde{\Sigma} < 1$. 
It is now clear that
this bound cannot be met by way of homodyne detection:
a small part of the information extracted from the atom is simply lost
in this particular procedure.
Constructing the optimal pointers on the photocurrent does yield 
$\Sigma \tilde{\Sigma} = 1.056$ 
however, a figure much closer to the bound than the 4.437 provided by 
the na\"ive choice of \eqref{naive choice}.

\chapter[Optimal Estimation of Qubit States]{Optimal Estimation\\ of Qubit States 
}
\label{ch:OEQS}
\setsubdir{hoofdstukken/hoofdstuk5/}

%
%
%
%
In this chapter,
we propose an adaptive, two step strategy for the estimation of mixed qubit 
states. We show that the strategy is optimal in a local minimax sense for the 
trace norm distance as well as other locally quadratic figures of merit. Local 
minimax optimality means that given $n$ identical qubits, there exists no 
estimator which can perform better than the proposed estimator on a 
neighbourhood of size $n^{-1/2}$ of an arbitrary state. 
In particular, it is asymptotically Bayesian optimal for a large class 
of prior distributions.

We present a physical implementation of the optimal 
estimation strategy based on continuous time 
measurements in a field that couples with the qubits, not unlike
the setup in the previous chapter.

The crucial ingredient of the result is the 
concept of local asymptotic normality (or LAN) for qubits. 
This means that, for large $n$, the statistical model described by
$n$ identically prepared qubits is locally equivalent to a 
model with only a classical Gaussian distribution and a 
Gaussian state of a quantum harmonic oscillator.

The term `local' refers to a shrinking neighbourhood around a 
fixed state $\rho_{0}$. 
An essential result is that the neighbourhood radius can be chosen 
arbitrarily close to $n^{-1/4}$. This allows us to use a two step 
procedure by which we first localize the state within a smaller 
neighbourhood of radius $n^{-1/2+\epsilon}$, and then use LAN to 
perform optimal estimation. 

%
%
\section{Introduction}

State estimation is a central topic in quantum  statistical inference 
\cite{Ho,Helstrom,Barndorff-Nielsen&Gill&Jupp,Hayashi.editor}.
In broad terms the problem can be formulated as follows: 
given a quantum system prepared in an unknown state $\rho$, 
one would like to reconstruct the state by performing a measurement 
$M$ whose random result $X$ will be used to build an estimator 
$\hat{\rho}(X)$ of $\rho$. The quality of the measurement-estimator 
pair is
given by the {\it risk}
\begin{equation}\label{eq.risk}
\smash{
R_{\rho}(M, \hat{\rho} ) = \mathbb{E} \left(d(\hat{\rho} (X) , \rho)^{2}\right),
}
\end{equation}
where $d$ is a distance on the space of states, for instance the fidelity
distance or the trace norm, and the expectation is taken with respect to the
probability distribution $\mathbb{P}^{M}_{\rho}$ of $X$, when the 
measured system is in state $\rho$. Since the risk depends on the 
unknown state
$\rho$, one considers a global figure of merit by either averaging 
with respect to a prior distribution $\pi$ (Bayesian setup)
\begin{equation}\label{eq.risk.bayes}
R_{\pi}(M, \hat{\rho}) = \int \pi(d\rho) R_{\rho}(M, \hat{\rho} ),
\end{equation}
or by considering a maximum risk (pointwise or minimax setup)
\begin{equation}\label{eq.risk.max}
R_{\rm max} (M, \hat{\rho}) = {\sup}_{\rho}  \, R_{\rho}(M, \hat{\rho} ).
\end{equation}
An optimal procedure in either setup is one which achieves the minimum risk.

Typically, one measurement result does not provide enough information in
order to significantly narrow down on the true state $\rho$. Moreover, if
the measurement is ``informative'' then the state of the system after the
measurement will contain little or no  information about the initial state
(cf. chapter \ref{ch:UDHP}) and one needs to repeat the preparation and 
measurement procedure 
in order
to estimate the state with the desired accuracy.

It is then natural to consider a framework in which we are given a number
$n$ of identically prepared systems and look for estimators $\hat{\rho}_{n}$ 
which are optimal, or become optimal in the limit of large $n$. This problem 
is the quantum analogue of the classical statistical 
problem \cite{vanderVaart} of estimating a parameter
$\theta$ from independent identically distributed random variables
$X_{1},\dots , X_{n}$ with distribution $\smash{\mathbb{P}_{\theta}}$, and 
some of the methods developed in this paper are inspired by the 
classical theory.

Various state estimation problems have been investigated in the
literature and the techniques may be quite different depending on
a number of factors: the dimension of the density matrix, the number of
unknown parameters, the purity of the states, and the complexity of 
measurements over which one optimizes.
A short discussion on these issues can be found in section \ref{sec.estimation}.

In this paper we give an asymptotically optimal measurement strategy for qubit
states that is based on the technique of {\it local asymptotic normality} 
introduced 
in \cite{Guta&Kahn,Guta&Jencova}. The technique is a quantum 
generalization of Le Cam's classical statistical result \cite{LeCam}, 
and builds on previous work of Hayashi and Matsumoto 
\cite{Hayashi.conference,Hayashi&Matsumoto}.
We use an adaptive two stage procedure involving
continuous time  measurements, which could in principle be implemented in
practice. The idea of adaptive estimation methods, which has a long history 
in classical 
statistics, was introduced in the quantum set-up by 
\cite{Barndorff-Nielsen&Gill}, 
and was subsequently used in \cite{Gill&Massar,Hayashi,Hayashi&Matsumoto2}.  
The aim there is similar: one wants to first localize the state, and then  
perform a suitably tailored measurement which performs optimally around a 
given state.  A different adaptive technique was proposed independently by 
Nagaoka \cite{Nagaoka} and further developed in \cite{Fujiwara}.

%
In the first stage,
the spin components $\sigma_x$, $\sigma_y$ and $\sigma_z$
are measured separately on a small portion
$\tilde{n} \ll n$ of the
systems, and a rough estimator $\tilde{\rho}_{n}$ is constructed.
By standard statistical arguments (see lemma \ref{lemma.small.probability}) we deduce that with high probability, the true state
$\rho$ lies within a ball of radius slightly larger than $n^{-1/2}$, say 
$n^{-1/2 + \epsilon}$ with $\epsilon>0$, centered at
$\tilde{\rho}_{n}$. The purpose of the first stage is thus
to localize the state within a small neighbourhood 
(up to a unitary rotation), as illustrated in 
the figure below
using the Bloch sphere representation 
of qubit states.\\[-1 cm] 

\begin{center}
\begin{tabular}{p{8cm}}
\hspace{5.2 mm}
\includegraphics[width=6.6cm]{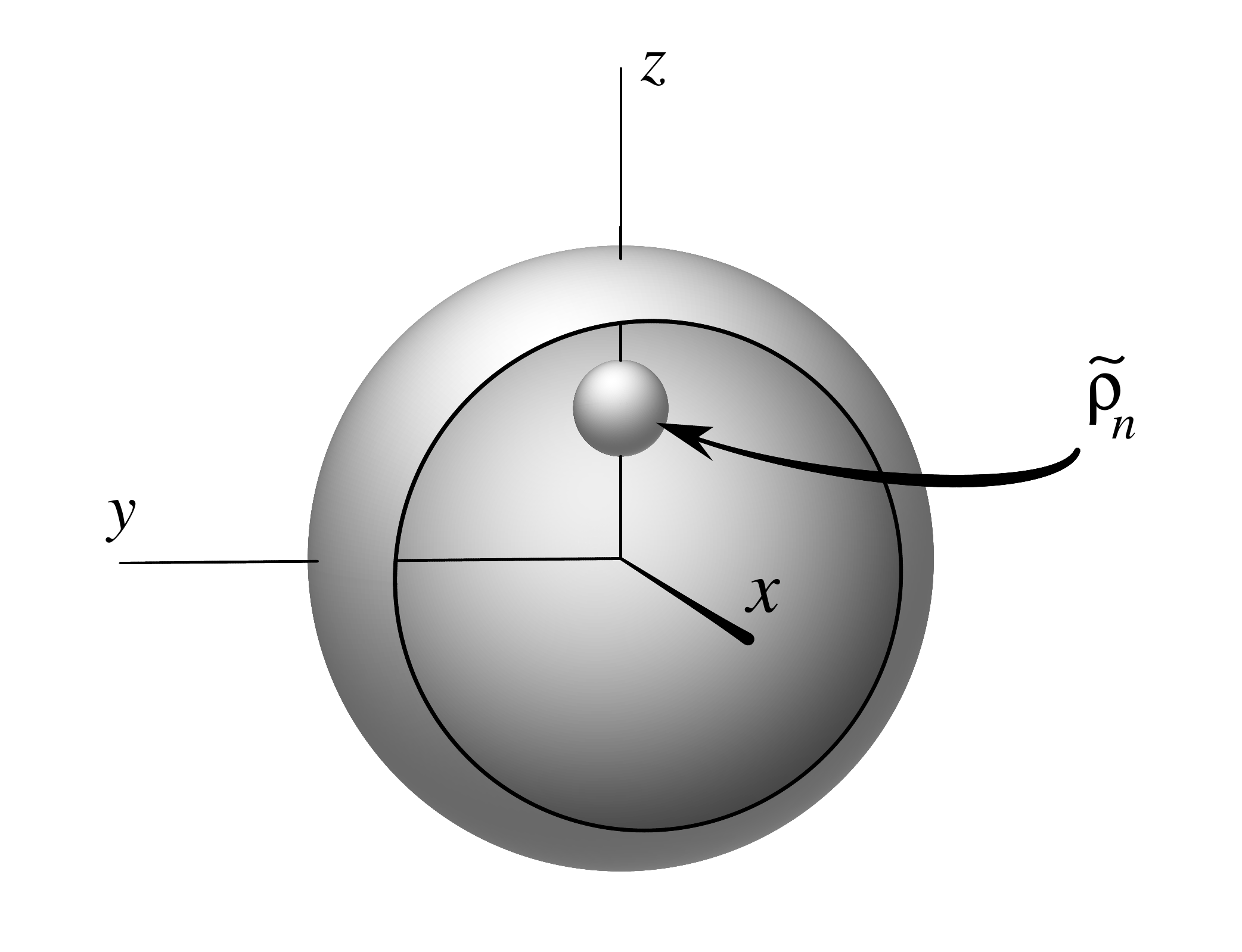}
\\[-0.3 cm]
{\small \fig 
After the first measurement stage, the state $\rho$ lies in a small ball 
centered at $\tilde{\rho}_{n}$.}
\end{tabular}\\[-0.5 cm]
\end{center}


This information is then used in the second stage, which is a
\emph{joint} measurement on the remaining $n-\tilde{n}$ systems.
This second measurement is implemented physically by two consecutive
couplings, each to a bosonic field. The qubits are first 
coupled to the field via a spontaneous emission interaction and a continuous 
time heterodyne detection measurement is performed in the field. This yields 
information on the eigenvectors of $\rho$. Then the interaction is changed, 
and a continuous time homodyne detection is performed in the field. 
This yields information on the eigenvalues of $\rho$.
 

We prove that the second stage of the measurement
is asymptotically optimal for all states in a ball of radius $n^{-1/2 + \eta}$  
around $\tilde{\rho}_{n}$. Here $\eta$ can be chosen to be bigger than 
$\epsilon>0$ implying that the two stage procedure as a whole is 
asymptotically optimal for any state, as depicted in 
the figure below.
\begin{center}
\begin{tabular}{p{8cm}}
\hspace{15 mm}
\includegraphics[width = 4.8cm]
{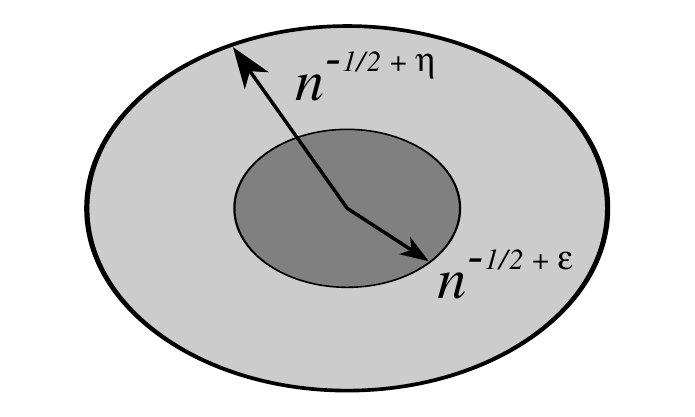}\\[-0.2 cm]
{\small \fig \setcounter{tempfig}{\value{figuur}}
The smaller domain is the localization region of the first step.
The second stage estimator is optimal for all states in the bigger domain.}
\end{tabular}
\end{center}


The optimality of the second stage
relies heavily on the principle of \emph{local asymptotic normality} or  LAN, 
see \cite{vanderVaart}, which we will briefly explain below,
and in particular on the fact that it holds in a ball of radius 
$n^{-1/2 + \eta}$
around $\smash{\tilde{\rho}_{n}}$ rather than just $\smash{n^{-1/2}}$, as was the case 
in \cite{Guta&Kahn}. 

Let $\rho_{0}$ be a fixed state.
We parameterize the neighbouring states as $\rho_{{\bf u}/\sqrt{n}}$,
where $\smash{{\bf u}=(u_{x},u_{y},u_{z})\in\mathbb{R}^{3}}$
is a certain set of local parameters around $\rho_0$. Then LAN entails 
that the joint state
$\smash[t]{\rho^{\bf u}_{n}: =\rho_{{\bf u}/\sqrt{n}}^{\otimes n}}$
of $n$ identical qubits
converges in the limit $n \rightarrow \infty$
to a Gaussian state of the form $N^{\bf u}\otimes\phi^{\bf u}$,
in a sense explained in theorem \ref{th.qlan}. By $N^{\bf u}$ we denote a 
{\it classical} one-dimensional
normal distribution centered at $u_{z}$. The second term $\phi^{\bf u}$ 
is a Gaussian state of a harmonic
oscillator, i.e. a displaced thermal equilibrium state with
displacement proportional to $(u_{x},u_{y})$. We thus have the convergence
$$
\rho^{\bf u}_{n}\leadsto N^{\bf u}\otimes \phi^{\bf u},
$$
to a much simpler family of classical -- quantum states for which we 
know how
to optimally estimate the parameter {\bf u} \cite{Ho,Yuen&Lax}.

The idea of approximating a sequence of statistical experiments by a
Gaussian one goes back to Wald \cite{Wald2}, and was subsequently developed 
by Le Cam \cite{LeCam} who coined the term local asymptotic normality. 
In quantum statistics  the first ideas in the direction of local asymptotic 
normality for d-dimensional 
states appeared in the Japanese paper \cite{Hayashi.japanese}, as well as 
\cite{Hayashi.conference} and were subsequently developed in 
\cite{Hayashi&Matsumoto}. In theorem \ref{th.qlan} we strengthen 
these results for the case of qubits, by proving a strong version 
of LAN in the spirit of Le Cam's pioneering work. We then exploit 
this result to prove optimality of the second stage. A different 
approach to local asymptotic normality has been developed in 
\cite{Guta&Jencova} to which we refer for a more general exposition 
on the theory of quantum statistical models. A short discussion on 
the relation between the two approaches is given in the remark 
following theorem \ref{th.qlan}.


\begin{center}
\begin{tabular}{p{8cm}}
\hspace{10 mm}
\includegraphics[width = 6cm]
{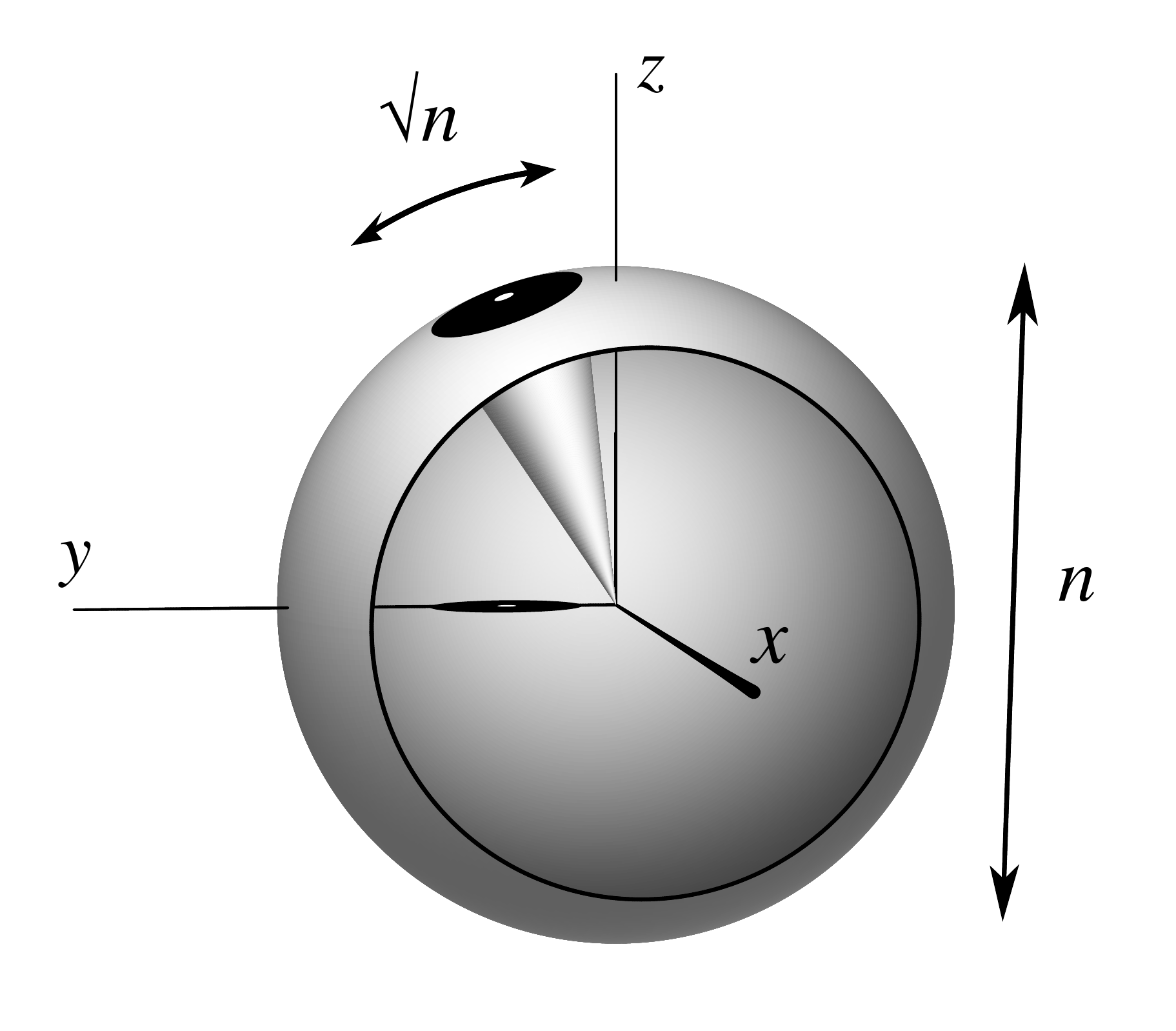}\\[-0.2cm]
{\small \fig 
Total spin representation 
of the state of $n\gg 1$  spins. The quantum fluctuations of the $x$ and 
$y$ spin directions coincide with those of a coherent state of a harmonic 
oscillator.}
\end{tabular}
\end{center}

From the physics perspective, our results put on a more rigorous basis the 
treatment 
of collective states of many identical spins, the keyword here being 
{\it coherent spin states} \cite{Holtz}. Indeed, it has been known 
since Dyson \cite{Dyson1} that $n$ spin-$\frac{1}{2}$ particles 
prepared in the spin up state $|\!\uparrow\,\rangle^{\otimes n}$ behave 
asymptotically as the ground state of a quantum oscillator, when 
considering the fluctuations of properly normalized total spin components 
in the directions orthogonal to $z$. We extend this to spin directions 
making an ``angle'' of order $n^{-1/2 + \eta}$ with the $z$ axis, as 
illustrated in the figure above, as well as to mixed states. 
We believe that a similar approach can be followed in the case of spin 
squeezed states and continuous time measurements with 
feedback control \cite{GSM04}.
%

 In theorem \ref{th.unitary.evolution} we prove a dynamical version of LAN.
 The trajectory in time of the joint state of the qubits together with the field  converges for large  $n$ to the corresponding trajectory of the joint state of the oscillator and field.
 In other words, time evolution preserves local asymptotic normality.
 This insures that for large $n$  the state of the qubits ``leaks'' into a Gaussian state of the field, providing a concrete implementation of the convergence to the limit Gaussian experiment.

The punch line of the paper is theorem \ref{automobiel} which says that the
estimator $\hat\rho_{n}$ is optimal in local minimax sense, which is the
modern statistical formulation of optimality in the frequentist setup \cite{vanderVaart}.
Also, its asymptotic risk is calculated explicitly.


The paper is structured as follows: in section \ref{sec.estimation},
we show that the first stage of the measurement sufficiently
localizes the state. In section \ref{sec.lan}, we prove
that LAN holds with radius of validity $n^{-1/2 + \eta}$, and we bound
its rate of convergence.
Sections \ref{sec.timevolf} and \ref{finmes} are concerned with the 
second stage of the measurement, i.e. with the coupling to the bosonic 
field and the continuous time
field-measurements. Finally, in section \ref{sec.endresult},
asymptotic optimality of the estimation scheme is proven.

The technical details of the proofs are relegated to the appendices in order to give the reader a more direct access  to the ideas and results.

\section{State Estimation}\label{sec.estimation}

In this section we introduce the reader to a few general aspects of quantum
state estimation after which we concentrate on the qubit case.

State estimation is a generic name for a variety of results which may be classified
according to the dimension of the parameter space, the kind or family of states 
to be estimated and the preferred estimation method. For an introduction to 
quantum statistical 
inference we refer to the books by Helstrom \cite{Helstrom} and Holevo \cite{Ho} 
and the more recent review paper \cite{Barndorff-Nielsen&Gill&Jupp}. The 
collection \cite{Hayashi.editor} is a good reference on quantum statistical 
problems, with many important contributions by the Japanese school.

For the purpose of this chapter, any quantum state representing a particular
preparation of a quantum system, is described by a density matrix 
(positive self-adjoint operator of trace one) on the Hilbert space 
$\mathcal{H}$ associated to the system. The algebra of observables 
is $\mathcal{B}(\mathcal{H})$, and the expectation of an observable 
$a\in\mathcal{B}(\mathcal{H})$ with respect to the state $\rho$ is 
$\tr(\rho a)$. A measurement $M$ with outcomes in a measurable space 
$(\mathcal{X},\Sigma)$ is completely determined by a $\sigma$-additive 
collection of positive self-adjoint operators $M(A)$ on $\mathcal{H}$, 
where $A$ is an event in $\Sigma$. This collection is called a positive 
operator valued measure (cf.~def.~\ref{povem}). 
The distribution of the results $X$ when the 
system is in state $\rho$ is given 
by $P_{\rho}(A) = \tr(\rho M(A))$.  

We are given $n$ systems identically prepared in state $\rho$ and we are 
allowed to perform a measurement $M_{n}$ whose outcome is the estimator
 $\hat{\rho}_{n}$ as discussed in the introduction.
 
The dimension of the density matrix may be finite, such as in the case of 
qubits or d-levels atoms, or infinite as in the case of the state of a 
monochromatic beam of
light. In the finite or parametric case one expects that the risk converges 
to zero as
$n^{-1}$ and the optimal measurement-estimator sequence $(M_{n},\hat{\rho}_{n})$
achieves the {\it best constant} in front of the $n^{-1}$ factor.
In the non-parametric case the rates of convergence are in general slower that
$n^{-1}$ because one has to simultaneously estimate an infinite number of 
matrix elements, each with rate $n^{-1}$.
An important example of such an estimation technique is that of quantum homodyne
tomography in quantum optics \cite{Vogel&Risken}. This allows the estimation 
with arbitrary precision \cite{D'Ariano.2, D'Ariano.3,Leonhardt.Munroe} of the 
whole density matrix  of a monochromatic beam of light by repeatedly measuring 
a sufficiently large number of identically prepared 
beams \cite{Smithey,Breitenbach&Schiller&Mlynek,Zavatta}. 
In \cite{Artiles&Guta&Gill,Butucea&Guta&Artiles} it is shown how to 
formulate the
problem of estimating infinite dimensional states without the need for 
choosing a cut-off in the dimension of the density matrix, and how to 
construct optimal minimax estimators of the Wigner function for a class 
of ``smooth'' states.

If we have some prior knowledge about the preparation procedure, we may encode
this by parameterizing the possible states as $\rho=\rho_{\theta}$ with
$\theta\in\Theta$ some unknown parameter. The problem is then to
estimate
$\theta$ optimally with respect to a distance function on $\Theta$.

Indeed, one of the main problems in the finite dimensional case is to find optimal estimation procedures for a given family of states. It is known that if the state $\rho$ is pure or belongs to a one parameter family, then separate measurements achieve
the optimal rate of the class of joint measurements \cite{Matsumoto}.
However for multi-dimensional families of mixed states this is no longer the case and joint measurements perform strictly better than separate ones \cite{Gill&Massar}.

In the Bayesian setup, one optimizes $R_{\pi}(M_{n},\hat{\rho}_{n})$ for some prior 
distribution $\pi$. We refer to 
\cite{Jones,Massar&Popescu,Latorre&Pascual&Tarrach,Fisher&Kienle&Freyberger,
Wunderlich,Bagan&Baig&Tapia,Narnhofer,Bagan&Monras&Tapia} for the pure state case, and to
 \cite{Cirac,Vidal,Mack, Keyl&Werner,Bagan&Baig&Tapia&Rodriguez,Sommers}, and \cite{Bagan&Gill} for the
mixed state case. The methods used here are based on group theory and can be applied 
only to invariant prior distributions and certain distance functions.
In particular, the optimal covariant measurement in the case of completely unknown
qubit states was found in \cite{Bagan&Gill,Hayashi&Matsumoto} but it has 
the drawback that it does not give any clue as to how it can be
implemented in a real experiment.

In the pointwise approach 
\cite{Hayashi,Hayashi&Matsumoto2,Gill&Massar,Barndorff-Nielsen&Gill, 
Fujiwara&Nagaoka,Matsumoto,Barndorff-Nielsen&Gill&Jupp,Hayashi&Matsumoto} 
one tries to minimize the risk for {\it each} unknown state $\rho$. 
As the optimal measurement-estimator pair cannot depend on 
the state itself, one optimizes the maximum risk 
$R_{\rm max}(M_{n},\hat{\rho}_{n})$, (see \eqref{eq.risk.max}), 
or a local version of this which will be defined shortly. 
The advantage of the pointwise approach is that it can be applied to 
arbitrary families of states and a large class of loss functions provided 
that they are locally quadratic in the chosen parameters.  The underlying 
philosophy is that as the number $n$ of states is sufficiently large, the 
problem ceases to be global and becomes a local one as the error in 
estimating the state parameters is of the order $n^{-1/2}$. 

The Bayesian and pointwise approaches can be compared
\cite{Gillunpub}, and in fact for large $n$ the prior distribution $\pi$ of the Bayesian approach becomes increasingly irrelevant and the optimal Bayesian estimator becomes asymptotically optimal in the minimax sense and vice versa.

\subsection{Qubit State Estimation: the Localization Principle}

Let us now pass to the quantum statistical model which will be the object of 
our investigations. Let $\rho \in M_{2}(\mathbb{C})$ be an arbitrary density 
matrix describing the state of a qubit. Given $n$ identically prepared qubits 
with joint state $\rho^{\otimes n}$, we would like to optimally estimate 
$\rho$ based on the result of a properly
chosen joint measurement $M_{n}$. For simplicity of the exposition we assume 
that
the outcome of the measurement is an estimator 
$\hat \rho_{n}\in M_{2}(\mathbb{C})$. In practice however, the result $X$ 
may belong to a complicated measure space (in our case the space of continuous 
time paths) and the estimator is a function of the ``raw'' data 
$\hat{\rho}_{n} := \hat{\rho}_{n}(X)$. The quality of the estimator at the 
state $\rho$ is quantified by the risk
$$
R_{\rho}(M_{n}, \hat{\rho}_{n}) := \mathbb{E}_{\rho} (d (\rho, \hat{\rho}_{n})^{2}),
$$
where $d$ is a distance between states. The above expectation is taken with 
respect to
the distribution $P_{\rho}(dx):=\tr(\rho M(dx))$ of the measurement results, 
where $M(dx)$ represents the associated positive operator valued measure of 
the measurement $M$. 
In our exposition $d$ will be the trace norm (see p.~\pageref{treesafstand})
$$
\| \rho_{1} - \rho_{2}\|_{1} :=  \tr(|\rho_{1} - \rho_{2}|),
$$
but similar results can be obtained using the fidelity distance. The aim is to 
find a sequence of measurements and estimators $(M_{n}, \hat{\rho}_{n})$ which 
is asymptotically optimal in the {\it local minimax} sense: for any given 
$\rho_{0}$
$$
\limsup_{n\to \infty}\sup_{\|\rho - \rho_0\|_{1} \leq n^{-1/2 + \epsilon}}  n R_{\rho}(M_{n}, \hat{\rho}_{n}) \leq
\limsup_{n\to \infty}\sup_{\|\rho - \rho_0\|_{1} \leq n^{-1/2 + \epsilon}}  n R_{\rho}(N_{n}, \check\rho_{n}),
$$
for any other sequence of measurement-estimator pairs $(N_{n},\check{\rho}_{n})$. The factor $n$ is inserted because typically $R_{\rho}(M_{n}, \hat{\rho}_{n})$ is of the order $1/n$ and the optimization is about obtaining the smallest constant factor possible.
The inequality says that one cannot find an estimator which performs better that
$\hat\rho_{n}$ over a ball of size $n^{-1/2 + \epsilon}$ centered at $\rho_{0}$,
even if one has the knowledge that  the state $\rho$ belongs to that ball! 

Here, and elsewhere in the paper $\epsilon$ will appear in different contexts, as a generic strictly positive number and will be chosen to be sufficiently small for each specific use. At places where such notation may be confusing we will use additional symbols to denote small constants.

As set forth in the Introduction,
our measurement procedure consists of two steps. The first one is to perform
separate measurements of $\sigma_x$, $\sigma_y$ and $\sigma_z$
on a fraction $\tilde{n}=\tilde{n}(n)$ of the systems.
In this way we
obtain a rough estimate $\tilde \rho_{n}$ of the true state $\rho$
which lies in a local neighbourhood around $\rho$ with
high probability.
The second step uses the information obtained in the first step
to perform a measurement which is optimal precisely for the states
in this local neighbourhood.
The second step ensures optimality and requires more sophisticated techniques
inspired by the theory of local asymptotic normality for qubit states \cite{Guta&Kahn}.
We begin by showing that the first step amounts to the fact that, without loss of
generality, we may assume that the unknown state is in a local neighbourhood of a known state. This may serve also as an a posteriori justification of the definition of local
minimax  optimality.
\begin{lemma}\label{lemma.small.probability}
Let $M_{i}$ denote the measurement of the $\sigma_{i}$ spin component of a
qubit with $i=x,y,z$. We perform each of the measurements $M_{i}$ separately
on $\tilde{n}/3$ identically prepared qubits and define
$$
\tilde\rho_{n} = \frac{1}{2}(\mathbf{1} + \tilde{\mathbf{r}} \sigma),\qquad {\rm if~}\quad |\tilde{r}|\leq 1,
$$
where $\tilde{\mathbf{r}} =(\tilde r_{x}, \tilde r_{y}, \tilde r_{z})$ is the vector average
of the measured components. If $| \tilde{r} |>1$ then we define $\tilde\rho_{n}$ as the state which has the smallest trace distance to the right hand side expression. 
Then for all
$\epsilon \in [0,2]$,
we have
$$
\mathbb{P} \left(\| \tilde\rho_{n} -\rho \|_{1}^{2} >
3 n^{2 \epsilon-1} \right) \leq
6 \exp(-\half \tilde{n} n^{2 \epsilon - 1} ), \qquad \forall \rho.
$$
Furthermore, for any
$0< \kappa < \epsilon/2$, if $\tilde{n}=n^{1-\kappa}$, the contribution to the risk 
$\mathbb{E} (\| \tilde{\rho}_{n}-\rho\|_{1}^{2})$ brought by the event 
$E= [\, \| \tilde\rho_{n} -\rho \|_{1} > \sqrt{3} n^{-1/2 +\epsilon}\,] $ satisfies
$$
\mathbb{E} \left(\, \| \tilde{\rho}_{n} - \rho \|_{1}^{2}\, \chi_{E}\, \right) \leq 24 \exp(-\half n^{2 \epsilon - \kappa}) = \mathcal{O}(1).
$$
\end{lemma}

\proof
For each spin component $\sigma_{i}$ we obtain i.i.d coin tosses $X_{i}$ with distribution
$\mathbb P(X_{i} = \pm 1) =(1\pm r_{i})/2$ and average $r_{i}$.

Hoeffding's inequality \cite{vanderVaart&Wellner} then states that for all real numbers $c > 0$, we have
$\mathbb{P}( |X_i - \tilde{X} |^2 > c ) \leq 2 \exp(- \half \tilde{n}c ) $.
By using this inequality  three times with
$c = n^{2 \epsilon -1}$,
once for each component, we get
$$
\mathbb{P} \left( \sum_{1}^{3} | \tilde r_{i} -r_{i}|^{2} > 3
n^{2 \epsilon-1}  \right) \leq
6 \exp(-\half \tilde{n} n^{2 \epsilon - 1} )
  \qquad \forall \rho,
$$
which implies the statement for the norm distance since 
$\|\tilde{\rho}_{n} -\rho \|_{1}^{2} =\sum_{i} |\tilde{r}_{i} -r_{i}|^{2}$.
The bound on conditional risk follows from the previous bound and the fact that
$\| \rho- \tilde{\rho}_{n} \|_{1}^{2} \leq 4$.
\qed

\noindent In the second step of the measurement procedure we rotate the
remaining $n-\tilde{n}$ qubits such that after rotation the vector $\tilde r$ is parallel to the $z$-axis. Afterwards,
we couple the systems to the field and perform certain measurements in the field which will
determine the final estimator $\hat{\rho}_{n}$.
The details of this second step are given in sections \ref{sec.timevolf} and
\ref{finmes}, but at this moment we can already prove that the effect of errors in the the first stage of the measurement is
asymptotically negligible compared to the risk of the second estimator. Indeed
by lemma
\ref{lemma.small.probability} we get that if
$\tilde{n}= n^{1 - \kappa}$,
then the probability that the first stage gives a
``wrong'' estimator (one which lies outside the local neighbourhood of the  true state)  is of the order
$\exp(-\half n^{2 \epsilon - \kappa})$ and so is the risk contribution. As the typical risk of
estimation is of the order $1/n$, we see that the first step is practically ``always''
placing the estimator in a neighbourhood of order
$n^{- 1/2 + \epsilon}$
of the true state $\rho$, as shown in figure \arabic{tempfig}.
In the next section we will show that for such neighbourhoods, the state of the remaining
$n-\tilde{n}$ systems behaves asymptotically as a Gaussian state. This will allow us to devise an optimal
measurement scheme for qubits based on the optimal measurement for Gaussian states.

\section{Local Asymptotic Normality}\label{sec.lan}

The optimality of the second stage of the measurement
relies on the concept of local asymptotic normality \cite{vanderVaart,Guta&Kahn}.
After a short introduction, we will prove that LAN holds for the qubit
case, with radius of validity $n^{-1/2 + \eta}$ for all $\eta \in [0,1/4)$.
We will also show that its rate of convergence is
$\mathcal{O}(n^{-1/4 + \eta + \epsilon})$ for arbitrarily
small $\epsilon$.

\subsection{Introduction to LAN and some Definitions}

Let $\rho_{0}$ be a fixed state, 
which by rotational symmetry can be chosen of the form
\begin{equation}
\label{rho0}
\rho_{0} =\left(
\begin{array}{cc}
\mu & 0\\
0 & 1-\mu
\end{array}
\right),
\end{equation}
for a given $\frac{1}{2}<\mu< 1$.
We parameterize the neighbouring states as $\rho_{{\bf u}/\sqrt{n}}$ where
${\bf u}=(u_{x},u_{y},u_{z})\in\mathbb{R}^{3}$ such that the first two
components  account for unitary rotations around $\rho_{0}$, while the
third one describes the change in eigenvalues
  \begin{equation}\label{eq.family}
\rho_{\bf v} := U\left({\bf v}\right)
\left(
\begin{array}{cc}
\mu + v_{z} & 0\\
0 & 1-\mu -v_{z}
\end{array}
\right)
U\left({\bf v} \right)^{\dagger}  , 
\end{equation}
with unitary
$U({\bf v}):= \exp(i(v_{x}\sigma_{x}  + v_{y}\sigma_{y}))$.
The ``local parameter'' ${\bf u}$ should be thought of, as having a
bounded range in $\mathbb{R}^{3}$ or may even ``grow slowly''
as $\|{\bf u}\|\leq n^{\eta}$.

Then, for large $n$, the joint state
$\rho^{\bf u}_{n}: =\rho_{{\bf u}/\sqrt{n}}^{\otimes n} $ of $n$ identical qubits
approaches a Gaussian state of the form $N^{\bf u}\otimes\phi^{\bf u}$ with
the parameter ${\bf u}$ appearing solely in the average of the two Gaussians.
By $N^{\bf u}$ we denote a {\it classical} one-dimensional
normal distribution centered at $u_{z}$ which relays
information about the eigenvalues of $\rho_{{\bf u}/\sqrt{n}}$.
The second term $\phi^{\bf u}$ is a Gaussian state of a harmonic
oscillator which is a displaced thermal equilibrium state with
displacement proportional to $(u_{x},u_{y})$.
It contains information on the eigenvectors of $\rho_{{\bf u}/\sqrt{n}}$.
We thus have the convergence
$$
\rho^{\bf u}_{n}\leadsto N^{\bf u}\otimes \phi^{\bf u},
$$
to a much simpler family of classical - quantum states for which we know how
to optimally estimate the parameter {\bf u}. The asymptotic splitting into
a classical estimation problem for eigenvalues and a quantum one for the
eigenbasis has been also noticed in \cite{Bagan&Gill} and in
\cite{Hayashi&Matsumoto}, the latter coming pretty close to our
formulation of local asymptotic normality.

The precise meaning of the convergence is given in theorem \ref{th.qlan}
below.
In short, there exist quantum channels (completely positive maps) $T_{n}$ which map the states
\smash{$\rho_{{\bf u}/\sqrt{n}}^{\otimes n}$}
into $N^{\bf u}\otimes\phi^{\bf u}$ with vanishing error in trace norm
distance, and uniformly over the local parameters ${\bf u}$. From the
statistical point of view the convergence implies
that a statistical decision problem concerning the model $\rho^{\bf u}_{n}$
can be mapped into a similar problem for the model
$N^{\bf u}\otimes\phi^{\bf u}$ such that the optimal solution for the
latter  can be translated into an asymptotically optimal solution
for the former. In our case the problem of estimating the state
$\rho$ turns into that of estimating the local parameter ${\bf u}$ around
the first stage estimator
$\tilde{\rho}_{n}$ playing the role of $\rho_{0}$. For the family of
displaced Gaussian states it is well known that the optimal estimation
of the displacement is achieved by the heterodyne detection
\cite{Ho,Yuen&Lax}, while for the classical part it sufficient to
take the observation as best estimator. Hence the second step will
give an optimal estimator $\hat{\bf u}$ of ${\bf u}$ and an optimal
estimator of the initial qubit state
$\hat{\rho}_{n}:= \rho_{\hat{\bf u}/\sqrt{n}}$. The precise result
is formulated in theorem \ref{automobiel}



\subsection{Convergence to the Gaussian Model}

We describe the state $N^{\bf u} \otimes \phi^{\bf u}$ in more detail.
$N^{\bf u}$ is simply the classical Gaussian distribution
\begin{equation} \label{gausje}
N^{\bf u}:=N(u_{z}, \mu(1-\mu)),
\end{equation}
with mean $u_z$ and variance $\mu(1-\mu)$.

The state $\phi^{\bf u}$ is a density matrix on
$\mathcal{H} = \mathcal{F}(\mathbb{C})$,
the representation space of the harmonic oscillator.
In general, for any Hilbert space $\mathfrak{h}$,
the {\it Fock space} over $\mathfrak{h}$ is defined 
(see sec.~\ref{focksec})
as
\begin{equation}\label{eq.Fock}
\mathcal{F}(\mathfrak{h}) := \bigoplus_{n=0}^{\infty}
\mathfrak{h} \otimes_{s}\dots
\otimes_{s} \mathfrak{h},
\end{equation}
with $\otimes_{s}$ denoting the symmetric tensor product.
Thus $\mathcal{F}(\mathbb{C})$ is the simplest example of a Fock space.
Let
\begin{equation}
\label{phi0}
\phi := (1-p) \sum_{k=0} p^{k} |k\rangle \langle k|
\end{equation}
be a thermal equilibrium state with $|k\rangle$ denoting the $k$-th 
energy level of
the oscillator and $p=\frac{1-\mu}{\mu}<1$.
For every $\alpha\in\mathbb{C}$ define the
displaced thermal state
\begin{equation*}
\label{phiu}
\phi(\alpha):= D( \alpha)  \,\phi\, D(-\alpha),
\end{equation*}
where $D(\alpha): = \exp( \alpha a^{\dagger} - \bar \alpha a)$
is the displacement operator (i.e.~the Weyl operator $W(\alpha)$,
cf.~sec.~\ref{hofock}), mapping
the vacuum vector $\ket{0}$ to the coherent vector
$$
| \alpha \rangle = \exp( - |\alpha|^2 \!/ 2)\sum_{k=0}^{\infty}
\frac{\alpha^k}{\sqrt{k!}} |k\rangle. 
$$
Here $a^{\dagger}$ and $a$ are the creation and annihilation operators on
$\mathcal{F}(\mathbb{C})$,
satisfying $[a, a^{\dagger}] = \mathbf{1}$.
The family $\phi^{\bf u}$ of states in which we are interested
is given by
\begin{equation}\label{eq.displacedthermal}
\phi^{\bf u} := \phi(\sqrt{2\mu -1} \alpha_{\bf u}) ,\qquad {\bf u} \in\mathbb{R}^{3},
\end{equation}
with $\alpha_{\bf u}:= -u_{y}+ iu_{x}$. Note that $\phi^{\bf u}$ does not depend on $u_z$.

We claim that the ``statistical information'' contained in the joint
state of $n$ qubits
\begin{equation}\label{eq.family.n}
\rho_n^{\bf u} :=  \rho_{{\bf u}/\sqrt{n}}^{\otimes n},
\end{equation}
is asymptotically identical to that contained in the couple
$(N^{\bf{u}} , \phi^{\bf u})$. More precisely:

\begin{theorem}\label{th.qlan}
Let $\rho_n^{\bf u}$ be the family of states \eqref{eq.family} on the
Hilbert space
$\left( \mathbb{C}^{2} \right)^{\otimes n}$, let $N^{\mathbf{u}}$ be the family
\eqref{gausje} of Gaussian distributions, and let $\phi^{\bf u}$ be the family
\eqref{eq.displacedthermal} of displaced thermal equilibrium states of a quantum
oscillator. Then for each $n$ there exist quantum channels (trace preserving CP maps)
\begin{equation*}
\begin{split}
T_{n} :
\mathcal{T}((\mathbb{C}^{2})^{\otimes n})
\to
           L^{1}(\mathbb{R})\otimes \mathcal{T}(\mathcal{F}(\mathbb{C})), \\
S_{n} :  L^{1}(\mathbb{R})\otimes\mathcal{T}(\mathcal{F}(\mathbb{C})) \to
\mathcal{T}((\mathbb{C}^{2})^{\otimes n})
\end{split}
\end{equation*}
with $\mathcal{T}(\mathcal{H})$ the trace-class operators on $\mathcal{H}$,
such
that, for any $0 \leq \eta<1/4$ and any $\epsilon> 0$,
\begin{eqnarray}\label{eq.channel.conv.}
&&
\sup_{\|{\bf u}\| \leq n^{\eta}} \,
\|  N^{\bf u}\otimes \phi^{\bf u} - T_{n} \left(  \rho^{\bf u}_{n}\right)
\|_{1} = \mathcal{O}(n^{-1/4+\eta + \epsilon} ), \\
&&
\sup_{\|{\bf u}\| \leq n^{\eta}} \,
\| \rho^{\bf u}_n - S_{n} \left(  N^{\bf u}\otimes \phi^{\bf u}\right)  \|_{1}
= \mathcal{O}(n^{-1/4+ \eta + \epsilon }). \label{eq.channel.conv.inverse}
\end{eqnarray}
Moreover, for each $\epsilon_2 > 0$ there exists a function $f(n)$ of order
$\mathcal{O}(n^{-1/4+\eta + \epsilon} )$ such that the above convergence rates are
bounded by $f(n)$, with $f$ independent of $\rho^{\mathbf{0}}$ as long as
$|\half - \mu| > \epsilon_2$.\vspace{-2mm}
\end{theorem}

\basremark{
Note that the equations \eqref{eq.channel.conv.} and 
\eqref{eq.channel.conv.inverse} imply that the expressions on the 
left side converge to zero as $n\to\infty$. 
Following the classical terminology of Le Cam \cite{LeCam}, we will call this 
type of result 
{\it strong convergence} of quantum statistical models (experiments). 
Another local asymptotic normality result has been derived in 
\cite{Guta&Jencova} 
based on a different concept of convergence, which is  an extension of the 
{\it weak convergence} of classical (commutative) statistical experiments. 
In the classical set-up it is known that strong convergence implies weak 
convergence for arbitrary statistical models, and the two are equivalent 
for statistical models consisting of a finite number of distributions. 
A similar relation is conjectured to hold in the quantum set-up, but for 
the moment this has been shown only under additional assumptions 
\cite{Guta&Jencova}. 

These two approaches to local asymptotic normality in quantum statistics are 
based on completely different methods and the results are complementary in the 
sense that the weak convergence of \cite{Guta&Jencova} holds for the larger 
class of finite dimensional states while the strong convergence has more 
direct consequences as it is shown 
in this paper for the case of qubits. Both results are part of a larger 
effort to develop a general theory of local asymptotic normality in quantum 
statistics. Several extensions are in order: from qubits to arbitrary 
finite dimensional systems (strong convergence), from finite dimensional 
to continuous variables systems, from identical system to correlated ones, 
and asymptotic normality in continuous time dynamical set-up. 

Finally, let us note that the development of a general theory of convergence 
of 
quantum statistical models will set a framework for dealing with other 
important 
statistical decision problems such as quantum cloning \cite{We2} and quantum 
amplification \cite{Caves}, which do not necessarily involve
measurements.\vspace{-1mm}
}

{{\bf \noindent Remark} \quad
The construction of the channels $T_{n}, S_{n}$ in the case of
fixed eigenvalues $(u_{z} =0)$ is given in theorem 1.1 of \cite{Guta&Kahn}. It
is also shown that a similar result holds uniformly over
$\| {\bf u}\|<C $ for any fixed finite constant $C$.
In \cite{Guta&Jencova}, it is shown that such maps also exist in the general
case, with unknown eigenvalues. A classical component then appears in the limit
statistical experiment. In the above result we extend the domain of validity of
these 
theorems from ``local'' parameters $\|{\bf u}\|<C $ to ``slowly growing'' local
neighbourhoods $\|{\bf u}\| \leq n^{\eta}$ with $\eta<1/4$.
Although this may be seen as merely a
technical improvement, it is in fact essential in order to insure
that the result of the first step of the estimation will, with high probability, fall
inside a neighbourhood $\|{\bf u}\| \leq n^{\eta}$ for which local
asymptotic normality still holds (see figure \arabic{tempfig}).\
}

\proof Following \cite{Guta&Kahn} we will first indicate how the channels
$T_{n}$ are constructed. The technical details of the proof can be found in
appendix \ref{sec.proof.qlan}.

The space $\left( \mathbb{C}^{2}\right)^{\otimes n}$ carries two 
unitary representations.
The representation $\pi_{n}$ of $\mathrm{SU}(2)$ is given by
$\pi_{n}(u) = u^{\otimes n}$ for any $u\in \mathrm{SU}(2)$,
and the representation $\tilde{\pi}_{n}$ of the symmetric group $S(n)$ is given
by the permutation of factors
$$
\tilde{\pi}_{n} (\tau) : v_{1} \otimes \dots \otimes  v_{n }
\mapsto v_{\tau^{-1}(1)} \otimes \dots \otimes v_{\tau^{-1}(n)}, 
\qquad\tau\in S(n).
$$
As $[\pi_{n} (u) , \tilde{\pi}_{n} (\tau)] =0$ for all $u\in \mathrm{SU}(2), 
\tau \in S(n)$,
we have the decomposition
\begin{equation}\label{eq.decomposition}
\left( \mathbb{C}^{2}\right)^{\otimes n} =
\bigoplus_{j=0, 1/2}^{n/2} \mathcal{H}_{j} \otimes \mathcal{H}^{j}_{n}.
\end{equation}
The direct sum runs over all positive (half)-integers $j$ up to $n/2$.
For each fixed $j$,
$\mathcal{H}_{j} \cong \mathbb{C}^{2j+1}$ is an irreducible representation 
$U_{j}$ of
$\mathrm{SU}(2)$ with total angular momentum \mbox{$J^{2} = j(j+1)$}, and
$\mathcal{H}^{j}_{n}\cong \mathbb{C}^{n_{j} }$ is the irreducible 
representation
of the symmetric group $S(n)$ with $n_{j}=\binom{n}{n/2 -j} - 
\binom{n}{n/2-j-1} $.
The density matrix $\rho_n^{\bf u}$ is invariant under
permutations and can be decomposed as a mixture of ``block'' density 
matrices
\begin{equation}
\label{blocks}
\rho_n^{\bf u} = \bigoplus_{j=0, 1/2}^{n/2}  p_{n,{\bf u}}(j) 
\,\rho^{\bf u}_{j,n} \otimes
\frac{\mathbf{1}}{n_{j}} \,.
\end{equation}
The probability distribution $p_{n,{\bf u}}(j)$ is given by 
\cite{Bagan&Gill}:
\begin{equation}
\label{pnj}
p_{n,{\bf u}}(j) := \frac{n_{j}}{2\mu_{\bf u} -1} 
\left(1-\mu_{\bf u}\right)^{\frac{n}{2} - j}
\mu_{\bf u}^{\frac{n}{2} + j +1}
\left(1-p_{\bf u}^{2j +1}\right),
\end{equation}
with $\mu_{\bf u}:= \mu + u_{z} /\sqrt{n}$, $p_{\bf u}:= 
\frac{1-\mu_{\bf u}}{\mu_{\bf u}} $.
We can rewrite $p_{n,{\bf u}}(j)$ as
\begin{equation}\label{eq.pnj}
p_{n,{\bf u}}(j) := B_{n,\mu_{\bf u}} (n/2+j)\times K(j,n,\mu,{\bf u}),
\end{equation}
where
$$
B_{n,\nu} (k) :=  \binom{n}{k}  \nu^{k} \left(1-\nu \right)^{n-k} , 
\qquad k= 0, \dots, n
$$
is a binomial distribution, and the factor $ K(j,n,\mu,{\bf u})$ is given by
$$
K(j,n,\mu,{\bf u}):=
\left(1 - p_{\bf u}^{2j +1}\right)
\frac{n+ (2(j-j_{n} -\sqrt{n}u_{z}) +1)/(2\mu_{\bf u} -1) }
{n + (j-j_{n} -\sqrt{n}u_{z}+1)/\mu_{\bf u}}\,,
$$
with $j_{n} := n(\mu-1/2)$.
Now $K(j,n,\mu,{\bf u}) = 1 + \mathcal{O}(n^{-1/2+\epsilon})$
on
the relevant values of $j$, i.e. the ones in an interval of order 
$n^{1/2+\epsilon}$ around
$j_{n}$, as long as $\mu_{\bf u}$ is bounded away from $1/2$, which 
is
automatically so for big $n$. As $B_{n,\mu_{\bf u}} (k) $ is the 
distribution of a sum of
i.i.d. Bernoulli random variables, we can now use standard local 
asymptotic normality
results \cite{vanderVaart} to conclude that if $j$ is distributed 
according to $p_{n,{\bf u}}$, then
the centered and rescaled variable
$$
g_{n} := \frac{j}{\sqrt n} - \sqrt{n}(\mu - 1/2),
$$
converges in distribution to a normal $N^{\bf u}$, after an additional 
randomization has been performed.
The latter is necessary in order to ``smooth''  the discrete distribution 
into a distribution which is 
continuous with respect to the Lebesgue measure, 
and will convergence to the Gaussian distribution in total variation norm.

The measurement ``which block'', corresponding to the decomposition
\eqref{blocks}, provides us with a result $j$ and a posterior state 
$\rho^{\bf u}_{j,n}$.
The function $g_{n}= g_{n}(j)$ (with an additional randomization) is 
the classical part of the channel $T_{n}$.
The randomization consists of ''smoothening'' with a Gaussian kernel of mean
$g_{n}(j)$ and variance $1/(2\sqrt{n})$, i.e. with
$\tau_{n,j} :=  (n^{1/4}/\sqrt{\pi})\exp\left(-\sqrt{n} (x-g_n(j))^2\right)$.

Note that this measurement is not disturbing the state
$\rho_n^{\bf u}$ in the sense that the average state after the measurement is 
the same as before. 

The quantum part of $T_{n}$ is the same as in \cite{Guta&Kahn} and consists of  embedding each block state $\rho^{\bf u}_{j,n}$ into the state space
of the oscillator by means of an isometry $V_{j} : \mathcal{H}_{j} 
\to\mathcal{F}(\mathbb{C})$,
$$
V_{j}: |j,m \rangle \mapsto |j-m \rangle,
$$
where $\{ |j,m \rangle : m=-j,\dots, j\}$ is the eigenbasis of the total 
spin component
$L_{z}:= \sum_{i}\sigma^{(i)}_{z}$, cf. equation (5.1) of \cite{Guta&Kahn}.
Then the action of the
channel $T_{n}$ is
\[
T_n: \bigoplus_{j} p_{n,{\bf u}}(j) \rho_{j,n}^{\bf u} \otimes \frac{{\bf 1}}
{
n_j} \mapsto \sum_{j} p_{n,{\bf u}}(j) \, \tau_{n,j} \otimes
V_{j}\rho_{j,n}^{\bf u} V^\dagger_{j}\,.
\]

The inverse channel $S_{n}$ performs the inverse operation with respect to
$T_{n}$. First the oscillator state is ``cut-off'' to the dimension of an 
irreducible representation
and then a block obtained in this way is placed into the decomposition
\eqref{eq.decomposition} (with an additional normalization from the remaining
infinite dimensional block which is negligible for the states in which we 
are interested).

The rest of the proof is given in appendix \ref{sec.proof.qlan}.
\qed

%

\section{Time Evolution of the Interacting System}    \label{sec.timevolf}

In the previous section, we have investigated the asymptotic
equivalence between the states $\rho_{n}^{\bf u}$ and
$N^{\bf u} \otimes \phi^{\bf u}$  by means of the channel $T_{n}$.
We now seek to implement this in a physical situation.
The $N^{\bf u}$-part will follow in section
\ref{sec.energy}, the $\phi^{\bf u}$-part will be treated in this section.

We couple the $n$ qubits to a bosonic field;
this is the physical implementation of LAN.
Subsequently, we perform a
measurement in the field  which will provide the information about the state
of the qubits; this is the utilization of LAN in order to
solve the asymptotic state estimation problem.

In this section we will limit ourselves to analyzing the joint
evolution of the qubits and field. The measurement on the field
is described in section
\ref{finmes}.

\subsection{Quantum Stochastic Differential Equations}

In the weak coupling limit, \cite{AFL,Gardiner}, 
the joint evolution of the qubits and field can
be described mathematically by quantum stochastic differential 
equations (QSDE). See section \ref{introqsdesec} for a short introduction,
and
\cite{HuP84, Par92}
for a thorough account. 
The basic notions here are the Fock space, the creation and annihilation 
operators and the quantum stochastic differential equation of the unitary 
evolution. The Hilbert space of the field is the Fock space 
$\mathcal{F}(L^{2}(\mathbb{R}))$ as defined in \eqref{eq.Fock}. 
An important linearly complete set in
$\mathcal{F}(L^{2}(\mathbb{R}))$ is that of the exponential vectors
\begin{equation}\label{eq.exponential}
e(f) := \bigoplus_{n=0}^{\infty} \frac{1}{\sqrt{n!}} f^{\otimes n} :=\bigoplus_{n=0}^{\infty} \frac{1}{\sqrt{n!}} |f \rangle_{n}  , \qquad f\in L^{2}(\mathbb{R}),
\end{equation}
with inner product $\inp{e(f)}{e(g)} = \exp(\inp{f}{g})$.
The normalized exponential states
$\ket{f} := e^{-\inp{f}{f} /2} e(f)$ are called
coherent states.
The vacuum vector is $|\Omega\rangle:=e(0)$ and we will denote the
corresponding density matrix $\ket{\Omega} \bra{\Omega}$ by $\Phi$.
The quantum noises are described by the creation and annihilation martingale 
operators
$A_{t}^{\dagger}: = a^{\dagger}(\chi_{[0,t]})$ and $A_{t}: = a(\chi_{[0,t]})$
respectively, where $\chi_{[0,t]}$ is the indicator function for $[0,t]$ and
$$
a(f) :e(g) \mapsto \inp{f}{g} e(g).
$$
The increments $dA_{t} := a(\chi_{[0,t+dt]})-a(\chi_{[0,t]})$ and 
$dA^{\dagger}_{t}$ 
play the role of non-commuting integrators in quantum stochastic 
differential equations, 
comparable to the role of Brownian motion in classical stochastic calculus.

We now consider the joint unitary evolution for qubits and field defined by 
the quantum stochastic differential equation \cite{HuP84,Bo}:
$$
dU_{n}(t) = (a_{n}dA^{\dagger}_{t} -a^{\dagger}_{n} dA_{t} -\frac{1}{2} 
a^{\dagger}_{n}a_{n} dt )U_{n}(t),
$$
where $U_{n}(t)$ is a unitary operator on 
$(\mathbb{C}^{2})^{\otimes n}\otimes \mathcal{F}(L^{2}(\mathbb{R}))$, and
$$
a_{n}: = \frac{1}{\sqrt{2j_{n}}}\sum_{k=1}^{n} \sigma_{+}^{(k)}\,,\qquad
 \sigma_{+}^{(k)} := \mathbf{1}\otimes \dots \otimes 
 (\sigma_{x} +i \sigma_{y})/2 \otimes \dots \otimes \mathbf{1}\,,
$$
with $j_{n} := (\mu-1/2)n$.
As we will see later, the ``coupling factor'' $1/\sqrt{j_n}$ of the order
$n^{-1/2}$, is necessary in order to obtain convergence to the unitary 
evolution of the quantum harmonic oscillator and the field.

We remind the reader that the $n$-qubit space can be decomposed into 
irreducible representations as in \eqref{eq.decomposition}, and the 
interaction between the qubits and field respects this decomposition\vspace{-1mm}
$$
U_{n}(t) = \bigoplus_{j=0,1/2}^{n/2} U_{j,n}(t) \otimes \mathbf{1},\vspace{-1mm}
$$
where $\mathbf{1} $ is the identity operator on the multiplicity space 
$\mathcal{H}_{n}^{j}$, and
$$
U_{j,n}(t) : \mathcal{H}_{j}\otimes \mathcal{F}(L^{2}(\mathbb{R})) \to
\mathcal{H}_{j}\otimes \mathcal{F}(L^{2}(\mathbb{R})) ,
$$
is the restricted cocycle
\begin{equation}\label{eq.cocycle.j}
d U_{j,n}(t) = (a_{j}dA^{\dagger}_{t} -a^{\dagger}_{j} dA_{t} -\frac{1}{2} 
a^{\dagger}_{j}a_{j} dt )U_{j,n}(t),
\end{equation}
with
$a_{j} $ acting on the basis $|j,m\rangle$ of $\mathcal{H}_{j}$ as
\begin{eqnarray*}
&&
a_{j} |j,m\rangle = \sqrt{j-m} \sqrt{ (j+m+1)/2j_{n}} \, |j,m+1\rangle ,\\
&&
a^{\dagger}_{j} |j,m\rangle = \sqrt{j-m+1} \sqrt{ j+m/2j_{n}}  \,  |j,m-1 
\rangle .\vspace{2mm}
\end{eqnarray*}

{\bf \noindent Remark}\quad
We point out that the {\it lowering} operator for $L_{z}$ acts 
as {\it creator} for our cut-off
oscillator since the highest vector $|j,j\rangle$ corresponds by $V_{j}$ to 
the vacuum of the oscillator. This choice does not have any physical meaning 
but is only related with our convention $\mu>1/2$. Had we chosen $\mu<1/2$, 
then the raising operator on the qubits would correspond to creation operator 
on the oscillator.

By \eqref{blocks} the initial state $\rho^{\otimes n}$ decomposes in the same 
way as the unitary cocycle, and thus the whole evolution decouples into 
separate ``blocks'' for each value of $j$. We do not have explicit solutions 
to these equations but based on the conclusions drawn from LAN we expect that 
as $n\to\infty$, the solutions will be well approximated by similar ones for 
a coupling between an oscillator and the field, at least  for the states in 
which we are interested. As a warm up exercise we will start with this 
simpler limit case where the states can be calculated explicitly.

\subsection{Solving the QSDE for the Oscillator}\label{sec.oscillator.field}

Let $a^{\dagger}$ and $a$ be the creation and annihilation operators of a 
quantum oscillator acting on $\mathcal{F}(\mathbb{C})$.  We couple the 
oscillator with the bosonic field
and the joint unitary evolution is described by the family of unitary 
operators $U(t)$ satisfying the quantum stochastic differential equation
$$
dU(t) =
(adA^{\dagger}_{t} -a^{\dagger} dA_{t} -\frac{1}{2} a^{\dagger}a dt )U(t).
$$
We choose the initial (unnormalized) state $\psi (0):= e({\bf z} )\otimes
|\Omega\rangle$, where ${\bf z}$ is any complex number, and we shall find the 
explicit form of the vector state of the system and field at time $t$:
$\psi(t):=U(t) \psi(0)$.

We make the following ansatz: $\psi(t) = e(\alpha_{t}) \otimes e(f_{t})$,
where $f_t$ is given by $f_t(s) := f(s) \chi_{[0,t]}(s)$ for some $f \in L^2(\RR)$.
For each $\beta \in \CCCC$, $g \in L^2(\RR)$, define
$I(t) := \inp{e(\beta) \otimes e(g)}{\psi(t)}$.
We then have $I(t) = \exp(\bar{\beta} \alpha(t) + \inp{g}{f_t})$,
so that it satisfies
\begin{equation}\label{sarcofaag}
d I(t) = \left( \bar{\beta} \dt \alpha(t) +
\bar{g}(t) f(t) \right)I(t) dt\,.
\end{equation}

We now calculate $\frac{d}{dt} I(t)$ with the help of the QSDE.
Since the annihilator acts as $A_{t} e(f) = \inp{\chi_{[0,t]}}{f} e(f)$,
we have, for continuous $g$,
$dA_t e(g) = g(t)e(g) dt$.
However, since $A_{s} e(f_t)$ is constant for $s \geq t$,
we have $dA_t e(f_t) = 0$.
Thus
$$
d I(t) =
\inp{e(\beta) \otimes e(g)}{(adA^{\dagger}_{t} -a^{\dagger} dA_{t} -
\half a^{\dagger}a dt ) \psi(t)} \,,
$$
so that
\begin{equation}\label{knaagdier}
dI(t) = (\bar{g}(t) \alpha(t)  -
\half \bar{\beta} \alpha(t)) I(t) dt\,.
\end{equation}

Equating \eqref{sarcofaag} with \eqref{knaagdier}
for all $t$, all $\beta$, and all continuous $g$,
we find $f(s) = \alpha(s)$,
$\frac{d}{dt} \alpha(t) = - \frac{1}{2} \alpha(t)$.
Thus $\alpha(t) = \alpha(0) e^{- \frac{1}{2} t}$,
$f_t(s) = \alpha(0)\chi_{[0,t]}(s) e^{- \frac{1}{2} s}$ with 
$\alpha(0)={\bf z}$.
In conclusion, $\psi(t) = e({\bf z} e^{-\frac{1}{2}t}) \otimes
e({\bf z}  e^{-\frac{1}{2} s} \chi_{[0,t]}(s))$. For later use, we 
denote by
$\psi_{\bf z}(t):= U(t) |{\bf z}\rangle \otimes |\Omega\rangle = 
e^{-|{\bf z}|^{2}/2}U(t) e({\bf z})\otimes | \Omega \rangle$
the {\it normalized} solution of the QSDE.

\subsection{QSDE for Large Spin}\label{sec.unitary.evol}

We consider now the unitary evolution for qubits and field:
$$
dU_{n}(t) = (a_{n}dA^{\dagger}_{t} -a^{\dagger}_{n} dA_{t} -\frac{1}{2} 
a^{\dagger}_{n}a_{n} dt )U_{n}(t).
$$
It is no longer possible to obtain an explicit expression for the joint
vector state $\psi_{n}(t)$ at time $t$.
However we will show that for the states in which we are interested, a 
satisfactory explicit {\it approximate} solution exists.

The trick works for an arbitrary family of unitary solutions of a quantum 
stochastic differential equation $dU(t)= G_{dt}U(t)$, and the general idea 
is the following. If $\psi(t)$ is the true state $\psi(t) = U(t)\psi$ and 
$\xi(t)$ is a vector describing an approximate evolution ($\psi(0)=\xi(0)$) 
then with $U^{t}_{t+dt} :=U(t+dt)U(t)^{-1}$,  we can write  
$\psi(t+dt) - \xi(t+dt)$ as
$$
\psi(t+dt) - U^{t}_{t+dt} \xi(t) + U^{t}_{t+dt} 
\xi(t) -\xi(t) +\xi(t) -\xi(t+dt)\,.
$$
This equals
$$
U^{t}_{t+dt} \left[ \psi(t) -\xi(t) \right] + [U(t+dt) - U(t)]U(t)^{-1} 
\xi(t) 
[\xi(t)  -\xi(t+dt)]\,,
$$
so that 
$$
\psi(t+dt) - \xi(t+dt)
=
U^{t}_{t+dt} \left[ \psi(t) -\xi(t) \right] + G_{dt} \xi(t) - d\xi(t)\,.
$$
By taking norms we finally obtain
\begin{equation}\label{eq.approx.qsde.solution}
d\| \psi(t) - \xi(t) \| \leq \| G_{dt} \xi(t) - d\xi(t)\|.
\end{equation}
The idea is now to devise a family $\xi(t)$ such that the right side is 
as small as possible.

We apply this technique block-wise, that is to each unitary $U_{j,n}(t)$ 
acting on
$\mathcal{H}_{j}\otimes \mathcal{F}(L^{2}(\mathbb{R}))$ (see equation 
\eqref{eq.cocycle.j})
for a ``typical'' $j\in \mathcal{J}_{n}$ (see equation  \eqref{eq.typical.j}).
By means of the isometry $V_{j}$ we can embed the space $\mathcal{H}_{j}$ 
into the first $2j+1$ levels of the oscillator and for simplicity we will 
keep the same notions as before for the operators acting on 
$\mathcal{F}(\mathbb{C})$. As initial states for the qubits we choose 
the block states $\rho^{\bf u}_{j,n}$.

\begin{theorem}\label{th.unitary.evolution}
Let $\rho^{\bf u}_{j,n} (t) = U_{j,n} (t)\, \left[ \rho^{\bf u}_{j,n} 
\otimes  \Phi\right] \, U_{j,n}^\dagger (t)$ be the $j$-th block of the 
state of qubits and field at time $t$.
Let $\phi^{\bf u}(t):= U(t) \,\left[\phi^{\bf u} \otimes \Phi \right]\,
U(t)^\dagger $ be the joint
state of the oscillator and field at time $t$.  For any $\eta <1/6$, for 
any
$\epsilon>0$,
\begin{equation}\label{eq.unitary.error}
 \sup_{j\in\mathcal{J}_{n}} \,
 \sup_{\|{\bf u} \| \leq n^{\eta}} \, \sup_{t} \| \rho^{\bf u}_{j,n} (t) -
\phi^{\bf u}(t) \|_{1} =\mathcal{O}(n^{-1/4+\eta + \epsilon}, n^{-1/2+ 3\eta + \epsilon}).
\end{equation}
\end{theorem}
\proof From the proof of the 
local asymptotic normality 
theorem \ref{th.qlan} we know that the initial states of the two 
unitary evolutions are asymptotically close to each other
\begin{equation}\label{eq.qlan.evolution}
\sup_{j\in\mathcal{J}_{n}} \,
 \sup_{\|{\bf u} \| \leq n^{\eta}} \| \rho^{\bf u}_{j,n}  - \phi^{\bf u}
\|_{1} = \mathcal{O}(n^{-1/4+\eta + \epsilon}).
 \end{equation}

The proof consists of two estimation steps. In the first one, we will devise 
another initial state $\tilde{\rho}^{\bf u}_{j,n}$ which is an approximation 
of $\phi^{\bf u}$ and thus also of $\rho^{\bf u}_{j,n}$:
\begin{equation}\label{eq.error.initial.condition}
\sup_{j\in\mathcal{J}_{n}} \,
 \sup_{\|{\bf u} \| \leq n^{\eta}} \| \tilde{\rho}^{\bf u}_{j,n}  -
\phi^{\bf u}  \|_{1} = \mathcal{O}( e^{-n^{\epsilon}}).
\end{equation}
In the second estimate we show that the evolved states 
$\tilde{\rho}^{\bf u}_{j,n}(t)$ and
$\phi^{\bf u}(t)$ are asymptotically close to each other
\begin{equation}\label{eq.error.evolved.states}
\sup_{j\in\mathcal{J}_{n}} \,
 \sup_{\|{\bf u} \| \leq n^{\eta}} \sup_{t} \| \tilde{\rho}^{\bf u}_{j,n}(t)  
 -
 \phi^{\bf u}(t)  \|_{1} =
 \mathcal{O}(n^{-1/4+\eta +\epsilon}, n^{-1/2 + 3\eta + \epsilon}).
\end{equation}
This estimate is important, because the two trajectories are driven by 
different Hamiltonians, and in principle there is no reason why they should 
stay close to each other.   
From \eqref{eq.qlan.evolution}, \eqref{eq.error.initial.condition} and  
\eqref{eq.error.evolved.states}, and using the triangle inequality, we get
\begin{equation*}
\sup_{j\in\mathcal{J}_{n}} \,
 \sup_{\|{\bf u} \| \leq n^{\eta}}
 \sup_{t}\| \rho^{\bf u}_{j,n} (t) - \phi^{\bf u}(t)  \|_{1} =
 \mathcal{O}(n^{-1/4+\eta + \epsilon},
n^{-1/2 + 3\eta + \epsilon}).
\end{equation*}

The following diagram illustrates the above estimates. 
The upper line concerns the
time evolution of the block state $\rho^{\bf u}_{j,n}$ and the field. 
The lower
line describes the time evolution of the oscillator and the field. 
The estimates show that the diagram is ``asymptotically commutative'' 
for large $n$.
\begin{center}
$
\CD
\CS(\mathcal{H}_{j}) @> {\rm Id}_{j} \otimes \Phi>> \CS(\mathcal{H}_{j} \otimes
\CF) @>U_{j,n}(t)>>
\CS(\mathcal{H}_{j} \otimes \CF) 
\\
@V {V_{j}\cdot V_{j}^{\dagger} }VV  @VVV  @VVV 
\\
\CS(\CF(\CCCC)) @>{\rm Id} \otimes \Phi >> \CS(\CF(\CCCC) \otimes \CF) @>U(t)>>
\CS(\CF(\CCCC) \otimes \CF) 
\endCD
$
\end{center}
For the rest of the proof, we refer to appendix \ref{kool}.
\qed

\noindent We have shown how the mathematical statement of LAN 
(the joint state of qubits converges to a Gaussian state of a 
quantum oscillator plus a classical Gaussian random variable) 
can in fact be physically implemented by coupling the spins to 
the environment and letting them ``leak'' into the field. In 
the next section, we will use this for the specific purpose 
of estimating $\bf u$ by performing a measurement in the field.


\section{The Second Stage Measurement}
\label{finmes}

We now describe the second stage of our measurement procedure.
Recall that in the first stage a relatively small part $\tilde{n}=n^{1-\kappa}, 1>\kappa>0,$
of the qubits is measured and a rough estimator $\tilde{\rho}_{n}$ is obtained. The purpose of 
this estimator is to localize the state within a small neighbourhood such
that the 
machinery of local asymptotic normality of 
theorem \ref{th.qlan} can be applied.

In theorem \ref{th.unitary.evolution} the local asymptotic normality was 
extended to the level of time evolution of the qubits interacting with a 
bosonic field. We have proven that at time $t$ the joint state of the 
qubits and field is
\begin{eqnarray*}
\rho_n^{\bf u}(t)
& := &\bigoplus_{j=0,1/2}^{n/2} p_{n,{\bf u}}(j) \frac1{2\pi s^2}
\int_{ \mathbb{C} } d{\bf z}\,
e^{-|{\bf z} -
\sqrt{2\mu -1} \alpha_{\bf u}|^2/2s^2} \exp(-|{\bf z}|^2)
\times
\\& &
|e({\bf z} e^{-t/2 })_j\rangle\langle e({\bf z}
e^{-t/2 })_j| \otimes
|e({\bf z} e^{-u/2 }\chi_{[0,t]}(u))\rangle \langle e({\bf z} e^{-u/2 }
\chi_{[0,t]}(u))|
\\& &
+ \mathcal{O}(n^{\eta -1/4 +\epsilon }, n^{3\eta -1/2 + \epsilon}),
\end{eqnarray*}
for $\|{\bf u}\|\leq n^{\eta}$.
The index $j$ serves to remind the reader that the first exponential states
live in different copies $\mathcal{F}(\mathbb{C})_j$ of the oscillator space,
corresponding to $\mathcal{H}_j$ via the isometry $V_j$. We will continue to
identify $\mathcal{H}_j$ with its image in $\mathcal{F}(\mathbb{C})_j$.

We can now approximate the above state by its limit for large $t$, since
\begin{align}
\label{approx.time}
\exp(-|{\bf z}|^2) 
\inp{ e({\bf z} e^{-t/2 })_j  }{ j,j \,} 
\inp
{
e({\bf z}e^{-u/2 }\chi_{[0,t]}(u))
}
{
e({\bf z} e^{-u/2})
}  =
\exp( - |{\bf z}|^2 e^{-t}).
\end{align}

As we are always working with $\|{\bf u}\|\leq n^{\eta}$, the only relevant
${\bf z}$
are bounded by $n^{\eta + \delta}$ for small $\delta$.
(The remainder of the Gaussian
integral has an exponentially decreasing norm, as discussed before).
Thus, for large enough time (i.e. for $t \geq \ln(n)$),
 we can write $\rho_n^{\bf u}(t) =
 \rho_n^{\bf u}(\infty)  + \mathcal{O}(n^{\eta -1/4 + \epsilon} ,
 n^{3\eta -1/2 + \epsilon})$ with
\begin{equation}
\begin{gathered}
 \rho_n^{\bf u}(\infty) :=
 \bigoplus_{j=0,1/2}^{n/2} p_{n,{\bf u}}(j) |j, j\rangle\langle j,j|
\otimes \\
\left[ \frac1{2\pi s^2}\int_{ \mathbb{C} } d{\bf z} \, e^{-|{\bf z} -
\sqrt{2\mu -1} \alpha_{\bf u}|^2/2s^2}
|e({\bf z} e^{-u/2 })\rangle \langle e({\bf z} e^{-u/2})|
\exp(-|{\bf z}|^2) \right]. \label{limstate}
\end{gathered}
\end{equation}

Thus, the field is approximately in the state 
$\phi^{\bf u}$ depending on $(u_{x}, u_{y})$, which is carried by the mode 
$(u\mapsto e^{-u/2} \chi_{[0,\infty)} (u))\in L^{2}(\mathbb{R})$ denoted 
for simplicity by $e^{-u/2}$. The atoms end up in a mixture of $|j,j\rangle$ 
states with coefficients 
$p_{n,{\bf u}}(j)$, which depend only on $u_{z}$, and are well approximated 
by the Gaussian random variable $N^{\bf u}$ as shown in 
theorem \ref{th.qlan}.  
Moreover since there is no correlation between atoms and field, the 
statistical problem decouples into one concerning the estimation of 
the displacement in a family of Gaussian states $\phi^{\bf u}$, and 
one for estimating the center of $N^{\bf u}$.

For the former problem, the optimal estimation procedure is known to be 
the heterodyne measurement \cite{Ho,Yuen&Lax}; for the latter, we 
perform a  ``which block'' measurement. These measurements are 
described in the next two subsections.

\subsection{The Heterodyne Measurement}
\label{subsec.hetero}

A heterodyne measurement is a ``joint measurement'' of the two quadratures 
${\bf Q}:= (a+a^{\dagger})/\sqrt{2}$ and ${\bf P}:= 
-i(a-a^{\dagger})/\sqrt{2}$ of a quantum harmonic oscillator which in 
our case represents a mode of light. Since the two operators do not 
commute, the price to pay is the addition of some ``noise'' which will 
allow for an approximate measurement of both operators. The light beam 
passes through a beamsplitter having a vacuum mode as the second input, 
and then one performs a homodyne (quadrature) measurement on each of the 
two emerging beams. If ${\bf Q}_{v}$ and ${\bf P}_{v}$ are the vacuum 
quadratures then we measure the following output quadratures 
${\bf Q}_{1} := ({\bf Q} + {\bf Q}_{v})/\sqrt{2}$ and 
${\bf P}_{2} := ({\bf P} - {\bf P}_{v})/\sqrt{2}$, with 
$[{\bf Q}_{1}, {\bf P}_{2}]
=0$. Since the two input beams are independent,  the distribution of 
$\sqrt{2}{\bf Q}_{1}$ is the convolution between the distribution of 
${\bf Q}$ and the distribution of ${\bf Q}_{v}$, and similarly for 
$\sqrt{2}{\bf P}_{2}$. 

In our case we are interested in the mode 
$e^{-u/2}$ which is in the state $\phi^{\bf u}$, up to a factor of order
$\mathcal{O}(n^{\eta -1/4 + \epsilon} , n^{3\eta -1/2 + \epsilon})$. From 
\eqref{eq.displacedthermal} we obtain that the distribution of ${\bf Q}$ 
is $N(\sqrt{2(2\mu-1)} u_{x}, 1/(2(2\mu-1)))$, 
that of ${\bf P}$ is $N(\sqrt{2(2\mu-1)} u_{y}, 1/(2(2\mu-1)))$, and the 
joint distribution 
of the rescaled output 
$$
\left( ({\bf Q}+ {\bf Q}_{v})/ \sqrt{2(2\mu-1)} \, , \, 
({\bf P}- {\bf P}_{v})/ \sqrt{2(2\mu-1)} \right) ,\vspace{-1mm}
$$
is\vspace{-1mm} 
\begin{equation}\label{gaus1}
N(u_{x} ,  \mu/(2(2\mu-1)^2)) \times N(u_{y} ,  \mu/(2(2\mu-1)^2)).
\end{equation}
We will denote by $(\tilde{u}_x, \tilde{u}_y)$ the result of the heterodyne 
measurement rescaled by the factor $\sqrt{2\mu-1}$ such that with good 
approximation 
$(\tilde{u}_x, \tilde{u}_y)$ has the above distribution and is an unbiased 
estimators of the parameters $(u_x, u_y)$.

Since we know in advance that the parameters $(u_{x}, u_{y})$ must be within 
the radius 
of validity of LAN we modify the estimators $(\tilde{u}_{x}, \tilde{u}_{y})$ 
to account for this information and obtain the final estimator $(\hat{u}_{x}, 
\hat{u}_{y})$:
\begin{align}
\label{mod1}
\hat{u}_i = \left\{\begin{array}{cc}
 \tilde u_i &\quad \textrm{if $|\tilde u_i|\leq 3 n^{\eta}  $} \\
  0 & \quad        \textrm{if $|\tilde u_i|>  3 n^{\eta}  $}
\end{array}
 \right.
\end{align}

Notice that if the true state $\rho$ is in the radius 
of validity of LAN around $\tilde \rho$, then $\|{\bf u}\|\leq n^{\eta}$, so 
that $|\hat{u}_i -
u_i |\leq |\tilde{u}_i - u_i|$. We shall use this when proving optimality of 
the
estimator. 

\subsection{Energy Measurement}        
\label{sec.energy}

Having seen the $\phi^{\bf u}$-part,
we now move to the $N^{\bf u}$-part of the equivalence
between $\rho_{n}^{\bf u}$ and $N^{\bf u} \otimes \phi^{\bf u}$.
This too is a coupling to a bosonic field, albeit a different coupling.
We also describe the measurement in the field which will provide the
information on the qubit states.

The final state of the previous measurement, restricted to the atoms alone
(without the field), is obtained by a partial trace of
equation \eqref{limstate} (for large time) over the field\vspace{-1mm}
$$
\tau_{n}^{\bf u} = \sum_{j=0,1/2}^{n/2} p_{n,{\bf u}}(j) |j,j\rangle
\langle j,j|
+ \mathcal{O}(n^{\eta -1/4 + \epsilon} , n^{3\eta -1/2 + \epsilon}) \,.\vspace{-1mm}
$$
We will take this as the initial state of the second measurement, which 
will
\mbox{determine j.}

A direct coupling to the $J^2$ does not appear to be physically available,
but a coupling to the energy$J_z$ is realizable. This suffices, because the 
above state satisfies
$j = m$ (up to order $\mathcal{O}(n^{\eta -1/4 + \epsilon}, n^{3\eta -1/2 + \epsilon})$).
We couple the atoms to a new field (in the vacuum state $\ket{\Omega}$)
by means of the interaction
$$
dU_t = \{ J_z (dA_{t}^{\dagger} - dA_{t}) - \half J_{z}^2 dt \} U_t \,,
$$
with $J_z := \frac{1}{\sqrt{n}}\sum_{k=1}^{n}\sigma_z$.
Since this QSDE is `essentially commutative', i.e. driven by a single
classical noise
$B_t = (\smash{A^\dagger_{t}} - A_{t})/i$, the solution is easily seen to be
$$
U_t = \exp(J_z \otimes (A^\dagger_{t} - A_t)) \, .
$$
Indeed, we have $d f(B_t) = f'(B_t)dB_t + \half f''(B_t)dt$ by the classical 
It\^o rule, so that
$$
d \exp (i J_z \otimes B_t) = \{ iJ_z dB_t - \half J_z^2 dt \}
\exp (iJ_z \otimes B_t)\, .
$$
For an initial state $\ket{j,m}\otimes \ket{\Omega}$, this evolution
gives rise to the final state
\begin{eqnarray*}
U_t \ket{j,m}\otimes\Omega &=& \ket{j,m}\otimes
\exp((m/\sqrt{n})(A^{\dagger}_t - A_{t}))\Omega \\
&=& \ket{j,m}\otimes \ket{(m/\sqrt{n}) \chi_{[0,t]}},
\end{eqnarray*}
where $\ket{f} \in \mathcal{F}(L^2(\mathbb{R}))$ denotes the normalized
vector $\exp(-\inp{f}{f}/2)e(f)$.
Applying this to the states $\ketbra{j,j}$ in $\tau_{n}^{\bf u}$
yields
\begin{eqnarray*}
U_{t} \,\tau_{n}^{\bf u} \otimes \Phi\, U^\dagger_{t} &=& 
\sum_{j=0,1/2}^{n/2} p_{n,{\bf u}}(j)
\ketbra{j,j} \otimes \ket{j/\sqrt{n} \chi_{[0,t]}} 
\bra{j/\sqrt{n} \chi_{[0,t]}}\\  
& &+ \mathcal{O}(n^{\eta -1/4 + \epsilon}, n^{3\eta -1/2 + \epsilon})\,.
\end{eqnarray*}
The final state of the field results from a partial trace over the atoms;
it is given by
\begin{equation}\label{smeerpoets}
\sum_{j=0,1/2}^{n/2} p_{n,{\bf u}}(j) \,
\ketbra{ (j/\sqrt{n}) \chi_{[0,t]}}
+ \mathcal{O}(n^{\eta -1/4 + \epsilon}, n^{3\eta -1/2 + \epsilon})\,.
\end{equation}

We now perform a homodyne measurement on the field, which amounts to a direct
measurement of $(A_t + A^\dagger_t)/2t$.
In the state $\ket{(j/\sqrt{n} \chi_{[0,t]}}$, this yields
the value of $j$ with certainty for large time (i.e. $t\gg \sqrt{n}$).
Indeed, for this state, $\mathbb{E}((A_t + A^\dagger_t)/2t) = j/\sqrt{n}$,
whereas $\var(A_t + A^\dagger_t)/2t) = 1/(4t)$.
Thus the probability distribution $p_{n,{\bf u}}$ is reproduced up to order
$\mathcal{O}(n^{\eta -1/4 + \epsilon} , n^{3\eta -1/2 + \epsilon} )$ in $L^1$-distance.

The following is a reminder from the proof of theorem \ref{th.qlan}. 
If we start with $j$ distributed according to $p_{n}(j)$ and we 
smoothen $\frac{j}{\sqrt{n}} - \sqrt{n} (\mu-1/2)$ with a Gaussian kernel, 
then we obtain a random variable $g_n$ which is continuously distributed 
on $\mathbb{R}$ and converges in distribution to $N(u_z, \mu(1 - \mu))$, 
the error term being of order $\mathcal{O}(n^{\eta - 1/2}) + \mathcal{O}(n^{\epsilon - 1/2})$. 
For $j$ distributed according to the actual distribution, as measured by the
homodyne detection experiment, we can therefore state that $g_n$
is distributed according to
\begin{equation}
\label{gaus2}
N(u_z,\mu(1-\mu)) + \mathcal{O}(n^{\eta -1/4 + \epsilon}, n^{3\eta -1/2 + \epsilon})
+ \mathcal{O}(n^{\eta - 1/2}) + \mathcal{O}(n^{\epsilon - 1/2}).
\end{equation}

As in the case of $(\hat{u}_{x}, \hat{u}_{y})$, we take into account the 
range of validity of LAN by defining the final estimator 
\begin{align}
\label{mod2}
\hat{u}_z = \left\{\begin{array}{cl}
 g_{n} &\quad \textrm{if $| g_{n}|\leq 3 n^{\eta}  $} \\
  0 & \quad        \textrm{if $| g_{n}|>  3 n^{\eta}  $\,.}
\end{array}
 \right.
\end{align}
Similarly, we note that if the true state $\rho$ is in the radius of 
validity of LAN around $\tilde \rho$, then $\|{\bf u}\|\leq n^{\eta}$, 
so that 
$|\hat{u}_z -u_z |\leq |\tilde{u}_z - u_z|$.

\section{Asymptotic Optimality of the Estimator}  \label{sec.endresult}

In order to estimate the qubit state,
we have proposed a strategy consisting of the following steps.
First, we use
$\tilde{n}:=n^{1-\kappa}$ copies of the state $\rho$ to get a rough estimate
$\tilde \rho_{n}$. Then we couple the remaining qubits with a field, and 
perform a
heterodyne measurement. Finally, we couple to a different field,
followed by homodyne measurement. From the measurement outcomes,
we construct an estimator $\hat{\rho}_{n} := \rho_{\hat{\bf u}_{n}/\sqrt{n}}$.

This strategy is asymptotically optimal in a global sense: for {\it any} 
true
state $\rho$ even if we knew beforehand that the true state $\rho$ is in a 
small
ball around a known state $\rho_0$, it would
be impossible to devise an estimator that could do better asymptotically, 
than our estimator $\hat{\rho}_{n}$ on a small ball around $\rho$. 
More precisely:
\begin{theorem}  \label{automobiel}
Let $\hat{\rho}_{n}$ be the estimator defined above.
For any qubit state $\rho_{0}$ different from the totally mixed state, for 
any sequence of
estimators $\hat{\varrho}_n$, the following local asymptotic minimax result
holds for any $0<\epsilon<1/12$:
\begin{equation}
\label{resultat}
\limsup_{n\to\infty}\sup_{\|\rho - \rho_{0}\|_1 \leq n^{-1/2+\epsilon}}
 n R(\rho,\hat{\rho}_{n}) \leq
 \limsup_{n\to\infty}\sup_{\|\rho - \rho_{0}\|_1 \leq n^{-1/2+\epsilon}}
n R(\rho,\hat{\varrho}_n).
\end{equation}
Let \mbox{$(\mu_{0}, 1-\mu_{0})$}  be the eigenvalues of $\rho_{0}$ with 
$\mu_{0}>1/2$. Then the local asymptotic minimax risk is
\begin{equation}
\label{risquefinal}
\limsup_{n\to\infty}\sup_{\|\rho - \rho_{0}\|_1 \leq n^{-1/2+\epsilon}}
n R(\rho,\hat{\rho}_{n}) = R_{\rm minimax}(\mu_{0}) =8\mu_{0}-4\mu_{0}^2  .
\end{equation}
\end{theorem}
\proof
We write the risk as the sum of two terms corresponding to the events
$E$ and $E^{c}$
that $\tilde{\rho}_{n}$ is inside or outside the ball of radius 
$n^{-1/2+\epsilon}$ around $\rho$. Recall that LAN is valid inside the ball. 
Thus
$$
R(\rho,\hat{\rho}_{n})  =
\mathbb{E} (\| \rho-\hat{\rho}_{n}\|_{1}^{2}\, \chi_{E^{c}} ) +
\mathbb{E} (\| \rho-\hat{\rho}_{n}\|_{1}^{2} \,\chi_{E} ),
$$
where the expectation comes from $\hat{\rho}_n$ being random. 
The distribution
of the result $\hat{\rho_n}$ of our measurement procedure applied 
to the true unknown
state $\rho$ depends on $\rho$.
We bound the first part by $R_{1}$ and the second part by $R_{2}$ 
as shown below.

$R_{1}$ equals $\mathbb{P}(E^{c})$ times the maximum error, which is $4$ 
since for any pair of density matrices $\rho$ and $\sigma$, we have 
$\| \rho-\sigma\|_{1}^{2}\leq 4$. Thus
$$
R_{1} = 4\mathbb{P} (\|\rho-\tilde{\rho}_{n} \|_{1} \geq n^{-1/2+\epsilon}).
$$
According to lemma \ref{lemma.small.probability} this probability goes to zero 
exponentially fast, therefore the contribution brought by this term can be 
neglected.

We can now assume that $\tilde \rho_{n}$ is in the range of validity of local 
asymptotic normality and we can write $\rho^{\otimes n}= \rho^{\bf u}_{n}$ 
with ${\bf u}$ the
local parameter around $\tilde{\rho}_{n}$. We get the following inequalities 
for the second term in the risk.
\begin{eqnarray}
\mathbb{E} (\| \rho-\hat{\rho}_{n}\|_{1}^{2}\, \chi_{E} )
\hspace{-1.5mm}&\leq&\hspace{-1.5mm}
\mathbb{E} \left[ \|\hat{\rho}_{n} -\rho\|_{1}^{2} ~\Big|~ \|\tilde 
\rho_{n} - \rho\|_1
\leq n^{-1/2+\epsilon} ~\right] \nonumber \\
\hspace{-1.5mm}&\leq&\hspace{-1.5mm}
\sup_{ \|\rho - \rho_{0}\| < n^{-1/2+\epsilon}}
\mathbb{E} \left[ \|\hat{\rho}_{n} -\rho\|_{1}^{2} ~\Big|~ \tilde 
\rho_{n} =\rho_{0} \right]
\nonumber\\
\hspace{-1.5mm}&\leq&\hspace{-1.5mm}
\sup_{ \|\rho - \rho_{0}\| < n^{-1/2+\epsilon}}
\mathbb{E}_{\rho^{\bf u}_{n}(\infty)}
\left[ \|\hat{\rho}_{n} -\rho\|_{1}^{2} ~\Big|~ \tilde \rho_{n} =
\rho_{0} \right]
 \nonumber \\
&&\hspace{-1.5mm}+\hspace{-3mm}
\sup_{ \|\rho - \rho_{0}\| < n^{-1/2+\epsilon}}
  \|\rho^{\bf u}_{n} (t) - \rho^{\bf u}_{n}(\infty)\|_{1} ~
  \sup_{ \hat{\bf u}_{n} } \|\hat{\rho}_{n} -\rho\|_{1}^{2}
  \nonumber  \\
\hspace{-1.5mm}&\leq&\hspace{-1.5mm}
\sup_{ \|\rho - \rho_{0}\| < n^{-1/2+\epsilon}}  
\mathbb{E}_{\rho^{\bf u}_{n}(\infty)}
\left[  \|\hat{\rho}_{n} -\rho\|_{1}^{2} ~\Big|~ 
\tilde \rho_{n} =\rho_{0}   \right]
\nonumber  \\
&&+
 c n^{-1+2\eta}\sup_{ \|\rho - \rho_{0}\| < n^{-1/2+\epsilon}}
 \|\rho_{n}^{\bf u}(t) - \rho^{\bf u}_{n}(\infty) \|_1 \nonumber\\
 \hspace{-1.5mm}&=& R_{2}.
 \label{error}
\end{eqnarray}
The first  two inequalities are trivial. In the third inequality we 
change the expectation from the one with respect to the probability 
distribution of our data
$\mathbb{P}_{\rho^{\bf u}_{n}(t)}$ to the probability distribution
$\mathbb{P}_{\rho^{\bf u}_{n}(\infty)}$. In doing so, an additional term
$\|\mathbb{P}_{\rho^{\bf u}_{n}(t)} - 
\mathbb{P}_{\rho^{\bf u}_{n}(\infty)}\|_{1}$ appears which is bounded 
from above by  $\|\rho^{\bf u}_{n} (t) - \rho^{\bf u}_{n}(\infty)\|_{1}$.
In the last inequality we can bound $ \|\hat{\rho}_{n} -\rho\|_{1}^{2}$ by
$c n^{-1+2\eta}$ for some constant $c$. Indeed from definitions \eqref{mod1} 
and \eqref{mod2} we know that $\| \hat{\rho}_{n} -\rho_{0} \|_{1} 
\leq c^{\prime} n^{-1/2 +\eta} $ and additionally we are under the assumption 
$\| \rho -\rho_{0}\|_{1}\leq  n^{-1/2 +\epsilon}$ with $\epsilon<\eta$.

For the following, recall that all our LAN estimates are valid uniformly
around any state $\rho^{\bf 0}=\tilde{\rho}$ as long as 
$\mu-1/2\geq \epsilon_2 > 0$. 
As we are working with $\rho$ different
from the totally mixed state and $\|\rho -\tilde \rho\| \leq
n^{-1/2+\epsilon}$, we know that for  big enough $n$, $\tilde \mu - 1/2\geq
\epsilon_2$ for any possible $\tilde \rho$. We can then apply the uniform 
results of the previous sections.

The second term in $R_{2}$ is $\mathcal{O}(n^{-5/4+ 3\eta +\delta},
n^{-3/2+5\eta+\delta })$ where $\delta>0$ can be chosen arbitrarily small. 
Indeed, in the end of section \ref{sec.timevolf}, we have proven that after 
time $t\geq \ln n$, we have
$  \|\rho_{n}^{\bf u}(t) -\rho^{\bf u}_{n} (\infty) \|_1 =  
\mathcal{O}(n^{-1/4+\eta+ \delta}, n^{-1/2+ 3\eta+\delta} )$.
The contribution to $nR(\rho,\hat{\rho}_{n})$ brought by this term will not 
count in the 
limit, as long as $\eta$ and $\epsilon$ are chose such that 
$1/12 > \eta>\epsilon$.

We now deal with the first term in $R_{2}$. We write $\rho$ in local 
parameterization around $\rho_{0}=\tilde{\rho}$ as $\rho_{{\bf u}_n/\sqrt{n}}$. 
We have
\begin{eqnarray}
\label{eq.dist.param}
\|\hat{\rho}_{n}-\rho\|_{1}^{2} &=&
 \| \rho_{{\bf u}/\sqrt{n}} - \rho_{{\bf \hat{u}_{n}}/\sqrt{n}}\|_1^2\notag\\
&=&
4\frac{(u_z-\hat{u}_z)^2 + (2\mu-1)^2 ((u_x-\hat{u}_x)^2 +
(u_y-\hat{u}_y)^2)}{n} \notag\\
&&+ \,\, \mathcal{O}(\|{\bf u} - \hat{\bf u}_{n} \|^3 n^{-3/2}).
\end{eqnarray}
The remainder term
$\mathcal{O}(\|{\bf u} - \hat{\bf u}_{n}\|^3 n^{-3/2})$ is negligible.
It is $\mathcal{O}(n^{3\eta-3/2})$ which does not contribute to $nR(\rho,\hat{\rho}_{n})$
 for $\eta<1/6$. This is because on the one hand we have asked for $
\|\tilde \rho_{n} - \rho\| < n^{-1/2+\epsilon} $, and on the other hand, 
we have
bounded our estimator $\hat{\bf u}_{n}$ by using \eqref{mod1} and \eqref{mod2}.

We now evaluate $\mathbb{E}_{\rho^{\bf u}_{n}(\infty)}
\left[ d({\bf u},\hat{\bf u}_{n})^{2} \right]$,
where $d({\bf u},\hat{\bf u}_{n})^{2}$ is defined as
\begin{equation}\label{eq.loss.function}
d({\bf u}, {\bf v})^{2} := 4\left[ (u_z- v_z)^2 + (2\mu-1)^2 ((u_x- v_x)^2 +
(u_y- v_y)^2)\right] .
\end{equation}
Note that the risk of $\hat{\bf u}_{n}$ is smaller than that of 
$\tilde{\bf u}_{n}$
(see discussion below \eqref{mod1} and \eqref{mod2}). Under the law 
$\mathbb{P}_{\rho^{\bf u}_{n}(\infty)}$ the estimator  $\tilde{\bf u}_{n}$ 
has a Gaussian distribution  as shown in \eqref{gaus1} and \eqref{gaus2} 
with fixed and known variance and unknown expectation. In statistics this 
type of model is known as a Gaussian shift experiment \cite{vanderVaart}. 
Using
\eqref{gaus1} and \eqref{gaus2}, we get
$\mathbb{E}_{\rho^{\bf u}_{n}(\infty)}\left[ (u_z - \hat{u}_z)^2\right] \leq
\mu(1-\mu)$ and $\mathbb{E}_{\rho^{\bf u}_{n}(\infty)}\left[ (u_i -
\hat{u}_i)^2\right] \leq
\mu/(2(2\mu-1)^2) $ for $i=x,y$. Substituting these bounds in 
\eqref{eq.dist.param},
we obtain \eqref{risquefinal}.

We will now show that the sequence $\smash{\hat{\rho}_{n}}$ is optimal in the local 
minimax sense: for any $\rho_{0}$ and any other sequence of estimators 
$\smash{\hat{\varrho}_{n}}$ we have
$$
R_{0} =
\limsup_{n\to \infty} \sup_{\|\rho -\rho_{0}\|_{1} \leq n^{-1/2+\epsilon}} 
nR(\rho, \hat \varrho_{n}) \geq
8\mu_{0} -4\mu_{0}^{2}.
$$
We will first prove that the right hand side is the  minimax risk 
$R_{\rm minimax} (\mu_{0})$ for the family of states 
$N^{\bf u}\otimes \phi^{\bf u}$ which is the limit of the local families 
$\rho^{\bf u}_{n}$ 
of qubit states centered around $\rho_{0}$. We then extend the result to 
our sequence of quantum statistical models $\rho^{\bf u}_n$.

The minimax optimality for $N^{\bf u}\otimes \phi^{\bf u}$ can be checked 
separately for
the classical and the quantum part of the experiment. For the quantum part
$\phi^{\bf u}$, the optimal measurement is known to be the heterodyne 
measurement.
A proof of this fact can be found in lemma 7.4 of \cite{Guta&Kahn}. 
For the classical part, which corresponds to the measurement of $L_{z}$, the
optimal estimator is simply the random variable $X\sim N^{\bf u}$ itself 
\cite{vanderVaart}. 


We now end the proof by using the other direction of LAN.
Suppose that there exists a better sequence of estimators $\hat{\varrho}_n$ 
such that
$$
R_{0}  < R_{\rm minimax} (\mu_{0}) = 8\mu_{0} - 4\mu_{0}^{2}.
$$ 
We will show that this leads to an estimator $\hat{u}$ of ${\bf u}$ for the 
family 
$N^{\bf u}\otimes \phi^{\bf u}$ whose maximum risk  is smaller than the  
minimax risk $R_{\rm minimax}(\mu_{0})$, which is impossible. 

By means of a beamsplitter one can divide the state 
$\phi^{\bf u}$ into two independent Gaussian modes, using a thermal state 
$\phi:=\phi^{0}$ as the second input. If $r$ and $t$ are the reflectivity 
and respective transmittivity of the beamsplitter ($r^{2}+t^{2}=1$), then 
the transmitted beam has state 
$\phi^{\bf u}_{\mathrm{tr}} = \phi^{t{\bf u}}$ and the reflected one 
$\phi^{\bf u}_{\mathrm{ref}} = \phi^{r{\bf u}}$.
By performing a heterodyne measurement on the latter, and observing the
classical part $N^{\bf u}$, we can localize ${\bf u}$ within a big ball around
the result $\tilde{\bf u}$ with high probability, in the 
spirit of lemma \ref{lemma.small.probability}. 
More precisely, for any small $\tilde\epsilon>0$ we can find $a>0$ big enough 
such that the risk contribution from unlikely $\tilde{\bf u}$'s is small
$$
\mathbb{E}(\|{\bf u} - \tilde{\bf u}\|^{2} \chi_{ \|{\bf u} - 
\tilde{\bf u}\|>a})  < \tilde\epsilon.
$$
Summarizing the localization step, we may assume that the parameter ${\bf u}$ 
satisfies $\|{\bf u}\| < a$ 
with an $\tilde{\epsilon}$ loss of risk, where $a= a(r,\tilde{\epsilon})$.

Now let $n$ be large enough such that $n^{\epsilon}> a$, then the parameter 
${\bf u}$ falls within the domain of convergence of the inverse map $S_{n}$ of 
theorem \ref{th.qlan} and by \eqref{eq.channel.conv.inverse} 
(with $\epsilon$ replacing $\eta$ and $\delta$ replacing $\epsilon$) we have 
$$
\| \rho^{t{\bf u}}_{n} - S(N^{t{\bf u}} \otimes \phi^{t{\bf u}} )\|_{1} 
\leq C n^{-1/4+\epsilon+ \delta},
$$
for some constant $C$.

Next we perform the measurement leading to the estimator $\hat{\varrho}_n$ and 
equivalently to an estimator $\hat{\bf u}_{n}$ of 
${\bf u}$. Without loss of risk we can implement the condition 
$\|{\bf u}\|<a$ into the estimator  $\hat{\bf u}_{n}$ in a similar fashion as 
in 
\eqref{mod1} and \eqref{mod2}. 
The risk of this estimation procedure for $\phi^{\bf u}$ is then bounded 
from above by the sum of three terms: the risk 
$n R_{\rho}(\hat{\varrho}_{n})/t^{2} $ coming from the qubit estimation, the 
error contribution from the map $S_{n}$ which is  
$a^{2} n^{-1/4 + \epsilon +\delta}$, and the localization risk contribution 
$\tilde{\epsilon}$. This risk bound uses the same technique as the third 
inequality of \eqref{error}. The second contribution can be made arbitrarily 
small by choosing $n$ large enough, for $\epsilon<1/4$. From our assumption 
we have $R_{0} < R_{\mathrm{minimax}}(\mu_{0})$ and we can choose $t$ close to one 
such that $R_{0} /t^{2}< R_{\mathrm{minimax}}(\mu_{0})$ and further choose 
$\tilde{\epsilon}$ such that  
$R_{0} /t^{2}+\tilde{\epsilon}< R_{\mathrm{minimax}}(\mu_{0})$.

In conclusion, we get that the risk for 
estimating ${\bf u}$ is asymptotically smaller that the risk of the 
heterodyne measurement combined with observing the classical part which is 
known to be 
minimax \cite{Guta&Kahn}. Hence no such sequence $\hat{\varrho}_{n}$ exists, 
and $\hat{\rho}_{n}$ is optimal.
\qed

{{\bf \noindent Remark}\quad
In theorem \ref{resultat}, we have used the risk function
$R(\rho , \hat{\rho}) = \mathbb{E}(d^2(\rho , \hat{\rho})) $,
with $d$ the $L_1$-distance
$d(\rho , \hat{\rho}) = \|\rho - \hat{\rho}\|_{1}$.
However,
the obtained results can easily
be adapted to \emph{any}
distance measure $d^2 (\rho_{\hat{\bf u}} , \rho_{\bf u} )$ which is locally
quadratic in $\hat{\bf u} - \bf{u}$, i.e.
$$
d^2(\rho_{\hat{\bf u}} , \rho_{\bf u} ) =
\sum_{\alpha, \beta = x,y,z} \gamma_{\alpha \beta}
(u_{\alpha} - \hat{u}_{\alpha}) (u_{\beta} - \hat{u}_{\beta})
+ \mathcal{O}(\|u - \hat{u}\|^{3})\,.
$$

For instance, one may choose
$d^2(\hat{\rho} , \rho) = 1 - F^2(\hat{\rho} , \rho)$
with the fidelity $F(\hat{\rho} , \rho) :=
\tr(\sqrt{\sqrt{\hat{\rho}} \rho \sqrt{\hat{\rho}}  })$.
For non-pure states, this is easily seen to be locally quadratic with
$$
\gamma =
\left(
\begin{array}{c c c}
(2 \mu_0 - 1)^2&0&0\\
0&(2 \mu_0 - 1)^2&0\\
0&0&\frac1{1 -(2 \mu_0 - 1)^2}\\
\end{array}
\right)     \,.
$$
For the corresponding risk function
$R_{F}(\rho,\hat{\rho}_{n}) :=
\mathbb{E}(1 - F^2(\rho,\hat{\rho}_{n}))$,
this yields
\begin{equation}
\label{risqueF}
\limsup_{n\to\infty}\sup_{\|\rho - \rho_{0}\|_1 \leq n^{-1/2+\epsilon}}
n R_{F}(\rho,\hat{\rho}_{n}) = \mu_0 + 1/4 \,,
\end{equation}
with the same asymptotically optimal $\hat{\rho}$.
The asymptotic rate $R_{F} \sim \frac{4\mu_0 + 1}{4n}$ was found earlier in
\cite{Bagan&Gill}, using different methods.
}
\section{Conclusions}

In this chapter, we have shown two properties of
quantum local asymptotic normality (LAN) for qubits.
First of all, we have seen that its radius
of validity is arbitrarily close to $n^{-1/4}$ rather than $n^{-1/2}$.
And secondly, we have seen how LAN
can be implemented physically, in a quantum optical setup.

We use these properties to construct an
asymptotically optimal estimator
$\hat{\rho}_n$ of the qubit state
$\rho$, provided that we are given $n$ identical copies
of $\rho$.
Compared with other optimal estimation methods
\cite{Bagan&Gill,Hayashi&Matsumoto}, our measurement technique
makes a significant step in the direction  of an experimental implementation.

%
The construction and optimality of $\hat{\rho}_n$
are shown in three steps.
\begin{itemize}
\item[I]
In the preliminary stage, we perform
measurements of $\sigma_x$, $\sigma_y$ and $\sigma_z$
on a fraction
$\tilde{n} = n^{1-\kappa}$ of the $n$ atoms.
As shown in section \ref{sec.estimation},
this yields a rough estimate $\tilde{\rho}_{n}$
which lies within a distance $n^{-1/2 + \epsilon}$
of the true state $\rho$ with high probability.
\item[II]
In section \ref{sec.lan}, it is shown that
local asymptotic normality holds within a
ball of radius $n^{-1/2 + \eta}$ around $\rho$ ($\eta>\epsilon$).
This means that locally, in the limit $n \rightarrow \infty$,
all statistical problems concerning the $n$
identically prepared qubits are equivalent
to statistical problems concerning a
Gaussian distribution $N^{\bf u}$
and its quantum analogue,
a displaced thermal state $\phi^{\bf u}$ of the harmonic oscillator.
\end{itemize}

Together, I and II imply that the principle of
LAN 
has been extended to a global setting.
It can now be used for
a wide range of asymptotic statistical problems,
including the global problem of state estimation.
Note that this hinges on the rather subtle extension
of the range of validity of LAN to neighbourhoods
of radius larger than $n^{-1/2}$.

\begin{itemize}
\item[III] LAN provides an abstract equivalence between
the n-qubit states $\rho^{\otimes n}_{{\bf u}/\smash{\sqrt{n}}}$
on the one hand, and on the other hand
the Gaussian states $N^{\bf u} \otimes \phi^{\bf u}$.
In sections \ref{sec.timevolf} and \ref{finmes} it is shown that this 
abstract equivalence
can be implemented physically by two consecutive couplings
to the electromagnetic field.
For the particular problem of state estimation,
homodyne and heterodyne detection on the
electromagnetic field then yield the data from 
which the optimal estimator $\hat{\rho}_{n}$ is computed.
\end{itemize}

Finally, in section \ref{sec.endresult}, it is shown that
the estimator $\hat{\rho}_n$, constructed above,
is optimal in a local minimax sense.
Local here means that optimality holds
in a ball of radius slightly bigger than $n^{-1/2}$ around
\emph{any} state $\rho_0$ except the tracial state. That is,
even if we had known
beforehand that the true state
lies within this ball around $\rho_0$,
we would not have been able
to construct a better estimator than $\hat{\rho}_n$,
which is of course
independent of $\rho_0$.

For this asymptotically optimal estimator, we have shown that
the risk $R$ converges to zero at rate
$R(\rho,\hat{\rho}_{n}) \sim \frac{8\mu_{0}-4\mu_{0}^2}{n}$,
with $\mu_0 > \half$ an eigenvalue of $\rho$.
More precisely, we have
$$
\limsup_{n\to\infty}\sup_{\|\rho - \rho_{0}\|_1 \leq n^{-1/2+\epsilon}}
n R(\rho,\hat{\rho}_{n}) = 8\mu_{0}-4\mu_{0}^2.
$$
The risk is defined as
$R(\rho,\hat{\rho}) = \mathbb{E}(d^2(\rho,\hat{\rho}))$,
where we have chosen
$d(\hat{\rho} , \rho)$ to be the
$L_1$-distance $\|\hat{\rho} - \rho\|_{1} := \tr(|\hat{\rho} - \rho|)$.
This seems to be a rather natural choice because of its
direct physical significance
as the worst case difference between the probabilities
induced by $\hat{\rho}$ and $\rho$ on a single event.

Even still, we emphasize that the same procedure can be applied to a
wide range of other risk functions.
Due to the local nature
of the estimator $\hat{\rho}_n$ for large $n$,
its rate of convergence in a risk $R$
is only sensitive
to the lowest order Taylor expansion of $R$
in local parameters $\hat{\bf u} - \bf{u}$.
The procedure can therefore easily be adapted to
other risk functions, provided that the distance measure
$d^2 (\rho_{\hat{\bf u}} , \rho_{\bf u} )$
is locally
quadratic in $\hat{\bf u} - \bf{u}$.\medskip

{{\bf \noindent Remark}\quad
The totally mixed state ($\mu=1/2$) is a singular point in the
parameter space, and theorem \ref{th.qlan} does not apply in this case. 
The effect of the 
singularity is that the family of states \eqref{eq.displacedthermal} 
collapses to a 
single degenerate state of infinite temperature. However this phenomenon 
is only due 
to our particular parameterization, which was chosen for its convenience in 
describing the local neighbourhoods around arbitrary states, with the 
exception of the totally mixed state. Had we chosen a different 
parameterization, e.g. in terms of the
Bloch vector, we would have found that local asymptotic normality holds 
for the totally mixed state as well, but the limit experiment is different: 
it consists of a three dimensional 
{\it classical} Gaussian shift, each independent component corresponding to 
the local change in the Bloch vector along the three possible directions. 
Mathematically, the optimal measurement strategy in this case is just to 
observe the classical variables. However this strategy cannot be implemented 
by coupling with the field since this coupling becomes singular 
(see equation \eqref{eq.cocycle.j}). 
These issues become more important for higher dimensional systems where the 
eigenvalues may exhibit more complicated multiplicities, and will be dealt 
with in that context.  
}

\let\thesectionsav\thesection
\let\thesubsectionsav\thesubsection
\renewcommand{\thesection}{\Alphminseven{section}}
\renewcommand{\thesubsection}{\Alphminseven{section}.\arabic{subsection}}

\section{Appendix: Proof of \texorpdfstring{Theorem \ref{th.qlan}}{a
 Theorem}}\label{sec.proof.qlan}

Here we give the technical details of the proof of local asymptotic normality
with ``slowly growing'' local neighbourhoods $\| {\bf u}\| \leq
n^{\eta}$, with $\eta<1/4$.
We start with the map $T_n$.

\subsection{\texorpdfstring{The Map $T_n$}{The First Map}}

 Let us define, for $0 < \epsilon < (1/4 -\eta) $ the interval
\begin{equation}\label{eq.typical.j}
 \mathcal{J}_{n} = \left\{ j \,;\, 
 (\mu-1/2)n -n^{1/2+\epsilon}\leq j\leq (\mu-1/2)n + n^{1/2+\epsilon}\right\}.
\end{equation}
Notice that $j\in \mathcal{J}_n$ satisfies $2j\geq \epsilon_{2} n$ for all
$\mu-1/2\geq \epsilon_2$ and $n$ big enough, independently of $\mu$.
Then $ \mathcal{J}_n$ contains the relevant values of $j$, uniformly for
$\mu-1/2\geq \epsilon_2$:
\begin{equation}
\label{tail}
\lim_{n\to\infty} p_{n,{\bf u}}(\mathcal{J}_{n})=1 - \mathcal{O}(n^{-1/2+\epsilon}).
\end{equation}
This is a consequence of Hoeffding's inequality applied to the 
binomial distribution, and
recalling that $p_{n,{\bf u}}(j) = B(n/2 + j) (1 + \mathcal{O}(n^{-1/2+\epsilon}))$ for
$j\in \mathcal{J}_n$.

We upper-bound $ \| T_n(\rho_n^{\bf u}) -N^{\bf u}\otimes \phi^{\bf u}\|$ by 
the sum
\begin{equation}
 3
\sum_{j\not\in \mathcal{J}_{n}} p_{n,j}^{\bf u}
+
\left \| N^{\bf u}  - \sum_{j\in\mathcal{J}_{n} } p_{n, {\bf u}}(j) 
\tau_{n,j}     \right\|_{1}
+
\sup_{j\in \mathcal{J}_n} \| V_{j} \rho_{j,n}^{\bf u} V_{j}^{\dagger} - 
\phi^{\bf u}\|_{1}.
\end{equation}
%
The first two terms are ``classical'' and converge to zero uniformly over $\|
{\bf u}\| \leq n^{\eta}$. For the first term, this is (\ref{tail}), while the
second term converges uniformly on $\mu-1/2\geq \epsilon_2$ at rate 
$n^{\eta-1/2}$ \cite{Kahn}. The third term can be analyzed as in 
proposition 5.1 of \cite{Guta&Kahn}:
\begin{eqnarray}\label{eq.two.terms}
\left\|V_{j} \rho^{\bf u}_{n,j} V_{j}^{\dagger} - \phi^{\bf u}\right\|_{1} & 
\leq &
\left\| \rho^{\bf u}_{n,j} - V_{j}^{\dagger} \phi^{\bf u}V_{j}\right\|_{1} +
\left\| \phi^{\bf u} - P_{j} \phi^{\bf u} P_{j} \right\|_{1},
\end{eqnarray}
where $P_{j} := V_{j} V_{j}^{\dagger}$ is the projection onto the image of
$V_{j}$. We will show that both terms on the right side go to zero uniformly at
rate $n^{-1/4+\eta+\epsilon}$ over
$j\in \mathcal{J}_{n}$ and
$\| {\bf u}\| \leq n^{\eta}$. The trick is to note that displaced thermal 
equilibrium states are Gaussian mixtures of coherent states
\begin{equation}\label{eq.displaced.thermalstate.coherent}
\phi^{{\bf u}} = \frac{1}{\sqrt{2\pi s^2}}
\int e^{- |{\bf z}-\sqrt{2\mu -1}\alpha_{\bf u}|^{2}/2s^{2}}
\left( | {\bf z} \rangle \langle {\bf z} | \right) d^{2} {\bf z},
\end{equation}
where $s^{2}:= (1-\mu)/(4\mu-2)$.

The second term on the right side of \eqref{eq.two.terms} is bounded from 
above by
$$
 \frac{1}{\sqrt{2\pi s^2}}
\int e^{- |{\bf z}-\sqrt{2\mu -1}\alpha_{\bf u}|^{2}/2s^{2}}
\|  | {\bf z} \rangle \langle {\bf z} |  - P_{j}  | {\bf z} \rangle 
\langle {\bf z} | P_{j} \|_{1} \,
d^{2} {\bf z},
$$
which after some simple computations can be reduced (up to a constant) to
\begin{equation}\label{zwaardvis}
\int e^{- |{\bf z}|^{2}/2s^{2}}
\| P_{j}^{\perp} |{\bf z}+\sqrt{2\mu -1}\alpha_{\bf u}\rangle\| \,
d^{2} {\bf z}.
\end{equation}

We now split the integral. the first part is integrating over $|{\bf z}|\geq
n^{\eta+\delta} $ with $0<\delta <1/4 -\eta/2 $. The integral is 
dominated by the Gaussian and its value is 
$\mathcal{O}(e^{-n^{2(\eta + \delta)}/(2s^2)}) $.
The other part is bounded by the supremum over
$|{\bf z}|\leq 2n^{\eta+ \delta}$ (as $\|{\bf u}\|\leq n^{\eta}$) of $\|
P_{j}^{\perp} |{\bf z}\rangle\| $. Now $\| P_{j}^{\perp} 
|{\bf z}\rangle\| \leq
|{\bf z}|^{j}/\sqrt{j!} =\mathcal{O}(e^{-  n (1/2-\eta- 2\delta) })$ uniformly on 
$j\in \mathcal{J}_n$,
for any
$\mu-1/2\geq \epsilon_2$ since then $2 j\geq \epsilon_2 n$.

The same type of estimates apply to the first term:
\begin{eqnarray}
\left\| \rho^{\bf u}_{n,j} -  V_{j}^{\dagger}\phi^{\bf u}V_{j} \right\|_{1} =
\left\|
\mathrm{Ad} \left[ U_{j} \left({ \textstyle \frac{\bf u}{\sqrt{n}}}\right) \right]
\left( \rho^{\bf 0}_{n,j} \right)
-V_{j}^{\dagger}  \phi^{\bf u} V_{j}\right\|_{1} &\leq& \nonumber \\
\left\|\rho^{\bf 0}_{n,j} - V_{j}^{\dagger} \phi^{\bf 0} V_{j} \right\|_{1} +
\left\|
 \mathrm{Ad} \left[ U_{j} \left({ \textstyle \frac{\bf u}{\sqrt{n}}}\right) \right]
\left( V^{\dagger}_{j} \phi^{\bf 0} V_{j}\right)
-V_{j}^{\dagger} \phi^{\bf u} V_{j}
\right\|_{1} &\hspace{-8mm}.&\label{eq.triangle.ineq.}
\end{eqnarray}
The first term on the right side does not depend on ${\bf u}$. From the proof
of lemma 5.4 of \cite{Guta&Kahn}, we know that
\begin{align*}
\left\|\rho^{\bf 0}_{n,j} - V_{j}^{\dagger} \phi^{\bf 0} V_{j} \right\|_{1} 
\leq
\left(\frac{p^{2j+1}}{1-p^{2j+1}}+ p^{2j+1}\right)
\end{align*}
with $p= (1-\mu)/\mu$. Now the left side is of the order $p^{2j+1}$ which
converges exponentially fast to zero uniformly on $\mu-1/2\geq \epsilon_2$ and 
$j\in
\mathcal{J}_{n}$.

The second term of \eqref{eq.triangle.ineq.} can be bounded again by a 
Gaussian integral
\begin{eqnarray}\label{eq.gaussian.integral}
&&
\frac{1}{\sqrt{2\pi s^{2}}}
\int e^{- |{\bf z}|^{2}/2s^{2}} \| \Delta({\bf u}, {\bf z}, j) 
\|_{1} d^{2}{\bf z},
\end{eqnarray}
where the operator $\Delta({\bf u}, {\bf z}, j) $ is given by
\begin{eqnarray*}
\Delta({\bf u}, {\bf z}, j) &:=&
\mathrm{Ad} \left[ U_{j} \left({ \textstyle \frac{\bf u}{\sqrt{n}}}\right) \right]
\left( V^{\dagger}_{j} |{\bf z}\rangle \langle{\bf z} | V_{j}\right)
\\ & & \hspace{-2.5mm}-  V_{j}^{\dagger} \,
| {\bf z} + \sqrt{2\mu - 1}\alpha_{\bf u} \rangle \langle  {\bf z} + 
\sqrt{2\mu -1}\alpha_{\bf u} |  \,V_{j}.
\end{eqnarray*}

Again, we split the integral along $\|{\bf z}\|\geq n^{\eta + \delta}$. 
The outer
part converges to zero faster than any power of $n$, as we have already seen.
The inner integral, on the other hand, can be bounded uniformly over $\|{\bf
u}\|\leq n^{\eta}$, $\mu-1/2\geq \epsilon_2$  and $j\in \mathcal{J}_n$ 
by the supremum of $\| \Delta({\bf
u}, {\bf z}, j)\|_1$  over
$|{\bf z}|\leq 2 n^{\eta + \delta}$, $\mu-1/2\geq \epsilon_2$,  
$j\in \mathcal{J}_n$ and $\|{\bf u}\|\leq n^{\eta}$.

Let the vector $\tilde{\bf z}\in\mathbb{R}^{2}$ be such that $\alpha_{\tilde{\bf z}}= 
{\bf z}/\sqrt{2\mu-1}$, and  introduce the notation $\psi(n,j, {\bf v}) =  V_{j}U_{j}
( {\bf v}/\sqrt{n}) |j,j\rangle$. Then, up to a $\sqrt{2}$ factor,  
$ \| \Delta({\bf u}, {\bf z}, j) \|_{1}$ is bounded from above by\vspace{-1mm}
\begin{eqnarray}
&
 \Big\|
\psi(n,j, \tilde{\bf z}) -
|{\bf z}\rangle
\Big\|\,\,+& \nonumber\\
&
\Big\|
\psi(n,j, {\bf u}+\tilde{\bf z} ) -
 | {\bf z} + \sqrt{2\mu -1 }\alpha_{\bf u}\rangle
 \Big\|\,\,+& \nonumber\\
 &
 \Big\|
U_{j} \left(\frac{\bf u}{\sqrt{n}}\right) U_{j} \left(\frac{\tilde{\bf z}}
{\sqrt{n}}\right)
|jj \rangle  -
 U_{j} \left(\frac{{\bf u}+ \tilde{{\bf z}}}{\sqrt{n}}\right)   | jj \rangle
\Big\|\,.&
\label{eq.3terms}
\end{eqnarray}
This is obtained by adding and subtracting the states
$|\psi( n,j,\tilde{\bf z})\rangle \langle \psi( n,j,\tilde{\bf z})|$ and
$|\psi( n,j, {\bf u}+\tilde{\bf z})\rangle \langle \psi( n,j, 
{\bf u}+\tilde{\bf z})|$ 
and using the fact that 
$\| |\psi \rangle\langle \psi | - |\phi \rangle \langle\phi | \|_{1} 
= \sqrt{2}\|\psi-\phi\|$ for normalized vectors $\psi,\phi$.

The first two terms are similar, we want to dominate
them uniformly. We replace ${\bf u} + \tilde{\bf z}$ by $\tilde{\bf z}$ with
$|{\bf z}|\leq 2n^{\eta + \delta}$, and then write:
\begin{eqnarray}\label{eq.splitsum}
\left\|
\psi(n,j, \tilde{\bf z}) - |{\bf z}\rangle \right\|^{2}
&=&
\sum_{k=0}^{\infty} |\langle k| \psi(n,j, \tilde{\bf z})\rangle -
\langle k|{\bf z}\rangle |^{2}\nonumber\\
&\leq&
\sum_{k=0}^{r-1} |\langle k| \psi(n,j, \tilde{\bf z})\rangle -
\langle k|{\bf z}\rangle |^{2}\nonumber \\
&& 
\hspace{-1.3em}\hspace{-2pt}+2
\sum_{k=r}^{\infty} \left(  |\langle k| \psi(n,j, \tilde{\bf z})\rangle |^{2} 
+
|\langle k|{\bf z}\rangle |^{2}\right)\,.
\end{eqnarray}
If ${\bf z} = |{\bf z}| e^{i\theta}$ then we have \cite{Hayashi&Matsumoto}
\begin{eqnarray*}
\langle k| \psi(n,j, \tilde{\bf z})\rangle
&=&
\sqrt{\binom{2j}{k}} 
\left(\sin(|{\bf z}|/\sqrt{n})e^{i\theta}\right)^k
\left(\cos(|{\bf z}|\sqrt{n})\right)^{2j-k},\\
\langle k|  {\bf z}\rangle &=&
\exp\left(-\frac{(2\mu-1)|{\bf z}|^2}{2}\right)
\frac{\left(e^{i\theta}|{\bf z}|\sqrt{2\mu-1}\right)^k}{\sqrt{k!}}.
\end{eqnarray*}
In \eqref{eq.splitsum}, we choose $r= n^{2\eta + \epsilon_3}$ with
$\epsilon_3$ such that it satisfies the conditions
$ 2\delta+2\eta+\epsilon <2\eta+\epsilon_3 +\epsilon< 1/2 $
and $\eta+\epsilon_3 <1/4$. The tail sums are then of the order
\begin{eqnarray*}
\sum_{k=r}^{\infty} |\langle k | {\bf z} \rangle |^{2}
&\leq&
\frac{ |{\bf z}|^{2 r} }{r !}
\leq
\frac{(2n^{(\eta+\delta)})^{2 n^{2\eta + \epsilon_3}}}{(n^{2\eta
+\epsilon_3})!}\\
&=&\mathcal{O}\left(\exp(-n^{2\eta+\epsilon_3})\right),\\
\sum_{k=r}^{\infty} |\langle k| \psi(n,j,\tilde{\bf z})\rangle |^{2} 
&\leq&
\sum_{k=r}^{j} \left( \frac{|{\bf z}|^{2}}{n} \right)^{k}
\frac{(2j)!}{(2j-k)! k!} \leq n \frac{|{\bf z}|^{2r}}{r!}\\
&=& \mathcal{O}\left(\exp(-n^{2\eta+\epsilon_3})\right).
\end{eqnarray*}
For the finite sums we use the following estimates which are uniform over
all $|{\bf z}|\leq 2n^{\eta+\delta}$, $k\leq r$, $j\in \mathcal{J}_{n}$:
\begin{align*}
\sqrt{\binom{2j}{k}} & =
\frac{((2\mu-1)n)^{k/2}}{\sqrt{k!}}(1+\mathcal{O}(n^{-1/2+\epsilon +2\eta +
\epsilon_3})), \\
(\sin(|{\bf z}|/\sqrt{n}))^k
& = (|{\bf z}|/\sqrt{n})^k(1+\mathcal{O}(n^{4 \eta +  \epsilon_3+ 2\delta -1} )) ,\\
(\cos(|{\bf z}|/\sqrt{n}))^{2j-k} &=
\exp\left(-\frac{(2\mu-1)|{\bf z}|^2}{2}\right) (1+ \mathcal{O}(n^{2\eta-1/2 + \epsilon
+2\delta})),
\end{align*}
where we have used on the last line that $(1+x/n)^n = \exp(x) (1+
\mathcal{O}(n^{-1/2}x))$ for $x\leq n^{1/2-\epsilon_4}$ (cf. \cite{Kahn}). This is enough
to show that the finite sum converges uniformly to zero at rate
$\mathcal{O}(n^{2\eta-1/2 + \epsilon
+\epsilon_3})$ (at least if $\epsilon_3$ is small enough). 
The first two terms in \eqref{eq.3terms} decrease as the square root of this,
that is $\mathcal{O}(n^{\eta-1/4+\epsilon/2 + \epsilon_3/2})$.

Notice that the error terms depend on $\mu$ only through $j$, and that $2 j\geq
\epsilon n$ for $\mu-1/2\geq \epsilon_2$. Hence they are uniform in $\mu$.

We now pass to the third term of \eqref{eq.3terms}. 
By direct computation, it can be shown that if we consider two general elements
$\exp(i  X_{1} )$ and $\exp(i X_{2})$ of $\mathrm{SU}(2)$, with $X_{i}$ 
self-adjoint elements of $M(\mathbb{C}^{2})$,
then
\begin{equation}
e^{-i(X_{1}+X_{2})} e^{iX_{1}} e^{iX_{2}} e^{[X_{1},X_{2}]/2} =
\mathbf{1} + \mathcal{O}( X_{i_{1}} X_{i_{2}} X_{i_{3}}),
\end{equation}
where the $\mathcal{O}(\cdot)$ contains only third order terms in $X_{1}, X_{2}$. If 
$X_{1},X_{2}$ are in the linear span of $\sigma_{x}$ and $\sigma_{y}$ then
all third order monomials are such linear combinations as well.

In particular we get that for ${\bf z},{\bf u}\leq n^{\eta+\epsilon_3}$:
\begin{eqnarray}
U(\beta)&:=&
U\left(-{ \textstyle \frac{{\bf u}+{\bf v}}{\sqrt{n}}}\right)
U\left({ \textstyle \frac{\bf u}{\sqrt{n}}}\right)
U\left({ \textstyle \frac{\bf v}{\sqrt{n}}} \right)
\exp(i (u_{x} v_{y}-u_{y}v_{x})\sigma_{z}/n)\nonumber
\\[2mm]
&=& \begin{bmatrix}
1 + \mathcal{O}(n^{-2+4\eta+4\epsilon_3}) & \mathcal{O}(n^{-3/2 + 3\eta +3\epsilon_3}) \\
\mathcal{O}(n^{-3/2 + 3\eta + 3\epsilon_3}) & 1 + \mathcal{O}(n^{-2+ 4\eta +4\epsilon_3})
\end{bmatrix}\label{eq.curvature}.
\end{eqnarray}

Finally, using the fact that $|j,j \rangle$ is an eigenvector of 
$L_{z}$, the third term in \eqref{eq.3terms} can be written as
$$
\| |j,j\rangle\langle j,j| - U_j(\beta) |j,j\rangle\langle j,j|
U_j(\beta)^{\dagger}\|
$$
and both states are pure, so it suffices to show that the scalar product 
converges 
to one uniformly. Using \eqref{eq.curvature} and the expression of
$\langle j| U_j (\beta)|j\rangle $ \cite{Hayashi&Matsumoto} we get, 
as $j\leq n$,
$$
\langle j,j| U_j (\beta)|j,j\rangle = 
\left[ U(\beta)_{1,1}\right]^j =1 + \mathcal{O}(n^{-1
+ 4\eta + 4\epsilon_3}),
$$
which implies that the third term in \eqref{eq.3terms} is of order
$\mathcal{O}(n^{-1+ 4\eta + 4\epsilon_3})$. By choosing $\epsilon_{3}$ and $\epsilon$
small enough, we obtain that all terms used in bounding 
(\ref{eq.gaussian.integral}) are uniformly
$\mathcal{O}(n^{-1/4+ \eta + \epsilon  })$ for any $\epsilon>0$.

This ends the proof of
convergence \eqref{eq.channel.conv.} from the $n$ qubit state to the 
oscillator.

\subsection{The \texorpdfstring{Map $S_{n}$}{Second Map}}\label{app.a2 }

The opposite direction \eqref{eq.channel.conv.inverse} does not require much
additional estimation, so we will only give an outline of the argument.

Given the state $N^{\bf u}\otimes \phi^{\bf u}$, we would like to map it 
into $\rho^{\bf u}_{n}$ by means of a 
completely positive map $S_{n}$, or at least come as close as possible.

Let $X$ be the classical random variable with probability distribution 
$N^{\bf u}$.
With $X$ we generate a random $j\in \mathbb{Z}$  as follows
$$
j(X) = [ \sqrt{n}X + n (\mu-1/2) ].
$$
This choice is evident from the scaling properties of the probability 
distribution
$p^{\bf u}_{n}$ which we want to reconstruct. Let $q^{\bf u}_{n}$ be the 
probability distribution of $j(X)$.
By classical local asymptotic normality results we have the convergence
\begin{equation} \label{wildzwijn}
\sup_{\|{\bf u} \|\leq n^{\eta}} \| q^{\bf u}_{n} - p^{\bf u}_{n}\|_{1} = 
\mathcal{O}(n^{\eta-1/2}).
\end{equation}

Now, if the integer $j$ is in the interval $\mathcal{J}_{n}$, then
we prepare the $n$ qu\-bits in block diagonal state with the only non-zero
block corresponding to the $j^{\mathrm{\,th}}$ irreducible representation of 
$\mathrm{SU}(2)$:
$$
\tau^{\bf u}_{n,j} := \left(V^{\dagger}_{j} \phi^{\bf u} V_{j} +
\tr(P_{j}^{\perp} \phi^{\bf u}) \mathbf{1} \right)  \otimes 
\frac{\mathbf 1}{n_{j}}.
$$
The transformation $\phi^{\bf u}\mapsto \tau^{\bf u}_{n,j} $ is trace
preserving and completely positive~\cite{Guta&Kahn}.

If $j\notin\mathcal{J}_{n}$ then we may prepare the qubits in an arbitrary 
state which
we also denote by $\tau^{\bf u}_{n,j}$.
The total channel $S_{n}$ then acts as follows
$$
S_{n}: N^{\bf u} \otimes \phi^{\bf u} \mapsto  \tau^{\bf  u}_{n} := 
\bigoplus_{j=0,1/2}^{n/2} \,
 q^{\bf u}_{n,j} \tau^{\bf u}_{n,j}.
$$

We estimate the error $\| \rho^{\bf u}_{n} -   \tau^{\bf  u}_{n} \|_{1}$ as
$$
\| \rho_{n}^{\bf u} -   \tau^{\bf  u}_{n} \|_{1} \leq
\| q^{\bf u}_{n} - p^{\bf u}_{n}\|_{1} +
2\mathbb{P}_{ p^{\bf u}_{n}} (j\notin\mathcal{J}_{n}) +
\sup_{j\in\mathcal{J}_{n}} \|\tau_{n,j}^{\bf u} - \rho_{n,j}^{\bf u}\|_{1}
$$
The first term on the r.h.s. is $\mathcal{O}(n^{\eta - 1/2})$ (see \eqref{wildzwijn}),
the second term
is $\mathcal{O}(n^{\epsilon - 1/2})$ (see \eqref{tail}). As for the third term, we use
the triangle inequality to write, for $j\in\mathcal{J}_{n}$,
$$
\| \tau_{n,j}^{\bf u} - \rho_{n,j}^{\bf u}\|_{1} \leq
\| \tau_{n,j}^{\bf u} - V_{j}^{\dagger} \phi^{\bf u} V_{j}^{\dagger}\|_{1} +
\| V_{j}^{\dagger} \phi^{\bf u} V_{j}^{\dagger} - \rho_{n,j}^{\bf u}\|_{1}\,.
$$
The first term is $\mathcal{O}(e^{-n(1/2 - \eta - 2\delta)})$, according to the 
discussion
following equation \eqref{zwaardvis}.
The second term on the right is $\mathcal{O}(n^{-1/4 + \eta + \epsilon})$
according to equations \eqref{eq.triangle.ineq.} through \eqref{eq.curvature}.

Summarizing, we have
$\|S_{n}(N^{\bf u} \otimes \phi^{\bf u}) - \rho^{\bf u}_{n} \|_{1}
= \mathcal{O}(n^{-1/4 + \eta + \epsilon})
$, which establishes the proof in the inverse direction.
\qed

\section{Appendix: 
Proof of 
\texorpdfstring{Theorem \ref{th.unitary.evolution}}{Another Theorem}
} \label{kool}
We split the proof of theorem \ref{th.unitary.evolution}
in two steps.

\subsection{First Estimate}
We build up the state $\tilde{\rho}^{\bf u}_{j,n}$ by taking linear
combinations of number states $|m\rangle$ to obtain an
approximate coherent state $|{\bf z}\rangle$, and finally mixing such states
with a Gaussian distribution to get an approximate displaced thermal state.
Consider the approximate coherent vector
$
P_{\tilde m} |{\bf z}\rangle ,
$
for some fixed ${\bf z}\in\mathbb{C}$ and $\tilde m = n^{\gamma}$, with 
$\gamma$ to be fixed later.
Define the normalized vector
\begin{equation}\label{eq.normalized.vectors}
|\psi^{n}_{{\bf z},j} \rangle : = \frac{1}{\| P_{\tilde m} |{\bf z}\rangle\| }
\sum_{m=0}^{\tilde m} \frac{|{\bf z}|^{m}}{\sqrt{m!}} |m\rangle,
\end{equation}
We mix the above states to obtain
$$
\tilde{\rho}^{\bf u}_{j,n} :=
\frac{1}{\sqrt{2\pi s^{2}}}\int e^{-|{\bf z -\sqrt{2\mu-1} 
\alpha_{\bf u}}|^{2}/2s^2}
\left(|\psi^{n}_{{\bf z},j} \rangle \langle \psi^{n}_{{\bf z},j}| \right)\,
d^{2}{\bf z}.
$$
Recall that $s^2 = (1-\mu)(4 \mu - 2)$, and
$$
\phi^{\bf u} = \frac{1}{\sqrt{2\pi s^{2}}}\int e^{-|{\bf z -\sqrt{2\mu-1} 
\alpha_{\bf u}}|^{2}/2s^2}
\left(|{\bf z} \rangle \langle {\bf z}| \right)\,
d^{2}{\bf z}.
$$
From the definition of $|\psi^{n}_{{\bf z},j} \rangle$ we have
\begin{equation}
\label{bound1}
\| |\psi^{n}_{{\bf z},j} \rangle- |{\bf z}\rangle\| \leq \sqrt{2} 
\frac{|{\bf z}|^{\tilde
m}}{\sqrt{\tilde m!}} \wedge 2,
\end{equation}
which implies that
$\| \tilde{\rho}^{\bf u}_{j,n} - \phi^{\bf u}\|_{1} = 
\mathcal{O}(e^{-n^{2(\eta+\epsilon)}})$
for any $\epsilon> 0$, for any $\gamma\geq 2(\eta + \epsilon)$.
Indeed, we bound it 
by
$$
 \frac{\sqrt{2}}{\sqrt{\pi s^{2}}}\int e^{-|{\bf z} |^{2}/2s^2}
\left( \frac{|{\bf z} + \sqrt{2\mu-1} \alpha_{\bf u}|^{\tilde m}}
{\sqrt{\tilde m!}}\wedge \sqrt{2} \right)  \,
d^{2}{\bf z}\,,
$$
and split the integral into two parts. The integral over the domain
$|{\bf z}|\geq n^{\eta+\epsilon}$ is dominated by the Gaussian factor and
is $\mathcal{O}(e^{-n^{2(\eta+\epsilon)}})$. The integral over the disk 
$|{\bf z}|\leq n^{\eta+\epsilon}$ is bounded by the supremum of (\ref{bound1}) 
since the Gaussian integrates to one,
and is $\mathcal{O}(e^{-(\gamma/2-\eta-\epsilon) n^{\gamma}})$. In the last step we 
use Stirling's formula to obtain
$\log \left[ (n^{\eta+\epsilon})^{n^{\gamma}} / \sqrt{n^{\gamma}!}\right] 
\approx  (\eta+\epsilon-\gamma/2) n^{\gamma} \log n$. Note that the estimate 
is uniform with respect
to $\mu-1/2> \epsilon_{2}$ for any fixed $\epsilon_{2}>0$.

\subsection{Second Estimate} We now compare the evolved qubits
state $\tilde{\rho}^{\bf u}_{j,n}(t)$ and the evolved oscillator state 
$\phi^{\bf u}(t)$.
Let $\ket{\psi_{m,j}^{n}(t)} = U_{j,n}(t)\, |m\rangle \otimes \ket{\Omega}$ 
be the joint state at time $t$ when the initial state of the system is 
$|m\rangle$ corresponding to $|j, j-m\rangle$ in the $L_{z}$ basis notation. 
We choose the following approximation of $\ket{\psi_{m,j}^{n}(t)}$,
\begin{equation}\label{eq.approx.m}
\ket{\xi_{m,j}^{n}(t)} :=\sum_{i=0}^{m} c_{n}(m,i) \alpha_{i}(t) | m-i\rangle 
\otimes
|e^{-1/2 u} \chi_{[0,t]}(u) \rangle_{i},
\end{equation}
where $\alpha_{i}(t) = \exp((-m+i)t/2)$,
and we define the numbers $c_{n}(m,i)$ inductively by $c_{n}(m,0):=1$ and
$$c_{n}(m,i) := c_{n}(m,i-1) \sqrt{{ \textstyle \frac{2j-m+i}{2j_{n}}}} 
\sqrt{{ \textstyle \frac{m-i+1}{i}}}\,.$$
The $n$-particle vector $|f\rangle_{n} := f^{\otimes n}$ was defined in 
\eqref{eq.exponential}.
Note that
in particular, for $\mu-1/2>\epsilon_2$ and $j\in\mathcal{J}_{n}$, we have
$$
c_{n}(m,i)\leq  \sqrt{{ \textstyle \binom{m}{i} (1+ \frac{2}{\epsilon_2} 
n^{-1/2+\epsilon} )^{i}}}
\,.$$
We now apply the estimate \eqref{eq.approx.qsde.solution}. By direct 
computations, we get
\begin{equation}
\begin{gathered}
d\ket{\xi_{m,j}^{n}(t)}
=
 -\frac{1}{2}\sum_{i=0}^{m}
 c_{n}(m,i) \alpha_{i}(t)(m-i) |m-i \rangle \otimes |e^{-1/2 u} 
 \chi_{[0,t]}(u) \rangle_{i} dt
\\
+
 \sum_{i=1}^{m}  c_{n}(m,i) \alpha_{i - 1} (t) |m-i \rangle \otimes
 |e^{-1/2 u} \chi_{[0,t]}(u) \rangle_{i-1} \otimes_{s} |\chi_{[t,t+dt]}
 \rangle,
\label{eq.dxi}
\end{gathered}
\end{equation}
where
$$
f^{\otimes i} \,\otimes_{s}g := \sum_{k=1}^{i+1} f\otimes f \otimes \dots 
\otimes g\otimes \dots \otimes f.
$$
From the quantum stochastic differential equation, we obtain
\begin{equation*}
\begin{gathered}
G_{dt} \, \ket{\xi_{m,j}^{n}(t)} =
\sum_{i=0}^{m} c_{n}(m,i) \alpha_{i}(t) {\textstyle \sqrt{\frac{(m-i)(2j-m+i+1)}
{2j_{n} (i+1)}}} \, \times
\\
| m-i-1\rangle \otimes
|e^{-1/2 u} \chi_{[0,t]}(u) \rangle_{i} \otimes_{s} |\chi_{[t,t+dt]}\rangle
\\
-\half \sum_{i=0}^{m} c_{n}(m,i) \alpha_{i}(t) (m-i) 
{\textstyle \frac{2j-m+i +1} {2j_{n}}}  | m-i\rangle \otimes
|e^{-1/2 u} \chi_{[0,t]}(u) \rangle_{i} dt\,. 
\end{gathered}
\end{equation*}
%
%
In the first term of the right hand side of the above equation, 
we can replace
$c_{n}(m,i) \sqrt{\frac{(m-i)(2j-m+i+1)}{2j_{n} (i+1)}}$ by 
$c_{n}(m,i+1)$ and
thus we obtain the same sum as in the second term of the left side of 
\eqref{eq.dxi}. Thus
\begin{eqnarray*}
&\displaystyle
G_{dt} \ket{\xi_{m,j}^{n}(t)} - d \ket{\xi_{m,j}^{n}(t)} =& \\
&
\displaystyle
\half  \sum_{i=0}^{m-1} c_{n}(m,i) \alpha_{i}(t) (m-i)
{\textstyle \frac{2(j_{n}-j) +m-i-1}{2j_{n}}} |m-i \rangle \otimes
 |e^{-1/2 u} \chi_{[0,t]}(u) \rangle_{i} \, dt.&
\end{eqnarray*}
Then, using  
$c_{n}(m,i)\leq  
\sqrt{\binom{m}{i}  (1+  (2/\epsilon_2) n^{-1/2 + \epsilon} )^{i}}
$, we find that the expression $
\|G_{dt} \xi_{m,j}^{n}(t) - d\xi_{m,j}^{n}(t)\|$ is bounded from above by
\begin{equation*}
\half \Big[
\sum_{i=0}^{m-1}  {\textstyle \binom{m}{i}}  ((1+ n^{-1/2+\epsilon})(1-e^{-t}))^{i} 
e^{-(m-i)t} 
\left({\textstyle \frac{(2(j_{n}-j) +m-i-1)(m-i)}{2j_{n}}} \right)^{2} \Big]^{1/2}dt.
\end{equation*}
We have
$$
{\textstyle \frac{(2(j_{n}-j) +m-i-1)(m-i)}{2j_{n}}} = \mathcal{O}( m( n^{-1/2 +\epsilon} +
n^{-1}m))\,.
$$
Inside the sum, we recognize the binomial terms with the $m^{\mathrm{th}}$ term missing.
Thus the sum is
\begin{eqnarray*}
&&
\left(1+ n^{-1/2+\epsilon} -e^{-t} n^{-1/2+\epsilon}  \right)^{m} -
\left((1-e^{-t}) (1+n^{-1/2+\epsilon})\right)^{m} \\
&&
\leq
(1+n^{-1/2+\epsilon})^{m} ( 1 - (1-e^{-t})^{m} )\leq
(1+n^{-1/2+\epsilon})^{m} \, m e^{-t}.
\end{eqnarray*}
Then there exists a constant $C$ (independent of $\mu$ if $\mu-1/2\geq
\epsilon_2$) such that
$$
\|G_{dt} \xi_{m,j}^{n}(t) - d\xi_{m,j}^{n}(t)\| \leq  \frac{C}{2} e^{-t/2}
 m^{3/2} (n^{-1/2+\epsilon} +m n^{-1})   \left(1+\frac{2}{\epsilon_{2}} 
 n^{-1/2+\epsilon}\right)^{m/2}
$$
By integrating over $t$, we finally obtain
\begin{equation}\label{eq.approx.error.m}
\|\psi_{m,j}^{n}(t)- \xi_{m,j}^{n}(t)\| \leq  C
 m^{3/2} (n^{-1/2+\epsilon} +m n^{-1})   \left(1+ \frac{2}{\epsilon_{2}}
 n^{-1/2+\epsilon}\right)^{m/2}\hspace{-4mm}.
\end{equation}
Note that under the assumption $\gamma < 1/3-2\epsilon/3$, the right side 
converges to zero at rate $n^{3\gamma/2- 1/2+\epsilon}$ for all 
$m\leq\tilde{m}= n^{\gamma}$. Summarizing, the assumptions which we 
have made so far on $\gamma$ are
$$
2\eta+2\epsilon <\gamma<1/3-2\epsilon/3.
$$

Now consider the vector
$
\ket{\psi^{n}_{{\bf z},j}} 
$
as defined in \eqref{eq.normalized.vectors}, and introduce the notation
$\ket{\psi_{{\bf z},j}^{n}(t)} = 
U_{j,n}(t) \ket{\psi^{n}_{{\bf z},j} } \otimes \ket{\Omega}$.
Then based on \eqref{eq.approx.m}, we choose the approximate solution
$$
\ket{\xi_{{\bf z},j}^{n}(t)}= e^{-|{\bf z}|^{2}/2} \sum_{m=0}^{\tilde m} 
{\textstyle \frac{|{\bf z}|^{m}}{\sqrt{m!}}}
\sum_{i=0}^{m} c_{n}(m,i) \alpha_{i}(t) | m-i\rangle \otimes
|e^{-1/2 u} \chi_{[0,t]}(u) \rangle_{i}.
$$
Note that the vectors $\ket{\psi_{k,j}^{n}(t)}$ and $\ket{\xi_{k,j}^{n}(t)}$ 
live in the ``$k$-particle'' subspace of $\mathcal{H}_{j}\otimes 
\mathcal{F}(L^{2}(\mathbb{R}))$ and thus are orthogonal to all vectors 
$\ket{\psi_{p,j}^{n}(t)}$ and $\ket{\xi_{p,j}^{n}(t)}$ with $p\neq k$.
By \eqref{eq.approx.error.m},  the error is
\begin{equation}
\begin{gathered}
\big\|\psi_{{\bf z},j}^{n}(t) - \xi_{{\bf z},j}^{n}(t) \big\|
\leq 
\\
C e^{-|{\bf z}|^{2}/2}
\left( \sum_{m=0}^{\tilde m} {\textstyle \frac{|{\bf z}|^{2m}}{m!}}
 m^{3} (n^{-1/2+\epsilon} +m n^{-1})^{2}   \left(1+ {\textstyle \frac{2}{\epsilon_{2}}}
 n^{-1/2+\epsilon}\right)^{m}\right)^{1/2}
 +  {\textstyle \frac{|{\bf z}|^{2\tilde{m}}}{\tilde{m}!}} \leq \\
C\tilde{m}^{3/2}  (n^{-1/2+\epsilon} +\tilde{m} n^{-1})   
\left(1+ {\textstyle \frac{2}{\epsilon_{2}}}n^{-1/2+\epsilon}\right)^{\tilde m/2}
+ {\textstyle \frac{|{\bf z}|^{2\tilde{m}}}{\tilde{m}!}}\,. \label{eq.xi-psi.j}
\end{gathered}
\end{equation}
We now compare the approximate solution $\xi_{{\bf z} ,j}^{n}(t)$ with the 
``limit'' solution
$\psi_{{\bf z}}(t)$ for the oscillator coupled with the field as described 
in section \ref{sec.oscillator.field}. We can write
$$
\psi_{{\bf z}}(t) = e^{-|{\bf z}|^{2}/2} \sum_{m=0}^{\infty} 
{ \textstyle \frac{|{\bf z}|^{m}}{\sqrt{m!}}}
\sum_{i=0}^{m} \sqrt{{ \textstyle \binom{m}{i}}} e^{-(m-i)t/2} | m-i\rangle \otimes
|e^{-1/2 u} \chi_{[0,t]}(u) \rangle_{i}.
$$
Then
\begin{equation*}
\begin{gathered}
\| \xi_{{\bf z},j}^{n}(t) - \psi_{\bf z}(t) \|^{2}  = 
e^{-|{\bf z}|^{2}} \sum_{m=\tilde m}^{\infty} { \textstyle \frac{|{\bf z}|^{2m}}{m!}} \\
+\,\,
e^{-|{\bf z}|^{2}}  \sum_{m=0}^{\tilde m} { \textstyle \frac{|{\bf z}|^{2m}}{m!}} 
\sum_{i=0}^{m}
e^{-(m-i)t} \left|c_{n}(m,i) - \sqrt{{ \textstyle\binom{m}{i}}} 
\right|^{2} (1-e^{-t})^{i} \,.
\end{gathered}
\end{equation*}
Now
\begin{eqnarray*}
\left|c_{n}(m,i) -\sqrt{ { \textstyle\binom{m}{i}}} \right|^{2}
&\leq&
\left| c_{n}(m,i)^{2} - { \textstyle \binom{m}{i}}\right|\\
&\leq&
 { \textstyle\binom{m}{i}}
\left|1- \prod_{p=1}^{i} \left(1+ { \textstyle \frac{ 2(j- j_n) - m +p }{2j_{n}}} \right)
\right|
\\
&\leq&
C_{2} { \textstyle\binom{m}{i}}  m n^{-1/2+\epsilon},
\end{eqnarray*}
where $C_2$ does not depend on $\mu$ as long as $\mu-1/2\geq\epsilon_2$ 
(recall
that the dependence on $\mu$ is hidden in $j_n = (2\mu-1)n$). Thus
\begin{eqnarray}\label{eq.xi-psi}
\| \xi_{{\bf z},j}^{n}(t) - \psi_{\bf z}(t) \|^{2} &\leq&
C_{2} n^{-1/2+\epsilon}
e^{-|{\bf z}|^{2}}  \sum_{m=0}^{\tilde m} \frac{ m|{\bf z}|^{2m}}{m!}  +  
\frac{|{\bf z}|^{2\tilde m}}{\tilde m!} \nonumber \\
& \leq &
C_{2} n^{-1/2+\epsilon} |{\bf z}|^{2} + \frac{|{\bf z}|^{2\tilde m}}
{\tilde m!} .
\end{eqnarray}

From \eqref{eq.xi-psi.j} and \eqref{eq.xi-psi} we get
\begin{eqnarray*}
\|\psi_{{\bf z},j}^{n} (t) -  \psi_{\bf z}(t)\| & \leq & 
2 \wedge \Bigg[C\tilde{m}^{3/2}  (n^{-1/2+\epsilon} +\tilde{m} n^{-1})
\left(1+ \frac{2}{\epsilon_{2}}n^{-1/2+\epsilon}\right)^{\tilde m/2} \\
& &
+ 
\frac{|{\bf z}|^{2\tilde{m}}}{\tilde{m}!} +
 \left[C_{2} n^{-1/2+\epsilon} |{\bf z}|^{2} + \frac{|{\bf z}|^{2\tilde
m}}{\tilde m!}   \right]^{1/2} \Bigg]\,, 
\end{eqnarray*}
and we denote the r.h.s. by $E(\tilde{m}, n, {\bf z})$.
We now integrate the coherent states over the displacements ${\bf z}$ as we 
did in the case of local asymptotic normality in order to obtain the thermal 
states in which we are interested
$$
\tilde{\rho}^{\bf u}_{j,n} :=
\frac{1}{\sqrt{2\pi s^{2}}}\int e^{-|{\bf z -\sqrt{2\mu-1} 
\alpha_{\bf u}}|^{2} /2s^{2}}
\left(|\psi^{n}_{{\bf z},j} \rangle \langle \psi^{n}_{{\bf z},j}| \right)\,
d^{2}{\bf z}.
$$
If we define the evolved states
$$
\tilde{\rho}^{\bf u}_{j,n} (t) := U_{j,n}(t) \tilde{\rho}^{\bf u}_{j,n} 
U_{j,n}(t)^{\dagger}
\qquad
\mathrm{and}
\qquad
\phi^{\bf u} (t):= U(t) \phi^{\bf u} U(t)^{\dagger},
$$
then $\sup_{j\in\mathcal{J}_{n}} \sup_{\|{\bf u}\|\leq  n^{\eta}} \| 
\tilde{\rho}^{\bf u}_{j,n} (t) - \phi^{\bf u} (t)\|_{1}$
is bounded by
$$
\sup_{\|{\bf u}\|\leq  n^{\eta}}
 \frac{1}{\sqrt{\pi s^{2}}}\int e^{-|{\bf z -\sqrt{2\mu-1} 
 \alpha_{\bf u}}|^{2}/2s^{2}}
 E(\tilde{m}, n, {\bf z})\,
d^{2}{\bf z}.
$$

Here again we cut the integral in two parts. On $|{\bf z}|\geq n^{\eta +
\epsilon}$, the Gaussian dominates, and this outer part is less
than $e^{-n^{\eta+\epsilon}}$. The inner part is then dominated by $\sup_{|{\bf
z}|\leq n^{\eta + \epsilon}} E(\tilde{m},n,{\bf z})$. Now 
on the one hand, we want (\ref{eq.xi-psi.j}) to be small,
so that 
$\tilde m$ cannot be too big.
On the other hand, we want
${\bf z}^{2\tilde m}/ {\tilde m}!$ to go to zero. A choice which satisfies the 
condition
is $\gamma = 2\eta + 3\epsilon$. By renaming $\epsilon$, we then get
\[
 E(\tilde{m}, n, {\bf z}) = \mathcal{O}(n^{\eta-1/4+ \epsilon }, 
 n^{3\eta -1/2 +\epsilon}),
\]
for any small enough $\epsilon>0$.
Hence we obtain \eqref{eq.unitary.error}. 
\hfill$\Box$
\renewcommand{\thesection}{\thesectionsav}
\renewcommand{\thesubsection}{\thesubsectionsav}
\chapter{Bundles with a Lift of Infinitesimal Diffeomorphisms}
\label{ch:BLID}
\setsubdir{hoofdstukken/hoofdstuk6}
This chapter is not about algebra and quantum mechanics,
but about differential geometry and classical field theory.
We slightly extend the notion of a natural fibre bundle by requiring
diffeomorphisms of the base to lift to automorphisms of the bundle
only infinitesimally, i.e.\ at the level of the Lie algebra of vector fields.
Spin structures are natural only in this extended sense.
We classify the fibre bundles with this property, 
under the additional assumption of
a finite dimensional structure group.
This includes all spin structures, but only some $\mathrm{Spin}^{c}$
and $\spg$-structures. This classification links the gauge group $G$
to the topology of space-time.\medskip

{{\bf \noindent Remark} \quad
Many of the results in this chapter, including the main theorem 
\ref{grondsmakreefknol}, can be found in Lecomte's 
paper \cite{Lecomte}.
Our results were obtained independently.
}
\section{Introduction}\label{een}
All bundles are equal, but some are
more natural than others. 
For example, the naturality of
tangent bundles,
frame bundles and jet bundles is beyond all question.
A bundle is called `natural' if diffeomorphisms of the base 
lift to automorphisms of the bundle in a local fashion. 

Let us phrase this more carefully.
Any principal fibre bundle $\pi : P \rightarrow M$ determines an exact 
sequence of
groups
\begin{eqnarray}\label{flapperdrop}
1 \rightarrow 
\Gamma_{c}(\mathrm{Ad}(P))
\rightarrow
\mathrm{Aut}_{c}(P)
\rightarrow
\mathrm{Diff}^{P}_{c}(M) \rightarrow 1 \,. 
\end{eqnarray}
In this expression,
$\mathrm{Aut}_{c}(P)$ is the group of automorphisms of 
$P$ that are trivial outside a compact subset of $M$.
Its image $\mathrm{Diff}^{P}_{c}(M)$, the group of compactly supported
liftable diffeomorphisms of $M$, is an open subgroup of
$\mathrm{Diff}_{c}(M)$. In particular, it is big enough to contain 
its connected component of unity.  
We have identified the gauge 
group of vertical automorphisms with the group of sections of the adjoint bundle 
$\mathrm{Ad}(P)$, and $\Gamma_{c}(\mathrm{Ad}(P))$ is the subgroup of
compactly supported ones. 
%
%
%
%
\label{kroep} 
\begin{basdef}\label{defvannatbun}
A `natural' principal fibre bundle is
a principal fibre bundle $P$ for which
all diffeomorphisms are liftable and
(\ref{flapperdrop})
is split, together with 
a distinguished splitting 
homomorphism
$\Sigma : \mathrm{Diff}_{c}(M) \rightarrow \mathrm{Aut}_{c}(P)$.
Moreover, the splitting homomorphism
$\Sigma$ is required to be local in the sense that 
it should lift each diffeomorphism 
$\phi : U_1 \rightarrow U_2$
between open subsets of $M$
to a bundle isomorphism 
$\Sigma(\phi) : \pi^{-1}(U_1) \rightarrow \pi^{-1}(U_2)$
in a functorial fashion.
\end{basdef}

\noindent Natural bundles have been classified.
A theorem of Salvioli, Palais, Terng, Ep\-stein and Thurston
\cite{Sa,PT,ET} states that 
any natural fibre bundle is associated to 
the $k^{\mathrm{th}}$ order frame bundle $F^{k}(M)$.

In this chapter, we seek to extend the notion of a natural fibre bundle
in two separate ways. First of all, we do not require locality, but prove it.
And secondly, we only require diffeomorphisms of the base 
to lift to automorphisms of the bundle infinitesimally, i.e.\
at the level of Lie algebras. We will call the principal fibre bundles 
which are natural in this extended sense `infinitesimally natural',
and our main objective 
will be their classification.
 
Let us formulate this more precisely. For any principal fibre bundle, 
the sequence of groups (\ref{flapperdrop}) gives rise 
to the exact sequence of Lie algebras
\begin{equation}\label{rijtje}
0 \rightarrow \Gamma_{c}(\mathrm{ad}(P)) \rightarrow \Gamma_{c}(TP)^{G}
\rightarrow \Gamma_{c}(TM) \rightarrow 0\,,
\end{equation}
where `$c$' stands for `$0$ outside a compact subset of $M$'.

\label{kroepstrup}
\begin{basdef}\label{defvaninnatbun}
An `infinitesimally natural' principal fibre bundle is 
a principal fibre bundle $P$, together with 
a distinguished 
Lie algebras homomorphism
$\sigma :\Gamma_{c}(TM) \rightarrow \Gamma_{c}(TP)^{G}$
that splits the exact sequence of Lie algebras~(\ref{rijtje}).
\end{basdef}

\noindent We emphasize that we will not require $\sigma$ to be local, 
continuous or $C^{\infty}(M)$-linear, and it need not come from a map
of bundles. 
We only require $\sigma$ to be a 
homomorphism of Lie algebras.

The price we pay for this level of generality is that we must 
confine ourselves to fibre bundles with 
a finite dimensional structure group. This allows us to
first study principal fibre bundles, and then generalize 
the results to more general fibre bundles with a finite dimensional structure 
group, such as vector bundles.

There are three main reasons for wanting to extend the notion
of a natural bundle. First of all, natural bundles are meant
to describe geometric objects, but not all geometry is local.
One might feel that the locality requirement on $\Sigma$ is therefore
inappropriate.
The universal cover $\tilde{M} \rightarrow M$ for instance is
not a natural bundle, although 
its ties with the global topology of $M$ are unmistakable.
It does however allow for a lift of infinitesimal
diffeomorphisms, and is therefore infinitesimally natural. 
 
The second reason comes from field theory. 
Fields are described by sections 
of a fibre bundle $F$ over space-time $M$, where $F$ is 
associated to some principal fibre bundle $P$. 
Fermions, for example, 
are sections of a spinor bundle $F$, 
associated to a spin structure $P$.
As spinors transform according to a projective representation of the 
Lorentz group, they acquire a minus sign upon a full rotation.
This means that there is no hope of 
lifting global diffeomorphisms, so that spinor bundles cannot be 
natural.
One might feel however that
it would be unfair to discard them as unnatural.
After all, they do occur in nature, or at least in our 
mathematical description of it. 
As one can lift diffeomorphisms infinitesimally, i.e.\ 
at the level of vector fields, a spinor bundle does constitute
an infinitesimally natural bundle.

The third reason is that in field theory, 
one would like to construct a 
stress-energy-momentum tensor, 
corresponding to infinitesimal
transformations of space-time $M$.
Unfortunately, $\mathrm{Diff}(M)$
does not act on the fields $\Gamma(E)$
directly. Since $\mathrm{Aut}(P)$ does, we
need an infinitesimal lift of  
$\mathrm{Diff}(M)$ into $\mathrm{Aut}(P)$,
i.e.\ a Lie algebra homomorphism $\sigma$
that splits (\ref{rijtje}).
This means that an infinitesimally natural bundle is precisely 
what one needs in order to construct a SEM-tensor from Noether's 
principle. (See chapter \ref{ch:intro}, or e.g. \cite{GM,FR}.)
The splitting $\sigma$ is to be interpreted as the transformation 
behaviour of the fields under
infinitesimal space-time transformations.

The outline of the chapter is as follows.
Sections \ref{twee}, \ref{drie} and \ref{vier} are devoted to 
the classification of infinitesimally natural principal fibre bundles.
The central result is theorem \ref{grondsmakreefknol}, which 
states the following.
\begin{center}
\emph{Any infinitesimally natural principal fibre bundle
is associated to the universal cover of the connected component of the 
$k^{\mathrm{th}}$ order frame bundle $\tf$.}
\end{center}
We extend this to fibre bundles
with a finite dimensional structure group in 
section \ref{zes}, with special attention for vector
bundles.
In section \ref{vijf}, 
we seek conditions under which 
a splitting of (\ref{rijtje}) gives rise to a flat connection.

Finally, in section \ref{zeven}, we study
spin structures in the presence of a gauge field.
Ordinary spin structures are infinitesimally natural, 
as they carry a canonical 
splitting of (\ref{rijtje}).
In the presence of a gauge group $G$ however, one should consider 
$\spg$-structures instead (see \cite{HawPop}, \cite{AI}). 
For example, if $G = U(1)$, the appropriate bundles are 
$\mathrm{Spin}^{c}$-structures. 
In contrast to ordinary spin structures, not all 
$\spg$-structures are infinitesimally natural.
We show that for a compact gauge group, 
infinitesimally natural $\spg$-structures 
correspond precisely to
homomorphisms
$
\pi_{1}(F(M)) \rightarrow G
$ that are injective on $\pi_{1}(\gln)$.

In a sense, this shows that complicated manifolds
call for complicated gauge groups.
Some manifolds, such as $CP^{2}$, do not allow for any 
infinitesimally natural $\spg$-structure at all.  
If one feels that fermions should have a well defined 
transformation behaviour under infinitesimal space-time
transformations, then 
these manifolds are disqualified as models for space-time.

\section{Principal Bundles as Lie Algebra Extensions}\label{twee}

We seek to classify infinitesimally natural principal fibre bundles.
It is to this end that we study Lie algebra homomorphisms
$\sigma$ that split (\ref{rijtje}).
In this section, we will prove that $\sigma$ must be a differential
operator of finite order.

The first step, to be taken in section \ref{ideaal}, is to show 
that maximal ideals in $\Gamma_{c}(TM)$ correspond precisely 
to points in $M$. Using this, we will see in section \ref{localmap} that 
$\sigma$ must be a local map. 
Once again, we stress that locality of $\sigma$ will be a theorem, 
not an assumption.
We will then prove, in section \ref{diffop}, that $\sigma$ is in fact a
differential operator.

Before we proceed, let us have a closer look at
the exact sequence of Lie algebras 
(\ref{rijtje}), derived 
from a smooth principal $G$-bundle $\pi : P \rightarrow M$.
The last term, $\Gamma_{c}(TM)$, is the Lie algebra of smooth, compactly supported 
vector fields on $M$.
The middle term, $\Gamma_{c}(TP)^{G}$, is the Lie algebra of 
$G$-invariant vector fields $v$ 
on $P$ such that $\pi(\mathrm{Supp}(v))$ is compact.
The pushforward defines a
Lie algebra homomorphism $\pi_{*} : \Gamma_{c}(TP)^{G} \rightarrow
\Gamma_{c}(TM)$ because of $G$-invariance, and
its kernel $\Gamma_{c}(TP)_{v}^{G}$ is the ideal of vertical
vector fields. 
It can be identified with $\Gamma_{c}(\mathrm{ad}(P))$, 
the Lie algebra of compactly supported
sections of the adjoint bundle $\mathrm{ad}(P) := P \times_{G}\mathfrak{g}$,
where $\mathfrak{g}$ is the Lie algebra of $G$.


\subsection{Ideals of the Lie Algebra of Vector Fields}\label{ideaal}
The following lemma, due to Shanks 
and Pursell \cite{SP},
constitutes the linchpin of the proof. 
It identifies the maximal ideals of $\Gamma_{c}(TM)$
as points in $M$.
The proof is taken from
\cite{SP}, with perhaps some minor clarifications. 
\begin{lemma} \label{Ranzig}
Let $\Gamma_{c}(TM)$ be the Lie algebra
of smooth compactly supported vector fields on $M$. 
Then the maximal ideals of $\Gamma_{c}(TM)$
are labelled by points $q$ in $M$, 
$I_{q}$ being the ideal of vector fields in $\Gamma_{c}(TM)$ which are 
zero and flat at $q$.
That is,  
$I_{q} = \{ v \in \Gamma_{c}(TM) \,|\, v(q)=0 \,\, \mathrm{and}\,\, 
(\mathrm{ad}(w_{i_1})\ldots \mathrm{ad}(w_{i_n}) v)(q)=0
\, \,\forall \,\,w_{i_1} , \ldots w_{i_n} ,\in \Gamma_{c}(TM)\}\,.$
\end{lemma}\vspace{-2mm}
\proof  
Let $I \subset \Gamma_{c}(TM)$ be an ideal. Suppose that 
there exists a point $q \in M$ such that $v(q) = 0$ for all $v \in
I$. Then $I \subseteq I_{q}$.
Indeed, for all $w_{i_1} , \ldots , w_{i_n} \in \Gamma_{c}(TM)$,
one has $(\mathrm{ad}(w_{i_1})\ldots \mathrm{ad}(w_{i_n}) v)(q)=0$ 
because the l.h.s. is in $I$.

Now suppose that $I$ does not have such a point. 
We will prove that this implies $I = \Gamma_{c}(TM)$.
From this, the lemma will follow. If $I$ is a maximal ideal, 
then surely $I \neq \Gamma_{c}(TM)$, so we must have $I \subseteq I_{q}$.
But it is clear from its definition that $I_{q}$ is an ideal, so, 
$I$ being maximal, we must have $I = I_{q}$. 
Conversely, the ideals
$I_{q}$ are all maximal, because any enveloping ideal must either be
all of $\Gamma_{c}(TM)$, or else be  
contained some $I_{\tilde{q}}$, which cannot be unless $\tilde{q} = q$.  

We are therefore left to prove the following: suppose that
$I$ is an ideal such that 
for all $q \in M$, there exists a $v \in I$ with $v(q) \neq 0$.
Then $I = \Gamma_{c}(TM)$. 

We will shortly prove the following statement:
for all $q \in M$, there exists a neighbourhood
$U_q $ such that all $w \in \Gamma_{c}(TM)$
with $\mathrm{Supp}(w) \subset U_q$ can be written as $w = [v,u]$
with $v \in I$ and $u \in \Gamma_{c}(TM)$.

Such $w$ are therefore elements of $I$. 
Membership of $I$ then extends to arbitrary 
$w \in \Gamma_{c}(TM)$ by covering its support with a finite 
amount of neighbourhoods 
$U_{q_1} \ldots U_{q_{N}}$, for which we construct a partition of unity
$\sum_{s=1}^{N} f_s = 1$, $\supp (f_s) \subseteq U_{q_s}$.
We may then write $w = \sum_{s=1}^{N} f_s w = \sum_{s=1}^{N} [v_s , u_s] \in I$.
This concludes the proof modulo lemma \ref{techniek}, which is next in line. \qed


\noindent Although one can easily find $u \in \Gamma_{c}(TM)$ and $v \in I$ such that 
$w = [u,v]$ locally, 
it is not always
clear how to extend these to global vector fields while simultaneously
satisfying $w = [u,v]$. For instance, with $M$ the circle 
$S^1$ and $v = w = \del_{\theta}$, the solution $u = \theta \del_{\theta}$
does not globally exist. 
We need to do some work in order to define a proper cutoff procedure.

\begin{lemma}\label{techniek}
Let $I$ be a maximal ideal of $\Gamma_{c}(TM)$ containing a 
vector field which does not vanish at $q \in M$.
Then $q$ has a neighbourhood
$U_q $ such that all $w \in \Gamma_{c}(TM)$
with $\mathrm{Supp}(w) \subset U_q$ can be written  $w = [v,u]$
with $v \in I$ and $u \in \Gamma_{c}(TM)$.
\end{lemma}\vspace{-2mm}
\proof
Choose $v \in I$ with $v(q) \neq 0$. There exist local co-ordinates
$x_1, \ldots, x_n$ 
and an open neighbourhood $W$ of $q$ such that $v|_{W} = \del_1$.
Choose $W$ to be a block centered around $q$, and
nest two smaller blocks (in local co-ordinates) 
inside, so that
$
W \supset V \supset U \,.
$
   
We take 
$W = (-\varepsilon,\varepsilon)^{n}$,
$V = (- { \frac{2}{3}}\varepsilon,\frac{2}{3}\varepsilon)^{n}$, and
$U = (- { \frac{1}{3}}\varepsilon,\frac{1}{3}\varepsilon)^{n}$. 
Choose a smooth function $g$ on $M$ such that $g|_U = x_1$ 
and $g|_{M\backslash V}=0$. 
Require also that $\del_i g(x) = 0$ for $i \neq 1$ and
$(x_2, \ldots ,x_n) \in (- { \frac{1}{3}}\varepsilon,\frac{1}{3}\varepsilon)^{n-1}$.
Define $h := \del_1 g$, and set 
$\tilde{v} := [v,g\del_1] = h \del_1$.
Then $\tilde{v}|_U = \del_1$, $\tilde{v}|_{M\backslash V} = 0$, 
and most importantly
$\tilde{v} \in I$.

Now let $w \in \Gamma_{c}(TM)$ with $\supp(w) \subset U$.
We will find a global vector field $u \in \Gamma_{c}(TM)$ that realizes 
$[\tilde{v},u]=w$. 

For $i \neq 1$, the $i^{\mathrm{th}}$ component of the above reads
$h \del_1 u^i = w^i$. In the region $x_1 \leq \frac{1}{3}\varepsilon$, we set
$u^i (x_1, \ldots , x_n) = 
\int_{-\infty}^{x_1} w^i(t,x_2 \ldots,x_n) \mathrm{d}t$. 
On $U$, where 
$\tilde{v} = \del_1$, we then have $\del_1 u^i = w^i$
and therefore $h \del_1 u^i = w^i$. This is obviously also correct for  
points outside $U$ with $x_1 \leq \frac{1}{3}\varepsilon$, as both $u^i$
and $w^i$ are zero.

For $\frac{1}{3}\varepsilon \leq x_1 \leq \frac{2}{3}\varepsilon$,
where $w=0$, 
we let $u^i (x_1, x_2, \ldots , x_n)$ be a constant function of
$x_1$, so that 
$u^i (x_1, x_2, \ldots , x_n) 
= u^i (\frac{1}{3} \varepsilon, x_2, \ldots , x_n)$.
We then have $\del_1 u^i (x_1, \ldots , x_n) = 0$, 
guaranteeing $h \del_1 u^i = w^i$.
Note that this does not effect the smoothness of $u^i$.

Finally, for $\frac{2}{3}\varepsilon \leq x_1 \leq \varepsilon$,
let $u^i$ tend to zero, and let $u^i$ be zero for 
$x_1 \geq \varepsilon$.
This can be done in such a way that $u^i$ remains smooth. 
Since both $h$ and $w$ are zero, we have $h \del_1 u^i = w^i$ on all of $M$.

The case $i = 1$ is handled similarly, the only difference being that 
the first component of 
$[\tilde{v},u]=w$ is now $h \del_1 u^1 - \sum_{j}u^j \del_j h = w^1$,
the term $u^1 \del_1 h$ of which cannot be dispensed with.
We have arranged for
$u^i(x)$, with $i \neq 0$, to be zero if   
$(x_2, \ldots ,x_n) \notin (- {
\frac{1}{3}}\varepsilon,\frac{1}{3}\varepsilon)^{n-1}$,
so that $\del_i h$ equals zero if $u^i$ does not.
This leaves us with 
$h \del_1 u^1 - u^1 \del_1 h = w^1$.

For $x_1 \leq \frac{1}{3}\varepsilon$,
one once again sets $u^1 (x_1, \ldots , x_n) = 
\int_{-\infty}^{x_1} w^1(t,x_2 \ldots,x_n) \mathrm{d}t$.
For $\frac{1}{3}\varepsilon \leq x_1 \leq \frac{2}{3}\varepsilon$ however,
one now has to define $u^i (x_1, \ldots , x_n)$ 
$=$ $h(x_1, \ldots , x_n)$ $u^i (\frac{1}{3}\varepsilon, x_2, \ldots , x_n)$
in order for $h \del_1 u^1 - u^1 \del_1 h = w^1$ to hold.
Since this renders $u^1$ zero on the boundary of $V$, one is then free 
to define $u^i$ to be zero on $M\backslash V$.

Thus, if $\supp(w) \subset U$, we see that 
$w = [\tilde{v},u]$ with $\tilde{v} \in I$ and
$u \in \Gamma_{c}(TM)$. This concludes the proof 
lemma \ref{techniek}, and thereby that of lemma \ref{Ranzig}. \qed
A maximal subalgebra $A$ of a Lie algebra $L$ is either self-normalizing 
or ideal. 
Indeed it is contained in its normalizer,
which therefore equals either $A$ or $L$. 
A theorem of Barnes \cite{Ba} states that a finite-dimensional 
Lie algebra is nilpotent if and only
if\footnote{
Actually, the `only if' part in Barnes' theorem is not written down in
\cite{Ba}, but this is immediately clear from the proof of Engel's theorem.
(See e.g. \cite{Hu}).
} every maximal subalgebra is an ideal. 
On the other extreme: 
\begin{proposition}\label{perfectisideaal}
Let $L$ be a Lie algebra over a field $\mathbf{K}$, and
let 
$\cal{S}$ be the set of subspaces $A \subset L$ such that $A$ is
both an ideal and a maximal subalgebra.
Then
$$
[L,L] = \bigcap_{A \in \cal{S}} A \,.
$$ 
In particular, since the r.h.s reads `$L$' if $\cal{S} = \emptyset$,
$L$ is 
perfect $([L,L] = L)$
if and only if every maximal subalgebra is self-normalizing.
\end{proposition}\vspace{-2mm}
\proof
Let $X \notin [L,L]$. Then choose 
$[L,L] \subseteq A \subsetneq L$ where $A$ has codimension 1 in $L$,
and $X \notin A$. A is an ideal maximal subalgebra,  
which does not contain $X$.
Thus $\bigcap_{A \in \cal{S}} A \subseteq [L,L]$.

Let $A$ be an ideal maximal subalgebra, and $X \notin A$.
Then $A + \mathbf{K} X$ is a subalgebra strictly containing $A$,
so that it must equal $L$. Thus $[L,L] = [A + \mathbf{K} X,A + \mathbf{K} X] 
\subseteq A$,
whence $[L,L] \subseteq \bigcap_{A \in \cal{S}} A$. \qed
As a corollary, we have the following well known statement: 
\begin{corollary} \label{pstorm} The Lie algebra $\Gamma_{c}(TM)$ is perfect; 
$$[\Gamma_{c}(TM),\Gamma_{c}(TM)] = \Gamma_{c}(TM)\,.$$
\end{corollary}\vspace{-2mm}
\proof 
According to lemma \ref{Ranzig}, the maximal ideals are precisely the ideals 
$I_{q}$ of vector fields in $\Gamma_{c}(TM)$ which are zero and
flat at $q$. 
$I_{q}$ is strictly contained in the subalgebra 
$A_{q}$ of vector fields which are zero at $q$, so that
no ideal is a maximal subalgebra. So every maximal subalgebra
is self-normalizing.\vspace{-2mm} \qed

\subsection{The Splitting as a Local Map}\label{localmap}
With the main technical obstacles out of the way, we turn our attention to 
the sequence (\ref{rijtje}).
We now prove that $\sigma$ is a local map.
%
\begin{lemma}\label{zwabbernoot}
Let $P \rightarrow M$ be a principal $G$-bundle over $M$, with $G$
any Lie group.
Let $\sigma: \Gamma_{c}(TM) \rightarrow \Gamma_{c}(TP)^{G}$
be a Lie algebra homomorphism splitting the exact sequence of Lie-algebras
\begin{equation}\label{srijtje}
0 \rightarrow \Gamma_{c}(TP)^{G}_{v} \rightarrow \Gamma_{c}(TP)^{G}
\rightarrow \Gamma_{c}(TM) \rightarrow 0 \,.
\end{equation}
Then $\sigma$ is local in the sense that 
$\pi (\supp(\sigma(v))) \subseteq \supp(v)$.
\end{lemma}\vspace{-2mm}
\proof 
Any principal fibre bundle possesses an equivariant connection 1-form 
$\omega \in \Omega^{1}(P,\mathfrak{g})$, which
enables one to lift vector fields.
More precisely, the lifting map 
$\gamma :\Gamma_{c}(TM) \rightarrow \Gamma_{c}(TP)^{G}$ is defined 
by the requirement that 
$\omega_{p}(\gamma(v_{\pi(p)}))$ equal zero for all $p \in P$,
and that $\pi_{*} \circ \gamma$ be the identity. 
The lifting map $\gamma$ splits the exact sequence (\ref{srijtje}) as a sequence of  
$C^{\infty}(M)$-modules, but generally \emph{not} as a sequence 
of Lie algebras, since it need not be 
a homomorphism.


Define $f := \sigma - \gamma$. 
Then $\pi_{*} \circ f = 0$, 
so that $f$ is a map $\Gamma_{c}(TM) \rightarrow \Gamma_{c}(TP)^{G}_{v}$.
Since $\gamma$ is local by definition, 
the lemma will follow if we show
$f$ to be local.

The homomorphism property of $\sigma$, 
$[\sigma(v),\sigma(w)] - \sigma([v,w]) = 0$,
translates to 
\begin{equation}\label{vlijmenfileermes}
\begin{gathered}
f([v,w]) - [f(v) , f(w)] =\\
[\gamma(v) , f(w)] - [\gamma(w), f(v)] 
+ 
[\gamma(v),\gamma(w)] - \gamma([v,w])
.
\end{gathered}
\end{equation}
The $G$-action effects a Lie algebra isomorphism 
$\Gamma_{c}(TP)^{G}_{v} \simeq \Gamma_{c}(\mathrm{ad}(P))$. 
From that perspective, we have for $s \in \Gamma_{c}(\mathrm{ad}(P))$ that
$[\gamma(v) , s] = \nabla_{v} s$, the covariant derivative
along $v$. The curvature of the connection is then given by
$R(v,w) = [\gamma(v),\gamma(w)] - \gamma([v,w])$.
We rewrite equation (\ref{vlijmenfileermes}) as
\begin{equation}\label{snelfornuis}
f([v,w]) - [f(v) , f(w)] = 
\nabla_v f(w) - \nabla_{w} f(v) + R(v,w)  \,.
\end{equation}
Pick $m \in M$, and identify the fibre of $\mathrm{ad}(P)$ 
over $m$ with the structure Lie algebra $\mathfrak{g}$.
The restriction
$r_m : \Gamma_{c}(\mathrm{ad}(P)) \rightarrow \mathfrak{g}$
is a homomorphism of Lie algebras.
Now define $A_m := \{v \in \Gamma_{c}(TM) \,|\, v(m) = 0 \}$ to be the 
maximal subalgebra of vector fields that vanish at $m$,
and consider the map
$\tilde{f}_{m} : A_m \rightarrow \mathfrak{g}$
defined by $\tilde{f}_{m} = r_m \circ f |_{A_m}$.
Equation (\ref{snelfornuis}) transforms into 
\begin{equation}\label{cofflofot}
\tilde{f}_m([v,w]) - [\tilde{f}_m(v) , \tilde{f}_m(w)] = 
\nabla_v (\tilde{f}(w))_m - \nabla_{w} (\tilde{f}(v))_m + R(v,w)_m  \,,
\end{equation}
the r.h.s. of which vanishes
because $\nabla$ and $R$ are 
$C^{\infty}(M)$-linear in $v$ and $w$, 
and $v(m) = w(m) = 0$.
This means that $\tilde{f}_m$ is a Lie algebra homomorphism. 

Let us restrict $\tilde{f}_m$ even further to the ideal 
$\Gamma_{c}(T(M-\{m\})) \subset A_m$ (which is non-maximal), and note that
$\hat{f}_m : \Gamma_{c}(T(M-\{m\})) \rightarrow \mathfrak{g}$
is a homomorphism.
%
Its kernel $\mathrm{Ker}(\hat{f}_m)$ is therefore an ideal
in $\Gamma_{c}(T(M-\{m\}))$, and one of finite codimension at that. 
Indeed,
$$\Gamma_{c}(T(M-\{m\})) / \mathrm{Ker}(\hat{f}_m) 
\simeq \mathrm{Im}(\hat{f}_m) \subseteq \mathfrak{g}\,.
$$ 

According to lemma \ref{Ranzig} however, all proper ideals
are of infinite codimension, forcing 
$\mathrm{Ker}(\hat{f}_m) = \Gamma_{c}(T(M-\{m\}))$. 
But the vanishing of $\hat{f}_m$ for all $m \in M$ is tantamount to
locality of $f$ in the sense that $\supp(f(v)) \subseteq \supp(v)$. \qed
We gather some notation that we introduced in the course 
of the proof. If $\sigma : \Gamma_{c}(TM) \rightarrow \Gamma_{c}(TP)^{G}$
is a Lie algebra homomorphism splitting $\pi_{*}$ and 
$\gamma : \Gamma_{c}(TM) \rightarrow \Gamma_{c}(TP)^{G}$ is the lift induced by
an equivariant connection, then we define $f := \sigma - \gamma$. 
It maps to the vertical vector fields, identified with $\Gamma(\mathrm{ad}(P))$.
We define $A_{m}$ to be the maximal subalgebra of vector fields that 
vanish at
$m$, and $\tilde{f}_{m} : A_{m} \rightarrow \mathfrak{g}$ to be the restriction
of $f$ to $A_{m}$, followed by the map 
$\Gamma(\mathrm{ad}(P)) \rightarrow \mathfrak{g}$ which picks out the fibre over 
$m$ and identifies it with $\mathfrak{g}$. 

\subsection{The Splitting as a Differential Operator}\label{diffop}

In this section, we will prove that $\sigma$ is a differential operator
of finite order. 
Since $f$ is local, it defines a map from the sheaf of smooth sections 
of $TM$ to the sheaf of smooth sections of $\mathrm{ad}(P)$.
An elegant theorem of Peetre (\cite{Pe}) then says that $f$, and therefore 
$\sigma$, must be
a differential operator of locally finite order. 
All we need to do is find a global bound on the order, which will occupy
us for 
the remainder of the section.

The fact that $f$ is a differential operator of locally finite order 
means that for each $m \in M$, there exists an $r$ such that 
$H^{r}_m = \{v \in A_m \,|\, j^{r}_{m}(v) = 0\}$ is contained in 
$\ker(\tilde{f}_m)$, where $j^{r}_{m}(v)$ is the $r$-jet of $v$ at $m$.
Consequently, $\tilde{f}_{m}$ factors through the jet Lie algebra
$J^{r,0}_m(TM) := A_{m}/H^{r}_{m}$. Since we do not know $r$, we define\vspace{-0.5mm} 
$$\vecn := \lim_{\longleftarrow} A_m / H^{r}_{m}\,,\vspace{-0.5mm}$$
and remark that $\tilde{f}_m$ induces a homomorphism 
$\vecn \rightarrow \mathfrak{g}$. 
The Lie algebra $\vecn$ depends on $M$
only through its dimension $n$. 


Local co-ordinates provide one with a basis $x^{\vec{\alpha}}\del_i$,
where $x^{\vec{\alpha}}$ is shorthand for
$x_1^{\alpha_{1}} \ldots x_{n}^{\alpha_{n}}$.  
With
$\vecn^{k} = \mathrm{Span}\{ x^{\vec{\alpha}}\del_i \,|\, 
|\vec{\alpha}| = k +1 \, ,\, i = 1 \ldots n\} $, 
one may write\vspace{-1mm}
$$
\vecn = \bigoplus_{k=0}^{\infty}\vecn^{k}\,.\vspace{-1mm}
$$
That is, each element of $\vecn$ can be uniquely written as a 
\emph{finite} sum of homogeneous vector fields.
Note that $\vecn^{k}$ is the $k$-eigenspace of 
the Euler vector field $E := \sum_{i=1}^{n} x^i \del_i$. 
If $I$ is an ideal containing $v = \sum_{k=0}^{N} v_{k}$, 
then $\mathrm{ad}(E)^{j}v  = \sum_{k=0}^{N} k^{j} v_{k} \in I$ for all 
$j$, so that $v_{k} \in I$.
Thus any ideal splits into homogeneous components\vspace{-1mm}
$$
I  = \bigoplus_{k=0}^{\infty} I^{k}\vspace{-1mm}
$$ 
with $I^{k} = I \cap \vecn^{k}$. This renders the ideal structure 
of $\vecn$ more or less tractable, so that we may prove the following 
bound on the order of $\tilde{f}_m$.

\begin{lemma}\label{hallekal}
The order of the differential operator $\sigma$ is at most 
$\mathrm{dim}(\mathfrak{g})$ unless $\mathrm{dim}(M) = 1$
and $\mathrm{dim}(\mathfrak{g}) = 2$, in which case the order is 
at most 3.\vspace{-2mm}
\end{lemma}
\proof We closely follow Epstein and Thurston, \cite{ET}.
One checks by hand that the
only ideals of
$\mathrm{Vec}_{1} = \mathrm{Span}\{x^{k} \del \,|\, k \geq 1\}$
are 
$\mathrm{Span}\{x^2\del , x^k \del \,| \, k \geq 4 \}$
and $\mathrm{Span}\{x^k \del \,| \, k \geq N \}$ with $N \geq 1$.

Consider
$\mathrm{Vec}_{1}$ as a subalgebra of $\vecn$,
define $K$ to be the kernel of $\tilde{f}_m$, and let 
$K_1 := K \cap \mathrm{Vec}_{1}$. 
We then have injective homomorphisms
$$
\mathrm{Vec}_{1}/K_1 \hookrightarrow \vecn/K \hookrightarrow \mathfrak{g}\,,
$$
so that $\mathrm{dim}(\mathrm{Vec}_{1}/K_1) \leq \mathrm{dim}(\mathfrak{g})$.
As $K_1$ is an ideal, it must be of the shape mentioned above.
This leads us to conclude that $x_1^{k} \del_1 \in K$ 
for all $k > \mathrm{dim}(\mathfrak{g})$ unless 
$\mathrm{dim}(\mathfrak{g}) = 2$, in which case 
$x_1^{k} \del_1 \in K$ for all $k > 3$, and $x_1^2\del_1 \in K$.

The following short calculation shows that if $\mathrm{dim}(\mathfrak{g}) = 2$ and
$\mathrm{dim}(M) > 1$, then also $x_1^3 \del_1 \in K$.
As $K$ contains $x_1^2\del_1$, it also contains 
$[x_1^2\del_1 , x_1 \del_2] = x_1^2 \del_2$,
and thus $[x_1^2 \del_2 , x_1 x_2\del_1] = x_1^3 \del_1 - 2 x_1^2 x_2 \del_2$.
But by bracketing with $x_1^2 \del_2$ and $x_2 \del_1$ respectively,
we see that $x_1^3 \del_1 - 3 x_1^2 x_2 \del_2$ is in $K$,
ergo $x_1^3 \del_1 \in K$. 

The next step is to show that if $x_1^{s}\del_1 \in K$, then 
$K$ also contains all $x^{\vec{\alpha}} \del_i$ with $|\vec{\alpha}| = s$.
First of all, we remain in $K$ if we repeatedly apply 
$\mathrm{ad}(x_i \del_1)$ to $x_1^{s}\del_1$, to the effect of 
replacing $x_1$ by $x_i$ up to a nonzero factor.
This shows that $x^{\vec{\alpha}}\del_1 \in K$.
Then the relation 
$
x^{\vec{\alpha}}\del_i =
[x^{\vec{\alpha}}\del_1 , x_1 \del_i]  + x_1 \del_i x^{\vec{\alpha}} \del_1 
$
transfers membership of $K$ from right to left.

In the generic case $\mathrm{dim}(\mathfrak{g}) \neq 2$, 
$\mathrm{dim}(M) \neq 1$, we may conclude that  
the order of $\sigma$
is at most $\mathrm{dim}(\mathfrak{g})$, because
$H_m^{\mathrm{dim}(\mathfrak{g})} \!\subset K$.
In the exceptional case that $\mathrm{dim}(\mathfrak{g}) = 2$ and
$\mathrm{dim}(M) = 1$, the order of $\sigma$ is at most 3.\qed
In particular, $\sigma$ is a differential operator of finite rather
than locally finite order. 
Let us summarize our progress so far.
\begin{proposition} \label{alhetvet}
Let $P$ be an infinitesimally natural principal $G$-bundle.
Then $\sigma : \Gamma_{c}(TM) \rightarrow \Gamma_{c}(TP)^{G}$ 
factors through the bundle of $k$-jets, where 
$k=3$ if $\mathrm{dim}(M) = 1$, $\dg = 2$ and
$k = \dg$ otherwise.  
If we identify 
equivariant vector fields on $P$ with
sections of 
the Atiyah bundle $TP/G$, we can therefore
define a bundle map 
$\nabla : J^{k}(TM)\rightarrow TP/G$ by
$\nabla (j_{m}^{k}(v)) := \sigma(v)_m$.
It makes the following diagram commute.
\begin{center}
\begin{tikzpicture}
\pgfsetzvec{\pgfpoint{0.385cm}{-0.385cm}}
\node (links) at (0, 0, 0) {$\Gamma(TM)$};
\node (rechts) at (2.5, 0, 0) {$\Gamma(TP/G)$};
\node (onder) at (1.3, -1.3, 0) {$\Gamma(J^{k}(TM))$};
\node (nab) at (2.2,-0.8,0) {$\nabla$};
\draw [->] (links) -- node[above] {$\sigma$}(rechts);
\draw [->] (links) -- (onder);
\draw [->] (onder) -- (rechts);
\end{tikzpicture}
\end{center}
\end{proposition}
The point is that although $\sigma$ is defined only on sections, 
$\nabla$ comes from a veritable bundle map 
$J^{k}(TM) \rightarrow TP/G$.  
Note that although $\sigma$ was only defined on $\Gamma_{c}(TM)$, 
it extends to $\Gamma(TM)$ by locality. 

\section{Lie Groupoids and Algebroids of Jets}\label{drie}
The bundles $J^{k}(TM)$ and $TP/G$ are Lie algebroids, and it will be 
essential for us to prove that $\nabla : J^{k}(TM) \rightarrow TP/G$  
is a homomorphism of Lie algebroids. In order to do this, we will
first have a closer look at $J^{k}(TM)$ and $TP/G$, and at their 
corresponding Lie groupoids.

Let us first set some notation.
The jet group $G^{k}_{0,0}( \mathbb{R}^{n})$ is the group of 
$k$-jets of diffeomorphisms 
of $\mathbb{R}^n$ that fix $0$. It is the semi-direct 
product of $\gln$
and the connected, simply connected, unipotent Lie group of $k$-jets
that equal the identity to first order.

The subgroup $\smash{G^{+ k}_{0,0}(\mathbb{R}^{n})}$ of orientation preserving 
$k$-jets is  
connected, but not simply connected. 
As $\smash{G^{+ k}_{0,0}(\mathbb{R}^{n})}$ retracts to $\mathrm{SO}(\mathbb{R}^n)$, 
its homotopy group is isomorphic to 
$\{1\}$ if $n=1$, to
$\mathbb{Z}$ if $n=2$, and to
$\mathbb{Z} / 2\mathbb{Z}$ if
$n>2$. For brevity, we introduce the following notation.
\begin{basdef}
If $k > 0$, we denote $\pi_{1}(G^{+ k}_{0,0}(\mathbb{R}^{n}))$ by $Z$.
\end{basdef}
Thus for $n > 2$, the universal cover 
$\tilde{G}^{+ k}_{0,0}(\mathbb{R}^{n}) \rightarrow 
G^{+ k}_{0,0}(\mathbb{R}^{n})$ is $2 : 1$,
and restricts to the spin group over $\mathrm{SO}(\mathbb{R}^n)$.

%

\subsection{The Lie Groupoid of \texorpdfstring{$k$}{k}-Jets}

In this section, we define the Lie groupoid $G^{k}(M)$ of $k$-jets, 
its maximal source-connected Lie subgroupoid $G^{+k}(M)$,
and the $k^{\mathrm{th}}$ order 
frame bundle $F^{k}(M)$.
   
Denote by $G^{k}_{m',m}(M)$ the manifold of $k$-jets at $m$ of diffeomorphisms
of $M$ which map $m$ to $m'$, and denote 
by $\smash{G^{k}(M) = \cup_{M \times M} G^{k}_{m',m}(M)}$ the groupoid
of $k$-jets. If $\smash{j^{k}_{m}}(\alpha)$ is a $k$-jet at $m$ of a 
diffeomorphism $\alpha$, 
then its source is 
$\smash{s(j^{k}_m(\alpha)) = m}$, its target is $\smash{t(j^{k}_m(\alpha))} = \alpha(m)$,
and multiplication is given by concatenation.
  
Denote by $G^{k}_{*,m}(M)$ the manifold 
$s^{-1}(m)$ of $k$-jets with source $m$. The target map 
$t: \smash{G^{k}_{*,m}(M)} \rightarrow M$ endows it with a structure of 
principal fibre bundle, the structure
group $G^{k}_{m,m}(M) \simeq G^{k}_{0,0}(\mathbb{R}^{n})$ acting 
freely and transitively on the right. 
As $G^{1}(M)_{*,m}$ is isomorphic to the frame bundle $F(M)$, 
one calls $G^{k}(M)_{*,m}$ the $k^{\mathrm{th}}$ order frame bundle, 
sometimes denoted $F^{k}(M)$.

\begin{lemma}
Let $M$ be connected and let $k \geq 1$. Then $G^{k}(M)$ is
source-connected if and only if $M$ is not orientable.\vspace{-2mm}
\end{lemma}
\proof
We may as well consider $k=1$, because the fibres of 
$G^{k}(M) \rightarrow G^{1}(M)$
can be contracted.
Each source fibre $G^{1}(M)_{*,m}$ of $G^{1}(M)$ is isomorphic to
the frame bundle. By definition, $M$ is oriented precisely when
the frames can be grouped 
into positively and negatively oriented ones.\vspace{-1mm}\qed

\begin{basdef} We define $G^{+k}(M)$ to be the maximal source-connected
Lie subgroupoid of $G^{k}(M)$, and denote its source fibre by 
$F^{+k}(M)$.\vspace{-1mm}
\end{basdef}

\noindent In the light of the previous lemma, this means that 
$G^{+ k}(M)$ is
the Lie groupoid of $k$-jets of
orientation preserving diffeomorphisms
if $M$ is orientable, and simply $G^{k}(M)$ if $M$ is not.

Note that the map $D : \mathrm{Diff}(M) \rightarrow \mathrm{Diff}(G^{k}(M))$
defined by
$D\alpha : j_m^{k}(\gamma) \mapsto j_m^{k}(\alpha \gamma)$
is a homomorphism of groups.
We will call it the $k^{\mathrm{th}}$ order derivative.
Because $D\alpha$ is
source-preserving 
and right invariant, 
it defines a homomorphism 
$\mathrm{Diff}(M) \rightarrow 
\mathrm{Aut}^{G^{k}_{m,m}(M)}(G^{k}_{*,m}(M))$, 
splitting the exact sequence of groups (\ref{flapperdrop}).
This makes $G_{*,m}^{k}(M) = F^{k}(M)$ into
a natural bundle.
Note that $F^{+k}(M)$ is infinitesimally natural.




\subsection{The Lie Algebroid of \texorpdfstring{$k$}{k}-Jets}
The bundle $J^{k}(TM)$ possesses a structure of Lie algebroid, 
induced by the
Lie groupoid $G^{k}(M)$. 
We now describe
the Lie bracket on $\Gamma(J^{k}(TM))$ explicitly.
Later, in section \ref{loidhom}, we will use this to  
show that $\nabla$
is a Lie algebroid homomorphism. 

The Lie algebroid of $G^{k}(M)$ is a vector bundle
$A \rightarrow M$.
Its fibre 
$A_m$ is by definition the subspace of
the tangent space of $G^{k}(M)$ at $j_{m}^{k}(\id)$ 
which is annihilated by $ds$.
Sections of $A$ therefore correspond 
to right-invariant vector fields on $G^{k}(M)$ parallel to the source fibres.

Each curve in $G^{k}_{*,m}(M)$ through $j^{k}_{m}(\id)$ takes the shape
$c(t) = j^{k}_{m}(\alpha_{t})$ with $\alpha_{0} = \id$,
so that its tangent vector $a\in A_{m}$ takes the shape $a = j_{m}^{k}(v)$, 
with 
$v = {\del_{t}}|_{0} \alpha_{t}$. 
This shows that $A \simeq J^{k}(TM)$. 

The anchor $d t : J^{k}(TM) \rightarrow TM$ is easily seen to be the canonical
projection, so we shall denote it by $\pi$. The Lie bracket on 
$\Gamma(J^{k}(TM))$ however, 
which is defined as the restriction of 
the commutator bracket on $\Gamma(TG^{k}(M))$ to the right invariant 
source preserving vector fields, 
perhaps deserves 
some comment.


Define $J^{k,0}(TM)$ to be the kernel of $\pi$, and consider
the exact sequence of Lie algebras\vspace{-1mm}
$$
0 \rightarrow 
\Gamma(J^{k,0}(TM)) 
\rightarrow \Gamma(J^{k}(TM)) 
\stackrel{\pi}{\rightarrow} \Gamma(TM) 
\rightarrow 0\,.
$$ 
It is split by 
$j^{k} : 
\Gamma(TM) \rightarrow 
\Gamma(J^{k}(TM))
$, the infinitesimal version of the $k^{\mathrm{th}}$ order derivative. 
(The sequence of Lie algebroids of course
does not split, as differentiation is not linear over 
$C^{\infty}(M)$.)

We will describe the Lie bracket on 
$\Gamma(J^{k}(TM))$ by giving it on $\Gamma(TM)$
and $\Gamma(J^{k,0}(TM))$ separately, and then giving the action of 
$\Gamma(TM)$
on $\Gamma(J^{k,0}(TM))$.

\begin{proposition}
Let $u$ and $u'$ be sections of $TM$, and 
let
$\tau : m \mapsto j_m^{k}(v_m)$ and 
$\tau' : m \mapsto j_m^{k}(v'_m)$ be sections of $J^{k,0}(TM)$. 
Then\vspace{-1mm} 
\begin{eqnarray*}
{}[j^{k}(u) , j^{k}(u')]_{m} &=& j_{m}^{k}([u,u'])\,,\\
{}[\tau , \tau']_{m} &=& j^{k}_{m}([v_{m},v_{m}'])\,,\\
{}[j^{k}(u) , \tau]_{m} &=& j^{k}_{m}([u, v_{m}]) + 
j_{m}^{k}( d_{u}|_{m}(x \mapsto v_x) )\,,\vspace{-1mm}
\end{eqnarray*}
where 
$d_{u}|_{m}(x \mapsto v_x)$ 
is the ordinary derivative at $m$ along $u$ of a map from
$M$ to $\Gamma(TM)$. 
Although both terms on the right hand side depend on the choice of 
$m \mapsto v_{m}$, their sum does not.
\end{proposition}\vspace{-2mm}
\proof
The first equality is clear, as $j^{k}$ is a homomorphism of Lie algebras.
The second equality can be seen as follows. 
Consider the bundle of groups
$G^{k}(M)_{*,*} := \{j^{k}_{m}(\alpha) \in G^{k}(M)\, | \,\alpha(m) = m \}$,
with bundle map $s = t$. 
Its sections $\Gamma(G^{k}(M)_{*,*})$ form a group under 
pointwise multiplication,
the Lie algebra of which is $\Gamma(J^{k,0}(TM))$,
with the pointwise bracket.
As $\Gamma(G^{k}(M)_{*,*})$ acts from the left on 
$G^{k}(M)$ by 
$j^{k}_m(\gamma) \mapsto j_{m}^{k}(\alpha_m) \circ j_{m}^{k}(\gamma)$, 
respecting
both the source map and right multiplication, the inclusion
$\Gamma(J^{k,0}(TM)) \rightarrow \Gamma(J^{k}(TM))$ is a 
homomorphism of Lie algebras.
This proves the second line.

To verify the third line, we must choose a smooth
map $x \mapsto v_{x}$ from $M$ to $\Gamma(TM)$ such that
$\tau_{x} = j_{x}^{k}(v_{x})$ in a neighbourhood of $m$.
Each $v_{x}$ necessarily has a zero at $x$.  
If we denote $m(s) := \exp(su)m$, then the bracket 
$[j^{k}(u) , \tau]_{m}$ is by 
definition\footnote{The 
Lie algebra of the diffeomorphism group is
the Lie algebra of vector fields,
but the exponential map is given by
$v \mapsto \exp(-v)$, where $\exp$ denotes the
unit flow along $v$.
This is why the groupoid commutator might seem odd at first sight.} 
minus the mixed second derivative along $s$ and $t$ at $0$ of
the groupoid commutator
$$
j^{k}_{m(s)}(\exp(-su)) 
j^{k}_{m(s)}(\exp(-tv_{m(s)})) 
j^{k}_{m}(\exp(su)) 
j_{m}^{k}(\exp(t v_{m}))
\,,$$
which is just 
$
j^{k}_{m} \left(
\exp(-su) 
\exp(-tv_{m(s)}) 
\exp(su) 
\exp(t v_{m})
\right)\,.
$
The terms not involving derivatives of $s \mapsto m(s)$ yield
$j_{m}^{k}([u,v_{m}])$, and the terms which do provide the extra
$j^{k}_{m}(d_{u}|_{m}(x \mapsto v_{x}))$. \qed

\subsection{The Gauge Groupoid and its Algebroid}
Given a principal $G$-bundle $\pi: P \rightarrow M$, one can define the
gauge groupoid $(P \times P)/G$, that is the pair groupoid modded out by the
diagonal action. 
Source and target 
come from projection on the second and first term respectively, 
$M \hookrightarrow (P \times P)/G$ as $\id_{\pi(p)} = [(p,p)]$, and 
multiplication is well defined by $[(r,q)] \circ [(q,p)] = [(r,p)]$.
An element $[(q,p)]$ corresponds precisely to a $G$-equivariant 
diffeomorphism $\pi^{-1}(p) \rightarrow \pi^{-1}(q)$, and the product to
concatenation of maps.

Its Lie algebroid $TP/G$ is sometimes called the Atiyah algebroid.
Indeed, 
the subspace of $T_{\id_{m}}((P \times P)/G)$ which annihilates $ds$
is canonically $(TP/G)_{m}$.
A section of $TP/G$ can be identified with a $G$-equivariant
section of $TP$, endowing $\Gamma(TP/G)$ with the Lie bracket
that comes from $\Gamma(TP)^{G}$.

\subsection{The Splitting as a Homomorphism of Lie Algebroids}\label{loidhom}

The point of considering the Lie algebroid structure of $J^{k}(TM)$
was of course to prove the following.
\begin{lemma}
The map $\nabla : J^{k}(TM) \rightarrow TP/G$ is a homomorphism of
Lie algebroids.
\end{lemma}
 
\proof The fact that $\nabla$ respects the anchor is immediate.
As $\pi_{*} \circ \nabla (j_{m}^{k}(v_m)) = 
\pi_{*} \circ \sigma(v_{m})(m) = v_{m}(m)$, it equals 
$\pi \circ j^{k}(v_{m})$ in the point $m$.

We now show that $\nabla : \Gamma(J^{k}(TM)) \rightarrow \Gamma(TP/G)$
is a homomorphism of Lie algebras. First of all, the restriction 
of $\nabla$ to $j^{k}(\Gamma(TM))$ is a homomorphism. Indeed,
$[\nabla(j^{k}(v)) , \nabla(j^{k}(w))] = [\sigma(v), \sigma(w)]$,
which equals $\sigma([v,w])$ because $\sigma$ is a homomorphism.
This in turn is just $\nabla(j^{k}([v,w]))$, so that
$[\nabla(j^{k}(v)) , \nabla(j^{k}(w))] = \nabla([j^{k}(v),j^{k}(w)] )$.

Secondly, its restriction to $\Gamma(J^{k,0}(TM))$ is a homomorphism.
If $\tau_x = j^{k}_{x}(v_x)$ and $\upsilon_x = j^{k}_{x}(w_x)$ 
are sections of $J^{k,0}(TM)$, then 
$\nabla \tau$ and $\nabla \upsilon$ are in the kernel of the anchor.
This implies that
their commutator at a certain point $m$ depends only on their
values at $m$, not on their derivatives.
To find the commutator at $m$, we may therefore
replace $j^{k}_{x}(v_x)$ by $j^{k}_{x}(v_{m})$ and likewise 
$j^{k}_{x}(w_x)$ by $j^{k}_{x}(w_{m})$. We then see that
$[\nabla j^{k}_{x}(v_x) , \nabla j^{k}_{x}(w_x)]_{m}$
$=$ $[\nabla j^{k}_{x}(v_{m}) , \nabla j^{k}_{x}(w_{m})]_{m}$.
We already know that this is
$\nabla( j^{k}_{x}([v_{m},w_{m}]))_{m}$,
so that
$[\nabla(\tau),\nabla(\upsilon)]_{m} = \nabla([\tau,\upsilon])_{m}$. 

The last step is to show that $\nabla$ respects the bracket between 
$j^{k}(\Gamma(TM))$ and $\Gamma(J^{k,0}(TM))$.
Again, let
$j^{k}(v)$ be an element of the former and $\tau_x = j^{k}_{x}(w_x)$ 
of the latter.
Considered as an equivariant vector field on $P$, the vertical 
vector field $\nabla(\tau)$ takes the value $\sigma_{p}(w_{\pi(p)})$
at $p \in P$.
Then $[\nabla(j^{k}(v)) , \nabla(\tau)]$ is the Lie derivative along 
$\sigma(v)$ of the vertical vector field $\sigma_{p}(w_{\pi(p)})$.
Differentiating along $\sigma(v)$ is done by considering
$(p,p') \mapsto \sigma_{p}(w_{\pi(p')})$, differentiating w.r.t. 
$p$ and $p'$ separately, and then putting $p=p'$.
This results in
$[\sigma(v), \sigma_{p}(w_{\pi(p)})]_{p_0} = 
[\sigma(v), \sigma(w_{\pi(p_0)})]_{p_0} + 
\sigma(d_{v}|_{\pi(p_0)} (x \mapsto w_{x}))$, 
which is in turn the same as 
$\sigma_{p_0}([v,w_{\pi(p_0)}] + d_{v}|_{\pi(p_0)} (m \mapsto w_{m})))$,
so that $[\nabla(j^{k}(v)) , \nabla(\tau)] = \nabla([j^{k}(v),\tau])$
as required. 
Therefore $\nabla$ must be a homomorphism on all of $\Gamma(J^{k}(TM))$. \qed

\begin{basdef}
A connection $\nabla$ of a Lie algebroid $A$ on a vector bundle $E$
is by definition a bundle map of $A$ into $\mathrm{DO}^{1}(E)$,
the first order differential operators on $E$, which respects the anchor.
If moreover it is a morphism of Lie algebroids, then the connection is 
called flat. A flat connection of $A$ on $E$ is also called a representation
of $A$ on $E$.
\end{basdef}

\noindent This explains our notation for the map $\nabla$ induced by $\sigma$.
Given a representation $V$ of $G$, one may form the 
associated vector bundle $E := P \times_{G} V$. The map $\nabla$
then defines a Lie algebroid homomorphism of $\Gamma(J^{k}(TM))$ 
into the Lie algebroid of first order differential operators on $E$.
(Simply consider a section of $E$ as a $G$-equivariant function
$P \rightarrow V$, and let $\Gamma^{G}(TP)$ act by Lie derivative.)
By definition, this is a flat connection, or equivalently a 
Lie algebroid representation.

\section{The Classification Theorem}\label{vier}
We use the fact that $\nabla : \Gamma(J^{k}(TM)) \rightarrow \Gamma(TP/G)$ 
is a homomorphism of Lie algebroids
to find a corresponding homomorphism of Lie groupoids.
This will give us the desired classification of 
infinitesimally natural principal fibre
bundles.

\subsection{Integrating a Homomorphism of Lie Algebroids}
The following theorem states that homomorphisms of Lie algebroids
induce homomorphisms of Lie groupoids if the initial groupoid is 
source-simply connected.\vspace{-3mm}
\begin{theorem}[Lie II for algebroids]
Let $G$ and $H$ be Lie groupoids, with corresponding Lie algebroids 
$A$ and $B$ respectively. Let $\nabla : A \rightarrow B$
be a homomorphism of Lie algebroids. If $G$ is source-simply 
connected, then there exists a unique homomorphism
$G \rightarrow H$ of Lie groupoids which integrates $\nabla$.
\end{theorem}
{{\bf \noindent Remark} \quad
The result was probably announced first in \cite{Pr}, and proofs have 
appeared e.g. in \cite{MX} and \cite{MM}.
We follow \cite{CF}, which the reader may consult for details.}

\sketchproof
The idea is that $\nabla$ allows one to lift 
piecewise smooth paths 
of constant source in $G$ to piecewise smooth paths of constant source 
in $H$. Source-preserving piecewise smooth homotopies in $G$ of course 
do not affect the endpoint of the path in $H$, so that,
if $G$ is source-simply connected, 
one obtains a map $G \rightarrow H$ by identifying 
elements $g$ of $G$ with equivalence classes of source preserving
paths from $\id_{s(g)}$ to $g$. 
One checks that this is 
the unique homomorphism of Lie groupoids integrating $\nabla$. 
\qed

\noindent Unfortunately, $G^{k}(M)$ is not always source connected, 
let alone source-simply connected.
Recall that $G^{+k}(M)$ is the maximal source-connected Lie subgroupoid of 
$G^{k}(M)$, and therefore has the same Lie algebroid $J^{k}(TM)$. 

We define $\tilde{G}^{+k}(M)$ to be the 
set of piecewise smooth, source preserving paths in $G^{+k}(M)$
beginning at an identity, 
modulo piecewise smooth, source preserving homotopies.
It is a smooth manifold because $G^{+k}(M)$ is, and a Lie groupoid
under the unique structure making the projection on the endpoint
$\tilde{G}^{+k}(M) \rightarrow G^{+k}(M)$ into a morphism of 
groupoids.
Explicitly, the multiplication is given as follows. 
If $g(t)$ a path from $\id_{m}$ to $g(1)_{m'm}$,
and $h(t)$ a path from $\id_{m'}$ to $h(1)_{m''m'}$,
then the product $[h] \circ [g]$ is 
$[(h \cdot g(1)) * g]$, where the dot denotes groupoid multiplication 
and the star concatenation of paths. The proof of associativity is the usual
one. 

Note that the source fibre $\tilde{G}^{+k}(M)_{*,m}$ is precisely
the universal cover of the connected component of
the $k^{\mathrm{th}}$ order frame bundle 
$G^{+k}(M)_{*,m} = F^{+k}(M)$.
In order to cut down on the subscripts, we introduce new notation 
for $\tilde{G}^{+k}(M)_{*,m}$ and its structure group 
$\tilde{G}^{+k}(M)_{m,m}$.
\begin{basdef} \label{zwalksmurf}
We denote the universal cover of the connected component of the 
$k^{\mathrm{th}}$ order frame bundle by $\tf$, and
its structure group by $G(k,M)$.  
\end{basdef}
It is an infinitesimally natural bundle because 
$F^{+k}(M)$ is. 
Note that $G(k,M)$ is not the universal cover of 
$G_{m,m}^{+k}(M)$,
but rather its extension by $\pi_{1}(F^{k}(M))$. 
As $\pi_{1}(F^{k}(M)) = \pi_{1}(F(M))$, we have the exact 
sequence of groups
$$
1 \rightarrow \pi_{1}(F(M)) \rightarrow G(k,M) \rightarrow 
G_{m,m}^{+k}(M) \rightarrow 1\,.
$$
The group $G_{m,m}^{+k}(M)$ in turn is isomorphic to 
$G_{0,0}^{+k}(\mathbb{R}^{n})$ if $M$ is orientable, and to
$G_{0,0}^{k}(\mathbb{R}^{n})$ if it is not.

\subsection{Classification}\label{Classification}

Now that we've found a source-simply connected Lie groupoid with 
$J^{k}(TM)$ as Lie algebroid, we can finally apply 
Lie's second theorem for algebroids to obtain the following.
\begin{proposition}\label{knolletje}
If $\sigma$ splits the exact sequence of Lie algebras (\ref{rijtje}),
then it induces a morphism of groupoids 
$\exp\nabla : \tilde{G}^{+k} \rightarrow (P \times P)/G$ such that the following diagram
commutes, with $\exp_{m}$ the flow along a vector field starting at 
$\id_{m}$.
\end{proposition}
\begin{center}
\begin{tikzpicture}
\pgfsetzvec{\pgfpoint{0.385cm}{-0.385cm}}
\node (links) at (0, 0) {$\Gamma(TM)$};
\node (midden) at (3, 0) {$\Gamma(TP/G)$};
\node (onderlinks) at (1.5, -1.5) {$\Gamma(J^{k}(TM))$};
\node (rechts) at (6,0) {$(P \times P )/G$};
\node (onderrechts) at (4.5, -1.5) {$\tilde{G}^{+k}(M)$};
\node at (0.3,-0.8) {$j^{k}$};
\node at (4.7,-0.7) {$\exp\nabla$};
\node at (1.95,-0.7) {$\nabla$};
\draw [->] (links) -- node[above] {$\sigma$}(midden);
\draw [->] (links) -- (onderlinks);
\draw [->] (onderlinks) -- (midden);
\draw [->] (midden) -- node[above] {$\exp_{m}$} (rechts);
\draw [->] (onderlinks) -- node[above] {$\exp_{m}$} (onderrechts);
\draw [->] (onderrechts) -- (rechts);
\end{tikzpicture}
\end{center}
\proof As $\tilde{G}^{+k}(M)$ is a source-simply connected Lie groupoid with 
$J^{k}(TM)$ as Lie algebroid, we can apply 
Lie's second theorem for algebroids. \qed
It is perhaps worth wile to formulate this
for general transitive Lie groupoids, as it clarifies the link with the recent
work of Grabowski, Kotov and Poncin \cite{GKP}. 
For the Atiyah algebroid $A$ of a principal
fibre bundle with a connected, reductive structure group,
they classify Lie algebra 
isomorphisms of $\Gamma(A)$ in terms of the 
Lie algebroid isomorphisms of $A$.    
\begin{proposition}
Let $\mathcal{G}\rightrightarrows M$ be a transitive Lie groupoid, 
with Lie algebroid $A$. The kernel of the anchor $K$ is then a 
bundle of Lie algebras with fixed dimension $d$.  Suppose that the sequence 
$$
0 \rightarrow \Gamma(K) \rightarrow \Gamma(A) \rightarrow 
\Gamma(TM) \rightarrow 0\,,
$$
with $K$ the kernel of the anchor, splits as a sequence of Lie algebras. 
Then this splitting is induced by a morphism of Lie algebroids 
$\nabla : J^{k}(TM) \rightarrow A$,  
and there is a corresponding morphism of Lie groupoids 
$\tilde{G}^{+k}(M) \rightarrow \mathcal{G}$. The number $k$ is at most $3$
if $d$ is $2$ and $\mathrm{dim}(M)=1$, and at most $d$ otherwise.
\end{proposition}
\proof (Or rather a flimsy sketch thereof.) 
Analogous to the case of the gauge groupoid. \qed
We have paved the way for a classification of
infinitesimally natural principal fibre bundles.
\label{fruukkiop}
\begin{theorem}\label{grondsmakreefknol}
Let $\pi : P \rightarrow M$ be an infinitesimally natural principal 
$G$-bundle with splitting $\sigma$ of (\ref{rijtje}).
Then there exists a group homomorphism 
$\rho : G(k,M) \rightarrow G$ such that 
the bundle $P$ is 
associated to $\tf$ through $\rho$, i.e.\ 
$$P \simeq \tf \times_{\rho} G\,.$$
Moreover, $\sigma$ is induced by the canonical one for
$\tf$.
\end{theorem}
\proof Fix a base point $m$ on $M$. The map $\exp \nabla$ 
yields a homomorphism of groups
$\rho : \tilde{G}_{m,m}^{+k}(M) \rightarrow ((P\times P)/G)_{m,m}$, 
the latter isomorphic to $G$, 
the former to 
$G(k,M)$.

The map
$\tilde{G}_{*,m}^{+k}(M) \times_{\rho}  
((P\times P)/G)_{m,m} \rightarrow ((P\times P)/G)_{*,m}$
which is given by $(g_{m' , m} , p_{m,m}) \mapsto 
(\exp \nabla (g_{m' , m})) \cdot p_{m,m}$
is well defined and injective because two pairs share the same 
image if and only
if they are equivalent modulo $\tilde{G}_{m,m}^{+k}(M)$.
It is also surjective and $G$-equivariant, and hence an
isomorphism of principal $G$-bundles. As 
$((P \times P)/G)_{*,m} \simeq P$
and
$\tilde{G}_{*,m}^{+k}(M) \times_{\rho}  
((P\times P)/G)_{m,m} \simeq \smash{\tf \times_{\rho} G}$,
the equivalence is proven. The remark on $\sigma$ 
follows from the construction.
\qed
This classifies the infinitesimally natural principal fibre bundles. 
They are all associated (via a group homomorphism) 
to the bundle $\tilde{G}^{+k}_{*,m} = \tf$.

The classification of natural principal fibre bundles is now an easy
corollary. The following well known result (\cite{PT}, \cite{Te}) states that 
they are precisely the
ones associated to $G^{k}_{*,m}(M) = F^{k}(M)$.
\begin{corollary}\label{PalaisTerng}
Let $\pi : P \rightarrow M$ be a natural principal $G$-bundle
with local splitting $\Sigma$ of (\ref{flapperdrop}). 
Then $P$ is associated to 
$F^{k}(M)$.
That is, there exists a 
homomorphism $\rho : G^{k}_{0,0}(\mathbb{R}^{n}) \rightarrow G$ such that 
$$
P \simeq F^{k}(M) \times_{\rho} G\,.
$$
Moreover, $\Sigma$ is induced by the canonical one for $F^{k}(M)$.
\end{corollary}
\proof
As the homomorphism 
$\Sigma : \mathrm{Diff}_{c}(M) \rightarrow \mathrm{Aut}_{c}(P)$
is local, it induces
a homomorphism of groupoids $\Sigma : 
\mathrm{Germ}(M) \rightarrow (P \times P)/G$,
with $\mathrm{Germ}(M)$ the groupoid of germs of diffeomorphisms of $M$.
We need but show that $\Sigma$ factors through 
$j^{k} : \mathrm{Germ}(M) \rightarrow G^{k}(M)$ for some $k > 0$,
cf. the proof of theorem \ref{grondsmakreefknol}.
 
The Lie algebra homomorphism 
$\sigma : \Gamma_{c}(TM) \rightarrow TP/G$ defined by
$\sigma(v) := \del_{t}|_{0} \Sigma(\exp(tv))$ 
is local by assumption, and according to proposition \ref{alhetvet} 
it factors through 
the $k$-jets for some $k>0$.
It suffices to show that $\Sigma(\phi)_{m,m} = \id_{m,m}$
for any $\phi \in \mathrm{Germ}_{m,m}(M)$ 
that agrees with the identity to $k^{\mathrm{th}}$ order at $m$.

In local co-ordinates $\{x^{i}\}$,
we write $\phi^{i}(x) = x^{i} + v^{i}(x)$, where 
$v : \mathbb{R}^{n} \rightarrow \mathbb{R}^{n}$
vanishes to $k^{\mathrm{th}}$ order.
We define the one parameter family of germs of diffeomorphisms 
$\phi^{i}_{t}(x) := x^{i} + t v^{i}(x)$.
Then 
$\del_{t}|_{\tau} \Sigma(\phi_{\tau})^{-1} \Sigma(\phi_{t})_{m} = 0$,
as it equals $\sigma_{m}( \del_{t}|_{\tau} \phi^{-1}_{\tau} \phi_{t})$,
the image of a vector field that 
vanishes to 
order $k$ at $m$. Therefore  
$t \mapsto \Sigma(\phi_{t})_{m,m}$ is constant, and
$\Sigma(\phi)_{m,m} = \id_{m,m}$ as required.
\qed 
To summarize: natural principal fibre bundles are associated to a
higher frame bundle, whereas infinitesimally natural 
principal fibre bundles are associated to the universal cover of 
a higher frame bundle.  



\subsection{The \texorpdfstring{Bundle $\tf$}{Frame Bundle}} \label{dekhetkframe}
The above considerations prompt a few remarks on the 
universal cover of the connected component of the frame bundle
$\tf$, and on its (disconnected) structure group $G(k,M)$.
Recall that they are just the source fibre $\smash{\tilde{G}^{+k}_{*,m}(M)}$
and isotropy group $\smash{\tilde{G}^{+k}_{m,m}(M)}$ of $\smash{\tilde{G}^{+k}(M)}$.
 
\subsubsection{General Manifolds}
If $\pi_{1}(M)$ is the homotopy groupoid of $M$, define
the homomorphism of groupoids $\mathrm{Pr} : \pi_{1}(M)_{m',m} \rightarrow 
\pi_{0}(G_{m',m}(M))$
by lifting a path in $M$ to a path in $G^{k}(M)$ with fixed source,
and taking the connected component of its end point.
It makes\vspace{-1mm}
\begin{center}
\begin{tikzpicture}
\pgfsetzvec{\pgfpoint{0.385cm}{-0.385cm}}
\node (links) at (0, -0.7, 0) {$\tilde{G}_{m',m}^{+k}(M)$};
\node (rechts) at (2.8, 0, 0) {$\pi_{1}(M)_{m',m}$};
\node (linksonder) at (2.8, -1.4, 0) {$G_{m',m}^{k}(M)$};
\node (rechtsonder) at (5.8,-0.7, 0) {$\pi_{0}(G_{m',m}^{k}(M))$};
\draw [->] (links) -- (rechts);
\draw [->] (links) -- (linksonder);
\draw [->] (linksonder) -- (rechtsonder);
\draw [->] (rechts) -- (rechtsonder);
\end{tikzpicture}\vspace{-1mm}
\end{center}
into a commutative diagram.

Define 
$(G^{k}(M) \times \pi_{1}(M))^{\mathrm{Pr}}$
to be the groupoid of pairs $(g,[f])$
such that $\pi_{0}(g) = \mathrm{Pr}([f])$. 
If $M$ is orientable, this is simply 
$G^{+k}(M) \times \pi_{1}(M)$.
The map of groupoids
$\tilde{G}^{+k}(M) \rightarrow 
(G^{k}(M) \times \pi_{1}(M))^{\mathrm{Pr}}$
is well defined and surjective. 
It restricts to a covering map of principal fibre bundles
\begin{equation}\label{tauprojectie}
\tau : \tilde{G}_{*,m}^{+k}(M) \rightarrow 
(G^{k}(M) \times \pi_{1}(M))_{*,m}^{\mathrm{Pr}}\,.
\end{equation}
The kernel of the corresponding cover of groups 
is precisely
$i_{*}\pi_{1}(\smash{G_{m,m}^{+k}(M)})$, with 
$i : G_{m,m}^{+k}(M) \rightarrow G_{*,m}^{+k}(M)$
the inclusion. 
Note that $i_{*}$ has a 
nonzero kernel precisely when a vertical loop is contractible
in $\smash{G_{*,m}^{+k}(M)}$, but not by a homotopy
which stays inside the fibre. Denoting 
$\pi_{1}(G_{m,m}^{+k}(M))$ by $Z$, we obtain the exact sequence\vspace{-1mm}
\begin{equation}\label{centraalgroep}
1 \rightarrow Z/\mathrm{Ker}(i_{*}) 
\stackrel{i_{*}}{\rightarrow} 
\tilde{G}_{m,m}^{+k}(M)
\stackrel{\tau}{\rightarrow}
(G^{k}(M) \times \pi_{1}(M))_{m,m}^{\mathrm{Pr}}
\rightarrow
1\,.
\end{equation}
A moment's thought reveals that this extension is central: 
if $g(t)$ is a path in
$G_{m,m}^{+k}(M)$ and $h(t)$ one in 
$G_{*,m}^{+k}(M)$, then both $h * (i \circ g)$ and 
$(i \circ g)\cdot h(1) * h$
can be homotoped into $t \mapsto h(t)g(t)$.

We may as well restrict attention to the case $k=1$, 
in which $\smash{G_{*,m}^{1}(M)}$ is the frame bundle $F(M)$.
Indeed, as 
$\smash[t]{G^{+k}_{*,m}(M) \rightarrow G^{1,+}_{*,m}(M)}$ has contractible fibres,
$\smash{\tilde{G}^{+k}_{*,m}(M)}$ is just
the pullback of $G^{+k}_{*,m}(M)$ along 
$\smash{\tilde{G}_{*,m}^{1,+}(M) \rightarrow G_{*,m}^{1,+}(M)}$.

\subsubsection{Orientable Manifolds}
For orientable manifolds, the situation simplifies.
If we identify the connected component of $G_{m,m}^{1}(M)$ with $\glnp$, we 
obtain a homomorphism 
$i_{*}$ of $\smash{\tglnp}$ into $\smash[b]{\tilde{G}^{1,+}_{m,m}(M)}$.
There is a second homomorphism
$\pi_{1}(F(M)) \rightarrow \tilde{G}^{1,+}_{m,m}(M)$.
Their images intersect in $Z/\mathrm{Ker}(i_{*})$, and commute
by an argument similar to the one on centrality of (\ref{centraalgroep}).
If we define $(\smash{\tglnp} \times \pi_{1}(F(M)))_{Z}$ to be the 
quotient of $\smash{\tglnp} \times \pi_{1}(F(M))$ by the equivalence
$(gz,h) \sim (g,zh)$,
we can regard it as a subgroup of $\tilde{G}^{1,+}_{m,m}(M)$.
Note that if $\mathrm{Ker}(i_{*})$ is nonzero, the above equivalence
relation sets it to 1.

If $M$ is orientable, we may restrict our attention to $F^{+}(M)$,
which has connected fibres. Any path in
$F^{+}(M)$ which starts and ends in the same fibre can therefore be obtained 
by combining a closed loop with a path in $\glnp$.
For orientable manifolds, we thus have 
$
\tilde{G}_{m,m}^{1,+}(M) \simeq (\smash{\tglnp} \times \pi_{1}(F^{+}(M)))_{Z}
$, and in the same vein
\begin{equation}\label{orient}
G(k,M) \simeq 
(G(k,\mathbb{R}^{n}) \times \pi_{1}(F^{+}(M)))_{Z}\,.
\end{equation}

\subsubsection{Spin Manifolds}\label{Spinmanifolds}
Let $M$ be an orientable manifold, equipped with a 
pseudo-Riemannian metric $g$ of signature 
$\eta \in \mathrm{Bil}(\mathbb{R}^{n})$.
Then $OF^{+}_{g} := \{f \in F^{+}(M) \,|\, f^* g = \eta\}$
is the bundle of positively oriented orthogonal frames. 
A spin structure is then by definition an
$\widetilde{\mathrm{SO}}(\eta)$-bundle\footnote{
There is a subtlety here. Suppose 
$\eta$ has indefinite signature, say $(3,1)$.
The group $\mathrm{SO}(3,1)$
has 2 connected components,  
so that a universal cover does not exist. 
As it is a subgroup of the simply connected group 
$\glnvp$,
we simply define
$\widetilde{\mathrm{SO}}(3,1)$ to be $\kappa^{-1}(\mathrm{SO}(3,1))$ 
with $\kappa : \tglnvp \rightarrow \mathrm{GL}^{+}(\mathbb{R}^{4})$ the covering map.
Thus $\widetilde{\mathrm{SO}}(3,1)$ is, perhaps surprisingly, 
not isomorphic to 
the 2-component spin group $\mathrm{Spin}(3,1)$. 
Indeed,
if
$T$ is time inversion and $P$
is the inversion of $3$ space co-ordinates, then
$(PT)^{2} = \one$
in $\widetilde{\mathrm{SO}}(3,1)$,
as opposed to
$(PT)^{2} = -\one$ in $\mathrm{Spin}(3,1)$.
Therefore 
$\pi^{-1}(\pm \one) \simeq 
\mathbb{Z}/2\mathbb{Z} \times \mathbb{Z}/2\mathbb{Z}$
in $\widetilde{\mathrm{SO}}(3,1)$,
whereas
$\pi^{-1}(\pm \one) = \mathbb{Z}/4\mathbb{Z}$ in $\mathrm{Spin}(3,1)$
(see \cite{BDWGK}).  

Of course 
the connected component of unity of $\widetilde{\mathrm{SO}}(3,1)$
and that of $\mathrm{Spin}(3,1)$ are both isomorphic 
to $\mathrm{SL}(\mathbb{C}^{2})$, so that
none of this is relevant
if $M$ is both orientable and time-orientable,
i.e.\ if the structure group of 
the frame bundle reduces to $\mathrm{SO}^{\uparrow}(3,1)$.
}
$Q$ over $M$, plus a map $u : Q \rightarrow \smash{OF^{+}_{g}}$ 
such that the following diagram 
commutes, with $\kappa$ the canonical homomorphism 
$\smash{\widetilde{\mathrm{SO}}(\eta)} \rightarrow \mathrm{SO}(\eta)$.
A manifold is called spin if it admits a spin structure.
\begin{center}
\begin{tikzpicture}
\pgfsetzvec{\pgfpoint{0.385cm}{-0.385cm}}
\node (lg) at (-1,0.7) {$\widetilde{\mathrm{SO}}(\eta)$};
\node (rg) at (1,0.7) {$\mathrm{SO}(\eta)$};
\node  (rpijl) at (1,0.3) {$\curvearrowleft$};
\node  (lpijl) at (-1,0.3) {$\curvearrowleft$};
\node (lb) at (-1,0) {$Q$};
\node (rb) at (1,0) {$OF^{+}_{g}$};
\node (lo) at (0,-1.2) {$M$};
\draw [->] (lg) --node[above]{$\kappa$} (rg);
\draw [->] (lb) --node[above]{$u$} (rb);
\draw [->] (lb) -- (lo);
\draw [->] (rb) -- (lo);
\end{tikzpicture}
\end{center}

Define $\hat{Q} := Q \smash{\times_{\widetilde{\mathrm{SO}}(\eta)}} \tglnp$, 
and let us again denote the induced
map $\hat{Q} \rightarrow F^{+}(M)$ by $u$. 
As any cover of $F^{+}(M)$ by a $\smash{\tglnp}$-bundle
can be obtained in this way, there is a 1:1 correspondence
between spin covers of $OF^{+}_{g}(M)$ and $F^{+}(M)$.
In particular, whether or not $M$ is spin does not depend on the metric.

The Serre spectral sequence gives rise to the exact sequence
\begin{equation} \label{piraterij}
1 \rightarrow Z/\mathrm{Ker}(i_{*}) \rightarrow 
\pi_{1}(F^{+}(M)) \rightarrow
\pi_{1}(M) \rightarrow 1\,.
\end{equation}
The following proposition is well known.
\begin{proposition}\label{klopboor?}
A spin structure exists
if and only if $i_{*} : Z \rightarrow \pi_{1}(F^{+}(M))$
is injective and
(\ref{piraterij}) splits 
as a sequence of groups. If spin structures exist, 
then 
equivalence classes of spin covers correspond to
splittings of (\ref{piraterij}).
\end{proposition}
\proof This will follow from theorem \ref{naturalspin} later
on, but
see e.g. \cite{Mo} for an independent proof. Our criterion for 
$M$ to be spin is equivalent to the vanishing of the second Stiefel-Whitney
class, see e.g. \cite{LM}.
\qed

{{\bf \noindent Remark} \quad
In terms of group cohomology, one can consider the sequence 
(\ref{piraterij}) as an element
$[\omega] \in H^{2}(\pi_{1}(M) , Z/\mathrm{Ker}(i_*))$.
Spin bundles exist if and only if both $\mathrm{Ker}(i_*)$  
and $[\omega]$ are trivial, in which case they are indexed by
$H^{1}(\pi_{1}(M) , Z)$.\vspace{-2mm}
}\medskip

\noindent If a spin structure exists, then $\smash{\widetilde{F}^{+}}$
is simply the pullback along
the universal cover $\tilde{M} \rightarrow M$ of  
$\hat{Q} \rightarrow M$. 
The picture then becomes\vspace{2mm} 
\begin{center}
\begin{tikzpicture}
\pgfsetzvec{\pgfpoint{0.385cm}{-0.385cm}}
\node (mb) at (0,0) {$u_{*}\big( F^{+k}(M) \big)$};
\node (mo) at (0,-1.5) {$\hat{Q}$};
\node (lb) at (-3,0) {$\tilde{F}^{+k}(M)$};
\node (rb) at (3,0) {$F^{+k}(M)$};
\node (lo) at (-3,-1.5) {$\tilde{F}^{+}(M)$};
\node (ro) at (3,-1.5) {$F^{+}(M)$};
\node (moo) at (0,-3) {$M$};
\node (loo) at (-3,-3) {$\tilde{M}$};
\draw [->] (lb) -- (mb);
\draw [->] (mb) -- (rb);
\draw [->] (lb) -- (lo);
\draw [->] (lo) -- (mo);
\draw [->] (mo) -- node[above]{$u$} (ro);
\draw [->] (rb) -- (ro);
\draw [->] (mb) -- (mo);
\draw [->] (mo) -- (moo);
\draw [->] (ro) -- (moo);
\draw [->] (lo) -- (loo);
\draw [->] (loo) -- (moo);
\end{tikzpicture}\vspace{0mm}
\end{center}
with each of the three squares a pullback square.

\section{More General Fibre Bundles}\label{zes}

In this section, we will prove 
a version of theorem \ref{grondsmakreefknol} 
for fibre bundles which are not principal.
It would however be overly optimistic to expect an analogue
of of theorem \ref{grondsmakreefknol} to hold for arbitrary 
smooth fibre bundles, so we will restrict ourselves to those bundles 
that carry a
sufficiently rigid structure on their fibres.




\subsection{Structured Fibre Bundles}
We start by making this statement more precise.

\begin{basdef} Let $\mathbf{C}$ be a subcategory of the category of smooth manifolds such that 
the group of automorphisms of each object of $\mathbf{C}$  
is a finite dimensional Lie group. 
Then a `structured fibre bundle'
with structure $\mathbf{C}$ and fibre $F_0 \in \mathrm{ob}(\mathbf{C})$ is 
by definition a smooth fibre bundle $\pi : F \rightarrow M$ where the fibres
are objects of $\mathbf{C}$. We also require
each point to possess a neighbourhood $U$ and a local trivialization 
$\phi :  \pi^{-1}(U) \rightarrow F_0 \times U$ such that $\phi$ restricted to
a single fibre is a $\mathbf{C}$-isomorphism $\pi^{-1}(m) \rightarrow F_0$.  
\end{basdef}
For example, a structured fibre bundle in the category of finite dimensional
vector spaces is a vector bundle.

If $\pi : F \rightarrow M$ is any smooth fibre bundle, then an 
automorphism of $\pi$ 
is by definition a 
diffeomorphism $\alpha$ of $F$ such that  
$\pi(f) = \pi (f')$ implies $\pi(\alpha(f)) = \pi(\alpha(f'))$.
It is called vertical if it maps each fibre to itself.

\begin{basdef} We define an automorphism of a structured fibre bundle to be an automorphism
of the smooth fibre bundle such that its restriction to each single fibre 
is an isomorphism in $\mathbf{C}$.      
\end{basdef}    
One can then construct a sequence of groups 
\begin{equation} \label{structuurgroep}
1 \rightarrow \mathrm{Aut}^{\mathbf{C}}_c(F)^{V} 
\rightarrow \mathrm{Aut}^{\mathbf{C}}_c(F)
\rightarrow \mathrm{Diff}^{F}_c(M)
\rightarrow 1 
\end{equation}
and its corresponding exact sequence of Lie algebras
\begin{equation} \label{structuurvezel}
0 
\rightarrow \Gamma^{\mathbf{C}}_c(TF)^{V} 
\rightarrow \Gamma^{\mathbf{C}}_c(TF)^{P}
\rightarrow \Gamma_c(TM)
\rightarrow 
0 \,,
\end{equation}
where `$V$' is for vertical, `$P$' for projectable, and $c$
again stands for `$0$ outside a compact subset of $M$'.

The proof of the following corollary of theorem \ref{grondsmakreefknol} 
is now a formality.
\begin{corollary}\label{structuurcor}
Let $\pi : F \rightarrow M$ be a structured fibre bundle with fibre $F_0$
such that 
(\ref{structuurvezel}) splits as a sequence of Lie algebras.
Then there exists an action $\rho$ of $G(k,M)$ 
by $\mathbf{C}$-automorphisms on a single fibre $F_m$ such that
$$
F = \tf \times_{\rho} F_0 \,.
$$
\end{corollary}

\proof
Construct the principal $\mathrm{Aut}^{\mathbf{C}}(F_0)$-bundle 
$\pi : P \rightarrow M$, the fibre over $x$
of which is precisely the set of $\mathbf{C}$-isomorphisms 
$\phi : F_{m} \rightarrow F_{x}$.
Then there is a natural isomorphism 
$\mathrm{Aut}^{\mathbf{C}}_{c}(F) \simeq \mathrm{Aut}_{c}(P)^{G}$ under which 
the vertical subgroups of the two correspond, so that the exact sequence 
(\ref{structuurgroep}) is isomorphic to (\ref{flapperdrop}),
and therefore (\ref{structuurvezel}) to (\ref{rijtje}).
As $F = P \times_{\mathrm{Aut}(F_m)} F_{m}$,
we can now apply theorem \ref{grondsmakreefknol} to $P$ in order to substantiate our
claim. \qed

\subsection{Vector Bundles}
We specialize to the case of vector bundles.  
These are precisely structured fibre bundles in the category of finite
dimensional vector spaces.

The exact sequence of Lie algebras (\ref{structuurvezel}) for a vector bundle
$E$ with fibre $V$ is then\vspace{-0mm}
\begin{equation}\label{veclie}
1 \rightarrow \mathrm{DO}^{0}_{c}(E) \rightarrow \mathrm{DO}^{1}_{c}(E) 
\rightarrow \mathrm \Gamma_{c}(TM) \rightarrow 0\,,
\end{equation}
where $\mathrm{DO}_{c}^{1}(E)$ is the Lie algebra of compactly supported 
$1^{\mathrm{st}}$ order differential operators on $\Gamma(E)$, and
$\mathrm{DO}_{c}^{0}(E)$ the ideal of $0^{\mathrm{th}}$ order ones, 
that is to say 
$\mathrm{DO}_{c}^{0}(E) \simeq \Gamma_{c} (E \otimes E^{*})$.  

Corollary \ref{structuurcor} then says that \smash{(\ref{veclie})}
splits as a sequence of Lie algebras if and only if there
is a representation $\rho$ of $G(k,M)$ on $V$ 
such that
$
E \simeq \tf \times_{\rho} V\,.
$

But thanks to the fact that all finite dimensional representations
of the universal cover of $\mathrm{GL}^{+}(\mathbb{R}^{n})$
factor through $\mathrm{GL}^{+}(\mathbb{R}^{n})$ itself, we can 
even say something slightly stronger.

\begin{proposition} \label{structuurvector}
Let $E \rightarrow M$ be a vector bundle 
for which (\ref{veclie}) splits as a sequence of Lie algebras. 
Then there exists a representation
$\rho$ of the group $(G^{k} \times \pi_{1}(M))_{m,m}^{\mathrm{Pr}}$ on V such that 
$$
E \simeq (G^{k}(M) \times \pi_{1}(M))_{*,m}^{\mathrm{Pr}} 
\times_{\rho} V\,.
$$
\end{proposition}
{{\bf \noindent Remark} \quad
If $M$ is orientable, this reads   
$
E \simeq \pi^{*} F^{+k}(M) 
\times_{\rho} V
$. In this expression, $\pi^{*} F^{+k}(M)$ is the pullback of $F^{+k}(M)$ along
$\pi : \tilde{M} \rightarrow M$, considered as a principal
$G_{0,0}^{+k}(\mathbb{R}^{n}) \times \pi_{1}(M)$-bundle
over $M$.}

\proof
%
Consider the restriction of the map $\tau$ in equation 
(\ref{tauprojectie}) to the group $\tilde{G}_{m,m}^{+k}(M)$.
In order to prove the proposition, we need but show that 
its kernel $Z$ 
acts trivially on $V$.
For $k=0$, this is clear.

If $k$ is at least 1, 
the homomorphism 
$
\tglnp
\rightarrow
\tilde{G}_{m,m}^{+k}(M)
$
makes $V$ into a finite dimensional representation
space for $\tglnp$.  
But it is known (see \cite[p.\,311]{Kn}) that all finite dimensional
representations of its cover factor through $\glnp$
itself. This implies that
the subgroup 
$Z$ which covers the identity
must act trivially on $V$,
and we may consider \vspace{-1mm}
$$\tilde{G}_{*,m}^{+k}(M) / Z 
\simeq 
(G^{k}(M) \times \pi_{1}(M))_{*,m}^{\mathrm{Pr}}\vspace{-1mm}
$$
to be the underlying bundle, as announced.\qed

\noindent This reduces the problem of classifying vector bundles with split 
sequence (\ref{veclie})
to the representation
theory of $\smash{(G^{k} \times \pi_{1}(M))_{m,m}^{\mathrm{Pr}}}$.
%

The above extends a result \cite{Te} of Terng, in which she
classifies vector bundles which allow for a local splitting of the sequence of 
groups (\ref{structuurgroep}). 
It is an extension first of all in the sense that 
we prove, rather than assume, that the splitting is local.
Secondly, we have shown that in classifying 
vector bundles
with split sequence (\ref{veclie}) of Lie algebras 
rather than groups, one encounters only slightly more.  
Intuitively speaking, the extra bit is the representation theory
of $\pi_{1}(M)$.
We refer to \cite{Te} for a thorough exposition of the representation theory
of $G^{k}_{0,0}(\mathbb{R}^{n})$.

\section{Flat Connections}\label{vijf}

Having concluded our classification of bundles
in which (\ref{rijtje}) is split as a sequence
of Lie algebras, the time has come to
apply our newly acquired knowledge.

In this section, we will investigate splittings
that come from a flat equivariant connection on a 
principal $G$-bundle $P \rightarrow M$.
We will prove that if the Lie algebra 
$\mathfrak{g}$ of $G$
does not contain $\mathfrak{sl}(\mathbb{R}^{n})$
as a subalgebra, 
then 
the sequence of Lie algebras (\ref{rijtje}) splits
if and only if $P$ admits a flat equivariant connection.
In other words, the sequence (\ref{rijtje}) then splits
as a sequence of Lie algebras if and only if it
splits as a sequence of Lie algebras and $C^{\infty}(M)$-modules.

Note that this is certainly not the case for general groups $G$.
The frame bundle for example always allows for a splitting of 
(\ref{rijtje}), but usually not for a flat connection.

\subsection{Lie Algebras That Do Not Contain 
\texorpdfstring{$\mathfrak{sl}(\mathbb{R}^{n})$}{the Special Linear Lie Algebra}}

Although lemma \ref{hallekal} exhibits $\sigma$
as a differential operator of finite order,  
the bound on the order is certainly not optimal. 
With full knowledge of the Lie algebras
at hand, sharper restrictions can be put 
on the kernel of $\sigma$.
In particular, if $\mathfrak{g}$
does not contain $\mathfrak{sl}(\mathbb{R}^{n})$,
there is only a single relevant ideal, and $\sigma$ is of order 
at most 1.
For notation, see section \ref{localmap}.
%


\begin{lemma}\label{grondkusrelmuis}
Let $n=1$, and let $\mathfrak{g}$ be such that it does not 
contain two nonzero elements such that $[X,Y] = Y$. 
Or
let $n \geq 2$, and let $\mathfrak{g}$ be such that it does not
admit $\mathfrak{sl}(\mathbb{R}^n)$ as a subalgebra.
Then the kernel of the homomorphism 
$\check{f}_m : \vecn \rightarrow \mathfrak{g}$ 
contains 
$\{v \in \vecn \, | \, \mdiv_m(v) = 0\}$.
\end{lemma}
\proof
We start with the case $n=1$. Again, we note that the
only ideals of
$\mathrm{Vec}_{1} = \mathrm{Span}\{x^{k} \del_x \,|\, k \geq 1\}$
are 
$\mathrm{Span}\{x^2\del_x , x^k \del_x \,| \, k \geq 4 \}$,
and for each $N\geq1$ an ideal
$\mathrm{Span}\{x^k \del_x \,| \, k \geq N \}$.  
The corresponding quotients all contain two elements $X$ and $Y$ with $[X,Y]=Y$,
except the ideals corresponding to $N=1,2$.
This means that also $\mathfrak{g}$, containing
$\mathrm{Vec}_{1} / \ker(\check{f}_m)$ as the image of $\check{f}_m$,  
will possess $X$ and $Y$ such that $[X,Y] = Y$ unless 
the kernel of $\check{f}_m$ contains the ideal
$\mathrm{Span}\{x^k \del_x \,| \, k \geq 2 \} = 
\{v \in \mathrm{Vec}_1 \, | \, \mdiv_m(v) = 0\}$.

Now for $n \geq 2$.
Under the identification $\vecn^{0} \simeq \mathfrak{gl}(\mathbb{R}^n)$
 given by $x_{i}\del_j \mapsto e_{i j}$,
the Euler vector field is the identity $\mathbf{1}$
and $\mdiv_{m}$ becomes the trace. 
As $\ker(\check{f}_m)^{0}$ is an ideal in $\mathfrak{gl}(\mathbb{R}^n)$,
it can be either $0$, $\mathbb{R} \mathbf{1}$, $\mathfrak{sl}(\mathbb{R}^n)$
or $\mathbb{R} \mathbf{1} \oplus \mathfrak{sl}(\mathbb{R}^n)$. 
In the former two cases, 
$\mathrm{Im}(\check{f}_m) \simeq \vecn / \ker(\check{f}_m)$,
and hence $\mathfrak{g}$,  
would contain $\mathfrak{sl}(\mathbb{R}^n)$ as a subalgebra, 
contradicting the hypothesis.
Hence $\mathfrak{sl}(\mathbb{R}^n) \subseteq \ker(\check{f}_m)$.
If we now show that 
$[\vecn , \mathfrak{sl}(\mathbb{R}^n)] = 
\mathfrak{sl}(\mathbb{R}^n) \bigoplus_{k=1}^{\infty} \vecn^{k}$,
the proof will be complete.

Let $i \neq j$. We then have 
$[x_{i}\del_j , x_{j}x^{\vec{\alpha}}\del_{j}] = (\alpha_j + 1) x_i
x^{\vec{\alpha}} \del_j$, showing that 
$x_i x^{\vec{\alpha}} \del_j \in [\vecn , \mathfrak{sl}(\mathbb{R}^n)]$.
The only basis elements not of this shape are of the form $x_{j}^{k}\del_j$.
But $[x_{j}\del_i , x_i x_{j}^{k-1}\del_{j}] = 
x_{j}^{k}\del_j - x_i x_{j}^{k-1}\del_{i}$. If $k \geq 2$, the latter part 
was just shown to be in $[\vecn , \mathfrak{sl}(\mathbb{R}^n)]$, so that also 
$x_{j}^{k}\del_j \in [\vecn , \mathfrak{sl}(\mathbb{R}^n)]$. 
If $k=1$, the elements 
$x_j \del_j - x_i \del_i$ join $x_i \del_j$ to form a basis of 
$\mathfrak{sl}(\mathbb{R}^n)$.
\qed

\noindent This rather limits the possibilities. 
Not only can we restrict to first order, but also 
the Lie algebroid map 
$\nabla : J^{1}(TM) \rightarrow TP/G$ vanishes on
the trace-zero jets
$K_{m} = \{\, j_{m}^{1}(v) \in J^{1}(TM) \, | \, v(m) = 0 \, \, 
\mathrm{and} \,\,
\tr(v)=0\, \}$,
so that it factors through the `trace Lie algebroid' 
$\tr_m(M) := J_{m}^{1}(TM)/K_{m}$.

This in turn is the Lie algebroid of the `determinant groupoid'
$\mathrm{Det}(M)$. An element $[\alpha]_{m',m}$ of $\mathrm{Det}(M)_{m',m}$
is by definition an equivalence class of diffeomorphisms mapping 
$m$ to $m'$, with $\alpha \sim \beta$ if and only if
$\mathrm{Det}(\beta^{-1}\alpha) = 1$. 

As $[\alpha]_{m',m}$ identifies $\wedge^{n}(T^{*}_{m}M)$ 
with $\wedge^{n}(T^{*}_{m'}M)$,
the source fibre $\mathrm{Det}(M)_{*,m}$ is 
isomorphic to the determinant line bundle $\wedge^{n}(T^{*}M) \rightarrow M$.
Its\footnote{An isomorphism 
$\wedge^{n}(T^{*}M) \simeq \mathrm{Det}(M)_{*,m}$ is only given after 
a choice of $\lambda_{0} \in \wedge^{n}(T_{m}^{*}M)$. This determines 
the connected component.} 
connected component $\wedge^{n,+}(T^{*}M)$ is the 
the bundle of positive top forms
if $M$ is orientable, and the whole bundle otherwise.
 
Its universal covering space is the bundle $\wedge^{n,+}(T^{*}\tilde{M})$
of positive top forms on $\tilde{M}$.
Indeed, 
$\tilde{M}$ is always orientable, regardless of whether or not $M$ is.
This means that 
$\wedge^{n}(T^{*}\tilde{M})$ is a trivial bundle,
and that its connected component 
$\wedge^{n,+}(T^{*}\tilde{M}) \simeq \tilde{M} \times \mathbb{R}^{+}$ 
is simply connected. The covering map is induced by the map 
$\tilde{M} \rightarrow M$. This leads to the following version of 
theorem \ref{grondsmakreefknol}.


\begin{proposition}\label{habbelsnap}
Let $P$ be a principal $G$-bundle over an $n$-dimensional manifold $M$.
Let $G$ be such that its Lie algebra $\mathfrak{g}$ does not 
contain $\mathfrak{sl}(\mathbb{R}^{n})$ if $n > 1$,
or $[X,Y] = Y$ if $n=1$. 
Then there is a homomorphism $\pi_{1}(M) \times \mathbb{R}^{+} \!\rightarrow\! G$
associating $P$ to the principal $\pi_{1}(M) \times \mathbb{R}^{+}$-bundle
$
\wedge^{n,+}(T^{*}\tilde{M}) \rightarrow M 
$.
$$
P \simeq \wedge^{n,+}(T^{*}\tilde{M}) \times_{\pi_{1}(M) \times \mathbb{R}^{+}}G\,.
$$
\end{proposition}
We may even classify the possible splittings. 
\begin{corollary}\label{grondsmakgu}
Under the hypotheses of proposition \ref{habbelsnap}, any 
Lie-algebra homomorphism 
$\sigma: \Gamma_{c}(TM) \rightarrow \Gamma_{c}(TP)^{G}$
which splits the sequence of Lie algebras (\ref{rijtje}) 
can be written
\begin{equation}\label{raldaldalgu}
\sigma = \nabla^{\mu} + \Lambda \mdiv_{\mu}\,,
\end{equation}
where $\nabla^{\mu}$ is a flat equivariant connection on $P$,
and $\Lambda$ a section of $\mathrm{ad}(P)$ which is 
constant w.r.t. the connection induced on $\mathrm{ad}(P)$ by $\nabla^{\mu}$. 
\end{corollary}
{{\bf \noindent Remark} \quad
In particular, this shows that there exists a flat
connection which splits (\ref{rijtje}), even though 
most splittings are not flat connections.}

\proof First, we prove the case $P = \wedge^{n,+}(T^{*}\tilde{M})$. 
Pick a nonzero \mbox{(pseudo-)} density $\mu$ on $M$. This induces 
an honest density $\tilde{\mu}$ on $\tilde{M}$, which in turn
identifies 
$
\wedge^{n,+}(T^{*}\tilde{M}) 
$
with
$\tilde{M} \times \mathbb{R}^{+} 
$. 
The local trivializations of
$\wedge^{n,+}(T^{*}\tilde{M}) \rightarrow \tilde{M}$
and 
$\tilde{M} \rightarrow M$
combine to locally trivialize 
$\wedge^{n,+}(T^{*}\tilde{M}) \rightarrow M$.
This yields a flat 
equivariant
connection $\nabla^{\mu}$ on $\wedge^{n,+}(T^{*}\tilde{M})$,
which annihilates $\tilde{\mu}$.

The splitting $\sigma$ is uniquely determined by the action of $\sigma(v)$ 
on local sections $\tilde{\nu}$, 
which reads 
$\sigma(v) (\tilde{\nu}) = \pi^{*} \circ \mathcal{L}_{v} \circ 
\pi^{* -1}
\tilde{\nu}$,
where $\pi$ is the map from $\tilde{M}$ to $M$. 
If we define the divergence w.r.t. $\mu$
by the requirement that the Lie derivative
$\mathcal{L}_v \mu $ equal $\mathrm{Div}_{\mu}(v) \mu$,
then we have
\begin{eqnarray*}
\sigma(v)(f\tilde{\mu}) &=& \pi^{*} (\mathcal{L}_{v}(f\mu))\\
&=& \pi^{*}(v(f) \mu + \mathrm{Div}_{\mu}(v) f \mu)\\
&=& v(f) \tilde{\mu} + \mathrm{Div}_{\mu}(v) f \tilde{\mu}\\
&=& \nabla^{\mu}_{v} (f\tilde{\mu}) + \mathrm{Div}_{\mu}(v) f \tilde{\mu}\,.
\end{eqnarray*}
This shows that $\sigma(v) = \nabla^{\mu}_{v} + \Lambda \mathrm{Div}_{\mu}(v)$, 
with $\Lambda = \del_{r}$, the equivariant vertical vector field defined by the 
action of $\mathbb{R}^{+}$ on $\wedge^{n}(T^{*}\tilde{M})$.
(Equivariant vertical vector fields on $P$ 
correspond to sections of $\mathrm{ad}(P)$.)

The general case follows by proposition \ref{habbelsnap}. \qed

\subsection{Lie Algebra Cohomology}

If we specialize to the case of a trivial bundle over an 
Abelian group $G$, we find ourselves in the realm of Lie algebra cohomology.
The continuous cohomology of the Lie algebra of vector 
fields with values in the functions has already been 
unravelled in all degrees \cite{Fu}.
Corollary \ref{grondsmakgu} describes this cohomology only in
degree 1, but now with all cocycles rather than just the continuous ones.
\begin{corollary}\label{homologerengu}
Let $H_{LA}$ denote Lie algebra cohomology and $H_{dR}$ de Rham cohomology.
Let $\mathfrak{g}$ be Abelian, and consider the representation
$C_{c}^{\infty}(M,\mathfrak{g})$ of $\Gamma_{c}(TM)$ where 
a vector field $v$ acts by the Lie 
derivative $\mathcal{L}_v$. Then 
$$
H^1_{LA}(\Gamma_{c}(TM) ,C_{c}^{\infty}(M,\mathfrak{g}))
\simeq
H^1_{dR}(M,\mathfrak{g}) \oplus \mathfrak{g}\,.
$$
\end{corollary}
\proof 
Consider the trivial
bundle $M \times G \rightarrow M$ over an Abelian Lie group $G$, 
which comes equipped with a 
flat connection $\nabla^{0}$,
which acts as Lie derivative.  
Note that Abelian $\mathfrak{g}$ certainly satisfy 
the conditions of propositions \ref{habbelsnap} and \ref{grondsmakgu}. 
View $\Gamma_{c}(\mathrm{ad}(P)) \simeq C^{\infty}_{c}(M,\mathfrak{g})$
as a representation of $\Gamma_{c}(TM)$, and consider its
Lie algebra cohomology. An $n$-cochain is an alternating linear map 
$\Gamma_{c}(TM)^{n} \rightarrow C_{c}^{\infty}(M,\mathfrak{g})$.
For $f^{1} \in C^{1}$, closure $\delta f^{1} = 0$ amounts to
$$
\mathcal{L}_{v}f^{1}(w) - \mathcal{L}_{w}f^{1}(v) - f^{1}([v,w]) 
= 0\,. 
$$
Due to this cocycle condition, $\sigma = \nabla^{0} + f^{1}$ 
is once again a Lie algebra homomorphism splitting $\pi_{*}$.
According to corollary \ref{grondsmakgu},
it must therefore take the shape 
$\sigma = \nabla^{\mu} + \Lambda \mdiv_{\mu}$,
where $\Lambda \in \mathfrak{g}$ is constant.
One can write $\nabla^{\mu} = \nabla^{0} + \omega^{1}$
for some closed 1-form $\omega^{1}$, so that
$f^{1} = \omega^1 + \Lambda \mdiv_{\mu}$.
This classifies the closed $1$-cocycles.


Exact 1-cocycles satisfy $f^1(v) = \delta f^{0} (v) = 
\mathcal{L}_{v}f^{0} = df^{0}(v)$, with $f^0$ a 
$0$-cocycle, that is an element of
$C^{\infty}(M,\mathfrak{g})$.

Note that a change of density $\mu' = e^{h}\mu$
alters $f^1$ by a mere coboundary $\Lambda dh$,
so that the choice of $\mu$ is immaterial. 
The class of $\omega^1 + \Lambda \mdiv_{\mu}$
modulo $\delta C^0$ is therefore determined 
by $[\omega^1] \in H^1_{dR}(M,\mathfrak{g})$
and $\Lambda \in \mathfrak{g}$. 
\qed

\noindent Continuity turns out to be implied by the closedness-condition.
A similar situation was encountered by Takens in \cite{Takens73}, when proving that
all derivations of $\Gamma_{c}(TM)$ are inner, i.e.\
$H^1_{LA}(\Gamma_{c}(TM) ,\Gamma_{c}(TM)) = 0$.

\section{General relativity, Gauge theory, and Spinors}\label{zeven}
In this section, we briefly reflect on the relationship
between (generalized) spin structures and principal fibre
bundles with a splitting of (\ref{rijtje}). 

We formulate general relativity and gauge theory, 
including the fermionic fields, 
in terms of principal fibre bundles over the 
manifold $M$ which describes space-time.
For convenience, we will take
$M$ to be a smooth and orientable manifold of dimension at least 2.

\subsection{General Relativity}

The fundamental degrees of freedom in general relativity
are a pseudo-Rie\-mann\-ian metric $g$ on space-time $M$, 
and a connection $\nabla$ on $TM$. 

We identify the metric
with a section of $F^{+}(M)/\mathrm{SO}(\eta)$ by 
associating to
$g_{x} : T_{x}M \times T_{x}M \rightarrow \mathbb{R}$ the coset 
of all frames $f : \mathbb{R}^{n} \rightarrow T_{x}M$ 
such that $f^{*}g = \eta$.
We identify the connection on $TM$ with an equivariant connection on 
$F^{+}(M)$, that is a section of
$J^{1}(F^{+}(M))/\mathrm{GL}^{+}(\mathbb{R}^{n}) \rightarrow M$.
(Its value at $x$ is the class of a
1-jet at $x$ of a section $\phi \in \Gamma(F^{+})$ with 
$\nabla_{x} \phi = 0$ \cite{Sr}.) 
These two, the metric $g$ and the connection $\nabla$, 
are conveniently combined into
a single section of the fibre bundle $J^{1}(F^{+}(M))/\mathrm{SO}(\eta)$. 

The dynamics of the theory are then governed by 
the Einstein-Hilbert action 
$S_{EH} : \Gamma(J^{1}(F^{+}(M))/\mathrm{SO}(\eta)) \rightarrow \mathbb{R}$,
defined in terms of the Ricci scalar $R$ by
$$
S_{EH}(g,\nabla) = \int_{M} R(g,\nabla) \sqrt{|g|}dx_{0}\ldots dx_{n}\,.
$$
In the Einstein-Hilbert approach, the connection is constrained 
to equal the Levi-Civita connection, and only the metric is varied.
In the Palatini approach \cite{As},
the connection varies 
independently, and the fact that $\nabla$ is the Levi-Civita
connection of the metric is a consequence of the field equations.




Either way, it is clear that the fields transform in a 
natural fashion under diffeomorphisms of $M$.
Any orientation-preserving diffeomorphism $\alpha$ of $M$ lifts to an 
automorphism $\Sigma(\alpha)$ of 
$J^{1}(F^{+}(M))/\mathrm{SO}(\eta)$, defined
by $\Sigma(\alpha)(j^{1}_{m}(\phi)) = j^{1}_{m}(\alpha_{*}\circ \phi \circ
\alpha^{-1})$.
This splits the sequence of groups
$$
1 
\rightarrow
\mathrm{Aut}^{V}( J^{1}(F^{+}(M))/\mathrm{SO}(\eta)) )
\rightarrow
\mathrm{Aut}(J^{1}(F^{+}(M))/\mathrm{SO}(\eta))
\rightarrow
\mathrm{Diff}^{+}(M)
\rightarrow
1\,.
$$

This splitting is central to the theory of general relativity.
The requirement that the action be invariant under co-ordinate 
transformations, $S_{EH}(\Sigma(\alpha) \phi) = S_{EH}(\phi)$,
cannot even be formulated without providing $\Sigma$ explicitly.

Note that as the above sequence of groups splits, so does
the corresponding sequence of Lie algebras of vector fields.





\subsection{Fermions and Spin Structures}
We wish to describe fermions.
As these are known to transform under Lorentz transformations 
by a projective representation rather than a linear one,
we must extend our framework.\vspace{-1.5mm}
\subsubsection{The Bundle}\vspace{-1.5mm}
Suppose that we have a spin structure $Q$ w.r.t
a background metric $g$. 
Let $u : \hat{Q} \rightarrow F^{+}(M)$ be as in section \ref{Spinmanifolds}, and
let $V$ be a unitary spinor representation\footnote{
The indefinite article is appropriate since there is a choice involved here.
The connected component of $\one$ of $\widetilde{\mathrm{SO}}(3,1)$ 
is $\mathrm{Spin}^{\uparrow}(3,1) \simeq \mathrm{SL}^{2}(\mathbb{C})$. 
A spinor representation 
for the connected component can then be unambiguously derived from 
a Clifford algebra representation \cite{Hermann74}.
But as $\widetilde{\mathrm{SO}}(3,1)$ is not isomorphic to
$\mathrm{Spin}(3,1)$, the action of the order 2 central elements
covering $PT$ will have to be specified `by hand'.
Again, this is not relevant if $M$ is time-orientable as well as 
orientable.
} 
of 
$\smash{\widetilde{\mathrm{SO}}}(\eta)$.
Then one has the composite bundle\vspace{-1mm}
$$\hat{Q} \times_{\widetilde{\mathrm{SO}}(\eta)} V 
\rightarrow F^{+}(M)/\mathrm{SO}(\eta)
\rightarrow M\,.\vspace{-1.5mm}
$$
A section $\tau : M \rightarrow \hat{Q} \times_{\widetilde{\mathrm{SO}}(\eta)} V$
can then be interpreted as a metric $g$
along with a fermionic field $\psi$.
Consider $g$ as the induced section of
$F^{+}(M)/\mathrm{SO}(\eta)$.  
Use $g$ to construct the spinor bundle
$u^{-1}(OF^{+}_{g}) 
\times_{\widetilde{\mathrm{SO}}(\eta)} V$,
and obtain a section $\psi$  
by simply restricting the image of $\tau$.
 
In the same vein, we will describe physical fields by sections of
the fibre bundle
\mbox{$J^{1}(\hat{Q})\smash{\times_{\widetilde{\mathrm{SO}}}(\eta)}V$}. 
This is equivalent to providing a triple of sections: one of
$F^{+}(M)/\mathrm{SO}(\eta)$, one of
$J^{1}(F^{+}(M))/\mathrm{GL}^{+}(\mathbb{R}^{n})$,
and one of
$u^{-1}(OF^{+}_{g}) \times_{\widetilde{\mathrm{SO}}(\eta)} V$.
These correspond to  
the metric $g_{\mu \nu}$,
the Levi Civita-connection $\smash{\Gamma^{\alpha}_{\mu \beta}}$,
and the 
fermionic field $\smash{\psi^{i}}$ respectively.\vspace{-1.5mm} 
\subsubsection{Transformation Behaviour}\vspace{-1.5mm}
Let us investigate its transformation behaviour.
As a spinor changes sign under a $2\pi$-rotation, there is no hope 
of finding an interesting homo\-morphism of groups 
$\mathrm{Diff}^{+}(M) \rightarrow \mathrm{Aut}
(J^{1}(\hat{Q})\times_{\widetilde{\mathrm{SO}}(\eta)}V)$.
There is however a canonical homomorphism of Lie algebras.


Because $u : \hat{Q} \rightarrow F^{+}(M)$ has discrete fibres, it
has a unique flat equivariant connection $\nabla$.
This means that the exact sequence of Lie algebras\vspace{-1mm}
$$
0\rightarrow
\Gamma_{c}(\mathrm{ad}(\hat{Q}))
\rightarrow
\Gamma_{c}(T\hat{Q})^{\tglnp}
\rightarrow
\Gamma_{c}(TM)
\rightarrow 0\vspace{-1mm}
$$
is split by $\sigma := \nabla \circ D$, with
$D : \Gamma_{c}(TM) \rightarrow \Gamma_{c}
(TF^{+}(M))^{\glnp}$
the first order derivative $D(v)_f = \del/\del t|_{0} \exp(tv)_{*}\circ f$.
This induces a splitting for $J^{1}(\hat{Q})$ by prolongation 
(see e.g. \cite{FR}), and consequently also one for
$J^{1}(\hat{Q})\times_{\smash[b]{\widetilde{\mathrm{SO}}(\eta)}}V$.
Note that $u_{*}\circ \sigma $ equals $D$.

We see that $\hat{Q}$ is an infinitesimally natural bundle,
and that the canonical splitting $\sigma$ of (\ref{rijtje})  
does not come from a splitting of groups.
We would like to emphasize that even if a splitting at the level of 
groups does exist, it will not be physically relevant, because 
it cannot reproduce the minus sign under a full rotation
that one knows and loves in fermions.

Take for example
the spin structure 
$Q = \mathbb{R}^{n} \times \widetilde{\mathrm{SO}}(\eta)$ 
over $\mathbb{R}^{n}$,
and 
lift $\alpha \in \mathrm{Diff}^{+}(M)$ 
to $\mathrm{Aut}(Q)$ by $\Sigma(\alpha)(m,q) = (\alpha(m),q)$.
This the wrong thing to do: if we restrict attention to 
$\mathrm{SO}(\eta) < \mathrm{Diff}^{+}(M)$, we see that
sections of
$Q \times_{\widetilde{\mathrm{SO}}(\eta)}V$ transform under the trivial
representation of the Lorentz group. Our fermions are Lorentz scalars
rather than spin-$1/2$ particles.
This is a very real problem: 
using the wrong splitting will generally
result in an incorrect energy-momentum tensor \cite{GM}.

We conclude that not only the bundle $Q$ and
the covering map $u : Q \rightarrow OF_{g}$ are
relevant, but 
also the splitting 
$\sigma : \Gamma_{c}(TM) \rightarrow \Gamma_{c}(T\hat{Q})^{\tglnp}$.
It must satisfy $u_{*} \circ \sigma = D$
in order for the metric $g \in \Gamma(TF^{+}(M)/SO(\eta))$ to transform
properly. Luckily, 
such a $\sigma$ is naturally associated to any ordinary
spin structure $Q$. For generalized spin structures however,
this is no longer so.

\subsection{Gauge Fields and \texorpdfstring{$\spgtitel$-}{Spin }structures}

In the presence of gauge fields, the topological conditions 
on $M$ in order to support a spin structure are more relaxed.
Intuitively, this is because the gauge group $G$ can absorb some of the 
indeterminacy that stems from the 2:1 cover of the Lorentz group.\vspace{-0.7mm}

\subsubsection{Generalized Spin Structures}\vspace{-0.7mm}
This is made more rigorous by the notion of a 
generalized spin structure or 
$\spg$-structure \cite{AI}.
\begin{basdef} Let $G$ be a Lie group with a central subgroup isomorphic to $Z$. 
Let\footnote{This notation is convenient but slightly misleading.
Be ware that 
if $\eta$ is of signature $+++-$, then
$\mathrm{Spin}^{Z}$ is isomorphic to 
$\widetilde{\mathrm{SO}}(\eta)$, not to $\mathrm{Spin}(3,1)$. 
} $\spg := \widetilde{\mathrm{SO}}(\eta) \times_{Z} G$.
A $\spg$-structure is a
$\spg$-bundle $Q$ over $M$,
together with a map  $u : Q \rightarrow OF^{+}_{g}$
that makes\vspace{-1mm} 
\begin{center}
\begin{tikzpicture}
\pgfsetzvec{\pgfpoint{0.385cm}{-0.385cm}}
\node (lg) at (-1.2,0.7) {$\spg$};
\node (rg) at (1.2,0.7) {$\mathrm{SO}(\eta)$};
\node  (rpijl) at (1.2,0.3) {$\curvearrowleft$};
\node  (lpijl) at (-1.2,0.3) {$\curvearrowleft$};
\node (lb) at (-1.2,0) {$Q$};
\node (rb) at (1.2,0) {$OF^{+}_{g}$};
\node (lo) at (0,-1.2) {$M$};
\draw [->] (lg) --node[above]{$\kappa$} (rg);
\draw [->] (lb) --node[above]{$u$} (rb);
\draw [->] (lb) -- (lo);
\draw [->] (rb) -- (lo);
\end{tikzpicture}\vspace{-1mm}
\end{center}
commute. We again denote the map $(x , g) \mapsto \kappa(x)$ by $\kappa$.    
\end{basdef}
This gives rise to 
the principal 
$\widetilde{\mathrm{GL}}(\mathbb{R}^{n}) \times_{Z} G$-bundle
$\hat{Q} := Q \times_{\widetilde{\mathrm{SO}}(\eta)}
\widetilde{\mathrm{GL}}(\mathbb{R}^{n})
$. 
If $G = Z$, we recover the notion of a spin structure.
Apart from spin structures, 
the best known examples of $\spg$-structures are probably
the spin$^{c}$-structures, which are precisely
$\spg$-structures for $G = U(1)$.

Let $V$ be a representation of $\spg$.
The `physical bundle' is then the fibre
$J^{1}(\hat{Q})\times_{\spgtitel}V$.
A single section of $J^{1}(\hat{Q}) \times_{\spgtitel}V$
represents  
a metric,
a Levi Civita-connection,
a gauge field  
and a fermionic field.

The metric is the induced section of
$F^{+}(M)/\mathrm{SO}(\eta)$,
and the the Levi-Civita connection 
that of
$J^{1}(F^{+}(M))/\mathrm{GL}^{+}(\mathbb{R}^{n})$.
One constructs the 
principal $G/Z$-bundle 
$P := \hat{Q} / \smash{\tglnp}$, and the gauge field
is the induced equivariant connection on $P$,
a section of $J^{1}(P)/(G/Z)$.
The fermionic field is the induced section of
$\smash{\pi^{-1}(OF^{+}_{g}) 
\times_{\smash[b]{\widetilde{\mathrm{SO}}(\eta) \times_{Z}G}} V}$,
where one should note that the bundle
itself depends on $g$.

\subsubsection{Infinitesimally Natural \texorpdfstring{$\spgtitel$-}{Spin }Structures}

We argue that it only makes sense to consider $\spg$-structures 
which admit an appropriate transformation
law under infinitesimal diffeomorphisms of
space-time. We will call these $\spg$-structures
infinitesimally natural.

\begin{basdef}
A $\spg$-structure $u : Q \rightarrow FO_{g}$ will be called
`infinitesimally natural'
if 
$\hat{Q} \rightarrow M$ is infinitesimally natural as a principal fibre
bundle. Moreover, we require that
the splitting 
$\sigma : \Gamma_{c}(TM) \rightarrow 
\Gamma_{c}(T\hat{Q})^{\widetilde{\mathrm{GL}}(\mathbb{R}^{n})\times_{Z}G}$ 
of (\ref{rijtje}) satisfy $u_{*} \circ \sigma = D$.
\end{basdef}
The mathematical requirement that 
$\hat{Q} \rightarrow M$ be infinitesimally natural as a principal fibre
bundle
corresponds to the physical requirement that
fields should have a definite transformation behaviour under 
infinitesimal co-ordinate transformations.

The requirement $u_{*} \circ \sigma = D$ 
corresponds to the fact that we need to interpret a section of
$\hat{Q}/\widetilde{\mathrm{GL}}(\mathbb{R}^{n})\times_{Z}G 
\simeq F^{+}(M)/\mathrm{SO}(\eta)$ as a metric, and we know that 
its transformation behaviour is governed by $D$.

We view infinitesimally natural $\spg$-structures as the underlying principal
fibre bundles in any classical field theory combining gravity, 
fermions and gauge fields.
Let us work towards 
their classification.

\subsubsection{Classification of Infinitesimally Natural 
\texorpdfstring{$\spgtitel$-}{Spin }Structures}
Theorem \ref{grondsmakreefknol} is of course the main tool
when classifying infinitesimally natural $\spg$-structures. 
It provides a homomorphism
$
\rho :
G(k,M)
\rightarrow
\widetilde{\mathrm{GL}}(\mathbb{R}^{n})\times_{Z}G
$
such that
$\hat{Q}$ 
is isomorphic to
$
\tf
\times_{\rho} 
(\widetilde{\mathrm{GL}}(\mathbb{R}^{n})\times_{Z}G)
$, and a map $\exp\nabla : \tf \rightarrow \hat{Q}$.
The splitting $\sigma$ is induced by the lift of the $k^{\mathrm{th}}$
order derivative
$
\tilde{D} :
\Gamma(TM) \rightarrow 
\Gamma(T \tf) ^{G(k,M)}
$.
In summary, we have the following commutative diagram.\vspace{-0.4cm}
\begin{center} 
\begin{tikzpicture}
\pgfsetzvec{\pgfpoint{0.385cm}{-0.385cm}}
\node (linkser) at (0, 0) {$0$};
\node (links) at (2.3, 0) 
{$\Gamma_{c}(T\hat{Q})_{v}^{\widetilde{\mathrm{GL}}(\mathbb{R}^{n})\times_{Z}G}$};
\node (midden) at (6, 0) 
{$\Gamma_{c}(T\hat{Q})^{\widetilde{\mathrm{GL}}(\mathbb{R}^{n})\times_{Z}G}$};
\node (rechts) at (9.5,0) {$\Gamma_{c}(TM)$};
\node (rechtser) at (11, 0) {$0$};
\draw [->] (linkser) -- (links);
\draw [->] (links) -- (midden);
\draw [->] (midden) -- node[above] {$\pi_{*}$} (rechts);
\draw [->] (rechts)..controls(8,0.7)..node[above]{$\sigma$}(midden); 
\draw [->] (rechts) --  (rechtser);
\node (2linkser) at (0, -1.2) {$0$};
\node (2links) at (2.3, -1.2) 
{$\Gamma_{c}(TF(M))_{v}^{\mathrm{GL}(\mathbb{R}^{n})}$};
\node (2midden) at (6, -1.2) 
{$\Gamma_{c}(TF(M))^{\mathrm{GL}(\mathbb{R}^{n})}$};
\node (2rechts) at (9.5,-1.2) {$\Gamma_{c}(TM)$};
\node (2rechtser) at (11, -1.2) {$0$};
\draw [->] (2linkser) -- (2links);
\draw [->] (2links) -- (2midden);
\draw [->] (2midden) -- node[above] {$\pi_{*}$} (2rechts);
\draw [->] (2rechts)..controls(8,-0.5)..node[above]{$D$}(2midden); 
\draw [->] (2rechts) --  (2rechtser);
\node (3linkser) at (0, 1.2) {$0$};
\node (3links) at (2.3, 1.2) 
{$\Gamma_{c}(T\tf)_{v}^{G(k,m)}$};
\node (3midden) at (6, 1.2) 
{$\Gamma_{c}(T\tf)^{G(k,m)}$};
\node (3rechts) at (9.5,1.2) {$\Gamma_{c}(TM)$};
\node (3rechtser) at (11, 1.2) {$0$};
\draw [->] (3linkser) -- (3links);
\draw [->] (3links) -- (3midden);
\draw [->] (3midden) -- node[above] {$\pi_{*}$} (3rechts);
\draw [->] (3rechts)..controls(8,1.9)..node[above]{$\tilde{D}$}(3midden); 
\draw [->] (3rechts) --  (3rechtser);
\draw [->] (links) --node[left]{$u_{*}$} (2links);
\draw [->] (midden) --node[left]{$u_{*}$} (2midden);
\draw [->] (rechts) --node[right]{id} (2rechts);
\draw [->] (3links) --node[left]{$\exp\nabla_{*}$} (links);
\draw [->] (3midden) --node[left]{$\exp\nabla_{*}$} (midden);
\draw [->] (3rechts) --node[right]{id} (rechts);
\end{tikzpicture}
\end{center}

The classification theorem for  
infinitesimally natural $\spg$-structures
will take the following form. 
\label{naturalspinpage}
\begin{theorem}\label{naturalspin}
Let $M$ be an orientable smooth manifold of dimension $n \geq 2$. 
Let $G$ be a Lie group with a central 
subgroup $Z$ isomorphic to $\pi_{1}(\gln)$, and let $\mathrm{Lie}(G)$
be such that it does
not contain $\mathfrak{sl}(\mathbb{R}^{n})$. 
Finally, let $(Q,u)$ be an infinitesimally natural $\spg$-structure
over $M$.
Then $i_{*} : Z \rightarrow \pi_{1}(F^{+}(M))$
is injective, and
there exists a homomorphism
$
\tau : \pi_{1}(F^{+}(M)) \rightarrow G
$
which fixes $Z$, and makes $(Q,u)$ isomorphic to the bundle
$$
\widetilde{OF}^{+}_{g} \times_{\tau} G\,,
$$
with $u$ the natural projection map onto $OF^{+}_{g}$.
\end{theorem} 
{{\bf \noindent Remark} \quad We have already seen that an ordinary spin structure is always 
infinitesimally natural, so that proposition \ref{klopboor?} 
is just the special case $G = Z$.
}\medskip

\noindent Theorem \ref{naturalspin} rather simplifies the data needed to construct 
the relevant bundle $J^1(\hat{Q})\times_{\spgtitel}V$.
Indeed, it suffices to have a triple $(M,G,(\rho,V))$ of:\label{fluitekruid} 
\begin{itemize}
\item[-]An orientable manifold $M$ such that 
$i_* : \pi_{1}(\mathrm{SO}(\eta)) \rightarrow \pi_1(F(M))$ is injective.
Its image $Z= \{\pm 1\}$ is then central in $\pi_1(F(M))$.
\item[-]A representation $\rho$ of 
$\mathrm{\widetilde{SO}}(\eta)\times_{Z}\pi_1(F(M))$ on $V$ which is unitary
when restricted to $\pi_1(F(M))$, and faithful on $Z$.
\item[-]A subgroup $G < U(V)$ that commutes with the image of 
$\mathrm{\widetilde{SO}}(\eta)$ under $\rho$, and contains
the image of $\pi_1(F(M))$. 
\end{itemize}
We can then form $\hat{Q} = \tilde{F^{+}}(M) \times_{\pi_1(F)}G$, 
with structure group $\spg = \smash{\mathrm{\widetilde{SO}}(\eta) \times_{Z} G}$.
The relevant bundle is $J^1(\hat{Q}) \times_{\spgtitel} V$, 
a single section providing the metric, Levi-Civita connection, gauge field and
fermions. 

Note that the `gauge bundle' 
is the trivial bundle $P = M \times G/Z$, and that the bundle 
describing fermions and metric 
is $\smash{\widetilde{F}^{+}(M)} \times_{H} V$, where
$H$ is defined as  
$H := \widetilde{\mathrm{SO}}(\eta) \times_{Z} \pi_{1}(F)$.

According to theorem \ref{naturalspin}, the above setting exhausts the
possibilities-- at least under the (natural) assumptions that $Q$ is 
infinitesimally natural, and that $V$ is a faithful unitary representation
for $G$. 
Indeed, $G$ must then be compact, 
so that its Lie-algebra cannot allow $\mathfrak{sl}(\mathbb{R}^{n})$
as a subalgebra.

\subsubsection{Two Lemmas}

We set out to prove theorem \ref{naturalspin}. We start with 
two lemmas designed to explicate the homomorphism
$
\rho :
G(k,M)
\rightarrow
\widetilde{\mathrm{GL}}(\mathbb{R}^{n})\times_{Z}G
$.


\begin{lemma}
Denote the natural projection map
$\tf \rightarrow F^{+}(M)$ by $\pi$,
and write $\nu : G(k,M) \rightarrow \glnp$ for the corresponding homomorphism of groups.
Then each infinitesimally natural $\spg$-structure is isomorphic to one
for which 
$\kappa \circ \rho = \nu$
and
$u \circ \exp\nabla = \pi$.
\end{lemma}
\proof 
Denote $u \circ \exp\nabla$ by $\gamma$. 
Pick $\tilde{f}  \in \tfx$,
and let $\pi(\tilde{f}) = f$.
There exists a $c \in \glnp$
such that $\gamma(\tilde{f}) = f c$.
Equivariance of $\gamma$ then implies
$\gamma(\tilde{f} \tilde{y}) = f c \kappa(\rho (\tilde{y}))$
for all $\tilde{y} \in G(k,M)$.

We required $\gamma_{*} \circ \tilde{D} = D$, but 
clearly we also have $\pi_{*} \circ \tilde{D} = D$. 
Therefore,
if $\exp(\tilde{D}(v)) \tilde{f} = \tilde{f} \tilde{y}$
and 
$\exp({D}(v)) f = f y$,
we must have both 
$\nu(\tilde{y}) = y$
and
$
\kappa(\rho(\tilde{y})) = c^{-1} y c
$.

We may obtain all $\tilde{y}$ in $G(k,M)^{0}$,
the connected component of $G(k,M)$,
by choosing an appropriate $v$, so that
we have 
$
\kappa(\rho(\tilde{y})) = c^{-1} \nu(\tilde{y}) c
$
for all 
$\tilde{y} \in G(k,M)^{0}$.
Although $c$ depends on $\tilde{f}$ a priori,
it turns out to be constant up to scaling.
Indeed, as both $\kappa$ and $\rho$ are constant, 
so is the adjoint action of $c$.
Moreover, $\tilde{f} \mapsto c(\tilde{f})$ is invariant under
$G(k,M)^{0}$,
making it a function on 
$\tf / G(k,M)^{0} \simeq
\tilde{M}$.
All in all, we have established that there exist
$c_{0} \in \mathrm{SL}(\mathbb{R}^{n})$
and
$h : \tilde{M} \rightarrow \mathbb{R}^{+}$
such that
$c(\tilde{f}) = c_{0} h(\tilde{x})$, with $\tilde{x}$
the projection of $\tilde{f}$ to $\tilde{M}$.
We may write 
$u \circ \exp\nabla(\tilde{f}) = \pi(\tilde{f})c_{0}h(\tilde{x})$.


We show that we may as well take $c$ and $h$ to be $1$. 
Pick a $\tilde{c} \in \widetilde{\mathrm{SL}}(\mathbb{R}^{n})$ 
which covers $c$, and
construct the bundle $\smash{\hat{Q}_{c} := 
\hat{Q} \times_{\mathrm{Ad}(\tilde{c})}
\widetilde{\mathrm{GL}}(\mathbb{R}^{n})}$.
It is isomorphic to $\hat{Q}$,
with isomorphism $\hat{Q}_{c} \rightarrow \hat{Q}$ given by 
$[q,y] \mapsto q\tilde{c}^{-1}h^{-1}(\tilde{x})y$.

If we simply pull back the covering map on $Q$, we obtain
$u_c : \hat{Q}_{c} \rightarrow F^{+}(M)$
given by $[q,y] \mapsto u(q\tilde{c}^{-1}h^{-1}(\tilde{x})y)$.
This makes $Q_{c} := u_{c}^{-1}(OF_{g})$
into a $\spg$-structure 
isomorphic to $Q$, but with the desired properties.
\qed

\noindent Recall from (\ref{orient}) that 
$G(k,M) \simeq G(k,\mathbb{R}^{n}) \times_{Z} \pi_{1}(F^{+}(M))$,
and that we have $G(k,\mathbb{R}^{n}) \simeq \tglnp \ltimes G^{>1}$,
with $G^{>1}$ the subgroup of $k$-jets 
that are the identity to first order.
We unravel $\rho$, considering it as a map
$$
\tsln \times \mathbb{R}^{+} \ltimes G^{>1} 
\times_{Z} \pi_{1}(F^{+}(M))
\rightarrow
\tsln \times \mathbb{R}^{+} \times_{Z} G\,.
$$

\begin{lemma}
Under the assumptions of theorem \ref{naturalspin}, the map
$\rho$ 
is completely determined by a homomorphism
$\bar{\rho}_{1} : \pi_{1}(F^{+}(M)) \rightarrow  G$ fixing $Z$,
a homomorphism $\bar{\rho}_{2} : \pi_{1}(F^{+}(M)) \rightarrow \mathbb{R}^{+}$,
and an element $\Lambda$ of $\mathrm{Lie}(G)$ which
commutes with $\mathrm{Im}(\bar{\rho}_{1})$. 
We have
$$
\rho : \,\,(\tilde{x} , e^{t} ,g, [p]) \mapsto 
(\tilde{x} , e^{t} \bar{\rho}_{2}([p]), e^{t\Lambda} \bar{\rho}_{1}([p]))\,.
$$
\end{lemma}
\proof
Consider $\dot{\rho}$ as a Lie algebra homomorphism 
from
$\lsln \times \mathbb{R} \ltimes \mathfrak{g}^{>1}$
to
$\lsln \times \mathbb{R} \times \mathrm{Lie}(G)$,
and let $\dot{\rho}_{ij}$ be its $(i,j)$
component for $i,j \in \{1,2,3\}$. 

Because $\lsln$ is simple and not contained in $\mathrm{Lie}(G)$,
we must have $\dot{\rho}_{13} = 0$.
Due to the previous lemma, $\dot{\rho}_{12} = 0$ and 
$\dot{\rho}_{11} = \mathrm{id}$. As $\smash{\tsln}$ is simply connected, 
we must have $(\tilde{x},1,1,1) \mapsto (\tilde{x},1,1)$.
In particular, this forces the map $i_{*} : Z \rightarrow \pi_{1}(F^{+}(M))$
to be injective.

Again due to the previous lemma, $\dot{\rho}_{21} = 0$
and $\dot{\rho}_{22} = \mathrm{id}$. 
Define the `scaling element' $\Lambda := \dot{\rho}_{23}(1)$.
Then $(1,e^{t},1,1) \mapsto (1,e^{t},e^{t\Lambda})$.
The image of $\pi_{1}(F^{+}(M))$ must commute with 
$(\tilde{x} , e^{t} , e^{t\Lambda})$, so that 
$\bar{\rho} : \pi_{1}(F^{+}(M)) \rightarrow \mathbb{R^{+}} \times G$
is well defined, and $\Lambda$ commutes with its image.

We now show that 
$\dot{\rho}(\mathfrak{g}^{> 1}) = 0$.
First of all, as $[\lsln , \lsln + \mathfrak{g}^{> 1}]$
equals
$\lsln + \mathfrak{g}^{> 1}$ if $n$ is at least 2 
(see Lemma \ref{grondkusrelmuis}),
we must have
$$
\dot{\rho}(\mathfrak{g}^{>1}) \subset 
[\dot{\rho}(\lsln) , \dot{\rho}(\lsln + \mathfrak{g}^{>1})]
\subset \lsln \oplus 0 \oplus 0\,.
$$
But on the other hand, as 
$[\mathbb{R} , \mathfrak{g}^{>1}] = \mathfrak{g}^{>1}$,
(recall that $\mathbb{R}$ represents the Euler vector field),
we have 
$$
\dot{\rho}(\mathfrak{g}^{>1}) = 
[\dot{\rho}(\mathbb{R}), \dot{\rho}(\mathfrak{g}^{>1})]
\subset 0 \oplus 0 \oplus \mathrm{Lie}(G)\,.
$$
The intersection being zero, we have $\dot{\rho}(\mathfrak{g}^{>1}) = 0$.
But then $\rho(G^{>1}) = 1$, because $G^{>1}$ is simply connected.\qed

\subsubsection{Proof of Theorem \ref{naturalspin}}

Clearly, $\mathbb{R}^{+}$ cannot cause any topological obstruction, so we should
be able to eliminate both $\bar{\rho}_{2}$ and $\Lambda$ from the story.
The former is easy. recall that we seek $Q$, not $\hat{Q}$. 
As $Q$ is a subbundle of $\hat{Q}/\mathbb{R}^{+}$, we will focus on
the latter from now on, allowing us to simply disregard $\bar{\rho}_{2}$.

In order to remove $\Lambda$, we choose a 
volume form $\lambda$ on $M$.
This endows each frame $f \in F_{x}(M)$
with a volume $\mathrm{vol}_{\lambda}(f)$. 
A frame has volume 1 precisely when it is the jet of a diffeomorphism 
which preserves $\lambda$. Denote by $F^{\lambda}(M)$
the $\mathrm{SL}(\mathbb{R}^{n})$-bundle of frames with volume 1.
Its universal cover $\widetilde{F}^{\lambda}(M)$ 
is its inverse image under $\pi$, a principal
$\smash{\widetilde{\mathrm{SL}}}(\mathbb{R}^{n}) \times_{Z} \pi_{1}(F^{+}(M))$-bundle.

%
%

Define the isomorphism 
$$
\hat{Q}
\rightarrow 
\widetilde{F}_{\lambda}(M) \times_{\pi_{1}(F^{+}(M))} (G \times \mathbb{R}^{+})\vspace{-2mm}
$$
by 
$$
(\tilde{f} , g) \mapsto 
\Big(\tilde{f} \, \mathrm{vol}^{-1}_{\lambda}(\pi(\tilde{f})) \,,\,
\exp\big(\log(\mathrm{vol}_{\lambda}(\pi(\tilde{f}))) \,\Lambda \big) g 
, \mathrm{vol}_{\lambda}(\pi(\tilde{f}))\Big)
$$
where we consider $\hat{Q}$ as 
$(\tf \times_{\bar{\rho}_1 \bar{\rho}_2} G)$.
One can see that it is well defined, and that it intertwines the 
natural maps to $F/\mathbb{R}^{+}$.

This shows that $\hat{Q}/\mathbb{R}^{+}$, and therefore the spin structure 
$(Q,u)$, is completely determined 
by the homomorphism $\bar{\rho}_{1} : \pi_{1}(F^{+}(M)) \rightarrow G$.
We denote it by $\tau$ from now on.

If we choose $\lambda$ to be the volume form induced by the metric $g$,
then $\smash{OF^{+}_{g}}$ is a subbundle of $F_{\lambda}(M)$. Since 
$\pi_{1}(F^{+}(M)) = \pi_{1}(\smash[b]{OF_{g}^{+}})$, we have
$Q \simeq \smash[b]{\smash{\widetilde{OF}}_{g}^{+}} \times_{\tau} G$.
The spin map $u$ is simply the projection 
$\widetilde{OF}{}^{+}_{g} \rightarrow \smash[b]{OF^{+}_{g}}$.
The principal $G$-bundle $Q \rightarrow OF_{g}$ has a flat equivariant
connection induced by the one on \raisebox{-0.3mm}{$\widetilde{OF}{}^{+}_{g} \rightarrow OF^{+}_{g}$}.
This concludes the proof of theorem \ref{naturalspin}. \qed

\noindent Tracking back through the isomorphisms, 
we can formulate the following.
\begin{corollary}
Under the assumptions of theorem \ref{naturalspin},
there exists an element $\Lambda \in \mathrm{Lie}(G)$
which com\-mutes with the image of $\tau$,
such that the splitting
$\sigma : \Gamma(TM) \rightarrow
\Gamma(T\hat{Q}/\mathbb{R}^{+})^{\widetilde{\mathrm{SL}}
(\mathbb{R}^{n})\times_{Z}G}$
is given by 
$$
\sigma(v) = \nabla \circ \bar{D}(v) + \mathrm{Div}_{\lambda}(v) \Lambda\,,
$$
where $\lambda$ is the volume form induced by $g$,
$\bar{D}$ is the natural lift
from $\Gamma(TM)$ to $\Gamma(\widetilde{F}^{+}(M))/\mathbb{R}^{+})$,
and we have identified $\Lambda$ with the vector 
field on $\hat{Q}/\mathbb{R}^{+}$
induced by the action of the Lie algebra element.
\end{corollary}
It is clear that two different homomorphisms 
$\tau_1$ and $\tau_{2}$ yield isomorphic
$\spg$-structures if one is obtained from the other by conjugation
within $G$.

\subsection{Induced Unitary Representations 
\texorpdfstring{of $\smash{\tdiff}$}{}}
Before we explore the ramifications of theorem \ref{naturalspin}
on some theories of physical interest, we pause for a second 
to examine the representations of $\tdiff$ that live on 
infinitesimally natural bundles of the type described on page
\pageref{fluitekruid}. 

Recall that 
a metric is precisely a section of the bundle of symmetric spaces
$\pi : R(M) \rightarrow M$, with $R(M) := F^{+}(M)/\mathrm{SO}(\eta)$.
The bundle describing fermions and metric is 
$E(V) := \tilde{F}^{+}(M) \times_{\mathrm{Spin}^{Z}} V$, with $V$ 
a representation
of $\mathrm{Spin}^{Z} = \smash{\widetilde{\mathrm{SO}}}(\eta)\times_{Z}\pi_{1}(F(M))$
of spinor type.
The projection $E(V) \rightarrow R(M)$ makes it into a 
vector bundle over $R(M)$, and the projection $E(V) \rightarrow M$
makes it into a fibre bundle, unfortunately not a vector bundle, over $M$.

If we fix $M$ to be $\mathbb{R}^{n}$, and $g$ to be the Minkowski
metric $\eta$, then the action of the connected Poincar\'e group 
$\smash{\mathrm{SO}^{0}(\eta) \ltimes \mathbb{R}^{n}}$ 
on $\smash{OF^{+}_{\eta}(M)}$ induces an action (by automorphisms)
of its cover on $\smash[b]{\widetilde{OF}{}^{+}_{\eta}(M)}$. This in turn
induces
a representation on the sections of the associated vector bundle
$\smash{\widetilde{OF}{}^{+}_{\eta}(M)
\times_{\widetilde{\mathrm{SO}}{}^{0}(\eta)}} V$, the space of spinors.
The Minkowski metric gives rise to a unique (up to a scalar)
Poincar\'e-invariant volume form $\vol_{\eta}$ on $M$, 
which we use to turn the
nondegenerate Lorentz-invariant form on
$V$ into the nondegenerate Poincar\'e-invariant form
$\inp{\psi}{\psi'} = \int_{M}(\psi,\psi')\vol_{\eta}$
on the space of sections. The group action
commutes with the Dirac operator, and the representation
therefore decomposes according to mass.

We would like to do something similar for the diffeomorphism group
instead of the Poincar\'e group. 
The action of $\mathrm{Diff}{}^{0}(M)$ on $F^{+}(M)$ induces an action 
by automorphisms of 
$\tdiff$ on $\widetilde{F}{}^{+}(M)$, and therefore also on 
$R(M)$ and $E(V)$. Unfortunately, the sections $\Gamma(E(V) \rightarrow M)$
of physical relevance
do not form a representation, or even a vector space for that matter.

We do however have a representation of $\tdiff$ on 
$\Gamma(E(V) \rightarrow R(M))$. 
We will show that $R(M)$ carries a unique (up to a scalar)
$\mathrm{Diff}^{0}(M)$-invariant volume form\footnote{It 
would be tempting to interpret the volume form
$(\psi , \psi)\vol$ on $R(M)$ as the probability density of observing
a fermion at a certain point of space-time \emph{and}
with a certain value of the metric, but as far as I am aware,
this is not what physicists do. }
$\vol$,
so that $\Gamma(E(V) \rightarrow R(M))$ has a nondegenerate 
$\tdiff$-invariant form $\inp{\phi}{\psi} = 
\int_{R(M)}(\phi,\psi)\vol$.

\subsubsection{An Invariant Volume}
The only special asset of $\mathrm{SO}(\eta)$ that we will use 
is the fact that it is an open subgroup of the 
fixed point set 
of the involution
$g \mapsto \eta (g^{t})^{-1}\eta$ of $\mathrm{SL}(\mathbb{R}^{n})$, so 
we may as well work in that
setting.
\begin{proposition}
Let $M$ be a connected $n$-dimensional smooth manifold, and 
let
$H < \mathrm{SL}(\mathbb{R}^{n})$ be an open subgroup of the 
fixed point set of an involution.
Then the bundle $R(M) = F^{+}(M) / H$ of symmetric spaces 
admits a 
diffeomorphism-invariant volume form $\vol$.
It is everywhere nonzero, and unique up to $\mathbb{R}^{\times}$. 
%
\end{proposition}
\proof 
We write $\mathfrak{h}$ for the Lie algebra of $H$. The decomposition
$\lgln = \mathfrak{h} \oplus \mathfrak{m}$ of $\lgln$
into a positive and a negative eigenspace of the Lie algebra 
involution is
invariant under the adjoint action of $\mathfrak{h}$.
In particular, $\mathfrak{m}$ is an $\mathfrak{h}$-representation.
We construct a canonical map
$$
\wedge^{n}(\mathbb{R}^{n\, *}) \otimes 
\wedge^{k}(\lgln^*)^{\mathfrak{h}}
\stackrel{\sim}{\rightarrow} \Omega^{n + k}(R(M))^{\mathrm{Diff}^{0}(M)},
$$
with $k = \mathrm{dim}(\mathfrak{m})$.
Choose $\lambda \in \wedge^{n}(\mathbb{R}^{n\, *})$
and $\mu \in \wedge^{k}(\lgln^*)^{\mathfrak{h}}$. 
The former is a volume form on $\mathbb{R}^{n}$, 
the latter can be thought of as a volume form on 
$\lgln/\mathfrak{h} \simeq \mathfrak{m}$.
We construct the volume form on $R(M)$ over each point
$[f] \in R(M)$ separately, using the isomorphism 
$T_{f}F^{+}(M)/i_{f}(\mathfrak{h}) \stackrel{\sim}{\rightarrow}  T_{[f]}R(M)$
induced by the quotient map $F^{+}(M) \rightarrow R(M)$.

Each tangent space $T_{f}F^{+}(M)$ fits into 
the exact sequence of vector spaces
$$
0 \longrightarrow \lgln \stackrel{i_{f}}{\longrightarrow} T_{f}F^{+}(M)
\stackrel{\theta_{f}}{\longrightarrow} \mathbb{R}^{n} \longrightarrow 0\,,
$$
where $i_{f}(\xi) := \frac{d}{dt}|_{0} f \circ e^{t\xi}$, and
$\theta_{f} = f^{-1} \circ \pi_{*}$. 

The $n$-form $\theta_{f}^{*} \lambda$
annihilates the image of $i_{f}$, so in particular, it defines an 
$n$-form on on $T_{[f]}R(M)$.
This is in fact independent of the choice 
of $f$ representing $[f]$. 
The two identifications 
$T_{[f]}R(M) \simeq T_{f}F^{+}(M) / i_{f}(\mathfrak{h})$ and
$T_{[f]}R(M) \simeq T_{fh}F^{+}(M) / i_{fh}(\mathfrak{h})$
are intertwined by $R_{h\,*}$, the pushforward of right 
multiplication by $h$.
We have 
$
\theta_{fh} \circ R_{h\,*} = 
h^{-1} f^{-1} \pi_{*} R_{h \, *} = 
h^{-1} f^{-1} \pi_* =
h^{-1} \circ \theta_{f}$,
so that 
$(R_{h\,*})^{*} \theta_{fh}^{*} \lambda = \theta_{f}^* h^{-1\,*} \lambda$. 
Because $h^{-1\,*} \lambda = \lambda$, ($h \in H < \mathrm{SL}(\mathbb{R}^{n})$ 
is volume preserving),
we even have 
$(R_{h\,*})^{*} \theta_{fh}^{*} \lambda = \theta_{f}^*\lambda$. 
Since the two forms $\theta_{f}^*\lambda$
and $\theta_{fh}^{*} \lambda$ 
are related by $R_{h\,*}$, 
they
define the same $n$-form on
$T_{[f]}R(M)$. 

In order to transport $\mu$ to $T_{f}F^{+}(M)$, we need to choose 
a splitting $s_{f} : T_{f}F^{+}(M) \rightarrow \lgln$ such that
$s_{f} \circ i_{f} = \id$. Consider the $n+k$-form
$s_{f}^* \mu \wedge \theta_{f}^* \lambda$.
Although $s_{f}^* \mu$ depends on the choice of splitting, 
the form $s_{f}^* \mu \wedge \theta_{f}^* \lambda$
does not.
Indeed, since 
$(s_{f}-s'_{f}) \circ i_{f} = 0$ for every two 
sections $s_{f}$ and $s'_{f}$,
we have $i_{f}^{*} (s_{f}^{*} \mu - s'{}_{\!\!f}^{*}\mu) = 0$, which implies that
$(s_{f}^*\mu - s'{}_{\!\!f}^{*}\mu) \wedge \theta_{f}^* \lambda = 0$.

It is also independent of the choice of $f$. Since 
$
R_{h \,*} \circ i_{f}(\xi) = 
\smash{\frac{d}{dt}|_{0} f e^{t\xi} h} =
\frac{d}{dt}|_{0} f h (h^{-1}e^{t\xi} h) = 
i_{fh} \circ \mathrm{Ad}_{h^{-1}} (\xi)$,
we obtain 
$s_{f} = \mathrm{Ad}_{h} \circ s_{fh} \circ R_{h\,*}$, 
so that $s_{f}^* \mu = (R_{h\,*})^{*} s_{fh}^* {\mathrm{Ad}_{h}}^* \mu$.

Note that $\mathrm{Ad}_{h}$
respects 
the nondegenerate bilinear form
$\kappa(\xi, \chi) = \tr(\xi\chi)$ on $\lgln$, 
and hence its restriction to $\mathfrak{m}$.
Since
$\mathrm{Ad}_{h} : \mathfrak{m} \rightarrow \mathfrak{m}$
is orthogonal w.r.t.\ $\kappa$, it certainly preserves the
volume form $\mu$.
%
%

Now since $\mathrm{Ad}_{h}^* \mu = \mu$, we have
$s_{f}^* \mu = (R_{h\,*})^{*} s_{fh}^* \,\mu$, as required.
This means that $\vol_{[f]} = s_{f}^* \mu \wedge \theta_{f}^* \lambda$
is a well-defined $(n + k)$-form on $T_{[f]}R(M)$ for any $[f]$.
Since our construction varies smoothly with $[f]$,  
they combine into a smooth volume form $\vol$ on $R(M)$,
which clearly depends only on $\lambda \otimes \mu$.
 
A diffeomorphism $\alpha \in \mathrm{Diff}^{0}(M)$ induces an 
automorphism $\alpha_*$
of $F^{+}(M)$ by pushforward, and therefore also one of $R(M)$. 
The diffeomorphism invariance of the volume form
is a straightforward consequence of the diffeomorphism invariance 
of $f \mapsto \theta_{f}$ and $f \mapsto i_{f}$.
Indeed,
invariance of the volume, $(\alpha_*)^* \vol = \vol$, is equivalent to
$(\alpha_{**})^* (s_{\alpha_* f}^{*} \mu \wedge \theta_{\alpha_* f}^* \lambda)
= s_{f}^* \mu \wedge \theta_{f}^* \lambda$.
This is true because
$\theta_{\alpha_* f} \circ \alpha_{**} = 
\theta_{f}
$
and 
$s_{\alpha_* f} \circ \alpha_{**} = s_{f}$.
The first can be seen 
by expanding the definitions and using 
$\alpha^{-1} \circ \pi \circ \alpha_{*} = \pi$, 
the second is a consequence of 
$\alpha_{**} \circ i_{f}(\xi) = \frac{d}{dt}\alpha_{*}f e^{t\xi} 
= i_{\alpha_*f}(\xi)$.
%

If $M$ is connected, then $\mathrm{Diff}^{0}(M)$ 
acts transitively on $R(M)$, so that the invariant volume 
is determined by its value in a single point. It is therefore 
unique up to a real scalar. The fact that it is nonzero is clear from the
construction.
\qed

{{\bf \noindent Remark} \quad
The space $R(M)$ has more structure than just a canonical volume.
In fact, its 
frame bundle has a canonical $\mathrm{Diff}^{0}(M)$-invariant 
subbundle with structure group
$H \ltimes 
\mathrm{Lin}(\mathbb{R}^{n} ,\mathfrak{m})$,
with multiplication given by 
$(h,v)(h',v')= (hh', \mathrm{Ad}_{h}\circ v' + v \circ h' )$.  
This implies the existence of an invariant volume, because
the structure group is contained in
$\mathrm{SL(\mathbb{R}^{n} 
\times \mathfrak{m})}$.}\medskip

\noindent Now let $\widetilde{H}$ be the preimage of $H$ under 
$G(1,M) \rightarrow \mathrm{GL}^{+}(\mathbb{R}^{n})$, an extension
of $H$ by $\pi_{1}(F(M))$. (Cf.\ def.\ \ref{zwalksmurf}.)
If $V$ is a representation of $\widetilde{H}$, we construct 
the vector bundle $E(V) \rightarrow R(M)$ by 
$E(V) := \widetilde{F}{}^{+}(M) \times_{\widetilde{H}}V$.

The action of $\mathrm{Diff}{}^{0}(M)$ on $F^{+}(M)$ induces an action 
by automorphisms of 
$\smash{\tdiff}$ on $\smash{\widetilde{F}}{}^{+}(M)$, and therefore also on 
$R(M)$ and $E(V)$.
(It is just the integrated version of the natural lift of $\Gamma(TM)$.)
This makes $\Gamma(E(V))$, the space of smooth sections of 
$E(V) \rightarrow R(M)$,
into a representation of $\tdiff$. 

If $V$ has an $\widetilde{H}$-invariant 
(bilinear or sesquilinear)
form $(\,\cdot\,,\,\cdot\,)$, 
then this induces a form on each fibre of $E(V) \rightarrow R(M)$.
This defines a form 
$\inp{\phi}{\psi} = \int_{R(M)} (\phi,\psi) \vol$
on $\Gamma(E(V))$.
Because $\vol$ is everywhere nonzero, the form 
$\inp{\,\cdot \,}{\,\cdot \,}$ 
is positive definite and/or nondegenerate if and only if
$(\,\cdot\,,\,\cdot\,)$ is. 
Moreover, the form $\inp{\,\cdot \,}{\,\cdot \,}$ 
is $\tdiff$-invariant because
the volume is.
%
It is perhaps worth while to summarize 
this in a corollary.  

\begin{corollary}Let $H < \mathrm{SL}(\mathbb{R}^{n})$ be an open subgroup of 
the fixed point set of an
involution. Let
$\widetilde{H}$ be its preimage in $G(1,M)$, and
let $V$ be a re\-pre\-sen\-ta\-tion of $\widetilde{H}$
with invariant  
form $(\,\cdot\,,\,\cdot\,)$. 
Then 
the induced re\-pre\-sen\-ta\-tion $\Gamma(E(V))$ of $\tdiff$
has an invariant form
$\inp{\phi}{\psi} = \int_{R(M)} (\phi,\psi) \vol$, with $\vol$
the $\mathrm{Diff}^{0}$-invariant volume on $R(M)$. 
It is nondegenerate and/or positive definite if and only if 
$(\,\cdot\,,\,\cdot\,)$ is.
In particular, any unitary representation of $\widetilde{H}$ 
induces a unitary representation of $\tdiff$.
\end{corollary}
In the case that $V$ is a spinor type representation of 
$\widetilde{\mathrm{SO}}(\eta) \times_{Z}\pi_{1}(F)$,
the induced representation on $\Gamma(E(V))$
is only unitary if $\eta$ is positive definite.
On the other hand, an invariant sesquilinear form
exists for any signature.

\subsubsection{Decomposition into Irreducibles}

These induced representations will in general not be irreducible.
We sketch the decomposition into irreps, in the particular case
that $M$ is a parallelizable, compact, connected manifold
of dimension $\geq 2$, and $V$ is a unitary representation
of $\widetilde{H}$ that factors through $H$. (So that the 
induced representation will factor through $\mathrm{Diff}^{0}(M)$.)
The problem then essentially reduces to 
harmonic analysis on the symmetric space 
$\Delta = \mathrm{SL}(\mathbb{R}^{n})/H$.

We fix once and for all a global section $f : M \rightarrow F^{+}(M)$ of the
frame bundle. It identifies $F^{+}(M)$ with 
$M \times \mathrm{GL}^{+}(\mathbb{R}^{n})$, 
by way of $(x,g) \mapsto f_{x}g$.
If $\alpha_* : T_{x}M \rightarrow T_{\alpha(x)}M$
is the pushforward, then we define 
$D_{x}\alpha$ to be 
the matrix of $\alpha_{*}$ w.r.t.\ the frames at $x$ and
$\alpha(x)$, i.e.\
$D_{x}\alpha := f_{\alpha(x)}^{-1} \circ \alpha_* \circ f_{x}$.
The action of $\alpha$ on $M \times \mathrm{GL}^{+}(\mathbb{R}^{n})$
is then $\alpha : (x,g) \mapsto (\alpha(x), D_{x}\alpha \cdot g)$.

We identify $\mathbb{R}^{+} \times \mathrm{SL}(\mathbb{R}^{n})$
with $\mathrm{GL}^{+}(\mathbb{R}^{n})$ by 
$(\lambda,g) \mapsto \sqrt[n]{\lambda}\, g$. 
With $\Delta = \mathrm{SL}(\mathbb{R}^{n})/H$, we then have 
$R(M) = M \times \mathbb{R}^{+} \times \Delta$.
The action of $\alpha$ on $R(M)$ is  
$\alpha : (x,\lambda,\omega) \mapsto 
(\alpha(x) , \mathrm{det}(D_{x}\alpha) , \check{D}_{x}\alpha \cdot \omega)$,
with $\check{D}_{x}\alpha := D_{x}\alpha / \sqrt[n]{\mathrm{det}(D_{x}\alpha)}$.

If $\vol_{M}$ is the volume on $M$ defined by the frame, and $\vol_{\Delta}$ is
the unique (up to $\mathbb{R}^{\times}$) 
$\mathrm{SL}(\mathbb{R}^{n})$-invariant volume on $\Delta$, then 
it is not hard to see that the unique (up to $\mathbb{R}^{\times}$) 
$\mathrm{Diff}^{0}$-invariant volume
on $R(M)$ is $\vol = \vol_{M} \wedge d\lambda/\lambda^{2} \wedge \vol_{\Delta}$. 

The bundle $E(V)$ is simply $M \times \mathbb{R}^{+} \times E_{\Delta}(V)$,
with $E_{\Delta}(V)$ the vector bundle 
$\mathrm{SL}(\mathbb{R}^{n}) \times_{H}V$ over $\Delta$. Its space $\CH$ 
of $L^{2}$-sections 
constitutes a unitary representation
of $\mathrm{SL}(\mathbb{R}^{n})$, the induced representation of $V$.

We identify the space of $L^{2}$-sections of $E(V)$ with
$L^{2}(M , L^{2}(\mathbb{R}^{+})\otimes\CH)$,
where $\mathbb{R}^{+}$ is equipped with $d\lambda/\lambda^{2}$,
and $M$ with $\vol_{M}$.
The representation of $\mathrm{Diff}^{0}(M)$ is then given by
$$(\alpha \psi) (x) = \rho(\mathrm{det}(D_{x}\alpha^{-1}))^{-1} 
\otimes \pi (\check{D}_{x}\alpha^{-1})^{-1} \psi(\alpha^{-1}(x))\,,$$
with $\rho$ 
the nonunitary $\mathbb{R^{+}}$-representation 
$(\rho(c)f)(\lambda) = f(c^{-1}\lambda)$ on the Hilbert space 
$L^{2}(\mathbb{R}^{+}, d\lambda/\lambda^{2})$, and $\pi$ the 
unitary $\mathrm{SL}(\mathbb{R}^{n})$-representation on $\CH$.

We decompose $\rho$
into irreducible representations. 
The unitary transformation 
$W : L^{2}(\mathbb{R}^{+}, d\lambda/\lambda^{2})
\rightarrow 
L^{2}(\mathbb{R}, dt)$
defined by $(Wf)(t) = e^{-t/2} f(e^{t})$
takes $\rho$ into $\rho'(c) = W\rho W^{-1}(c)$, with
$(\rho' f')(t) = c^{-1/2}f'(t - \log(c))$.
The Fourier transform 
$\CF : L^{2}(\mathbb{R}, dt) \rightarrow L^{2}(\mathbb{R}, dk)$
then takes $\rho'$ into $\hat{\rho}(c) = \CF \rho'(c) \CF^{-1}$,
with $(\hat{\rho}(c) \hat{f})(k) = c^{-1/2 - ik} \hat{f}(k)$.
We conclude that the overall transformation 
$
\CF \circ W : 
L^{2}(\mathbb{R}^{+}, d\lambda/\lambda^{2})
\rightarrow 
L^{2}(\mathbb{R} ,dk),
$
given by
$$
\hat{f}(k) = \frac{1}{\sqrt{2\pi}} 
\int_{0}^{\infty} \lambda^{1/2 - ik} f(\lambda) d\lambda/\lambda^{2}\,,
$$
identifies $\rho$ with the representation 
$(\hat{\rho}(c) \hat{f})(k) = c^{-1/2 - ik} \hat{f}(k)$
on $L^{2}(\mathbb{R},dk)$. We write 
$L^{2}(\mathbb{R^{+}},d\lambda/\lambda^{2}) 
\simeq \int^{\oplus}_{\mathbb{R}} \mathbb{C}_{k}dk$ and
$\rho \simeq \int_{\mathbb{R}}^{\oplus} \rho_{k} dk$, with 
$(\rho_{k},\mathbb{C}_{k})$ the one-dimensional representation 
$\rho_{k}(c)z =c^{-\frac{1}{2} - ik} z$.
(See e.g.\ \cite[ch.\ II]{DixmierNEU} on direct integrals.)

Consequently, the $\mathrm{Diff}^{0}(M)$-representation
$L^{2}(M \times \mathbb{R}^{+} , \CH)$,  
with measure $\vol_{M}\wedge d\lambda/\lambda^2$, is unitarily
equivalent to
$L^{2}(M \times \mathbb{R}, \CH)$, with measure $\vol_{M}\wedge dk$ 
and
representation
\begin{equation}\label{llaatt}
(\alpha \psi)(x,k) =        
\mathrm{det}(D_{x}\alpha^{-1})^{1/2 + ik}
\pi(\check{D}_{x}\alpha^{-1})^{-1} 
\psi(\alpha^{-1}(x),k)\,.
\end{equation}
We thus have a decomposition
$
L^{2}(M , L^{2}(\mathbb{R}^{+})\otimes \CH) \simeq
\int_{\mathbb{R}}^{\oplus} 
L^{2}(M , \mathbb{C}_{k}\otimes \CH)dk\,
$ into subrepresentations,
where the unitary $\mathrm{Diff}^{0}(M)$-representation on the Hilbert space
$L^{2}(M , \mathbb{C}_{k}\otimes \CH)$ is given by
equation (\ref{llaatt}), but now with fixed $k$.
%

In general, 
the representations $L^{2}(M , \mathbb{C}_{k}\otimes \CH)$ will 
still be reducible. Now
suppose that we can decompose $\CH$ into unitary irreps as
$\CH = \int_{B}^{\oplus} \CH_{\beta} \mu(d\beta)$
for some measure space $(B,\mu)$, with
$\pi = \int_{B}^{\oplus}\pi_{\beta} \mu(d\beta)$.
Then 
\begin{equation}\label{irrepsdif}
L^{2}(M , L^{2}(\mathbb{R}^{+}) \otimes \CH)
\simeq 
\int^{\oplus}_{\mathbb{R} \times B}
L^{2}(M , \mathbb{C}_{k} \otimes \CH_{\beta}) dk \, \mu(d\beta)\,,
\end{equation}
with the representation on $L^{2}(M , \mathbb{C}_{k} \otimes \CH_{\beta})$
given by
\begin{equation}
(\alpha \psi)(x) =        
\mathrm{det}(D_{x}\alpha^{-1})^{1/2 + ik}
\pi_{\beta}(\check{D}_{x}\alpha^{-1})^{-1} 
\psi(\alpha^{-1}(x))\,.
\end{equation}

According to a theorem stated (but not proven) in
\cite{Ismagilov72}, these representations 
of $\mathrm{Diff}^{0}(M)$ are 
irreducible, at least under the assumption that
$M$ is a compact manifold of dimension at least 2.

This means that the problem of decomposing the $L^{2}$-closure of 
$\Gamma(E(V))$ into irreducible representations of 
$\mathrm{Diff}^{0}(M)$ reduces to that  
of decomposing the $L^{2}$-closure of 
$\Gamma(E_{\Delta}(V))$ into irreps of $\mathrm{SL}(\mathbb{R}^{n})$.
Although this is not an easy problem, it is one that has been studied 
in considerable detail.

\subsubsection{Example: The Hyperbolic Plane}
We explicitly carry out the decomposition in the special case that 
$M$ is the 2-torus $T^{2}$
with global frame $(\frac{d}{dx},\frac{d}{dy})$, the subgroup is 
$H$ is $\mathrm{SO}(\mathbb{R}^{2})$, and
$V = \mathbb{C}$ is the trivial representation. 

We identify the symmetric space
$\mathrm{SL}(\mathbb{R}^{2})/\mathrm{SO}(\mathbb{R}^{2})$
with the open unit disk 
$\Delta = \{z \in \mathbb{C}\,;\,|z|<1\}$.
It comes
equipped with the 
$\mathrm{SL}(\mathbb{R}^{2})$-action
$g \cdot z  =  A(g) z + B(g) / (\bar{B}(g)z + \bar{A}(g))$,
where $A(g)$ and $B(g)$ are defined by
\begin{equation}\label{sunaarsl}
\begin{pmatrix}A & B\\\bar{B} & \bar{A} \end{pmatrix} = 
\frac{1}{2}\begin{pmatrix}1 & -i\\1 & i \end{pmatrix}
\,g\,
\begin{pmatrix}1 & 1\\i & -i \end{pmatrix}
\,.
\end{equation}
(This is the identification of
$\mathrm{SU}(1,1)$ with $\mathrm{SL}(\mathbb{R}^{2})$ that comes
from the complex isomorphism of $\Delta$ with the upper half plane.)
The $\mathrm{SL}(\mathbb{R}^{2})$-invariant volume on $\Delta$ is 
$$\vol_{\Delta} = \frac{i dz \wedge d\bar{z}}{2(1-|z|^{2})^{2}}\,.$$
Since $E_{\Delta}(V)$ is the trivial bundle $\Delta \times \mathbb{C}$,
its
space of $L^{2}$-sections is $L^{2}(\Delta,\vol_{\Delta})$. 
It decomposes into the spherical principal series representations
$\pi_{\beta}$, which we will define shortly. 

The unit circle $S^{1} = \{\omega \in \mathbb{C}\,;\, |\omega| = 1\}$,
equipped with the rotation invariant measure
$d\omega/(2\pi i \omega)$,
carries an $\mathrm{SL}(\mathbb{R}^{2})$-action 
$\omega \mapsto g\cdot \omega$ because it is the boundary
of $\Delta$.
Define the
Fourier transform
$\CF : L^{2}(\Delta)\stackrel{\sim}{\rightarrow}L^{2}(\mathbb{R}^{+}\times S^{1})$
by
$$
(\CF f)(\beta,\omega) = \int_{\Delta}f(z) \left
(\frac{1-|z|^{2}}{|\omega - z|^{2}}\right)^{(1-i\beta)/2}\vol_{\Delta}\,.
$$
If we equip $\mathbb{R}^{+}$ with the measure
$\mu(d\beta) = (2\pi)^{-1} \beta \tanh(\half\pi\beta)d\beta$, 
then it 
is a unitary isomorphism \cite[p.\ 33]{HelgasonGGA}.
By direct calculation, one verifies that
$$\left(\CF \pi(g^{-1}) \CF^{-1} \hat{f}\right)(\beta,\omega) = 
|\bar{B}(g)\omega + \bar{A}(g)|^{i\beta - 1} \hat{f}(\beta,g\cdot \omega)\,,
$$ 
with $A(g)$ and $B(g)$ as in equation (\ref{sunaarsl}).
In other words, the Fourier transform identifies the 
representation of $\mathrm{SL}(\mathbb{R}^2)$ on $L^{2}(\Delta)$  
with the direct integral 
$\int^{\oplus}_{\mathbb{R}^{+}} \pi_{\beta} \mu(d\beta)$, where 
$\pi_{\beta}$ is the spherical principal series representation of
$\mathrm{SL}(\mathbb{R}^2)$
on $L^{2}(S^{1})$ defined by 
$\left(\pi_{\beta}(g^{-1})f\right)(\omega) = 
|\bar{B}(g)\omega + \bar{A}(g)|^{i\beta - 1} \, f( g\cdot \omega )\,.$

In this particular case, equation
(\ref{irrepsdif}) means that the $\mathrm{Diff}^{0}$-representation 
on the $L^{2}$-sections of
$E(V)$ is unitarily equivalent to 
$L^{2}(T^{2} \times S^{1} \times \mathbb{R}^{+} \times \mathbb{R}^{+})$,
with the measure
$d\phi \wedge d\theta \wedge d\omega/(2\pi i \omega) \wedge dk \wedge
(2\pi)^{-1} \beta \tanh(\pi\beta/2)d\beta
$, and the representation\vspace{-1.5mm}
$$
(\pi(\alpha)\psi)(x,\omega,k,\beta) = 
c_{k,\beta}(D_{x}\alpha^{-1})
\psi(\alpha^{-1}(x) , D_{x}\alpha^{-1} \cdot \omega , k, \beta)\,,\vspace{-0.5mm}
$$
with the cocycle ($A$ and $B$ as in (\ref{sunaarsl}))\vspace{-0.5mm}
$$
c_{k,\beta}(D_{x}\alpha^{-1}) = 
\mathrm{det}(D_{x}\alpha^{-1})^{1/2 + ik}
|\bar{B}(D_{x}\alpha^{-1})\omega + \bar{A}(D_{x}\alpha^{-1})|^{i\beta - 1}\,.\vspace{-0.5mm}
$$
(There is no need for the projection 
$\check{D}_{x}\alpha^{-1}$, as the scalars act trivially on $\Delta$.)

The unitary representations of $\mathrm{Diff}^{0}(T^{2})$ 
that one gets when fixing $k$ and $\beta$ are 
irreducible according to \cite{Ismagilov72}.


\subsection{Some Physical Theories} 
If we accept that any physically relevant 
$\spg$-structure must be
infinitesimally natural, and that $G$ should be 
the gauge group of the theory, then
theorem \ref{naturalspin}
provides a link between
the spectrum of elementary particles  
and the global topology of space-time.

The fact that the mere existence of a generalized spin structure
may place restrictions on the space-time manifold was recognized 
by Hawking and Pope \cite{HawPop}.   
Generalized spin structures were classified \cite{AI},
and it was found that 
if the Lie group contains $\mathrm{SU}(2)$, 
then `universal spin structures' exist 
\cite{BFF}, irrespective of the topology of $M$.
In particular, there are no topological obstructions to the existence of
a $\spg$-structure as soon as $\mathrm{SU}(2) < G$.

But according to theorem \ref{naturalspin}, this changes 
if one requires the $\spg$-structure to be infinitesimally natural. 
Universal spin structures then exist only for certain noncompact 
groups.
For compact $G$, the requirement that there exist a 
homomorphism $\pi_{1}(F^{+}(M)) \rightarrow G$ fixing $Z$
provides an obstruction on the space-time manifold $M$
in terms of the group of internal symmetries $G$. 
Let us see what this means in some specific cases.\vspace{-1mm}
\subsubsection{Weyl Spinors}\vspace{-1mm}
Consider a single massless charged Weyl spinor coupled to 
a $U(1)$ gauge field.
For simplicity, let us assume that $M$ is 4-dimensional, 
oriented, and time-oriented, so that $Z = \mathbb{Z}/2\mathbb{Z}$
and we may use $\mathrm{SL}(\mathbb{C}^{2})$ instead of 
$\widetilde{\mathrm{SO}}(\eta)$. 

This means that $G = U(1)$ and $V = \mathbb{C}^{2} \otimes \mathbb{C}_{q}$,
the two-dimensional
defining representation of 
$\mathrm{SL}(\mathbb{C}^{2})$ tensored with
the one dimensional defining representation 
of $U(1)$. 

Note that the 
representation descends to
$\mathrm{Spin}^{c} = \mathrm{SL}(\mathbb{C}^{2}) \times_{Z} U(1)$.
($Z$ is just $\pm\one$ in $U(1)$.)
This means that upon choosing a $\mathrm{Spin}^{c}$-structure $Q$,
the configuration space is
$\Gamma(J^{1}(\hat{Q})\times_{ \mathrm{ Spin }^{c} }V)$.

%
%
%

Let us now impose the requirement that $Q$ be infinitesimally natural.
Theorem \ref{naturalspin} tells us that there must then be 
a homomorphism $\pi_{1}(F^{+}(M)) \rightarrow U(1)$
sending $\mathbb{Z}/2\mathbb{Z}$ to $\pm 1$. 
If $\pi_{1}(M)$ is finitely generated,
then the image of $\pi_{1}(F^{+}(M))$ in $U(1)$ is a finitely generated
subgroup
containing $\pm\one$. 
It must be isomorphic to $\mathbb{Z}^{n} \times (\mathbb{Z}/2m\mathbb{Z})$
for some $n, m \in \mathbb{N}$. 
If we then send $\mathbb{Z}^{n}$ to $1$, we obtain a homomorphism
$\pi_{1}(F^{+}(M)) \rightarrow \mathbb{Z}/2m\mathbb{Z}$.
The sequence
$1 \rightarrow
\mathbb{Z} / 2\mathbb{Z} \rightarrow
\mathbb{Z}/2m\mathbb{Z} \rightarrow
\mathbb{Z}/m\mathbb{Z}
\rightarrow 1$
splits precisely when $m$ is odd.

We conclude that an infinitesimally natural $\mathrm{Spin}^{c}$-structure
exists on $M$ if and only if there is a surjective homomorphism
$\pi_{1}(F^{+}(M)) \rightarrow \mathbb{Z}/2m\mathbb{Z}$
which preserves $\mathbb{Z}/2\mathbb{Z}$. 
Only if $m$ is odd does this give rise to a spin
structure.\vspace{-1mm}

\subsubsection{Dirac Spinors}\vspace{-1mm}
Next, consider the case of a Dirac spinor. 
That is, $V = \mathbb{C}^{4} \otimes \mathbb{C}_{q}$,
with $\mathbb{C}^{4}$ the representation 
of $\mathrm{Cl}(\mathbb{C}^{4})$ which splits into
two identical irreps
$\mathbb{C}^{2} \oplus \mathbb{C}^{2}$
under $\mathrm{SL}(\mathbb{C}^{2})$, the left handed
and right handed spinors.

The fact that $V$ is reducible under $\mathrm{Spin}^{c}$
makes us re-examine our assumption that
the group $G$ in theorem \ref{naturalspin}
should be the gauge group $U(1)$.
Indeed, the unitary commutant of $\mathrm{SL}(\mathbb{C}^{2})$
in $V$ is $U(2)$ rather than $U(1)$. 
If we take any discrete subgroup $H < U(2)$ and form
the group $U(1)_{H}$ generated by $H$ and $U(1)$,  
can we take $Q$ to
be a $\spg$-structure with structure group $U(1)_{H}$?

As far as only the kinematics is concerned, the answer is yes.
The generic fibre of 
$J^{1}(\hat{Q}) \times_{\spgtitel} V$ is the same for 
$G = U(1)$ as it is for $G = U(1)_{H}$, so 
adding $H$ will not change the space of local sections.

But if we take into account the dynamics, the answer becomes:
`only if $H$ is a group of symmetries of the Lagrangian'.
The reason is that Lagrangians are usually defined 
in local co-ordinates,
yielding a local action $S_{U}$ for each co-ordinate patch $U \subset M$.
However, as the dynamics of the theory should be governed by a global 
action functional $S_{M}$,
it is necessary for $S_{U}$ and $S_{V}$ to agree on $U \cap V$.
This means that if $H$ is part of the structure group of the bundle,
then it must leave the Lagrangian invariant. 
If $H$ is a global symmetry, then the transition functions
must be constant. This is automatic if $H$ is discrete.

For instance, in the case of 
a massive Dirac fermion, the 
subgroup of $U(2)$ which preserves the Lagrangian is precisely the diagonal
$U(1)$.
This means that the relevant 
$\spg$-structures are precisely the $\mathrm{Spin}^{c}$-structures 
classified above.

For massless Dirac spinors, the left and right Weyl spinors decouple,
so that the relevant symmetry group is $U_{L}(1) \times U_{R}(1)$. 
Although the requirement on a manifold to carry 
a $\spg$-structure does not change, this does give us more 
$\spg$-structures for the same manifold.

This illustrates that we may enlarge the gauge group $G$
by any group of discrete symmetries of the Lagrangian in order
to obtain $\spg$-structures.
In particular, this means that infinitesimally natural $\mathrm{Spin}^{c}$-structures
are allowed even for uncharged Weyl spinors. (The image of 
$\pi_{1}(F^{+}(M))$ is automatically discrete.)\vspace{-1mm}

\subsubsection{The Standard Model}\vspace{-1mm}
In the the standard model
of elementary particle physics,
the gauge group $G$ is 
$(\mathrm{SU}(3) \times \mathrm{SU}(2)_{L} \times \mathrm{U}(1)_{Y})/ N$,
with $N$ the cyclic subgroup of order 6 generated by
$(e^{2\pi i/3}\one ,  -\one , e^{2\pi i/6})$. 
It is isomorphic to $S(U(3) \times U(2))$, a subgroup of $\mathrm{SU}(5)$,
and it has a unique central subgroup of order 2 
generated by $\mathrm{\textsl{diag}}\,(1,1,1,-1,-1)$.


The fermion representation $V$ for a single family 
can be conveniently described 
(see e.g.\ \cite{BaezHuerta})
by
$\mathbb{C}^{2} \otimes \wedge^{\bullet}\mathbb{C}^{5}$,
the tensor product of the defining representation of 
$\mathrm{SL}(\mathbb{C}^{2})$ and the exterior algebra 
of the defining representation of $\mathrm{SU}(5)$.
Under $\mathrm{SL}(\mathbb{C}^{2}) \times S(U(3) \times U(2))$,
this decomposes into 12 irreps corresponding to left and right handed 
electrons, neutrinos, up and down quarks and their antiparticles.

Unfortunately, $\mathrm{\textsl{diag}}\,(1,1,1,\!-1,\!-1) \in G$ acts
by $+1$ on right-handed fermions,
whereas $-\one \in \mathrm{SL}(\mathbb{C}^{2})$ acts by $-1$.
This means that $V$ does not define a representation of
$\mathrm{SL}(\mathbb{C}^{2}) \times_{Z} S(U(3) \times U(2))$
if one were to identify the central order 2 elements on both sides.

As the gauge group alone is of no use when trying to find
a $\spg$-structure,  
one has to involve the group 
of global symmetries
of the standard model Lagrangian.
It contains the gauge group $G$, but also (at least on the classical level)
the global $U(1)_{B} \times U(1)_{L}$-symmetries
that rotate quarks and leptons 
independently.
(These are connected to baryon and lepton number.)

We conclude that the only infinitesimally natural $\spg$-structures 
relevant to the standard 
model are the ones associated to  homomorphisms\vspace{-1mm}
\begin{equation} \label{homomo}
\pi_{1}(F^{+}(M)) \rightarrow 
\hat{G}\vspace{-1mm}
\end{equation}
that preserve $Z$, the subgroup of 
$U(1)_{B} \times U(1)_{L}$ generated by
$(-1,-1)$. 
In this expression, 
$\hat{G}$ is the group of global symmetries of the standard
model Lagrangian, which at least contains 
$S(U(3) \times U(2)) \times U(1)_{B} \times U(1)_{L}$.\vspace{-1mm}

\subsubsection{Infinitesimally Natural \texorpdfstring{$\spgtitel$-}{}Structures for 
the Standard Model}\vspace{-1mm}
It seems that
any manifold which possesses an infinitesimally natural 
$\spg$-structure for the standard model 
automatically permits an infinitesimally natural
$\mathrm{Spin}^{c}$-structure.
(This is at least the case if the group $\hat{G}$ of global symmetries 
is $S(U(3) \times U(2)) \times U(1)_{B} \times U(1)_{L}$.)
On the other hand, there do exist $\spg$-structures for the
standard model which are not $\mathrm{Spin}^{c}$.
We construct an example.

Consider de Sitter space 
$H = \{\vec{x} \in \mathbb{R}^{5} | 
- x^{2}_{0} + x^{2}_{1} + x^{2}_{2} + x^{2}_{3} + x^{2}_{4} = 1 \}$,
which has a pseudo-Riemannian metric $g$ with constant curvature 
induced by the Minkowski metric 
in the ambient $\mathbb{R}^{5}$.
Its group of orientation preserving isometries is $\mathrm{SO}(1,4)$,
and $H \simeq \mathrm{SO}(1,4)/ \mathrm{SO}(1,3)$.
Denote by $OF^{+\uparrow}_{g}(H)$ the bundle of orthogonal
frames with positive orientation and time-orientation.
By viewing $OF^{+\uparrow}_{g}(H)$ as a submanifold of 
$\mathbb{R}^{5} \times \mathrm{SO}(1,4)^{0}$, one can see that
$\mathrm{SO}(1,4)^{0}$ acts freely and transitively 
by $x : f \mapsto x_{*}f$. Therefore $OF^{+\uparrow}_{g}(H)$
is diffeomorphic to $\mathrm{SO}(1,4)^{0}$.

Now let $\Gamma < \mathrm{SO}(4)$ 
be a discrete group which acts freely, isometrically 
and properly discontinuously on $S^{3}$.
The manifold $\Gamma \backslash S^{3}$ is called 
a spherical space form. (See \cite{Wo} for a classification.)
As $\Gamma$ includes into $\mathrm{SO}(1,4)^{0}$, it acts
on $H$, making $M = \Gamma \backslash H$ into a pseudo-Riemannian manifold 
with constant curvature.

We immediately see that $\pi_{1}(M) = \Gamma$, 
because $H$ is simply connected.
We calculate the homotopy group of the frame bundle.
Because $OF_{g}^{+\uparrow}(M)$ is just 
$\Gamma \backslash OF_{g}^{+\uparrow}(H)$,
it is isomorphic to
$\Gamma \backslash \mathrm{SO}(1,4)^{0}$. 
Going to the universal cover, we see that
$OF_{g}^{+}(M) = \tilde{\Gamma} \backslash \smash{\widetilde{\mathrm{SO}}}(1,4)^{0}$.
As $\Gamma < \mathrm{SO}(4)$, we may consider $\tilde{\Gamma}$ to be 
the preimage of $\Gamma$ in $\mathrm{Spin}(4)$.
As the universal cover is simply connected, it is now clear that
$\pi_{1}(OF_{g}^{+\uparrow}(M)) = \tilde{\Gamma}$.
We get for free a homomorphism
$\smash{\tilde{\Gamma}} \rightarrow \mathrm{Spin}(4)
\simeq
\mathrm{SU}(2)_{l} \times \mathrm{SU}(2)_{r}$, which maps the 
noncontractible loop in the fibre to $(-\one, -\one)$.

Triggered by the WMAP-data on cosmic background radiation, 
there has been some interest in the case where $\Gamma$ is $I^{*}$,
the binary icosahedral group \cite{Lea}. 
Some evidence against \cite{Kea} as well as in favour of \cite{Rea}
this hypothesis
appears to have been found. 
We carefully avoid choosing sides 
in the debate, 
and simply point out that $I^{*} \backslash H$
gives rise to an interesting infinitesimally natural $\spg$-structure.

Under the identification
$\mathrm{Spin}(4) \simeq \mathrm{SU}(2)_{l} \times \mathrm{SU}(2)_{r}$,
we see that $\Gamma = I^{*} \times \one$ lives only in $\mathrm{SU}(2)_{l}$,
so that $\tilde{\Gamma}$ is the direct product of $I^{*} \times \one$
and the $\mathbb{Z}/2\mathbb{Z}$ generated by $(-\one, -\one)$.
One can therefore define a homomorphism (\ref{homomo})
by identifying $\mathrm{SU}(2)_{l}$ with 
$\mathrm{SU}(2)_{L}< G$, and mapping $(-\one, -\one)$ to 
$(-\one, -\one) \in U(1)_{B} \times U(1)_{L}$.
This yields an infinitesimally natural $\spg$-structure which uses the 
noncommutativity of the gauge group in an essential fashion. 
Note however that `ordinary'
$\mathrm{Spin}^{c}$-structures also exist on $M$.\vspace{-1mm}

\subsubsection{Extensions of the Standard Model}\vspace{-1mm}
The fact that $S(U(3) \times U(2))$ 
does not contribute to the obstruction of finding $\spg$-structures
on $M$
is due to the fact that
it never acts by $-1$ on $V$. This is not true for some GUT-type extensions of 
the standard model, such as the Pati-Salam
$\mathrm{SU}(2)_{L} \times \mathrm{SU}(2)_{R} \times \mathrm{SU}(4)$
model and anything which extends it, for example $\mathrm{Spin}(10)$.

If $N$ is the group of order 2 generated by $(-\one,-\one,-\one)$, then
infinitesimally natural $\spg$-structures in the Pati-Salam model 
correspond, neglecting global symmetries, to homomorphisms 
$\pi_{1}(F^{+}(M)) \rightarrow 
\mathrm{SU}(2)_{L} \times \mathrm{SU}(2)_{R} \times \mathrm{SU}(4)/N$
which take $Z$ to $\langle (-\one,-\one,\one)\rangle$.
This has the rather intriguing consequence that there may well
exist space-time manifolds $M$ which are compatible with the Pati-Salam model,
but not with the standard model.
A manifold $M$ would have this property
if 
the smallest quotient of $\pi_{1}(F(M))$ containing $Z$
is a nonabelian subgroup of 
$\mathrm{SU}(2)_{L} \times \mathrm{SU}(2)_{R} \times \mathrm{SU}(4)$
containing $(-\one,-\one,\one)$.

\section{Discussion}\label{acht}
A natural bundle is one in which diffeomorphisms of the base 
lift to automorphisms of the bundle in a local fashion.
In corollary \ref{PalaisTerng}, we have re\-derived the well known 
\cite{PT} 
result that any 
natural principal fibre bundle is associated to the $k^{\mathrm{th}}$
order frame bundle
$F^{k}(M)$. 

Bundles associated to spin structures are almost never natural, 
but nonetheless indispensable to physics.
We would therefore like to extend the notion of `naturality'.
One way to do this is to introduce
`gauge natural bundles' \cite{Eck} \cite{KMS}, 
which do not transform according to diffeomorphisms of the base $M$,
but according to automorphisms of a bundle over $M$ which has to be
specified.

We here propose to accommodate spin structures in a different fashion.
Although the Lorentz group does
not act on a spin structure, its associated Lie algebra does.
Analogously, we define an
`infinitesimally natural' bundle over $M$ to be a bundle
in which it is not the group of diffeomorphisms of $M$
that lifts, but its Lie algebra of vector fields.
%
This has the advantage that the link between the base and the 
fibres is not lost, so that a stress-energy-momentum tensor can be
constructed from Noether's principle \cite{No}. 

We assume only that the lift is a homomorphism of Lie algebras.
However,  
the careful analysis by Gotay and Marsden \cite{GM} 
reveals that in order to define a SEM-tensor from Noether's principle, 
one needs a lift that is given by a differential operator.
Assuming only that the lift is a homomorphism,
we prove in proposition \ref{alhetvet} that it is actually a differential operator .
This shows that 
one can construct a SEM-tensor from Noether's 
principle precisely when the bundle is infinitesimally natural.  
We therefore propose to describe fields on space-time 
by sections of infinitesimally natural bundles. 

We have classified the infinitesimally natural principal fibre bundles.
According 
theorem \ref{grondsmakreefknol},   
they are associated to the universal cover $\tf$ 
of the $k^{\mathrm{th}}$
order frame bundle, and their transformation behaviour is governed by
the disconnected group $G(k,M)$. 
This group depends on the manifold. It generalizes the 
$\mathrm{Pin}$ group in the sense that it
regulates parity, time reversal and Lorentz transformations. 

Theorem \ref{grondsmakreefknol} is
originally due to Lecomte \cite{Lecomte}.
Although our proof was obtained independently,
it does rather resemble Lecomte's.
We pause to highlight a few differences.
First of all, Lecomte does not work in the compactly supported
setting: he bases his version of lemma~\ref{zwabbernoot}
on results of Amemiya \cite{Amemiya} rather than on 
Pursell and Shanks' lemma~\ref{ideaal}.
Secondly, our uniform bound, in 
proposition~\ref{alhetvet}, is different
from Lecomte's.
But the most significant difference is probably that Lecomte bases
his proof on integration of foliations rather than on integration of
Lie algebroids, which probably makes our proof easier to generalize.
Strikingly enough, Lecomte's work does not appear to be motivated by physics;
he mentions spin structures only briefly. 



The consistent description of fermions in the presence of a 
gauge group $G$ requires 
a $\spg$-structure rather than a spin structure \cite{HawPop}, \cite{AI}. 
All spin structures are infinitesimally natural, but some
$\spg$-structures are not. 
As described in theorem \ref{naturalspin},  
infinitesimally natural $\spg$-structures for compact $G$
correspond to homomorphisms $\pi_{1}(F^{+}(M)) \rightarrow G$
that 
are injective on $\pi_{1}(\glnp)$.
 

The existence of infinitesimally natural $\spg$-structures therefore
provides 
an obstruction on space-time $M$ in terms of the symmetry group
$G$.
Some manifolds are even excluded for any (compact) choice of $G$. 
For example, $\mathbb{C}P^{2}$ does not admit
any infinitesimally natural $\spg$-structure 
because its frame bundle is simply connected. 
In our eyes, this disqualifies it as a model for space-time.


If all global symmetries are gauged, 
then theorem \ref{naturalspin} constitutes the
``connexion between the topology of space-time and the
spectrum of elementary particles'' alluded to in \cite{HawPop}.

\pagestyle{plain}
  \phantomsection
  	\addcontentsline{toc}{chapter}{Bibliography}

\newcommand{\etalchar}[1]{$^{#1}$}




\chapter*{Samenvatting} \label{ch:vatbaar}
\phantomsection
\addcontentsline{toc}{chapter}{Samenvatting}
\setsubdir{hoofdstukken/samenvatting/}
\selectlanguage{dutch}
{\frenchspacing 
Het meetproces in de kwantummechanica is een wat bizarre aangelegenheid. 
Iets wat in het dagelijks leven volstrekt vanzelfsprekend is, 
zoals het gelijktijdig meten van twee verschillende dingen, 
kan op kleine schaal opeens volstrekt onmogelijk blijken.

In hoofdstuk \ref{ch:ITISC} en \ref{ch:UDHP} is
een aantal van dit soort fundamentele beperkingen
zo scherp 
mogelijk geformuleerd. Je kunt immers moeilijk je grenzen opzoeken, 
als je niet weet waar die liggen. Vervolgens is voor een tweetal 
meetmethoden vastgesteld hoeveel verbetering er theoretisch nog mogelijk 
is. 
Het zijn allebei methoden om 
de inwendige toestand van 
een atoom vast te stellen aan de 
hand van het licht dat het uitstraalt.
De methode
van hoofdstuk \ref{ch:OEQS} blijkt 
optimaal te zijn, die van hoofdstuk \ref{ch:OPJM} laat 
nog een klein beetje ruimte voor verbetering.

Hoofdstuk \ref{ch:BLID} gaat over een ander deel van het onderzoek.
Alle velden die in de natuur voorkomen (elektrische velden, magnetische velden,
\ldots) 
bewegen mee met kleine vervormingen van de tijdruimte. Het belangrijkste 
resultaat in dit deel van het proefschrift is een classificatie van alle 
`meebewegende' velden die theoretisch mogelijk zijn. 
Verrassend genoeg zijn dat er in zekere 
zin tamelijk weinig, en bovendien hangen de mogelijkheden sterk af van 
de vorm van de tijdruimte.

In de rest van deze samenvatting wil ik graag 
\emph{zonder wiskundige formules}
een preciezer beeld geven van de inhoud van dit proefschrift. 
Ik beperk me daarbij tot 
de wiskundige idee\"en in
hoofdstuk \ref{ch:BLID}, omdat die het best uit te drukken zijn in woorden
en vooral ook plaatjes.



\subsection*{Elektrische Velden en de Raakbundel}

Een elektrisch veld is iets wat trekt aan geladen deeltjes. Hoe sterker het elektrisch veld, hoe harder 
de deeltjes in de richting van dat veld getrokken worden.  

\begin{center}
\begin{tabular}{c}
\includegraphics[width = 8 cm]{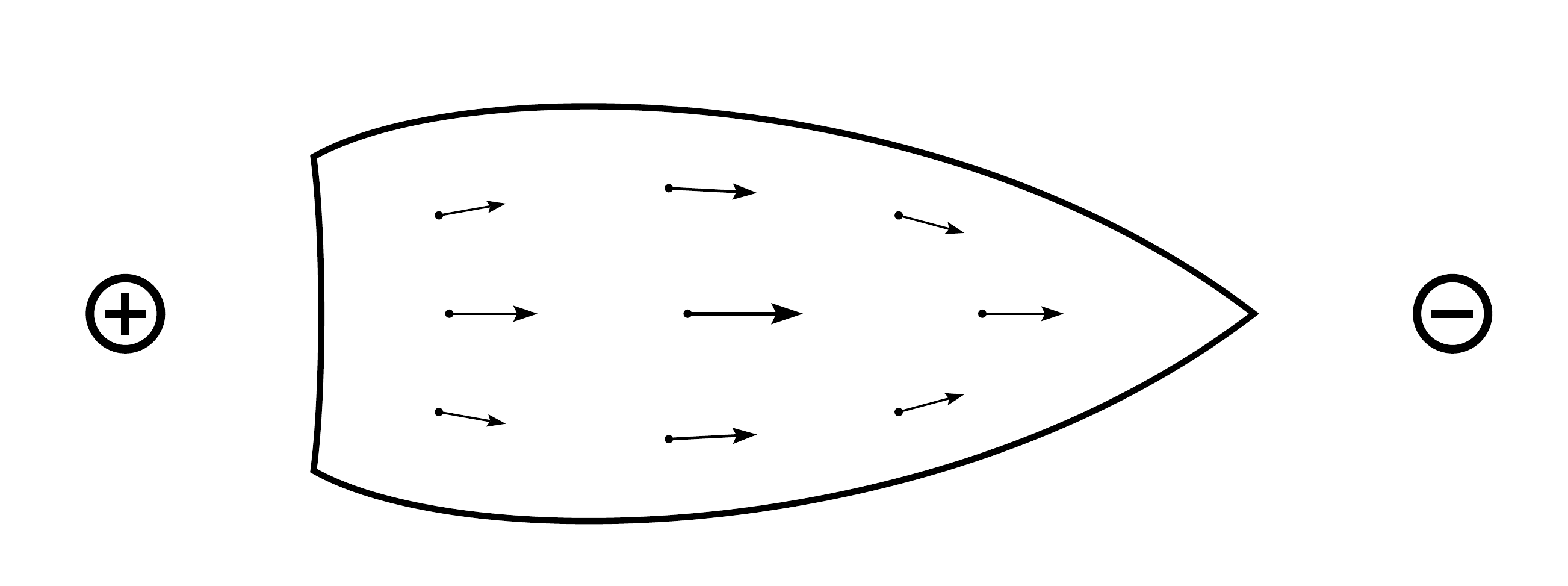}\\
{\small \quad \fig \label{strijkijzer}
Een elektrisch veld zet geladen deeltjes in beweging.}
\end{tabular}
\end{center}

Dit is bijzonder handig, bijvoorbeeld als je een broek wilt strijken. 
Je neemt dan een dunne metalen plaat, 
die vol zit met geladen deeltjes (elektronen). Plaats deze bovenop de te 
strijken broek.
Zorg vervolgens voor een elektrisch veld in de metalen plaat,  
bijvoorbeeld m.b.v.\ een stopcontact.
De geladen deeltjes worden in de richting van het veld door de plaat
getrokken, die door alle wrijving opwarmt.
Door vervolgens de plaat zachtjes heen en weer over de broek te bewegen 
kan men vouwen 
en kreukels verwijderen, of deze er juist in strijken 
(in het geval van de auteur vooral dat laatste).

\subsubsection{De Ruimte van Punten}
Laten we proberen om het elektrisch veld wiskundig te beschrijven. 
Stap \'e\'en is een beschrijving van de ruimte $M$ 
van punten waar het veld kan leven.
%
Dit vraagt om een zekere mate van algemeenheid.
Wat $M$ precies is, hangt immers af van de situatie. 

In een dunne metalen plaat kan ieder punt bijvoorbeeld
beschreven worden door 2 co\"ordinaten, lengte 
$(x)$ en breedte $(y)$. 
De ruimte $M$ heet dan 2-dimensionaal.
Maar elektrische velden leven ook
in de aarde, in de lucht en zelfs in de lege ruimte.
Je hebt daar
3 co\"ordinaten nodig om een 
punt te beschrijven: lengte $(x)$, breedte $(y)$ en hoogte $(z)$.
De ruimte $M$ is dus 3-dimensionaal.
Wonderlijk genoeg houdt het hier niet op. Volgens Einstein zijn tijd 
(1-dimensionaal) en ruimte (3-dimensionaal)
`van hetzelfde spul gemaakt'. Er zijn dus 4 co\"ordinaten nodig 
om een punt in de tijdruimte te beschrijven:
lengte $(x)$, breedte $(y)$, hoogte $(z)$ en tijd $(t)$.
Met andere woorden: 
de tijdruimte $M$ is 4-dimensionaal.

\begin{center}
\begin{tabular}{c}
\includegraphics[width = 11 cm]{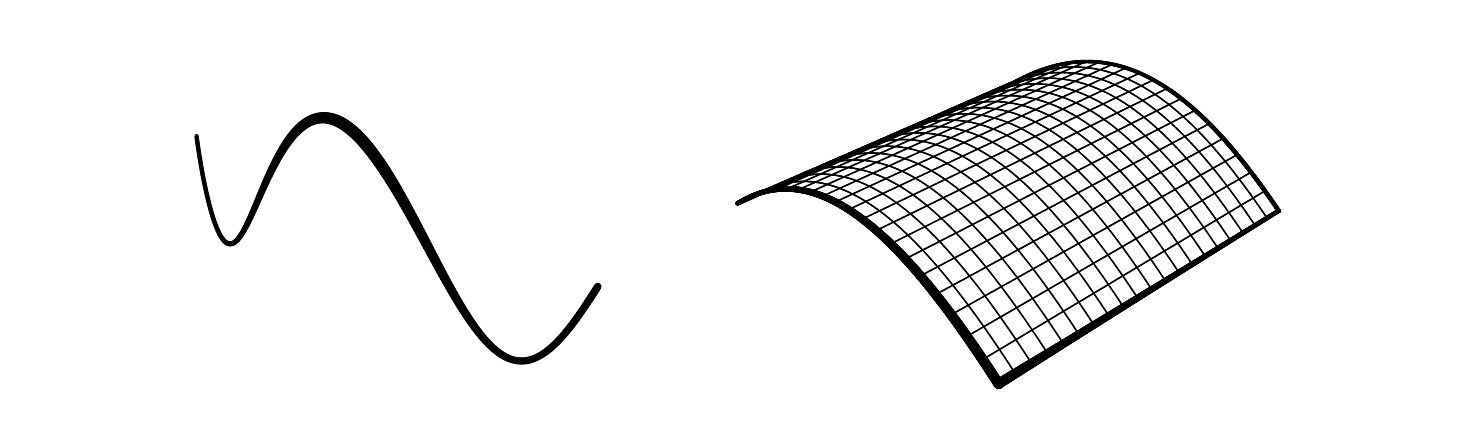}\\[5mm]
{\small \quad \fig \label{menigvoud}
Gladde vari\"{e}teiten van dimensie 1 (links) en 2 (rechts).}
\end{tabular}
\end{center}

Om dit alles in \'e\'en klap te vangen,
beschrijven we de ruimte $M$ van punten 
als een \emph{gladde vari\"eteit} 
(smooth manifold).
Gladde vari"eteiten heb je in elke gewenste dimensie.
Een gladde vari"eteit van dimensie 1 is niets anders dan een 
(kromme) lijn. Een gladde vari"eteit van dimensie 2 is een (gekromd)
oppervlak. 

Zoals gezegd is onze tijdruimte
$M$ een gladde vari"eteit van dimensie 4. Dit is lastig
te tekenen, zeker aangezien de tijdruimte ook nog eens gekromd is. 
Maar wiskundig gezien zijn 4 dimensies niet veel moeilijker dan 1 of 2. 
Bij woorden als `gladde vari\"eteit', `ruimte' of `tijdruimte' kan 
men dus gerust figuur \ref{menigvoud} in het achterhoofd houden.

\subsubsection{De Raakbundel}
Stap \'e\'en was het beschrijven van de ruimte $M$ van punten. 
Stap twee is het veld
in een vast punt\footnote{We 
schrijven kortweg $x$ in plaats van
$(x,y)$, $(x,y,z)$ of $(x,y,z,t)$. Op een heel proefschrift 
scheelt dat alras een bladzij.} 
$x$ van $M$.  
Het veld geeft aan hoe hard er aan een geladen 
deeltje wordt getrokken. Hoe sterker het veld, hoe groter de trekkracht.
Maar daarnaast bepaalt het veld ook in welke richting er wordt getrokken. 
Het veld in $x$ is dus iets wat naast grootte ook richting heeft: een pijl of \emph{vector}.
Hij moet natuurlijk raken aan $M$, omdat de deeltjes anders de ruimte zouden verlaten.
Kortom: het elektrisch veld in een punt $x$ beschrijven we wiskundig als een vector
die aangrijpt in $x$ en raakt aan $M$.

\begin{center}
\begin{tabular}{p{12cm}}
\hspace{0.5 cm}\includegraphics[width = 11 cm]{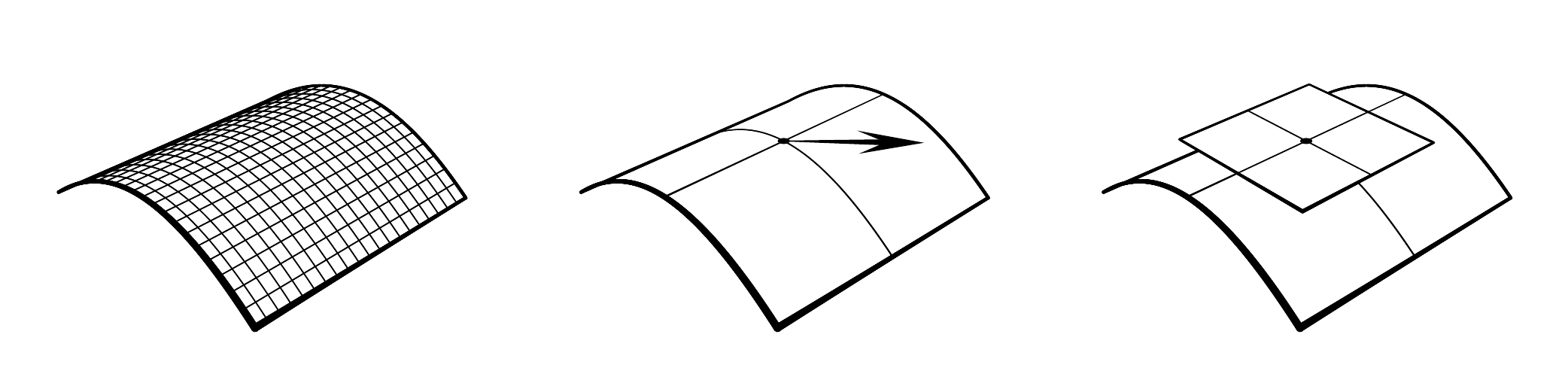}\\
{\small \quad \fig \label{raakzooi}
Links een gladde vari\"eteit $M$ van dimensie 2. In het midden een raakvector 
die aangrijpt in het punt $x$. 
Rechts het raakvlak in dit punt. Alle vectoren die 
raken aan $M$ en aangrijpen in $x$ 
hebben hun kop in dit vlak.}
\end{tabular}
\end{center}

Voor de derde en laatste stap in de beschrijving van het veld
moeten we ons realiseren dat
het elektrisch veld niet louter leeft in \'e\'en punt. 
Voor ieder mogelijk punt
in de ruimte bepaalt het hoe hard er aan geladen deeltjes getrokken wordt.
We beschrijven het veld dan ook als 
een hele familie vectoren, \'e\'en raakvector die aangrijpt in $x$ 
voor elk punt $x$ van de ruimte $M$.

\begin{center}
\begin{tabular}{ p{11 cm} }
\hspace{30 mm}\includegraphics[width = 5 cm]{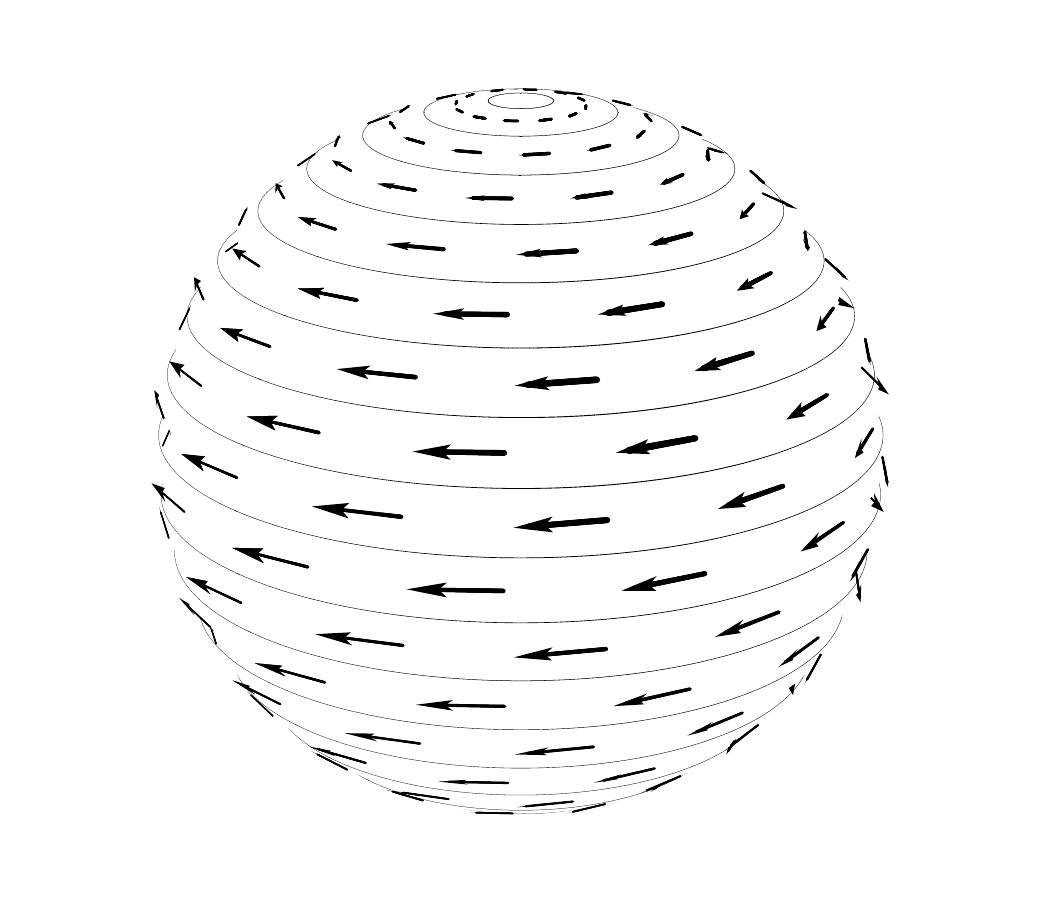}\\
{\small \quad \fig \label{harigebol}
De bolschil is een 2-dimensionale gladde vari\"{e}teit.
Het elektrisch veld hierop wordt beschreven door
een snede van zijn raakbundel. 
Dit is iets wat in ieder punt $x$ \'{e}\'{e}n raakvector uitkiest: 
de waarde van het veld in $x$.}
\end{tabular}
\end{center}

Alles bij elkaar geeft dit een wiskundige beschrijving van het elektrisch veld.
De ruimte waarop het veld leeft, beschrijven we met een gladde vari\"eteit $M$.
Bij ieder punt $x$ van $M$ hoort een \emph{raakvlak} (tangent space) $T_{x}M$, 
de ruimte van alle vectoren die raken aan $M$ en aangrijpen in $x$ 
(zie figuur \ref{raakzooi}).
Het raakvlak in $x$ beschrijft precies alle mogelijke waarden die het elektrisch veld 
in $x$ kan aannemen.
We defini\"{e}ren de \emph{raakbundel} (tangent bundle)
als de familie van alle mogelijke raakvlakken $T_{x}M$, waar
$x$ loopt over de punten van $M$.

We beschrijven het elektrisch veld als een \emph{snede} van de raakbundel.
Dit is iets wat in ieder raakvlak $T_{x}M$ precies 
\'{e}\'{e}n vector uitkiest: de waarde van het veld in $x$.


\subsection*{Algemene Velden en Vezelbundels}

Het elektrisch veld is niet het enige veld.
Zo bestaat er bijvoorbeeld een drukveld, dat als het ware bijhoudt in hoeverre
een rubber balletje 
zou worden samengedrukt als het de moed had op een bepaalde plek in de ruimte
zijn bol gezicht te vertonen.
Een trilling in het drukveld ervaren wij als geluid.
%
Ook het elektrisch veld kan trillen:
dit doet het samen met zijn tweelingbroer, het magnetisch veld.
Rond 1864 ontdekte J.C.~Maxwell\footnote{Hij publiceerde zijn bevindingen in 
"A Dynamical Theory of the Electromagnetic Field", 
\emph{Phil.\,Trans.\,R.\,Soc.\,Lond.}
{\bf 155}, blz. 459--512, 1865. 
Maxwell is veel geprezen voor zijn `ontdekking van het licht', iets wat voor 
een man met die initialen nochtans geen grote opgaaf zal zijn geweest.} dat
licht precies bestaat uit dit soort gecombineerde  
elektromagnetisch trillingen.  


Maar ook in de moderne natuurkunde spelen velden een hoofdrol.
Het hele idee van een `fundamenteel deeltje' is daar zelfs
vervangen door dat van een `fundamenteel veld'!
Deeltjes zijn dan niets anders dan trillingen in dit veld.

Het loont dus de moeite om de voorgaande bespiegelingen iets algemener 
te trekken,
en een wiskundige beschrijving te geven voor velden in het algemeen.
Dit hebben we gedaan in hoofdstuk \ref{bunklasveld}.
Een sleutelrol is daarbij weggelegd voor het begrip \emph{vezelbundel}
(fibre bundle), dat de raakbundel generaliseert.

\begin{center}
\begin{tabular}{p{11 cm}}
\hspace{2 cm} \includegraphics[width = 6 cm]{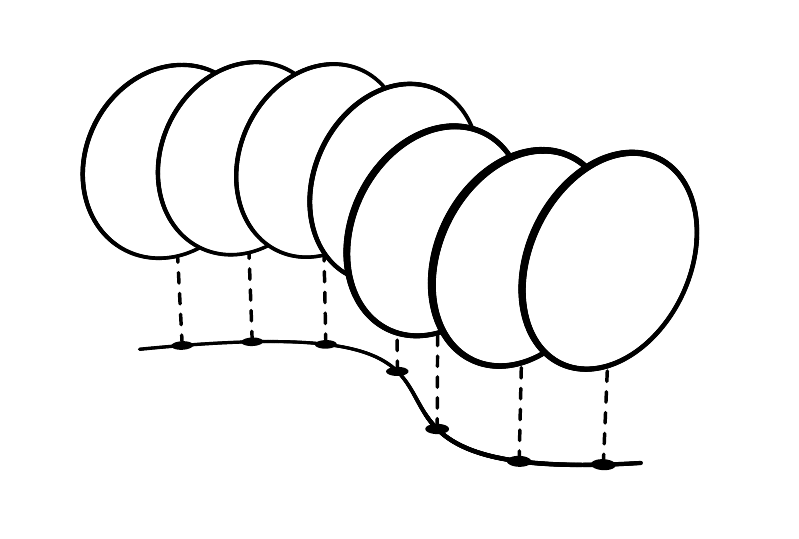}\\[-1mm]
{\small \quad \fig \label{bundleofjoy}
Schets van een vezelbundel. 
De basisruimte $M$ is de kromme beneden.
Bij ieder punt $x$ op de kromme hoort een vezel $F_{x}$.
Alle vezels hebben hier de vorm van een schijf.}
\end{tabular}
\end{center}

Een vezelbundel over een gladde vari\"eteit $M$ is per definitie een familie 
van \emph{vezels}, \'e\'en vezel $F_{x}$ bij ieder punt $x$ van $M$.
De vezel bij het punt $x$ 
stelt de verzameling voor van waarden 
die het veld mag aannemen in $x$.
Die vezels mogen van alles zijn, maar we eisen wel dat 
alle vezels 
er `hetzelfde uitzien' (in jargon: isomorf zijn).

Een belangrijk voorbeeld van een vezelbundel is de raakbundel. 
De vezel $F_{x}$ bij een punt $x$ wordt hier gevormd door 
het raakvlak $T_{x}M$, de ruimte van alle 
vectoren die raken aan $M$
en aangrijpen in het punt $x$. 
Alle vezels zien er hetzelfde uit, in dit geval als vlakken.


In het algemeen beschrijven we een veld als een \emph{snede} (section) 
van de vezelbundel.
Een snede is iets wat in iedere vezel $F_{x}$ precies \'e\'en
waarde uitkiest.\vspace{-2mm}

\begin{center}
\begin{tabular}{p{11 cm}}
\hspace{2 cm}\includegraphics[width = 7 cm]{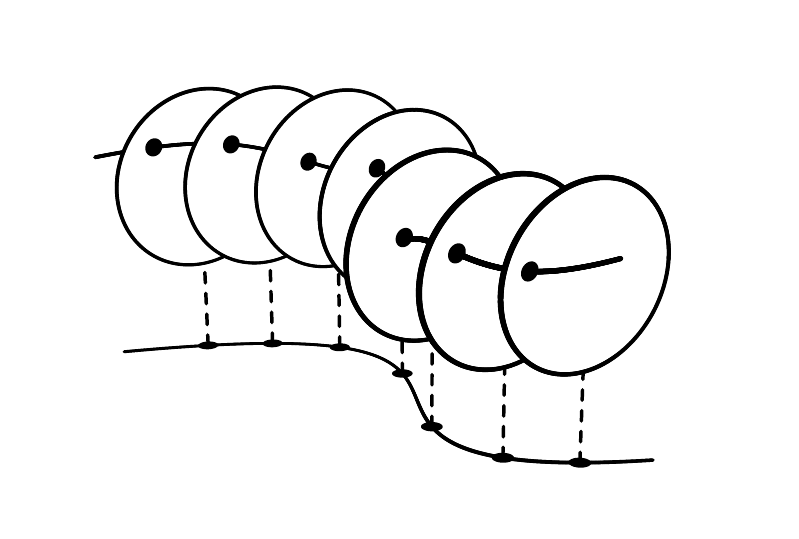}\\[-3mm]
{\small \quad \fig \label{Wesley}
Schets van een snede. In iedere vezel is \'e\'en waarde
aangewezen.}
\end{tabular}
\end{center}

De wiskundige beschrijving van een veld werkt dus als volgt. 
De ruimte waarop het veld leeft, beschrijven we met een gladde vari\"eteit $M$.
Het soort veld (bijvoorbeeld: het elektrisch veld) wordt bepaald door 
een vezelbundel (bijvoorbeeld: de raakbundel).
Bij ieder punt $x$ van $M$ hoort dan een vezel $F_{x}$ van waarden die het veld 
daar aan kan nemen (voor de raakbundel is dit het raakvlak $T_{x}M$).
Het veld zelf wordt beschreven door een snede van de vezelbundel.
Die wijst in iedere vezel $F_{x}$ \'e\'en waarde aan:
de waarde van het veld in $x$. (Een snede van de raakbundel wijst 
bij ieder punt dus \'e\'en raakvector aan.)

\subsubsection*{Natuurlijke Bundels}
Alle vezelbundels zijn gelijk, maar sommige vezelbundels zijn gelijker dan
andere. Een vezelbundel heet \emph{natuurlijk} als iedere vervorming
(jargon: diffeomorfisme) van de ruimte $M$ een welbepaald 
effect heeft op zijn vezels.
(Zie definitie \ref{defvannatbun} op blz.~\pageref{kroep} voor een 
preciezere formulering.)
Je kunt de natuurlijke bundels zien als een soort `elite' onder 
de vezelbundels.

\begin{center}
\begin{tabular}{c}
\includegraphics[width = 10 cm]{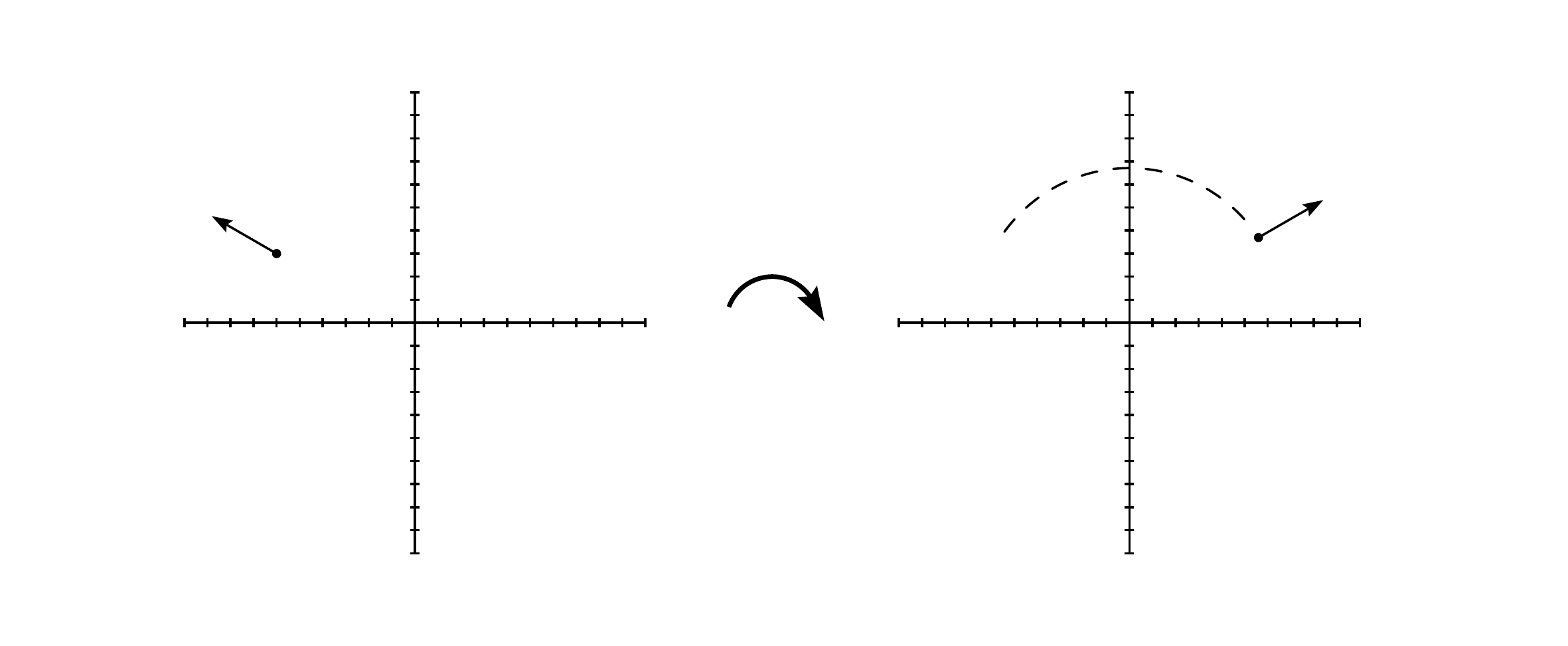}\\[-2mm]
{\small \quad \fig \label{figdraai}
Een raakvector draait mee als je het vlak draait.}
\end{tabular}
\end{center}

Elke raakbundel is natuurlijk.
%
Bekijk bijvoorbeeld 
het platte vlak $M$, met een raakvector
die aangrijpt in een punt $x$ van $M$. 
Vervolgens draaien we het vlak. Wat gebeurt er met de vector?
Het antwoord is duidelijk: \emph{de vector beweegt mee
als je $M$ roteert}, zie figuur \ref{figdraai}.
De raakbundel heet `natuurlijk'
omdat elke vervorming 
van de ruimte $M$ een welbepaald effect heeft op de raakvectoren.

Dit heeft belangrijke gevolgen voor het elektrisch veld, dat wordt
beschreven door een snede van de raakbundel. 
Bij een vervorming van de ruimte $M$ is er een welbepaald effect
op de raakvectoren in \emph{ieder} afzonderlijk punt van $M$.
Omdat een elektrisch veld niets anders is dan een manier om
aan \emph{alle} punten van $M$ een raakvector in dat punt toe te kennen,
moeten we concluderen dat \emph{het gehele veld meeverandert}
als $M$ wordt vervormd.
Zie figuur~\ref{figvelddraai}. 

\begin{center}
\begin{tabular}{c}
\includegraphics[width = 10 cm]{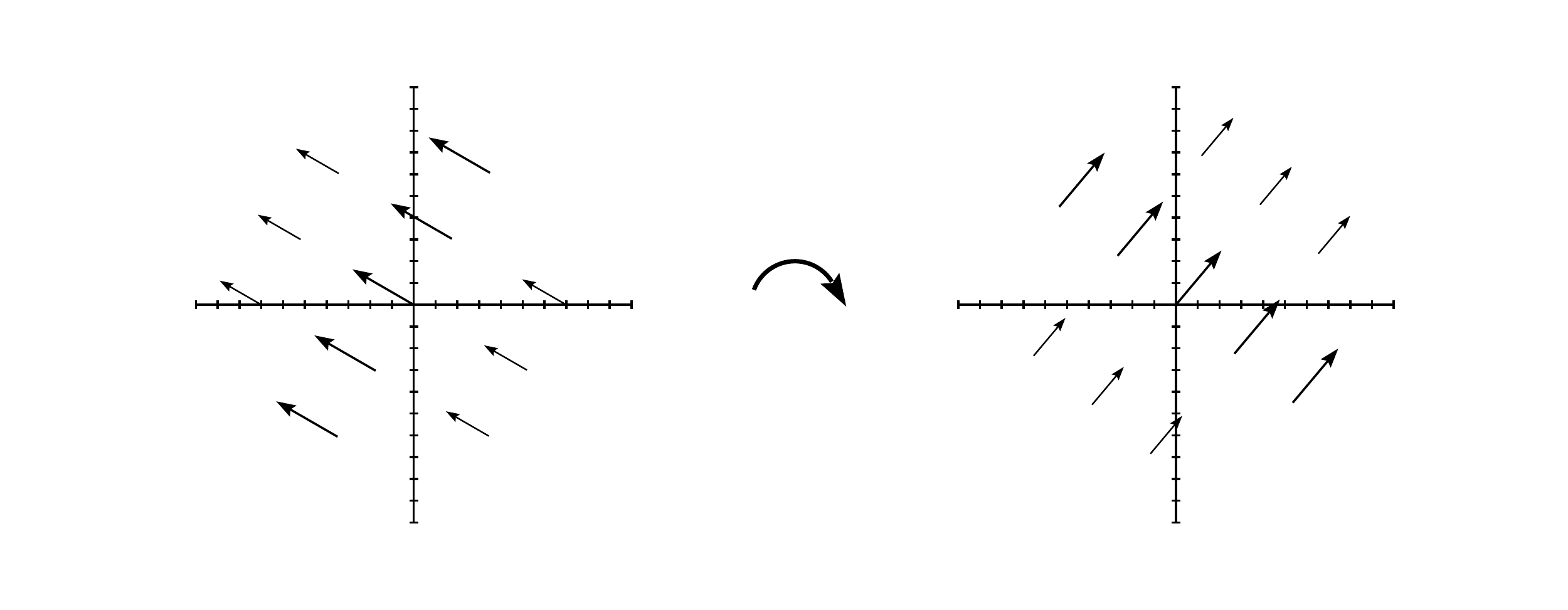}\\[-2mm]
{\small \quad \fig \label{figvelddraai}
Ook het veld draait mee als je het vlak draait.}
\end{tabular}
\end{center}

Aan het einde van hoofdstuk \ref{bunklasveld} wordt uitgelegd
dat het meebewegen van het veld bij kleine vervormingen van $M$
noodzakelijk is om
de \emph{energie} 
(jargon: de `stress-energy-momentum tensor')
van een veld te defini\"eren.

Het feit dat
elektrische en magnetische velden energie kunnen bevatten,
hoeft niet te verbazen.
Licht is immers niets anders dan een golf in het elektromagnetisch veld.
Iedereen die wel eens roodverbrand onder de zonnebank vandaan is gekomen, 
heeft het energiedragend vermogen van het elektromagnetisch veld 
dus experimenteel vastgesteld, om niet te zeggen aan den lijve ondervonden. 

Energietransport door het elektromagnetisch veld heeft trouwens ook
bijverschijnselen die prettiger zijn dan eerstegraads brandwonden. 
Leven op aarde, om eens wat te noemen.
De energie die het leven hier mogelijk maakt, komt direct 
of indirect van de zon. 
Maar als het elektromagnetisch veld geen energie zou dragen,
dan zou al die 
zonne-energie vrolijk op de zon achterblijven. 
De aarde zou dan koud en levenloos achterblijven.



Enfin, elektrische velden kunnen dus, samen met magnetische velden,
energie overdragen. Dit komt doordat de raakbundel natuurlijk is.  
Hoe zit dat met andere velden? 

We weten uit de natuurkunde dat velden eigenlijk altijd
in staat zijn energie te dragen.
Het ligt dus voor de hand om velden te beschrijven met
\emph{natuurlijke} vezelbundels.
Het eerste voorbeeld van een natuurlijke vezelbundel 
hebben we dus al gezien: de raakbundel.

Natuurlijke bundels zijn witte raven onder de vezelbundels. 
Een algemene vezelbundel is zomaar een familie vezels $F_x$, 
zonder enige verdere structuur.
In het algemeen heb je geen flauw idee hoe de vezels 
zouden moeten meeveranderen als de ruimte $M$ vervormt,
en is het zelfs nog maar de vraag of dit hoe dan ook wel op een 
consistente manier kan. Als een vezelbundel daarentegen 
natuurlijk is, heeft
elke vezel $F_x$ als het ware een band met zijn bijbehorend punt $x$.
Als de ruimte $M$ vervormt, `voelen' de vezels $F_x$ hoe hun 
eigen punt beweegt, en bewegen mee, vgl.~figuur \ref{figdraai}
en \ref{figvelddraai}.  
 

%


\subsubsection*{Infinitesimaal Natuurlijke Bundels}

Een bundel verdient dus de eretitel `natuurlijk' als 
elke vervorming van de ruimte $M$ een welbepaald effect heeft op de vezels $F_{x}$.
Die eis is soms een beetje al te zwaar.
We noemen een bundel daarom \emph{infinitesimaal natuurlijk} als elke
\emph{kleine} vervorming van $M$ een welbepaald effect heeft op de vezels
$F_{x}$. (Zie ook definitie \ref{defvaninnatbun} op blz.~\pageref{kroepstrup}.)
Dat is een soort troostprijs. Als je een bundel bent, is 
infinitesimaal natuurlijk zijn 
niet zo goed als \'echt natuurlijk zijn, 
maar het is toch al heel wat.  
Voor het defini\"eren van energie zijn bijvoorbeeld
alleen \emph{kleine} vervormingen nodig.
Met andere woorden, een veld kan energie dragen precies dan als
hij hoort bij een infinitesimaal
natuurlijke bundel.


De zogenaamde \emph{spinorbundels} zijn
infinitesimaal natuurlijke vezelbundels bij uitstek.
De spinorbundel bij het platte vlak $M$ is precies de
raakbundel, behalve dat vectoren \emph{op halve snelheid} met de ruimte 
meedraaien, zie figuur~\ref{figspindraai}.
Zo een halfslachtige vector heet ook wel een \emph{spinor}.\vspace{-0.5cm} 

\begin{center}
\begin{tabular}{c}
\includegraphics[width = 11 cm]{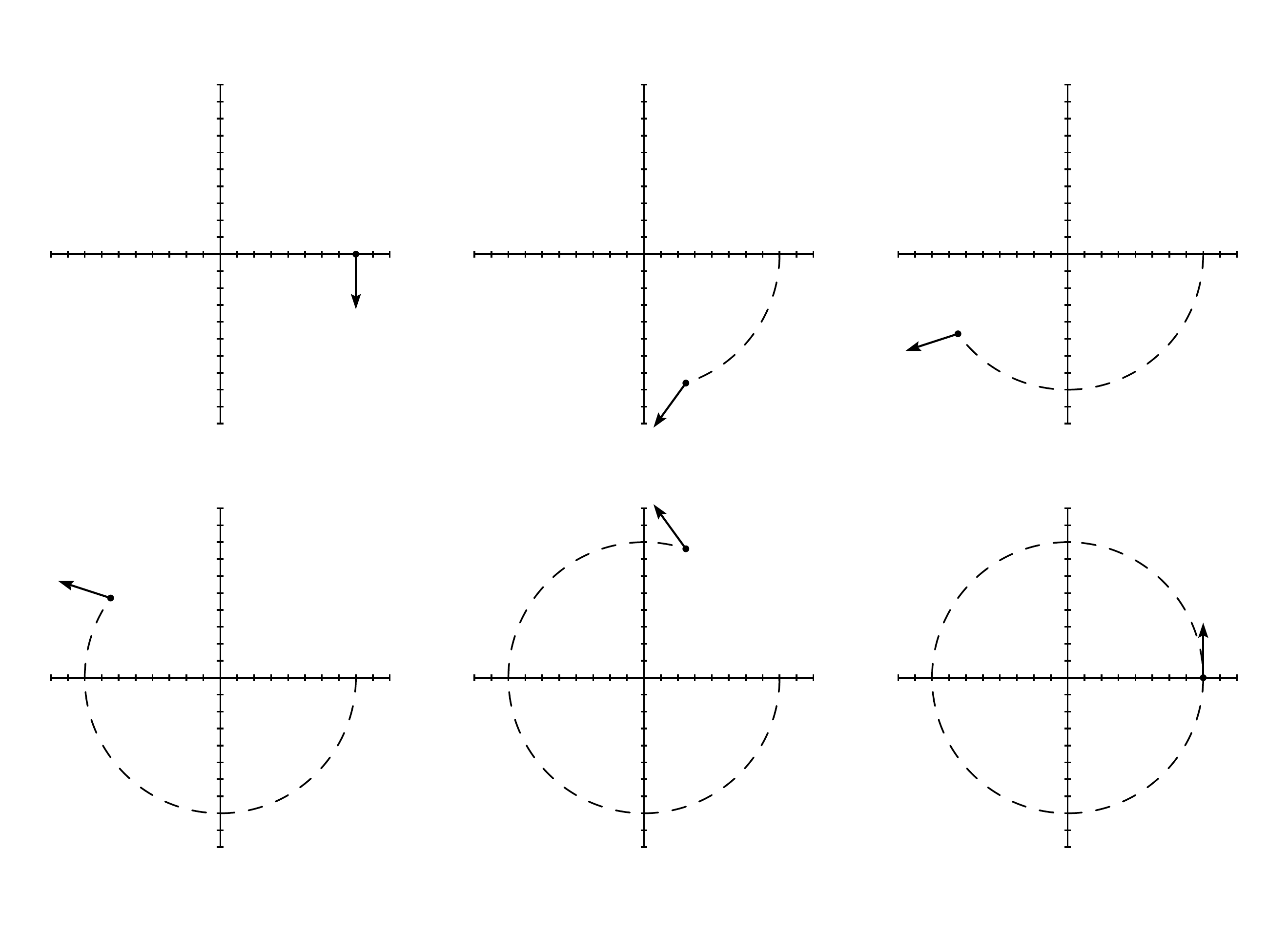}\\
{\small \quad \fig \label{figspindraai}
Een spinor draait mee met de ruimte, maar half zo snel.}
\end{tabular}
\end{center}

Als de ruimte 
$360^{\circ}$ om haar as draait, is er netto helemaal niets 
gebeurd. Als je immers een heel rondje draait, is daarna 
alles weer precies zoals het geweest is.
In figuur \ref{figspindraai} valt echter te zien dat
de spinor desondanks is omgeklapt, van $\downarrow$ naar $\uparrow$.
Dat is een beetje merkwaardig, en deze kleine
ongerijmdheid komt de arme spinorbundel dan ook op 
verlies van het predicaat `natuurlijk' te staan. 

Dit soort problemen ontstaat echter alleen bij \emph{grote} vervormingen
van de ruimte. Bij \emph{kleine} vervormingen
bewegen de spinoren w\'{e}l netjes mee.
We strijken daarom over ons hart en reiken de troostprijs uit:
spinorbundels mogen dan niet `natuurlijk' zijn, ze zijn ten minste
nog wel `infinitesimaal natuurlijk'.



Een spinorveld is een snede van de spinorbundel.
Het wijst in ieder punt van de ruimte $M$
\'e\'en spinor aan, de waarde van het veld in dat punt.
%
%
Spinorvelden worden veel gebruikt in de natuurkunde. 
Je zou zeggen dat de rare flip van 
$\downarrow$ naar $\uparrow$ onder een rotatie van $360^{\circ}$
problemen oplevert, maar dat blijkt mee te vallen.
De reden hiervoor is wat subtiel; het 
komt er uiteindelijk op neer dat een spinorveld nooit direct kan worden 
waargenomen, we zien alleen zijn kwadraten.
En in het kwadraat van een spinorveld heb je twee flips, 
die elkaar opheffen.


Het is ook zeker niet zo dat spinorvelden alleen gebruikt worden
om bizarre effecten te beschrijven in obscure hoekjes van de 
moderne fysica.
In het `standaard model', de natuurkundige theorie
die al sinds de jaren '70 vrijwel algemeen geaccepteerd wordt,
wordt alle huis-tuin-en-keukenmaterie beschreven door
spinorvelden.
Net als licht een trilling is in het elektromagnetisch veld,
zijn sommige
elementaire deeltjes 
trillingen in een spinorveld. 
Twee van dat soort deeltjes, de 
\emph{elektronen} en de \emph{quarks}, klonteren via een aantal
tussenstappen aaneen tot atomen.
Atomen op hun beurt zijn weer de bouwstenen van alle materie
om ons heen.
Welbeschouwd is dus iedere steen, iedere appel, 
iedere hond, mens of neushoorn een serie trillingen in \'e\'en 
groot spinorveld.

\subsection*{De Hamvraag}
We hebben gezien dat
verschillende soorten velden worden beschreven door
verschillende soorten vezelbundels. 
Een vezelbundel
komt echter alleen in aanmerking om een veld te beschrijven
als hij \emph{infinitesimaal natuurlijk} is.
Zijn vezels moeten dan als het ware meebewegen
met kleine vervormingen van de ruimte, 
zoals in figuur~\ref{figdraai}, \ref{figvelddraai} en \ref{figspindraai}.

Dit laat de natuur 
relatief weinig vezelbundels om uit te kiezen.
Als we ze allemaal kunnen beschrijven,
dan geeft ons dit
een 
idee welke velden w\'el in de natuur kunnen voorkomen, en welke niet.
De hamvraag is dus: 
\begin{center}
\emph{
`Welke vezelbundels zijn infinitesimaal
natuurlijk?'
}
\end{center}

Het is natuurlijk niet zo dat, als een veld in de natuur \emph{mag}
voorkomen, hij dat ook automatisch \emph{doet}.
Het infinitesimaal natuurlijk zijn van zijn vezelbundel moet je dan ook zien
als een minimumvereiste voor een veld om in de natuur voor te komen.
%
%

Hoofdstuk \ref{ch:BLID} is gewijd aan het 
beantwoorden van de hamvraag. 
Het belangrijkste resultaat is de classificatiestelling 
\ref{grondsmakreefknol}
op bladzijde \pageref{fruukkiop}.
Zij zegt dat alle infinitesimaal natuurlijke vezelbundels
afstammen van \'e\'en overzichtelijke oerfamilie, 
de zogenaamde `universele overdekkingen van bundels van
$k$-frames'.

Hoe die beesten er precies uitzien, is vrij goed bekend 
(zie sectie \ref{dekhetkframe} voor een algemene beschrijving en 
figuur \ref{kframe} voor een eenvoudig voorbeeld).
Van belang is vooral dat het er niet al te veel zijn, en dat hun
structuur sterk af\-hangt van de gladde vari\"{e}teit $M$.
Voor iedere $M$ geeft ons dit 
een `grijze' lijst van velden die infinitesimaal natuurlijk zijn
en dus in de natuur
kunnen voorkomen, en 
een `zwarte' lijst van velden die 
dat zeker niet kunnen.

De grap is nu dat we eigenlijk al wel een heel aardig beeld hebben van de
velden die in het wild voorkomen.
De laatste decennia hebben fysici (CERN, Fermilab etc.) 
enorme deeltjesversnellers gebouwd, 
waarmee een `witte' lijst is opgesteld van velden die echt gemeten zijn,
en die dus zeker in de natuur voorkomen.

Aan de andere kant hebben we eigenlijk geen flauw benul van de
vorm van onze 4-di\-men\-sio\-na\-le tijdruimte $M$.
Zij is op grote schaal recht voorzover wij kunnen zien, maar dat is 
helaas niet ver genoeg om uitsluitsel te kunnen geven.

\begin{center}
\begin{tabular}{p{11 cm}}
\hspace{3 cm} \includegraphics[width = 5 cm]{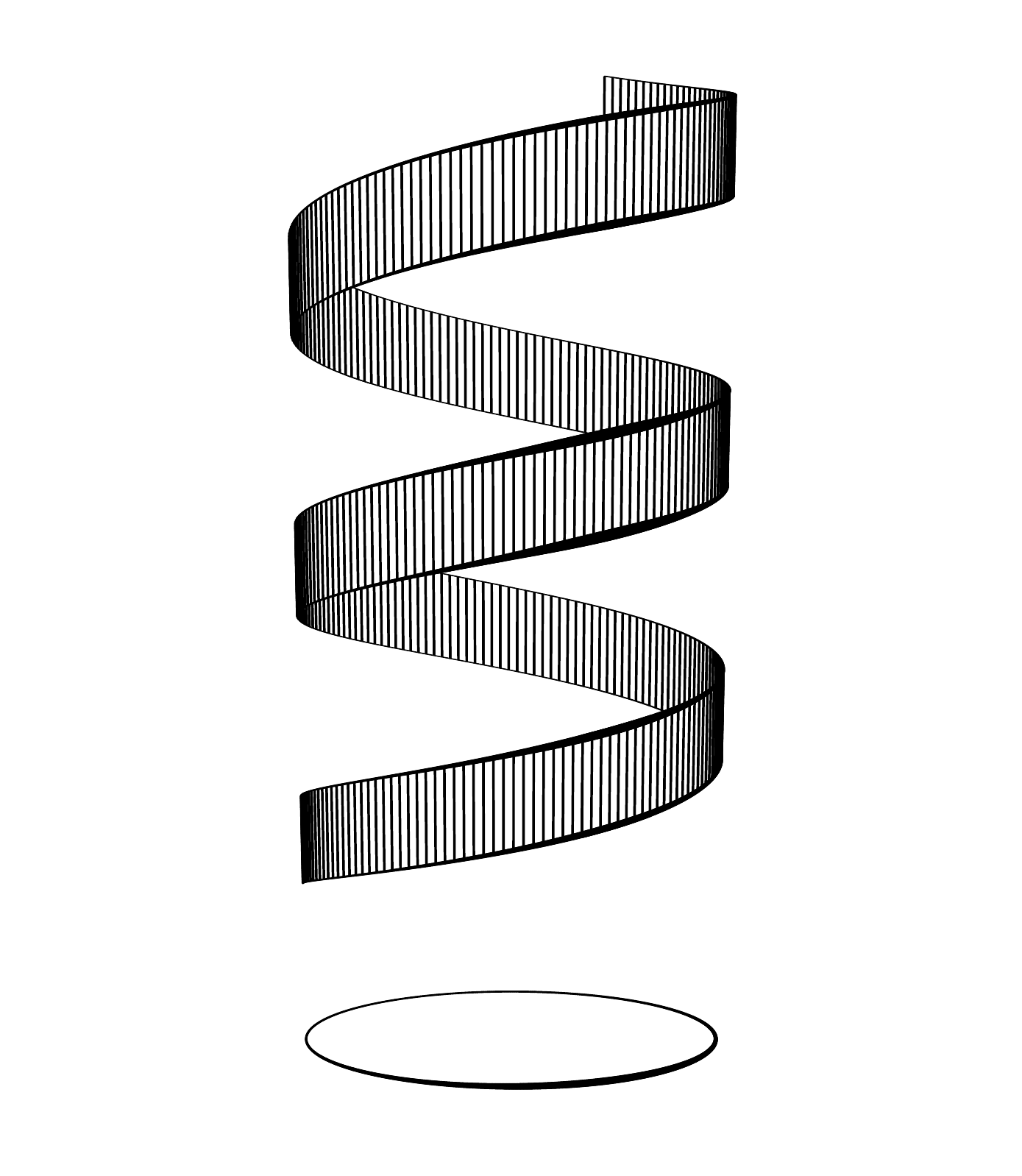}\\
{\small \quad \fig \label{kframe}
Schets van de universele overdekking van de bundel van $k$-frames
over de cirkel, voor het geval $k=1$.}
\end{tabular}
\end{center}

We willen dus eigenlijk stelling \ref{grondsmakreefknol} ondersteboven houden!
In plaats van alle mogelijke velden te bepalen voor een gegeven 
tijdruimte $M$, bepalen we nu alle mogelijke tijdruimten,
gegeven de `witte' lijst van velden waarvan we al weten dat ze voorkomen. 
Dit is de portee van stelling
\ref{naturalspin} op bladzijde
\pageref{naturalspinpage}.

Een aantal mogelijkheden kun je nu al direct afstrepen.
De gladde 4-di\-men\-sio\-na\-le vari\"eteit $M$ met de naam 
$\mathbb{C}P^{2}$
valt bijvoorbeeld af, omdat zijn `zwarte lijst' volgens 
stelling \ref{grondsmakreefknol} overlapt met de 
`witte lijst' van waargenomen velden. 
Als het beeld dat we hebben van de elementaire velden
correct en compleet is, kan de tijdruimte $M$ dus nooit van deze 
vorm zijn.

Welbeschouwd is dit een wonderlijke conclusie:
door goed te kijken naar de kleinst mogelijke structuren, 
de elementaire deeltjes waaruit alle materie is opgebouwd, 
krijgen we informatie over de vorm van ons heelal.

}
\selectlanguage{english}


\chapter*{Dankwoord}
\phantomsection
\addcontentsline{toc}{chapter}{Dankwoord}
{\frenchspacing
\selectlanguage{dutch}
Voor de inhoud van dit proefschrift draagt de auteur de ongedeelde 
ver\-ant\-woor\-de\-lijk\-heid, maar zeker niet de ongedeelde eer. 
Graag wil ik hier iedereen van harte danken die er een bijdrage aan heeft geleverd.

In de eerste plaats zijn dat natuurlijk mijn copromotoren, Hans Maassen en 
Johan van de Leur,
die mij vertrouwd hebben gemaakt met respectievelijk de kwantumkansrekening 
en de Lie-theorie. De vele plezierige discussies met Hans tijdens het
laatste jaar van mijn studie in Nijmegen hebben de basis gevormd
voor hoofdstuk \ref{ch:ITISC} en \ref{ch:UDHP}.
De eerste daarvan is dan ook
gepubliceerd onder ons beider naam. 
Johan heeft mij begeleid bij mijn promotietraject in Utrecht en ik wil hem 
graag bedanken voor alle hulp en raad die ik van hem heb gehad, 
al dan niet wiskundig van aard.
Voor zijn enthousiaste en essenti\"ele hulp bij hoofdstuk \ref{ch:BLID}
ben ik ook Marius Crainic veel dank verschuldigd.

Verder wil ik graag mijn coauteurs bedanken voor de fijne en leerzame
samenwerking. Hoofdstuk \ref{ch:OPJM} is geschreven samen met
Luc Bouten,
en hoofdstuk \ref{ch:OEQS} samen
met M\u{a}d\u{a}lin Gu\c{t}\u{a} en Jonas Kahn. 
Graag dank ik ook de res\-pec\-tie\-ve\-lij\-ke thuisbases van Luc en M\u{a}d\u{a}lin
voor hun gastvrijheid: 
het Mabuchilab, toen nog aan het California Institute
of Technology, en de School of Mathematical Sciences aan de 
University of Nottingham.

Hartelijk dank ook aan mijn 
promotor, Roberto Fern\'andez, en
aan de leescommissie:
Erik van den Ban,
Nilanjana Datta,
Mark Fannes,
Eduard Looijenga en
Christoph Wockel.

Ten slotte wil ik graag mijn medewiskundigen bedanken voor alle interessante
gesprekken, waar ik veel van geleerd heb. 
Bij naam wil van hen slechts diegenen noemen die ik helaas niet meer persoonlijk
bedanken kan: Hanneke Janssen en Hans Duistermaat.

} 
\selectlanguage{english}

\chapter*{Curriculum Vitae}
\phantomsection
\addcontentsline{toc}{chapter}{Curriculum Vitae}

\selectlanguage{dutch}
{\frenchspacing

Bas Janssens werd op 20 mei 1981 geboren te Maastricht en behaalde in 1999
cum laude zijn gymnasiumdiploma aan het Stella Maris College in Meerssen.
Hij studeerde natuurkunde (1999-2004) en wiskunde (2000-2005) aan de 
Katholieke Universiteit Nijmegen (nu Radboud Universiteit Nijmegen), 
waar hij
onder begeleiding van Hans Maassen 
cum laude afstudeerde in de richting van de kwantumkansrekening.
Vervolgens bracht hij korte tijd door aan het California Institute 
of Technology (waar hij samenwerkte met Luc Bouten) en
aan de University of Nottingham (waar hij samenwerkte met 
Jonas Kahn en M\u{a}d\u{a}lin Gu\c{t}\u{a}). 
In 2006 begon hij zijn promotieonderzoek in de richting van 
oneindigdimensionale 
Lie-algebras aan de Universiteit Utrecht,
onder begeleiding van Johan van de Leur. 
Dit proefschrift is niet enkel een verslag van dit promotieonderzoek,
maar van al zijn bevindingen tot nu toe.
  
}
\selectlanguage{english}


 \newpage
 \thispagestyle{empty} 
 \mbox{}
 \newpage
 \thispagestyle{empty} 
 \mbox{}
 
\end{document}